# A class of finite dimensional spaces and $H(\mathrm{div})$ – conformal elements on general polygons

*Rémi Abgrall* *, *Élise Le Mélédo*†, *Philipp Öffner*‡

*Institute of Mathematics*

*Universität Zürich, Winthherturerstrasse 190, Zürich, Switzerland*

#### Abstract

We present a class of discretisation spaces and $H(\mathrm{div})$ – conformal elements that can be built on any polytope. Bridging the flexibility of the Virtual Element spaces towards the element's shape with the divergence properties of the Raviart – Thomas elements on the boundaries, the designed frameworks offer a wide range of $H(\mathrm{div})$ – conformal discretisations. As those elements are set up through degrees of freedom, their definitions are easily amenable to the properties the approximated quantities are wished to fulfil. Furthermore, we show that one straightforward restriction of this general setting share its properties with the classical Raviart – Thomas elements at each interface, for any order and any polytopial shape. Then, we investigate the shape of the basis functions corresponding to particular elements in the two dimensional case.

## 1 Introduction

Several classes of schemes adapted to hyperbolic problems, like Discontinuous Galerkin [6], Flux Reconstruction [12, 21, 15] or Residual Distribution use the Finite Elements framework. The approximated solution is considered piecewise and is split onto a tessellation of the computational domain, using usually simplicial or quadrangular tiles. A finite space of approximation is then set up on every tile, offering a discretisation framework from which schemes can be built.

*remi.abgrall@math.uzh.ch
†elise.lemeledo@math.uzh.ch
‡philipp.oeffner@math.uzh.ch

In order to provide an adequate spatial approximation of complex shape domains, extensive efforts have been made on designing strategies allowing various element shapes [18]. Also, some variational formulations, spaces included, have been redesigned to handle non simplicial elements, as by example [20]. Our contribution is along similar lines, going beyond simplicial elements by providing a framework for Flux Reconstruction methods on arbitrary polytopes.

The classical Flux Reconstruction starts from a point-wise approximation of the flux and modify it in order to take care of local conservation. This amounts to introduce a modification of the centred approximation. When dealing with triangles, quad, hex or tets, this correction can be written as adding an element of the Raviart – Thomas [16] space associated to this element. Going further and using the Residual Distribution framework, we have shown in [1] how to construct schemes similar to those occurring in the Flux Reconstruction framework for arbitrary polytopes, convex or non-convex. This method enjoys a non-linear entropy stability property, but this is only seen algebraically. It is possible to recover a clean variational setting if one is able to construct space sharing the properties of the Raviart-Thomas setting on the boundaries of arbitrary polygons. This is the purpose of this report.

The theory of $H(\text{div})$ – conformal element has already been studied by Raviart, Thomas [16] and later generalized by Nédélec [13] or Brezzi, Douglas and Marini in the context of mixed finite element method [3]. More recently, a mixed Petrov-Galerkin scheme using Raviart – Thomas elements has also been investigated in [8]. However, – up to the authors knowledge – those elements are limited to simplicial and quadrangular shapes.

Several attempts to use general polygons have been made [10, 17], but they usually make use of generalized barycentric coordinates and are delicate to handle in distorted non-convex elements. A first polygonal $H(\text{div})$ – conformal element has been proposed in [7] using gradient reconstruction and pyramidal submeshes tessellation. However, their construction requires some shape regularity within the mesh, and the parallel with Raviart – Thomas spaces is limited to the lowest order space. Some other approaches as Virtual Elements Method [20] introduced approximation spaces based on Poisson's solutions. Although more flexible towards the element's shape, those are scalar and not H(div) – conformal.

In this report, we propose a construction that inherits the interface properties of the Raviart – Thomas elements and benefits from the shape flexi-bility of the Virtual Element discretisation. Moreover, more than defining a basis on which the correction functions can be decomposed, this new setting offers a new element class that can be used as such in the construction of further numerical schemes.

After briefly recalling classical algebra results in the *Section 2*, we detail the simplicial and quadrilateral the Raviart – Thomas elements respectively in the *Sections 3 and 4*. Then, we introduce our new class of discretisation spaces in the *Section 5* before providing two dimensional elements examples, along with a method to build the corresponding basis functions in the *Section 6*. Lastly, we test the discussed $H(\mathrm{div})$ – conformal elements through the behaviour of their corresponding basis functions in the numerical results given in the *Section 7*. For interested readers, an application of those elements can be found in [1].

# Contents

## Notations

We present here for reference the main notations used in this report.

### Geometrical notations

$d$ dimension of the coordinate space
$K$ polytope on which elements will be built
$n$ number of edges of the polytope
$f$ generic face / edge of an element $K$
$f_i$ face $i$ of the element $K$
$n$ generic normal to the face $f$ of the element $K$
$n_i = (n_{ix}, n_{iy})$ normal to the face $f_i$ of the element $K$
$\mathrm{x} = (x, y)$ coordinate vector in the two dimensional case
$x = (x_1, \cdots, x_d)$ coordinate vector in the $d$ - dimensional case
$\gamma$ path of an hyper dimensional face within the coordinate space

### Functional space notations

$\mathbb{P}_k(K)$ polynomial space $\mathbb{P}_k$ built on $\mathbb{R}^d$
$\mathbb{P}_k(f)$ polynomial space $\mathbb{P}_k$ built on $\mathbb{R}^{d-1}$
$\mathbb{Q}_k(K)$ polynomial space $\mathbb{Q}_k$ built on $\mathbb{R}^d$
$\mathbb{Q}_{[k]}(K)$ polynomial space $\mathbb{Q}$ of polynomials having degree exactly k
$\mathbb{Q}_k(f)$ polynomial space $\mathbb{Q}_k$ built on $\mathbb{R}^{d-1}$
$RT_k(K)$ Raviart – Thomas space
$f = e\left\{ \quad (\cdot) \quad \right\}$ vector formed by each coordinate of the vector's $e$ multiplied by the scalar solution of the problem $(\cdot)$.

### Operator notations

$\zeta_i$, $\xi_i$ permutation operators
$\circ$ Hadamard product
$\times$ usual scalar product when used between two numbers
$\bigtimes$ Cartesian product
$\mathbb{1}$ indicator function

### Lexicography

Free set: linearly independent set
Poisson function: solution to a Poisson's problem

# 2 Some algebraic properties

To begin with, we recall some well-known algebraic results that will be in use in the following sections.

## 2.1 Multivariate polynomial spaces

Approximations and interpolations of solutions living in some domain $K \subset \mathbb{R}^d$ are usually based on polynomials as finite dimensional polynomial spaces of approximation are easy to handle.

Classically being restrictions of the space $\mathbb{P}(K)$ of all $d$-variate polynomials, their characterisation is usually based on the polynomials' degree whose definition is not unique. Therefore, we enjoy several straightforward definitions of finite dimensional polynomial spaces. We detail here three of them.

***Note.*** Any polynomial $p$ can be written as a finite sum of monomial elements $p_l$, $l \in [\![1, m]\!]$, $m$ being any natural number. The polynomial degree is then always defined by $\deg(p) = \max\left(\{\deg(p_l)\}_{l\in[\![1,m]\!]}\right)$. Therefore, we only discuss without loss of generality the definition of the degree on monomials, transferring immediately to polynomials when defining the subsequent polynomial spaces. ▲

Let us set some monomial $p = \prod_{i=1}^{d} x_i^{\alpha_i}$ belonging to the space of monomials $\mathcal{M}(K) \subset \mathbb{P}(K)$. A natural definition of the degree arises when one considers the monomial as a single entity and interprets the degree as a mixed growing rate. We then have

$$\deg_1(p) = \sum_{i=1}^{d}\{\alpha_i\},$$

leading to the space

$$\mathbb{P}_k(K) = \{p \in \mathbb{P}(K),\, \deg_1(p) \leq k\}.$$

Using algebraic results [14] one can derive its dimension, being

$$\dim \mathbb{P}_k(K) = \binom{k+d}{k}.$$

An other view would rather be to define the degree as

$$\deg_2(p) = \max\left(\{\alpha_i\}_{i\in[\![1,\, d]\!]}\right).$$

There, it would represent the maximal unidirectional growth rate of the monomial regardless the variable it applies to. The subsequent finite dimensional space would read:

$$\mathbb{Q}_k(K) = \{p \in \mathbb{P}(K),\, \deg_2(p) \leq k\},$$

and would enjoy the dimension

$$\dim \mathbb{Q}_k(K) = (k+1)^d.$$

***Example.*** Let us focus on the two dimensional case and consider the monomial $(x,\, y) \mapsto x^3 y$. There, when seeing the monomial as a single entity, the two coordinates contribute conjointly to the behaviour. This leads to a combined growth rate of order four ($\deg_1(x^3y) = 4$).

Considering however the monomial's behaviour as being a collection of unidirectional growths as if we were uncoupling the domain, the growing rates of $x$ and $y$ would be compared with each other to determine the coordinate that supports the dominant variations. The steepest rate achieved in any direction is then seen as the order, being here three ($\deg_2(x^3y) = 3$). ♦

In the same spirit, one can prefer to prescribe the behaviour coordinate-wise and set a variant of $\mathbb{Q}_k(K)$ by imposing a maximal power on each variable. It then reads

$$\mathbb{P}_{k_1,\cdots,k_d}(K) = \{p \in \mathbb{P}(K),\, p = \sum_{l=1}^{m}\prod_{i=1}^{d} x_i^{\alpha_{il}} \text{ s.t. } \forall i \in [\![1,\, d]\!],\ \max_{l\in[\![1,\, m]\!]}\{\alpha_{il}\} \leq k_i\}.$$

Note that in the above definition, any $p \in \mathbb{P}(K)$ is decomposed using $m > 0$ monomials. Furthermore, we directly get from algebra that

$$\dim \mathbb{P}_{k_1,\cdots,k_d} = \prod_{l=1}^{d}(k_l + 1).$$

This space is suitable for discretising an anisotropic problem, or when the quality of the discretisation projections is wished to be detailed variable-wise.

***Example.*** The projection quality being detailed variable-wise, the space $\mathbb{P}_{k+1,\,k}$ see its functions robust towards the $\partial_{x_1}$ operator.

Indeed, the partial derivative in $x_1$ maps $\mathbb{P}_{k+1,\,k}$ onto the space $\mathbb{Q}_k$, space that provides an homogeneous discretisation framework in the two variables $x_1$ and $x_2$. This will be useful when turning to divergence considerations. ⧫

***Remark.*** One can derive inclusion relations between those spaces. By example, it holds:

$$\mathbb{P}_{k_1,\,\cdots,\,k_d} \subset \mathbb{Q}_{\max\limits_{i\in[\![1,\,d]\!]} k_i}$$
$$\mathbb{Q}_k = \mathbb{P}_{k,\,\cdots,\,k}$$

▲

***Remark.*** In practice, when working with a polytopial shape $K$, one generally uses the $\mathbb{P}_k(K)$ spaces on simplicial shapes and the $\mathbb{Q}_k(K)$ spaces on higher polytopes. ▲

Some spaces based on polynomial approximation can also be designed specifically for the boundary $\partial K$ in order to allow discontinuities from one face to an other. In particular, the two spaces that will be used later on are described by

$$\mathcal{R}_k(\partial K) = \{p \in L^2(K),\, p|_{f_i} \in \mathbb{P}_k(f_i) \text{ for every face } f_i \in \partial K\}$$

and

$$\mathcal{T}_k(\partial K) = \{p \in L^2(K),\, p|_{f_i} \in \mathbb{Q}_k(f_i) \text{ for every face } f_i \in \partial K\}.$$

As no continuity is required at the polytope's vertices, each of those two spaces can be reduced to a collection of $n := |\{f_i\}_i|$ independent polynomials of order $k$, each of them living on one specific face $f_i \in \partial K$. Thus, assuming that none of the polygons' faces is degenerated, their dimension automatically reads

$$\dim \mathcal{R}_k(\partial K) = n \dim \mathbb{P}_k(f) \text{ for any generic face } f \in \partial K \tag{1}$$
$$\dim \mathcal{T}_k(\partial K) = n \dim \mathbb{Q}_k(f) \text{ for any generic face } f \in \partial K \tag{2}$$

***Note.*** In the two dimensional case, the spaces $\mathcal{R}_k(\partial K)$ and $\mathcal{T}_k(\partial K)$ coincide. Furthermore, it comes $\dim \mathcal{R}_k(\partial K) = \dim \mathcal{T}_k(\partial K) = |\{f_i\}_i|(k+1)$. ▲

## 2.2 Symmetric polynomial spaces

When the considered object is multi-dimensional, the discretisation space is commonly formed from a Cartesian product of some polynomial spaces introduced in the previous section. In particular, as the discretised quantity may have distinct behaviours in each coordinate, one may choose to work with a different type of scalar polynomial space per coordinate.

***Example.*** If one would like to discretise the quantity $u = (u_1, u_2)$ living on a subset of $\mathbb{R}^2$ where it is known that the behaviour of $u_1$ is comparable to $(x, y) \mapsto xy^6$ and the one of $u_2$ to $(x, y) \mapsto x^3y$, then the polynomial discretisation space may be set as $\mathbb{P}_{1,6} \times \mathbb{P}_{3,1}$. ♦

However, in order to have an homogeneous representation of each variable it may be interesting to restrict the choice of those Cartesian polynomial spaces to those that provide a similar discretisation setting in every coordinate. There, the similarity would be defined through circular permutations among the coordinates.

***Example.*** If one would like to discretise the quantity $u = (u_1, u_2)$ living on a subset of $\mathbb{R}^2$ where it is known that the behaviour of $u^1$ is comparable to $(x, y) \mapsto xy^6$ and the one of $u_2$ to $(x, y) \mapsto x^3y$, then the polynomial discretisation space may be set as $\mathbb{P}_{3,6} \times \mathbb{P}_{6,3}$ or as $\mathbb{P}_{1,3} \times \mathbb{P}_{3,1}$, depending on the quality one would like to achieve. ♦

If the object of the discretisation plays the role of a flux or appears under some differential operator, one can also take into account the impact of the operator onto the variables in order to retrieve an homogeneous discretisation after the operator has been applied.

***Example.*** If one would like to approximate the variable of interest $u = (u_1, u_2)$ with the help of polynomials of degree $k$, but that the considered equation is of the type $\nabla_x \cdot u = f$ for any function $f$ regular enough, then one may use a polynomial space of the type $\mathbb{P}_{k+1,k} \times \mathbb{P}_{k,k+1}$, as presented in the example of the last paragraph. ♦

The above presented restricted discretisations make use of Cartesian spaces which are so-called symmetric polynomial spaces.

**Definition 2.1** Symmetric polynomial spaces

We call symmetric polynomial space a vectorial polynomial space in which every coordinate is approximated in the same manner. More precisely, a vectorial space is called symmetric if it can be recast into the following shape.

Let $\{b_j\}_j$ be a basis of some scalar polynomial space. Then, a space generated by the set

$$\bigcup_{i=1}^{d}\{(0,\,\cdots,\,0,\,b_j(\zeta_i(x_i)),\,0,\,\cdots,\,0)^T\}_j$$

where the term $b_j(\zeta_i(x_i))$ lies in the $i^{\text{th}}$ position and where $\zeta_i$ represents the circular permutation $(x_1,\,\cdots,\,x_d)\mapsto(x_i,\,\cdots,\,x_d,\,x_1,\,\cdots)$ is symmetric.

***Example.*** *Two and three dimensional cases* Let us consider the canonical basis of the polynomial spaces $\mathbb{P}_{k_1,k_2}\times\mathbb{P}_{k_3,k_4}$ for any non-negative integers $k_1$, $k_2$, $k_3$ and $k_4$. Then, the following assertions hold.

- $\mathbb{P}_{2,0}\times\mathbb{P}_{0,2}$ is a symmetric space. Indeed, it is generated by

$$\left\{\begin{pmatrix}1\\0\end{pmatrix},\begin{pmatrix}x\\0\end{pmatrix},\begin{pmatrix}x^2\\0\end{pmatrix},\begin{pmatrix}0\\1\end{pmatrix},\begin{pmatrix}0\\y\end{pmatrix},\begin{pmatrix}0\\y^2\end{pmatrix}\right\},$$

where the permutation is clearly identifiable.

- $\mathbb{P}_{2,0}\times\mathbb{P}_{0,1}$ is not a symmetric space. Indeed, it is generated by

$$\left\{\begin{pmatrix}1\\0\end{pmatrix},\begin{pmatrix}x\\0\end{pmatrix},\begin{pmatrix}x^2\\0\end{pmatrix},\begin{pmatrix}0\\1\end{pmatrix},\begin{pmatrix}0\\y\end{pmatrix}\right\}.$$

where the used transformation $\zeta$ is not even a permutation as it omits the term $(0,\,y^2)^T$.

- The space

$$\left(\mathbb{P}_{2,0}\times\mathbb{P}_{0,2}\setminus\left\{\begin{pmatrix}1\\0\end{pmatrix},\begin{pmatrix}0\\1\end{pmatrix}\right\}\right)\bigcup\left\{\begin{pmatrix}1\\1\end{pmatrix}\right\}$$

is symmetric. Indeed, it is generated by

$$\left\{\begin{pmatrix}1\\1\end{pmatrix},\begin{pmatrix}x\\0\end{pmatrix},\begin{pmatrix}x^2\\0\end{pmatrix},\begin{pmatrix}0\\y\end{pmatrix},\begin{pmatrix}0\\y^2\end{pmatrix}\right\}.$$

There, a circular permutation between the first rows components and the second ones is identifiable. Indeed, by the definition of the circular permutation we have $\zeta((1,\,1)^T) = (1,\,1)$. Thus, when defining the discretisation of the second coordinate through the circular permutation of the three first vectors corresponding to the discretisation of the first coordinate, the obtained set reduces to

$$\left\{\begin{pmatrix}1\\1\end{pmatrix},\begin{pmatrix}1\\1\end{pmatrix},\begin{pmatrix}x\\0\end{pmatrix},\begin{pmatrix}x^2\\0\end{pmatrix},\begin{pmatrix}0\\y\end{pmatrix},\begin{pmatrix}0\\y^2\end{pmatrix}\right\}.$$

And as the definition of symmetric spaces we consider asks only to have a generating set, not necessarily free, this last space is compliant with the definition. Thus, it is a symmetric space.

• In three dimensions, the space $\mathbb{P}_{0,\,2,\,1}\times\mathbb{P}_{2,1,0}\times\mathbb{P}_{1,\,2,\,0}$ is symmetric. However, for the same reason as developed in the second case of this example, the space $\mathbb{P}_{0,\,2,\,1}\times\mathbb{P}_{2,0,1}\times\mathbb{P}_{1,\,2,\,0}$ is not. Its slightly modified version

$$\{\mathbb{P}_{0,\,2,\,1}\times\mathbb{P}_{2,0,1}\times\mathbb{P}_{1,\,2,\,0}\}\setminus\left\{\begin{pmatrix}z^2\\0\\0\end{pmatrix},\begin{pmatrix}0\\x^2\\0\end{pmatrix},\begin{pmatrix}0\\0\\y^2\end{pmatrix}\right\}$$

would however be symmetric. ♦

## 2.3 Degrees of freedom

**Definition of degrees of freedom** Let us consider some finite dimensional space $P$. Seeing it as a discretisation space, we can approximate any quantity of interest by projecting it onto $P$ and decompose it on some basis.

***Example.*** Let us consider the finite dimensional space $\mathbb{P}_2(\mathbb{R})$. There, any projected quantity will be polynomial of degree at maximum two and can be decomposed on any basis of $\mathbb{P}_2$. By example, the polynomial $p = 2x^2+3x+2$ can be decomposed on the canonical basis $\{x^2,\, x,\, 1\}$ and represented through the corresponding values 2, 3 and 2. ♦

One can see the coefficients of the decomposition as liberties in the definition of the function lying in the space $P$. Those liberties are expressed in terms of the behaviours represented by each of the basis functions.

Another way to represent the quantity of interest is to fix those liberties not by projecting on the basis functions, but by specifying values of some quantifiers that are characterising the functions living in the discretisation space.

***Example.*** Knowing that it belongs to $\mathbb{P}_2(\mathbb{R})$, the polynomial $p\colon x \mapsto 2x^2 + 3x + 2$ can be completely determined by either of the following relations.

$$(a)\begin{cases} p(0) = 2 \\ p(1) = 7 \\ p(-1) = 1 \end{cases} \qquad (b)\begin{cases} p(0) = 2 \\ p(1) = 7 \\ \displaystyle\int_0^1 x\,p\,\mathrm{d}x = 3 \end{cases} \qquad (c)\begin{cases} \displaystyle\int_0^1 p\,\mathrm{d}x = \frac{25}{6} \\ \displaystyle\int_0^1 x\,p\,\mathrm{d}x = \frac{5}{2} \\ \displaystyle\int_0^1 x^2\,p\,\mathrm{d}x = \frac{115}{60} \end{cases}$$

The relations $(a)$ are characterising the polynomial through so-called *point values*, while the configuration $(c)$ involves only *moments*. One can also mix the characterisation types and retrieve a *mixed* configuration as shown in the configuration $(b)$.

To see this characterisation, we can derive the linear system corresponding to a generic function of $\mathbb{P}_2(\mathbb{R})$, quantified by any set of degrees of freedom $(a)$, $(b)$ or $(c)$. Indeed, given any basis $\mathcal{B} = \{b_1, b_2, b_3\}$ of $\mathbb{P}_2(\mathbb{R})$, each of the relations $(a)$, $(b)$ and $(c)$ provides a determination of the three coefficients associated to the basis $\mathcal{B}$.

By example, we can look for a function $p$ of $\mathbb{P}_2(\mathbb{R})$ under the form $p\colon x \mapsto a\,x^2 + b\,x + c$, using the canonical basis for $\mathcal{B}$ and real coefficients $a$, $b$ and $c$. Seeing the configuration $(b)$ as a quantifier set, we can derive the following determination system.

$$\begin{cases} a \times 0^2 + b \times 0 + c = 2 \\ a \times 1^2 + b \times 1 + c = 7 \\ \displaystyle\int_0^1 x\,(a\,x^2 + b\,x + c)\,\mathrm{d}x = 3 \end{cases}$$

By linearity of the integral, it reduces to the linear system

$$\begin{pmatrix} 0 & 0 & 1 \\ 1 & 1 & 1 \\ \frac{1}{3} & \frac{1}{2} & 1 \end{pmatrix} \begin{pmatrix} a \\ b \\ c \end{pmatrix} = \begin{pmatrix} 2 \\ 7 \\ 3 \end{pmatrix},$$

and the coefficients $a = 2$, $b = 3$ and $c = 2$ follow. Note that at this point there is no particular relation between the basis $\mathcal{B}$ and the used linear forms.

♦

The shape of the quantifiers is up to one taste, but they should be linear forms over the coordinate space. Furthermore, to entirely determine any function of the discretisation space, one needs as many quantifiers as the dimension of the space, each of them emphasising the unique types of behaviours that can be found within $P$. One can then determine any function of $P$ through the values taken by those specific quantities. Those quantifiers are called degrees of freedom and live in the space that is dual to $P$.

**Definition 2.2** Degrees of freedom

A degree of freedom is any real valued linear form over the coordinate space that acts on $P$ and characterises some property of the functions living in $P$.

In other words, its value settles one liberty in the definition of the functions living in $P$.

Let us now consider as many degrees of freedom as the dimension of the space $P$, and explore the conditions ensuring that the obtained set entirely determines any function of the space. The condition that each quantifier characterises a unique type of behaviour translates into the linear independence of the equations retrieved when the degrees of freedom are applied to any function of $P$. This property is called *unisolvence*.

**Definition 2.3** Unisolvence

When a set of degrees of freedom values completely determines any element living in $P$, it is said to be *unisolvent* for the space $P$.

***Example.*** *Continued* The rows of the above matrix are linearly independent. Thus, as the quantifiers act on $P$ and as we have as many rows as the dimension of the space $P$, the matrix is fully ranked, and invertible. The set of degrees of freedom $(b)$ is then unisolvent for the space $\mathbb{P}_2(\mathbb{R})$. ⧫

The definition of the degrees of freedom will emphasise some properties or behaviour of functions one may wish to analyse. In particular, for any set of degrees of freedom there exists a basis of $P$ that is dual to its own. In other words, there exists a couple (degrees of freedom – dual basis functions): $(\{\sigma_i\}_i, \{\varphi_j\}_j)$, $(i, j) \in [\![1, \dim P]\!]^2$ such that

$$\sigma_i(\varphi_j) = \delta_{ij}.$$

This dual basis represents all the behaviours the functions of $P$ are decomposed on, while the degrees of freedom gives a direct access the corresponding quantified properties.

***Example.*** The following degrees of freedom are dual to the canonical basis of $\mathbb{P}_2(\mathbb{R})$.

$$\begin{cases} \sigma_1 \colon p \mapsto p(0) \\ \sigma_2 \colon p \mapsto \dfrac{p(2) - p(1)}{3} \\ \sigma_3 \colon p \mapsto p(1) - p(0) - \dfrac{p(2) - p(1)}{3} \end{cases}$$

One can notice that computing the value of the degrees of freedom for a specific function furnishes the coefficients of its projection onto the dual basis. ⧫

Though this dual basis exists, its knowledge is not necessary to express the discretised quantities. Indeed, as shown above the quantifiers' values are enough to characterise any function living in $P$ in their own. Thus, we only have to know the value of the degrees of freedom to represent a function in $P$, provided that the set of degrees of freedom is unisolvent for $P$.

***Example.*** As the matrix of linear forms applied to any basis of $\mathbb{P}_2(\mathbb{R})$ is invertible, we can characterise any function of $\mathbb{P}_2(\mathbb{R})$ only by looking at the values of its second member.

In the above example, the polynomial $p\colon x \mapsto 2\,x^2+3\,x+2$ is characterised by the values $(2,\ 7,\ 3)^T$ through the set of degrees of freedom $(b)$, independently from the basis functions used to represent the polynomials. Indeed, changing the basis of $\mathbb{P}_2(\mathbb{R})$ would only transfer the matrix to an other base and impact the values of $a$, $b$, $c$, which leaves the second member unchanged. Therefore, the representation of the quantity through the degrees of freedom values is not impacted. ⧫

It is possible to have several representations of the same quantity. Indeed, as it is possible to have several basis for the space $P$, it is equally possible to define several sets of set of degrees of freedom. As those two sets will be dual to two different basis of $P$, we can link them through a transfer matrix. We detail the procedure in the *Section 2.4*.

In practice, the two common types of degrees of freedom are so-called *point-wise values* and *moments*. The formers consists in evaluating the discretised quantities at a precise point, while the latter is averaging the quantity tested against some function characterizing a target behaviour on a subset of the coordinate space (as used *e.g.* in the virtual elements method). In particular, when using moment based degrees of freedom one can test functions living in $P$ against some subspace contained in $P$. However, testing against a full subspace is computationally unrealistic. Therefore, practically one only computes the degrees of freedom from a basis of the subspace, as shown below.

**Decomposition of moment based degrees of freedom on the basis of their kernel** Let us consider in this paragraph a set of moment based degrees of freedom under the form

$$\begin{aligned} \sigma\colon Q &\longrightarrow \mathbb{R} \\ q &\longmapsto \int\limits_K q \cdot p \,\mathrm{d}x, \quad \forall p \in \mathcal{P} \end{aligned} \tag{3}$$

for some open domain $K \subset \mathbb{R}^d$, $d \in \mathbb{N}^*$ and some finite dimensional variational spaces $Q$ and $\mathcal{P}$. There, instead of computing the expression of the degrees of freedom for any function of $\mathcal{P}$, it is enough to generate the set of internal degrees of freedom from any basis of $\mathcal{P}$. Indeed, any $p \in \mathcal{P}$ can be decomposed on any basis of $\mathcal{P}$. Therefore, any degree of freedom can be decomposed on the set of degrees of freedom obtained from any basis of $\mathcal{P}$ thanks to their linearity.

More precisely, let $\{b_i\}_{i\in[\![1,\,\dim\mathcal{P}]\!]}$ be some basis for $\mathcal{P}$ and $\{\sigma_i\}_{i\in[\![1,\,\dim\mathcal{P}]\!]}$ their corresponding degree of freedom by duality. We have by definition; for all $i \in [\![1,\,\dim\mathcal{P}]\!]$:

$$\sigma_i\colon q \mapsto \int\limits_K q \cdot b_i \,\mathrm{d}x.$$

Let us now take any $p \in \mathcal{P}$. Then, by definition of the basis, we have that there exists a set of constants $\{\alpha_i\}_{i\in[\![1,\,\dim\mathcal{P}]\!]}$ such that $p = \sum_{i=1}^{\dim\mathcal{P}} \alpha_i b_i$. Its corresponding degree of freedom then writes

$$\sigma\colon q \mapsto \int\limits_K q \cdot p \,\mathrm{d}x = \int\limits_K p \cdot \left(\sum_{i=1}^{\dim\mathcal{P}} \alpha_i b_i\right) \mathrm{d}x.$$

By the linearity of the dot product and of the degrees of freedom, we have:

$$\sigma\colon q \mapsto \sum_{i=1}^{\dim \mathcal{P}} \alpha_i \int_K q \cdot b_i \,\mathrm{d}x = \sum_{i=1}^{\dim \mathcal{P}} \alpha_i \sigma_i.$$

Therefore, generating a set of degrees of freedom from a basis of $\mathcal{P}$ is enough to completely define the set of degrees of freedom (3). Note, however, that the degrees of freedom are not necessarily dual to the basis that have been used in their definition.

***Example.*** In $\mathbb{P}_1(\mathbb{R})$, let us consider the degrees of freedom

$$\sigma_1\colon q \mapsto \int_0^1 q\,x\,\mathrm{d}x$$

$$\sigma_1\colon q \mapsto \int_0^1 q\,\mathrm{d}x.$$

There, the moment $\int_0^1 p\,x\,\mathrm{d}x$ does not determine directly the coefficient applied to the basis function $x \mapsto x$. Indeed, let us take by example $q \equiv 1$, a polynomial of degree zero. We have

$$\int_0^1 1 \times x\,\mathrm{d}x = \frac{1}{2}[x^2]_0^1 = \frac{1}{2}$$

$$\int_0^1 1\,\mathrm{d}x = [x]_0^1 = 1,$$

which obviously prevents this set of degrees of freedom to be dual to the basis $\{1,\, x\}$. ♦

## 2.4 Basis tuning on degrees of freedom

Two different sets of basis functions that generate a same vector space do not necessarily share the same properties. However, it is possible to find an application mapping each other. As a consequence, one can design a specific basis from any generic basis and a set of quantification functions that express some wished properties through a duality relationship.

One just has to pay attention that those wished properties are compatible with the considered space. Typically, they should not contradict each other neither be equivalent, and for each property there should exist functions that belong to the vector space enjoying it. Furthermore, the set of relations enforcing them should be of duality type, and the set of the corresponding quantification functions should have the dimension of the vector space. We assume in this paragraph that all the above requirements hold, and detail the mapping.

Let us consider some vector space $V$ endowed with a basis $\{\phi_j\}_j$ and a set of linear applications $\{\sigma_i\}_i$ on $K$ from $V$ to $\mathbb{R}$ whose dimension matches the one of $V$ (such applications may be called degrees of freedom and can be seen as quantifiers. We refer to the above *Section 2.3* for more details. It is then possible to define a new basis $\{\varphi_j\}_j$ that will be dual to $\{\sigma_i\}_i$.

Indeed, since $\{\phi_j\}_j$ and $\{\varphi_j\}_j$ span the same vector space, we can write

$$\varphi_j = \sum_{m=1}^{N} \alpha_{jm}\phi_m \tag{4}$$

for some coefficients $\{\alpha_{jm}\}_{(j,\,m)\in[\![1,\,N]\!]^2}$, and where $N$ is the dimension of $V$. Furthermore, by definition of the degrees of freedom and of their related basis $\varphi$, we have for any $i$, $j$ both belonging to $[\![1,\,N]\!]$:

$$\sigma_i(\varphi_j) = \delta_{ij}. \tag{5}$$

We can then compose (5) from the left and retrieve

$$\forall i,\, j \in [\![1,\,N]\!],\ \sigma_i(\varphi_j) = \delta_{ij} = \sigma_i\left(\sum_{m=1}^{N} \alpha_{jm}\phi_m\right).$$

Then, by the linearity of degrees of freedom we get

$$\forall i,\, j \in [\![1,\,N]\!],\ \sigma_i(\varphi_j) = \delta_{ij} = \sum_{m=1}^{N} \alpha_{jm}\sigma_i\left(\phi_m\right).$$

When rewrote into a matrix form, the above relation writes

$$I_N = \Lambda A^T,$$

where $(\Lambda)_{ij} = \sigma_i(\phi_j)$ and $(A)_{jm} = \alpha_{jm}$. We then retrieve the coefficients $\alpha_{im}$ by computing

$$A = \left(\Lambda^{-1}\right)^T, \tag{6}$$

and get the basis functions $\{\varphi_j\}_j$ directly through the relation (4).

***Note.*** If the triplet $\{K, V, \{\sigma_i\}_i\}$ is unisolvent then the matrix $\Lambda$ is invertible. Therefore, under the setting of the paragraph the definition of $A$ through the equation (6) is always possible. ▲

We can also retrieve similarly the following relations on the derivatives;

$$\begin{aligned} \mathrm{D}_x\phi_{ij} &= \sigma_i(\partial_x\varphi_j), \\ (\mathrm{D}_x\varphi)_{ij} &= (\sigma_i(\partial_x\varphi_j))_{ij}, \\ \mathrm{D}_x\phi A^T &= \mathrm{D}_x\phi\Lambda^{-1} = \mathrm{D}_x, \varphi, \end{aligned}$$

holding for all $i,\ j \in [\![1,\ N]\!]$ and where $\varphi = (\varphi_1,\ \cdots,\ \varphi_n)^T$. Gathering all the the above transformations, the basis tuning reduces to the following relations.

**Definition 2.4** Tuning of basis functions

$$\begin{cases} (\alpha_{jm})_{jm} = ((\sigma_i(\phi_j))_{ij})^{-T} \\ \varphi_j(x) = \displaystyle\sum_{m=1}^{N} \alpha_{jm}\phi_m(x) \\ \partial_x\varphi_j(x) = \displaystyle\sum_{m=1}^{N} \alpha_{jm}\partial_x\phi_m(x). \end{cases} \tag{8}$$

***Note.*** The above relations emphases that the choice of the linear applications $\{\sigma_i\}_i$ matches the choice of the properties one wants to fulfil. Indeed, applying them to the corresponding set of basis functions $\{\phi_j\}_j$ will impose by duality the properties (5) formalising the following behaviour.

$$\underbrace{\overset{\text{Quantificaton}}{\sigma_i} \quad \overset{\text{Basis function}}{(\varphi_j)} \quad = \quad \overset{\text{Duality relation}}{\delta_{ij}}}_{\text{Set of desired properties}}$$

Let us give some example on a more practical side, establishing a connection with the previous section.

Setting a moment based quantification through the application $\sigma\colon q \mapsto \frac{1}{|K|}\int_K q\,\mathrm{d}x$ quantifies $q$ by its mean value on the cell $K$, and set $N = \dim V$ properties to the set of basis functions. Indeed, for some basis function $\varphi_i$, $i \in [\![1,\ N]\!]$, we then impose: $\frac{1}{|K|}\int_K \varphi_i \mathrm{d}x = 1$ and for any $j \in [\![1,\ N]\!]$ such that $j \neq i$, $\frac{1}{|K|}\int_K \varphi_i \mathrm{d}x = 0$.

Prescribing pointwise values is also possible. Assuming that one would like to prescribe the midpoint value to be zero for one specific basis function and one for the other $N-1$ ones, one can simply set $\sigma\colon q \mapsto 1 - q(\frac{1}{2})$ as a quantification function. During the tuning, it will match each basis function with its desired value at the midpoint.

All the quantifiers together will automatically tune the basis functions towards all the desired $N^2$ properties by the above duality relations. ▲

## 2.5 Basis functions adjustment to geometry distortion

As the later considerations will apply on some reference element $K$, we need to be able to deal with any of its distortion. Yet, to compile with the strive for $H(\mathrm{div})$ – conformity, we require an orientation preserving map.

### 2.5.1 Orientation preserving map for simplicial elements

In the case of simplicial elements, the Piola transform can be used straightforwardly. Although first defined from a reference element to a target element, the backward transformation can be done similarly. We succinctly recall it for the sake of completeness, and adopt the notations introduced in the *Figure 1.*

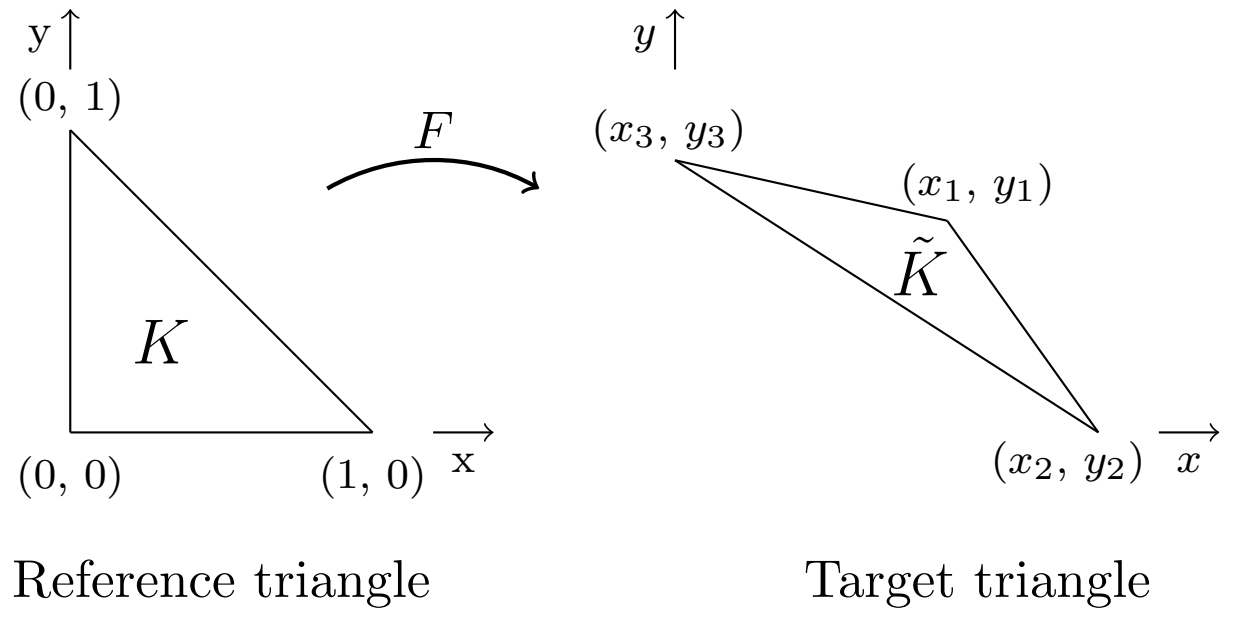

Fig. 1: Illustration of notations used to define the Piola transform

***Note.*** The Piola transform is valid in any dimension. However, for the sake of legibility the underlying change of coordinates will only be presented in two dimensions. A generalisation can immediately be drawn. ▲

First, we define the change of coordinates through the application $F$ by

$$F\colon (\mathrm{x},\, \mathrm{y}) \mapsto \begin{cases} x = x_1 + (x_2 - x_1)\mathrm{x} + (x_3 - x_1)\mathrm{y} \\ y = y_1 + (y_2 - y_1)\mathrm{x} + (y_3 - y_1)\mathrm{y} \end{cases}$$

and derive its Jacobian matrix to get

$$J_F = \begin{pmatrix} \frac{\partial F_x}{\partial \mathrm{x}} & \frac{\partial F_x}{\partial \mathrm{y}} \\ \frac{\partial F_y}{\partial \mathrm{x}} & \frac{\partial F_y}{\partial \mathrm{y}} \end{pmatrix} = \begin{pmatrix} x_2 - x_1 & x_3 - x_1 \\ y_2 - y_1 & y_3 - y_1 \end{pmatrix}.$$

Setting $|J_F| = |\det(J_F)|$, the Piola transform then reads as follows.

**Definition 2.5** Piola transform

$$\begin{aligned} \mathcal{P}\colon L^2(K) &\longrightarrow L^2(\tilde{K}) \\ \phi &\longmapsto \tilde{\phi} = \frac{1}{|J_F|} J_F \phi. \end{aligned} \tag{9}$$

***Remark.*** The transformation is orientation and regularity preserving. In particular, for any $\phi$ in $\mathcal{C}^1(K)$, $\tilde{\phi}$ belongs to $\mathcal{C}^1(\tilde{K})$. ▲

The same map also transfers to the gradient. In particular it holds

$$\begin{aligned} \nabla_{(x,y)}\tilde{\phi} &= \frac{1}{|J_F|} J_F (\nabla_{(\mathrm{x},\mathrm{y})}\phi) J_F^{-1} \\ \nabla_{(x,y)} \cdot \tilde{\phi} &= \frac{1}{|J_F|} \nabla_{(\mathrm{x},\mathrm{y})} \cdot \phi. \end{aligned}$$

Furthermore, we have for any $\phi \in H(\mathrm{div},\, K)$;

$$\begin{aligned} \int_{\tilde{K}} \tilde{\phi}_i \cdot \tilde{\phi}_j \,\mathrm{d}x &= \int_K \frac{1}{|J_F|} J_F \phi_i \cdot J_F \phi_j \,\mathrm{dx} \\ \int_{\tilde{K}} \tilde{\phi} \cdot \nabla_x \tilde{v} \,\mathrm{d}x &= \int_K \phi \cdot \nabla_{\mathrm{x}} v \,\mathrm{dx} \\ \int_{\tilde{K}} \nabla_x \cdot \tilde{\phi}\, \tilde{v} \,\mathrm{d}x &= \int_K \nabla_{\mathrm{x}} \cdot \phi\, v \,\mathrm{dx} \\ \int_{\partial\tilde{K}} \tilde{\phi} \cdot n\, \tilde{v} \,\mathrm{dx} &= \int_{\partial K} \phi \cdot n\, v \,\mathrm{d}x, \end{aligned}$$

where for any $v \in H^1(K)$, $\tilde{v} = v \circ F^{-1}$. We refer to [19] for a proof, whose key idea lies in the relation $\mathrm{dx} = |J_F| \mathrm{d}x$.

As a consequence, the degrees of freedom and basis functions can be defined through the reference element $K$ on any distorted element $\tilde{K}$.

### 2.5.2 Orientation preserving map for quadrilateral elements

For quadrilateral elements, the Piola transform can also be used, though a change in the underlying map of coordinates is required. Indeed, contrarily to simplicial shapes, mapping one quadrilateral element to another requires the mixed term $xy$ to be taken into account. We succinctly point out the necessary changes using the notations of the *Figure 2*.

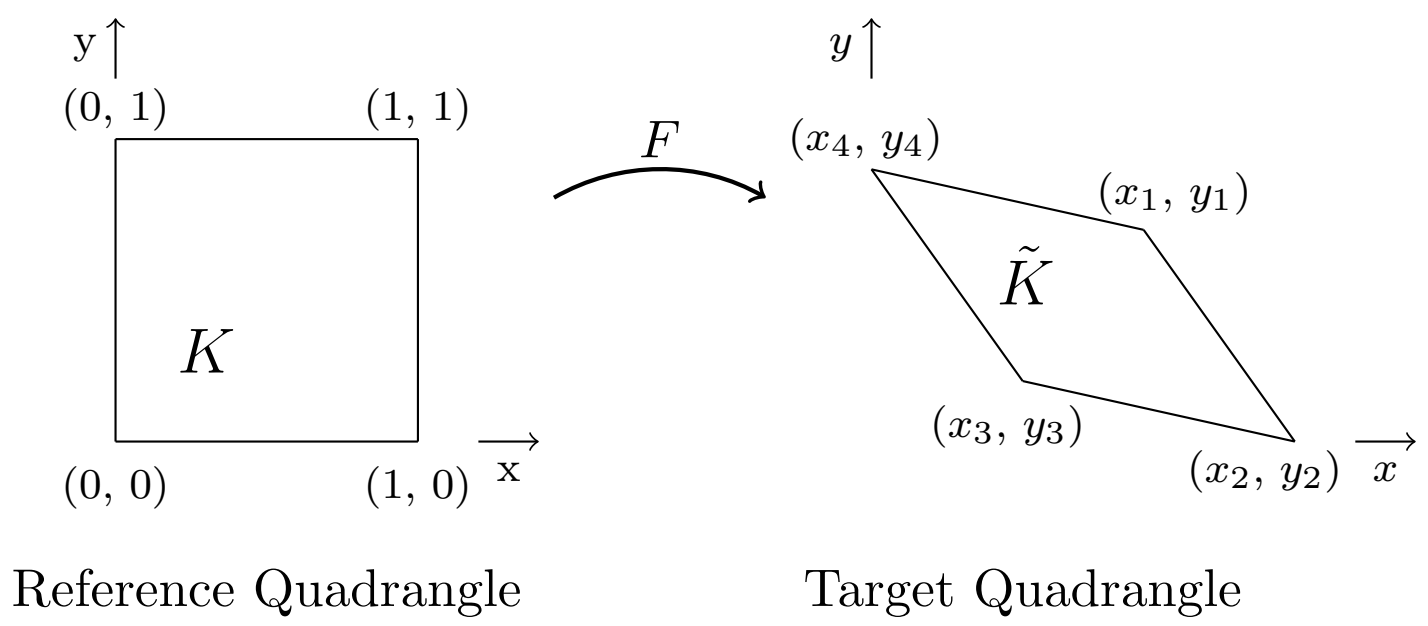

Fig. 2: Illustration of notations used to define the Piola transform in the quadrangle case

***Note.*** As in the simplicial case, though the Piola transform applied on quadrilateral shapes is valid in any dimension, the underlying change of coordinates will only be presented in two dimensions for the sake of legibility. ▲

Here, the used coordinates mapping from a reference quadrangle to a general one reads

$$F\colon (\mathrm{x},\, \mathrm{y}) \mapsto \begin{cases} x = a_1 + b_1\mathrm{x} + c_1\mathrm{y} + d_1\mathrm{xy} \\ y = a_2 + b_2\mathrm{x} + c_2\mathrm{y} + d_2\mathrm{xy}, \end{cases}$$

where the coefficients $\{a_i,\, b_1,\, c_1,\, d_i\}_{i1=1,2}$ are determined by the relations

$$A \begin{pmatrix} a_1 \\ b_1 \\ c_1 \\ d_1 \end{pmatrix} = \begin{pmatrix} 0 \\ 1 \\ 1 \\ 0 \end{pmatrix} \quad \text{and} \quad A \begin{pmatrix} a_2 \\ b_2 \\ c_2 \\ d_2 \end{pmatrix} = \begin{pmatrix} 0 \\ 0 \\ 1 \\ 1 \end{pmatrix}, \tag{10}$$

with

$$A = \begin{pmatrix} 1 & x_1 & y_1 & x_1y_1 \\ 1 & x_2 & y_2 & x_2y_2 \\ 1 & x_3 & y_3 & x_3y_3 \\ 1 & x_4 & y_4 & x_4y_4 \end{pmatrix}.$$

Note that the second members of the equations defined in (10) are respectively given by the $x$ and $y$ coordinates values from the vertices of the reference element. We then can compute the Jacobian of the transformation $F$, and plug it into the definition of the Piola transform given in (9). The same properties as in the simplicial case hold, and will not be repeated.

### 2.5.3 A word on orientation preserving map for polytopial elements

Defining an orientation preserving map for polytopial elements is much more delicate, in particular when one wants to allow the use of non-convex polytopes.

In addition to finding a mapping valid for any number of faces, one has to determine whether or not the normals orientation should be preserved, as depicted in the *Figure 3*. This choice is highly dependent both on the global shape and on the position of the distorted element with respect to the axes.

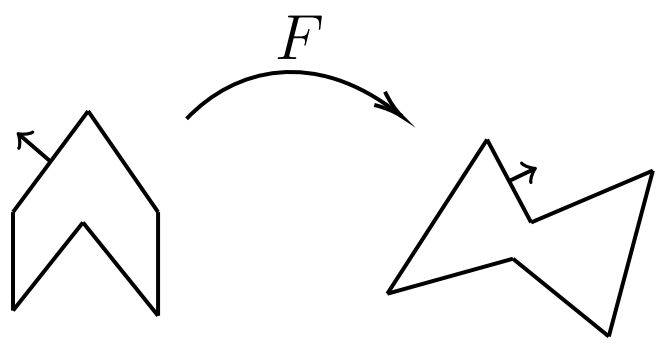

Fig. 3: Illustration of a possible transformation in the hexagonal case, where the shape's inflexion implies a change in the normal's definition

In particular, if the distorted element has a different convexity that the reference one, the normal components of the impacted faces should not be preserved (as pictured in the *Figure 3*). However, if the element is also rotated so that the impacted edge is comparable to the one of reference, then it should be preserved (*see the Figure 4*).

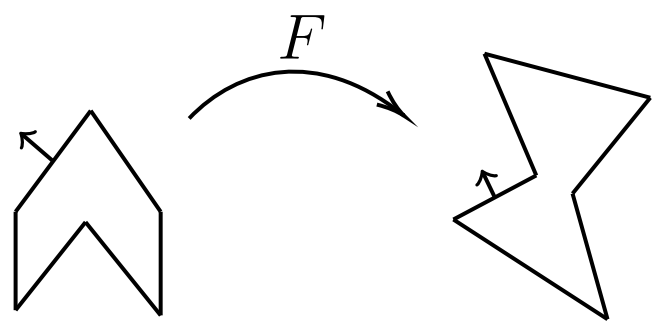

Fig. 4: Illustration of a transformation in the hexagonal case, where the normal of the edge impacted by the inflexion change should be preserved.

We then need to redefine the concept of orientation preserving, and determine a transformation which is dependent on the shape of the polytope, at least locally. This will be considered in a future work.

# 3 Definition of the simplicial Raviart – Thomas elements

We start by recalling the origin of the Raviart – Thomas elements [16] for simplicial shapes. The basic concept consists in working in a vectorial polynomial discrete subspace of $H(\mathrm{div}, K)$ and to enforce the conformity with proper polynomial projection of the derivatives. We progressively derive this idea paraphrasing [13] for the sake of clarity.

***Note.*** As in this section only simplicial shapes are considered, the involved finite dimensional polynomial space is set as $\mathbb{P}_k$. ▲

## 3.1 Framework of Nédélec

Let us progressively build the roots that lead to the classical definition of Raviart – Thomas spaces. All the pictured results are taken from [13] while the presented proofs are adaptations of the work of [11].

### 3.1.1 Structure of the wished discretisation

Let $K$ be a general simplicial reference shape contained in $\mathbb{R}^d$, where $d \in \mathbb{N}$ is the dimension of the space where the problem's domain lies. The aim is to find a polynomial approximation space $V$ of $H(\mathrm{div},\, K)$, where

$$H(\mathrm{div}, K) = \{u \in \left(L^2(K)\right)^d,\, \mathrm{div}\, u \in L^2(K)\}.$$

Indeed, as any space $V \subset H(\mathrm{div},\, K)$ is by definition $H(\mathrm{div})$ – conformity friendly, one can naturally design a set of degrees of freedom $\{\sigma\}$ such that the triplet $(K,\, V,\, \{\sigma\})$ defines a conformal element.

In the perspective of designing such a space, we start by noticing that a mean to ensure $(\mathrm{div}\, u) \in L^2(K)$ for any $u \in V$ is to prescribe some $k^{\text{th}}$ - order derivative of $u$ to vanish. Then, following the work of Nédélec [13], a symmetrised linear differential operator $\varepsilon^k$ is introduced. When applied to vectorial polynomial functions its kernel generates $V$ as a finite dimensional subspace of $H(\mathrm{div}, K)$, defining the simplicial Raviart – Thomas space. Let us now develop this construction.

### 3.1.2 Expression of the symmetrised differential operator

As a preamble to the construction of the Raviart – Thomas space, the differential operator $\varepsilon^k$ is introduced in its original definition before being reduced for the sake of later convenience.

**Construction of the symmetrised differential operator** It is set

$$\begin{aligned} \varepsilon^k \colon \mathbb{R}^d \times \cdots \times \mathbb{R}^d &\longrightarrow \mathbb{R} \\ (\chi_1, \cdots, \chi_{k+1}) &\longmapsto \frac{1}{(k+1)!} \sum_{\zeta \in \zeta_{k+1}} \mathrm{D}^k u(\chi_{\zeta(1)}, \cdots, \chi_{\zeta(k+1)}). \end{aligned} \tag{11}$$

There, the term $\zeta_{k+1}$ denotes the set of all permutations of $(1, \cdots, k+1)$ and $\mathrm{D}^k$ denotes the usual $(k+1)$ - differential form described below. Moreover, $u \colon \mathbb{R}^d \to \mathbb{R}^d$ is a $d$ - valued function whose components are given by $u = (u_1, \cdots, u_d)^T$ and whose expression along any directional vector $\chi_1 = (\chi_1^1, \cdots \chi_1^d)^T \in \mathbb{R}^d$ reads $u(x)(\chi_1) = u_1(x)\chi_1^1 + \cdots + u_d(x)\chi_1^d$ for any $x \in \mathbb{R}^d$.

***Note.*** Though implicit from the definition, we point out that for any directional vector $\chi_1 \in \mathbb{R}^d$, the function

$$\begin{aligned} \mathbb{R}^d &\longrightarrow \mathbb{R} \\ x &\longmapsto u(x)(\chi_1) \end{aligned}$$

is scalar. Also note that in the following, the dependency of $u$ on the variable $x$ is dropped in the notation, so that the solely term $u(\chi_1)$ refers to the above function. ▲

Let us now explicit the operator $\mathrm{D}^k$. First of all, the zero-th order differential operator $\mathrm{D}^0$ immediately reads as

$$\mathrm{D}^0 u(\chi_1) = u(\chi_1) = u \cdot \chi_1$$

by definition of $u$. For higher order derivatives, a set $\chi = (\chi_1, \cdots, \chi_{k+1})$ of $(k+1)$ directional vectors $\{\chi_i\}_i$ belonging to $\mathbb{R}^d$ is defined, and the $(k+1)$ - differential form $\mathrm{D}^k$ is set to

$$\mathrm{D}^k u(\chi_1, \cdots, \chi_{k+1}) = \bigcirc_{j=2}^{k+1} \Bigg( \underbrace{\left( \frac{\partial}{\partial x} u(\chi_1) \right) \cdot \chi_j}_{:= f_j(u(\chi_1))} \Bigg). \tag{12}$$

Here, the term $\left(\frac{\partial}{\partial x} u(\cdot)\right)$ represents the gradient vector $\left(\frac{\partial}{\partial x_1} u(\cdot), \cdots, \frac{\partial}{\partial x_d} u(\cdot)\right)^T$ and the operator $\bigcirc_{j=2}^{k+1} f_j \chi_j$ stands for the composition operator

$$f_1 \longmapsto f_{k+1}(f_k(\cdots(f_2(f_1) \cdot \chi_2) \cdots) \cdot \chi_k) \cdot \chi_{k+1},$$

which for any function $f_1 \colon \mathbb{R}^d \mapsto \mathbb{R}$ maps to another function of $\mathbb{R}^d$ taking its values in $\mathbb{R}$.

***Example.*** In two dimensions, the function $u$ is of the shape

$$\begin{aligned} u\colon\ \mathbb{R}^2 &\longrightarrow \mathbb{R}^2 \\ \begin{pmatrix} x \\ y \end{pmatrix} &\longmapsto \begin{pmatrix} u_1(x,\, y) \\ u_2(x,\, y) \end{pmatrix}. \end{aligned}$$

Thus, tested against any directional vector $\chi_1 \in \mathbb{R}^2$ it comes

$$u(\chi_1) = u_1\chi_1^1 + u_2\chi_1^2,$$

which gives the expression of $\mathrm{D}^0 u(\chi_1)$. Then, for any $\chi_2 \in \mathbb{R}^2$, the first order derivative $\mathrm{D}^1 u(\chi_1,\, \chi_2)$ reads

$$\begin{aligned} \mathrm{D}^1 u(\chi_1, \chi_2) &= \frac{\partial}{\partial x}\left(u(\chi_1)\right) \cdot \chi_2 \\ &= \begin{pmatrix} \frac{\partial}{\partial x_1}\left(u_1\chi_1^1 + u_2\chi_1^2\right) \\ \frac{\partial}{\partial x_2}\left(u_1\chi_1^1 + u_2\chi_1^2\right) \end{pmatrix} \cdot \chi_2 \\ &= \frac{\partial u_1}{\partial x_1}\chi_1^1\chi_2^1 + \frac{\partial u_2}{\partial x_1}\chi_1^2\chi_2^1 + \frac{\partial u_1}{\partial x_2}\chi_1^1\chi_2^2 + \frac{\partial u_2}{\partial x_2}\chi_1^2\chi_2^2. \end{aligned}$$

Going to the second order with the same spirit, the differential operator $\mathrm{D}^2$ is expressed following (12) as

$$\mathrm{D}^2 u(\chi_1,\, \chi_2,\, \chi_3) = \frac{\partial}{\partial x}\left(\mathrm{D}^1 u(\chi_1,\, \chi_2)\right) \cdot \chi_3$$

for any $(\chi_1,\, \chi_2,\, \chi_3) \in \mathbb{R}^{2\times 3}$. ♦

As working directly from the quantity $u(\chi_1)$ turns to be delicate when expressing the kernel of $\varepsilon^k$, we derive an alternative form of (12) that can be expressed solely in terms of $u$. Using the definition of $u(\chi_1)$ and the linearity of the operator $\mathrm{D}^k$ to swap $\chi_1$ and $\chi_j$, we first reorder the directional vectors to get

$$\begin{aligned} \mathrm{D}^k u(\chi_1,\, \cdots,\, \chi_{k+1}) &= \bigcirc_{j=2}^{k+1} \left(\frac{\partial}{\partial x} u(\chi_1)\right) \cdot \chi_j \\ &= \bigcirc_{j=2}^{k+1} \left(\frac{\partial u}{\partial x} \cdot \chi_1\right) \cdot \chi_j \\ &= \bigcirc_{j=2}^{k+1} \left(\frac{\partial u}{\partial x} \cdot \chi_j\right) \cdot \chi_1. \end{aligned}$$

Then, using once again this linearity, the last expression becomes

$$\begin{aligned}
\mathrm{D}^k u(\chi_1, \cdots, \chi_{k+1}) &= \bigcirc_{j=2}^{k+1} \left( \frac{\partial}{\partial x} \begin{pmatrix} u_1 \\ \vdots \\ u_d \end{pmatrix} \cdot \chi_j \right) \cdot \chi_1 \\
&= \begin{pmatrix} \bigcirc_{j=2}^{k+1} \frac{\partial}{\partial x} u_1 \cdot \chi_j \\ \vdots \\ \bigcirc_{j=2}^{k+1} \frac{\partial}{\partial x} u_d \cdot \chi_j \end{pmatrix} \cdot \chi_1,
\end{aligned}$$

which is nothing else than the $k^{\text{th}}$ - order Fréchet derivatives' vector of each component of $u$ tested along the directional vector $\chi_1$. Thus, by defining for any $l \in [\![1, d]\!]$ the Fréchet derivative as

$$\mathrm{d}_l^k u(\chi_2, \cdots, \chi_{k+1}) = \bigcirc_{j=2}^{k+1} \left( \frac{\partial}{\partial x} u_l \right) \cdot \chi_j$$

and by denoting the corresponding vector

$$\mathrm{d}^k u(\chi_2, \cdots, \chi_{k+1}) = \begin{pmatrix} \mathrm{d}_1^k u(\chi_2, \cdots, \chi_{k+1}) \\ \vdots \\ \mathrm{d}_d^k u(\chi_2, \cdots, \chi_{k+1}) \end{pmatrix},$$

the definition of $\mathrm{D}^k$ becomes

$$\mathrm{D}^k u(\chi_1, \cdots, \chi_{k+1}) = \left( \bigcirc_{j=2}^{k+1} \mathrm{d}^k u(\chi_2, \cdots, \chi_{k+1}) \right) \cdot \chi_1.$$

As a consequence, the definition of $\varepsilon^k$ comes down to

$$\begin{aligned}
\varepsilon^k \colon \mathbb{R}^d \times \cdots \times \mathbb{R}^d \quad &\longrightarrow \quad \mathbb{R} \\
(\chi_1, \cdots, \chi_{k+1}) \quad &\longmapsto \quad \frac{1}{(k+1)!} \sum_{\zeta \in \zeta_{k+1}} \mathrm{d}^k u(\chi_{\zeta(2)}, \cdots, \chi_{\zeta(k+1)}) \cdot \chi_{\zeta(1)}.
\end{aligned}$$

Lastly, as the permutations annihilate the importance of the index ordering, we can shuffle arbitrarily the vector $\chi = (\chi_1, \cdots, \chi_{k+1})$. Thus, up to shifting the indices of the Fréchet derivative so that

$$\mathrm{d}_l^k u(\chi_1, \cdots, \chi_k) = \bigcirc_{j=1}^{k} \left( \frac{\partial}{\partial x} u_l \right) \cdot \chi_j,$$

we finally obtain a more intuitive expression of $\varepsilon^k$ that first treats the derivatives of the functional $u$ in a vectorial way before testing them against the directional vector $\chi_1$.

**Definition 3.1** Expression of $\varepsilon^k$

$$\begin{aligned} \varepsilon^k\colon \mathbb{R}^d \times \cdots \times \mathbb{R}^d &\longrightarrow \mathbb{R} \\ (\chi_1, \cdots, \chi_{k+1}) &\longmapsto \frac{1}{(k+1)!} \sum_{\zeta\in\zeta_{k+1}} \mathrm{d}^k u(\chi_{\zeta(1)}, \cdots, \chi_{\zeta(k)}) \cdot \chi_{\zeta(k+1)}. \end{aligned} \tag{13}$$

***Example.*** *Continued* Let us take the same case as in the previous example where $d = 2$. Using the previous computations of $\mathrm{D}^k$, we derive the expression of $\varepsilon^k$ from this last formulation by considering any directional vectors $\chi_1$, $\chi_2$ and $\chi_3$ belonging to $\mathbb{R}^2$.

• In the case $k = 1$, the definition of $\mathrm{D}^1$ reads

$$\mathrm{D}^1 u(\chi_1, \chi_2) = \begin{pmatrix} \mathrm{d}_1^1 u(\chi_1) \\ \mathrm{d}_2^1 u(\chi_1) \end{pmatrix} \cdot \chi_2,$$

where

$$\begin{aligned} \mathrm{d}_1^1 u(\chi_1) &= \left(\frac{\partial u_1}{\partial \mathrm{x}}\right) \cdot \chi_1 = \begin{pmatrix} \frac{\partial u_1}{\partial x} \\ \frac{\partial u_1}{\partial y} \end{pmatrix} \cdot \begin{pmatrix} \chi_1^1 \\ \chi_1^2 \end{pmatrix} = \left(\frac{\partial u_1}{\partial x}\chi_1^1 + \frac{\partial u_1}{\partial y}\chi_1^2\right) \\ \mathrm{d}_2^1 u(\chi_1) &= \left(\frac{\partial u_2}{\partial \mathrm{x}}\right) \cdot \chi_1 = \begin{pmatrix} \frac{\partial u_2}{\partial x} \\ \frac{\partial u_2}{\partial y} \end{pmatrix} \cdot \begin{pmatrix} \chi_1^1 \\ \chi_1^2 \end{pmatrix} = \left(\frac{\partial u_2}{\partial x}\chi_1^1 + \frac{\partial u_2}{\partial y}\chi_1^2\right). \end{aligned}$$

By developing it comes

$$\mathrm{D}^1 u(\chi_1, \chi_2) = \left(\frac{\partial u_1}{\partial x}\chi_1^1 + \frac{\partial u_1}{\partial y}\chi_1^2\right)\chi_2^1 + \left(\frac{\partial u_2}{\partial x}\chi_1^1 + \frac{\partial u_2}{\partial y}\chi_1^2\right)\chi_2^2.$$

As the set of permutations of $\{\chi_1, \chi_2\}$ is only given by $\{(\chi_1, \chi_2), (\chi_2, \chi_1)\}$, we obtain by definition of $\varepsilon^1$;

$$\begin{aligned} \varepsilon^1 u(\chi_1, \chi_2) &= \frac{1}{2}\Bigg[\left(\frac{\partial u_1}{\partial x}\chi_1^1 + \frac{\partial u_1}{\partial y}\chi_1^2\right)\chi_2^1 + \left(\frac{\partial u_2}{\partial x}\chi_1^1 + \frac{\partial u_2}{\partial y}\chi_1^2\right)\chi_2^2 \\ &\quad + \left(\frac{\partial u_1}{\partial x}\chi_2^1 + \frac{\partial u_1}{\partial y}\chi_2^2\right)\chi_1^1 + \left(\frac{\partial u_2}{\partial x}\chi_2^1 + \frac{\partial u_2}{\partial y}\chi_2^2\right)\chi_1^2\Bigg] \\ &= \frac{1}{2}\Bigg[\left(\frac{\partial u_1}{\partial x} + \frac{\partial u_1}{\partial x}\right)\chi_1^1\chi_2^1 + \left(\frac{\partial u_1}{\partial y} + \frac{\partial u_2}{\partial x}\right)\chi_1^2\chi_2^1 \\ &\quad + \left(\frac{\partial u_2}{\partial y} + \frac{\partial u_2}{\partial y}\right)\chi_1^2\chi_2^2 + \left(\frac{\partial u_1}{\partial y} + \frac{\partial u_2}{\partial x}\right)\chi_1^1\chi_2^2\Bigg]. \end{aligned}$$

Furthermore, by its linearity the operator $\mathrm{D}^1$ is symmetric with respect to the vectors $\{\chi_i\}_i$. The above expression then reduces to

$$\begin{aligned}\varepsilon^k u(\chi_1,\,\chi_2) &= \frac{\partial u_1}{\partial x}\chi_1^1\chi_2^1 + \frac{\partial u_2}{\partial y}\chi_1^2\chi_2^2 \\ &\quad + \frac{1}{2}\left(\frac{\partial u_1}{\partial y} + \frac{\partial u_2}{\partial x}\right)\left(\chi_1^2\chi_2^1 + \chi_1^1\chi_2^2\right).\end{aligned}$$

- In the case $k = 2$ it holds in a similar way;

$$\begin{aligned}\mathrm{D}^2 u(\chi_1,\,\chi_2,\,\chi_3) &= \begin{pmatrix}\frac{\partial}{\partial x}\mathrm{D}^1 u(\chi_1,\,\chi_2) \\ \frac{\partial}{\partial y}\mathrm{D}^1 u(\chi_1,\,\chi_2)\end{pmatrix}\cdot\chi_3 \\ &= \begin{pmatrix}\frac{\partial^2 u_1}{\partial x^2}\chi_1^1\chi_2^1 + \frac{\partial^2 u_1}{\partial x\partial y}\chi_1^2\chi_2^1 + \frac{\partial^2 u_2}{\partial x^2}\chi_1^1\chi_2^2 + \frac{\partial^2 u_2}{\partial x\partial y}\chi_1^2\chi_2^2 \\ \frac{\partial^2 u_1}{\partial y\partial x}\chi_1^1\chi_2^1 + \frac{\partial^2 u_1}{\partial y^2}\chi_1^2\chi_2^1 + \frac{\partial^2 u_2}{\partial y\partial x}\chi_1^1\chi_2^2 + \frac{\partial^2 u_2}{\partial y^2}\chi_1^2\chi_2^2\end{pmatrix}\cdot\begin{pmatrix}\chi_3^1 \\ \chi_3^2\end{pmatrix} \\ &= \frac{\partial^2 u_1}{\partial x^2}\chi_1^1\chi_2^1\chi_3^1 + \frac{\partial^2 u_2}{\partial y^2}\chi_1^2\chi_2^2\chi_3^2 \\ &\quad + \frac{\partial^2 u_2}{\partial x^2}\chi_1^1\chi_2^2\chi_3^1 + \frac{\partial^2 u_1}{\partial x\partial y}\chi_1^2\chi_2^1\chi_3^1 + \frac{\partial^2 u_1}{\partial y\partial x}\chi_1^1\chi_2^1\chi_3^2 \\ &\quad + \frac{\partial^2 u_1}{\partial y^2}\chi_1^2\chi_2^1\chi_3^2 + \frac{\partial^2 u_2}{\partial x\partial y}\chi_1^2\chi_2^2\chi_3^1 + \frac{\partial^2 u_2}{\partial y\partial x}\chi_1^1\chi_2^2\chi_3^2\end{aligned}$$

Here, by definition, the operator $\varepsilon^2$ is generated from the set of permutations $\{(\chi_1,\,\chi_2,\,\chi_3),\,(\chi_1,\,\chi_3,\,\chi_2),\,(\chi_2,\,\chi_1,\,\chi_3),\,(\chi_2,\,\chi_3,\,\chi_1),\,(\chi_3,\,\chi_2,\,\chi_1),\,(\chi_3,\,\chi_1,\,\chi_2)\}$, which has for cardinal six. Furthermore, it can be seen that all the three vectors $\{\chi_i\}_i$ involved in the permutations appear conjointly in each term of the above expression. The symmetry therefore only redefines the directional vectors without interconnecting the derivative terms. In addition, as each of the two first terms is involving only the same spatial coordinate $\chi_i^1$ or $\chi_i^2$ of each directional vector, they are totally insensitive to the permutation. Thus, the operator $\varepsilon^k$ reads

$$\begin{aligned}\varepsilon^2 u(\chi_1,\,\chi_2,\,\chi_3) = \frac{1}{6}\Bigg[&6\,\frac{\partial^2 u_1}{\partial x^2}\chi_1^1\chi_2^1\chi_3^1 + 6\,\frac{\partial^2 u_2}{\partial y^2}\chi_1^2\chi_2^2\chi_3^2 \\ &+ \left(\frac{\partial^2 u_1}{\partial y^2} + 2\frac{\partial^2 u_2}{\partial x\partial y}\right)\left(\chi_1^2\chi_2^2\chi_3^1 + \chi_2^2\chi_3^2\chi_1^1 + \chi_3^2\chi_1^2\chi_2^1\right. \\ &\qquad \left. + \chi_2^2\chi_1^2\chi_3^1 + \chi_3^2\chi_2^2\chi_1^1 + \chi_1^2\chi_3^2\chi_2^1\right) \\ &+ \left(\frac{\partial^2 u_2}{\partial x^2} + 2\frac{\partial^2 u_1}{\partial x\partial y}\right)\left(\chi_1^2\chi_2^1\chi_3^1 + \chi_2^2\chi_3^1\chi_1^1 + \chi_3^2\chi_1^1\chi_2^1\right. \\ &\qquad \left. + \chi_2^2\chi_1^1\chi_3^1 + \chi_3^2\chi_2^1\chi_1^1 + \chi_1^2\chi_3^1\chi_2^1\right)\Bigg],\end{aligned}$$

where we have used the fact that $u \in L^2(K)$, ensuring that $\frac{\partial u}{\partial x\partial y} = \frac{\partial u}{\partial y\partial x}$. ♦

**Reduction of the expression of $\varepsilon^k$** As it appears on the above example, the form of $\varepsilon^k$ can be reduced by using the symmetry property of the underlying operator $\mathrm{D}^k$. The function $u$ being in $L^2(K)$, a further reduction using the commutativity of the partial differential operators is equally possible, leading to the following simpler form [11] of the operator (3.1).

**Definition 3.2** Reduced expression of $\varepsilon^k$

$$\begin{aligned} \varepsilon^k \colon \mathbb{R}^d \times \cdots \times \mathbb{R}^d &\longrightarrow \mathbb{R} \\ (\chi_1, \cdots, \chi_{k+1}) &\longmapsto \frac{1}{k+1} \sum_{j=1}^{k+1} \chi_j \cdot \mathrm{d}^k u(\chi_1, \cdots, \chi_{j-1}, \chi_{j+1}, \cdots, \chi_{k+1}). \end{aligned} \tag{14}$$

***Example.*** Let us consider again the case $d = 2$ and $k = 1$ to familiarise with this new expression. With $\chi_1$, $\chi_2$ and $\chi_3$ three directional vectors of $\mathbb{R}^2$, reading the expression (14) out leads us to

$$\varepsilon^k u(\chi_1, \chi_2) = \frac{1}{2} \left( \chi_1 \cdot \begin{pmatrix} \mathrm{d}_1^1 u(\chi_2) \\ \mathrm{d}_2^1 u(\chi_2) \end{pmatrix} + \chi_2 \cdot \begin{pmatrix} \mathrm{d}_1^1 u(\chi_1) \\ \mathrm{d}_2^1 u(\chi_1) \end{pmatrix} \right),$$

which matches exactly the classical form. In the case $k = 2$, it reads however

$$\varepsilon^k u(\chi_1, \chi_2, \chi_3) = \tfrac{1}{3} \left( \chi_1 \cdot \begin{pmatrix} \mathrm{d}_1^2 u(\chi_2, \chi_3) \\ \mathrm{d}_2^2 u(\chi_2, \chi_3) \end{pmatrix} + \chi_2 \cdot \begin{pmatrix} \mathrm{d}_1^2 u(\chi_1, \chi_3) \\ \mathrm{d}_2^2 u(\chi_1, \chi_3) \end{pmatrix} + \chi_3 \cdot \begin{pmatrix} \mathrm{d}_1^2 u(\chi_1, \chi_2) \\ \mathrm{d}_2^2 u(\chi_1, \chi_2) \end{pmatrix} \right),$$

which naturally reduces to the expression found in the previous example, with less efforts to put into the computation. ♦

### 3.1.3 A first definition of the simplicial Raviart – Thomas space

Having defined the operator $u \mapsto \varepsilon^k u$ by last expression (14), we can now set up an approximation space by defining its application domain, *i.e.* the space where the functions $u$ will lie, and expressing its corresponding kernel.

As we target the use of polynomial spaces, we start by restricting the space where the functions $u$ are defined to the finite dimensional polynomial space $(\mathbb{P}_k(K))^d$. Then, to force the $k^{\text{th}}$ order derivative to vanish while preserving the $H(\mathrm{div})$ – conformity friendly property, we restrict the space $(\mathbb{P}_k(K))^d$ to a subspace $W \subset (\mathrm{div}, K)$ such that

$$W = \{u \in (\mathbb{P}_k(K))^d,\ \mathrm{div}\, u \in \mathbb{P}_{k-1}(K)\}.$$

Lastly, as we want to drive the extinction of the $k^{\text{th}}$ order derivative along the kernel of the operator $\varepsilon^k$, we select the discretisation subspace $V \subset W$ as follows.

**Definition 3.3** Simplicial Raviart – Thomas space

$$V = \{u \in (\mathbb{P}_k(K))^d,\ \varepsilon^k u = 0\}$$

***Note.*** This early definition of the Raviart – Thomas space is not convenient to handle in practice. However, it immediately allows to witness the implication $V \subset (\mathbb{P}_k(K))^d$ and $\mathbb{P}_{k-1}(K) \subset V$. ▲

### 3.1.4 Determination of the kernel of $\varepsilon^k$

Let us now retrieve the classical definition of the simplicial Raviart – Thomas space by making the kernel of $\varepsilon^k$ explicit. We first parametrise $V$ from the derivatives of the functions $u \in (\mathbb{P}_k(K))^d$ before drawing a proper expression for the generated space.

**Parametrisation of the kernel of** $\varepsilon^k$ In order to parametrise the kernel of $\varepsilon^k$, we look for a characterisation of the functions $u \in (\mathbb{P}_k(K))^d$ satisfying

$$\forall(\chi_1, \cdots, \chi_{k+1}) \in \mathbb{R}^{d\times(k+1)}, \quad \varepsilon^k u(\chi_1, \cdots, \chi_{k+1}) = 0. \tag{15}$$

First of all, let us point out that $\varepsilon^k$ is a $(k+1)$ - multilinear form. In particular, for any set of directional vectors $\{\chi_j^{i_j}\}_{\substack{j\in[\![1,\, k+1]\!]\\ i_j\in[\![1,\, d]\!]}}$ belonging to $\mathbb{R}^d$, it holds

$$\begin{aligned}
\varepsilon^k u\left(\sum_{i=1}^{d} \alpha_1^{i_1}\chi_1^{i_1}, \cdots, \sum_{i=1}^{d} \alpha_{k+1}^{i_{k+1}}\chi_{k+1}^{i_{k+1}}\right) &= \sum_{i=1}^{d} \alpha_1^{i_1}\varepsilon^k u\left(\chi_1^{i_1}, \cdots, \sum_{i=1}^{d} \alpha_{k+1}^{i_{k+1}}\chi_{k+1}^{i_{k+1}}\right)\\
&= \sum_{i_1,\cdots,i_{k+1}=1}^{d} \alpha^i \varepsilon^k u\left(\chi_1^{i_1}, \cdots, \chi_{k+1}^{i_{k+1}}\right),
\end{aligned}$$

where $\alpha_i = \Pi_{j=1}^{k+1}\alpha_j^{i_j}$. Note that the set $\{\alpha_j^{i_j},\ j \in [\![1,\, k+1]\!] i_j \in [\![1,\, d]\!]\}$ contains real coefficients and that the term $\chi_j = \sum_{i=1}^{d} \alpha_j^{i_1}\chi_j^{i_j}$ stands for any directional vector of $\mathbb{R}^d$.

***Note.*** The vectors $\{\chi_j^{i_j},\ j \in [\![1,\, k+1]\!],\ i_j \in [\![1,\, d]\!]\}$ are not necessarily linearly independent. ▲

As a consequence, it is enough to consider $\{\chi_j^{i_j}\}_{i_j\in[\![1,\, d]\!]}$ as the canonical basis of $\mathbb{R}^d$ for all $j \in [\![1,\, k+1]\!]$ to represent the impact of any directional vector set $(\chi_1, \cdots, \chi_{k+1}) \in \mathbb{R}^{d\times(k+1)}$ on $\varepsilon^k$. Indeed, by letting the coefficients of the set $\{\alpha_j^{i_j}\}_{ij}$ vary in $\mathbb{R}$, the relation (15) reduces to find all the functions $u \in (\mathbb{P}_k(K))^d$ such that

$$\sum_{i_1,\,\cdots,\,i_{k+1}=1}^{d} \alpha^i \varepsilon^k u\left(b^{i_1},\,\cdots,\,b^{i_{k+1}}\right) = 0, \tag{16}$$

where $\mathcal{B} = \{b^1,\,\cdots,\,b^d\}$ denotes the canonical basis of $\mathbb{R}^d$ and for any $j \in [\![1,\,k+1]\!]$, $\{b^{i_j}\}_i = \mathcal{B}$.

***Example.*** We still consider the case $d = 2$ and $k = 1$. There, the problem (15) reduces to find all the $u \in (\mathbb{P}_1)^2$ such that for all $\chi_1,\,\chi_2 \in \mathbb{R}^2$, it holds

$$\varepsilon^k u(\chi_1,\,\chi_2) = 0.$$

Denoting $\mathcal{B} = \{b^1,\,b^2\}$ the canonical basis of $\mathbb{R}^2$ and defining four real coefficients $\alpha_1^1$, $\alpha_1^2$, $\alpha_2^1$ and $\alpha_2^2$, it comes $\chi_1 = \alpha_1^1 b^1 + \alpha_1^2 b^2$ and $\chi_2 = \alpha_2^1 b^1 + \alpha_2^2 b^2$. The problem (15) then reduces to

$$\begin{aligned}
& \varepsilon^1 u(\chi_1,\,\chi_2) = 0 \\
\Leftrightarrow \quad & \varepsilon^1 u(\alpha_1^1 b^1 + \alpha_1^2 b^2,\,\alpha_2^1 b^1 + \alpha_2^2 b^2) = 0 \\
\Leftrightarrow \quad & \alpha_1^1 \varepsilon^1 u(b^1,\,\alpha_2^1 b^1 + \alpha_2^2 b^2) + \alpha_1^2 \varepsilon^k u(b^2,\,\alpha_2^1 b^1 + \alpha_2^2 b^2) = 0 \\
\Leftrightarrow \quad & \alpha_1^1 \alpha_2^1 \varepsilon^1 u(b^1,\,b^1) + \alpha_1^1 \alpha_2^2 \varepsilon^1 u(b^1,\,b^2) + \alpha_1^2 \alpha_2^1 \varepsilon^1 u(b^2,\,b^1) + \alpha_1^2 \alpha_2^2 \varepsilon^1 u(b^2,\,b^2) = 0,
\end{aligned}$$

which is finally leading to

$$\sum_{i_1,\,i_2=1}^{2} \left( \prod_{j=1}^{2} \alpha_j^{i_j} \right) \varepsilon^1 u(b^{i_1},\,b^{i_2}) = 0,$$

where here $\{b^{i_1}\}_i = \{b^{i_2}\}_i = \{b^1,\,b^2\}$. It thus reduces to the expression (16). ♦

As the equation (16) should hold for any set of real coefficients $\{\alpha_j^{i_j}\}_{ij}$, we can express it in a stronger form by defining mappings $\{p_m\}_m$ such that

$$\begin{aligned}
p_m \colon \{1,\,\cdots,\,k+1\} &\longrightarrow \{1,\,\cdots,\,d\} \times \cdots \times \{1,\,\cdots,\,d\} \\
(1,\,\cdots,\,k+1) &\longmapsto (p_m(1),\,\cdots,\,p_m(k+1)).
\end{aligned}$$

Seeing $b_j^{i_j}$ through the relation $b_j^{i_j} = b^{p_m(j)}$, we can define the set $\mathcal{S} = \{p_m\}_m$ of cardinal $|\mathcal{S}| = (k+1)!$ as being the set of all the possibilities for defining the vector $(b_1^{i_1},\,\cdots,\,b_{k+1}^{i_{k+1}})$. The strong version of (16) thus comes down to finding all the functions $u$ belonging to $(\mathbb{P}_k(K))^d$ such that

$$\forall p \in \mathcal{S}, \quad \varepsilon^k u(b^{p(1)},\,\cdots,\,b^{p(k+1)}) = 0. \tag{17}$$

***Remark.*** For any $m \in [\![1, (k+1)!]\!]$, the set $\{p_m(j)\}_j$ do not necessarily contains unique elements. Indeed, by definition of the differential operator $\mathrm{D}^k$, each directional vector $\chi_i$, $i \in [\![1, k+1]\!]$ can be any vector of $\mathbb{R}^d$. Thus, for any vector $\chi_i$, $i \in [\![1, k+1]\!]$, any selection of a basis vector $b^j$, $j \in [\![1, d]\!]$ is possible regardless the other basis that have been selected for the other directional vectors $\chi_j$, $j \neq i$. ▲

***Example.*** *Continued* The functions $\{p_m\}_m$ read

$$\begin{aligned} p_i \colon \{1,\, 2\} &\longrightarrow \{1,\, 2\} \times \{1,\, 2\} \\ (1,\, 2) &\longmapsto (p_i(1),\, p_i(2)), \end{aligned}$$

and the set $\mathcal{S}$ is given by

$$\begin{aligned} \mathcal{S} =& \{p_m(1,\, 2)\}_m \\ =& \{(1,\, 2) \mapsto (1,\, 2), \\ & (1,\, 2) \mapsto (2,\, 1), \\ & (1,\, 2) \mapsto (1,\, 1), \\ & (2,\, 2) \mapsto (2,\, 2)\}\,. \end{aligned}$$

Thus, the problem (17) reduces to find all the functions $u \in (\mathbb{P}_1(K))^2$ such that

$$\forall i \in [\![1,\, 4]\!], \quad \varepsilon^1 u(b^{p_m(1)},\, b^{p_m(2)}) = 0.$$

Plugging the expressions of $\{p_m\}_m$ in the above equation, we obtain the system

$$\begin{cases} \varepsilon^1 u(b^1,\, b^2) = 0 \\ \varepsilon^1 u(b^2,\, b^1) = 0 \\ \varepsilon^1 u(b^1,\, b^1) = 0 \\ \varepsilon^1 u(b^2,\, b^2) = 0, \end{cases}$$

which coincides with the weakly expressed form obtained in the second example of this section. As a side comment, note that by the insensitivity of the differential operator $\mathrm{d}^1(\chi_1,\, \chi_2)$ on the order of $(\chi_1,\, \chi_2)$, the two first relations are identical. ♦

Therefore, there is only left to express the terms $\varepsilon^1 u(b^{p(1)}, \cdots, b^{p(k+1)})$ for obtaining a kernel parametrisation of the operator $\varepsilon^k$. Let us start by simplifying the Fréchet differentials $\mathrm{d}^k u(b^{p(1)}, \cdots, b^{p(k)})$ for any mapping $p \in \{\mathcal{S}\}$.

To this sake, we denote by $\alpha_i^m$ the number of occurrences of the vector $b^i$ in the set $(b^{p_m(1)}, \cdots, b^{p_m(k)})$ and define the directional vector

$$\chi^m = (b^{p_m(1)}, \cdots, b^{p_m(k)}).$$

***Note.*** By construction, $|\alpha^m| = |(\alpha_1^m, \cdots, \alpha_d^m)| = \sum_{i=1}^d \alpha_i^m = k$. ▲

Up to reordering the directional vectors, we can reduce $\chi^m$ to the expression

$$\chi^m = (\underbrace{b^1, \cdots, b^1}_{\alpha_1^m \text{ times}}, \cdots, \underbrace{b^d, \cdots, b^d}_{\alpha_d^m \text{ times}}),$$

leading to

$$\mathrm{d}_l^k u(b^{p(1)}, \cdots, b^{p(k)}) = \mathrm{d}_l^k u(\underbrace{b^1, \cdots, b^1}_{\alpha_1^m \text{ times}}, \cdots, \underbrace{b^d, \cdots, b^d}_{\alpha_d^m \text{ times}}).$$

Taking benefit of $\{b^i\}_i$ being the canonical basis functions, we have $\left(\frac{\partial u_l}{\partial x}\right) \cdot b^i = \frac{\partial u_l}{\partial x_i}$ and we can further derive

$$\begin{aligned}
\mathrm{d}_l^k u(\chi) &= \mathrm{d}_l^k u(\underbrace{b^1, \cdots, b^1}_{\alpha_1^m \text{ times}}, \cdots, \underbrace{b^d, \cdots, b^d}_{\alpha_d^m \text{ times}}) \\
&= \bigcirc_{j=1}^{k} \frac{\partial}{\partial x} (u_l) \cdot b^{p(m)} \\
&= \bigcirc_{j=2}^{k} \frac{\partial}{\partial x} \left( \frac{\partial u_l}{\partial x} \cdot b^{p(1)} \right) \cdot b^{p(m)} \\
&= \bigcirc_{j=2}^{k} \frac{\partial}{\partial x} \left( \frac{\partial u_l}{\partial x_{p(1)}} \right) \cdot b^{p(m)} \\
&= \bigcirc_{j=3}^{k} \frac{\partial}{\partial x} \left( \frac{\partial^2 u_l}{\partial x_{p(2)} \partial x_{p(1)}} \right) \cdot b^{p(m)},
\end{aligned}$$

which ultimately reduces to

$$\mathrm{d}_l^k u(\chi) = \frac{\partial^{|\alpha^m|} u_l}{\partial x_1^{\alpha_1^m} \cdots \partial x_d^{\alpha_d^m}}.$$

Thus, considering all the components of the function $u$ it comes

$$\mathrm{d}^k u = \frac{\partial^{|\alpha^m|}}{\partial x_1^{\alpha_1^m} \cdots \partial x_d^{\alpha_d^m}} \begin{pmatrix} u_1 \\ \vdots \\ u_d \end{pmatrix},$$

finally leading to the following expression of $\varepsilon^k$.

**Definition 3.4** Expression of $\varepsilon^k$ along the canonical directional vectors

Given some mapping $p$: $\{1, \cdots, k+1\} \mapsto \{1, \cdots, d\}^{k+1}$, it holds

$$\begin{aligned} \varepsilon^k \quad &: \mathbb{R}^d \times \cdots \times \mathbb{R}^d &&\longrightarrow \quad \mathbb{R} \\ &\left(b^{p(1)}, \cdots, b^{p(k+1)}\right) &&\longmapsto \quad \frac{1}{k+1} \sum_{j=1}^{k+1} \frac{\partial^k u_{p(j)}}{\partial x_1^{\alpha_1} \cdots \partial x_d^{\alpha_d}}, \end{aligned} \tag{18}$$

where $\alpha_i$ represents the number of occurrences of $i$ in the set $\mathfrak{s} = \{p(l),\, l \in \{1, \cdots, j-1, j+1, \cdots, k+1\}\}$, with $|\alpha| = k$.

Thus, expressing the kernel of $\varepsilon^k$ comes down to find every $u \in (\mathbb{P}_k(K))^d$ verifying

$$\sum_{j=1}^{k+1} \frac{\partial^k u_{p_m(j)}}{\partial x_1^{\alpha_1^m} \cdots \partial x_d^{\alpha_d^m}} = 0 \tag{19}$$

for any of the mappings $p_m \in \mathcal{S}$. Let us now describe the space determined by this parametrisation.

***Example.*** Let us consider the case $d = 2$, $k = 1$ and express the parametrisation of the kernel of $\varepsilon^k$ from the expression (18).

As shown before, the set $\mathcal{S}$ contains four mappings $\{p_i\}_i$, leading to four relations. We detail the expression (18) for the second of them and only state the computational outcome for the three others.

Let us first express the set $\mathfrak{s}_j$ and the terms $\{\alpha_i^j\}_{ij}$. We have $p_2(1, 2) = (2, 1)$. Therefore, it comes by definition

$$\mathfrak{s}_1 = \{p(2)\} = \{1\} \qquad \text{and} \qquad \mathfrak{s}_2 = \{p(1)\} = \{2\}.$$

Counting the respective occurrences of 1 and 2 within those sets leads to

$$\begin{aligned} \alpha_1^1 &= 1 \qquad \text{and} \qquad \alpha_2^1 = 0, \\ \alpha_2^1 &= 0 \qquad \text{and} \qquad \alpha_2^2 = 1. \end{aligned}$$

The relation (18) expressed from $p_2$ then becomes

$$\begin{aligned} \varepsilon^k \quad &: \mathbb{R}^d \times \cdots \times \mathbb{R}^d &&\longrightarrow \quad \mathbb{R} \\ &\left(b^2, \cdots, b1\right) &&\longmapsto \quad \frac{1}{2} \frac{\partial u_{p(1)}}{\partial x_1^{\alpha_1^1} \partial x_2^{\alpha_2^1}} + \frac{1}{2} \frac{\partial u_{p(2)}}{\partial x_1^{\alpha_1^2} \partial x_2^{\alpha_2^2}}, \end{aligned}$$

and leads to the relation

$$\frac{1}{2}\frac{\partial u_2}{\partial x_1} + \frac{1}{2}\frac{\partial u_1}{\partial x_2} = 0.$$

In a similar way, we retrieve the three other relations leading to the parametrisation of the kernel of $\varepsilon^1$. Presented respectively in the order corresponding to the mappings $p_1$, $p_2$, $p_3$, $p_4$, it reads

$$\begin{cases} \frac{1}{2}\frac{\partial u_1}{\partial x_2} + \frac{1}{2}\frac{\partial u_2}{\partial x_1} = 0 \\ \frac{1}{2}\frac{\partial u_2}{\partial x_1} + \frac{1}{2}\frac{\partial u_1}{\partial x_2} = 0 \\ \frac{\partial u_1}{\partial x_1} = 0 \\ \frac{\partial u_2}{\partial x_2} = 0. \end{cases}$$

We retrieve the system corresponding to the strong form of the expression of $\varepsilon^k$ obtained in the second example of this section. ♦

**Space generated from the kernel of $\varepsilon^k$** Now that we could determine the kernel's parametrisation, let us find a convenient expression for the generated space.

In all of this section, we keep in mind that the only functions $u$ that are considered belong to the vectorial space $(\mathbb{P}_k(K))^d$, which can also be written under the direct sum

$$(\mathbb{P}_k(K))^d = (\mathbb{P}_{k-1}(K))^d \oplus \left(\mathbb{P}_{[k]}(K)\right)^d \quad (20)$$

where $\mathbb{P}_{[k]}(K)$ is the space of polynomials having for degree exactly $k$. The kernel of the operator $\varepsilon^k$ acting on functions $(\mathbb{P}_k(K))^d$ can then be characterised by the union of the kernels of the operator $\varepsilon^k$ acting either on functions of $(\mathbb{P}_{k-1}(K))^d$ or on functions of $\left(\mathbb{P}_{[k]}(K)\right)^d$.

For characterising the former sub - kernel, it can be observed that for any integers $i$, $l \in [\![1,\, d]\!]$ there will always exist a mapping $p \in \mathcal{S}$ such that

$$p(1, \cdots, k+1) = (i, \cdots, i, l).$$

There, the multi - index $\alpha$ writes $\alpha = (0, \cdots, 0, k, 0, \cdots, 0)$ with $k$ being in the $i^{\text{th}}$ position, reducing the denominator of (19) to $\partial x_1^{\alpha_1} \cdots \partial x_d^{\alpha_d} = \partial x_i^k$.

The structure of the equation

$$\sum_{j=1}^{k+1} \frac{\partial^k u_{p(k+1)}}{\partial x_1^{\alpha_1} \cdots \partial x_d^{\alpha_d}} = 0$$

then comes down to

$$\frac{\partial^k u_l}{\partial x_i^k} = 0.$$

Thus, as this observation is valid for any $k \in \mathbb{N}$, any integer $d \geq 2$ and any integer $i \in [\![1,\, d]\!]$, there always exists $d$ mappings $p \in \mathcal{S}$ leading to the relations

$$\frac{\partial^k u_l}{\partial x_i^k} = 0 \quad \forall i,\, l \in [\![1,\, d]\!].$$

As a consequence, we get that all the functions $u \in (\mathbb{P}_{k-1}(K))^d \subset (\mathbb{P}_k(K))^d = (\mathbb{P}_{k-1}(K))^d \oplus \big(\mathbb{P}_{[k]}(K)\big)^d$ belong to the kernel of $\varepsilon^k$.

To characterise the latter sub - kernel, we should find the functions $u$ living in $\big(\mathbb{P}_{[k]}(K)\big)^d$ such that $\varepsilon^k u = 0$. To this aim, we paraphrase the work of Nédélec [13] and define the space $S^k$ of homogeneous polynomials of degree $k$ contained in $V$, that is

$$S^k = \left\{ u \in \big(\mathbb{P}_{[k]}(K)\big)^d \text{ such that } \sum_{i=1}^{d} x_i u_i \equiv 0 \right\}.$$

So defined, it holds [13] the following sub - kernel characterisation.

**Theorem 3.5** Latter sub - kernel characterisation

For any $u \in \big(\mathbb{P}_{[k]}(K)\big)^d$, it holds

$$u \in S^k \quad \Leftrightarrow \quad \varepsilon^k u \equiv 0$$

***Proof.*** *We briefly recall the arguments of [13].*

• We first show that $\{u \in \big(\mathbb{P}_{[k]}(K)\big)^d \text{ and } \varepsilon^k u \equiv 0\} \Rightarrow u \in S^k$. There, by definition of $u$, it is enough to show that $\varepsilon^k u \equiv 0 \Rightarrow \sum_{i=1}^d x_i u_i = 0$.

By assumption, we have for any $\chi = (\chi_1,\, \cdots,\, \chi_{k+1}) \in \mathbb{R}^{d\times(k+1)}$ and for any $x \in \mathbb{R}^d$,

$$\varepsilon^k u(x)(\chi) = 0.$$

Thus, taking in particular $\chi = (x, \cdots, x)$ for any $x \in \mathbb{R}^d$, it comes

$$\varepsilon^k u(\underbrace{x, \cdots, x}_{k+1 \text{ times}}) = 0,$$

which by definition reads

$$\frac{1}{k+1} \sum_{j=1}^{k+1} x \cdot \mathrm{d}^k u(\underbrace{x, \cdots, x}_{k \text{ times}}) = 0$$

and reduces to

$$x \cdot \mathrm{d}^k u(\underbrace{x, \cdots, x}_{k \text{ times}}) = 0.$$

In addition, since $u$ is homogeneous, we derive by the Euler's identity

$$\mathrm{d}^k u(\underbrace{x, \cdots, x}_{k \text{ times}}) = k!\, u.$$

Plugging it in the relation on $\varepsilon^k$, we retrieve for any $x \in \mathbb{R}^d$

$$\varepsilon^k u(\underbrace{x, \cdots, x}_{k+1 \text{ times}}) = 0 \quad \Rightarrow \quad \mathrm{d}^k(\underbrace{x, \cdots, x}_{k \text{ times}}) \cdot x = 0,$$

and finally

$$k!\, x \cdot u = 0.$$

Thus, $\varepsilon^k u \equiv 0 \;\Rightarrow\; \sum_{i=1}^d x_i u_i = 0$ for any $x \in \mathbb{R}^d$. And as $u \in \big(\mathbb{P}_{[k]}(K)\big)^d$, we get by definition that $u \in S^k$.

• Conversely, let us now show that $u \in S^k \;\Rightarrow\; \{u \in \big(\mathbb{P}_{[k]}(K)\big)^d$ and $\varepsilon^k u = 0\}$. As before, by definition of $s^K$ it is enough to show that $u \in S^k \;\Rightarrow\; \varepsilon^k u = 0\}$.

By assumption, we have that for any $x \in \mathbb{R}^d$, $\sum_{i=1}^d x_i u_i = x \cdot u = 0$. Therefore, deriving $\varepsilon^k u$ in the particular direction $\chi = (x, \cdots, x)$ yields

$$\begin{aligned} \varepsilon^k u(\underbrace{x, \cdots, x}_{k+1 \text{ times}}) &= \frac{1}{k+1} \sum_{j=1}^{k+1} x \cdot \mathrm{d}^k u(\underbrace{x, \cdots, x}_{k \text{ times}}) \\ &= x \cdot \mathrm{d}^k u(x, \cdots, x), \end{aligned}$$

which as $u$ is a polynomial of order exactly $k$ reduces to

$$\varepsilon^k u(\underbrace{x, \cdots, x}_{k+1 \text{ times}}) = 0.$$

Thus, for all $x \in \mathbb{R}^d$, $\sum_{i=1}^d x_i u_i = 0 \Rightarrow \varepsilon^k u(x, \cdots, x) = 0$. Expressing $x$ on the canonical basis of $\mathbb{R}^d$ so that $x = \sum_{i=1}^d x_i b^i$, we can develop the term $x \cdot u$ as $x \cdot u = \sum_{i=1}^d x_i b^i u_i$. Being is an homogeneous polynomial for which the term $\varepsilon^k u$ is identically vanishing, we can extend the above equality to any direction $\chi$, concluding the proof. We refer to [13] for further details.

■

Let us now merge the characterisation of the two sub - kernels. By considering $u \in (\mathbb{P}_k(K))^d = (\mathbb{P}_{k-1}(K))^d \oplus \left(\mathbb{P}_{[k]}(K)\right)^d$, it has been shown above that

$$\begin{cases} (\mathbb{P}_{k-1}(K))^d \subset \operatorname{Ker}(\varepsilon^k) \\ \forall\, u \in \left(\mathbb{P}_{[k]}(K)\right)^d_, \ \displaystyle\sum_{i=1}^d x_i u_i \equiv 0 \Leftrightarrow \varepsilon^k u \equiv 0. \end{cases}$$

Thus, it comes the following kernel characterisation.

**Property 3.6** Kernel of $\varepsilon^k$

$$\operatorname{Ker}(\varepsilon^k) = (\mathbb{P}_{k-1}(K))^d \oplus S^k \tag{21}$$

**The simplicial Raviart – Thomas space** With respect to the above computations and the *Definition 3.3*, the discretisation space is defined as follows.

$$V = (\mathbb{P}_{k-1}(K))^d \oplus S^k$$

Recalling that $S^k$ is the space of polynomials of degree k whose gradient are homogeneous polynomials of degree $k-1$, we get

$$u \in S^k \quad \Rightarrow \quad u = (k+1) x \nabla u$$

by Euler's identity. And since $\nabla p$ belongs to $\mathbb{P}_{[k-1]}(K)$, we can identify $S^k$ with $x\,\mathbb{P}_{[k-1]}$. Therefore, up to shifting the indices it comes the following classical definition of the Raviart – Thomas simplicial spaces.

**Definition 3.7** Raviart – Thomas simplicial spaces.

$$RT_k = (\mathbb{P}_k)^d \oplus x\,\mathbb{P}_{[k]}$$

We conclude this section by recalling the examples of [13].

***Example.*** *Two dimensional cases* ($d = 2$) We consider the two lowest orders differential operators $\varepsilon^k$, $k = 1, 2$ in the two dimensional case and derive the space generated by the kernel of the operator $\varepsilon^k$.

• The kernel of the $\varepsilon^1$ operator is generated by the following relations.

$$\varepsilon^1 = 0 \quad \Leftrightarrow \quad \begin{cases} \frac{\partial u_1}{\partial x_1} = 0 \\ \frac{\partial u_2}{\partial x_2} = 0 \\ \frac{\partial u_1}{\partial x_2} + \frac{\partial u_2}{\partial x_1} = 0 \end{cases}$$

When coupled with the constraint $u \in \mathbb{P}_1^2(K)$ coming from the wish of defining a finite subset of $H(\mathrm{div}, K)$ and formulated through (3.1.1), those equations describe the space

$$V = \begin{pmatrix} a \\ b \end{pmatrix} \oplus c \begin{pmatrix} -x_2 \\ x_1 \end{pmatrix},$$

where $a$, $b$, $c$ are constants. In particular, $\big(a,\, b\big)^{\mathrm{T}} \subset (\mathbb{P}_0)^2$, space which is generated by the two first equations. The third equation generates $S^1$, which is reformulated as $x\,\mathbb{P}_0(K)$. Thus, we retrieve directly the classical Raviart – Thomas space. Meaningly,

$$V = (\mathbb{P}_0)^2 \oplus S^1 = (\mathbb{P}_0)^2 \oplus x\,\mathbb{P}_0(K) = RT_0(K).$$

• Similarly, the kernel of the $\varepsilon^2$ operator is generated by the following relations.

$$\varepsilon^2 = 0 \Leftrightarrow \begin{cases} \frac{\partial^2 u_1}{\partial x_1^2} = 0 \\ \frac{\partial^2 u_2}{\partial x_2^2} = 0 \\ \frac{\partial^2 u_1}{\partial x_2^2} + \frac{2\partial^2 u_2}{\partial x_1 \partial x_2} = 0 \\ \frac{\partial^2 u_2}{\partial x_1^2} + \frac{2\partial^2 u_1}{\partial x_1 \partial x_2} = 0 \end{cases}$$

Coupled with the condition $u \in (\mathbb{P}_1(K))^2$, those equations generate $RT_1(K)$. Thus we can write:

$$V = (\mathbb{P}_1)^2 \oplus S^2 = RT_1(K).$$

♦

***Remark.*** The space $RT_k$ is conformity - ready. Provided a good choice of degrees of freedom , the elements defined on this space are conformal due to the constraint on the $\varepsilon^k$ operator, which drives the conformity behaviour on the divergence. ▲

## 3.2 Structure of the RT – space

When working with the $\mathbb{P}_k$ polynomial spaces, we recall that the Raviart – Thomas space on a reference element $K \subset \mathbb{R}^{[d]}$ can be written

$$RT_k(K) = (\mathbb{P}_k(K))^d \oplus x\mathbb{P}_{[k]}(K), \tag{22}$$

with $x \in \mathbb{R}^d$ and where $\mathbb{P}_{[k]}$ is the polynomials space of degree exactly $k$.

As this space is designed to provide a $H(\mathrm{div}, K)$ – conformal finite element discretisation framework, we also need to consider a discretisation space on the boundary of $K$ so that degrees of freedom can be designed directly on the boundary. In the case of simplicial elements, this space is set as $\mathcal{R}_k(\partial K)$.

## 3.3 Properties of the $RT_k(K)$ space

### Dimension

Since $RT_k(K)$ is built as a direct sum between two spaces, we can write

$$\dim RT_k(K) = d \dim \mathbb{P}_k(K) + \dim \mathbb{P}_{[k]}(K). \tag{23}$$

And by the previous algebra results, we have

$$\dim \mathbb{P}_k(K) = \binom{k+d}{k}.$$

Thus,
$$\begin{aligned}\dim \mathbb{P}_{[k]}(K) &= \dim \mathbb{P}_k(K) - \dim \mathbb{P}_{k-1}(K) \\ &= \binom{k+d}{k} - \binom{k+d-1}{k-1} \\ &= \binom{k+d}{k} - \left[\binom{k+d}{k} - \binom{k+d-1}{k}\right],\end{aligned}$$

by Pascal's formula, reducing to

$$\dim \mathbb{P}_{[k]}(K) = \binom{k+d-1}{k}.$$

Therefore, the dimension reads; $\dim RT_k(K) = d\binom{k+d}{k} + \binom{k+d-1}{k}$.

**Property 3.8** Raviart – Thomas space's dimension for simplicial shapes

$$\dim RT_k(K) = \frac{d(k+d)!}{k!d!} + \frac{(k+d-1)!}{k!(d-1)!} \tag{24}$$

From there, we have two ways to develop; either on $k$, writing $(k+d)! = k!(k+1)\cdots(k+d)$ or on $d$, writing $(k+d)! = d!(d+1)\cdots(d+k)$. To be able to use the formula for any $k$, we choose to develop on $k$. It yields

$$\begin{aligned}\dim RT_k(K) &= \frac{d(1)\cdots(k)(k+1)\cdots(k+d)}{(1)\cdots(k)(1)\cdots(d)} + \frac{(1)\cdots(k-1)(k)\cdots(k+d-1)}{(1)\cdots(k)(1)\cdots(d-1)} \\ &= \frac{d}{d!}(k+1)\cdots(k+d) + \frac{1}{(d-1)!}(k+1)\cdots(k+d-1).\end{aligned}$$

In practice, the dimension for lowest orders Raviart – Thomas spaces reads

| k \ d | 2 | 3 | 4 |
|---|---|---|---|
| 0 | 2 + 1 = 3 | 3 + 1 = 4 | 4 + 1 = 5 |
| 1 | 6 + 2 = 8 | 12 + 3 = 15 | 20 + 4 = 24 |
| 2 | 12 + 3 = 15 | 30 + 6 = 36 | 60 + 10 = 70 |
| 3 | 20 + 4 = 2 | 60 + 10 = 70 | 140 + 20 = 160 |

where we left in black the dimension of $\mathbb{P}_k$ times the domain's spatial dimension and set in grey the dimension of $\mathbb{P}_{[k]}$. The dimension of $RT_k(K)$ is their sum. In view of later use, we recall the dimensions formulas of the $RT_k(K)$ spaces for the two and three dimensional cases.

***Example.*** *Dimension of simplicial $RT_k(K)$ spaces for $d = 2$ and $d = 3$.*

- In the case $d = 2$ we have

$$\begin{aligned}\dim RT_k(K) &= \frac{2}{2}(k+1)(k+2) + (k+1) \\ &= (k+1)(k+2+1) \\ &= (k+1)(k+3),\end{aligned}$$

- while for $d = 3$ it comes

$$\begin{aligned}\dim RT_k(K) &= \frac{3}{3!}(k+1)(k+2)(k+3) + \frac{1}{2}(k+1)(k+2) \\ &= \frac{1}{2}(k+1)(k+2)(k+3+1) \\ &= \frac{1}{2}(k+1)(k+2)(k+4).\end{aligned}$$

♦

### Divergence properties

The Raviart – Thomas simplicial space enjoy the following property [4].

**Property 3.9** Divergence properties

For any element $K$, we have for any $q \in RT_k(K)$:

$$\begin{cases} \operatorname{div} q \in \mathbb{P}_k(K) \\ q \cdot n|_{\partial K} \in \mathcal{R}_k(\partial K) \end{cases} \tag{25}$$

Moreover, the divergence is onto from $RT_k(K)$ onto $\mathbb{P}_k(K)$.

***Proof.*** Any $q \in RT_k(K)$ can be written $q = q_0 + xp_k$ with some $q_0 \in (\mathbb{P}_k(K))^d$ and $p_k \in \mathbb{P}_{[k]}(K)$. Then, still by Euler's identity we have that $\operatorname{div} q = \operatorname{div} q_0 + \operatorname{div}(x\, p_k)$ is a polynomial of degree $k$. Thus, $\operatorname{div} q \in \mathbb{P}_k(K)$.

Let us now consider any side $f$ of the element $K$ and set $n_f$ its normal. By the properties of the normal vector, $x \cdot n_f$ is constant for any $x$ belonging to the side $f$. In other words, we have that for any face $f$ of the boundary $\partial K$, the application

$$\begin{aligned} f &\longrightarrow \mathbb{R} \\ x &\longmapsto x \cdot n_f \end{aligned}$$

is constant. Thus, for any point $x$ on the side $f$, $q \cdot n = q_0 \cdot n + p_k(x \cdot n_f)$ and $x \cdot n_f \equiv c$ for some constant $c \in \mathbb{R}$. Recalling that $q_0 \in \mathbb{P}_k(K)$, $p_k \in \mathbb{P}_{[k]}(K)$ and that the definition of $\mathcal{R}_k(\partial K)$ does not require continuity across the edges, we obtain that $q \cdot n|_{\partial K} \in \mathcal{R}_k(\partial K)$.

A proof of the surjectivity of the divergence operator from $RT_k(K)$ onto $\mathbb{P}_k(K)$ can be found in [2, *Paragraph 2.3.3*] and will not be repeated.

■

***Note.*** The structure of the $RT_k(K)$ space has been specifically designed so that the divergence relations (25) hold. ▲

Now that as we have the property 3.9, we can tune the degrees of freedom on $K$, $RT_k(K)$ and $\mathcal{R}_k(K)$ to form a H(div,$K$) conformal element.

## 3.4 Construction of $RT_k(K)$ elements

We start by defining the classical set of degrees of freedom of $RT_k(K)$ elements and explain a way to construct their corresponding basis functions.

### 3.4.1 Set up of degrees of freedom

To completely define the element we need to set up degrees of freedom, a set $\{\sigma_i\}_i$ of linear forms on $K$ that uniquely determines the element. The dimension of the set must match the one of $RT_k(K)$, and each of the degrees of freedom should map $RT_k(K)$ to $\mathbb{R}$. Their values will then entirely characterise some function living in a $RT_k(K)$ space.

As we want to build H(div) – conformal elements, we will select degrees of freedom whose shape will enforce conformity. With this in mind, we start by noticing the following.

- For d=2 it holds from (23)

$$\begin{aligned}
\mathbf{dim}\,\boldsymbol{RT_k(K)} &= 2\,\dim(\mathbb{P}_k(K)) + \dim \mathbb{P}_{[k]}(K)\\
&\boldsymbol{= \mathbf{dim}\,(\mathbb{P}_k(K)^2) + \mathbf{dim}\,\mathbb{P}_{[k]}(K)}\\
&= (k+1)(k+2) + (k+1)\\
&= (k+1)(k+3)\\
&= k(k+1) + 3(k+1)\\
&\boldsymbol{= \mathbf{dim}\,(\mathbb{P}_{k-1}(K)^2) + 3\,\mathbf{dim}\,\mathbb{P}_k(f)}\\
&= \dim(\mathbb{P}_{k-1}(K)^2) + \dim \mathcal{R}_k(\partial K), \quad \text{by the relation (1),}
\end{aligned}$$

- whereas for d=3 we have:

$$\begin{aligned}
\mathbf{dim}\,\boldsymbol{RT_k(K)} &= 3\dim(\mathbb{P}_k(K)) + \dim \mathbb{P}_{[k]}(K)\\
&\boldsymbol{= \mathbf{dim}\,(\mathbb{P}_k(K)^3) + \mathbf{dim}\,\mathbb{P}_{[k]}(K)}\\
&= \frac{1}{2}(k+1)(k+2)(k+3) + \frac{1}{2}(k+1)9k+20\\
&= \frac{1}{2}(k+1)(k+2)(k+4)\\
&= \frac{1}{2}k(k+1)(k+2) + \frac{4}{2}(k+1)(k+2)\\
&= \frac{3}{3!}k(k+1)(k+2) + 4\frac{(k+1)(k+2)}{2}\\
&= 3\binom{k-1+3}{k-1} + 4\binom{k+3-1}{k}\\
&\boldsymbol{= \mathbf{dim}\,(\mathbb{P}_{k-1}(K)^3) + 4\,\mathbf{dim}\,\mathbb{P}_k(f)}.
\end{aligned}$$

- More generally for a dimension $d$ we have:

$$\begin{aligned}
\mathbf{dim}\, \boldsymbol{RT_k(K)} &= d \dim(\mathbb{P}_k(K)) + \dim \mathbb{P}_{[k]}(K) \\
&\boldsymbol{= \mathbf{dim}\, (\mathbb{P}_k(K)^d) + \mathbf{dim}\, \mathbb{P}_{[k]}(K)} \\
&= \frac{d(k+d)!}{k!d!} + \frac{(k+d-1)!}{k!(d-1)!} \\
&= \frac{1}{k!d!}(d(k+d)! + d(k+d-1)!) \\
&= \frac{1}{k!(d-1)!}((k+d-1)!(1+k+d)) \\
&= \frac{d}{dk!(d-1)!}k(k+d-1) + \frac{1}{k!(d-1)!}(1+d)(k+d-1)! \\
&= d\binom{k-1+d}{k-1} + (d+1)\binom{k+d-1}{k} \\
&\boldsymbol{= \mathbf{dim}\, (\mathbb{P}_{k-1}(K)^d) + (d+1)\, \mathbf{dim}\, \mathbb{P}_k(f)}. \qquad (26)
\end{aligned}$$

We also notice that as $f$ is a hypersurface of $K$ contained in $\mathbb{R}^{d-1}$, $(1+d)\dim \mathbb{P}_k(f) = \dim \mathcal{R}_k(\partial K)$ for simplicial elements in dimension $d$. That means that for those elements, we can split the degrees of freedom up and build them from the two distinct spaces $\mathbb{P}_{k-1}(K)^d$ and $\mathcal{R}_k(\partial K)$. Both sets together will still provide a determination of an element over the $RT_k(K)$ space, provided that they are free from each other.

***Remark.*** Provided that the total dimension of the two sets matches the dimension of $RT_k(K)$, checking the unisolvence of the total set of degrees of freedom is equivalent to checking that they are two by two free. ▲

Keeping this split in mind, we can use the degrees of freedom defined on $\mathcal{R}_k(\partial K)$ to enforce the H(div) – conformity. Moment-based degrees of freedom are then designed.

**Definition 3.10** Degrees of Freedom in the simplicial case

For any $q \in RT_k(K)$, we set as degrees of freedom:

$$q \longmapsto \int_{\partial K} q \cdot n\, p_k \,\mathrm{d}\gamma(x), \quad \forall p_k \in \mathcal{R}_k(\partial K) \tag{27a}$$

$$q \longmapsto \int_{K} q \cdot p_{k-1} \,\mathrm{d}x, \quad \forall p_{k-1} \in (\mathbb{P}_{k-1}(K))^d \tag{27b}$$

where $\gamma$ represents the paths skimming the faces.

***Remark.*** In practice, it is enough to provide those degrees of freedom by taking $\{p_k\}$ as a basis of $\mathcal{R}_k(\partial K)$, and $\{p_{k-1}\}$ as a basis of $(\mathbb{P}_{k-1}(K))^d$. Classically, one can stay canonical and take:

$$q \longmapsto \int_f q \cdot n\, x^l \,\mathrm{d}x, \quad \forall l,\ 0 \le |l| \le k \text{ and for every } f \text{ face in } \partial K \tag{28a}$$

$$q \longmapsto \int_K q \cdot \bigodot_{i=1}^{d} \left(0, \cdots, 0, x_i^{l_i}, 0, \cdots, 0\right)^{\mathrm{T}} \mathrm{d}x, \quad {\scriptstyle \forall \{l_i\}_{i \in [\![1,\, d]\!]},\, 0 \le l_i \le k-1, \atop \sum_{i=1}^{d} l_i \le k-1,\, \forall j \in [\![1,\, d]\!]}, \tag{28b}$$

where $l = (l_1, \cdots, l_d)$ is a multi-index and where in the last row $x_i$ is in the $j^{\text{th}}$ position. The symbol $\bigodot$ denotes the Hadamard product. ▲

***Example.*** *Two dimensional case* For $d = 2$, the relations (28a) reduce to

$$q \longmapsto \int_f q \cdot n\, x^m \,\mathrm{d}x, \quad \forall m,\ 0 \le m \le k \text{ and for every } f \text{ edge in } \partial K \tag{29}$$

where here $i$ is a natural integer and $x$ lives in one dimension, while the relations (28b) reads

$$q \longmapsto \int_K q \cdot \begin{pmatrix} x^i\, y^j \\ 0 \end{pmatrix} \mathrm{d}x \quad \forall (i,\, j) \text{ s.t. } {\scriptstyle i,j \le k-1 \atop i+j \le k-1}$$

$$q \longmapsto \int_K q \cdot \begin{pmatrix} 0 \\ x^i\, y^j \end{pmatrix} \mathrm{d}x \quad \forall (i,\, j) \text{ s.t. } {\scriptstyle i,j \le k-1 \atop i+j \le k-1.}$$

♦

Those degrees of freedom are linear forms on $K$. Furthermore, it is well-known [16, 4] that it holds the following.

**Proposition 3.11**

$$\begin{cases} \int_{\partial K} q \cdot n\, p_k \,\mathrm{d}\gamma(x) = 0, & \forall p_k \in \mathcal{R}_k(\partial K) \\ \int_K q \cdot p_{k-1} \,\mathrm{d}x = 0, & \forall p_{k-1} \in (\mathbb{P}_{k-1}(K))^d \end{cases} \quad \Rightarrow \quad q \equiv 0.$$

Thus, since the dimension of $\{\sigma\}$ matches the one of $RT_k(K)$, and since $\operatorname{Ker}(\{\sigma\}) = \{0\}$, we obtain straightforwardly the following proposition.

**Proposition 3.12**

The set $\{\sigma\}$ is unisolvent for $RT_k(K)$ and $(K, RT_k(K), \{\sigma\})$ completely defines an element on $K$.

***Remark.***
- The function $q \in RT_k(K)$ is vectorial.
- The applications $q \longmapsto \int_{\partial K} q \cdot n\, p_k \mathrm{d}x$ are called normal degrees of freedom, whereas $q \longmapsto \int_K q \cdot p_{k-1} \mathrm{d}x$ are called internal degrees of freedom. ▲

***Note.*** Endowed with the set of degrees of freedom (27), the $RT_k(K)$ space is the smallest polynomial space within which the divergence maps $RT_k(K)$ to $\mathbb{P}_k(K)$. ▲

To conclude this paragraph, let us draw the classical representation of the Raviart – Thomas elements' degrees of freedom. In the simplicial case, it is given in the *Figure 5*, where the arrows represent the normal moments and the number the internal moments.

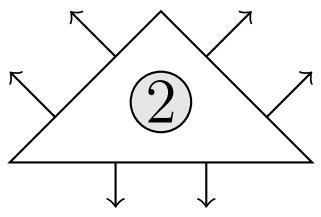

Fig. 5: $RT_2(K)$ element for $d = 2$ on triangle

### 3.4.2 Conformity

**Proposition 3.13**

The above-defined element $(K, RT_k(K), \{\sigma\})$ is H(div)–conformal.

***Proof.*** Here, the $H(\mathrm{div})$ conformity is equivalent to the continuity of the normal component at each interface. By definitions of the degrees of freedom and $\mathcal{R}_k(\partial K)$, the unknowns that lie on the faces are

$$\int_f q \cdot n\, p_k \mathrm{d}\gamma(x), \quad \forall f \in \partial K, \forall p_k \in \mathbb{P}_k(f).$$

Since $q \cdot n$ also belongs to $\mathbb{P}_k(f)$ as shown before, this element is H(div)–conformal.

■

***Remark.*** When working with a reference element, one has to be careful with the mesh orientation. Indeed, usual transformation to reference element will fail to preserve the $H(\text{div}, K)$ – conformity because of the normal impact on the degrees of freedom. One requires in general the Piola transformation, as described in the section 2.5.1. ▲

## 3.5 Associated basis functions for the 2D case

Even if the element $(K, RT_k(K), \{\sigma\})$ is already well defined through the degrees of freedom, it can be convenient to know the corresponding basis functions. We present here a construction method slightly adapted from [9].

The aim is to compute a basis of $RT_k(K)$ that matches the introduced degrees of freedom. We know that

- The basis functions are dual to the degrees of freedom.
- All the basis functions must be free from each other.
- Each one must lie in the $RT_k(K)$ space.
- The dimension of their set must match the one of $RT_k(K)$.

Thus, we first define any basis over the $RT_k(K)$ space, and then tune it to the preset degrees of freedom through their duality relationship.

### 3.5.1 Construction of a basis of the $RT_k(K)$ space

Although the method presented here works only for simplicial elements, it is described with general notations to make their later generalisation easier. Therefore, keep in mind that in this paragraph, $d = 2$.

For simplicity, we will work on the reference element presented in the *Figure 6.* Let us denote by $n$ the number of edges of the element (here $n = 3$). In all the following, the edges' normals $n_f$ are normalized with respect to the norm of the corresponding edge $f$, *i.e.* $n_f \leftarrow \frac{n_f}{||f||_2}$.

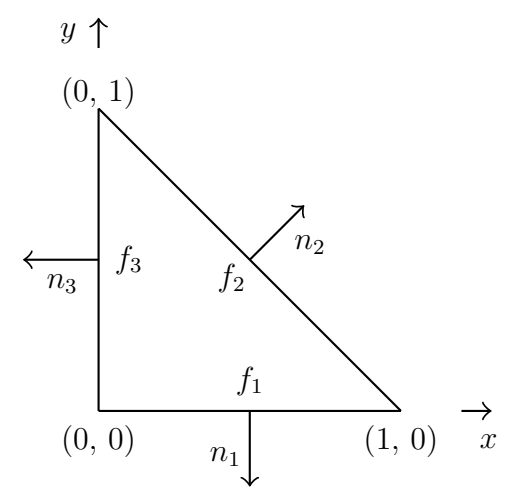

Fig. 6: Triangle of reference

As we are in the particular case of two dimensional simplicial elements, we can denote $\mathrm{x} = (x,\, y)^T$ and reduce the definition of $RT_k(K)$ to

$$RT_k(K) = (\mathbb{P}_k(K))^2 \oplus \mathrm{x}\,\mathbb{P}_{[k]}(K).$$

We want to define $\dim RT_k(K)$ functions that lie in $RT_k(K)$ and that are free from each other. As we aim to match them with either internal degrees of freedom or normal degrees of freedom, those basis functions will be split into two categories.

- $\mathcal{N}$, a free set of functions whose dimension is $\dim \mathcal{R}_k(\partial K) = (1 + d)\dim(\mathbb{P}_k(f))$ and where the functions lie in $RT_k(K)$. They should enforce the $H(\mathrm{div},\, K)$ – conformity.
- $\mathcal{I}$, a set of functions which is free strictly within $K$ and whose dimension is $d\dim(\mathbb{P}_{k-1}(K))$ and where the functions lie in $RT_k(K)$. They should preserve the $H(\mathrm{div},\, K)$ – conformity.

Each of $\mathcal{N}$ and $\mathcal{I}$ should be a free subset of $RT_k(K)$. Note that those functions do not lie exclusively in either $(\mathbb{P}_k(K))^d$ or $x\mathbb{P}_{[k]}(K)$, but rather in the full $RT_k(K)$ space. Indeed, we need $\mathcal{N} \cup \mathcal{I}$ to be a free set, and the imposed dimension on $\mathcal{N}$ and $\mathcal{I}$ prevents a direct split on the two subspaces $(\mathbb{P}_k(K))^d$ and $\mathrm{x}\mathbb{P}_{[k]}(K)$. We then have to build the functions conjointly on $RT_k(K)$, preserving the desired properties of dimension and freedom.

**Let us first generate $\mathcal{N}$.** We want to generate $\dim \mathcal{R}_k(\partial K)$ functions that are free from each other and that lie in $RT_k(K)$. As those functions will later be tuned to be dual to the set (27a), we follow the layout of those degrees of freedom and generate $(k+1) = \dim \mathbb{P}_k(f_i) = \dim(\mathcal{R}_k(\partial K))/n$ functions from each edge $f_i$, $i \in [\![1,\, n]\!]$.

To this end, we define the vector $e_i$ corresponding to any edge $f_i$ as

$$e_i = \begin{pmatrix} x + n_{ix} \\ y + n_{iy} \end{pmatrix},$$

where $n_{f_i} = (n_{ix},\, n_{iy})^{\mathrm{T}}$ is the $x$ and $y$ components vector of the normalized normal vector to the edge $f_i$. Then, we set one - dimensional Lagrangian functions of degree $k$ on each edge so that each Lagrangian set is generated from $(k+1)$ sample points lying strictly within the corresponding edge (*see the Figure 7*). This Lagrangian functions set will be denoted by $\{l_{im},\, i \in [\![1,\, n]\!],\, m \in [\![1,\, k+1]\!]\}$.

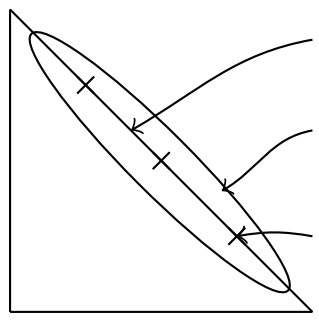

Fig. 7: Example of sampling points distribution

Multiplying each of those functions coordinate - wise with the previously introduced vectors $e_i$ forms a set that fulfil our requirements. Therefore, we set

$$\mathcal{N} = \{(x,\, y)^{\mathrm{T}} \mapsto l_{i,m}(x,\, y) e_i(x,\, y)\}_{\substack{i \in [\![1,\, n]\!] \\ m \in [\![1,\, k+1]\!]}}\,. \tag{30}$$

***Remark.*** In practice, the $(k+1)$ points are generated by a Gauss - Legendre distribution. However, any other sampling method is admissible and may change from one edge to another. ▲

***Remark.*** If we fix an homogeneous distribution method across the edges for the points generating the local Lagrangian functions, we can generate the functions only once and transform them to fit the other edges.

Indeed, we can orientate our domain and set a transformation from one Lagrangian set to another. To do so, we can define $x \mapsto l_m(x)$ the $m^{\text{th}}$ Lagrangian function generated by $(k+1)$ points on the segment $[0,\, 1]$, that is:

$$l_m(x) = \prod_{\substack{l = [\![1,\, k+1]\!] \\ l \neq m}} \frac{x - x_l}{x_l - x_m}$$

Then, by an affine transformation, the Lagrangian functions associated with the edge $f_i$ are given by:

$$\begin{cases} l_{i,m}(x,\, y) = l_m\left(\frac{1}{|\cos(\alpha_i)|} x\right), & \text{if } n_{x_i}, n_{y_i} \neq 0 \\ l_{i,m}(x,\, y) = l_m(y) & \text{if } n_{x_i} = 0 \\ l_{i,m}(x,\, y) = l_m(x) & \text{if } n_{y_j} = 0. \end{cases} \tag{31}$$

***Note.*** If $n_{x_i}, n_{y_i} \neq 0$, then $l_m\left(\frac{1}{|\cos(\alpha_i)|} x\right) = l_m\left(\frac{1}{|\sin(\alpha_i)|} y\right)$. Furthermore, note that for numerical integration reasons given in the *Remark 3*, the constant $1/|\cos(\alpha_i)|$ or $1/|\sin(\alpha_i)|$ could be dropped (*see the Figure 8 for an illustration*).

However, if one prefers to preserve the homogeneity of the sampling points distribution method across the edges, it should be kept.

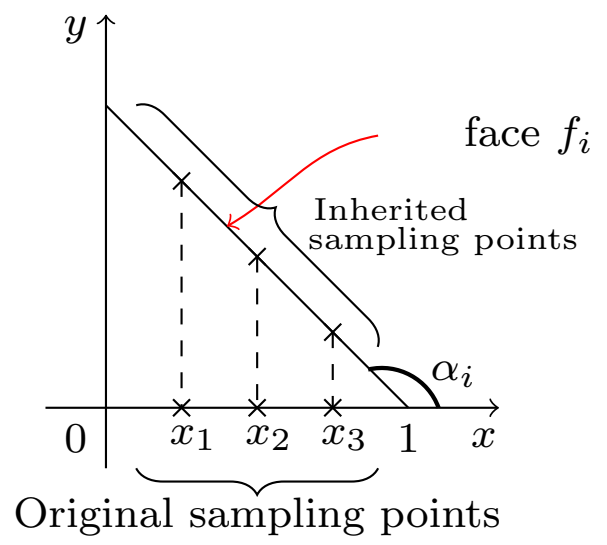

Fig. 8: Generating the functions $l_{im}$ at once as in (31) for $k = 2$ when the constant $1/|\cos(\alpha_i)|$is dropped. There, the functions $l_i(x)$ are generated on the horizontal edge. The Lagrangian functions emerging from the vertical edges are computed using the function defined on the horizontal edge with flipped variables, while the Lagrangian functions built from the last edge are just using the $x$ - variable spanning on the horizontal edge before using the definition of $l_i$. The behaviour of the retrieved function along the last edge is therefore given with respect to the $x$ - projection. ▲

***Remark 1.*** It may be convenient to order the indexing of the sampling points (and thus of the basis functions) with respect to the orientation of the element. Indeed, once tuned towards the degrees of freedom the basis functions will enjoy a global Lagrangian properties on the edges. Ordering the sampling points so that their indices are contiguous with respect to the element's orientation makes easier to recognize which function will value one on which point. Setting a single global permutation $\zeta$ over the sampling points in the definition of the basis function is enough to ensure a correct indexing (*see the Figure 9*). The set of all functions (30) would then turn into

$$\mathcal{N} = \{(x,\, y)^{\mathrm{T}} \mapsto l_{i,\zeta(m)}(x,\, y) e_i(x,\, y)\}_{\substack{i\in[\![1,\,n]\!]\\ m\in[\![1,\,k+1]\!]}}.$$

Furthermore, in the case of an homogeneous sampling across the edges, one can just set the permutation operator to $\zeta\colon i \mapsto k+1-i$ and change the first and second relations in (31) respectively to $l_{i,m}(x,\, y) = l_{k+1-m}\left(\frac{1}{|\cos(\alpha_i)|}x\right)$ and $l_{i,m}(x,\, y) = l_{k+1-m}(y)$. ▲

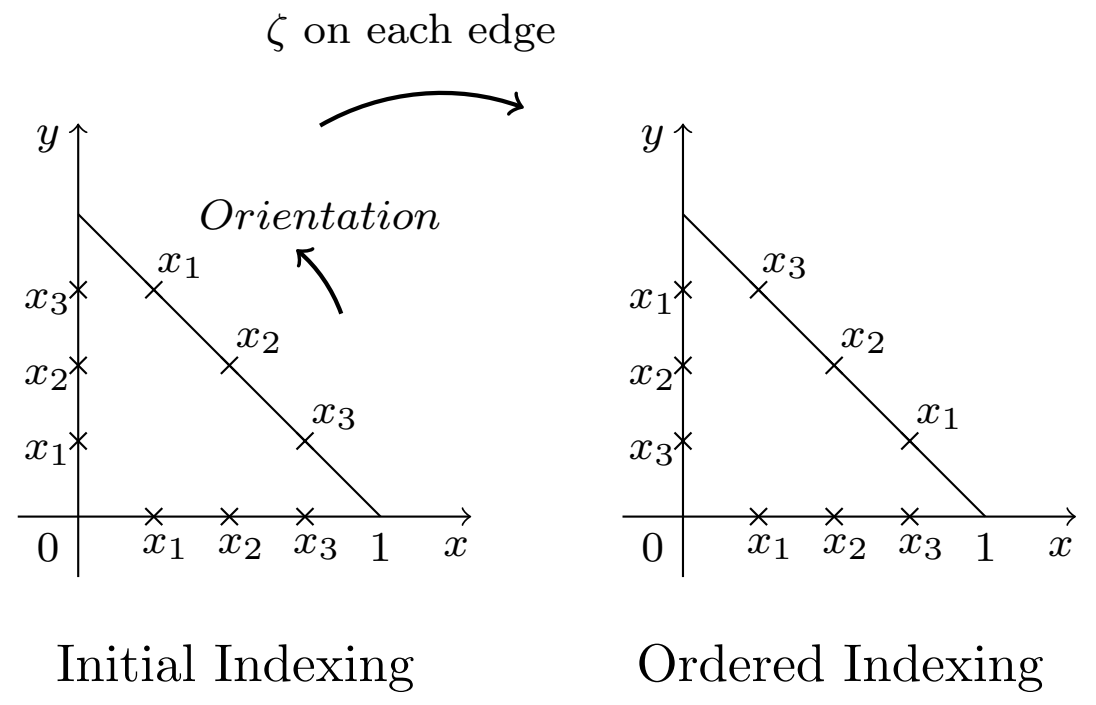

Fig. 9: Sampling points ordering in local indexing

We can immediately derive some properties.

**Property 3.14**

$\dim \mathcal{N} = (1+d)\dim \mathbb{P}_k(f) = \dim \mathcal{R}_k(\partial K)$

***Proof.*** We have $\dim \mathcal{N} = n \dim\{l_{i,m}\}_{m\in[\![1,\,k+1]\!]} = n(k+1)$. As in the simplicial case $n=3$, we retrieve $\dim \mathcal{N} = 3(k+1)$. Furthermore, since we are in the two-dimensional case, we have $3(k+1) = (d+1)(k+1)$ and $\dim \mathbb{P}_k(f) = (k+1)$. Thus, in the case of two dimensional simplicial elements, $\dim \mathcal{N} = (d+1)(k+1) = (k+1)\dim \mathbb{P}_k(f) = \dim \mathcal{R}_k(\partial K)$.

■

**Property 3.15**

If $q$ belongs to $\mathcal{N}$, then $q$ belongs to $RT_k(K)$.

***Proof.*** If $q$ belongs to $\mathcal{N}$, then one can write for all $(x,\,y)$ in $K$ that $q(x,\,y) = l_{i,m}(x,\,y)\begin{pmatrix} x+n_{ix} \\ y+n_{iy}\end{pmatrix}$ for given $i\in[\![1,\,n]\!]$ and $m\in[\![1,\,k+1]\!]$.

Thus, $q(x,\,y) = \begin{pmatrix} x \\ y\end{pmatrix}\underbrace{l_{i,m}(x,\,y)}_{\in\mathbb{P}_k(f_i)} + \underbrace{\begin{pmatrix} n_{ix} \\ n_{iy}\end{pmatrix}}_{\in\mathbb{R}^2}\underbrace{l_{i,m}(x,\,y)}_{\in\mathbb{P}_k(f_i)}$.

We note here that the term $(x,\,y)^{\mathrm{T}}$ only shifts by one the polynomial degree in each direction. Therefore, $q$ lies in a subspace of $(\mathbb{P}_k)^2 \cup \mathrm{x}\mathbb{P}_k$ and up to rewriting $\mathrm{x}\mathbb{P}_k = (x\mathbb{P}_{[k]}\cup\mathbb{P}_k)\times(y\mathbb{P}_{[k]}\cup\mathbb{P}_k)$, $(\mathbb{P}_k)^2\cup\mathrm{x}\mathbb{P}_k$ is contained in $(\mathbb{P}_k)^2\cup\mathrm{x}\mathbb{P}_{[k]}$. Thus, by definition $q\in RT_k(K)$.

■

**Proposition 3.16**

Let us assume that the sampling points method is injective. If the normals $\{n_f\}_f$ are two by two different, then the set $\mathcal{N}$ is free (linearly independent).

***Proof.*** For any $q \in \mathcal{N}$, $q$ can be written under the form $q = l_{i,m}e_i$ for some $i \in [\![1,\, n]\!]$ and $j \in [\![1,\, k+1]\!]$. We also know that the set $\{l_{i,m}\}_{i,m}$ is free and composed of functions mapping $\mathbb{R}^2$ to $\mathbb{R}$. Therefore, the set $\{al_{i,m}\}_{i,m}$ is also free for any constant $a \in \mathbb{R}$.

Furthermore, as $n = 3$ and as the normals are two by two different, we have that $\{e_i\}_i$ is a free set of vectors. Thus, for constants $\alpha_i \in \mathbb{R}$, $i \in [\![1,\, n]\!]$, it holds for every point $(x,\, y) \in K$:

$$\sum_{i=1}^{n} \alpha_i e_i(x,\, y) = 0 \quad \Rightarrow \quad \forall i \in [\![1,\, n]\!],\, \alpha_i = 0.$$

Therefore, for any fixed point $(x,\, y) \in K$, we have for any given $m \in [\![1,\, k+1]\!]$ and for some constants $\{\alpha_{im}\}_{i\in[\![1,\, n]\!]} \in \mathbb{R}$:

$$\forall i \in [\![1,\, n]\!], \quad \sum_{i=1}^{n} \alpha_{im}\, l_{i,m}(x,\, y)e_i(x,\, y) = 0 \quad \Rightarrow \quad \alpha_{im} l_{im}(x,\, y) = 0.$$

The above relation is valid for any $(x,\, y) \in K$. Therefore, knowing that there are only $k+1$ vanishing points for $l_{im}$, we can write the relation at $n$ other points and get

$$\alpha_{i,m} = 0 \quad \forall i \in [\![1,\, n]\!].$$

Finally, recalling that $\{l_{i,m}\}_{i,m}$ is a free set we can infer that there is no possible combination between the elements of $\{l_{i,m}e_i\}_m$ for any $i \in [\![1,\, n]\!]$. Thus, writing the previous relations for every $m$ in $[\![1,\, k+1]\!]$, summing them and enjoying the freeness of $\{l_{i,m}e_i\}_m$ for any fixed $i \in [\![1,\, n]\!]$ leads to:

$$\sum_{m=1}^{k+1} \sum_{i=1}^{n} \alpha_{im}\, l_{i,m}(x,\, y)e_i(x,\, y) = 0 \;\Rightarrow\; \forall (i,\, m) \in [\![1,\, n]\!] \times [\![1,\, k+1]\!],\; \alpha_{im} = 0$$

and in the end, the set $\mathcal{N}$ is free.

∎

**Proposition 3.17**

In the above setting of simplicial elements, the set $\mathcal{N}$ is free.

***Proof.*** Since in dimension $d$ we have $(d+1)$ faces, the element is necessarily convex and even if the outward normals vectors may be collinear, they cannot be identical. The previous assumption is then naturally fulfilled, in particular for $d = 2$. The freeness of $\mathcal{N}$ follows.

■

We can also derive properties that are due to the Lagrangian construction on each edge of the element.

**Property 3.18** Local Lagrangian property

For any $i \in [\![1,\, 3]\!]$ and for any sampling point $\mathrm{x}_l \in f_i$ that was used to generate the Lagrangian function set $\{l_{i,\,m}\}_{m\in[\![1,\,k+1]\!]}$, the set of functions

$$\{(x,\, y)^T \mapsto l_{i,\,m}(x,\, y) e_i(x,\, y)\}_{m\in[\![1,\,k+1]\!]}$$

enjoy the following local Lagrangian property.

$$l_{i,\,m}(x_l,\, y_l) e_i(x_l,\, y_l) \cdot n_i = c_i \delta_{\zeta(m)l} \tag{32}$$

Here, $c_i \in \mathbb{R}^*$ represents a constant depending on the shape of the reference element while $\zeta$ stands for some permutation function emphasizing the possible re-indexing of the Lagrangian functions (*see the Remark 1*).

***Proof.*** Let us fix some edge $f_i \in \partial K$, $i \in [\![1,\, 3]\!]$ and consider any $m \in [\![1,\, k+1]\!]$. Then, as $l_{i,\,m}$ is a scalar function one has

$$(l_{i,\,m}(x_l,\, y_l) e_i(x_l,\, y_l)) \cdot n_i = l_{i,\,m}(x_l,\, y_l)(e_i(x_l,\, y_l) \cdot n_i). \tag{33}$$

And by the definition of the vector $e_i$, we get:

$$\begin{aligned} l_{i,\,m}(x_l,\, y_l) e_i(x_l,\, y_l) \cdot n_i &= \left[ \begin{pmatrix} x_l + n_{ix} \\ y_l + n_{iy} \end{pmatrix} \cdot \begin{pmatrix} n_{ix} \\ n_{iy} \end{pmatrix} \right] l_{i,\,m}(x_l,\, y_l) \\ &= (x_l n_{ix} + n_{ix}^2 + y_l n_{iy} + n_{iy}^2) l_{i,\,m}(x_l,\, y_l) \\ &= \underbrace{(x_l n_{ix} + y_l n_{iy} + 1)}_{=c_i \text{ as } \mathrm{x}_l \in f_i} l_{i,\,m}(x_l,\, y_l) \end{aligned}$$

$$l_{i,\,m}(x_l,\, y_l) e_i(x_l,\, y_l) \cdot n_i = c_i \delta_{\zeta(m)l}, \quad \text{as } \mathrm{x}_l \text{ is a sampling point on } f_i.$$

■

***Remark 2.***

• Even though each of the functions in $\mathcal{N}$ enjoy a local Lagrangian property with respect to a set of $(k+1)$ aligned points on one specific edge, they do not enjoy a global Lagrangian property over the set of all the sampling points dispatched on every edge of the element. To see this, let us derive the same equations as previously without restricting ourselves to one specific edge.

When applying the relation (33) to a sampling point $\mathrm{x}_l$ belonging to the edge $f_j$, $j \neq i$, it comes equally $(l_{i,m}(x_l,\, y_l)e_i(x_l,\, y_l)) \cdot n_j = (x_l n_{jx} + y_l n_{jy} + (n_{ix}n_{jx} + n_{iy}n_{jy}))l_{i,m}(x_l,\, y_l)$ for any $m \in [\![1,\, k+1]\!]$. However, this relation is not reducible anymore as on any edge $f_j, j \neq i$ the term $n_{jx}n_{ix} + n_{jy}y_{iy}$ does not necessarily equal one.

Furthermore, as $\mathrm{x}_l$ is not one of the sampling points that generated $l_{i,m}$, the term $l_{i,m}(x_l,\, y_l)$ will not necessarily vanish and its value is not known *a - priori*. Indeed, it will vanish only if one of the coordinates of $\mathrm{x}_l$ matches some one - dimensional sampling point of $l_i$. It could typically arise when one designs a reference element $K$ having parallel edges and uses a homogeneous sampling for the points $\mathrm{x}_l$ across the edges, or when one drops the constant in the automatic definition of the $l_{i,m}$ in (31). Otherwise, as the degree of the polynomial is $k$, it cannot vanish elsewhere.

Therefore, except in rare triangle layouts where the term $(n_{ix}n_{jx}+n_{iy}n_{jy})$ vanishes for each sampling point $x_l$ of any edge $f_j$, there is no global Lagrangian property.

• Even if having a global Lagrangian property is not necessary to generate a free set $\mathcal{N}$ (as shown above), it may be computationally interesting to directly construct such a set. In the case of the presented triangle of reference, a simple modification of the vector $e_2$ is enough to ensure a global Lagrangian property. We then fall back in the setting of [9] by defining

$$\tilde{e_2} = \sqrt{2}\begin{pmatrix} x \\ y \end{pmatrix}.$$

One can easily check that the set $\{e_1,\, \tilde{e_2},\, e_3\}$ is free and that they enjoy the same properties as the one previously detailed. Furthermore, if we set $\{\mathrm{x}_l\}_l$ the set of all the sampling points of every edge, we have

$$l_{i,m}(\mathrm{x}_l)e_i(\mathrm{x}_l) \cdot n_i := l_s(\mathrm{x}_l)e_i(\mathrm{x}_l) \cdot n_i = \delta_{sl}$$

for any $s = (i-1)(k+1) + m \in [\![1,\, 3(k+1)]\!]$ representing the global index of the considered function. Indeed, for any function that have been generated from the first edge one can derive the following.

For any $\mathrm{x}_l \in f_1$, we have $\mathrm{x}_l = (x_l,\, 0)$ and

$$\begin{aligned} l_{1,m}(\mathrm{x}_l)e_1(\mathrm{x}_l)\cdot n_1 &= \left[\begin{pmatrix} x_l \\ y_l - 1 \end{pmatrix} \cdot \begin{pmatrix} 0 \\ -1 \end{pmatrix}\right] l_{1,m}(x_l,\, y_l) \\ &= (1 - y_l) l_{1,m}(x_l,\, y_l) \\ &= \delta_{\zeta(m)l}. \end{aligned}$$

For any $\mathrm{x}_l \in f_2$, we have $\mathrm{x}_l = (x_l,\, 1 - x_l)$ and

$$\begin{aligned} l_{1,m}(\mathrm{x}_l)e_1(\mathrm{x}_l)\cdot n_2 &= \left[\begin{pmatrix} x_l \\ y_l - 1 \end{pmatrix} \cdot \begin{pmatrix} \frac{\sqrt{2}}{2} \\ \frac{\sqrt{2}}{2} \end{pmatrix}\right] l_{1,m}(x_l,\, y_l) \\ &= \left(-\frac{\sqrt{2}}{2} + \frac{\sqrt{2}}{2}x_l + \frac{\sqrt{2}}{2}y_l\right) l_{1,m}(x_l,\, y_l) \\ &= 0. \end{aligned}$$

And lastly, for any $\mathrm{x}_l \in f_3$, we have $\mathrm{x}_l = (0,\, y_l)$ and

$$\begin{aligned} l_{1,m}(\mathrm{x}_l)e_1(\mathrm{x}_l)\cdot n_3 &= \left[\begin{pmatrix} x_l \\ y_l - 1 \end{pmatrix} \cdot \begin{pmatrix} -1 \\ 0 \end{pmatrix}\right] l_{1,m}(x_l,\, y_l) \\ &= -x_l l_{1,m}(x_l,\, y_l) \\ &= 0 \end{aligned}$$

Similarly, one has for the second edge

$$\begin{aligned} l_{2,m}(\mathrm{x}_l)\tilde{e_2}(\mathrm{x}_l)\cdot n_1 &= \sqrt{2}\left[\begin{pmatrix} x_l \\ y_l \end{pmatrix} \cdot \begin{pmatrix} 0 \\ -1 \end{pmatrix}\right] l_{2,m}(x_l,\, y_l) \\ &= -\sqrt{2} y_l l_{2,m}(x_l,\, y_l) \\ &= 0 \end{aligned}$$

for any $\mathrm{x}_l \in f_1$,

$$\begin{aligned} l_{2,m}(\mathrm{x}_l)\tilde{e_2}(\mathrm{x}_l)\cdot n_2 &= \sqrt{2}\left[\begin{pmatrix} x_l \\ y_l \end{pmatrix} \cdot \begin{pmatrix} \frac{\sqrt{2}}{2} \\ \frac{\sqrt{2}}{2} \end{pmatrix}\right] l_{2,m}(x_l,\, y_l) \\ &= \sqrt{2}\left(\frac{\sqrt{2}}{2}x_l + \frac{\sqrt{2}}{2}y_l\right) l_{2,m}(x_l,\, y_l) \\ &= (x_l + y_l) l_{2,m}(x_l,\, y_l) \\ &= l_{2,m}(x_l,\, y_l) \\ &= \delta_{\zeta(m)l} \end{aligned}$$

for any $\mathrm{x}_l \in f_2$, and for any $\mathrm{x}_l \in f_3$;

$$\begin{aligned} l_{2,m}(\mathrm{x}_l)\tilde{e_2}(\mathrm{x}_l) \cdot n_3 &= \sqrt{2}\left[\begin{pmatrix} x_l \\ y_l \end{pmatrix} \cdot \begin{pmatrix} -1 \\ 0 \end{pmatrix}\right] l_{2,m}(x_l,\, y_l) \\ &= -x_l l_{2,m}(x_l,\, y_l) \\ &= 0. \end{aligned}$$

Lastly, for the third edge we have

$$\begin{aligned} l_{3,m}(\mathrm{x}_l)e_3(\mathrm{x}_l) \cdot n_1 &= \left[\begin{pmatrix} x_l - 1 \\ y_l \end{pmatrix} \cdot \begin{pmatrix} 0 \\ -1 \end{pmatrix}\right] l_{3,m}(x_l,\, y_l) \\ &= -y_l l_{3,m}(x_l,\, y_l) \\ &= 0 \end{aligned}$$

for any $\mathrm{x}_l \in f_1$,

$$\begin{aligned} l_{3,m}(\mathrm{x}_l)e_3(\mathrm{x}_l) \cdot n_2 &= \left[\begin{pmatrix} x_l - 1 \\ y_l \end{pmatrix} \cdot \begin{pmatrix} \frac{\sqrt{2}}{2} \\ \frac{\sqrt{2}}{2} \end{pmatrix}\right] l_{3,m}(x_l,\, y_l) \\ &= \left(-\frac{\sqrt{2}}{2} + \frac{\sqrt{2}}{2}x_l + \frac{\sqrt{2}}{2}y_l\right) l_{3,m}(x_l,\, y_l) \\ &= 0 \end{aligned}$$

for any $\mathrm{x}_l \in f_2$, and for any $\mathrm{x}_l \in f_3$;

$$\begin{aligned} l_{3,m}(\mathrm{x}_l)e_3(\mathrm{x}_l) \cdot n_3 &= \left[\begin{pmatrix} x_l - 1 \\ y_l \end{pmatrix} \cdot \begin{pmatrix} -1 \\ 0 \end{pmatrix}\right] l_{3,m}(x_l,\, y_l) \\ &= (1 - x_l) l_{3,m}(x_l,\, y_l) \\ &= \delta_{\zeta(m)l}. \end{aligned}$$

Note that this global Lagrangian property also implies in itself the linearity independence of the above defined functions generating $\mathcal{N}$. ▲

***Remark.*** If one uses the definition of the vectors $e_1$, $\tilde{e_2}$ and $e_3$, the obtained set of functions (30) is already partly dual to the set of normal degrees of freedom (27a), in the sense that the corresponding transfer matrix is block diagonal. Indeed, all the normal moments are integral forms whose support matches exactly one edge. Therefore, when they are applied to functions vanishing on their support, the retrieved degree of freedom's value is zero.

By construction, this situation arises whenever the integral's supporting edge is different from the one being the support of the Lagrangian functions generating the tested basis function (*see the Figure 10*).

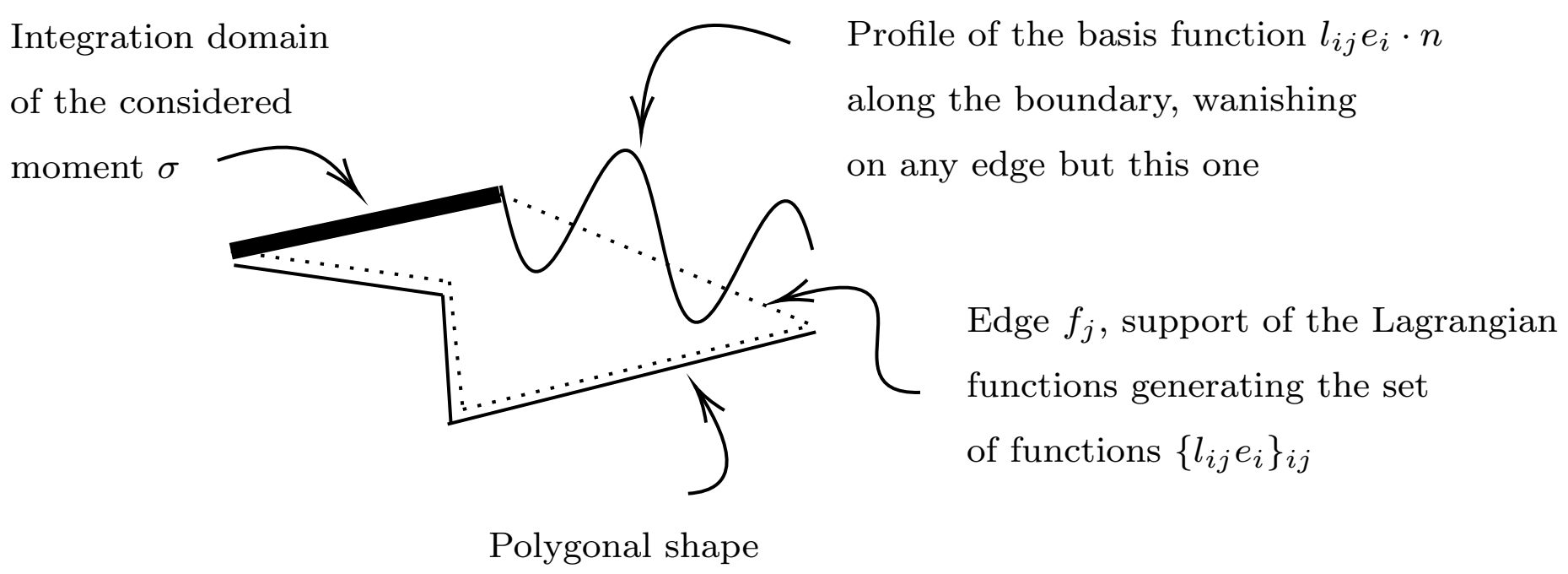

Fig. 10: Edge integration of basis functions vanishing on any edge but one.

Thus, writing the transfer matrix as explained in the *Section 2.4* for the described couple basis functions – set of degrees of freedom leads to a matrix of the following shape.

$$\begin{pmatrix} * & 0 & \cdots & 0 \\ 0 & * & 0 & \vdots \\ \vdots & 0 & \ddots & 0 \\ 0 & \cdots & 0 & * \end{pmatrix}$$

Moments whose supports match the $i^{\text{th}}$ edge

Basis functions not vanishing on the $i^{\text{th}}$edge

Here, each block of columns corresponds to the basis functions generated from a same edge and each block of rows corresponds to the moments whose support matches a specific edge. With an intuitive ordering of the rows and columns along the considered edges, the vanishing property detailed above leads to this block diagonal matrix. Furthermore, by definition of the space $\mathcal{R}_k(\partial K)$ the number of rows and column associated to one edge is $(k+1)$. Thus, the non-vanishing blocks are of size $(k+1) \times (k+1)$.

Though this matrix does not reduce to the identity (therefore the basis functions are not properly dual to the set of degrees of freedom), each non - vanishing block corresponds to all the normal information available at one specific edge only, information not being impacted by any other edge. Therefore, when one wants to transform the set of basis functions in order to retrieve the dual set to the degrees of freedom, the tuning is only done by blocks of $k+1$ functions and only changes the functions scaling, not their intrinsic behaviour. Indeed, functions of $\mathbb{P}_k$ stay in $\mathbb{P}_k$ by linear transformations and the support of the non - vanishing part on the boundary is unchanged.

To write it more precisely, as for any $q \in \mathcal{N}$ $q$ can be written under the form $q = l_{i,m}e_i$ for some $(i,\, m) \in [\![1,\, n]\!] \times [\![1,\, k+1]\!]$, we obtain by definition of the degrees of freedom (29)

$$\sigma(q) = \int_{f_i} l_{i,m}(x,y)e_i(x,y) \cdot n\, x^r \,\mathrm{d}\gamma(x,\, y)$$

for some path $\gamma$ skimming the edge and for any $r \in [\![0,\, k]\!]$. By example, in the case of the reference element presented in *Figure 6* we can derive in particular the degrees of freedom whose support is the third edge, $f_3$. We get for any $r \in [\![0,\, k]\!]$ and any $m \in [\![1,\, k+1]\!]$:

$$\sigma(q) = \int_0^1 l_{3,m}(0,\, y)e_3(0,y) \cdot n_3\, y^r \,\mathrm{d}y.$$

It comes

$$\begin{aligned}\sigma(q) &= \int_0^1 l_{3,m}(0,\, y)e_3(0,\, y) \cdot n_3 y^r \,\mathrm{d}y \\ &= \int_0^1 l_{3,m}(0,\, y) \begin{pmatrix} 0-1 \\ y \end{pmatrix} \cdot \begin{pmatrix} -1 \\ 0 \end{pmatrix} y^r \,\mathrm{d}y \\ &= \int_0^1 l_{3,m}(0,\, y) y^r \,\mathrm{d}y.\end{aligned}$$

Furthermore, we know that $y \mapsto l_{3,m}(0,\, y)$ is a polynomial of degree at most $k$ in the variable $y$, not depending on $x$. Thus, we can assume without lost of generality that the sampling points are allowing an exact quadrature rule. As there is only one of those points where $l_{3,\,m}$ is not null, and as it is associated with some real coefficient $c_k$ depending only on the number of points, and thus only on $k$, we have

$$\sigma(q) = c_k y^r_{\zeta(m)}. \tag{34}$$

For the other $2(k+1)$ functions that are generated from the edges $f_1$ and $f_2$, it comes for $i = 1,\, 2$:

$$\begin{aligned}\sigma(q) &= \int_0^1 l_{i,m}(0,\, y) e_i(0,\, y) \cdot n_i y^r \,\mathrm{d}y \\ &= \int_0^1 l_{i,m}(0,\, y) \times 0 \times y^r \,\mathrm{d}y \\ &= 0,\end{aligned}$$

where we have used the computations of the previous remark. One can derive similar properties when the support of the degrees of freedom matches the other three edges.

Thus, the tuning is only made by block of functions on one edge, having only a scaling role. The heuristic behaviour is not altered. If this may be more stable than tuning from the previous full system obtained when defining the basis functions from the vectors $e_1$, $e_2$ and $e_3$, we can observe that the amplitude of the obtained scaled functions may be high due to the relatively small constants (34) defining the matrix (6) to be inverted. Indeed, the sampling points should lie strictly inside the edges, and the reference element is contained within the unit circle, making the values $\sigma(q)$ sensibly small. ▲

***Remark.*** One can notice through the computation of $\sigma(q)$ that the highest $k$ is, the closest the matrix rows will be. Indeed, $y^r_{\zeta(m)} \to 0$ when $r \to \infty$, and thus $\sigma(q) \to 0$. Thus, the matrix rows corresponding to moments involving high $r$ are more and more similar and the conditioning of the matrix becomes worse and worse. ▲

***Note.*** As we already enjoy the desired conformity property, it may be advised to work directly on the non-tuned functions even if they do not correspond to the classical definition of the Raviart – Thomas elements. Therefore, no tuning is required to enjoy the desired property of $H(\mathrm{div},\, K)$ – conformity. However, if one prefers to use the classically defined basis functions and tune them after with the relations detailed in the *Section 3.5.2*, the tuned functions will equally enjoy a global Lagrangian property. ▲

***Remark 3.*** As the degrees of freedom are integration based, it is convenient to have an exact quadrature rule based on the sampling points only. Therefore, connecting with the *Remark 1* it may be better to generate the Gaussian points with respect to the $x$ or $y$ projection coordinates of the edges rather than to distribute the points on the edges themselves. Indeed, as the degrees of freedom are only one dimensional integrals, we can parametrize the edges and use only Gaussian quadrature points to compute it exactly, which would not be the case if one would have distributed it directly on the edge (*see the Figure 11*). Fortunately, in our reference element $\alpha_1 = 3\pi/4$ and therefore no particular care is required.

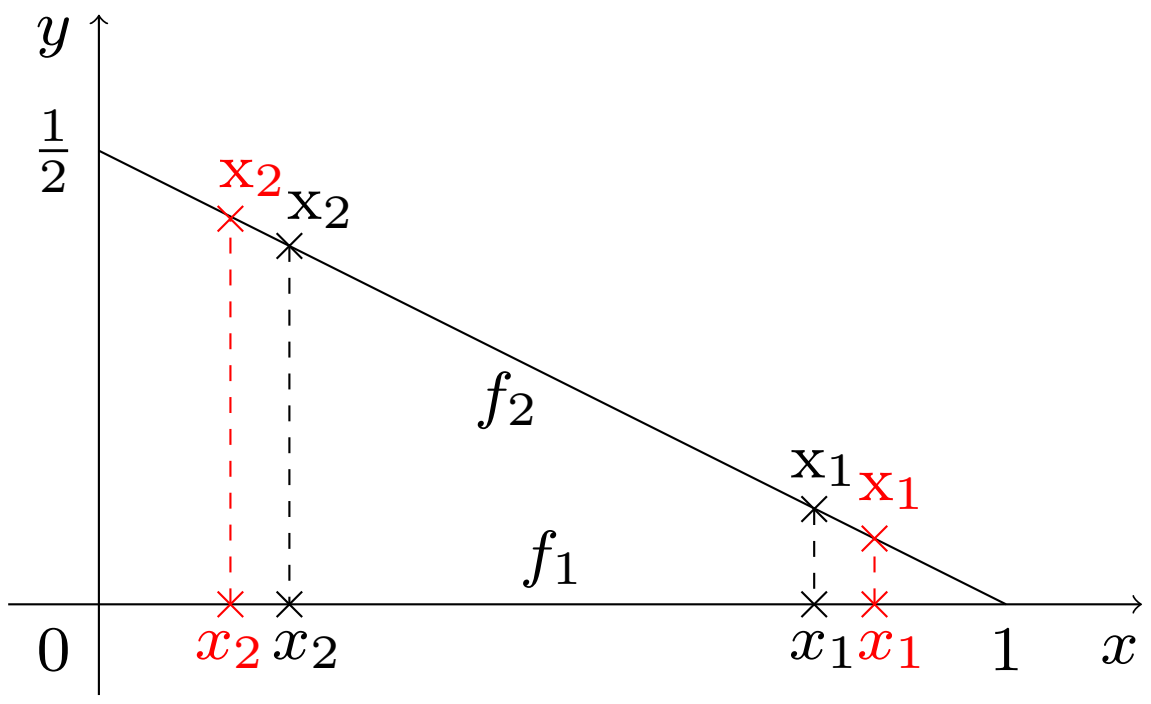

Fig. 11: Gaussian points projection: In red, Gauss-Legendre points generated from $f_2$ projected on $f_1$, reversely for the black labels. When integrating over $f_2$ through an edge parametrization by the variable $x$, only the black points will allow an exact quadrature.

▲

**Let us now generate $\mathcal{I}$.** We want to generate $d \dim \mathbb{P}_{k-1}(K)$ functions that are free from each other and whose set is free with $\mathcal{N}$. Those functions will later be tuned to form the set of internal basis functions dual to (27b).

***Note.*** For $k = 0$, $\mathbb{P}_{k-1}(K)$ does not exist. The set $\mathcal{I}$ is then automatically designed as the empty set and there exists no internal basis function. ▲

To this end, we first generate a subspace of $\mathbb{P}_{k-1}(K)$ that lies in $RT_k(K)$ and whose dimension is $d \dim \mathbb{P}_{k-1}(K) = \dim((\mathbb{P}_{k-1}(K))^d)$. In the same time, we also have to make sure that all of those function are free with $\mathcal{N}$. As the weakest sufficient conditions require to prescribe precise sampling points, it is nearly impossible to draw a general setting. However, as we enjoy the *Property 3.18* leading to the *Remark 2*, we inherit the strong sufficient condition $q \cdot n_f|_f = 0$ for all functions $q$ of $\mathcal{I}$ for any edge $f$ of $\partial K$.

On the same spirit as before, we start by defining the $d = 2$ vectors

$$e_{n+1}(x, y) = x \begin{pmatrix} x-1 \\ y \end{pmatrix} \quad \text{and} \quad e_{n+2}(x, y) = y \begin{pmatrix} x \\ y-1 \end{pmatrix}$$

enforcing the above vanishing property, where we used the general notation $n = 3$ to make easier later generalisation. We then take any basis $B_{k-1}(K)$ of $\mathbb{P}_{k-1}(K)$ and define our basis function set $\mathcal{I}$ as:

$$\mathcal{I} = \{p_{k-1}e_{n+i},\ i = 1, 2,\ p_{k-1} \in B_{k-1}(K)\}.$$

***Remark.*** Any definition of $e_{n+1}$ and $e_{n+2}$ would work here provided that $e_{n+1}$ and $e_{n+2}$ are free from each other, that they are of degree exactly two and verify $e_{n+1} \cdot n_f = 0$ and $e_{n+2} \cdot n_f = 0$ for every edge $f$ belonging to $\partial K$. ▲

**Proposition 3.19**

For all $p$ belonging to $\mathcal{I}$, $p$ belongs to $RT_k(K)$.

***Proof.*** If $p$ belongs to $\mathcal{I}$ then $p$ can be written as $p = p_{k-1}e_{n+1}$ or $p = p_{k-1}e_{n+2}$ for some $p_{k-1} \in B_{k-1}(K)$. We assume here that $p$ is written on the former form, $p_{k-1}e_{n+1}$. The other case can be treated similarly.

Let us assume that $p(x, y) = p_{k-1}(x, y) \begin{pmatrix} x^2 - x \\ yx \end{pmatrix}$.

Then,
$$\begin{aligned} p(x, y) &= p_{k-1}(x, y) \begin{pmatrix} x^2 \\ yx \end{pmatrix} + p_{k-1}(x, y) \begin{pmatrix} -x \\ 0 \end{pmatrix} \\ &= p_{k-1}(x, y) \begin{pmatrix} x \\ x \end{pmatrix} \circ \begin{pmatrix} x \\ y \end{pmatrix} + p_{k-1}(x, y) \begin{pmatrix} -x \\ 0 \end{pmatrix}. \end{aligned}$$

As $(x, y)^{\mathrm{T}} \mapsto (x, x)^{\mathrm{T}} p_{k-1}(x, y)$ belongs to $(\mathbb{P}_k(K))^2$ and since the last term also belongs to $(\mathbb{P}_k(K))^2$, up to separate the monomials we retrieve

$$\begin{aligned} p(x, y) &= \underbrace{\tilde{p}_{k-2}(x, y) \begin{pmatrix} x \\ x \end{pmatrix}}_{\in (\mathbb{P}_{k-1}(K))^2} \circ \begin{pmatrix} x \\ y \end{pmatrix} + \underbrace{\tilde{p}_{[k-1]} \begin{pmatrix} x \\ x \end{pmatrix}}_{\in (\mathbb{P}_{[k]})^2} \circ \begin{pmatrix} x \\ y \end{pmatrix} + \underbrace{p_{k-1}(x, y) \begin{pmatrix} -x \\ 0 \end{pmatrix}}_{\in (\mathbb{P}_k(K))^2} \\ &= \begin{pmatrix} x \\ y \end{pmatrix} \underbrace{p_{[k]}(x, y)}_{\in \mathbb{P}_{[k]}} + \underbrace{\left( \tilde{p}_{k-2}(x, y) \begin{pmatrix} x \\ x \end{pmatrix} \circ \begin{pmatrix} x \\ y \end{pmatrix} \right.}_{\in (\mathbb{P}_k(K))^2} + \underbrace{\left. \begin{pmatrix} -x \\ 0 \end{pmatrix} p_{k-1}(x, y), \right)}_{\in (\mathbb{P}_k(K))^2} \end{aligned}$$

showing that $p$ belongs to $RT_k(K)$.

■

***Remark.*** The above proof only works by the choice of the polynomial space as being $\mathbb{P}_k$. Indeed, there cannot be some function $p_{k-1}$ in $B_{k-1}(K)$ such that $(x, y) \mapsto (x, x)^{\mathrm{T}} p_{k-1}(x, y) \in \mathbb{P}_k(K) \times \mathbb{P}_{k-1}(K)$ or $\mathbb{P}_{k-1}(K) \times \mathbb{P}_k(K)$ allowing the direct split on the two first terms. This is due to the definition of the polynomial degree as being the product of the multiplicities independently on which coordinate takes the highest multiplicity. This will not be the case anymore when working on Quads, where the space $\mathbb{Q}_k$ is used. ▲

**Property 3.20**

$\dim \mathcal{I} = d \dim(\mathbb{P}_{k-1}(K))$

***Proof.*** Here we are restricted to the case $d = 2$. We have two linearly independent vectors $e_{n+1}$ and $e_{n+2}$. Each of them is used to generate functions through a basis of $\mathbb{P}_{k-1}(K)$. Thus, $\dim \mathcal{I} = 2 \dim \mathbb{P}_{k-1}(K) = 2\frac{k(k+1)}{2} = k(k+1)$.

■

**Property 3.21**

For all $p$ in $\mathcal{I}$, for all edges $f$ of $\partial K$ and for any point $(x, y) \in f$, we have $p(x, y) \cdot n_f = 0$.

***Proof.*** We recall that we are working on our reference element designed in *Figure 6*. Thus, we have for all $(x, y) \in f_1$;

$$
\begin{aligned}
p_1(x, y) \cdot n_1 &= p_{k-1}(x, y) \begin{pmatrix} x^2 - x \\ yx \end{pmatrix} \cdot n_1 \\
&= p_{k-1}(x, y) \begin{pmatrix} x^2 - x \\ 0 \end{pmatrix} \cdot \begin{pmatrix} 0 \\ -1 \end{pmatrix} \\
&= 0 \\
p_2(x, y) \cdot n_1 &= p_{k-1}(x, y) \begin{pmatrix} yx \\ y^2 - y \end{pmatrix} \cdot n_1 \\
&= p_{k-1}(x, y) \begin{pmatrix} 0 \\ 0 \end{pmatrix} \cdot \begin{pmatrix} 0 \\ -1 \end{pmatrix} \\
&= 0,
\end{aligned}
$$

and for all $(x, y) \in f_3$;

$$\begin{aligned} p_1(x, y) \cdot n_2 &= p_{k-1}(x, y) \begin{pmatrix} 0 \\ 0 \end{pmatrix} \cdot \begin{pmatrix} -1 \\ 0 \end{pmatrix} \\ &= 0 \\ p_2(x, y) \cdot n_2 &= p_{k-1}(x, y) \begin{pmatrix} 0 \\ y^2 - y \end{pmatrix} \cdot \begin{pmatrix} -1 \\ 0 \end{pmatrix} \\ &= 0. \end{aligned}$$

Lastly for all $(x, y) \in f_2$ we get;

$$\begin{aligned} p_1(x, y) \cdot n_0 &= \frac{1}{2} p_{k-1}(x, y) \begin{pmatrix} x^2 - x \\ x(1-x) \end{pmatrix} \cdot \begin{pmatrix} \sqrt{2} \\ \sqrt{2} \end{pmatrix} \\ &= \frac{\sqrt{2}}{2} p_{k-1}(x, y)(x^2 + x - x - x^2) \\ &= 0 \\ p_2(x, y) \cdot n_0 &= \frac{1}{2} p_{k-1}(x, y) \begin{pmatrix} x(1-x) \\ (1-x)^2 - (1-x) \end{pmatrix} \cdot \begin{pmatrix} \sqrt{2} \\ \sqrt{2} \end{pmatrix} \\ &= 0. \end{aligned}$$

■

**Proposition 3.22**

$\mathcal{I}$ is linearly independent strictly within the element $K$.

***Proof.*** *We repeat here the arguments of [9].* We consider a linear combination of functions in $\mathcal{I}$, and look for collapsing combinations. For any $\alpha_i$ and $\beta_i$, $i = [\![1, \frac{1}{2}k(k+1)]\!]$, we write:

$$\left( \sum_{i=1}^{N_k} \alpha_i p_i(x, y) \right) e_{n+1}(x, y) + \left( \sum_{i=1}^{N_k} \beta_i p_i(x, y) \right) e_{n+2}(x, y) = 0,$$

where we set $N_k = \frac{k(k+1)}{2}$ for legibility. We can assume without loss of generality that the chosen basis of $B_{k-1}(K)$ enjoys Lagrangian properties, *i.e.* there exists some nodes $\mathrm{x}_l$, $l \in [\![1, N_k]\!]$ strictly within the element $K$ such that $p_i(x_l, y_l) = \delta_{il}$. Taking in particular the above combination at those nodes, we get:

$$\forall l \in [\![1, N_k]\!],\ \alpha_l e_{n+1}(x_l, y_l) + \beta_l e_{n+2}(x_l, y_l) = 0.$$

This reduces to

$$\beta_l \begin{pmatrix} x_l y_l \\ y_l^2 - y_l \end{pmatrix} + \alpha_l \begin{pmatrix} x_l^2 - x_l \\ x_l y_l \end{pmatrix} = 0$$

$$\Leftrightarrow \quad \begin{cases} \beta_l x_l y_l + \alpha_l x_l^2 - \alpha_l x_l & = 0 \\ \beta_l y_l^2 - \beta_l y_l + \alpha_l x_l y_l & = 0 \end{cases} \tag{35}$$

We notice that if $(x_l,\, y_l)$ belongs to the segments $[\,0,\, 1\,] \times \{0\}$, $\{0\} \times [\,0,\, 1\,]$ or $(x,\, 1-x)_{x\in[0,\,1]}$, then there exists non-trivial solutions to (35). Otherwise, the only solution is given by:

$$(\alpha_l,\, \beta_l) = (0,\, 0) \quad \forall l \in [\![1,\, N_k]\!].$$

Thus, $\mathcal{I}$ is free strictly within the element.

■

**Let us gather the sets $\mathcal{N}$ and $\mathcal{I}$.** Let us finally show that both sets $\mathcal{N}$ and $\mathcal{I}$ together generates $RT_k(K)$.

**Proposition 3.23**

We have $RT_k(K) = \text{span}\,\{\mathcal{N}, \mathcal{I}\}$.

Let us prove the underlying properties leading to the above proposition.

**Proposition 3.24**

The set $\{\mathcal{N},\, \mathcal{I}\}$ is free.

***Proof.*** The case $k = 0$ is trivial as $\mathcal{I}$ is an empty set. Let us consider the case $k > 0$. We saw previously that both $\mathcal{N}$ and $\mathcal{I}$ are free sets. It is left to show that no interconnection happens within the gathered set.

Proving that there exists no interconnection between the sets $\mathcal{N}$ and $\mathcal{I}$ reduces to see that for a combination of functions to be constant on one edge, it should be constant everywhere in $K$. By extension it would be zero everywhere, which prevents for creating a function outside of $\text{span}\,\{\mathcal{N}\}$. Let us detail it more precisely.

• Let us show first that one cannot generate a function of $\mathcal{I}$ from functions belonging to $\mathcal{N}$ and the other functions of $\mathcal{I}$.

We know that $\mathcal{I}$ is a free set and that its functions vanish with respect to the outward normal on every edge, which is not the case for any function of $\mathcal{N}$. The only possibility would then be to find some combination of functions in $\mathcal{N}$ vanishing on every edge, which is not identically zero and that combined to functions in $\mathcal{I}$ generates an other function of $\mathcal{I}$. We show that this never happens.

To start with, we have by construction that any function of $\mathcal{N}$ verifies

$$\begin{aligned} l_{i,m}e_i(x,\, y)\cdot n_j|_{f_j} &= (xn_{jx} + yn_{jy} + (n_{ix}n_{jx} + n_{iy}n_{jy}))l_{i,m}(x,\, y) \\ &= C\, l_{i,j}(x,\, y) \end{aligned} \tag{36}$$

for a given constant $C_j$ depending on the edge $f_j$, some $i \in [\![1,\, n]\!]$ and some $j \in [\![1,\, k+1]\!]$. Note that the last relation is obtained as the term $xn_{ix} + yn_{iy}$ is constant on the edge $f_j$. We also have to keep in mind that by definition $l_{i,m}(x,\, y) = l_{i,m}(x)$ or $l_{i,m}(x,\, y) = l_{i,m}(y)$. No variable coupling is involved here. Therefore, we can assume without loss of generality that

$$\begin{cases} l_{1,m}(x,\, y) = l_{1,m}(x) \\ l_{2,m}(x,\, y) = l_{2,m}(x) \\ l_{3,m}(x,\, y) = l_{3,m}(y). \end{cases} \tag{37}$$

We can now show that there exists no non - trivial combination of the functions defined in (36) vanishing on every edge $f_j$ of $\partial K$.

Let us first show that for any $(x,\, y) \in \partial K$ no linear combination of the form

$$(x,\, y) \mapsto A\sum_{m=1}^{k+1} l_{1,m}(x) + B\sum_{m=1}^{k+1} l_{2,m}(x) + C\sum_{m=1}^{k+1} l_{3,m}(y), \quad {\scriptstyle A,B,C\in\mathbb{R}, \atop (A,B,C)\neq(0,0,0),} \tag{38}$$

of $\{l_{i,m}\}_{i,m}$ can vanish with respect to the normal component simultaneously on each edge. We consider for now the edge $f_1$. Using (37), one can observe that $l_{2,m}(x,\, y)|_{f_1} \equiv l_{2,m}(1)$ and $l_{3,m}(x,\, y)|_{f_1} \equiv l_{3,m}(0)$. Thus, those functions are constants on the edge $f_1$ (*see the Figure 12*).

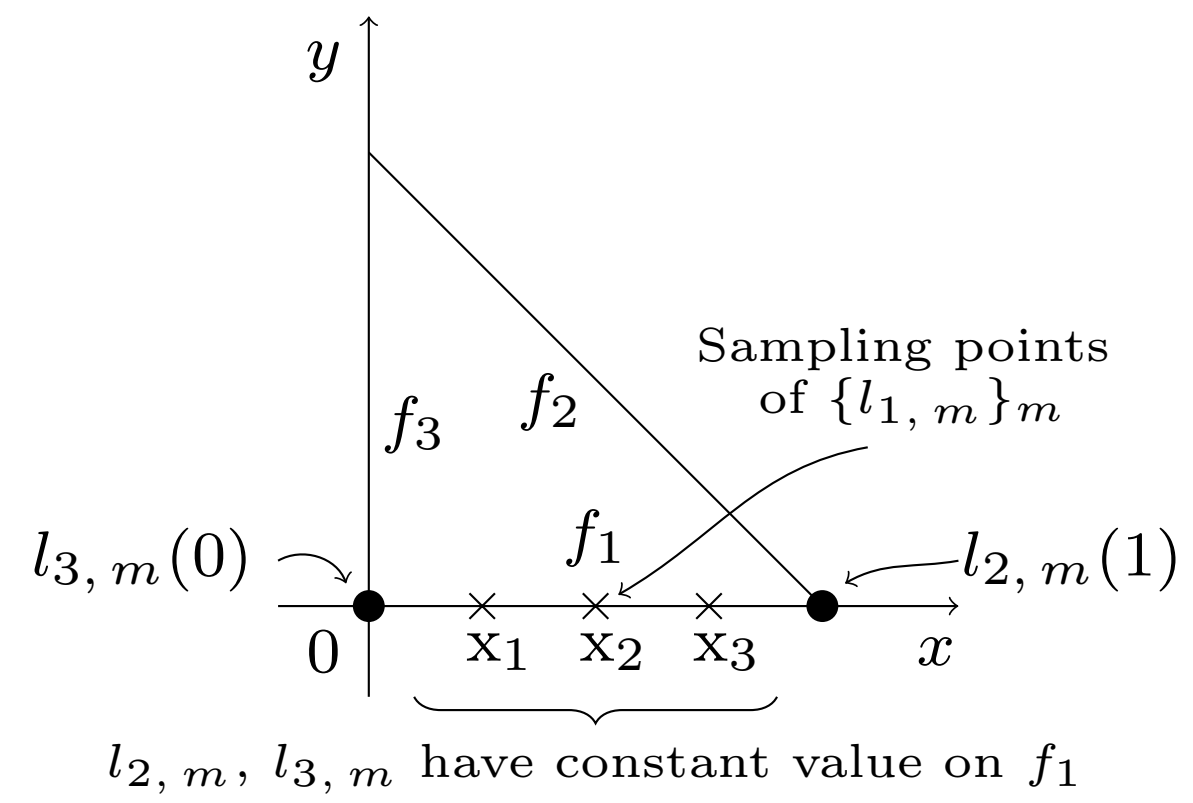

Fig. 12: Coordinate - wise values of Lagrangian functions.

Furthermore, the function $l_{1,m}$ is nothing else on $f_1$ that the $m^{\text{th}}$ Lagrangian function associated to $(k+1)$ sampling points lying on that edge. Therefore, if one considers the Lagrangian set $\{l_{1,m}\}$, the only constant value one can achieve on $f_1$ when $k > 0$ is when all the local Lagrangian functions are considered with the same weight (*see the Figure 13*). The constant will then necessarily value one.

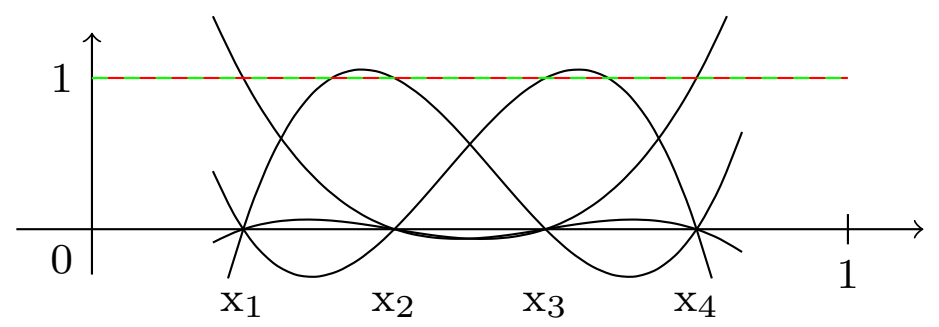

Fig. 13: Example of local Lagrangian functions on one edge for four sampling points $\{x_l\}_{l\in[\![1,4]\!]}$. In green is the constant value 1, in red is the sum of the four Lagrangian functions.

Thus, this set of local Lagrangian functions cannot vanish by themselves. One has to drag the constant down by the help of the two other sets of $k+1$ local Lagrangian functions that are constant on $f_1$, meaningly $l_{2,m}$ and $l_{3,m}$. Note also that the same discussion applies to the sets $\{l_{2,m}\}_m$ and $\{l_{3,m}\}_m$ on the edges $f_2$ and $f_3$ respectively. Therefore, we need to drag down those three sets of functions simultaneously (*see the Figures 14 and 15*). However, each edge is connected to the other two. Therefore, knowing that we can only play on the global offset of each equally weighted combinations of $\{l_{i,m}\}_m$, the only solution to see the expression (38) vanishing on every edge is to fulfil the equation

$$A + B + C = 0.$$

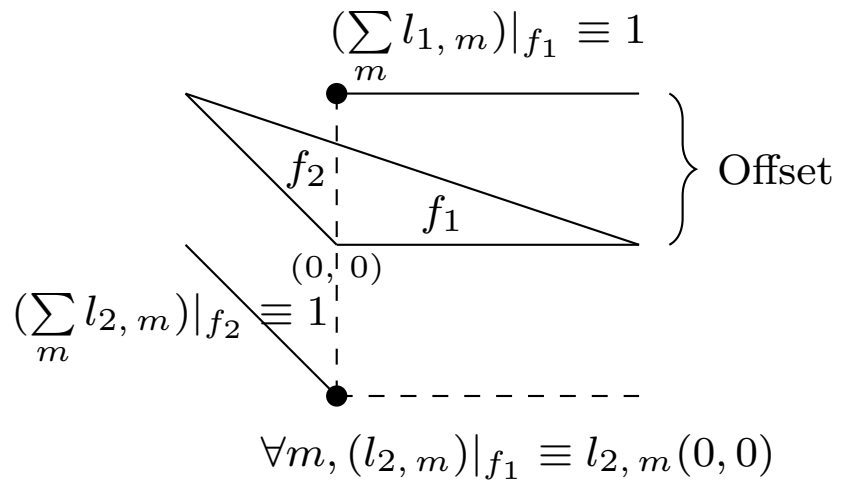

Fig. 14: Offset modulation on triangle while respecting the Lagrangian constraints.

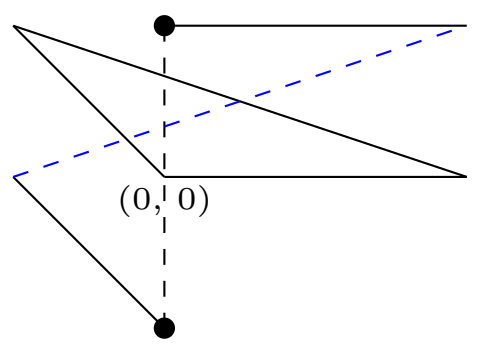

Fig. 15: If not null, the last offset constraint imposes the use of a non - constant Lagrangian function on the third edge, which is not allowed.

However, respecting this condition implies that the combination will be vanishing everywhere, and therefore, no linear combination of functions in $\{l_{i,m}\}_{i,m}$ can vanish on every edge at the same time without being identically null.

**Note:** This incompatibility appears to be more evident in the case where global Lagrangian functions are used. The property is less visible here as we have terms coupling some functions' features, but the fundamental behaviour is the same. A simple change of basis make it pop-up.

Let us now show that for any fixed edge $f_i$ the term $(xn_{ix} + yn_{iy} + n_{ix}n_{jx} + n_{iy}n_{jy})$ cannot vanish on each edge $f_j$ simultaneously. To do so, let us recall that we are in a convex case and that for any x $\in f_i$ the term $c_1 := xn_{ix} + yn_{iy}$ is always constant. We start by looking for the cases where $c_1 = -n_{ix}n_{jx} - n_{iy}n_{jy}$.

We consider first the case where $c_1 = 0$. The only cases where $c_1$ vanishes is when the submits of the edge $e_i$ are aligned with the origin (*see the Figure 16*). As we are on triangles, this can happen at maximum on two edges. On the other one(s) the combination of Lagrangian functions will not be vanishing for the reasons evoked above. Therefore, the terms (36) are not vanishing for the three edges simultaneously. Note that the case $-n_{ix}n_{jx} - n_{iy}n_{jy} = 0$, where the two edges $f_j$ and $f_i$ are orthogonal, merges with this one.

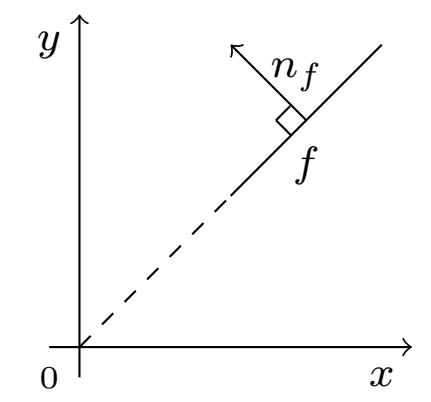

Fig. 16: Example of cases where $c_1$ vanishes (For any $(x,\, y)^T \in f$, $(x - 0,\, y - 0)^T \cdot n_f = 0$).

If $c_1 \neq 0$, then one gets immediately the relation $y = \frac{-n_{jx}}{n_{jy}}(x + n_{ix}) + n_{iy}$. Simultaneously we know by definition of the element $K$ that on that edge, $y = \frac{-n_{ix}}{n_{iy}}x + b$ for some real constant $b$. Therefore, the only vanishing possibility leads to the equation

$$\begin{aligned} \frac{-n_{ix}}{n_{iy}}x + b &= \frac{-n_{ix}}{n_{iy}} - \frac{n_{jx}n_{ix}}{n_{jy}} + n_{iy} \\ \Rightarrow n_{iy} + b &= \frac{-n_{ix}n_{jx}}{n_{jy}}. \end{aligned}$$

Up to translating the element of reference so that $b = 0$, we get

$$\begin{aligned} n_{iy}^2 &= -n_{ix}n_{jx} \\ \Rightarrow n_{jx} \leq 0 &\text{ or } n_{ix} \leq 0. \end{aligned}$$

However, this can happen at maximum on two edges since if $i = j$ then $n_{ix}n_{jx} \geq 0$. On the third one, there will be left a non - vanishing contribution of non - vanishing local Lagrangian functions that cannot be compensated, and the normal component of combinations of functions in $\{l_{i,\,m}\}_{i,\,m}$ cannot vanish on every edge. Furthermore, note that if $n_{iy} = 0$ we have $n_{ix} = 0$ at maximum only once as we work here with triangles.

We can then finally conclude that no combination of functions of $\mathcal{N}$ can have their normal component vanishing on every edge without being identically null. And therefore, no function of $\mathcal{I}$ can be generated. To be more precise, one can see it computationally by writing

$$\sum_{i=1}^{n} A_i e_i \cdot n_j = 0 \quad \forall j \in [\![1,\, n]\!] \Rightarrow A_i = 0 \quad \forall i \in [\![1,\, n]\!].$$

Indeed, one can see that the system

$$\begin{cases} \left(A_1 \begin{pmatrix} x+0 \\ y-1 \end{pmatrix} + A_2 \begin{pmatrix} x+\frac{\sqrt{2}}{2} \\ y+\frac{\sqrt{2}}{2} \end{pmatrix} + A_3 \begin{pmatrix} x-1 \\ y-0 \end{pmatrix}\right) \cdot \begin{pmatrix} 0 \\ -1 \end{pmatrix} = 0 \\ \left(A_1 \begin{pmatrix} x+0 \\ y-1 \end{pmatrix} + A_2 \begin{pmatrix} x+\frac{\sqrt{2}}{2} \\ y+\frac{\sqrt{2}}{2} \end{pmatrix} + A_3 \begin{pmatrix} x-1 \\ y-0 \end{pmatrix}\right) \cdot \begin{pmatrix} \frac{\sqrt{2}}{2} \\ \frac{\sqrt{2}}{2} \end{pmatrix} = 0 \\ \left(A_1 \begin{pmatrix} x+0 \\ y-1 \end{pmatrix} + A_2 \begin{pmatrix} x+\frac{\sqrt{2}}{2} \\ y+\frac{\sqrt{2}}{2} \end{pmatrix} + A_3 \begin{pmatrix} x-1 \\ y-0 \end{pmatrix}\right) \cdot \begin{pmatrix} -1 \\ 0 \end{pmatrix} = 0 \end{cases}$$

reduces to

$$\begin{cases} -\frac{\sqrt{2}}{2}A_2 - A_3 & = 0 \\ -A_1 - \frac{\sqrt{2}}{2}A_2 & = 0 \\ \frac{\sqrt{2}}{2}A_1 \times 0 + \frac{\sqrt{2}}{2}A_2(1+\frac{\sqrt{2}}{2}) + \frac{\sqrt{2}}{2} \times 0 & = 0 \end{cases}$$

which has for unique solution

$$A_1 = A_2 = A_3 = 0.$$

• We also have that no combination of normal basis functions can generate an internal basis function and *vice - versa.* Indeed, for any $(x,\, y) \in \partial K$ it holds

$$p_{i,m}(x,\, y) \cdot n|_{\partial K} \not\equiv 0$$

for all the normal basis functions $\{p_{i,\,m}\}_{i,\,m}$. Furthermore, we saw previously that no combination of normal functions can have their normal component vanishing on all the boundaries without being identically null. And in the same time, on all the boundaries the internal basis functions satisfy

$$p(x,\, y) \cdot n_f|_f = 0.$$

Therefore, no standalone interconnection can be built between the two subsets.

• Lastly, one cannot either create a function of $\mathcal{N}$ with the help of functions lying in $\mathcal{I}$ and the others functions of $\mathcal{N}$. To show this, it is enough to consider a function of $\mathcal{N}$,

$$f_i = \begin{pmatrix} x - n_{ix} \\ y - n_{iy} \end{pmatrix} l_{i,\,m}(x,\, y) \tag{39}$$

for some $m \in [\![1,\, k+1]\!]$, a function of $\mathcal{I}$ that reads

$$p_j = x \begin{pmatrix} x-1 \\ y \end{pmatrix} b_j(x,\, y)$$

for some $b_j \in \mathbb{P}_{k-1}(K)$ and to show that no function of the form

$$q_m = \alpha x \begin{pmatrix} x-1 \\ y \end{pmatrix} b_j(x,\, y) + \beta \begin{pmatrix} x - n_{ix} \\ y - n_{iy} \end{pmatrix} l_{i,\,m}(x,\, y) \tag{40}$$

can actually belong to the set $\mathcal{N}$. Indeed, as all functions in $\mathcal{N}$ are free from each other and share the same shape, their combination would enjoy the same structure as the one presented in (40) and the arguments presented hereafter would hold in the same way. Note that using $e_{n+2}$ instead of $e_{n+1}$ in the definition of $p_m$ also results in the same argumentation.

Let us first point out that any vanishing combination of $\mathcal{N}$ cannot construct a new function belonging to $\mathcal{N}$ as the functions $\mathcal{I}$ themselves cannot do it. This case will therefore be excluded, allowing us to only show that no combination of functions of $\mathcal{N}$ and $\mathcal{I}$ can create another function of $\mathcal{N}$. We start by reordering the terms of $q_m$ as

$$\begin{aligned} q_m &= \begin{pmatrix} \alpha\, x l_{i,\,m}(x,\, y) + \beta\; x(x-1) b_j(x,\, y) - \alpha\, n_{ix}\, l_{i,\,m}(x,\, y) \\ \beta\,(y - n_{iy})\, l_{i,\,m}(x,\, y) + \alpha\, xy b_j(x,\, y) \end{pmatrix} \\ &= \begin{pmatrix} x(\alpha\, l_{i,\,m} - \beta\, b_j(x,\, y)) + \beta\, x^2 b_j(x,\, y) - n_{ix}\, \alpha\, l_{i,\,m}(x,\, y) \\ y\,(\beta\, l_{i,\,m}(x,\, y) + \alpha\, x b_j(x,\, y)) - \beta\, n_{iy} l_{i,\,m}(x,\, y) \end{pmatrix}. \end{aligned}$$

To generate a function of $\mathcal{N}$, we would like to find some $a \in [\![1,\, n]\!]$ and $m \in [\![1,\, k+1]\!]$ such that

$$q_m = \begin{pmatrix} x - n_{ax} \\ y - n_{ay} \end{pmatrix} l_{a,\,m}(x,\, y) \in \mathcal{N}.$$

The only opportunity to build such a vector is to ask $(x,\, y) \mapsto x\, b_j(x,\, y) \in \mathbb{P}_k(K)$. Thus, by the definition of the polynomial's degree, $b_j$ must belong to $\mathbb{P}_{k-1}(K)$ so that the term $x^2 b_j$ can be catch up in the form $x l_{l,\,m}(x,\, y)$. However doing so imposes $(x,\, y) \mapsto x b_j(x,\, y) \in \mathbb{P}_k(K)$ also in the second coordinate, which makes the term $(x,\, y) \mapsto y x b_j(x,\, y)$ not belonging to $\in \mathbb{P}_k(K)$. Therefore, the second coordinate cannot be factorized in the shape $(y - n_{ay}) p_k$ for some $p_k \in \mathbb{P}_k(K)$, and *a fortiori* not in the shape $(y - n_{ay}) l_{l,\,m}$. The use of $e_{n+2}$ based internal functions to compensate would not help either as it would bring the same issue on the $x$ - coordinate.

Thus, $q_m$ cannot be written in the shape (39) and the construction (40) consequently does not belong to $\mathcal{N}$. Gathering all the discussed cases together, the claim of the property 3.24 follows.

■

***Remark.*** When one uses the vector $\tilde{e_2}$, a simpler proof of freedom can be derived. Indeed, in that case one has a global Lagrangian property on the sampling boundary points $x_l$;

$$p_{i,m}(x_l) \cdot n_i = \delta_{ml}$$

for the normal basis functions $\{p_{i,\,m}\}_{i,\,m}$, and

$$p(x,\,y) \cdot n_f|_f = 0$$

on all the boundaries for the internal basis functions. As the second relation holds in particular for the sampling points distributed on the boundary $\partial K$ used to generate $p_{im}$, $p(x_l,\,y_l) \cdot n_i = 0$. Thus, the set $\mathcal{N}$ will be free with respect to any $p$ in $\mathcal{I}$, and the set $\{\mathcal{N},\,\mathcal{I}\}$ is then free. ▲

Proposition 3.25

$\dim\{\mathcal{N},\,\mathcal{I}\} = \dim RT_k(K)$

***Proof.*** By the fact that the sets $\mathcal{N}, \mathcal{I}$ are free from each other and that we consider simplicial elements, we directly have from the previous sections that

$$\begin{aligned}\dim\{\mathcal{N},\,\mathcal{I}\} &= 3(k+1) + (k+1)(k+2)\\ &= 3\dim\mathbb{P}_k(f) + 2\dim\mathbb{P}_{k-1}(K)\\ &= \dim RT_k(K).\end{aligned}$$

■

Proposition 3.26

If $p$ belongs to span $\{\mathcal{N},\,\mathcal{I}\}$, then $p$ belongs to $RT_k(K)$.

***Proof.*** We saw previously that

$$\begin{cases} p \in \mathcal{N} & \Rightarrow \quad p \in RT_k(K)\\ p \in \mathcal{I} & \Rightarrow \quad p \in RT_k(K)\end{cases}$$

As $RT_k(K)$ is a vectorial space, the result immediately follows by linear combination of the two previous quantities.

■

***Remark.*** So constructed, the full set of basis functions $\{\mathcal{N}, \mathcal{I}\}$ is not partly dual to the full set of internal and normal degrees of freedom. Indeed, neither those normal basis functions are vanishing within the tetrahedra, nor is their contribution through the internal degrees of freedom. Therefore, the full transfer matrix from the canonical basis functions of the space $RT_k(K)$ to the dual set of the full degrees of freedom (27) is not diagonal by block anymore. Indeed, the values corresponding to the internal degree of freedom tested against normal functions are no zero. An example is given in the numerical results of the *Section 7*. ▲

**Summary of the construction** The following basis functions form a basis for $RT_k(K)$ on simplicial elements in dimension two.

*Normal basis functions:*

$$\phi_j = l_{i,m} e_i, \quad \substack{1\leq i\leq 3,\, 1\leq m\leq k+1 \\ 1\leq j\leq n(k+1)}$$

*Internal basis functions:*

$$\phi_{3(k+1)+l} = p_{k-1,l} e_{n+1}, \quad l = [\![1, \frac{1}{2}k(k+1)]\!]$$

$$\phi_{3(k+1)+N_k+l} = p_{k-1,l} e_{n+2}, \quad l = [\![1, \frac{1}{2}k(k+1)]\!],$$

with $p_{k-1,l}$ the $l^{\text{th}}$ basis function of the space $B_{k-1}(K)$.

### 3.5.2 Basis function tuning on the degrees of freedom

In the previous paragraph we generated a generic basis of $RT_k(K)$ that does not necessarily describe the element $(K, RT_k(K), \{\sigma_i\}_i)$. In order to find a corresponding basis, it is enough to adapt the generic basis functions towards the set of degrees of freedom $\{\sigma_i\}_i$ so that their set becomes dual to the moments (27).

***Note.*** Even if this two-steps basis functions building may seem tedious, it can be numerically convenient. In particular, it allows any arbitrary choice in the definition of the set $\{b_j\}_j$ and of the sampling points on each edge. ▲

A direct application of the algebra results presented in the *Section 2.4* allows us to retrieve our tuned basis functions as $\varphi_i = \sum_{m=1}^{N} \alpha_{im}\phi_m$, $i \in [\![1, N]\!]$ for some coefficients depending on the value of the degrees of freedom computed from the raw basis functions over the element $K$. As the basis change transfers to the derivatives, we also get as a by-product that $\partial_x\varphi_i = \sum_{m=1}^{N} \alpha_{im}\partial_x\phi_m$, $i \in [\![1, N]\!]$. We refer to the *Section 2.4* for details.

### 3.5.3 Transformations to and from the reference elements

All the discussions above apply to a reference element $K$. To be able to deal with any distortion of this reference element, we need an orientation-preserving map. Indeed, due to the definition of the normal moments and the wish of H(div) – conformity, we need the normal's orientation to be preserved through the transformation to or from the reference element. In the case of triangles, the Piola transform fulfil all of our requirements. We refer to the *Section 2.5.1* for details.

## 3.6 Numerical implementation

We implemented the proposed construction method for the global Lagrangian basis functions before tuning it against the set of degrees of freedom to retrieve the classical Raviart – Thomas setting. We report here the results for the third order space.

We detail first the results obtained for the build of the basis from the global Lagrangian functions set on the boundary. There, one can observe that we already have the same behaviour as the one prescribed by the Raviart – Thomas degrees of freedom, that is retrieving polynomials of degree three on the boundaries. Furthermore, the basis functions are bounded by one within the reference element. In particular, for the internal basis functions, one can verify that indeed, even if coordinate wise they are not vanishing on the boundary, their normal component is. Thus, the canonical basis functions are already compliant with the Raviart – Thomas setting.

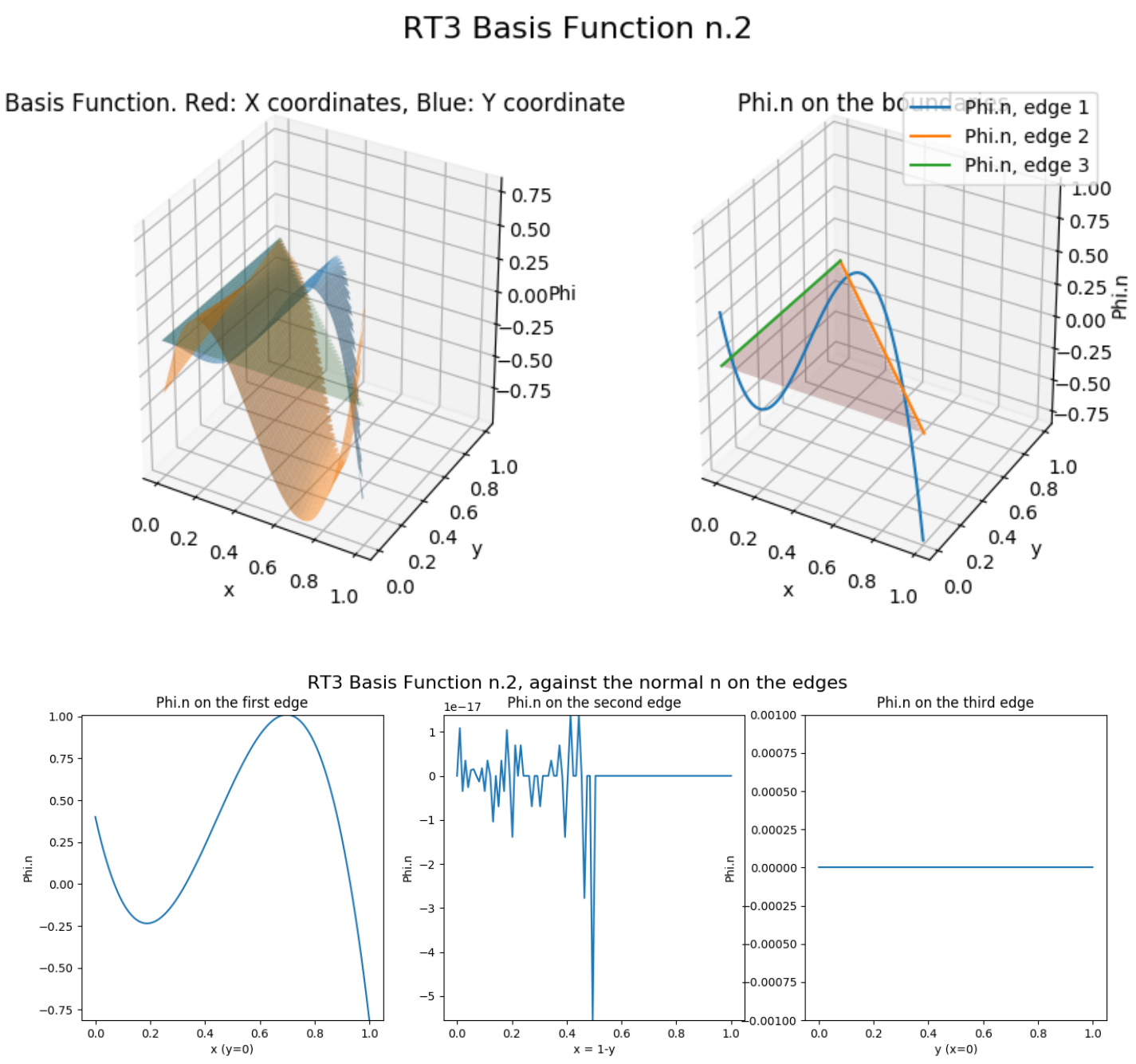

Fig. 17: Plot of the third canonical normal basis function in the case $k = 3$.

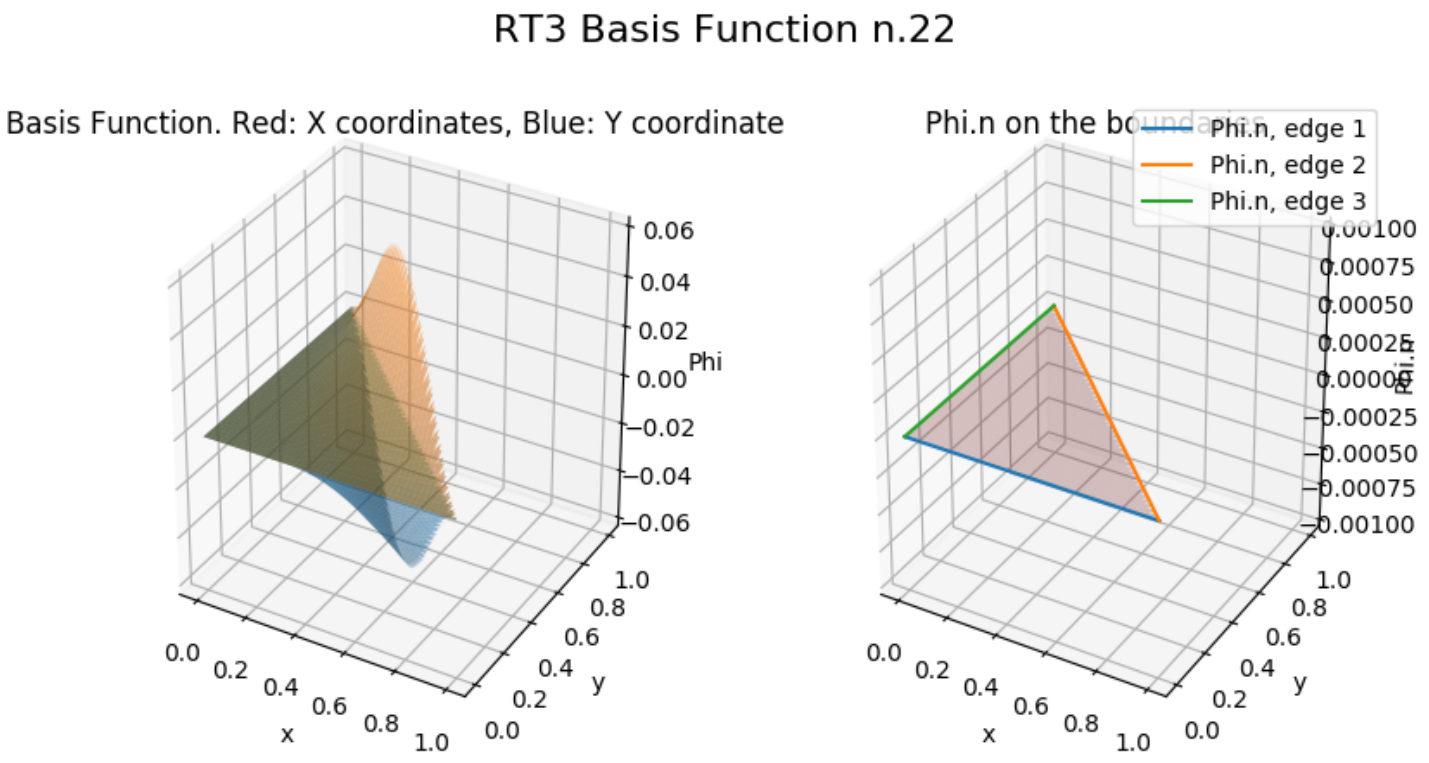

Fig. 18: Plot of one internal canonical basis function in the case $k = 3$.

The tuned set of local Lagrangian functions against the classical set of degrees of freedom of $RT_k(K)$ also observe the same behaviour. However, one can notice that the amplitude of the basis functions is much higher. This is due to the smallness of the coefficients in the transfer matrix.

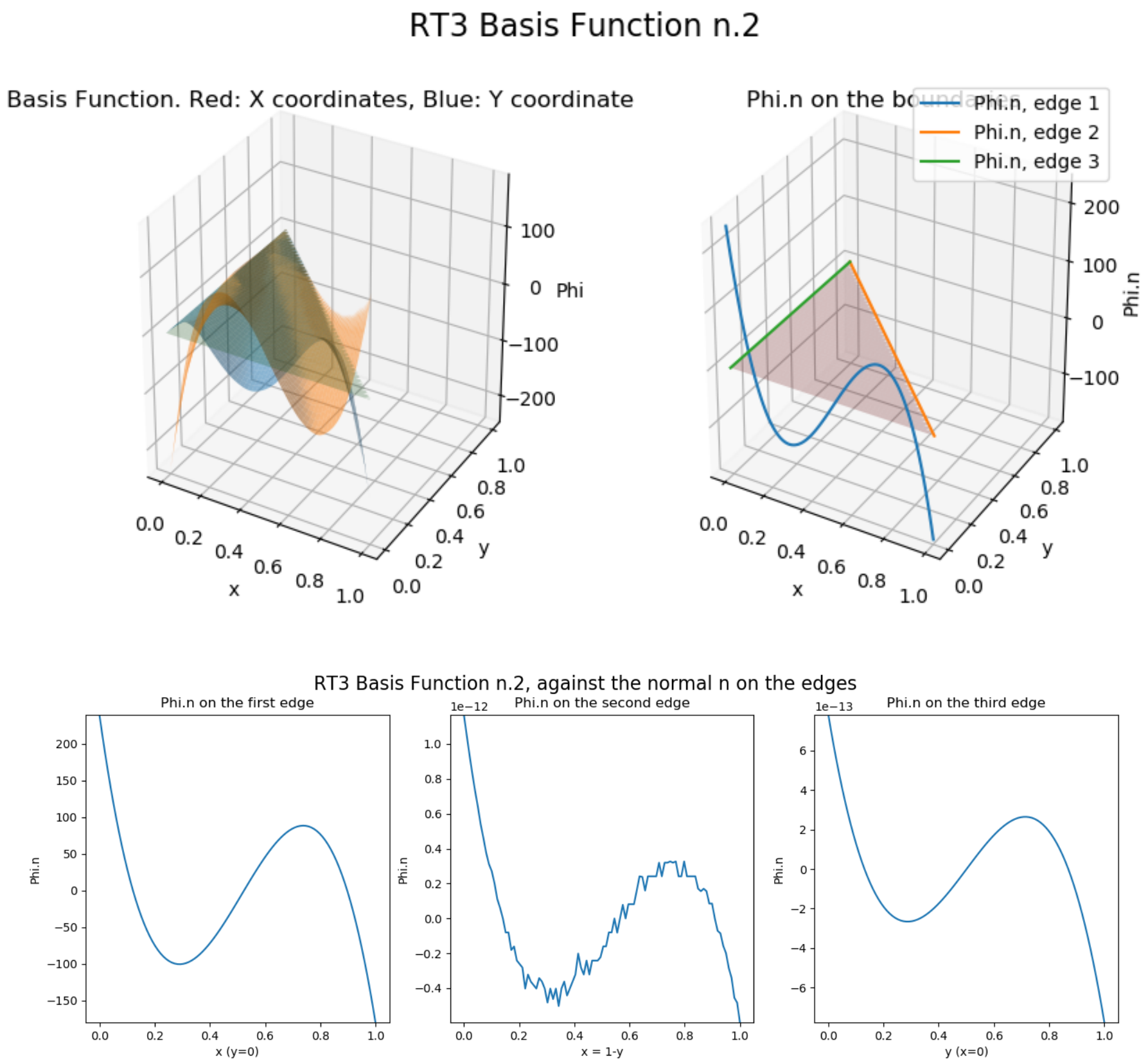

Fig. 19: Plot of the third tuned normal basis function in the case $k = 3$.

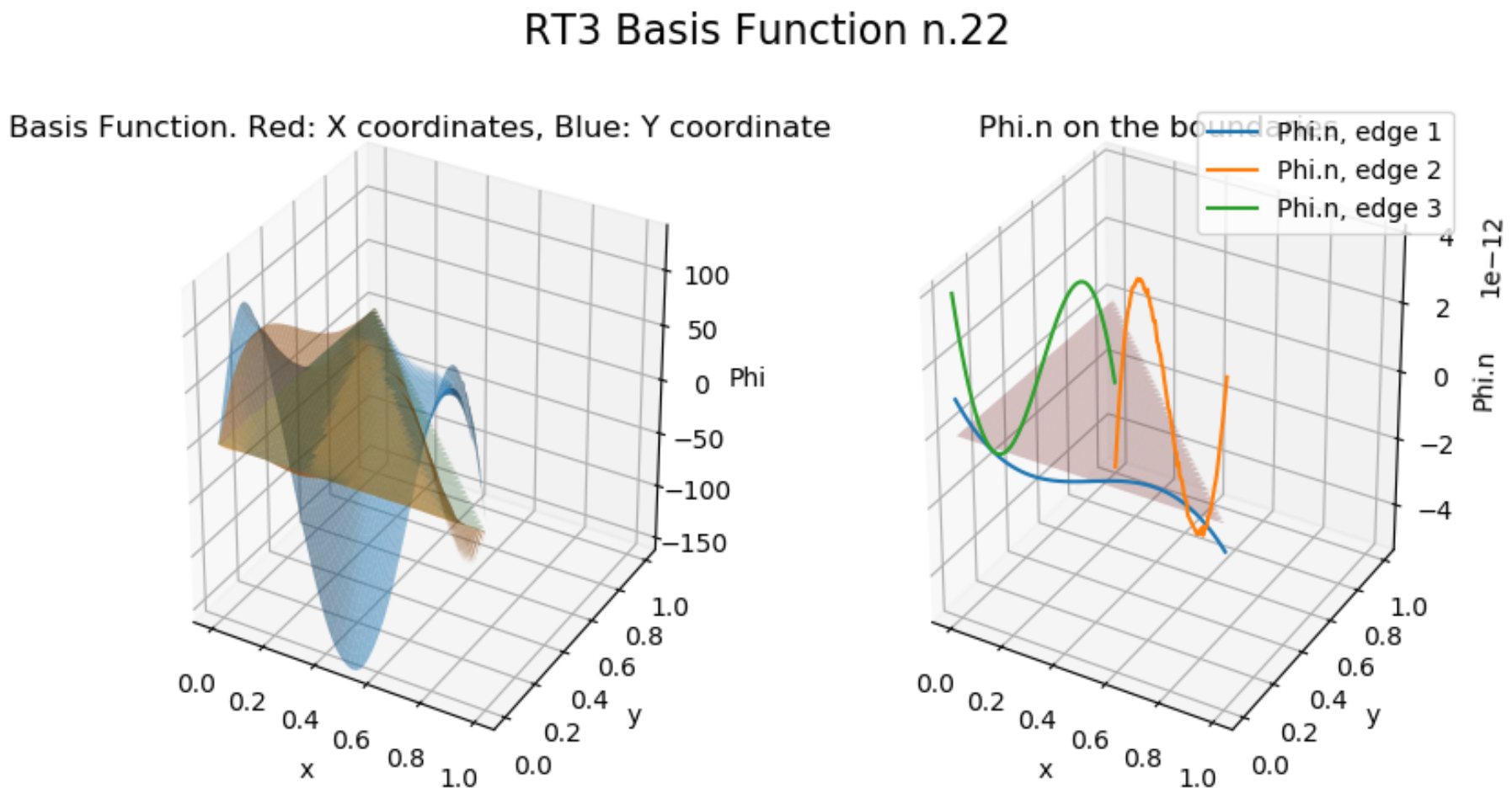

Fig. 20: Plot of one tuned canonical basis function in the case $k = 3$.

# 4 A classical extension to quadrilateral elements

When extending the definition of Raviart–Thomas elements to quadrilateral shapes, the same construction spirit applies up to some minor changes. We present here the resulting Raviart – Thomas quadrilateral elements before detailing a possible construction method for their corresponding basis functions.

## 4.1 A definition of Raviart – Thomas spaces for quads

**A first limitation arising from $\mathbb{P}^k$** The key property of the previously presented simplicial Raviart – Thomas elements is the adequacy between their dimension and the number of faces of the element. Indeed, the dimension of the discretisation space allows a natural split between internal and normal degrees of freedom. Furthermore, as the number of normal degrees of freedom is a multiple of the number of faces, their definition follow immediately.

Making use of the same definition (22) in the quadrilateral case would fail as the dimension split (26) does not make the number of faces of the quadrilateral emerge (here, $n = 2d \neq d + 1$).

However, changing the meaning given to the polynomial degree to work from the polynomial space $\mathbb{Q}_k$ as it is usually the case for quadrilateral elements solves this issue. There, the split between normal and internal degrees of freedom occurs naturally again. Let us detail the resulting several changes.

**Natural changes in the definition** Along with the change in the definition of a multivariate polynomial's degree, the definition of $|\alpha|$ in the differential operator (11) changes. As a consequence, the relation $|\alpha| = k+1$ does not generate the same set and the definition of the operator $\varepsilon^k$ changes accordingly. However, by analogy with the simplicial elements case one can still define a subspace of $H(\mathrm{div},\, K)$ as:

$$V = \{u \in (\mathbb{Q}_k(K))^d\,,\, \varepsilon^k u = 0\}.$$

It can be shown [13] that expliciting the kernel of $\varepsilon^k$ leads to the quadrilateral Raviart – Thomas space definition

$$RT_k(K) = (\mathbb{Q}_k(K))^d + x\,\mathbb{Q}_{[k]}(K), \tag{41}$$

where $x \in \mathbb{R}^d$ and $K \subset \mathbb{R}^{[d]}$ is a quadrilateral reference element.

In the same vein, as we use the polynomial space $\mathbb{Q}_k$ as part of the discretisation, we need to consider the discretisation $\mathcal{T}_k(\partial K)$ instead of $\mathcal{R}_k(\partial K)$ on the boundary of the element.

***Note.*** In the two dimensional case, we have $\mathcal{R}_k(\partial K) = \mathcal{T}_k(\partial K)$ as the boundary lies in one dimension, where the spaces $\mathbb{Q}_k$ and $\mathbb{P}_k$ merge. ▲

**An equivalent definition** Although we enjoy the same architecture of the space as in the simplicial case, using the definition (41) to build $H(\text{div})$ – conformal elements is not straightforward as the right hand side is not a direct sum anymore. Indeed, by the definition of the polynomial degree in the $\mathbb{Q}_k(K)$ space, we have $x\,\mathbb{Q}_{[k]} \cap (\mathbb{Q}_k(K))^d \neq \varnothing$ .

***Example.*** *Example of a non - empty intersection* To give an example, let us take the polynomial $\mathbb{Q}_2 \ni p\colon (x,\, y) \mapsto xy^2$. By definition, we have $\mathrm{x}\,p = (x^2y^2,\, xy^3)^T \in \mathrm{x}\mathbb{Q}_{[2]}(K)$, where we recall that $\mathbb{Q}_{[2]}(K)$ denotes the space of polynomial of degree exactly two; that is $\mathbb{Q}_{[2]}(K) = \{x^2,\, x^2y,\, x^2y^2,\, y^2,\, y^2x\}$. However, in the same time $\mathrm{x}\,p$ belongs to $\mathbb{Q}_2(K) \times \mathbb{Q}_3(K)$, thus intersecting with $(\mathbb{Q}_2(K))^2$. ♦

To overcome the problematics arising from the above definition, it has been shown [4] that the space (41) can be expressed as a Cartesian product of polynomial spaces $\mathbb{P}_{k_1,\,\cdots,\,k_d}$ for some integers $\{k_i\}_{i\in[\![1,\,d]\!]}$ through to the following definition.

**Definition 4.1** Raviart – Thomas spaces on quadrilaterals

The classical definition of Raviart – Thomas spaces on quadrilateral shapes reads

$$RT_k(K) = \bigtimes_{i=1}^{d} \mathbb{P}_{\zeta_i(k+1,\,\dots,\,k,\,k)}, \tag{42}$$

where $\{\zeta_i\}_i$ is the set containing the $d$ cyclic permutations of $\{k+1,\, \underbrace{k,\, \dots,\, k}_{d-1 \text{ times}}\}$.

For the sake of convenience, this last definition will be considered in all the following.

***Example.*** *Definition of $RT_k(K)$ in low dimension for quadrilateral elements*
In two dimension, the definition (42) reduces to:

$$RT_k(K) = \mathbb{P}_{k+1,\,k} \times \mathbb{P}_{k,\,k+1}, \tag{43}$$

while in three dimension we get

$$RT_k(K) = \mathbb{P}_{k+1,\,k,\,k} \times \mathbb{P}_{k,\,k+1,\,k} \times \mathbb{P}_{k,\,k,\,k+1}. \tag{44}$$

♦

We can now state its properties.

## 4.2 Properties of the quadrilateral $RT_k(K)$ space

### Dimension

From the definition (42) and the classical algebra results recalled in the section 2.1, we can derive the space dimension as follows.

$$\dim RT_k(K) = \sum_{i=1}^{d} \dim(\mathbb{P}_{\zeta_i\{k+1,\,k,\,\dots,\,k\}})$$

As the dimension of $\mathbb{P}_{\zeta_i\{k+1,\,k,\,\dots,\,k\}}$ is independent of $i$, we get:

$$\dim RT_k(K) = \sum_{i=1}^{d} (k+2)(k+1)^{d-1}$$

to finally achieve the following property.

**Property 4.2** Dimension of $RT_k(K)$ on quadrilateral shapes

$$\dim RT_k(K) = d(k+2)(k+1)^{d-1}$$

***Example.*** In the particular case of two and three dimensions, we get:

- For d = 2;

$$RT_k(K) = \mathbb{P}_{k+1,\,k} \times \mathbb{P}_{k,\,k+1}$$
$$\dim RT_k(K) = 2(k+1)(k+2)$$

- For d = 3;

$$RT_k(K) = \mathbb{P}_{k+1,\,k,\,k} \times \mathbb{P}_{k,\,k+1,\,k} \times \mathbb{P}_{k,\,k,\,k+1}$$
$$\dim RT_k(K) = 3(k+2)(k+1)^2.$$

In order to make easier later computational verifications, we provide here a dimension table for various dimensions and orders.

| k \ d | 2 | 3 | 4 |
|---|---|---|---|
| 0 | 4 | 6 | 8 |
| 1 | 12 | 36 | 96 |
| 2 | 24 | 108 | 432 |
| 3 | 40 | 240 | 1280 |

Tab. 1: Dimension table for quadrilateral $RT_k(K)$ spaces.

♦

### Divergence properties

As in the simplicial case, we can derive the following divergence property.

**Proposition 4.3**

For any $q$ belonging to $RT_k(K)$, it holds:

$$\begin{cases} \operatorname{div} q \in \mathbb{Q}_k \\ q \cdot n|_{\partial K} \in \mathcal{T}(\partial K) \end{cases} \tag{45}$$

***Proof.*** The proof is analogous to the one of the simplicial case. See the *Section 3.3* for details.

■

***Remark.*** The second relation in the properties (45) also reads:

$$q \cdot n_{f_i}|_{f_i} \in \mathbb{Q}_k(K) \quad \text{for each face } f_i \in \partial K, \quad i = [\![1,\, n]\!].$$

▲

Now that as we enjoy those two properties, we can tune the degrees of freedom on $\{K,\, RT_k(K),\, \mathcal{T}_k(\partial K)\}$ to define a conformal element.

## 4.3 Construction of the RT element

To build the Raviart – Thomas element on a quadrilateral of reference, we start by defining degrees of freedom ensuring the $H(\text{div},\, K)$ – conformity.

### 4.3.1 Definition of the degrees of freedom

As done in the simplicial case, we start by noticing the following.

• In the two dimensional case ($d = 2$), we have:

$$\begin{aligned}\dim RT_k(K) &= 2(k+1)(k+2)\\ &= 2k(k+1) + 4(k+1)\\ &= 2k(k+1) + 2\times 2(k+1).\end{aligned}$$

Noticing that $(k+1) = \dim \mathbb{P}_k(f)$ for any edge $f$ of $\partial K$ we can provide a dimension split that will be suitable to the design of the moment based degrees of freedom. Indeed, since we are here in the two dimensional case, $\mathbb{P}_k(f)$ matches $\mathbb{Q}_k(f)$ and $\mathbb{Q}_k(f)$ can be seen as the $\mathbb{Q}_k$ polynomial space built on a hypersurface of $K$ matching one of its boundaries. Coupled with algebra results, we derive:

$$\begin{aligned}\dim RT_k(K) &= \dim\left(\mathbb{P}_{k,k-1}\right)^2 + 2\times 2\dim \mathbb{Q}_k(f)\\ &= \dim\left(\mathbb{P}_{k,k-1}\right)^2 + 2\times 2\dim \mathbb{Q}_k(f)\\ &= \dim\left(\mathbb{P}_{k-1,k}\right)^2 + 2\times 2\dim \mathbb{Q}_k(f)\\ &= \dim\left(\mathbb{P}_{k-1,k}\times\mathbb{P}_{k,k-1}\right) + 2\times 2\dim \mathbb{Q}_k(f),\end{aligned}$$

allowing us to reshape the space in several ways, depending on the type of discretisation one wants to achieve (useful if you do no wish an homogeneous or symmetric space). Lastly, since $4 = 2\times 2$ matches the number of edges of a square, we get

$$\dim RT_k(K) = \dim\left(\mathbb{P}_{k-1,k}\times\mathbb{P}_{k,k-1}\right) + \dim \mathcal{T}_k(\partial K). \tag{46}$$

• For the three dimensional case $d = 3$, we get similarly:

$$\begin{aligned}\dim RT_k(K) &= 3(k+1)^2(k+2)\\ &= (k+1)^2(3k+6)\\ &= 6(k+1)^2 + 3k(k+1)^2\\ &= 2\times 3(k+1)^2 + 3k(k+1)^2\\ &= 2\times 3\dim \mathbb{Q}_k(f) + \dim\left(\mathbb{P}_{k-1,\,k,\,k}\times\mathbb{P}_{k,\,k-1,\,k}\times\mathbb{P}_{k,\,k,\,k-1}\right)\end{aligned}$$

One can notice that the constant $2\times 3$ corresponds to the number of faces of a cube. Therefore we can write

$$\dim RT_k(K) = \dim\left(\mathbb{P}_{k-1,\,k,\,k}\times\mathbb{P}_{k,\,k-1,\,k}\times\mathbb{P}_{k,\,k,\,k-1}\right) + \dim \mathcal{T}_k(\partial K).$$

- More generally for a dimension $d$ we can derive:

$$\begin{aligned}\dim RT_k(K) &= d(k+1)^{d-1}(k+2)\\ &= (k+1)^{d-1}(dk+2d)\\ &= 2d(k+1)^{d-1} + dk(k+1)^{d-1}\end{aligned}$$

$$\dim RT_k(K) = 2d \dim \mathbb{Q}_k(f) + \dim \bigtimes_{i=1}^{d} \mathbb{P}_{\zeta_i(\{k-1,\cdots,k,k\})}.$$

As for quadrilateral elements in dimension $d$, $2d \dim \mathbb{Q}_k(f) = \dim \mathcal{T}_k(\partial K)$, we get

$$\dim RT_k(K) = \dim \bigtimes_{i=1}^{d} \mathbb{P}_{\zeta_i(\{k-1,\cdots,k,k\})} + \dim \mathcal{T}_k(\partial K). \tag{47}$$

Therefore we can set up degrees of freedom respectively constructed from $\mathcal{T}_k(\partial K)$ and $\bigtimes_{i=1}^{d} \mathbb{P}_{\zeta_i(\{k-1,k,\cdots,k\})}$. With this in mind, we set

$$\Psi_k(K) = \bigtimes_{i=1}^{d} \mathbb{P}_{\zeta_i(\{k-1,\underbrace{k,\cdots,k}_{d-1 \text{ times}}\})}.$$

Furthermore, as in the simplicial case, as long as we ensure that the generated sets are free and free from each other we will be able to use the degrees of freedom built from $\mathcal{T}_k(\partial K)$ to enforce the $H(\text{div})$ – conformity. The definition of the degrees of freedom is then naturally given by the following relations.

**Proposition 4.4** Degrees of freedom for quadrilateral elements

For any $q$ belonging to $RT_k(K)$, the following relations form degrees of freedom on $K$.

$$q \mapsto \int_{\partial K} q \cdot n\, p_k \,\mathrm{d}\gamma(x), \qquad \forall p_k \in \mathcal{T}_k(\partial K) \tag{48a}$$

$$q \mapsto \int_{K} q \cdot p_k \,\mathrm{d}x, \qquad \forall p_k \in \Psi_k(K) \tag{48b}$$

where $\gamma$ represents the paths skimming the faces.

In practice, it is enough to provide those degrees of freedom by taking any basis of $\mathcal{T}_k(\partial K)$ and any basis of $\Psi_k(K)$ generating all the possible $p_k$. A good choice of those basis furnishes an easily computable set of degrees of freedom.

Typically, one can take respectively for (48a) and (48b)

$$q \mapsto \int_f q \cdot n\, p_m \,\mathrm{d}\gamma(x), \quad \forall f \in \partial K, \forall m \in [\![1, (k+1)^{d-1}]\!]$$

$$q \mapsto \int_K q \cdot \bigodot_{l=1}^{d} (0, \,\dots,\, 0, \underbrace{x_l^{\alpha_l}}_{i^{\text{th}}\text{ position}}, 0, \,\dots,\, 0)^T \,\mathrm{d}x, \quad {}_{\substack{\forall i \in [\![1,n]\!], \forall l \in [\![1,d,]\!] \\ \forall \alpha_l \in [\![1, (\zeta_i\{k,\dots,k,k+1\})_l]\!]}},$$

where $\{p_m\}$ forms any basis of $\mathbb{Q}_k(f)$ for any face $f$ of $\partial K$ and $\bigodot$ represents the Hadamard product.

***Example.*** *Two dimensional case* In the two dimensional case, the relations (48a) reduce to

$$q \mapsto \int_f q \cdot n\, x^m \,\mathrm{d}\gamma(x), \quad \forall f \in \partial K, \forall j \in [\![1,\, k+1]\!],$$

whereas the relations (48b) read

$$q \mapsto \int_K q \cdot \begin{pmatrix} x^{\alpha_x} \\ 0 \end{pmatrix} \circ \begin{pmatrix} y^{\alpha_y} \\ 0 \end{pmatrix} \,\mathrm{dx}, \quad \forall \alpha_x \leq k-1 \text{ and } \alpha_y \leq k$$

$$q \mapsto \int_K q \cdot \begin{pmatrix} 0 \\ y^{\alpha_x} \end{pmatrix} \circ \begin{pmatrix} 0 \\ x^{\alpha_y} \end{pmatrix} \,\mathrm{dx}, \quad \forall \alpha_x \leq k \text{ and } \alpha_y \leq k-1.$$

♦

Those applications are linear forms on $K$ and define moment based degrees of freedom. Furthermore, it is well - known [4] that it holds:

**Proposition 4.5**

For any $q \in RT_k(K)$,

$$\begin{cases} \int_{\partial K} q \cdot n\, p_k \,\mathrm{d}\gamma(x) = 0, & \forall p_k \in \mathcal{T}_k(\partial K) \\ \int_K q \cdot p_k \,\mathrm{d}x = 0, & \forall p_k \in \Psi_k(K) \end{cases} \Rightarrow \quad q = 0$$

The above property also reads $\mathrm{Ker}\{(48)\} = 0$. Thus, knowing that the dimension of $RT_k(K)$ matches the one of the set of degrees of freedom (48), it comes that the set $\{(48)\}$ is unisolvent for $RT_k(K)$. Thus, the set $\{(48)\}$ completely determines an element on $K$ for the space $RT_k(K)$.

We set $\{\sigma\} = \{(48)\}$ and define Raviart – Thomas quadrilateral elements by the triplet – domain, discretisation space, degrees of freedom – $\{K,\, RT_k(K),\, \{\sigma\}\}$.

***Remark.***
• As before, note that $q \in RT_k(K)$ is a vectorial function.
• The applications (48a) will be called normal degrees of freedom, and (48b) internal degrees of freedom. ▲

The representation of the Raviart–Thomas element on quadrilaterals is given in *Figure 21*.

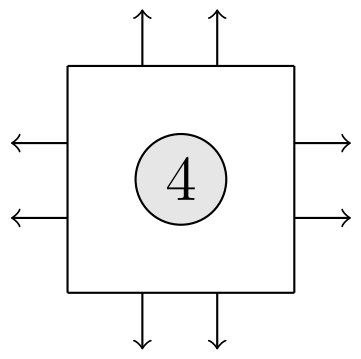

Fig. 21: Representation of $RT_k(K)$ on a square, here for $k = 2$

Here again, the arrows represent the normal degrees of freedom, and the circled number the number of internal moments.

### 4.3.2 Conformity

**Proposition 4.6**

The above-defined element $(K, RT_k(K), \{\sigma\})$ is $H(\text{div})$–conformal.

***Proof.*** The proof is analogous to the one given in the *Paragraph* 3.4.2, when replacing the space $\mathbb{P}_k$ by $\mathbb{Q}_k$. ■

### 4.3.3 Associated shape functions for the two dimensional case

The construction of the basis functions of the $RT_k(K)$ space for quadrilateral elements in two dimensions is done on the same spirit as for simplicial elements; that is find a generic basis for $RT_k(K)$ and tune it towards the degrees of freedom.

Therefore, as the adjustment of the basis functions applies identically here as in the *Section 3.5.2*, it will not be repeated. Likewise, the discussion on the transformation from a reference element to the considered one still holds. A method is recalled in the *Section 2.5.1* for interested readers. We just detail a possible construction method for some basis of $RT_k(K)$.

We consider the classical definition of $RT_k(K)$ over quadrangles;

$$\begin{aligned} RT_k(K) &= \mathbb{Q}_k(K) + x\mathbb{Q}_{[k]}(K) \\ &= \mathbb{P}_{k+1,\,k} \times \mathbb{P}_{k,\,k+1}. \end{aligned}$$

We want to define $\dim RT_k(K)$ functions that lie in $RT_k(K)$ and that are two by two free. To this end, let us work on the reference element presented in *Figure 22*. Note that as we restrict ourselves to the two dimensional case, the element has $n = 4$ edges.

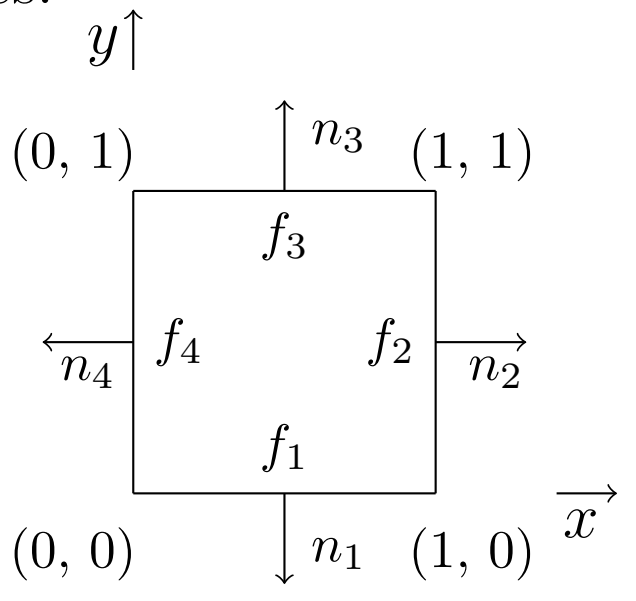

Fig. 22: Reference element for quadrangles

As before, the basis functions will be split into two categories.

- $\mathcal{N}$, a free set of functions of $RT_k(K)$ whose dimension is $\dim \mathcal{T}_k(\partial K) = (2d) \dim \mathbb{Q}_k(f)$. They should enforce the $H(\text{div})$ – conformity.
- $\mathcal{I}$, a free set of functions lying in $RT_k(K)$ whose dimension is $\dim \Psi_k(K)$. They should preserve the $H(\text{div})$ – conformity.

Note that here as well, each of the functions in $\mathcal{N}$ or $\mathcal{I}$ lie entirely in $RT_k(K)$. Doing so not only allows to cope with the two subspaces dimension, but also helps to guarantee the freeness of the set $\mathcal{N} \cup \mathcal{I}$.

**Let us first generate $\mathcal{N}$** We want to generate $\dim \mathcal{T}_k(\partial K)$ functions that are free from each other and that lie in $RT_k(K)$. They will be later tuned to be dual to the set of internal degrees of freedom (48a). Therefore, we follow the layout of those degrees of freedom and set $\dim \mathcal{T}_k(\partial K)/(2d) = \dim \mathcal{T}_k(\partial K)/n$ functions on each edge $f_i$, $i \in [\![1,\, n]\!]$.

***Note.*** Here, $\mathcal{T}_k(\partial K)$ is based on $\mathbb{Q}_k(f)$, which coincides with $\mathbb{P}_k(f)$ as we are on the boundaries of a two dimensional element. ▲

The construction is done in a similar spirit than in the simplicial elements case; that is defining one two-dimensional vector per edge belonging to $\mathbb{P}_{1,\,0} \times \mathbb{P}_{0,\,1}$ before stitching each to $(k+1)$ Lagrangian functions.

For the first three edges, the vectors are defined in the exact same way as in the simplicial case by setting

$$e_i = \begin{pmatrix} x + n_{ix} \\ y + n_{iy} \end{pmatrix}, \quad i = 1,\, 2,\, 3,$$

where $n_i = (n_{ix},\, n_{iy})^T$ is the normal vector to the edge $f_i$. However, defining the fourth vector in the same way would prevent the set $\{e_i\}_{i \in [\![1,\,4]\!]}$ to be free. Indeed, for any $n \geq 4$ it is always possible to find a linear combination between the elements of the set

$$\left\{ \begin{pmatrix} x + \alpha_i \\ y + \beta_i \end{pmatrix} \right\}_{i \in [\![1,\,n]\!]} \tag{49}$$

for any set $\{(\alpha_i,\, \beta_i)\}_{i \in [\![1,\,n]\!]}$ whose elements belong to $(\mathbb{R}^2)$. Therefore, we need to design the last vector as follows:

$$e_4 = \begin{pmatrix} |n_{4x}| + \mathrm{sign} n_{4x} x \\ |n_{4y}| + \mathrm{sign} n_{4y} y \end{pmatrix},$$

with the convention that the sign of 0 is positive. The importance here is to break the uniform sign coupling between the two variables $x$ and $y$ as for $n = 4$ vectors it is enough to prevent the linear interdependency of the vectors in (49). Therefore, according to the definition of the reference element this change in the vector definition may be done either for $e_1$ or for $e_4$.

Note however that any other vector built from the tuples $(x,\, -y)^T$ or $(-x,\, y)^T$ and belonging to $\mathbb{P}_{1,\,0} \times \mathbb{P}_{0,\,1}$ would be equally admissible. Furthermore, by a good choice of the vector to modify, this method is adaptable to any other quadrilateral reference element one may wish to use.

So defined, we enjoy the following property.

**Property 4.7**

The set $\{e_i\}_{i \in [\![1,\,4]\!]}$ is free within the element $K$.

***Proof.*** Let us show that no linear combination is possible within $\{e_i\}_{i \in [\![1,\,4]\!]}$ in the element $K$. We are looking for $\alpha$, $\beta$, $\gamma$ and $\tau$ four real numbers so that for any point $(x,\, y)$ in $K$ it holds;

$$\alpha e_1(x,\, y) + \beta e_2(x,\, y) + \gamma e_3(x,\, y) + \tau e_4(x,\, y) = 0.$$

When applied to the reference element (22) the set $\{e_i\}_{i\in[\![1,\,4]\!]}$ respectively reads

$$\left\{\begin{pmatrix}x+1\\ y\end{pmatrix},\,\begin{pmatrix}x\\ y-1\end{pmatrix},\,\begin{pmatrix}x\\ y+1\end{pmatrix},\,\begin{pmatrix}1-x\\ y\end{pmatrix}.\right\}.$$

It comes the system

$$\begin{cases}\alpha x+\beta x+\beta+\gamma x-\tau x+\tau=0\\ \alpha y-\alpha+\beta y+\gamma y+\gamma+\tau y=0,\end{cases}$$

which on the interior of $K$ reduces to

$$\begin{aligned}&\Leftrightarrow\begin{cases}-\tau+2\gamma-\tau=0\\ -\tau+2\gamma+\tau=0\\ \beta=-\tau\\ \gamma=\alpha\end{cases}\\ &\Leftrightarrow\quad \alpha=\beta=\gamma=\tau=0.\end{aligned}$$

Therefore the set is free within the element $K$.

■

We then set as previously $(k+1)$ functions on each edge with local Lagrangian property that are free from each other by defining the applications

$$(x,\,y)\mapsto l_{im}(x,\,y)e_i(e,\,y),\quad \forall i\in[\![1,\,2d]\!],\,\forall m\in[\![1,\,k+1]\!],$$

$l_{im}$ being the $m^{\text{th}}$ Lagrangian function generated from $(k+1)$ sampling points distributed along the edge $f_i$.

***Remark.*** Although the sampling points can be chosen freely on each edge, the Gauss-Legendre nodes are usually applied. ▲

All in all, we set:

$$\mathcal{N}=\{(x,\,y)\mapsto l_{i,\,m}(x,\,y)e_i(x,\,y)\}_{\substack{i\in[\![1,\,n]\!]\\ m\in[\![1,\,k+1]\!]}}\tag{50}$$

***Remark.*** *Important remark* With $\{e_i\}_{i\in[\![1,\,4]\!]}$ defined as such, the freeness of $\mathcal{N}$ is here guaranteed by the convexity of the element of reference. Indeed, when using homogeneous sampling rule across the edges, the edges $f_1$ and $f_3$ are similar in the view of generation of Lagrangian functions. Therefore, $\{l_{3,\,m}\}_{m\in[\![1,\,k+1]\!]}$ and $\{l_{1,\,m}\}_{m\in[\![1,\,k+1]\!]}$ are the same set, and thus not free from each other. However, by the convexity of the element no normal is identical, and in particular the set $\{e_1, e_3\}$ is free. Thus the set $\{l_{i,m} e_i\}_{i,m}$ contains two by two free functions as shown in the *Corollary 4.11*. ▲

***Remark 4.*** As in the simplicial case, if the sampling points generation rule is the same for each edge one can retrieve every $l_{ij}$ from a single function $l_i$ as in (31) but on the reference element given on *Figure 22* with the angles defined as in *Figure 23*.

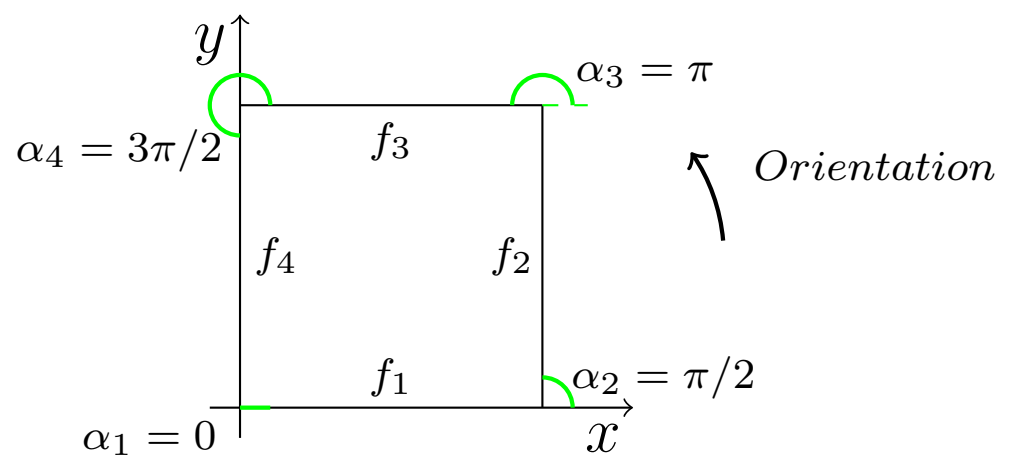

Fig. 23: Angle of the quadrilateral element

Note that here, only the third and fourth case of (31) occur. However, one would then have to be careful here when ordering the basis function to keep the coherence. Indeed, as here the angles $\alpha$ are all multiples of $\pi$ or $\pi/2$, if one wants to order the sampling points on the boundary element with a global indexing, the definition (31) should be modified as:

$$\begin{aligned} \tilde{l}_{i,m}(x, y) &= l_{\zeta(m)}(x), \quad \text{if } \alpha_i \equiv 0[\pi] \\ \tilde{l}_{i,m}(y) &= l_{\zeta(m)}(x), \quad \text{if } \alpha_i \equiv \pi/2[\pi] \end{aligned}$$

where $\zeta$ is a permutation of the sampling points $\{\mathrm{x}_m\}_m$ lying on the edge $f_i$. It will ensure that the points are correctly ordered with respect to the edges they belong to and the local orientation of the element (*see the Figure 24*). ▲

***Remark.*** In the case of a uniform sampling across the edges that is disregarding the orientation of the element (*see the Figure 26*), we have

$$\begin{aligned} l_{0,m}(x_l, y_l) &= l_j(x_l) = \delta_{lm} \\ l_{1,m}(x_l, y_l) &= l_j(y_l) = \delta_{lm} \\ l_{2,m}(x_l, y_l) &= l_j(x_l) = \delta_{lm} \\ l_{3,m}(x_l, y_l) &= l_j(y_l) = \delta_{lm} \end{aligned}$$

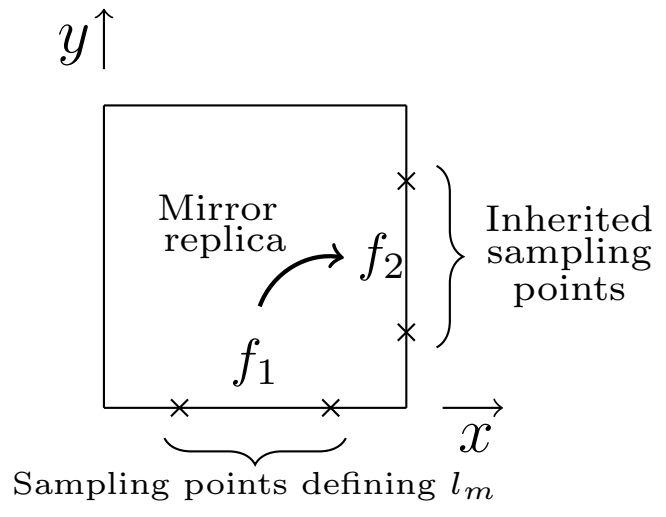

Fig. 25: Point mirroring example

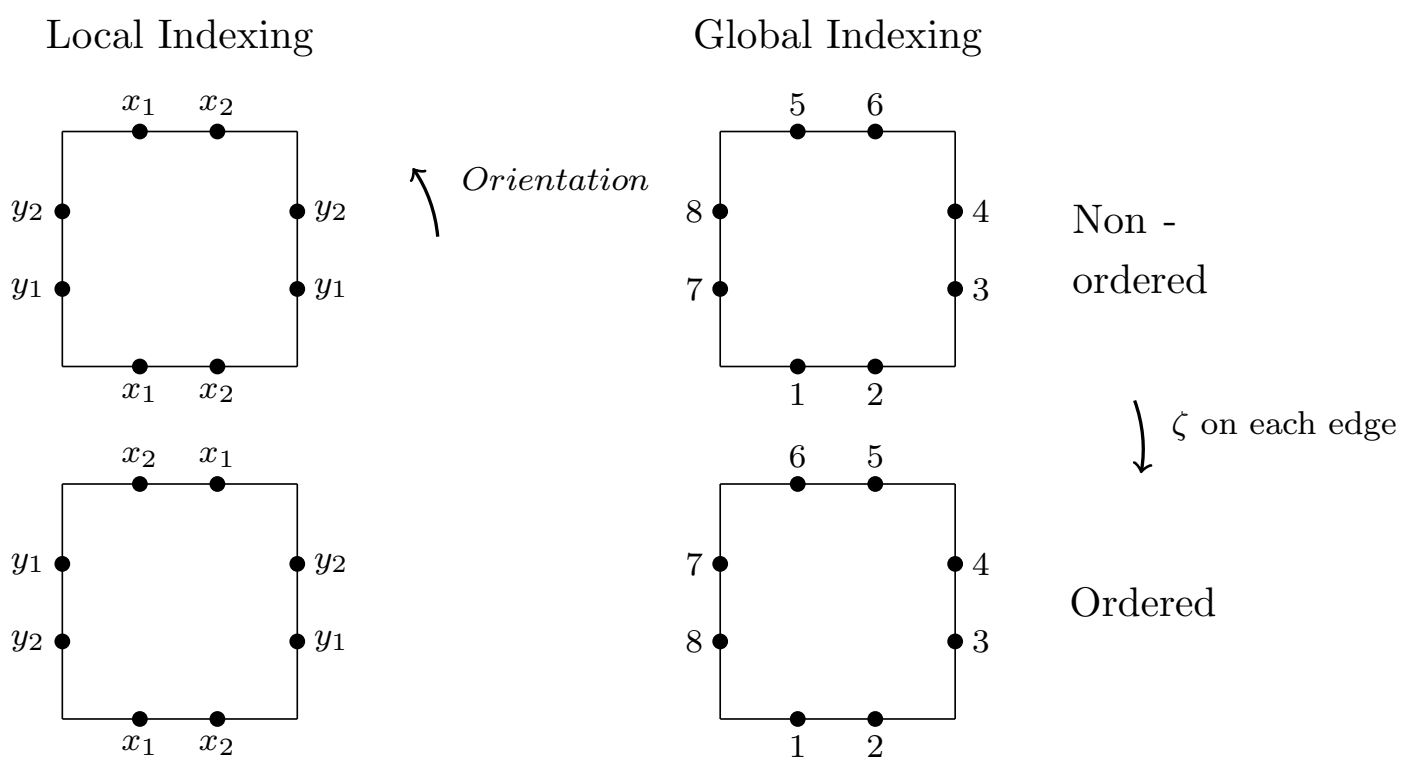

Fig. 24: Sampling point and labelling ordering

for some one - dimensional Lagrangian function $l_m$, $m \in [\![1,\, k+1]\!]$. Indeed, the edges of our reference element share the same length. Therefore, by the uniform distribution of the sampling points the Lagrangian functions will be identical, up to possibly flipping the $x$ and $y$ coordinates and reordering the points (*see the Figure 25*). The above relations follow immediately. ▲

We can now state some properties of $\mathcal{N}$.

**Property 4.8** Dimension of $\mathcal{N}$

$$\dim \mathcal{N} = 2d \dim \mathbb{Q}_k(f)$$

***Proof.*** By construction, we directly get:

$$\begin{aligned}\dim \mathcal{N} &= 4 \dim\{l_{i,\,m}\}_{m\in[\![1,\,k+1]\!]}\\ &= 2d \dim\{l_{i,\,m}\}_{m\in[\![1,\,k+1]\!]}\\ &= 2d(k+1)\\ \dim \mathcal{N} &= 2d \dim \mathbb{Q}_k(f),\end{aligned}$$

where the last line holds as we are in the two dimensional case.

■

**Proposition 4.9**

For any $q$ belonging to $\mathcal{N}$, $q$ belongs to $RT_k(K)$.

***Proof.*** Let us take any $q$ in $\mathcal{N}$. Then $q$ is of the form $q = l_{i,\,m}(x,\, y)e_i(x,\, y)$ for some $i \in [\![1,\, 4]\!]$ and some $m \in [\![1,\, k+1]\!]$.

If $i = 1$, 2 or 3, the same proof as in the simplicial case holds when simply replacing $\mathbb{P}_k(K)$ by $\mathbb{Q}_k(K)$ (up to some simplifications as the space is not made of a direct sum anymore).

For the case $i = 4$, one needs to consider the definition (42) of the $RT_k(K)$ space. Indeed, $q$ may be written as

$$q = \begin{pmatrix} |n_{4x}| + \text{sign} n_{4x} x \\ |n_{4y}| + \text{sign} n_{4y} y \end{pmatrix} l_{4,m}(x, y)$$

for some $m \in [\![1, k+1]\!]$. As the edge here is vertical we know that $l_{i,m} \in \mathbb{P}_{0,k} \subset \mathbb{P}_{k,k}$. Therefore, the application $(x, y) \mapsto (\text{sign}(n_{4x})x + |n_{4x}|)l_{4j}$ belongs to $\mathbb{P}_{k+1,k}$ and $(x, y) \mapsto (\text{sign}(n_{4x})y + |n_{4y}|)l_{4j}$ to $\mathbb{P}_{k,k+1}$. Finally, $q \in \mathbb{P}_{k+1,k} \times \mathbb{P}_{k,k+1} = RT_k(K)$.

■

As in the triangular case, we enjoy a local Lagrangian property.

**Property 4.10** Local Lagrangian property

For any $i \in [\![1, 4]\!]$ and for any sampling point $\mathrm{x}_l \in f_i$ that was used to generate the Lagrangian function set $\{l_{i,m}\}_{m \in [\![1, k+1]\!]}$, the set of functions

$$\{(x, y)^T \mapsto l_{i,m}(x, y) e_i(x, y)\}_{m \in [\![1, k+1]\!]}$$

enjoy the local Lagrangian property

$$l_{i,m}(x_l, y_l) e_i(x_l, y_l) \cdot n_i = c_i \delta_{\zeta(m)l}, \tag{51}$$

where $c_i \in \mathbb{R}$ is possibly vanishing and represents a constant depending on the shape of the reference element. The application $\zeta$ stands for some permutation function emphasizing the possible re-indexing of the Lagrangian functions (*see e.g. the Remark 4*).

***Proof.*** Let us fix some edge $f_i \in \partial K$, $i \in [\![1, 4]\!]$ and consider any $m \in [\![1, k+1]\!]$. Then, as $l_{i,m}$ is a scalar function one has:

$$(l_{i,m}(x_l, y_l) e_i(x_l, y_l)) \cdot n_i = l_{i,m}(x_l, y_l)(e_i(x_l, y_l) \cdot n_i). \tag{52}$$

And by the definition of the vector $e_i$, we get:

$$\begin{aligned} l_{i,m}(x_l, y_l) e_i(x_l, y_l) \cdot n_i &= \left[ \begin{pmatrix} x_l + n_{ix} \\ y_l + n_{iy} \end{pmatrix} \cdot \begin{pmatrix} n_{ix} \\ n_{iy} \end{pmatrix} \right] l_{i,m}(x_l, y_l) \\ &= (x_l n_{ix} - n_{ix}^2 + y_l n_{iy} - n_{iy}^2) l_{i,m}(x_l, y_l) \end{aligned}$$

$$= \underbrace{(x_l n_{ix} + y_l n_{iy} - 1)}_{= c_i \text{ as } \mathrm{x}_l \in f_i} l_{i,m}(x_l, y_l)$$

$$l_{i,m}(x_l, y_l) e_i(x_l, y_l) \cdot n_i = c_i \delta_{\zeta(m)l}, \quad \text{as } \mathrm{x}_l \text{ is a sampling point on } f_i.$$

■

***Remark.*** Even if we cannot directly enjoy a global Lagrangian property, we can still derive the following statement. For any $\mathrm{x} \in f_j$ and any sampling points $\{\mathrm{x}_l\}_l \in f_j$, it holds

$$l_{i,m}(x, y) e_i(x, y) \cdot n_j = c_i l_{i,m}(x, y) + c_{\zeta(j),\{\mathrm{x}_l\}_l} \tag{53}$$

for come constants $c_i$ and $c_{\zeta(j),\{x_l\}_l}$ depending respectively on the edge $f_i$ the function is generated from, and on the edge $f_j$ the function is tested against plus the layout of the sampling points lying on $f_j$.

When one uses a homogeneous sampling points distribution which disregards the element's local orientation (*see the Figure26*), it comes down to:

$$l_{i,m}(x_l, y_l) e_i(x_l, y_l) \cdot n_j|_{f_j} = c_i \delta_{\zeta(m)l} + c_{\zeta(j),\{x_l\}_l}. \tag{54}$$

Note that in both cases the constants may vanish. Indeed, we have for any $i, j \in [\![1, n]\!]$ and $m \in [\![1, k+1]\!]$:

$$l_{i,m}(x, y) e_i(x, y) \cdot n_j|_{f_j} = (x n_{jx} + y n_{jy} + (n_{ix} n_{jx} + n_{iy} n_{jy})) l_{i,m}(x, y).$$

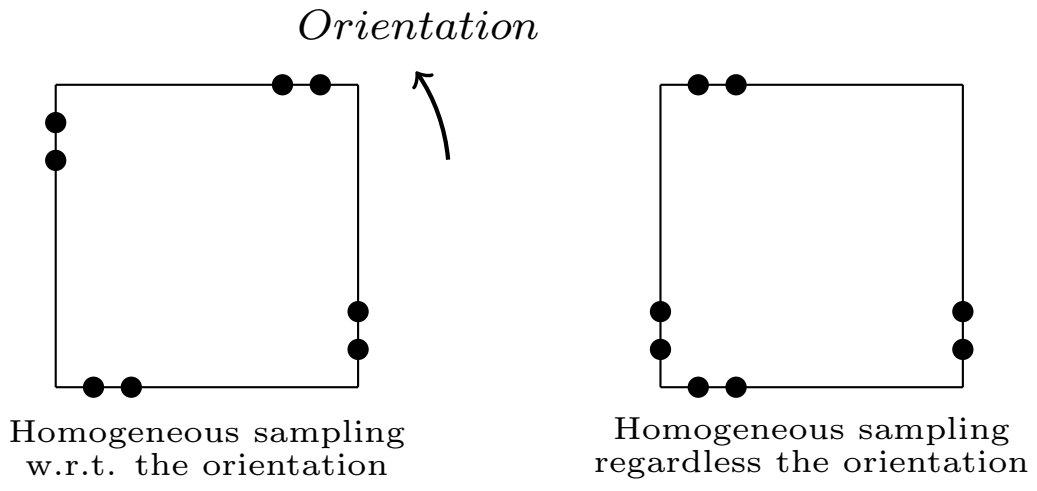

Fig. 26: Homogeneous sampling. Left: Taking care of the orientation. Right: Disregarding the orientation.

From there, we only detail the computations for the functions that have been generated from the first edge. The other ones can be derived in a similar way. Therefore we set $i = 1$ and get:

$$l_{1,m}(x, y) e_1(x, y) \cdot n_j|_{f_j} = (x n_{jx} + y n_{jy} + (0 - n_{jy})) l_{1,m}(x, y).$$

Thus, for $\mathrm{x} \in f_1$ it reads:

$$\begin{aligned} l_{1,m}(x,\, y)e_1(x,\, y) \cdot n_1|_{f_1} &= (0 + 0 + (0 - (-1)))l_{1,m}(x,\, y) \\ &= l_{1,m}(x,\, y), \end{aligned}$$

while when applied to the sampling points $\mathrm{x}_l$ lying on $f_1$ we get:

$$l_{1,m}(x_l,\, y_l)e_1(x_l,\, y_l) \cdot n_1|_{f_1} = \delta_{ml}.$$

Similarly, for any $\mathrm{x} \in f_2$ we have:

$$\begin{aligned} l_{1,m}(x,\, y)e_1(x,\, y) \cdot n_2|_{f_2} &= (x + 0 + (0 + 0))l_{1,m}(x,\, y) \\ &= l_{1,m}(x), \end{aligned}$$

which when applied to the sampling points $\mathrm{x}_l = (1,\, y_l)$ lying on $f_2$ comes down to:

$$\begin{aligned} l_{1,m}(x_l,\, y_l)e_1(x_l,\, y_l) \cdot n_2|_{f_2} &= l_{1,m}(x_l) \\ l_{1,m}(x_l,\, y_l)e_1(x_l,\, y_l) \cdot n_2|_{f_2} &= l_{1,m}(1) \in \mathbb{R}. \end{aligned}$$

For $\mathrm{x} \in f_3$, we have

$$\begin{aligned} l_{1,m}(x,\, y)e_1(x,\, y) \cdot n_3|_{f_3} &= (0 + y + (0 + 1))l_{1,m}(x,\, y) \\ &= 2l_{1,m}(x). \end{aligned}$$

Applied to the sampling points $\mathrm{x}_l$ lying on $f_3$, there holds:

$$l_{1,m}(x_l,\, y_l)e_1(x_l,\, y_l) \cdot n_3|_{f_3} = 2\delta_{ml} \in \mathbb{R}.$$

Lastly, for any $\mathrm{x} \in f_4$ it reads:

$$\begin{aligned} l_{1,m}(x,\, y)e_1(x,\, y) \cdot n_4|_{f_4} &= (0 + 0 + (0 + 1))l_{1,m}(x,\, y) \\ &= 0. \end{aligned}$$

Note that even if the functions emerging from the fourth edge is not of the similar construction type, its structure is the same. Therefore, the computations end in the same way. ▲

***Note.*** Even though at this point only a local Lagrangian property is satisfied, a global Lagrangian property will emerge after tuning against the set of degrees of freedom (48a). ▲

**Corollary 4.11**

The functions of $\mathcal{N}$ are linearly independent.

***Proof.*** By the local Lagrangian property the functions share, the set $\{l_{im}\}_{m\in[\![1,\,k+1]\!]}$ is free for each $i \in [\![1,\,4]\!]$. Furthermore, as shown before the functions $\{e_i\}_{i\in[\![1,\,4]\!]}$ are free strictly within $K$. Thus, the same discussion as in the simplicial case applies also here (*see the Proof 3.5.1*), and we have by construction that the set

$$\mathcal{N} = \{l_{im} e_i\}_{\substack{i\in[\![1,\,4]\!] \\ m\in[\![1,\,k+1]\!]}}$$

is free on $K$. Note that as the set $\{e_i\}_i$ is still free, no change in the demonstration is required despite the change in the definition of $e_4$.

■

***Remark.*** It is also possible to create directly a set of basis functions enjoying global Lagrangian properties, that is

$$l_{i,m}(x_l,\, y_l) e_i(x_l,\, y_l) \cdot n_j|_{f_j} = \delta_{sl}\delta_{ij},$$

where $\mathrm{x}_l$ is any sampling on $f_j$ with a global index and where $s = (i-1)(k+1) + m$ is the global index of the function. One just has to look for four vectors $\tilde{e_1}$, $\tilde{e_2}$, $\tilde{e_3}$ and $\tilde{e_4}$ of the form $(ax+b,\, cx+d)^T$ such that $\tilde{e_i} \cdot n_j = \delta_{im}$ for any $i \in [\![1,\,4]\!]$ and any $m \in [\![1,\,k+1]\!]$. The shape of the vectors being imposed by the conditions on the polynomial degree, we only have to select the constants $a$, $b$, $c$ and $d$.

We thus retrieve four systems of four equations. Let us detail the computations for the first vector. Its scalar product with the normal should be one on the first edge and zero on the others. All other vectors can be retrieved in a similar way. The systems one has to fulfil then reads:

$$\begin{cases} \begin{pmatrix} ax+b \\ cy+d \end{pmatrix} \cdot \begin{pmatrix} 0 \\ -1 \end{pmatrix} = 1 & \text{for the first edge } f_1, \\ \begin{pmatrix} ax+b \\ cy+d \end{pmatrix} \cdot \begin{pmatrix} 1 \\ 0 \end{pmatrix} = 0 & \text{for the second edge } f_2, \\ \begin{pmatrix} ax+b \\ cy+d \end{pmatrix} \cdot \begin{pmatrix} 0 \\ 1 \end{pmatrix} = 0 & \text{for the third edge } f_3, \\ \begin{pmatrix} ax+b \\ cy+d \end{pmatrix} \cdot \begin{pmatrix} -1 \\ 0 \end{pmatrix} = 0 & \text{for the fourth edge } f_4. \end{cases}$$

Using the values of $x$ and $y$ on each edge, the system reduces to:

$$\begin{cases} -d = 1 & \text{since on } f_1,\ y = 0 \\ a + b = 0 & \text{since on } f_2,\ x = 1 \\ c + d = 0 & \text{since on } f_3,\ y = 1 \\ -b = 0 & \text{since on } f_4,\ x = 0. \end{cases}$$

We finally get that $a = b = 0$, $c = 1$, $d = -1$. Repeating the process for the other vectors, we retrieve

$$\tilde{e_1} = \begin{pmatrix} 0 \\ y-1 \end{pmatrix} \quad \tilde{e_2} = \begin{pmatrix} x \\ 0 \end{pmatrix} \quad \tilde{e_3} = \begin{pmatrix} 0 \\ y \end{pmatrix} \quad \tilde{e_4} = \begin{pmatrix} x-1 \\ 0 \end{pmatrix}.$$

Fortunately, those vectors are still free from each other. We can then plug them in the above construction in place of the set $\{e_i\}_{i\in[\![1,\,n]\!]}$ and get the desired basis functions in a similar way. Furthermore, as they enjoy a global Lagrangian property, the proof of freeness of $\mathcal{N}$ is immediate. ▲

**Let us now generate $\mathcal{I}$** We now want to generate $\dim(\mathbb{P}_{k-1,\,k} \times \mathbb{P}_{k,\,k-1})$ functions lying in $RT_k(K)$ that are free from each other and free with $\mathcal{N}$. They will later be tuned to be dual to the set of degrees of freedom (48b).

***Note.*** For $k = 0$, $\dim(\mathbb{P}_{k-1,\,k} \times \mathbb{P}_{k,\,k-1}) = 0$. The set $\mathcal{I}$ is then automatically designed as the empty set and there exists no internal basis function. ▲

To build this set of functions for any $k > 0$, we use the following two facts.

- A sufficient condition to ensure that the set will be free with $\mathcal{N}$ is to ask every function $p$ in $\mathcal{I}$ to verify $p \cdot n|_{\partial K} = 0$.
- It holds $\dim(\mathbb{P}_{k-1,\,k} \times \mathbb{P}_{k,\,k-1}) = 2\dim(\mathbb{P}_{k-1,\,k})$.

Therefore, to construct $\mathcal{I}$ we first generate two free vectors of $\mathbb{P}_{2,0} \times \mathbb{P}_{0,2}$ that will enforce the desired property of $q \cdot n|_{\partial K} = 0$. Then, we stitch to them a set of $\dim(\mathbb{P}_{k-1,\,k})$ functions lying in $\mathbb{P}_{k-1,\,k} \times \mathbb{P}_{k,\,k-1}$ to cope with the definition of $RT_k(K)$ and the desired dimension of the subspace $\mathcal{I}$.

To define the two initial vectors, we select any vector function $(x,\,y) \mapsto (f(x,\,y),\,g(x,\,y))^T$ such that $f(x,\,y)n_{ix} + g(x,\,y)n_{iy} = 0$ for each edge $f_i$ of $\partial K$, and where $f$ and $g$ are two given scalar polynomial functions of $(x,\,y)$ that are regular enough.

As we are in the quadrilateral case and as the edges of the reference element do not couple the coordinates, we can derive it easily by setting

$$e_{n+1} = \begin{pmatrix} x(x-1) \\ y(y-1) \end{pmatrix} \quad \text{and} \quad e_{n+2} = \begin{pmatrix} x(1-x) \\ y(y-1) \end{pmatrix}. \tag{55}$$

**Property 4.12**

For any $i = 1, 2$, $e_{n+i} \in \mathbb{P}_{2,0} \times \mathbb{P}_{0,2}$ and $e_{n+i} \cdot n|_{\partial K} = 0$.

***Proof.*** We detail the proof only for $e_{n+1}$. The case $e_{n+2}$ can be done in the exact same way.

We have by construction that for any $k \geq 1$ $e_{n+1} \in \mathbb{P}_{2,0} \times \mathbb{P}_{0,2} \subset RT_k(K)$. Furthermore, it holds on the edges:

$$\begin{aligned} e_{n+1}(x, 0) \cdot \begin{pmatrix} 0 \\ -1 \end{pmatrix} &= \begin{pmatrix} x(x-1) \\ 0 \times (0-1) \end{pmatrix} \cdot \begin{pmatrix} 0 \\ -1 \end{pmatrix} = 0 \\ e_{n+1}(x, 1) \cdot \begin{pmatrix} 0 \\ 1 \end{pmatrix} &= \begin{pmatrix} x(x-1) \\ 1 \times (1-1) \end{pmatrix} \cdot \begin{pmatrix} 0 \\ 1 \end{pmatrix} = 0 \\ e_{n+1}(0, y) \cdot \begin{pmatrix} -1 \\ 0 \end{pmatrix} &= \begin{pmatrix} 0 \times (0-1) \\ y \times (y-1) \end{pmatrix} \cdot \begin{pmatrix} -1 \\ 0 \end{pmatrix} = 0 \\ e_{n+1}(1, y) \cdot \begin{pmatrix} 1 \\ 0 \end{pmatrix} &= \begin{pmatrix} 1(1-1) \\ y \times (y-1) \end{pmatrix} \cdot \begin{pmatrix} 1 \\ 0 \end{pmatrix} = 0, \end{aligned}$$

which immediately concludes the above assertion.

■

***Note.*** Note that so far, no coupling between $x$ and $y$ is allowed componentwise in the definition of the two vectors because of the normals layout. The above proof furnishes an illustration. ▲

***Remark.*** Two such vectors could have been designed here only because for a quadrilateral of the reference element's shape there exists a Gröbner reduction to a polynomial set vanishing on the edges of the element whose degree is less than two. This will not hold anymore for the cases $n \geq 5$. ▲

**Property 4.13**

The set $\{e_{n+1}, e_{n+2}\}$ is free on the interior of $K$.

***Proof.*** The only non-trivial real numbers $A$, $B$ such that

$$Ae_{n+1} + Be_{n+2} = 0$$

can be found if $x = 0$, $y = 0$, $y = 1$ or $x = 1$, which are exactly the edges of $K$. Indeed, let us look for two constants $A$ and $B$ such that

$$Ae_{n+1} + Be_{n+2} = 0$$

By definition we have:

$$\begin{aligned} &\begin{cases} Ax(x-1) + Bx(1-x) = 0 \\ Ay(y-1) + By(y-1) = 0 \end{cases} \\ \Rightarrow &\begin{cases} Ax^2 - Ax - Bx^2 + Bx = 0 \\ (A+B)y(y-1) = 0 \end{cases} \\ \Rightarrow &\begin{cases} (A-B)x^2 + (B-A)x = 0 \\ (A+B)y(y-1) = 0 \end{cases} \end{aligned}$$

Thus, the only possibilities to have non - trivial solutions for $A$ and $B$ arise when

$$\begin{cases} y = 0 \\ x = 1 \end{cases} \quad \text{or} \quad \begin{cases} y = 0 \\ x = 0 \end{cases} \quad \text{or} \quad \begin{cases} y = 1 \\ x = 0 \end{cases} \quad \text{or} \quad \begin{cases} y = 1 \\ x = 1 \end{cases},$$

which are exactly matching the edges of the reference element. Therefore, the set $\{e_{n+1}, e_{n+2}\}$ is free on the strict interior of $K$.

■

As mentioned at the top of this paragraph, *page 94*, we can now define the basis by stitching to those two vectors any basis $\mathcal{B}_{k-1}(K) = \{b_m\}_m$ of $\mathbb{P}_{k-1,k}(K)$ and define the functions

$$p_j\colon (x, y) \mapsto e_{n+1}(x, y) \circ \begin{pmatrix} b_m(x, y) \\ b_m(y, x) \end{pmatrix}, \quad {\scriptstyle m \in [\![1, k(k+1)]\!] \atop j \in [\![1, k(k+1)]\!]} \tag{56a}$$

$$p_j\colon (x, y) \mapsto e_{n+2}(x, y) \circ \begin{pmatrix} b_m(x, y) \\ b_m(y, x) \end{pmatrix}, \quad {\scriptstyle m \in [\![1, k(k+1)]\!] \atop j \in [\![1+k(k+1), 2k(k+1)]\!]}, \tag{56b}$$

where $\circ$ denotes the Hadamard product. By construction, we immediately get the following property.

**Proposition 4.14**

For any $j \in [\![1,\, 2k(k+1)]\!]$, the above defined function $p_j$ belongs to $RT_k(K)$.

***Proof.*** We have for any $j \in [\![1,\, k(k+1)]\!]$ that

$$p_j(x,\, y) = \begin{pmatrix} x(x-1)b_j(x,\, y) \\ y(y-1)b_j(y,\, x) \end{pmatrix},$$

with $(x,\, y) \mapsto b_j(x,\, y) \in \mathbb{P}_{k-1,\, k}$ and $(x,\, y) \mapsto b_j(y,\, x) \in \mathbb{P}_{k,\, k-1}$. Thus, $p_j(x,\, y) \in \mathbb{P}_{k+1,k} \times \mathbb{P}_{k,k+1} = RT_k(K)$ for all $j \in [\![1,\, k(k+1)]\!]$.

We get in the same way $p_j(x,\, y) \in \mathbb{P}_{k+1,k} \times \mathbb{P}_{k,k+1}$ for all $j \in [\![k(k+1),\, 2k(k+1)]\!]$. Therefore, $\{p_j\}_{j\in[\![1,\, 2k(k+1)]\!]} \subset RT_k(K)$.

■

***Note.*** By selecting any basis $\{b_j\}_{j\in[\![1,\, 2k(k+1)]\!]}$ of $\mathbb{P}_{k-1,\, k}(K)$, the set

$$\left\{(x,\, y) \mapsto \begin{pmatrix} b_j(x,\, y) \\ b_j(y,\, x) \end{pmatrix}\right\}_{j\in[\![1,\, 2k(k+1)]\!]}$$

is naturally free. ▲

We therefore set

$$\mathcal{I} = \{p_j\}_{j\in[\![1,\, 2k(k+1)]\!]}$$

and derive the following properties.

**Proposition 4.15**

For any $p$ belonging to span $\{\mathcal{I}\}$, $p$ belongs to $RT_k(K)$.

***Proof.*** For any $p$ belonging to span $\{\mathcal{I}\}$, $p$ can be written as a linear combination of elements in $\{p_j\}_{j\in[\![1,\, 2k(k+1)]\!]}$, which are all elements of $RT_k(K)$. Therefore, $p$ belongs to $RT_k(K)$.

■

**Property 4.16** Dimension

$\dim \mathcal{I} = \dim(\mathbb{P}_{k-1,k} \times \mathbb{P}_{k,k-1})$

***Proof.*** We directly have by construction and algebra result that $\dim \mathcal{I} = \dim(\{p_j\}_{j\in[\![1,\, 2k(k+1)]\!]}) = 2k(k+1) = \dim(\mathbb{P}_{k-1,k} \times \mathbb{P}_{k,k-1})$.

■

**Proposition 4.17**

For any $p$ belonging to $\mathcal{I}$, all edges $f$ of $\partial K$ and any point $(x, y)$ on $f$, it holds

$$p(x, y) \cdot n_f = 0.$$

***Proof.*** We recall that we are working on our reference element represented in *Figure 22*. We know that any $p$ belonging to $\mathcal{I}$ can be written under the form

$$p = \sum_{j=1}^{2k(k+1)} \alpha_j p_j$$

for some given real coefficients $\{\alpha_j\}_j$. By definition of $\{p_j\}_j$, the variables in each vector $e_{n+i}$, $i = 1, 2$ are uncoupled. The layout of the normals on the reference element does not involve a coupling of the two coordinates either and it is directly obtained that the properties of $\{e_{n+i}\}_{i=1,2}$ transfer to $\{p_j\}_{j\in[\![1,\,2k(k+1)]\!]}$. Thus $p_j \cdot n|_{\partial K} = 0$ for any $j \in [\![1, 2k(k+1)]\!]$.

More explicitly, one has for any generic edge $f$ and its corresponding normal vector $n_f = (n_x, n_y)^T$;

$$p_j \cdot n_f|_f = \begin{pmatrix} e_{n+i,1}(x, y) b_j(x, y) \\ e_{n+i,2}(x, y) b_j(y, x) \end{pmatrix} \cdot \begin{pmatrix} n_x \\ n_y \end{pmatrix}.$$

By example on the edge $f_1$ we get

$$p_j \cdot n_1|_{f_1} = \begin{pmatrix} e_{n+i,1}(x, 0) b_j(x, 0) \\ e_{n+i,2}(x, 0) b_j(0, x) \end{pmatrix} \cdot \begin{pmatrix} 0 \\ -1 \end{pmatrix}.$$

Here, the term $e_{n+i,2}(x, 0)$ vanishes for both $i = 1$ and $i = 2$. Therefore even if the term $b_j(0, x)$ actually involves variable coupling, the contribution is shut down.

Lastly, the definition of the normal vector to the edge $f_1$ does not bring variable coupling either and make the whole term vanish. As any other edge can be derived accordingly, the claim follows.

■

**Proposition 4.18**

The set $\mathcal{I}$ is free in the interior of $K$.

***Proof.*** It comes directly from the fact that $b_j$ belongs to $\mathbb{P}_{k-1,\,k}$, and that $(x,\,y) \mapsto b_j(y,\,x)$ belongs to $\mathbb{P}_{k,\,k-1}$. Indeed, we know that

$$e_{n+1} = \begin{pmatrix} x(x-1) \\ y(y-1) \end{pmatrix} \in \mathbb{P}_{2,\,0} \times \mathbb{P}_{0,\,2},$$

and that $e_{n+2} \in \mathbb{P}_{2,\,0} \times \mathbb{P}_{0,\,2}$ in the same way. Furthermore, we can assume without loss of generality that the set $\{b_j\}_{j \in [\![1,\,2k(k+1)]\!]}$ is the canonical basis of monomials $\{1,\, x,\, y,\, xy,\, x^2,\, y^2x,\, x^2,\, y^2\, \cdots,\, x^{k-1}y^k,\, x^{k-1},\, y^k\}$. The vectors $(x,\,y) \mapsto (b_j(x,\,y),\, b_j(y,\,x))^T$ can then be written under the form:

$$\left\{ \begin{pmatrix} 1 \\ 1 \end{pmatrix}, \begin{pmatrix} x \\ y \end{pmatrix}, \begin{pmatrix} xy \\ xy \end{pmatrix}, \begin{pmatrix} x^2y \\ y^2x \end{pmatrix}, \cdots, \begin{pmatrix} x^{k-1}y^k \\ x^ky^{k-1} \end{pmatrix}, \begin{pmatrix} x^{k-1} \\ y^{k-1} \end{pmatrix}, \begin{pmatrix} y^k \\ x^k \end{pmatrix} \right\},$$

and the set $\{p_j\}_j$ turns to write

$$\left\{ \begin{pmatrix} x(x-1) \\ y(y-1) \end{pmatrix}, \begin{pmatrix} x^2(x-1) \\ y^2(y-1) \end{pmatrix}, \begin{pmatrix} x^2y(x-1) \\ xy^2(y-1) \end{pmatrix}, \cdots, \right.$$
$$\left. \cdots, \begin{pmatrix} x^k(x-1) \\ y^k(y-1) \end{pmatrix}, \begin{pmatrix} y^k(x-1) \\ x^ky(y-1) \end{pmatrix} \right\}.$$

Since all the vectors have two by two distinct orders with respect to the $x$ and $y$ coordinates (*i.e.* it is impossible to find a same linear combination between the first coordinates of vectors belonging to $\{p_j\}_j$ and between their second coordinates), there exists a non-trivial solution on $\{\alpha_j\}_j \in \mathbb{R}$ for

$$\sum_j \alpha_j p_j = 0$$

only on the edges of the element $K$, where some vectors (or one of the coordinates of some vectors) vanishes. There, collinearity between the elements of the set is enabled because the order prevention coordinate-wise is broken. Thus, $\mathcal{I}$ is free in the strict interior of $K$.

■

**Let us gather $\mathcal{N}$ and $\mathcal{I}$** Let us show that the sets $\mathcal{N}$ and $\mathcal{I}$ generates $RT_k(K)$; as expressed in the following property.

**Proposition 4.19**

$$RT_k(K) = \operatorname{span}\{\mathcal{N},\, \mathcal{I}\}$$

The proof is immediate through the following assertions.

**Property 4.20** Dimension

$\dim \operatorname{span}\{\mathcal{N}, \mathcal{I}\} = \dim RT_k(K)$.

***Proof.*** We have by previous observations that

$$\begin{aligned}\dim\{\mathcal{N}, \mathcal{I}\} &= 2d \dim \mathbb{Q}_k(f) + \dim \mathbb{P}_{k-1,\,k} \times \mathbb{P}_{k,\,k-1} \\ &= \dim RT_k(K).\end{aligned}$$

■

**Proposition 4.21**

If $p$ belongs to $\{\mathcal{N}, \mathcal{I}\}$, then $p$ belongs to $RT_k(K)$.

***Proof.*** We already saw in the *Proposition 4.9* that if $p$ belongs to $\mathcal{N}$, then $p$ belongs to $RT_k(K)$. The same applies for the subset $\mathcal{I}$ by the *Proposition 4.15*. Thus, for any element $p$ of $\{\mathcal{N}, \mathcal{I}\}$, $p$ belongs to $RT_k(K)$.

■

**Proposition 4.22**

The set $\{\mathcal{N}, \mathcal{I}\}$ is free.

***Proof.*** We give a similar proof to the one of the simplicial elements case that can be found in the *Section 3.24*.

• Let us show first that one cannot create a function of the set $\mathcal{I}$ with the help of functions in $\mathcal{N}$ and the other functions of $\mathcal{I}$. Indeed, we should guarantee that $p \cdot n|_{\partial K} \equiv 0$. And on each edge, the only way to have a constant value on the boundary is to add all the local Lagrangian functions to each other with the same weight. We denote for any $i \in [\![1, n]\!]$;

$$l_i = \sum_{m=1}^{k+1} l_{i,\,m}.$$

Note that by construction, we retrieve $l_i|_{\partial K} = 1$. Indeed, the functions $\{l_i\}_i$ depend only on one variable and are thus constant along the other one.

Therefore, the only combinations of functions of $\mathcal{N}$ that may be considered when adding them to the other contributions of functions in $\mathcal{I}$ read

$$q = Al_1e_1 + Bl_2e_2 + Cl_3e_3 + Dl_4e_4$$

for some real constants $A$, $B$, $C$ and $D$. Let us look at the possible range of admissible functions that would allow a linear combination of functions in $\mathcal{I}$ and $\mathcal{N}$ to lie in $\mathcal{I}$. We require

$$q \cdot n|_{\partial K} = (Al_1e_1 + Bl_2e_2 + Cl_3e_3 + Dl_4e_4) \cdot n|_{\partial K} = 0.$$

Considering the constraints on the functions $l_i$ and the definition of the vectors $\{e_i\}$, we can split the above equation on each edge and retrieve the following system.

$$\begin{cases} (Ae_1 + Be_2 + Ce_3 + De_4) \cdot n_1|_{f_1} & = 0 \\ (Ae_1 + Be_2 + Ce_3 + De_4) \cdot n_2|_{f_2} & = 0 \\ (Ae_1 + Be_2 + Ce_3 + De_4) \cdot n_3|_{f_3} & = 0 \\ (Ae_1 + Be_2 + Ce_3 + De_4) \cdot n_4|_{f_4} & = 0 \end{cases}$$

Plugging the definitions of the vectors $e_i$ and $n_i$ before inserting the values of $x$ and $y$ on the respective edges, we get

$$\begin{cases} C - D & = 0 \\ 2B + C + D & = 0 \\ A + B + 2D & = 0 \\ A - B & = 0. \end{cases}$$

We finally obtain an infinite number of solutions trough the relation

$$A = -B = C = -D.$$

However, since by construction the functions $l_i$ are constants everywhere, this means that the sum of the pondered functions $l_i$ will value zero everywhere, and will not help to construct an other function of $\mathcal{I}$ (*see the Figure 27*). This concludes the first case.

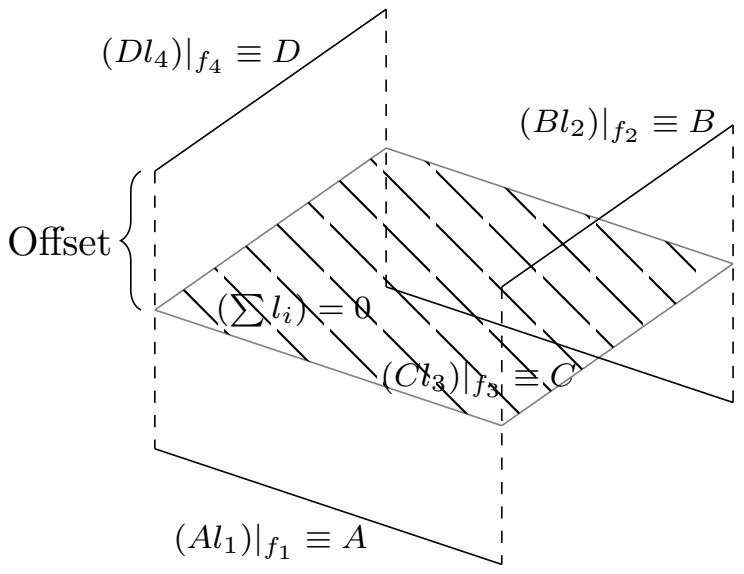

Fig. 27: Layout of the functions $l_i$ on the reference element.

• Secondly, we know that

$$\forall i \in [\![1,\, 4]\!],\, m \in [\![1,\, k+1]\!], \quad p_{i,\, m}(x,\, y) \cdot n|_{\partial K} \not\equiv 0$$

for the normal functions and

$$p_{i,\, m}(x,\, y) \cdot n|_{\partial K} \equiv 0$$

for the internal ones. Furthermore, we saw that the only function belonging to span $\mathcal{N}$ that can vanish on the edges with respect to the normal component is the identically null function. Therefore, no function created from $\mathcal{N}$ can fall in $\mathcal{I}$ and *vice - versa.*

• Lastly, let us show that it is also impossible to create a function of $\mathcal{N}$ with the help of other functions of $\mathcal{N}$ and functions in $\mathcal{I}$. We start by noticing that as $(x,\, y) \mapsto b_j(x,\, y) \in \mathbb{P}_{k-1,\, k}$ and $(x,\, y) \mapsto b_j(y,\, x) \in \mathbb{P}_{k,\, k-1}$, we have

$$(x,\, y) \mapsto \begin{pmatrix} x(x-1)b_j(x,\, y) \\ y(y-1)b_j(y,\, x) \end{pmatrix} \in \mathbb{P}_{k+1,\, k} \times \mathbb{P}_{k,\, k+1}.$$

Similarly, it comes $(x,\, y) \mapsto e_{n+2}(x,\, y) \circ (b_j(x,\, y),\, b_j(y,\, x))^T \in \mathbb{P}_{k+1,\, k} \times \mathbb{P}_{k,\, k+1}$. However in the same time we have

$$(x,\, y) \mapsto \begin{pmatrix} x - n_{ix} \\ y - n_{iy} \end{pmatrix} l_{i,\, m}(x,\, y) \in \mathbb{P}_{[k+1],\, [0]} \times \mathbb{P}_{[k],\, [1]} \text{ or } \mathbb{P}_{1,\, [k]} \times \mathbb{P}_{0,\, [k+1]},$$

as $l_{i,\, j} \in \mathbb{P}_{[k],\, 0}$ or $\mathbb{P}_{0,\, [k]}$. Therefore, when adding any combination of those two sets of functions we can only construct a function $q_m$ of the shape

$$q = \sum_{j=1}^{\dim\{\mathcal{I}\}} \alpha_j \begin{pmatrix} x(x-1)b_j(x,\, y) \\ y(y-1)b_j(y,\, x) \end{pmatrix} + \sum_{\substack{\{e_i l_{i,\, m}\} \in \tilde{\mathcal{N}} \subset \mathcal{N}, \\ \dim \tilde{\mathcal{N}} = \dim(\mathcal{N}) - 1}} \beta_{i,\, m} \begin{pmatrix} (x - n_{ix}) l_{i,\, m}(x,\, y) \\ (y - n_{iy}) l_{i,\, m}(x,\, y) \end{pmatrix}. \quad (57)$$

for some real constants $\{\alpha_j\}_j$ and $\{\beta_{i,\, m}\}_{i,\, m}$. As the case where the left sum vanishes has been treated previously, we assume here that the first term is not identically null. It follows that no function of that form can belong to $\mathcal{N}$. Indeed, in order to belong to $\mathcal{N}$ the function $q_m$ must be of the form

$$q = \begin{pmatrix} x - n_{ix} \\ y - n_{iy} \end{pmatrix} l_{i,\, m}(x,\, y) \quad (58)$$

for some given index $i \in [\![1,\, n]\!]$ and $m \in [\![1,\, k+1]\!]$.

However, no matter the used combination in (57) the function $b_j$ always contributes simultaneously to the two coordinates through the term

$$\begin{pmatrix} x^2 b_j - x b_j \\ y^2 b_j - y b_j \end{pmatrix} = \begin{pmatrix} x(x b_j) - x b_j \\ y(y b_j) - y b_j \end{pmatrix}.$$

It implies that $(x,\, y) \mapsto (x b_j(x,\, y),\, y b_j(x,\, y))^T$ must belong to the set $\{l_{i,j}\}_{i,j}^2$. Recalling that the functions $l_{i,j}$ are only dependent of one variable, we then infer that $q_j$ cannot be recast in the shape (58), and thus does not belong to $\mathcal{N}$.

• Finally, when merging all the above cases together, we get that the set $\{\mathcal{N},\, \mathcal{I}\}$ is free.

■

**Conclusion** The following functions form a basis of $RT_k(K)$ for quadrilateral elements in two dimensions.

*Normal basis functions*

$$\phi_j = l_{i,m} e_i, \qquad {}_{\substack{i \in [\![1,\, 4]\!],\, m \in [\![1,\, k+1]\!] \\ j \in [\![1,\, 4(k+1)]\!]}}$$

*Internal basis functions*

$$\phi_j = p_m, \qquad {}_{\substack{m \in [\![1,\, 2k(k+1)]\!] \\ j \in [\![4(k+1),\, 2(k+1)(k+2)]\!]}}$$

where the vectors $\{e_{n+i}\}_i$ are defined in (55) and the functions $\{p_j\}_j$ are defined in (56a).

## 4.4 Numerical implementation

As in the simplicial case, we implemented the proposed construction method for both global Lagrangian and local Lagrangian basis functions. The local Lagrangian basis function have then been tuned against the set of degrees of freedom to retrieve the classical Raviart – Thomas setting. We report here the results for the third order space.

When computing the locally Lagrangian basis functions, we do not fall directly in the setting of the Raviart – Thomas elements. Indeed, as one can see on the second edge, the normal component of the third basis function that is generated from the first edge is not vanishing as wished, but has a non - zero constant value.

However, the others specifications as the polynomial order on the boundary and the polynomial regularity of the two components within the space are preserved. The set of internal basis functions is also already matching the framework of the Raviart – Thomas elements.

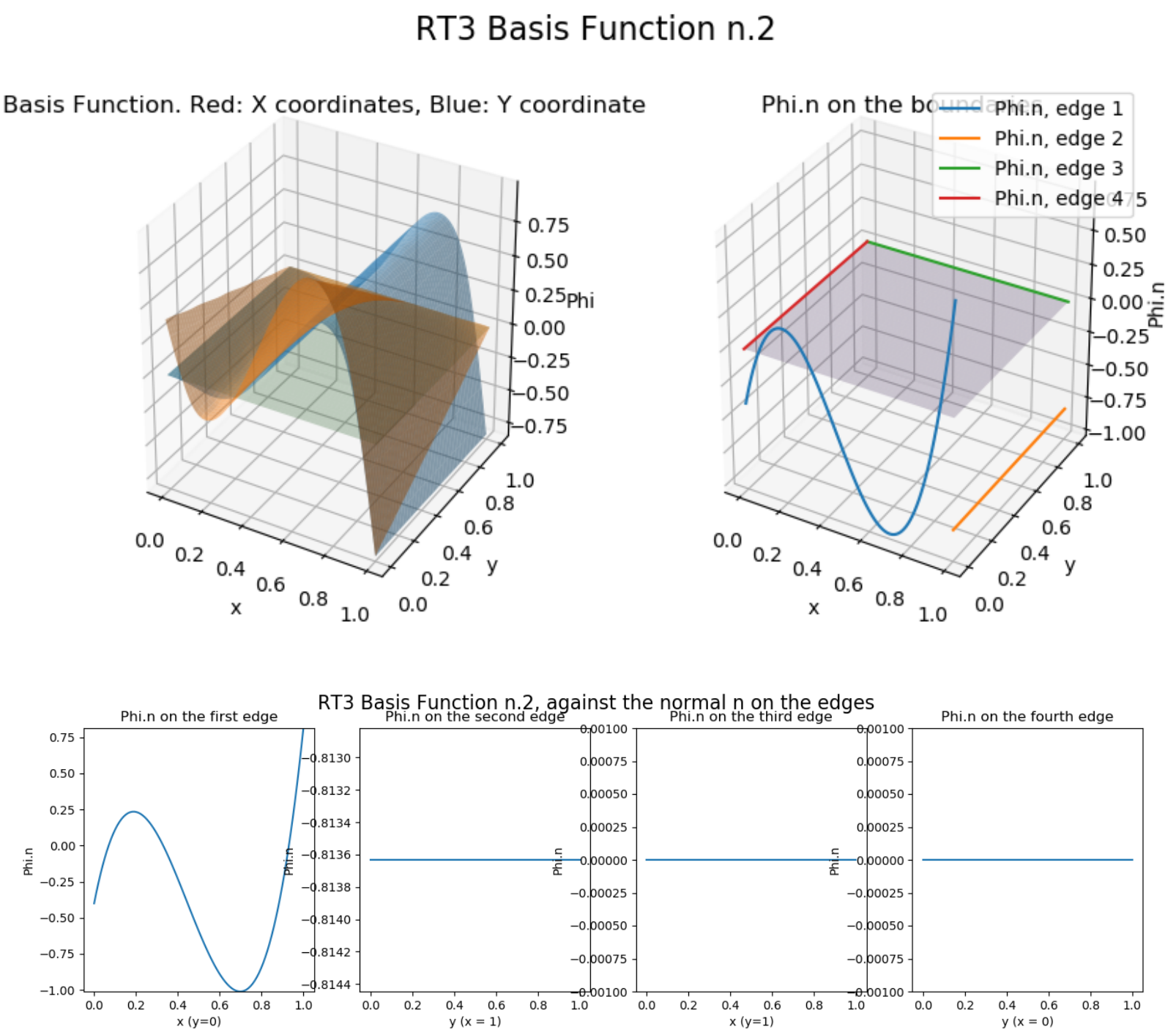

Fig. 28: Locally Lagrangian third normal basis function in the case $k = 3$.

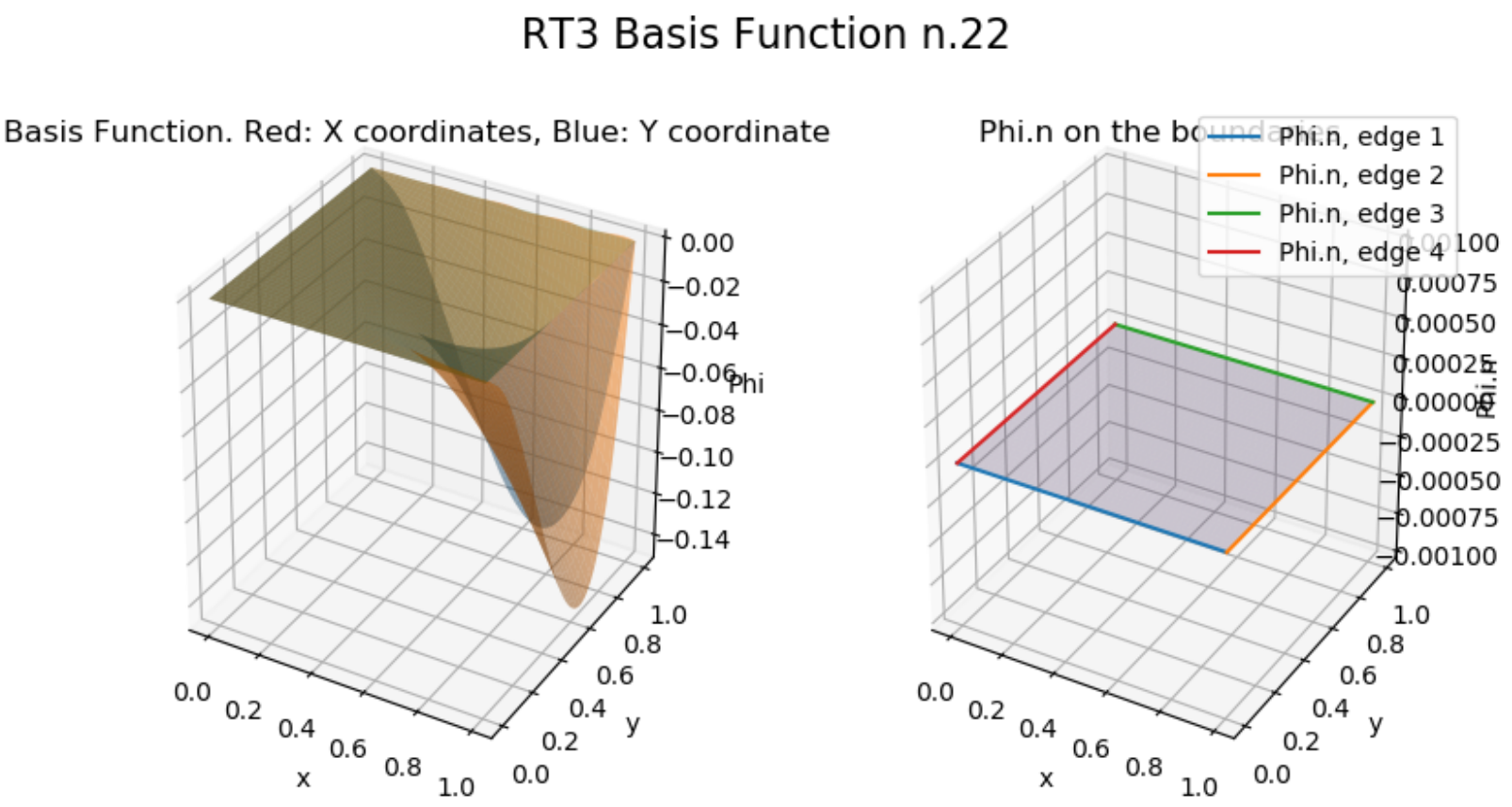

Fig. 29: Representative of a canonical internal function in the case $k = 3$

Changing the definition of the vectors $e_i$ allowed us to retrieve directly the expected behaviour. Indeed, for the basis functions enjoying a global Lagrangian properties, the properties of the Raviart – Thomas basis functions are directly verified both on the boundary and within the element.

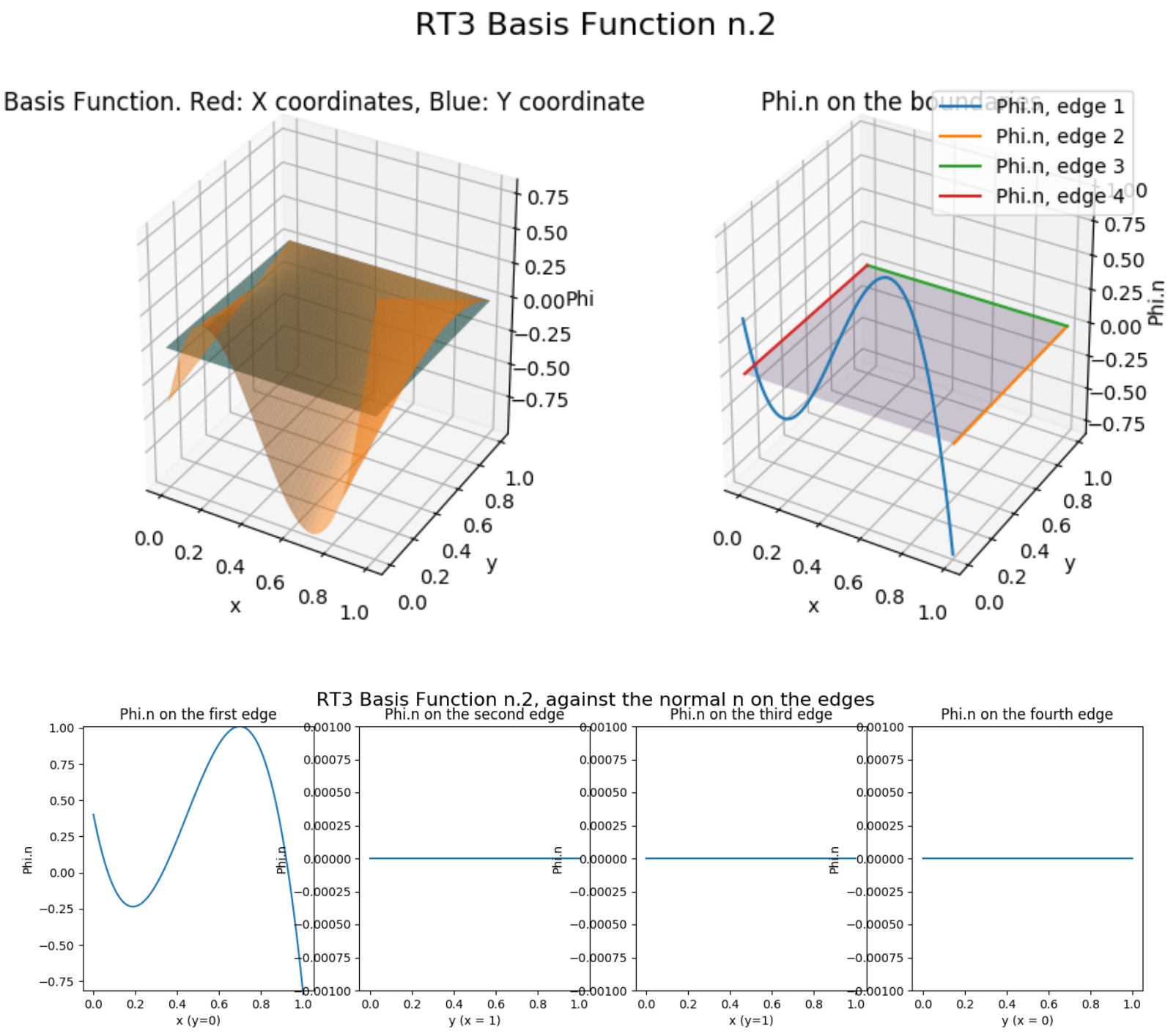

Fig. 30: Globally Lagrangian third normal basis function in the case $k = 3$.

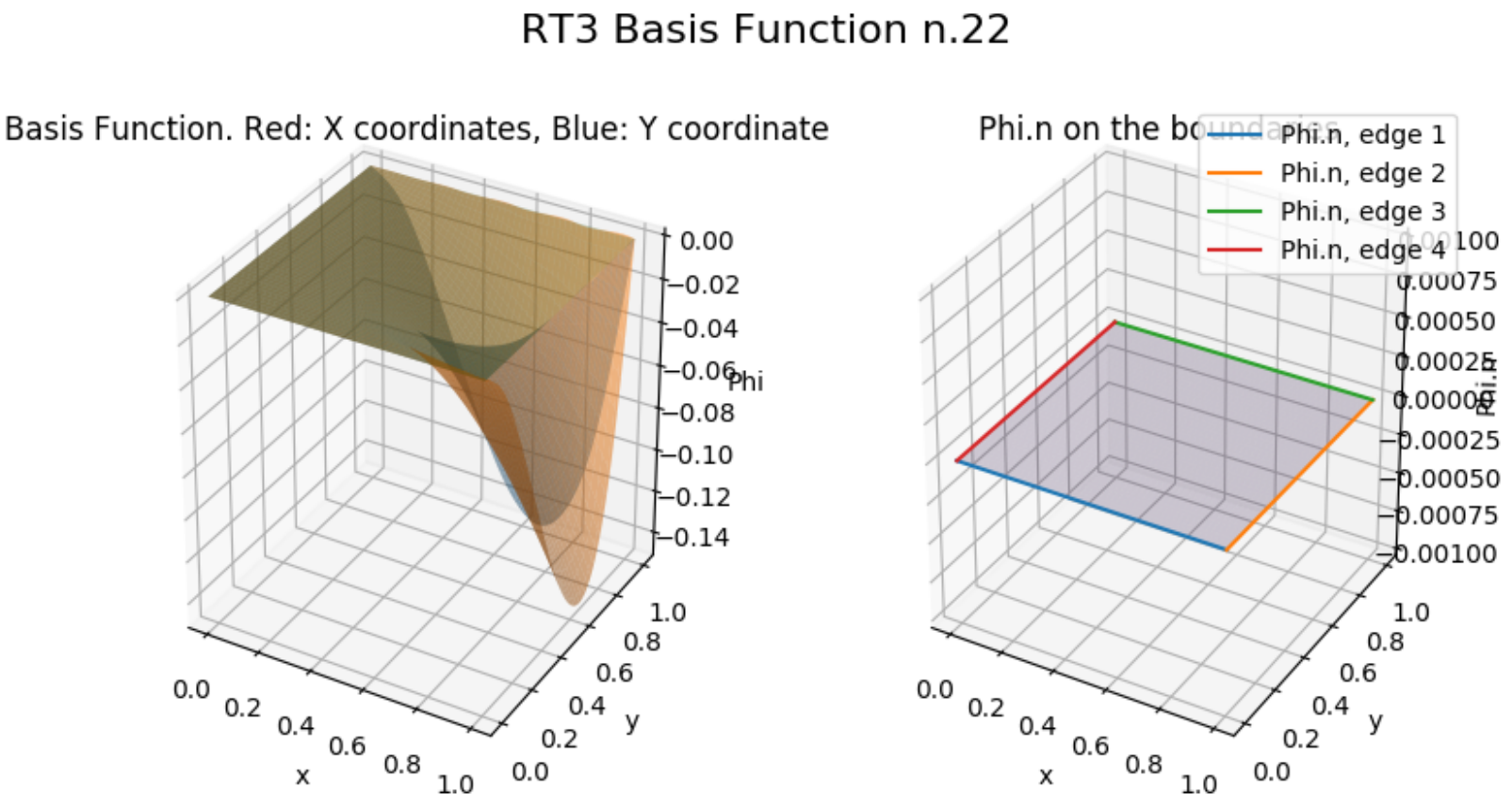

Fig. 31: Representative of one internal function for $k = 3$

We then tune the local Lagrangian basis functions set towards the classical set of degrees of freedom to retrieve the classical Raviart – Thomas basis functions. If as in the simplicial case the amplitude of the functions is huge, this time we can also observe dominated variations within the element for the normal basis functions.

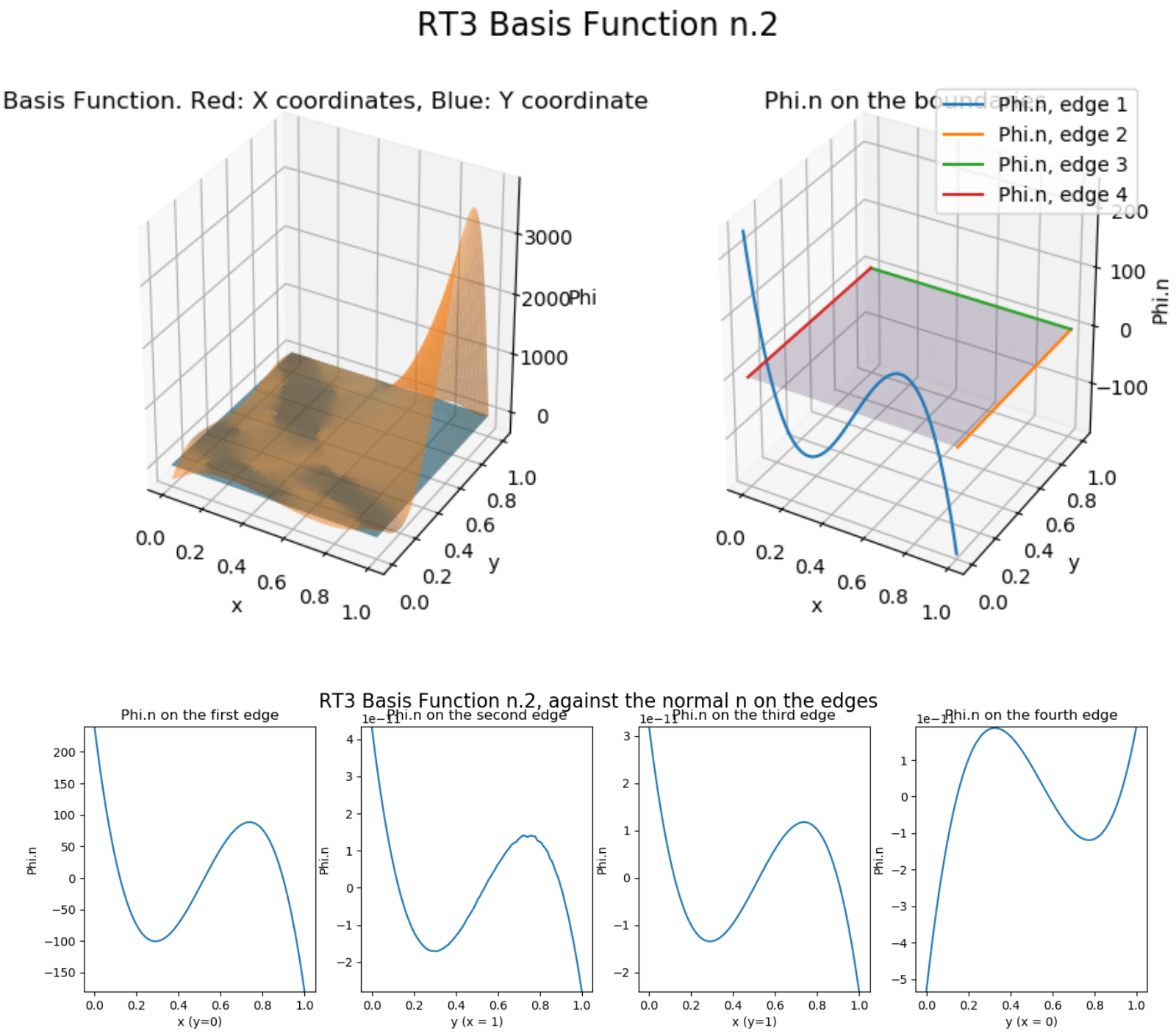

Fig. 32: Third normal basis function in the case $k = 3$.

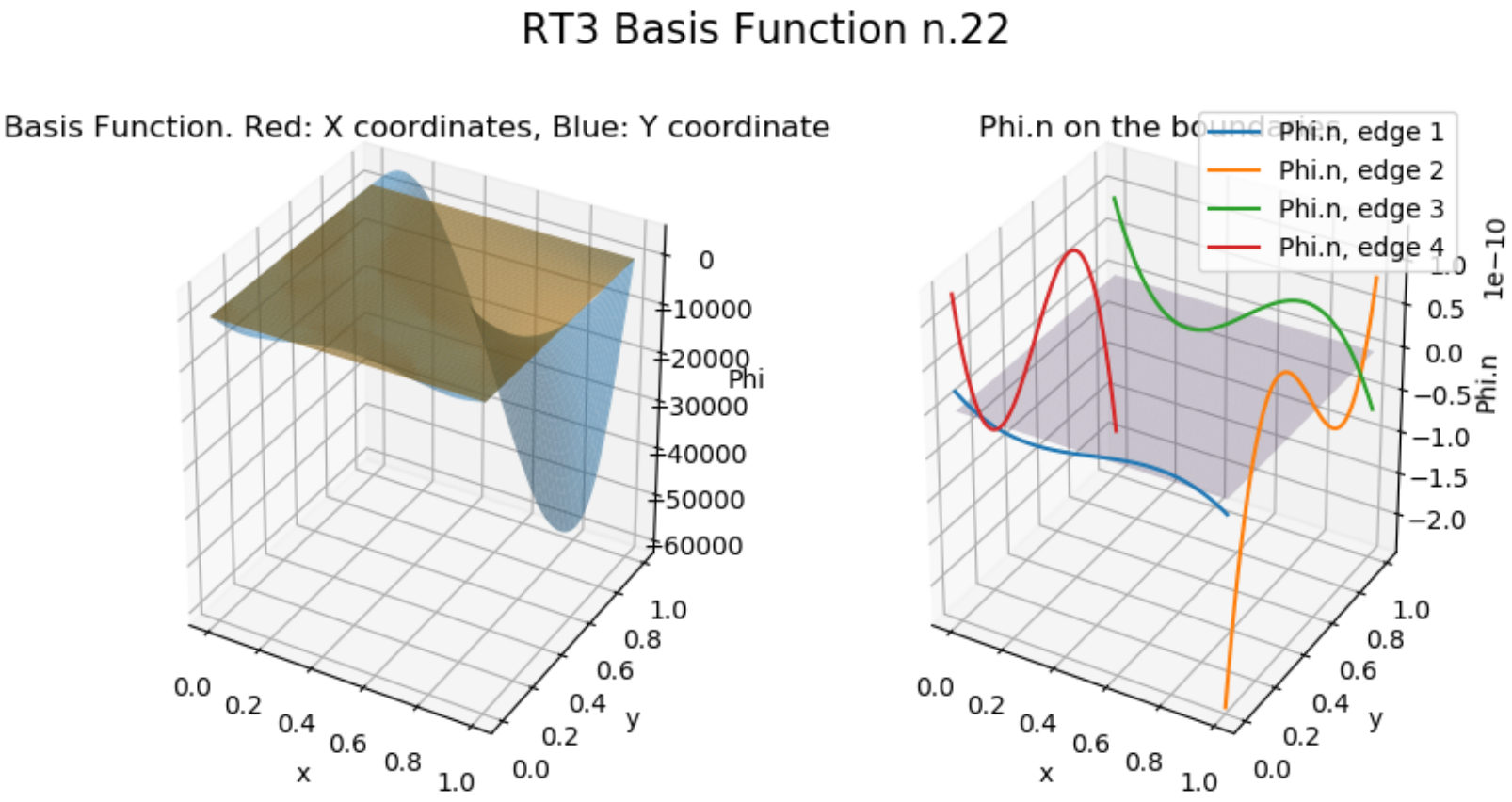

Fig. 33: Representative of one internal function for $k = 3$

# 5 A framework for arbitrary polytopes

We would like here to extend the previously described Raviart – Thomas elements by drawing a general framework that allows the use of any non - convex polytope. In other words, for a given polytope $K \subset \mathbb{R}^d$ that is reasonably squeezed, we would like to define a finite dimensional discretisation space $\mathbb{H}_k(K)$ of order $k$ enjoying the following property.

$\mathbb{H}_k(K)$ should be a subset of $H(\text{div}, K)$, that is; for any $p \in \mathbb{H}_k(K)$ one should have

$$\begin{cases} \operatorname{div} p \in L^2(K) \\ p \cdot n|_{\partial K} \in \mathcal{H}_k(\partial K) \end{cases}$$

for some space

$$\mathcal{H}_k(\partial K) = \{u|_{\partial K} \in L^2(\partial K),\, u|_f \in \mathcal{P}_{b,k}(f), \quad \forall f \in \partial K\} \tag{59}$$

where $P_{b,k}(f)$ represents any polynomial subspace of $\mathbb{P}(\mathbb{R}^{d-1})$ providing a discretisation of order $k$ on the face $f$. Note that as continuity across the faces is not required, the definition of the space $\mathcal{H}_k(\partial K)$ depends on the number of faces $n$ and its dimension is naturally given by $\dim \mathcal{H}_k(\partial K) = n\, \dim(P_{b,\,k}(f))$.

$H(\text{div}, K)$ – conformal elements can then be constructed in the spirit of Raviart – Thomas by enforcing the continuity of the normal components across the faces through degrees of freedom whose supports match the boundary of $K$. Therefore, we keep the idea of hardly separating normal and internal degrees of freedom and look for a set $\{\sigma\}$ such that the triplet $E_k(K) = (K, \mathbb{H}_k(K), \{\sigma\})$ fulfils the following three criterion.

- (*Layout of the degrees of freedom*) The set of degrees of freedom $\{\sigma\}$ can be split into two distinct subspaces; normal degrees of freedom that live only on the boundary of $K$ and (if any) internal degrees of freedom that live only in the interior of $K$.
- ($H(\text{div}, K)$ *– conformity*) The triplet $E_k(K) = (K, \mathbb{H}_k(K), \{\sigma\})$ forms a $H(\text{div}, K)$ – conformal element. In particular, the set of degrees of freedom $\{\sigma\}$ is unisolvent for $\mathbb{H}_k(K)$.
- (*Layout of the corresponding basis functions*) Let $\{p\}$ be the basis of $\mathbb{H}_k(K)$ that is dual to the set $\{\sigma\}$ in $\mathbb{H}_k(K)$. Then, for every $p \in \{p\}$ it holds:

$$\begin{aligned} &p \cdot n|_{\partial K} \in \mathcal{H}_k(\partial K) && \text{if } p \text{ is a normal function,} \\ &p \cdot n|_{\partial K} \equiv 0 && \text{if } p \text{ is an internal function.} \end{aligned}$$

***Note.*** We consider here only non-degenerated polytopes in dimension $d$, that is all polytopes $K$ whose faces $f$ share the same dimension $d-1$ (*see the Figure 34*). ▲

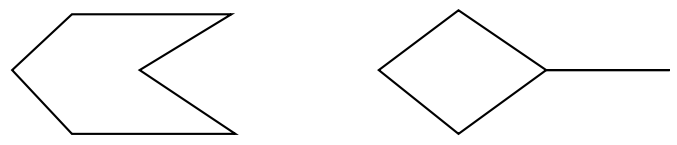

Fig. 34: Left: Considered non - degenerated polytope. Right: Excluded degenerated polytope.

We start by explaining why an analogous construction to the simplicial and quadrilateral elements ones fails before building admissible elements.

## 5.1 Limitations of the classical formulation of Raviart – Thomas elements

The definitions of simplicial and quadrilateral Raviart – Thomas elements cannot be directly extended to arbitrary polytopes as the dimension of the space $RT_k(K)$ is only dependent on the desired order $k$ and not adaptive towards the number of faces of the element.

In order to emphasize the arising issues, let us fork the definition of the classical quadrilateral $RT_k(K)$ elements to general polytopes. There, the space $RT_k(K)$ is still taken as

$$RT_k(K) = (\mathbb{Q}_k(K))^d + x\,\mathbb{Q}_{[k]}(K)$$

while the set of degrees of freedom $\{\sigma\}$ reads

$$\begin{cases} \sigma(q) = \int\limits_{\partial K} q \cdot n\, p_k \,\mathrm{d}\gamma(x), & \forall p_k \in \mathcal{T}_k(\partial K) \\ \sigma(q) = \int\limits_{K} q \cdot p_k \,\mathrm{d}x, & \forall p_k \in \mathcal{P}_k(K), \end{cases}$$

where if not empty $\mathcal{P}_k(K)$ represents the polynomial space against which the internal component of functions living in $RT_k(K)$ is projected (*see Section 2.3*). In the Raviart – Thomas setting, the $H(\mathrm{div}, K)$ – conformity is only enforced by the definition of the normal degrees of freedom whose supports are restricted to a single face. Therefore, one is required to have at least one normal degree of freedom per face and the dimension of $\mathcal{T}_k(\partial K)$ is naturally given by $\dim \mathcal{T}_k(\partial K) = n \dim(\mathbb{Q}_k(f))$ as the continuity across the face is not required.

*A - contrario,* the space $\mathcal{P}_k(K)$ is not directly dependent of the number of faces of the element as it is used to characterise the function on the strict inner part of the cell. However, a dimensional problematic emerges once observing the structural constraints of the space $RT_k(K)$. Indeed, as $\dim RT_k(K)$ does not vary on $n$ we get $\dim(\mathcal{P}_k) = \dim RT_k(K) - n\dim(\mathbb{Q}_k(f))$, which is very sensitive towards the number of faces and makes very difficult its general definition.

Therefore, even if this definition fits *a - priori* in the framework presented in the *Section 5* when taking $\mathcal{P}_{b,k}(f) = \mathbb{Q}_k(f)$, the corresponding triplet $E_k(K)$ does not fulfil the three criterion given in the *Paragraph 5.* In particular, even if the changes from the original definition of the space $RT_k(K)$ for quadrilaterals only take place in the dimension of $\mathcal{T}_k(\partial K)$ and in the definition of the space $\mathcal{P}_k(K)$, the above simple observation leads to several issues preventing the very definition of $H(\text{div}, K)$ – conformal elements in this way. We detail these issues below in order to help the construction of an other type of $H(\text{div}, K)$ – conformal elements avoiding those pitfalls.

To make easier the explanation, we consider the two dimensional case built on the quadrilateral Raviart – Thomas space, but the observations also hold when one considers higher dimensions or uses the definition of the simplicial elements.

### 5.1.1 Discretisation quality

The first consequence of the dimensional variability of $\mathcal{P}_k(K)$ impacts the quality of the discretisation in the inner cell. Indeed, it holds $\dim \mathcal{T}_k(\partial K) = n \dim \mathbb{Q}_k(f)$ for any generic face $f$ of $\partial K$. Therefore, as for any fixed order $k$ the dimension of $\mathbb{Q}_k(f)$ is positive, more there are faces higher the total number of normal degrees of freedom is. Thus, the amount of internal degrees of freedom logically decreases as $n$ increases.

Consequently, as those degrees of freedom take care of the discretisation quality in the interior of the cell, the quality of the representation within the cell drops as $n$ increases, and for a fixed order $k$ we have less an less control of the inner cell. Furthermore, we cannot rely on the order of the discrete space to infer the discretisation quality without knowing the type of element that is in use.

***Example.*** *Two dimensional case* In two dimensions, for any polygon with $n$ edges and for any fixed $k$, adding a supplementary edge to form a polygon with $(n+1)$ edges makes the number of internal moments drop by $(k+1) = \dim \mathbb{Q}_k(f)$. We picture this for $n = 5$ and $k = 3$ in the *Figure 35.* ⧫

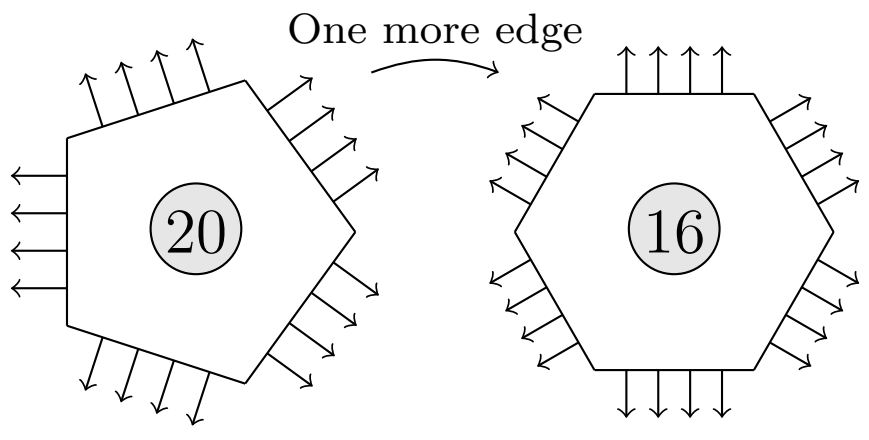

Fig. 35: Inner discretisation drops when $n$ increases at fixed order $k$.

More dramatically, when one works with a non - homogeneous mesh (*i.e.* when gathering various element's shapes) the quality of the discretisation is not homogeneous thorough the mesh. Indeed, for a same order $k$ but two different shapes, the number of internal degrees of freedom will not be the same and the discretisation will not be performed in the same way. In addition, we may have to consider $k$ large enough to guarantee that both elements can be properly defined (*see Section 5.1.2*). This impedes the use of low order spaces which would have been possible when considering only the elements with the smallest number of faces.

***Example.*** *Two dimensional case* We give two examples of an incompatible discretisation method for non - homogeneous meshes.

• Let be a mesh generated from elements having five and six edges, and consider $k = 3$. We represent the corresponding degrees of freedom in the *Figure 36*.

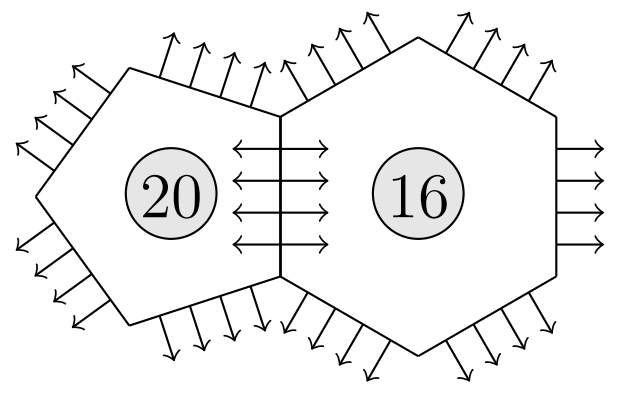

Fig. 36: Incompatible discretisation order within the mesh.

There, the number of internal degrees of freedom is 16 for the element with six edges and 20 for the element with only five edges. Although the $H(\text{div})$ – conformity is preserved by the definition of the normal degrees of freedom, the solution will not be represented with the same accuracy within each cell.

• For a mesh generated from elements having five or twelve edges, the situation is even more critical. We represent the layout of the corresponding degrees of freedom in the *Figure 37*.

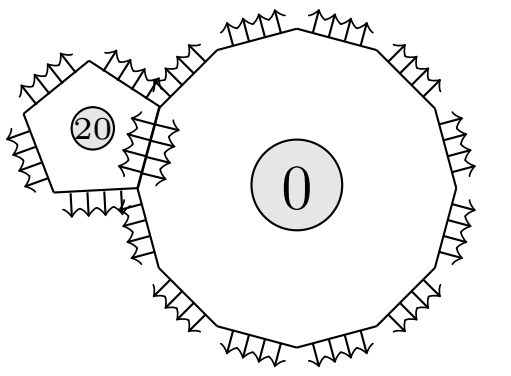

Fig. 37: Incompatible discretisation order within the mesh.

In this case, there is a dramatic disequilibrium for two neighbouring cells between a fine inner representation of the discretised quantity (35 internal degrees of freedom) and absolutely no representation of the inner part. Even if this case is extreme and happens only when two polygons with highly different shapes are next to each other, this is not acceptable in terms of the consistency of the discretisation quality through the domain. ♦

Therefore, as we aim to see the dimension of $\mathcal{P}_k(K)$ insensitive to $n$, we have to design the space $\mathbb{H}_k(K)$ so that for any fixed order $k$ its dimension increases by $\dim P_{b,k}(f)$ each supplementary face. That way, the space dimension increment will allow to distribute equally the normal degrees of freedom one every face while keeping unchanged the number of internal degrees of freedom. It would also be admissible to ask for a dimensional increment by a multiple of $n+1$ as the equal distribution of the normal degrees of freedom would also be possible (up to allow the definition of $P_{b,k}(f)$ to be dependent on $n$).

### 5.1.2 Range of admissible discretisation order

A simple fork of the definition of quadrilateral $RT_k(K)$ elements on polytopes also severely impacts the range of the discretisation orders one can enjoy. Indeed, following the logic of the *Section 5.1.1* one can see that for every order $k \in \mathbb{N}$ there exists $n$ - gons such that the dimension of the space $RT_k(K)$ will be lower than the dimension of $\mathcal{T}_k(\partial K)$. Therefore, when following the definition (5.1) the normal degrees of freedom cannot be set on every face.

For even more complex polytopes with very high number of faces, it can happen that the dimension of $RT_k(K)$ is lower than $n$. In those cases, the question of the very existence of $E_k(K)$ for early orders $k$ is assessed as it would be impossible to define the set of normal degrees of freedom either in their classical sense or through some space $\mathcal{H}_k(\partial K)$ built from any polynomial space $P_{b,k}(f)$.

***Example.*** *Two dimensional case* Let us consider the case $n = 5$ and $k = 0$. We have that $\dim RT_k(K) = 4$ and $\dim \mathcal{T}_k(\partial K) = 5 \dim \mathbb{Q}_k(f) = 5\,(0+1) = 5$. However, we cannot distribute equally four normal degrees of freedom over five edges (*see the Figure 38*).

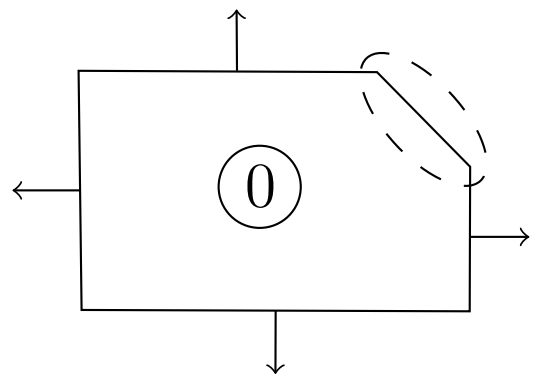

Fig. 38: The classical definition of the degrees of freedom is impossible to set for the lowest order space. ♦

This issue is only a dimensional problem. Therefore, changing the definition of the degrees of freedom while preserving the wish of $H(\mathrm{div}, K)$ – conformity will not help. However, we do have the following proposition.

**Proposition 5.1**

The application

$$\begin{aligned} \mathbb{N} &\longrightarrow \mathbb{N} \\ k &\longmapsto \dim RT_k(K) - \dim \mathcal{T}_k(\partial K) \end{aligned}$$

is increasing and onto on $\mathbb{N}$ for all $n \in \mathbb{N}^*$.

***Proof.*** By definition of $RT_k(K)$ and $\mathcal{T}_k(\partial K)$, we have

$$\begin{aligned} \dim RT_k(K) - \dim \mathcal{T}_k(\partial K) &= d(k+2)(k+1)^{d-1} - n(k+1)^{d-1} \\ &= (k+1)^{d-1}(d(k+2) - n) \end{aligned}$$

As $d > 0$, $(k+1)^{d-1} > 0$ for any $k$ and since for any fixed $n \in \mathbb{N}$ the term $(d(k+2) - n)$ is increasing along $k$, the claim immediately follows.

■

Thus, for any polytope of any number of faces $n$ there will be one $k_0$ from which the classical element $E_k(K)$ can be defined. The problem we are discussing here therefore only arises for early $k < k_0$ where the extent of "early" only depends on the number of faces.

However, even for $k \geq k_0$ one needs at least some rectification in the definition of the degrees of freedom when considering the fork of classical quadrilateral $RT_k(K)$ elements to polytopes. We detail three failing approaches that show the need of the properties detailed in the *Section* (5.2.1), each of them being incompatible with the classical definition of $RT_k(K)$.

**Use of arithmetic congruences** A natural approach to alleviate the problem of lack of low order spaces is to use arithmetic congruences. Indeed, one can always distribute some degrees of freedom on the boundary as early as possible up to switching $\mathbb{Q}_k(f)$ to some space $\mathcal{P}_{b,\,k}(f)$ in the definition of $\mathcal{T}_k(\partial K)$ for early orders $k$. For the sake of consistency with the classical definition of $RT_k(K)$ elements, the definition is switched back to the classical one as soon as it is possible. We then obtain a definition that is split into three cases.

First, the element $E_k(K)$ is left undefined when the dimension of the space $E_k(K)$ is less than the number of faces. Then, for the orders $k$ such that $n \leq \dim RT_k(K) < \dim \mathcal{T}_k(\partial K)$ the normal degrees of freedom are distributed modulo $n$ and the remaining is set as internal degrees of freedom. Lastly, for any order $k$ such that $\dim \mathcal{T}_k(\partial K) \leq \dim RT_k(K)$ the classical definition of $E_k(K)$ is used. Note that since once one could use the classical definition at some order $k_0$ it can be used for any subsequent $k \geq k_0$, those three cases happens only in the incremental order along $k$.

More precisely, we set $k_1$ the smallest non - negative integer such that for all $k \geq k_1$, $\dim RT_k(K) \geq n$ and define $k_0$ as the smallest non - negative integer such that for all $k \geq k_0$, $\dim RT_k(K) \geq \dim \mathcal{T}_k(\partial K)$. The definition then reads:

- For $0 \leq l < k_1$ the element $E_k(K)$ is undefined,
- for $k_1 \leq k < k_0$ the internal and normal degrees of freedom are set such that $\dim RT_k(K) \equiv \dim\{\mathcal{I}\}[n]$,
- for $k \geq k_0$ the classical quadrilateral definition of $E_k(K)$ is used.

Though this definitions allows us to enjoy $H(\text{div},\, K)$ – conformal elements for relatively early orders, one could notice that the number of normal degrees of freedom per face is not explicitly prescribed for $k_1 \leq k < k_0$. Even if by definition the number of degrees of freedom per face is non - decreasing, it highly depends on $n$ through a congruency relationship which *a - fortiori* is not linear in $k$. Therefore, when incrementing the discretisation order from $k$ to $k+1$ the refinement on each face is not necessarily performed in the same way, if any (see the example below).

Another consequence of this congruency relationship impedes the discretisation quality of the inner cell. As the number of normal degrees of freedom is not monotonous against $k$, the discretisation quality can decrease as $k$ increases before increasing again. Furthermore, this quality drop can arise several times.

Lastly, when switching to the classical definition at $k_1$, the inner discretisation of the cell through internal degrees of freedom may become coarser while the one given through normal degrees of freedom may suffer from some unexpected refinement jumps. Although rare, cases where the reverse phenomena happens may also be encountered. Thus, there is still no intelligible meaning for the early orders $k$ in terms of the discretisation quality.

***Example.*** *Two dimensional case* We give two examples of this construction.

• Let us first consider an element with $n = 12$ edges. In that case, the refinement logic is already broken for the elements built on the arithmetic definition.

For $k = 0$ we have as before $\dim RT_k(K) = 4$. And since $4 < 12 = \dim RT_k(K)$, the element $E_0$ remains undefined.

For $k = 1$, we have $\dim RT_k(K) = 12 = n$. However, $\dim RT_k(K) = 12 < 24 = n \dim \mathbb{Q}_k(f) = \dim \mathcal{T}_k(\partial K)$. Therefore, we use the arithmetic definition $12 \equiv 0[12]$, and set one normal degree of freedom per edge, without internal moment. We also set $k_1 = 1$ (*see the Figure 39*).

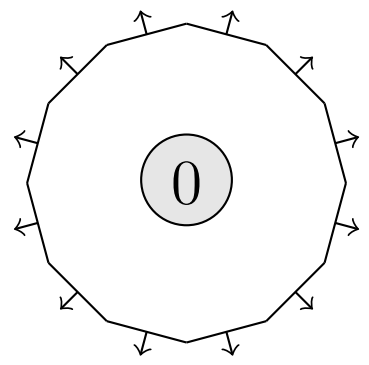

Fig. 39: Distribution of normal and internal degrees of freedom for $k = 1$ and $n = 12$.

For $k = 2$, we have $\dim RT_k(K) = 24 > 12 = n$. We still are in the arithmetic case as $\dim \mathcal{T}_2(\partial K) = 12 \times 3 = 36 > 24$. We then set $24 \equiv 0[12]$ and have two normal moments per edge, still without internal moment (*see the Figure 40*). We notice that we are already breaking the refinement logic. Indeed, if our normal moments are doubled and the discretisation quality on the edge is equally refined, no control of the interior of the cell appears. Thus, the discretisation is not homogeneously refined.

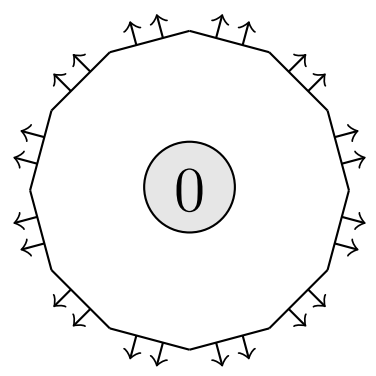

Fig. 40: Distribution of normal and internal degrees of freedom for $k = 2$ and $n = 12$.

For $k = 3$, $\dim RT_k(K) = 40$ and $\dim \mathcal{T}_3(\partial K) = 48$. We are still in the case of arithmetic congruences and define the number of normal moments following the relation $40 \equiv 4[12]$. We then have three normal moments per edge and four internal moments (*see the Figure 41*). The normal moments are still logically refined. However, the internal moments that pop up do not follow some refinement rate in their dimension with respect to $k$ only.

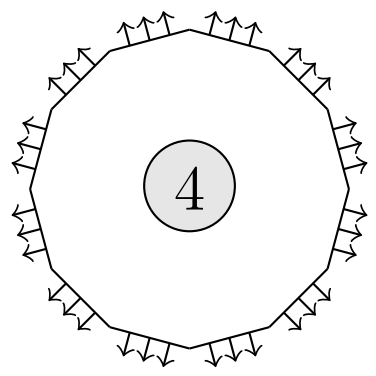

Fig. 41: Distribution of normal and internal degrees of freedom for $k = 2$ and $n = 12$.

For $k = 4$, we have $\dim RT_k(K) = 60$ and there, $\dim \mathcal{T}_4(\partial K) = 60$. Thus, we can switch to the classical definition of $E_4(K)$. We then set $5 = 4 + 1$ normal moments per edge and set the remaining, 0, as internal moments within the cell. We set $k_0 = 4$ (*see the Figure 42*). Here, the number of internal moments drops. Furthermore, following this algorithm, no element with four normal moments is designed.

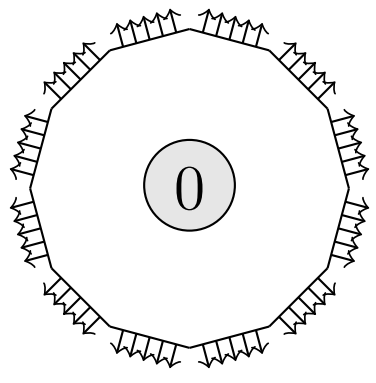

Fig. 42: Distribution of normal and internal degrees of freedom for $k = 2$ and $n = 12$.

• Let us now consider a polygon with $n = 10$ edges. There, the refinement process is only broken by the classical definition switch.

For $k = 0$, we have $\dim RT_k(K) = 4$. Since $4 < 10 = n$ there is no definition for $E_0(K)$.

For $k = 1$, we have $\dim RT_k(K) = 12$ and $n = 10 < 12$, but $12 < \dim \mathcal{T}_k(\partial K) = 20$. Therefore, no classical definition of the element is possible. We then define the normal degrees of freedom based on congruences through the relation $12 \equiv 2[10]$. There are two internal moments and ten normal moments. We set $k_1 = 1$ *(see Figure 43)*.

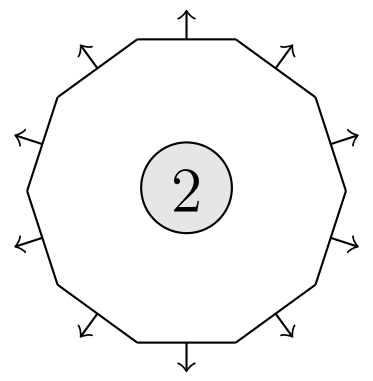

Fig. 43: Distribution of normal and internal degrees of freedom for $k = 1$ and $n = 10$.

For $k = 2$, we have $\dim RT_k(K) = 24$ and $\dim \mathcal{T}_k(\partial K) = 10(2+1) = 30 > 24$. Therefore, we are not yet allowed to define a classical element $E_2(K)$. Using again congruences, we have $24 \equiv 4[10]$ and get four internal degrees of freedom and two normal degrees of freedom per edge (*see the Figure 44*).

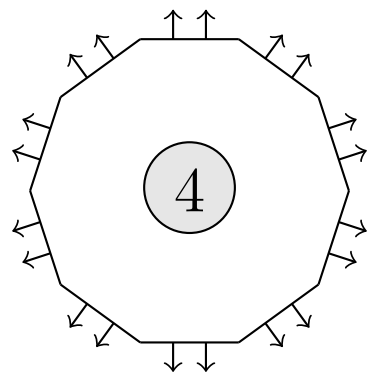

Fig. 44: Distribution of normal and internal degrees of freedom for $k = 2$ and $n = 10$.

So far, we are in a working example. The internal quality of discretisation and the quantity of information available on the edges are refined homogeneously. However, for $k = 3$ we have $\dim RT_k(K) = 40$. There, $\dim \mathcal{T}_k(\partial K) = 10(3 + 1) = 40 = \dim RT_k(K)$. We can then set $E_3(K)$ in its classical definition and get four normal degrees of freedom per edge without internal moment. We set $k_0 = 3$ (*see the Figure 45*).

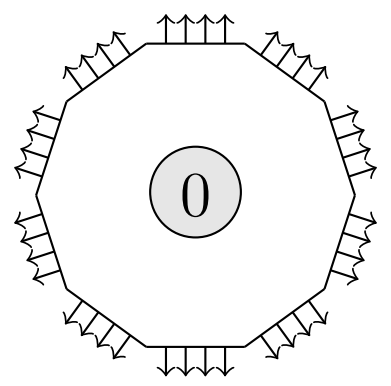

Fig. 45: Distribution of normal and internal degrees of freedom for $k = 3$ and $n = 10$.

In this example, the logic breaks in the switch to the classical definition. Indeed, there is no order $k$ where one can have exactly three moments per edge, creating a jump in the edges' discretisation refinement. Furthermore, the internal moments vanishes and the interior of the cell is not represented anymore. ♦

Here we could manage to have conformity as early as possible in $k$, but we pay it by an uncontrollable discretisation quality across the order $k$ for any polytope with a fixed number of faces $n$. We even get a non - monotonous, non-linear progression of the discretisation quality on the faces and severe fluctuations of the discretisation quality the within the cell for early $k < k_0$.

Therefore, when designing the elements $E_k(K) = (K, \mathbb{H}_k(K), \{\sigma\})$ we aim to set a definition of $\mathcal{P}_{b,k}$ such that the refinement on the faces is strictly monotonous (there is neither coarsening nor unexpectedly fine refinement of the discretisation on the faces when the order $k$ increases). Furthermore, the dimension of $\mathbb{H}_k(K)$ should increase only proportionally with $\dim \mathcal{P}_{b,k}$ when a face is added, not along with $n + 1$. That way it allows to keep the same definition of $\mathcal{P}_{b,k}$ for every $k$ while preventing the discretisation quality of the inner cell to drop.

**Use of internal moments** To avoid those jumps in the inner discretisation quality when increasing the order $k$, one could give up on the earliness of the definition of $H(\text{div}, K)$ – conformal element and define the moments exclusively as internal objects until $k \geq k_0$ where a classical definition of $E_k(K)$ is allowed. In other words, every moment would be set as internal for any $k$ such that $\dim RT_k(K) - \dim \mathcal{T}_k(\partial K) < 0$ (see the example below).

However, one has to be very careful here. Indeed, if this method provides a definition of $E_k(K)$ for early orders $k$, those spaces are not $H(\text{div})$ – conformal. Furthermore, as in the previous method when we switch to the classical form (*i.e.* $k = k_0$) the number of internal functions can decrease and the discretisation quality drop within the cell. Indeed, $k = k_0$ is the first time where normal degrees of freedom will be set up, and there will be the first time where a split between normal and internal functions will occur, impacting the largeness of the set of internal degrees of freedom. Note furthermore that this jump is more severe as $n$ grows. Therefore, even considering only one specific polygonal shape, the term "order" would again be meaningless.

***Example.*** *Two dimensional case* Let us consider the case $n = 5$.

- For $k = 0$, $\dim RT_k(K) = 4$. Since $4 < 5 = \dim \mathcal{T}_0(\partial K)$, we cannot use the classical definition of $E_0(K)$. Although we do not want to arbitrarily omit some information on the edges, we define the moments as internal ones despite the non $H(\mathrm{div})$ – conformity it implies (*see the Figure 46*).

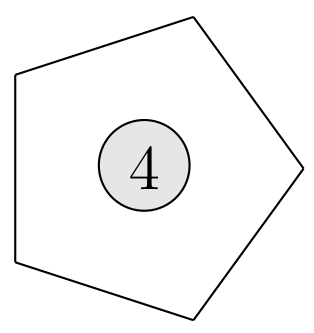

Fig. 46: Distribution of degrees of freedom for $k = 0$.

- For $k = 1$, we have $\dim RT_k(K) = 12$. There, $\dim \mathcal{T}_k(\partial K) = 5 \times 2 = 10 < 12$ and we have enough liberty to define normal degrees of freedom and $E_1(K)$ in its classical sense (*see the Figure 47*). We set $k_0 = 1$.

However, the number of internal degrees of freedom drops by two and the discretisation quality within the cell decreases. In particular, there exists cases where there is no internal moment, as for $n = 6$.

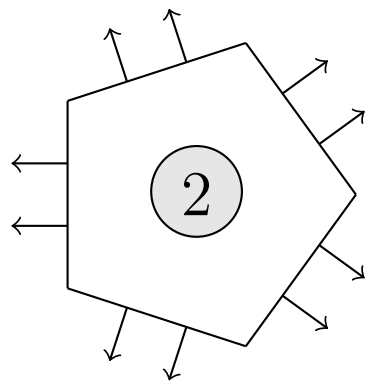

Fig. 47: Degrees of freedom for $k = 0$.

One can also observe that there is no case where there is only one moment per edge, which can be bothersome if one would like to prescribe only constants on the edges for example. Also, the term "order" is meaningless since it does not describe homogeneously the discretisation quality for both internal and normal degrees of freedom (at least for $k \leq k_0$). ⧫

***Remark.*** One could argue that instead of switching to the classical definition at the first $k_0$ where it is possible, one could wait an order $k$ high enough so that switching to the classical definition of $E_k(K)$ does not make the number of internal moment decrease. In this spirit, one could indeed start to define normal moments when the number of internal moments in the classically defined element $E_k(K)$ is bigger than the total dimension of $RT_{k-1}(K)$. That way, it preserves the discretisation quality within the cell, and by the increasing property of the quantity $\dim RT_k(K) - \dim \mathcal{T}_k(\partial K)$ (*see the Proposition 5.1*), the quality cannot drop later on either.

However, this is not always possible. In particular, in two dimensions for six edges, that would mean that the normal moments would appear from the smallest $k$ such that

$$\underbrace{2(k+1)(k+2)}_{\dim RT_k(K)} - \underbrace{6(k+1)}_{\substack{\text{number of}\\ \text{normal moments}}} \geq \underbrace{2k(k+1).}_{\substack{\dim(RT_{k-1}(K))\\ \text{Only internal moments for } k-1}}$$

This inequality reduces to $2(k+1)(k-2k-1) \geq 0$. The order $k$ being non - negative it is never fulfilled and waiting for any high order $k$ does not help. The same goes for any $n \geq 6$ in dimension two. ▲

In this case, the number of normal degrees of freedom are linearly expanding with $k$ from a certain order $k_0$. However, the early elements are not $H(\text{div}, K)$ – conformal and there is still one jump backwards in the number of internal moments when switching to the classical definition.

With this in mind, we aim to construct a space $\mathbb{H}_k$ allowing to set up a space $\mathcal{P}_k$ whose dimension is increasing with $k$. There, no backward jumps in the inner discretisation quality will occur.

**Shift the first admissible order** In the given setting, the only option that allows a homogeneous discretisation with respect to $k$ is to define $E_k(K)$ for any order $k \geq k_0$, where $k_0$ is the smallest integer so that $\dim RT_k(K) - \dim \mathcal{T}_k(\partial K) \geq 0$. That way, we simply omit the critical cases. This is possible as we have by the *Property 5.1* that if $E_{k_0}(K)$ can be defined, so can $E_k(K)$ for any $k \geq k_0$. However, even if such a $k_0$ can always be found please note that $k_0$ can be very large, especially when $n$ is large.

***Example.*** *Two dimensional example* Let us look at the first definable elements for $n = 6$ and $n = 12$.

| | $n = 6$ | $n = 12$ |
|---|---|---|
| $k = 0$ | 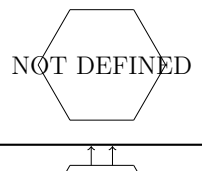  | 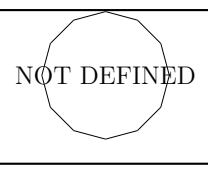  |
| $k = 1$ | 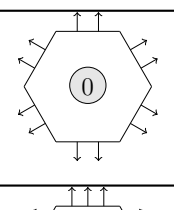  | 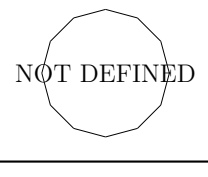  |
| $k = 2$ | 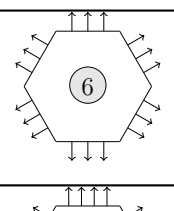  | 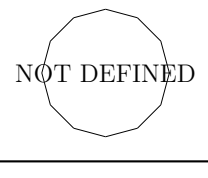  |
| $k = 3$ | 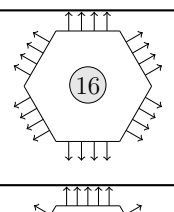  | 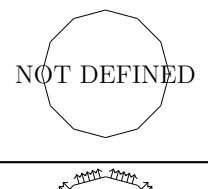  |
| $k = 4$ | 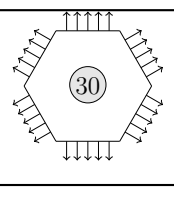  | 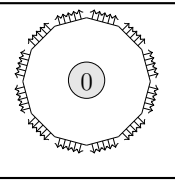  |

We observe that we still have the jump in the number of normal moments when the first element is defined. Therefore, assuming that the elements can be defined from $k_0$ on, there is no element that is enjoying $k+1$ moments on each edge for any $0 \leq k < k_0$. Furthermore, the relation between $n$ and $k_0$ is not linear. ♦

Even if it is very tempting, one should not shift hardly the element order by setting $\tilde{k} \leftarrow k + k_0$ in order to bypass artificially the elements that are not defined. Indeed, we are looking for $H(\text{div}, K)$ – conformal elements and shifting the order index would make impossible the definition of conformal meshes that are homogeneous in the discretisation order but built on elements that have various number of faces (*see the Figure 48*).

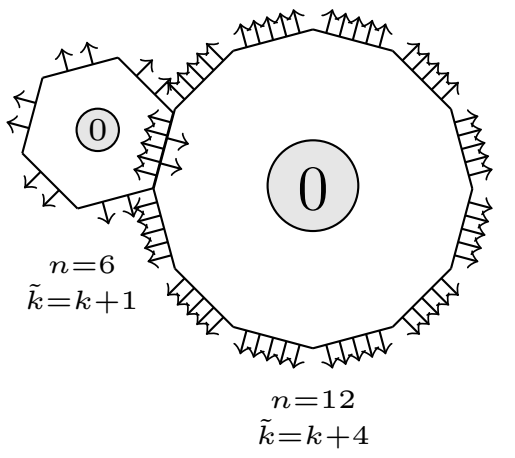

Fig. 48: Two - dimensional mesh with homogeneous order $\tilde{k} = 0$, built on elements that have six and twelve edges.

***Remark.*** Meshes that use different shapes of elements should only be defined from the highest starting order $k_0$ of the collection of reference element. In the *Figure 48*, the lowest possible discretisation order is $k = 4$, not $k = 1$. ▲

Once the first element $E_{k_0}$ is set up, then for any $k \geq k_0$ the discretisation requirement will be homogeneous for normal and internal moments. Furthermore, the number of added moments on the faces is linear with respect to $n$ (when adding an face) and $k$ (when increasing the order). Therefore, even if this method allows less available orders than the two previous ones, it preserves a meaning for $k$ when considering a specific polygon. Indeed, even if the order is somehow shifted, on can guess the discretisation quality just by knowing the order and the number of faces of the element without any advanced computation. Even if a low number of normal functions cannot necessarily be achieved, this is the best one can achieve when one wants to keep a classical $RT_k(K)$ definition.

However, it is not satisfying. Indeed, building meshes gathering several element types will again lead to the problem of a non - homogeneous discretisation of the inner part of the cell through the mesh (*see the Figure 37*) and you need to know the shape of the element (and thus $k_0$) to infer the discretisation quality of the inner cell through the value of $k$.

Therefore, when designing admissible elements for polytopes one has to pay attention that for a given order $k$ the number of internal and normal degrees of freedom is identical for any element having any number of faces $n$. For the sake of convenience, having a same starting element at $k = k_0 = 0$ for any element $K$ will also be target.

### 5.1.3 Definition of the lowest order space

To be able to sort our discretisation through the value of $k$ with respect to the quality of the approximation the space provides, we adopt temporarily the definition of $E_k(K)$ given in the last paragraph of the *Section 5.1.2* for the sake of the explanation. We denote by $RT_{k_0}(K)$ the lowest order space that can be designed.

Another limitation of using this simple fork of the quadrilateral elements definition occurs when considering the lowest order discretisation space. Indeed, as it may be convenient for some applications to have only normal moments we would like to have a space analogous to the quadrilateral $RT_0(K)$ one, where the corresponding element does not have internal moment.

However, while staying in the framework provided by the definition of $RT_k(K)$, there exists $n$ such that it is not possible to not have internal moments. As a consequence, for some $n$ no element similar to the quadrilateral $RT_0(K)$ element can be designed.

***Remark.*** By the *Property 5.1* it is enough to show that for some $n$, $E_{k_0}$ may have a non - empty set of internal functions. Indeed, $E_k(K)$ is not defined for earlier $k$ and the set of internal functions is monotonously expanding. ▲

***Example.*** *W*e detail two examples in two dimensions.

• (*Working example*) In two dimensions, if we consider the case $n = 6$, $E_1(K)$ is the lowest possible order element, and has no internal moments. It can then be assimilated to $RT_0(K)$ from this point of view (*see the Figure 49*).

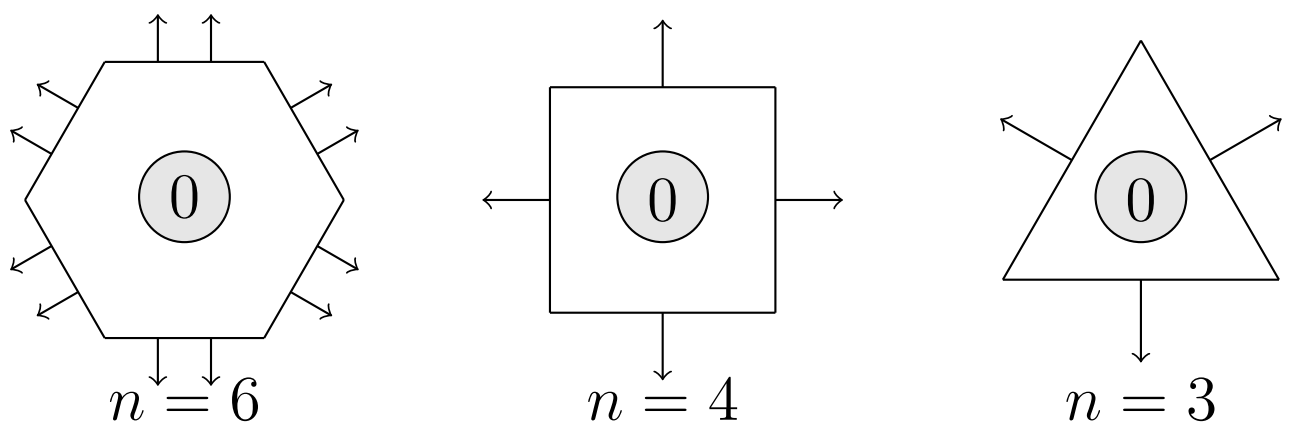

Fig. 49: Parallel with lowest order spaces

• (*Failing example*) When considering the case $n = 5$, no element can be assimilated with $RT_0(K)$. Indeed, the first element that can be designed is $E_1(K)$, $k_0 = 1$ being the smallest integer such that

$$\dim RT_k(K) = 2(k+1)(k+2) \geq 5(k+1) = \dim \mathcal{T}_k(\partial K).$$

Thus, defining the space $E_1(K)$, we have:

$$\dim RT_k(K) = 2(1+1)(1+2) = 12$$
$$\dim \mathcal{T}_k(\partial K) = 5(1+1) = 10.$$

Therefore, by distributing two normal functions on each edge we have a remaining of two internal basis functions. Thus, as the set of internal basis function is not empty, no parallel with the space $RT_0(K)$ can be done (*see the Figure 50*).

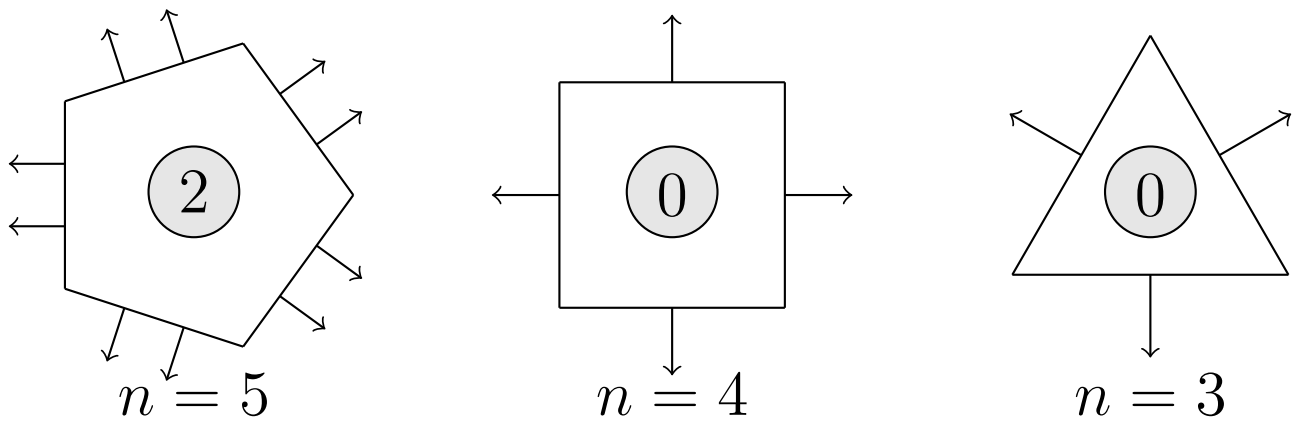

Fig. 50: Impossible parallel with the lowest quadrilateral $RT_k(K)$ space.
A more arithmetical way to see this impossibility is to notice that there exists no non - negative integer such that $2(k+1)(k+2) - 5(k+1) = 0$. And more generally, the equation $2(k+1)(k+2) - n(k+1) = 0$ does not always have a solution over the integers. ♦

Furthermore, note that no matter the number of faces we have (for $n \geq 2d+1$) the number of normal moments can never be one, and therefore can never be matched which is the number of normal degrees of freedom of $RT_0(K)$. Thus, the so - defined element $E_{k_0}(K)$ cannot always be linked or have the same properties as $RT_0(K)$, and when it is possible there is no particular meaning for the value of $k_0$.

We will then look for spaces $\mathbb{H}_k$ and $H(\text{div}, K)$ – conformal triplets so that the corresponding lowest order element does not have internal degrees of freedom.

### 5.1.4 Definition of the internal basis functions in two dimensions

As a consequence of the split between the internal and normal degrees of freedom, the direct generation of a basis of $RT_k(K)$ defining a conformal element built on polytopes is also challenged. We detail here specifically the two-dimensional case of polygons.

Indeed, when building the internal basis functions we have to ensure that their normal component will vanish on every edge to preserve the split. However, when using the classical definition of Raviart – Thomas spaces on general polytopes, the definition of $\mathcal{P}_k$ is delicate and its dimension is highly dependent on $n$. Therefore, we have to build $\dim \mathcal{P}_k$ functions that lie in the polynomial space $RT_k(K)$ and whose normal components are vanishing on the $n$ boundaries of $K$.

If it appears to be doable, constructing vanishing polynomials on the edges of arbitrarily polygons is not an easy task, especially when its degree is bounded by above. Indeed, when you try to define these functions from the product of linear functions vanishing on one edge of the polygon, the order of the obtained polynomial will prevent it to belong to the space $RT_k(K)$ for any $k \in \mathbb{N}$ lower than some specific order denoted $k_3 \in \mathbb{N}$.

***Example.*** *L*et us stay in two dimensions, with $k = 1$ and $n = 6$. Let $p_1$ be some polynomial that vanishes on the first edge, $p_2$ the one vanishing on the second edge and so on. Then by construction, the polynomial $\prod_{i=1}^{n} p_i$ would vanish on every edge. However, its degree would then be at least five, and would not lie in the space $RT_1(K)$. ♦

In some very specific cases, *e.g.* as when $n = 3$ or $n = 4$, it is possible to find an ideal of polynomials whose normal components vanish on the edges that are contained in $RT_k(K)$ for any $k$. One then can build on it a set of polynomials whose dimension matches the one of $\mathcal{P}_k$. However, this is far from being always possible and in those cases, the classical techniques such as the Groëbner basis reduction does not work either.

Therefore, this severe limitation prevents us to use this approach. We have to find another definition for $\mathbb{H}_k(K)$ so that basis functions whose normal component vanishes on every edge can be built for any polygons.

### 5.1.5 A last remark

As we aim to design a general setting, we leave down the option of designing one polytopial space per element to rather develop a discretisation space whose dimension is adaptive to the number of faces. Put together, the previous observations help up to design admissible elements, as presented in the next section.

## 5.2 Admissible spaces of discretisation

### 5.2.1 Conditions for admissible spaces

With respect to the limitations that were detailed previously and to the specifications that we would like to instil into $E_k$, we derive some conditions that the discrete space $\mathbb{H}^k$ has to fulfil and features that the element $E_k(K)$ must enjoy.

When building the space, we require the following conditions.

**Necessary conditions 5.2** Requirements on the discretisation space

1. The space $\mathbb{H}_k(K)$ should be vectorial and allow a definition of $H(\text{div}, K)$ – conformal elements for any discretisation order.
2. Its dimension has to be adaptive towards the number of faces of $K$ as well as toward the discretisation order $k$.

Furthermore, when building a $H(\text{div}, K)$ – conformal element $E_k(K)$ through the degrees of freedom, on has to fulfil the following criterion.

**Necessary conditions 5.3** Requirements on the elements

1. For any $k \in \mathbb{N}$, there exists a set of degrees of freedom unisolvent for $\mathbb{H}_k(K)$ that can be clearly subdivided into internal and normal subsets. In addition, the structure of the degrees of freedom should be identical for any order $k$.
2. The degrees of freedom are defined through intelligible definitions of $\mathcal{P}_{b,k}$ and $\mathcal{P}_k$ that quantify target aspects of the discretised quantity.
3. The order $k$ is meaningful regardless the shape of the element.

More precisely, for a given order $k$ we ask every element with any number of faces $n$ to share the same number of normal degrees of freedom per face and the same number of internal degrees of freedom (*see the Figure 51*).

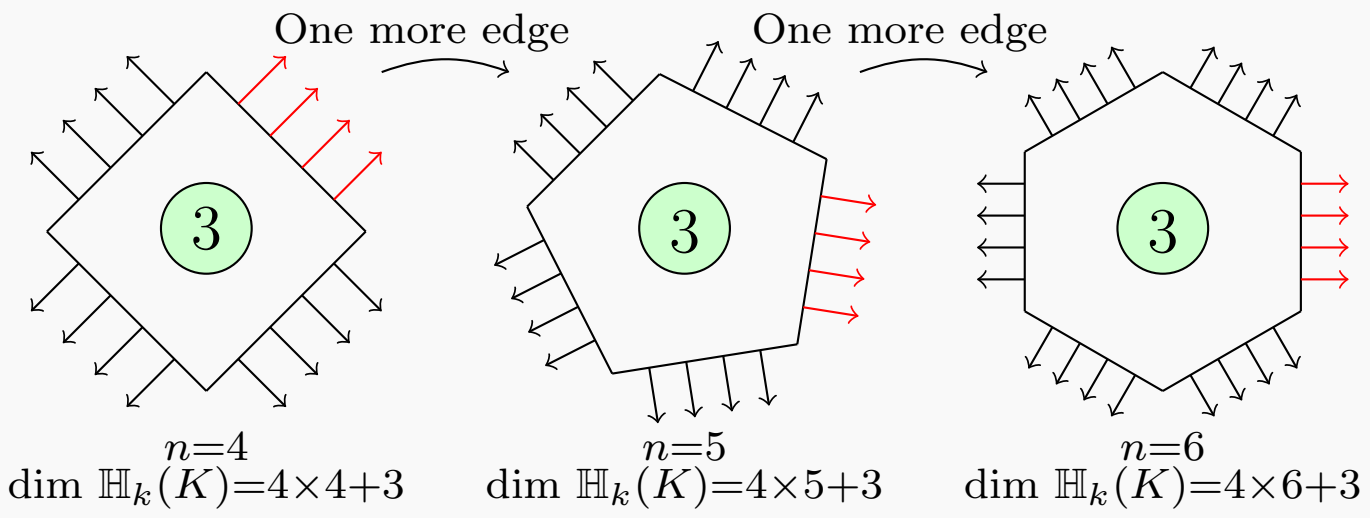

Fig. 51: : Adaptiveness of the number of degrees of freedom

4. For a fixed order $k$, the increment in the dimension of the set of normal degrees of freedom matches the increment in the dimension of $\mathbb{H}_k(K)$ each supplementary face.

5. For a fixed number of faces $n$, the number of normal moments per face increases strictly monotonously with $k$ in an understandable way (*i.e.* either linearly or following a predefined refinement sequence).

6. The refinement performed in the inner cell is strictly monotonous and is done in an understandable way with respect to the increment of the order $k$. In particular, there is no order $k$ where sudden or particularly fine refinement appears when compared to the discretisation obtained for the order $k-1$ (*see the Figure 52*).

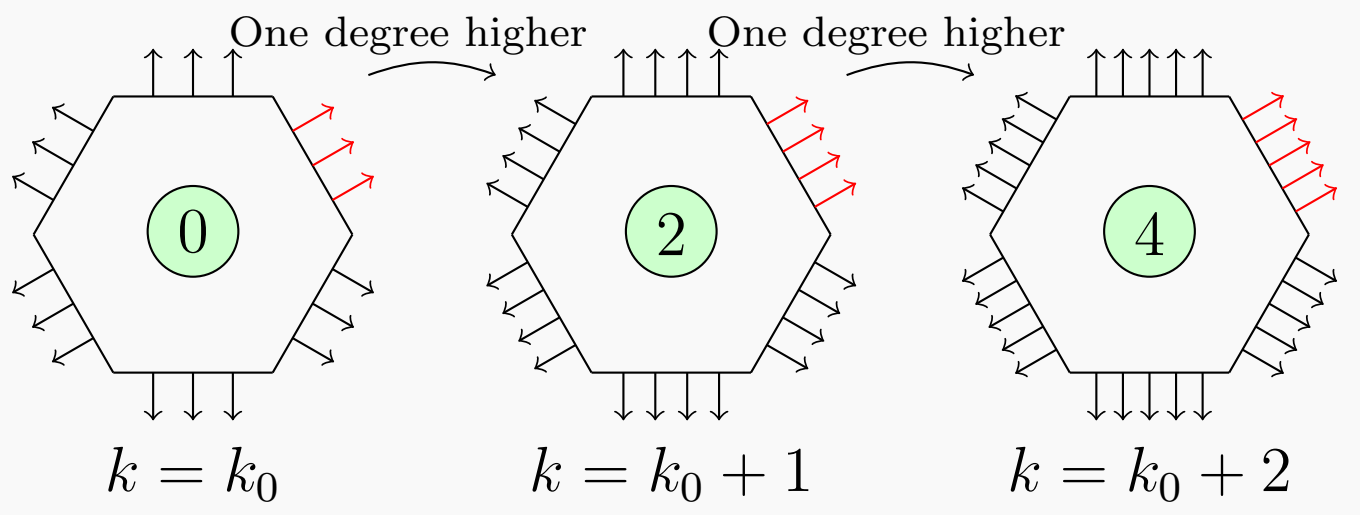

Fig. 52: :Intelligible discretisation refinement through the increase of the order

Lastly, when considering the dual space to the set of degrees of freedom, one has to make sure that the following requirements hold.

**Necessary condition 5.4** Requirements on the basis functions

The corresponding internal basis functions vanish on the faces of the element.

Additionally, when one wishes a parallel with the Raviart–Thomas elements from the lowest order on, one should require the following additional requirement.

**Optional condition 5.5** Requirements on the elements

The lowest order element has no internal degrees of freedom.

We recall that in order to obtain more easily $H(\text{div}, K)$ – conformity, we prescribe $P_{b,k}(f)$ to be any polynomial subspace of $\mathbb{P}(\mathbb{R}^{d-1})$ providing a discretisation order $k$ (*see Section 2.1*). Note also that though the framework we propose furnishes a discretisation setting for any polytope, it is specially thought for cells having more than five faces ($n \geq 5(d-1)$).

### 5.2.2 A class of admissible spaces

#### 5.2.2.1 Construction of the space

The wish of the possibility to design $H(\text{div}, K)$ – conformal elements expressed in second point of 5.2 leads us to define the space $\mathbb{H}_k(K)$ in the same spirit as for the classical Raviart – Thomas case. Meaningly we set

$$\mathbb{H}_k(K) = (A_k)^d + x\, B_k$$

for two given functional sets $A_k$ and $B_k$.

In addition, the eighth point of 5.3 makes us consider the spaces $A_k$ and $B_k$ based on solutions to Poisson's equation, in the same way as in the definition of the elements used in the VEM method [20]. Indeed, a way to force the internal basis functions to vanish on the boundary of any polygon is to define them from the set of solutions to the problems

$$\begin{cases} \Delta u = p_k, \\ u|_{\partial K} = 0. \end{cases} \tag{60}$$

There, $p_k$ stands for any function belonging to some functional space $\tilde{\mathcal{P}}_k$ that remains to be determined. Note however that independently of the definition of $\tilde{\mathcal{P}}_k$, the vanishing property of the solutions to (60) is naturally transmitted to functions belonging to the space $\mathbb{H}_k(K)$. Indeed, if $A_k$ and $B_k$ are built as being subsets of solutions to those equations, the functions in the space $A_k^d + x\,B_k$ enjoys the same vanishing property on the boundary.

Then, to allow the $H(\text{div}, K)$ – conformity to be enforced by normal quantities emerging only from the boundaries while proposing various qualities of discretisation on the faces, we also need to consider functions of the type

$$\begin{cases} \Delta u = 0 \\ u|_{\partial K} = p_k \mathbb{1}_f \end{cases} \tag{61}$$

for any face $f$ in $\partial K$ and any function $p_k$ belonging to the space $\mathcal{P}_{b,\,k}(f)$. Note that the discontinuity of the boundary conditions is required to tune the behaviour of the discretised quantity facewise. Thus, the $H(\text{div}, K)$ – conformity is still preserved as the normal component of the solution is facewise continuous. Furthermore, as the boundary functions are still in $L^2(\partial K)$, we still get uniqueness of the solution to the problem (61) (*see e.g. [5]*).

***Note.*** For the sake of concision, in all the following we denote by Poisson's function the solution to some Poisson's problem under some Dirichlet boundary conditions. ▲

As we want to enjoy some smoothness in the projection space, we will only consider spaces $\tilde{\mathcal{P}}_k$ and $\mathcal{P}_{b,\,k}$ that are polynomial. In particular, as we focus specifically on elements having more than $2d+1$ faces, we will consider $\tilde{\mathcal{P}}_k = \mathbb{Q}_k(K)$ and $\mathcal{P}_{b,\,k} = \mathbb{Q}_k(f)$, though at this point any other definition would be admissible. As a consequence, the definition (59) of $\mathcal{H}_k(\partial K)$ turns to read

$$\mathcal{H}_k(\partial K) = \{u|_{\partial K} \in L^2(\partial K),\, u|_f \in \mathbb{Q}_k(f), \quad \forall f \in \partial K\}. \tag{62}$$

Using the superposition principle, we can then define a space enjoying the desired properties 5.2 by writing

$$\begin{aligned}\mathbb{H}_k(K) =& \{u \in H^1(K),\, u|_{\partial K} \in \mathcal{H}_{l_1}(\partial K),\, \Delta u \in \mathbb{Q}_{m_1}(K)\}^d \\ &+ x\,\{u \in H^1(K),\, u|_{\partial K} \in \mathcal{H}_{l_2}(\partial K),\, \Delta u \in \mathbb{Q}_{[m_2]}(K)\}\end{aligned} \tag{63}$$

for some integers $l_1$, $l_2$, $m_1$ and $m_2$, and where $H^1$ is the minimum admissible regularity to fit our framework. Indeed, weakening it prevents from achieving $\mathbb{H}_k(K) \subset H(\text{div}, K)$ as the divergence of functions living in $\mathbb{H}_k(K)$ would be less controlled and the setting may allow inappropriate behaviours.

We adopt the convention $\mathbb{Q}_{-1} = \{0\}$, and set $A_k$ or $B_k$ to the empty space when one of the indices defining either the second members or the boundary conditions defining the functions of the subspace are strictly below $-1$.

***Remark.*** The spaces $\mathbb{H}_k(K)$ are constructed by blocks. Indeed, the liberty in their definition is restricted to the choices of spaces $\mathcal{P}_{b,l_1}$, $\mathcal{P}_{b,l_2}$ and by the use of $\mathbb{Q}_{m_1}$ and $\mathbb{Q}_{[m_2]}$. Furthermore, even if there exist combinations of choices that are more suitable than others (*see Section 5.5*) each of the definition of those four spaces does not impact the other ones, and the space $\mathbb{H}_k(K)$ can be defined for any of their combinations. ▲

Let us point out that the orders $l_1$, $l_2$, $m_1$ and $m_2$ of each block are not necessarily identical and that the space is based on solutions to (60) and (61). As a consequence, one cannot directly get the discretisation order of the space *a - priori* simply by reading out the orders that are in use in each block. Therefore, when speaking about the order of the space will from now on refer to the following definition.

**Definition 5.6** Order of the space $\mathbb{H}_k(K)$

The order of the space $\mathbb{H}_k(K)$ is the couple $(k_1,\, k_2)$ of the greatest integers such that

$$\{u \in H^1(K),\, u|_{\partial K} \in \mathcal{H}_{k_1}(\partial K),\, \Delta u \in \mathbb{Q}_{k_2-1}(K)\} \subset \mathbb{H}_k(K).$$

***Note.*** Linking the coefficients $m_1$, $m_2$ with $l_1$, $l_2$ through a linear relation allows to represent the order of the space as a scalar integer. By example, using $m_1 = m_2 = k-1$ and $l_1 = 0$, $l_2 = k$ for any $k \in \mathbb{N}$, the above definition leads to an order of $(k,\, k)$ that can be abbreviated by a scalar order $k$. ▲

***Remark 5.***

• One can notice that all of the spaces (63) prevent their functions to be all continuous. Indeed, asking continuity across the polytope's vertices reduces to use the space

$$\mathcal{H}_k(\partial K) = \{u|_{\partial K} \in \mathcal{C}_0(\partial K),\, u|_f \in \mathbb{Q}_k(f), \quad \forall f \in \partial K\}$$

instead of (62) in both $A_k$ and $B_k$. There, as the boundary of the polytope is closed, $\dim \mathbb{H}_k(K)$ would be reduced by $(n-2)\dim \mathbb{Q}_0(f) + (\dim \mathbb{Q}_1(f))$ if $k > 0$ or $(n-1)\dim \mathbb{Q}_0(f)$ if $k = 0$. The dependency of the space's dimension would then not be straightforwardly dependent on the number of faces, if it is at all (*see the Figure 53*). Therefore, a general setting is impossible to draw. This is easily observable for the case $k = 0$. There, by continuity only one constant remains to determine, which is obviously independent of the number of faces.

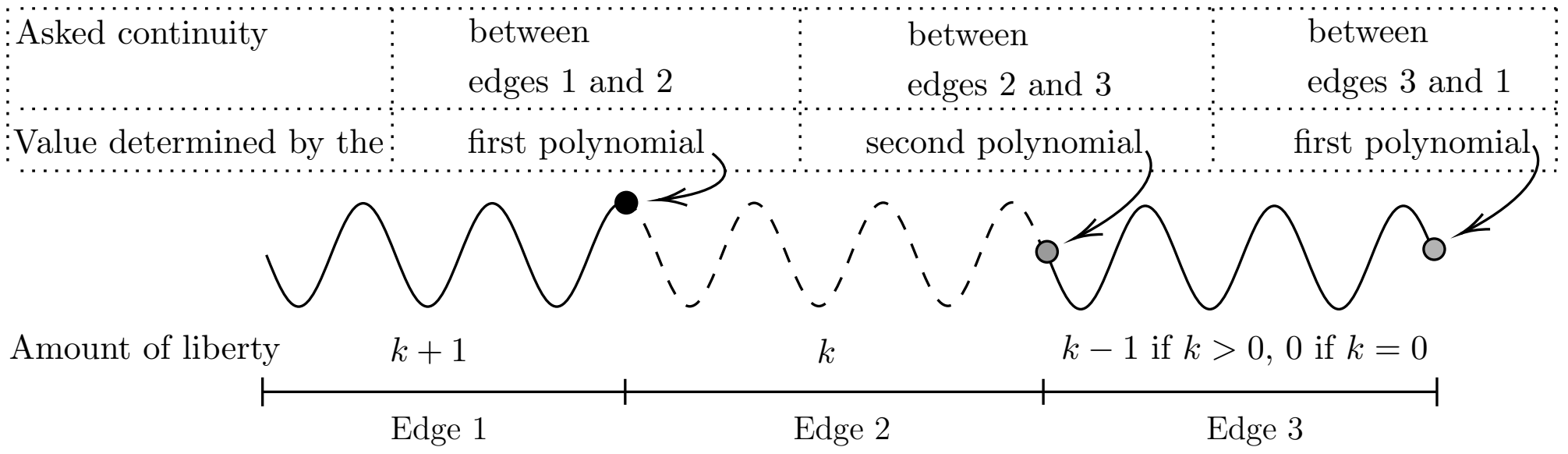

Fig. 53: Two dimensional example of the profile of a function in $\mathcal{H}_k(\partial K)$. Split view of a triangle's edges where continuity is asked across the vertices. The remaining freedom per edge is dependent on the neighbours. Here, the choices are initiated on the first edge.

• In the same spirit, no space based on

$$\mathcal{H}_k(\partial K) = \{u|_{\partial K} \in \mathbb{Q}_k(\partial K)\}$$

for both $A_k$ and $B_k$ can be designed, as here there is no dependency on the number of faces at all.

• A hybrid version would however be admissible, where the above space is used in either $A_k$ or $B_k$. It may be useful in particular when using

$$\mathcal{H}_k(\partial K) = \{u|_{\partial K} \equiv 1\}$$

as this term does not generate further freedom, the boundary behaviour being prescribed. The dependency on the number of faces is then only achieved through the other subspace of $\mathbb{H}_k(K)$, and becomes straightforward. As this case is a pure restriction of the general case designed here, it will not be discussed. Yet, an example will be given in the *Section 6.2*. ▲

By using this space's construction, we will see in the *Section 5.2.2.2* that the first point of the requirements 5.2 and the first to third points of the requirements 5.3 are naturally fulfilled. The last three requirements of 5.2 and the requirement 5.4 depend directly on the element, relying on the definition of the degrees of freedom. They will be discussed in the *Section 5.3*. The optional requirement 5.5 is achieved by a structural change, *i.e.* when using a hybrid version quickly evoked above. It is detailed as a reduced space in the *Section 6.2*.

#### 5.2.2.2 Common properties of $\mathbb{H}_k(K)$ spaces

We can already derive some properties that are common to those spaces regardless the values of $l_1$, $l_2$, $m_1$ and $m_2$.

### Divergence properties

**Property 5.7** $H(\mathrm{div},\, K)$ – conformity

For any function $q$ belonging to any space $\mathbb{H}_k(K)$, it holds:

$$\begin{cases} q \cdot n|_{\partial K} \in \mathbb{H}_{\max\{l_1,l_2\}}(\partial K) \\ \quad \operatorname{div} q \in L^2(K) \end{cases}$$

***Proof.***

• We start by deriving the first statement. By construction, any $q \in \mathbb{H}_k(K)$ can be decomposed into $q = q_0 + x\, q_1$ for some $q_0 \in (A_k)^d$ and $q_1 \in B_k$. Therefore, on the boundary of $K$ one has $q \cdot n|_{\partial K} = q_0 \cdot n|_{\partial K} + (x q_1 \cdot n)|_{\partial K}$. As the functions $q_1$ is scalar, this quantity can also read $q \cdot n|_{\partial K} = q_0 \cdot n|_{\partial K} + q_1 (x \cdot n)|_{\partial K}$ by linearity and commutativity of the dot product.

Since for every face $f$ of $K$ the term $x \cdot n|_f$ is constant, it reduces to $q \cdot n|_f = q_0 \cdot n|_f + c_f\, q_1|_f$ on each face $f$ for a constant $c_f \in \mathbb{R}$ depending only on the face layout and position with respect to the axes and origin. Therefore, since $q_0|_f \in (\mathbb{Q}_{l_1}(f))^d$ and $q_1|_f \in \mathbb{Q}_{l_2}(f)$, $q \cdot n|_f \in \mathbb{Q}_{\max\{l_1,l_2\}}(f)$. And since it is valid for any face $f \in \partial K$, we finally get that $q \cdot n|_{\partial K} \in \mathbb{H}_{\max\{l_1,l_2\}}(\partial K)$.

• Let us now derive the divergence property within the cell. Any $u \in \mathbb{H}_k(K)$ can be written under the form

$$u = \tilde{q} + x\, q$$

for some functions $q \in H^1(K)$ and $\tilde{q} = (\tilde{q}_1, \cdots, \tilde{q}_d)^T \in (H^1(K))^d$ such that

$$\begin{cases} \Delta q \in \mathbb{Q}_{[m_2]}(K) \\ q|_{\partial K} \in \mathbb{Q}_{l_2}(\partial K) \end{cases} \quad \text{and} \quad \begin{cases} \Delta \tilde{q}_i \in \mathbb{Q}_{m_1}(K) \\ \tilde{q}_i|_{\partial K} \in \mathbb{Q}_{l_1}(\partial K), \end{cases} \quad \forall i \in [\![1,\, d]\!]. \tag{64}$$

Thus,

$$\operatorname{div}(u) = \nabla \cdot u$$

$$= \sum_{i=1}^{d} \partial_{x_i}(x_i\, q) + \sum_{i=1}^{d} \partial_{x_i}\tilde{q}_i$$

$$= \sum_{i=1}^{d} (q + x_i\partial_{x_i} q) + \sum_{i=1}^{d} \partial_{x_i}\tilde{q}_i$$

$$= \underbrace{d\, q}_{\in L^2(K)} + \sum_{i=1}^{d} \big(x_i \underbrace{\partial_{x_i} q}_{\in L^2(K)}\big) + \sum_{i=1}^{d} \underbrace{\partial_{x_i}\tilde{q}_i.}_{\in L^2(K)}$$

Since by (64) we have $\nabla \cdot q \in L^2(K)$, it comes that for any $i \in [\![1,\, d]\!]$;

$$x_i\, \partial_{x_i} q \in L^2_{\text{loc}}(K).$$

As $K$ is compact and bounded, we have $L^2_{\text{loc}}(K) = L^2(K)$ and $\operatorname{div} q \in L^2(K)$. As a by-product, note that we can derive

$$\nabla \cdot (x \nabla q) = \nabla \cdot q + \underbrace{x \underbrace{\Delta q}_{\in \mathbb{Q}_{\max\,[m_1,\, m_2+1]}}}_{\in \mathcal{C}^\infty(K)} \in L^2(K).$$

■

We can directly state the following implication.

**Corollary 5.8** $H(\text{div},\, K)$ – conformity

For any $k \geq 0$,

$$\mathbb{H}_k(K) \subset H(\text{div},\, K)$$

By the above *Proposition 5.7*, the functions belonging to those spaces enjoy the continuity of their normal components through the faces of $K$. Therefore, as by the above statement $\mathbb{H}_k(K)$ is a subset of $H(\text{div},\, K)$, it allows to build $H(\text{div},\, K)$ – conformal elements. Observing furthermore that the space is vectorial by definition, the construction (63) verifies the second point of 5.2.

### Dimension

We start by noticing a property that will be useful to determine the dimension of the space $\mathbb{H}_k(K)$ and make the construction of basis functions for elements constructed on $\mathbb{H}_k(K)$ easier.

**Property 5.9**

The two subspaces arising in the definition of

$$\begin{aligned}\mathbb{H}_k(K) =&\{u \in H^1(K),\, u|_{\partial K} \in \mathcal{H}_{l_1}(\partial K),\, \Delta u \in \mathbb{Q}_{m_1}(K)\}^d \\ &\oplus x\{u \in H^1(K),\, u|_{\partial K} \in \mathcal{H}_{l_2}(\partial K),\, \Delta u \in \mathbb{Q}_{[m_2]}(K)\}\end{aligned}$$

are in direct sum.

***Proof.*** As we enjoy a block construction, the impact of $l_1$ and $l_2$ is restricted to driving the behaviour on the boundary of the functions living respectively in $\mathcal{A}_k$ or $B_k$. Their inner parts are only driven by the values of $m_1$ and $m_2$ through the elliptic equations (60). Thus, to prove that the two spaces $(\mathcal{A}_k)^d$ and $xB_k$ are in direct sum, it is enough to show that no function satisfying

$$\begin{cases} \Delta u = p_k \\ u|_{\partial K} = 0 \end{cases}$$

for any $p_k \in \mathbb{Q}_{m_1}(K)$ can be matched with a coordinate of any function living in $xB_k$, and *vice - versa*.

Therefore, we neglect the values of $l_1$ and $l_2$ and discuss without loss of generality the assertion depending on the relationship between $m_1$ and $m_2$ when considering only homogeneous Dirichlet conditions in the definitions of $\mathcal{A}_k$ and $B_k$.

The proof is split into two parts. We first treat the natural case $m_1 \leq m_2$ before treating the case $m_2 < m_1$. In order to make later explanations easier, let us first recall the Laplacian computation

$$\begin{aligned}\Delta(x_i u) &= \Delta(x_i)u + 2\nabla(x_i)\cdot\nabla(u) + x_i\Delta(u) \\ &= 2\nabla\cdot u + x_i\Delta(u), \end{aligned} \tag{65}$$

and state that two functions with different gradients cannot be identical.

• We consider first the case $m_1 \leq m_2$. Let us take some $v \in (A_k)^d \subset \mathbb{H}_k(K)$. Then, each coordinate of $v$ has its Laplacian in $\mathbb{Q}_{m_1}(K)$.

Let us assume that $v$ can be factorized similarly in each coordinate so that it can be recast into $v = xu$ for some scalar function $u$. We can now consider the coordinate wise functions $\{v_i\}_{i\in[\![1,\, d]\!]} = \{x_i u\}_{i\in[\![1,\, d]\!]}$. Then, one has by (65) that $\Delta(v_i) = \Delta(x_i u) = 2\nabla\cdot u + x_i\Delta(u)$ for any $i \in [\![1,\, d]\!]$. However, $\Delta(v_i) = \Delta(x_i u) \in \mathbb{Q}_{m_1}(K)$, and therefore so is $x_i\Delta(u)$. Thus, $\Delta(x_i u) = \Delta(v_i) \in \mathbb{Q}_{\zeta_i(m_1-1,\, m_1,\, \cdots,\, m_1)}(K)$ for any $i \in [\![1,\, d]\!]$. Lastly, since $u$ is a scalar function appearing identically for each coordinate, we get $u \in \mathbb{Q}_{m_1-1}(K) \not\subset \mathbb{Q}_{[m_2]}(K)$. Thus, if $v \in (A_k)^d$, $v$ cannot be recast into the form $v = x\,u$ with $u \in B_k$.

One could argue that this impossibility arises due to the strong assumption made on the shape of $v$. However, dropping the assumption that $v$ can be recast into $v = x\,u$ would prevent the structure of $v$ to match with the one of any functions in $xB_k$ and we end up with the same conclusion.

Reversely, if $v \in xB_k$, then still by (65) we get $\Delta v = 2\nabla \cdot v + x_i\Delta(v)$, where $\Delta(v)$ belongs to $\mathbb{Q}_{[m_2]}(K)$. Therefore, $x_i\Delta(v)$ belongs to $\mathbb{Q}_{\zeta_i([m_2+1],m_2,\cdots,m_2)}(K) \not\subset \mathbb{Q}_{m_1}(K)$ and no coordinate of $v$ can match the one of some function in $A_k$. Thus $v \notin (A_k)^d$.

Combining the two above facts, the space $\mathbb{H}_k(K)$ is in direct sum for any $m_1 \leq m_2$ .

• Let us now treat the case $m_2 < m_1$. There, we have $\mathbb{Q}_{m_2}(K) \subset \mathbb{Q}_{m_1}(K)$. Therefore, the set of Dirichlet homogeneous solutions of the set $B_k$ is included in the one of $A_k$. However, we can notice that $\mathrm{div}(x\,u) = du + \sum_{i=1}^{d}(x_i\partial_{x_i} u)$, which is not equal to $c\,\mathrm{div} u$ for any constant $c \in \mathbb{R}$ provided that $u$ is not entirely null and that the domain on which the equation live is not reduced to a point. As in $B_k$ the case $\Delta u \equiv 0$ is not possible, those two cases never happen and the two divergences $\mathrm{div}\, u$ and $\mathrm{div}\,(xu)$ always differ. Therefore, no function of $xB_k$ can be identified with one living in $(A_k)^d$ and *vice - versa*. Thus, those two subspaces are in direct sum, which concludes the proof.

■

***Note.*** The restriction of $\mathbb{H}_k(K)$ to the boundary, is not a direct sum of the two subspaces $A_k$ and $B_k$. Indeed, in the case where $l_2 < l_1$, one could form on the boundaries a function of $x\mathbb{Q}_{l_2}(f)$ from $d$ functions of $\mathbb{Q}_{l_1}(f)$. This is visible in the decomposition of the polynomials living on the boundary in the *Section 5.3.2.* ▲

Making use of this property and recalling that we enjoy a block construction in the definition of $\mathbb{H}_k(K)$, the dimension of the space $\mathbb{H}_k(K)$ can be easily derived. Indeed, as the two natural subspaces are in direct sum, we can simply add the dimension of the two main subspaces $(A_k)^d$ and $B_k$ to retrieve the dimension of $\mathbb{H}_k(K)$. Let us derive their respective dimensions.

First, we compute the dimension of $A_k$. In the way presented in [20], one can get it by using the superposition theorem. Indeed, for any second member belonging to $\mathbb{Q}_{m_1}$ there exists a unique solution to the Dirichlet problem (60). Similarly, for any face $f$ and any polynomial $p_k \in \mathbb{Q}_{l_1}(f)$ defining the boundary function $p_k \mathbb{1}_f \in L^2(K)$, there is a unique solution to (61) (*see e.g. [5]*). Thus, reading out the structure of the set $A_k$ implies the following relation.

$$\begin{aligned}\dim A_k &= \dim \mathcal{H}_{l_1}(\partial K) + \dim \mathbb{Q}_{m_1}(K)\\ &= n(l_1+1)^{d-1} + (m_1+1)^d\end{aligned}$$

Therefore, as $(A_k)^d$ is a simple Cartesian product of $d$ copies of $A_k$, we have immediately $\dim A_k = d(\dim A_k) = d(n(l_1+1)^{d-1} + (m_1+1)^d)$. In the exact same way, we retrieve the dimension of $B_k$ by

$$\begin{aligned}\dim B_k &= \dim \mathcal{H}_{l_2}(\partial K) + \dim \mathbb{Q}_{[m_2]}(K)\\ &= n(l_2+1)^{d-1} + (m_2+1)^d - m_2^d.\end{aligned}$$

Lastly, we recall that the space $x\,B_k$ simply corresponds to an identical $d$ - duplication of the space $B_k$ where each coordinate has been multiplied by the corresponding spatial variable. Thus, there is no liberty adjunction during its construction, and the dimension of $x\,B_k$ equals the one of $B_k$.

Bringing together the dimensions of the two subsets we finally retrieve

$$\begin{aligned}\dim \mathbb{H}_k(K) &= d \dim A_k + \dim B_k\\ &= d(n(l_1+1)^{d-1} + (m_1+1)^d) + n(l_2+1)^{d-1} + (m_2+1)^d - m_2^d.\end{aligned}$$

Reordering the terms, it comes the following property.

**Property 5.10** Dimension

The space $\mathbb{H}_k(K)$ enjoys the dimension

$$\begin{aligned}\dim \mathbb{H}_k(K) =& n\,\left(d(l_1+1)^{d-1} + (l_2+1)^{d-1}\right)\\ &+ \left(d(m_1+1)^d + (m_2+1)^d - m_2^d\right). \qquad (66)\end{aligned}$$

Let us now point out the following observations, emphasizing that the spaces $\mathbb{H}_k(K)$ fit in our desired framework.

First, we can directly notice that the dimension of the space is adaptive to both $n$ and $k$, which fulfils the second point of 5.2. Furthermore, a clear subdivision between internal and normal degrees of freedom in the spirit of the *Section 4.3.1* is possible when one prescribes $d(l_1+1)^{d-1} + (l_2+1)^{d-1}$ degrees of freedom on each face. The first point of 5.3 is then immediately verified.

On the side of the internal degrees of freedom, one can observe that adding one face to the element at a fixed order $k$ does not create supplementary internal degrees of freedom. As a consequence, the dimension of $\mathbb{H}_k(K)$ increases only by $d(l_1 + 1)^{d-1} + (l_2 + 1)^{d-1}$, which by construction matches the dimension of the normal degrees of freedom on one face. This fulfils the fifth point of 5.3.

As a by product, the number of normal moments per face is identical for any number of faces at $d$ and $k$ fixed. The same goes for the internal degrees of freedom as their number is totally independent from $n$. This fulfils the third point of 5.3. Furthermore, the polynomial discretisation on the boundary is of the same order on every face.

Lastly, the refinement of the normal degrees of freedom is linear with respect to $n$. Therefore, thanks to the definition of the spaces $\mathbb{Q}_{l_1}(f)$ and $\mathbb{Q}_{l_2}(f)$ the refinement process is monotonous and understandable in $k$ for any fixed $n$, fulfilling the sixth point of 5.3.

Enjoying the previous property helps us to design properly the degrees of freedom. However, a discussion on the admissible values of $l_1$, $l_2$, $m_1$ and $m_2$ is required to draw a general framework.

## 5.3 Definition of admissible elements

### 5.3.1 Admissible orders

In order to fulfil the three last points, one has to pay attention to the combination of the orders $l_1$, $l_2$ and $m_1$, $m_2$. Indeed, as it can be seen in the *Property 5.9*, their layout combined with the value of $d$ drives the distribution of the degrees of freedom into the normal and internal classification. In particular, as we work with polynomials on the boundary, the interaction between $l_1$ and $l_2$ may end to an amount of boundary degrees of freedom that would imply an over-determination. This typically occurs when the amount of liberties per face is higher than the dimension of the polynomial space in us in the definition of $\mathbb{H}_k(K)$. The definition of $m_1$ and $m_2$ may also prevents the definition of a lowest order space allowing a parallel with $RT_0(K)$. We then derive conditions ensuring that we will fall in an admissible case.

The first restriction comes from the wish of making a connection with $RT_0(K)$ for the lowest order element. As it should not have internal degrees of freedom to fulfil the sixth point of 5.3, the set $\tilde{\mathcal{P}}_0$ should be empty or reduced to one unique solution. Thus, we set $\mathbb{Q}_{-1} = \{0\}$ and allow $m_1$ and $m_2$ to take values only from $-1$ on. For the same purpose, we adopt as before the convention that for either $m_1 < -1$ or $m_2 < -1$ the two subspaces are respectively set to $A_k = \varnothing$ and $B_k = \varnothing$.

We also like to define spaces that provides several discretisation qualities built on the same spirit. Therefore, to make the discussions easier and allow an understandable refinement rule, we define the spaces starting from a lowest order space beginning at $k = 0$ and consider $l_1$, $l_2$, $m_1$ and $m_2$ through the affine relations $l_i = a_i k + b_i$, $m_i = c_i k + d_i$ for natural numbers $a_i$, $b_i$, $c_i$, $d_i$, $i \in \{1, 2\}$.

Then, for the sake of giving a meaning to the order of the space we forbid a stagnation in the refinement for both internal and normal discretisations. Therefore, we cannot set $m_2$ and $m_1$ both lower or equal than $k-2$. Indeed, one would have otherwise no internal moment for both $k = 0$ and $k = 1$, implying that no refinement would occur.

Furthermore, we wish to avoid holes in the projection space. By this, we mean that we would like to forbid $\mathbb{H}_k(K)$ to be described by functions that are solution to (60) for a polynomial space $\tilde{\mathcal{P}}$ gathering monomials whose degrees are not contiguous. This would indeed not be conceivable in terms of the categorization of the spaces.

***Example.*** *B*y example, spaces gathering the solutions of Laplacian problems of the form $\Delta u = 1$ and $\Delta u = x^2 y^2$ without providing any of other intermediate orders polynomials as a second member are not eligible. ♦

However, we accept sequences of contiguous degrees even if they are not complete with respect to the definition of the $\mathbb{Q}_k(K)$ spaces, as by example $\{1, xy, x^2y^2\}$. Fortunately, the spaces $\mathbb{Q}_{m_1}$ and $\mathbb{Q}_{[m_2]}$ satisfy this criteria. The restriction on the admissible order will come from the spaces $(A_k)^d$ and $xB_k$. Indeed, even if they do not share the same structure, the space $\mathbb{Q}_{m_1} \cup \mathbb{Q}_{[m_2]}$ should not contain holes either as we would like to have no hole in the polynomial part of the projection space representing the divergence of the discretised quantity. Since $\nabla \cdot (xu) = \sum(\frac{\partial(x_i u)}{\partial x_i}) = \sum(u + x_i \frac{\partial u}{\partial x_i})$ and $\Delta(xu) = 2\nabla \cdot u + (\sum_{i=1}^d x_i)u$, we need $x\mathbb{Q}_{[m_2]} \cup \mathbb{Q}_{m_1}$ not to have holes either. As we are in a direct sum for any $m_2 \geq m_1$, prescribing a too large $m_2$ with respect to $m_1$ will leave down all Poisson's functions based on the intermediate polynomials. Therefore, we will consider the spaces where $m_2 \leq m_1$. Furthermore, we would also like to keep the preservation of the discretisation quality through the divergence and will favour the spaces having $m_2 \geq m_1$. Therefore, though it is not mandatory, we strongly advise to set $m_2 = m_1$.

Lastly, no specific restriction is made *a - priori* on $l_1$ and $l_2$ besides asking them to be bigger or equal to $-1$ in order to define properly the boundary conditions to the Laplacian problems. Yet, on a practical side one must take $l_1 = 0$ or $l_1 = -1$. Indeed, asking higher variations would make issues arise in the definition of the degrees of freedom when looking for unisolvence. We detail it in the next paragraph.

### 5.3.2 Definition of admissible normal degrees of freedom

We now design degrees of freedom that totally defines a $H(\text{div}, K)$ – conformal element on any polytope.

Let us start by defining the normal degrees of freedom. Following the *Paragraph 5.2.2.2*, we should set $d(l_1+1)^{d-1} + (l_2+1)^{d-1}$ of those on each face to preserve the adaptivity of the element to arbitrary polygonal shapes.

Since we enjoy a block construction, the role of the normal degrees of freedom is restricted to determining the boundary parts of functions living in $\mathbb{H}_k(K)$. Therefore, their support is set to be the boundary of $K$, where functions of $\mathbb{H}_k(K)$ are reduced to polynomials. More precisely, they belong to a subspace of $\mathbb{H}_k(K)|_{\partial K} = (\mathcal{H}_{k+1}(\partial K))^d$ within each the directional components are polynomial functions of degree at most $\max\{l_1,\, l_2+1\}$. More precisely, any $p \in \mathbb{H}_k(K)|_f$ writes

$$p = \bigtimes_{j=1}^{d} \left( \sum_{\substack{|\alpha_{i,j}| \leq \\ \max\{l_1, l_2+1\}}} c_{i,j} x^{\alpha_{i,j}} \right) \tag{67}$$

for multi - indices $\{\alpha_{i,j}\}_{i,j}$, constants $\{c_{i,j}\}_{i,j}$ and where $|\alpha_{i,j}| = \max_i\{(\alpha_{i,j})_i\}$.

Therefore, determining the boundary part of a function living in $\mathbb{H}_k(K)$ comes down to determining $d$ polynomials of degree $\max\{l_1,\, l_2+1\}$ on each face. Thus, *a - priori* the set of normal degrees of freedom per face can be at maximum as large as the number of coefficients required to determine those polynomials living on a face $f$, that is $d \dim \mathbb{Q}_{\max\{l_1, l_2\}}(f) = d(\max\{l_1,\, l_2+1\} + 1.)^{(d-1)}$

***Remark.*** This limitation on the degree of the polynomial restrictions on the faces of functions in $\mathbb{H}_k(K)$ provides a first criteria for the admissibility of the couple $(l_1, l_2)$. Indeed, we should verify

$$d(l_1+1)^{d-1} + (l_2+1)^{d-1} \leq d(\max\{l_1, l_2+1\}+1)^{d-1}_{.}$$

In all the next, we will assume that this assumption is fulfilled. Note however that this is not yet a sufficient condition allowing to define properly the element $E_k(K)$. ▲

Taking a closer look on the structure of the space $\mathbb{H}_k(K)$, one can reduce the need of determination. Indeed, any function of $(A_k)^d$ reads

$$p \in (A_k)^d|_f \Rightarrow p = \bigtimes_{j=1}^{d} \left( \sum_{|\alpha_{i,j}| \leq l_1} a_{i,j} x^{\alpha_{i,j}} \right)$$

for a given set of multi - index $\{\alpha_{i,j}\}_{i,j}$ and coefficients $\{a_{i,j}\}_{i,j}$ depending on the coordinates. Therefore, one needs to prescribe the $(l_1+1)^{d-1}$ monomial coefficients $\{a_{i,j}\}_{i,j}$ in each direction to completely determine a function in $(A_k)^d|_f$. Thus, we have to find $|\{a_{i,j}\}_{i,j}| = d(l_1+1)^{d-1}$ coefficients. Any function of $xB_k$ reads however

$$p \in B_k \Rightarrow p = \bigtimes_{j=1}^{d} \left( x_j \sum_{|\beta_i| \leq l_2} b_i x^{\beta_i} \right)$$

for a given set of multi - indices $\{\beta_i\}_i$ and coefficients $\{b_i\}_i$independent of the coordinates $x_j$. There, it is enough to prescribe $(l_2+1)^{d-1}$ coefficients to determine the polynomial $\sum_{|\beta_i| \leq l_2} b_i x^{\beta_i}$. As $p$ is then constructed directly be shifting every component by the corresponding variable $x_j$, $|\{b_i\}|_i = (l_2+1)^{d-1}$ coefficients are enough to determine a function of $B_k|_f$. Thus, we get a first reduction of the required number of degrees of freedom from the general setting (67).

However, one has to be careful. Indeed, when restricted to the boundary of $K$ we are not in a direct sum anymore, it may happen that a function in $\mathbb{H}_k(K)|_f$ can be determined by only $d(l_1+1)^{d-1} + (l_2+1)^{d-1} - l_1^{d-1}$ monomial coefficients. To see it, we may derive the following. If $l_2 \geq l_1$, any polynomial

$p$ in $\mathbb{H}_k(K)|_f$ can be written under the form

$$p = \bigtimes_{j=1}^{d} \left( \sum_{|\alpha_i| \leq l_1} a_{ij} x^{\alpha_i} \right) + \bigtimes_{j=1}^{d} \left( x_j \sum_{|\beta_i| \leq l_2} b_i x^{\beta_i} \right) \tag{68}$$

$$= \bigtimes_{j=1}^{d} \left( \sum_{\substack{|\alpha_i| \leq l_1 \\ |\alpha_i| \neq 0}} (a_{ij} + b_{\xi_j(i)}) x^{\alpha_i} + x_j \sum_{l_1 \leq |\beta_i| \leq l_2} b_i x^{\beta_i} \right)$$
$$+ \bigtimes_{j=1}^{d} \alpha_{0j} x_j^0, \tag{69}$$

while if $l_1 \geq l_2 + 1$ it can be decomposed as

$$p \in \mathbb{H}_k(K)|_f \Rightarrow p = \bigtimes_{j=1}^{d} \left( \sum_{\substack{|\alpha_i| \leq l_1 \\ |\alpha_i \neq 0|}} (a_{ij} + b_{\xi_j(i)}) x^{\alpha_i} + a_{0j} x_j^0 \right). \tag{70}$$

There, the operators $\{\xi_j\}_{j \in [\![1, d]\!]}$ denote some coordinates permutation that reuptakes the shifted term generated by $x$ in the construction of $xB_k$.

***Example.*** *Example of a permutation $\xi_i$ in the two dimensional case*
Let us consider $l_1 = 2$ and $l_2 = 1$. Then, the intersection $(\mathcal{H}_{l_1})^2 \cap \mathrm{x}\mathcal{H}_{l_2}$ is not empty. In particular, any second order term involves the terms

$$b_1 \begin{pmatrix} x \\ y \end{pmatrix} xy, \quad b_2 \begin{pmatrix} x \\ y \end{pmatrix} x, \quad \text{and} \quad b_3 \begin{pmatrix} x \\ y \end{pmatrix} y,$$

living in $\mathrm{x}B_k$ while involving the term

$$\begin{pmatrix} a_{2,1} x^2 y \\ a_{2,2} y^2 \end{pmatrix}$$

lying in $A_k$. Thus, the coefficient associated to the term

$$\begin{pmatrix} x^2 y \\ y^2 \end{pmatrix}$$

in the decomposition of any element of the space $\mathbb{H}_k(K)|_{\partial K}$ gathers terms living in $(A_k)^2$ and a term living in $\mathrm{x}B_k$ through the factorised expression

$$\begin{pmatrix} (a_{2,1} + b_1) x^2 y \\ (a_{2,1} + b_3) y^2 \end{pmatrix}.$$

The permutation over the indices of the coefficients $\{b_i\}_i$ can be witnessed. ♦

The relations (68) and (70) yield

$$\dim \mathbb{H}_k(K)|_{\partial K} = \begin{cases} d(l_1+1)^{d-1} - d + (l_2+1)^{d-1} - l_1^{d-1} + d & \text{if } l_2 \geq l_1 \\ d(l_1+1)^{d-1} & \text{otherwise ,} \end{cases}$$

from where comes the following proposition.

**Proposition 5.11** Necessary condition for building conformal elements

A necessary bound that one get on the orders $l_1$ and $l_2$ can be obtained by deriving the inequality

$$\underbrace{\dim \mathcal{N}}_{\text{Wished number of normal moments per face}} \leq \underbrace{\dim \mathbb{H}_k(K)|_{\partial K}.}_{\text{Available tuning coefficients on one face}}$$

It our case, it reduces to the two following cases.

If $l_2 \geq l_1$,

$$\begin{aligned} d(l_1+1)^{d-1} + (l_2+1)^{d-1} &\leq d(l_1+1)^{d-1} + (l_2+1)^{d-1} - l_1^{d-1} \\ \Leftrightarrow \qquad l_1^{d-1} &\leq 0 \end{aligned}$$

while otherwise

$$\begin{aligned} d(l_1+1)^{d-1} + (l_2+1)^{d-1} &\leq d(l_1+1)^{d-1} \\ \Leftrightarrow \qquad (l_2+1)^{d-1} &\leq 0 \qquad (\Leftrightarrow l_2 = -1). \end{aligned}$$

***Remark.*** This is actually a necessary condition to be able to set the degrees of freedom. Indeed, asking more degrees of freedom would result into an over determination of the polynomials on the faces. Fulfilling this criterion does not however guarantee that a definition of a conformal element will be possible. It will depend on the value of $l_2$. Note also that in general, not achieving the bound gives us the liberty of choice in the projection space. ▲

We can now transpose the determination process of those coefficients in the dual space and define relevant degrees of freedom. By definition of $\mathbb{H}_k(K)|_f$, we will consider a polynomial projection space. Furthermore, as we want to enforce $H(\text{div}, K)$ – conformity through the normal degrees of freedom only, we limit ourselves to the three following shapes of degrees of freedom.

**Available types of degrees of freedom 5.12**

1. The face integral of coordinate - wise normal components tested against polynomials. They read

$$q \mapsto \int_{f_i} q_{x_i} n_{ix}\, p \,\mathrm{d}x_i \qquad \text{or} \qquad q \mapsto \int_{f_i} q_{x_i} n_{ix}\, p \,\mathrm{d}\gamma(x)$$

for any $q \in \mathbb{H}_k(K)$, and some polynomial $p$. Those degrees of freedom are suited to determine the vector polynomials living in $\mathbb{H}_k(K)$ in each direction by testing the coefficients of each polynomial component against a set of polynomials $\{p\}$ in which the maximum degree is higher that the maximum degree of the components of $q$. Typically, functions of $A_k$ would be represented that way. The $H(\mathrm{div}, K)$ – conformity will also be enforced coordinate - wise, which is a bit stronger than the usual case.

2. The face integral of the projection of a function in $\mathbb{H}_k(K)$ onto the face normal, tested against polynomials. They read

$$q \mapsto \int_{f_i} q \cdot n_i\, p \,\mathrm{d}\gamma(x)$$

for any $q \in \mathbb{H}_k(K)$ and some polynomial $p$. There, a function of $\mathbb{H}_k(K)$ is represented globally and all the directional polynomial will see their coefficients assigned homogeneously. In other words, polynomial vectors that are determined that way are a simple copy of one polynomial in every direction. This suits to the vectorial functions in $B_k$ as they are generated from a single polynomial. The $H(\mathrm{div}, K)$ – conformity is enforced in its classical sense.

3. The pointwise values of the discretised quantity against the normal of the face. They read

$$q \mapsto q(x_{im}) \cdot n_i$$

for any $q \in \mathbb{H}_k(K)$ and a set of sampling points $\{x_{im}\}_m$ lying on the face $f_i$. Those degrees of freedom are useful when one wants to enforce the variations of the polynomial with respect to val-ues at particular points, or when one wants to tune their offset. There, the tuning is global but it would also be possible to consider those values coordinate - wise. In any case, the $H(\text{div}, K)$ – conformity is achieved pointwise.

***Note.***

• Working with pointwise normal values requires special care. Indeed, the distribution of the sampling points along the face directly impacts the shape and amplitude of the corresponding polynomials in the dual space (*see e.g. Figure 54*).

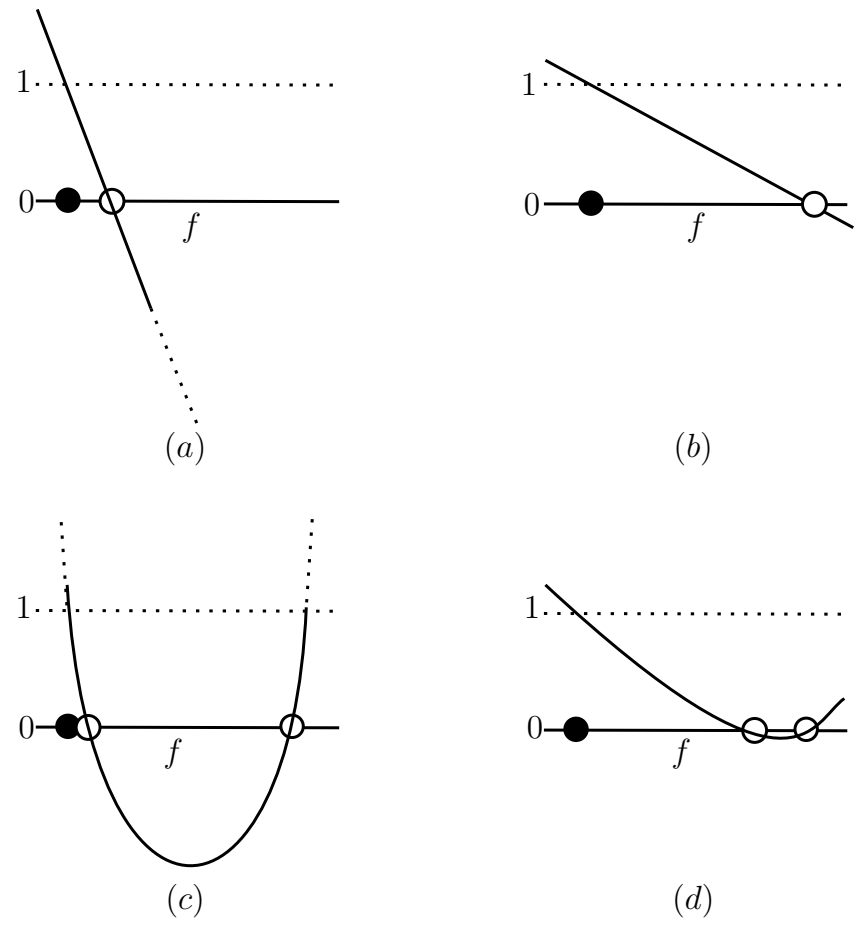

Fig. 54: (a) - (b): Amplitude impact. Using two close nodes increases the amplitude of the Lagrangian function, and thus the one of the basis function. (c) - (d): Shape impact. Distributing inhomogeneously the nodes may result in undesired polygons that are not suited to represent the solution with the same quality within the face. Even if the Lagrangian functions sum to one, the distribution of their behaviour is important.

In particular, using points that are close to each other weakens the discretisation reliability and performing computationally a change of basis may be delicate. Indeed, the prescribed duality relationship (5) reduces to ask the local Lagrangian relation on each face

$$p_j|_{f_i}(x_{im}) \cdot n_{ij} = \delta_{ij}.$$

for some vectorial functions $\{p_j\}$ in $\mathbb{H}_k(K)$. Thus, using close sampling points at a fixed polynomial order may end up to very high variations of the corresponding basis functions that are not recommended computationally thinking, arising in particular at the extremities of the faces (*see the Figure (54d)*). In addition, the conditioning number of the transfer matrix would be horrible. The distribution of the points must then be chosen carefully.

• When one uses only one pointwise value in the set of normal degrees of freedom and if no moments based degrees of freedom projects on constants, the above duality relation prescribes the global offset of the polynomials (*see the Figure 55*). Furthermore, when the pointwise values are considered coordinate wise, the tuning of the offset is also done coordinate wise.

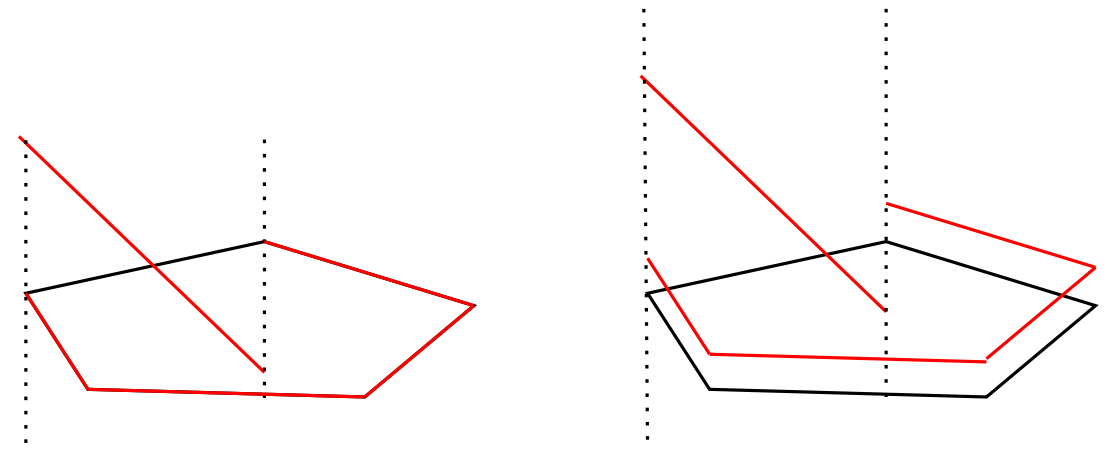

Fig. 55: Example of impact of the pointwise normal value on the offset of the determined polynomial when the space $\mathcal{H}_1$ is used. Left: non - shifted pointwise value. Right: Shifted pointwise value.

• Using several pointwise values prevents from prescribing the global offset as they would then describe higher variations of polynomials rather that their constant part. Furthermore, using several pointwise values in combination with moment based degrees of freedom makes the contribution of each degrees of freedom in the determination process hard to infer. ▲

Let us now present more specifically the sets of normal degrees of freedom that one may want to use. As we wish to keep the spirit of Raviart – Thomas, we will maximize the number of moment based degrees of freedom.

In our case, as the normal component of boundary functions are polynomials of degree $\max\{l_1, l_2\}$, we cannot define projections over polynomials whose degree is more than $\max\{l_1, l_2\}$ without breaking the unisolvence. Therefore, we can only select normal degrees of freedom among the following sets that are reduced from the types presented in the *Box 5.12*.

$$q \mapsto \int_f q_{x_i} n_{ix}\, p \,\mathrm{d}x_i, \quad \text{or} \quad q \mapsto \int\limits_{f_i} q_{x_i} n_{ix}\, p \,\mathrm{d}\gamma(x), \quad p \in \mathbb{Q}_{\max\{l_1, l_2\}}(f) \tag{71}$$

$$q \mapsto \int_f q \cdot n\, p \,\mathrm{d}\gamma(x), \quad p \in \mathbb{Q}_{\max\{l_1, l_2\}}(f) \tag{72}$$

$$q \mapsto \int_f q \cdot n\, p \,\mathrm{d}\gamma(x), \quad p \in \bigcup_{i=1}^{d-1} \{x_i\, q_i\},\ q_i \in \mathbb{Q}_{\zeta_i([l_2], l_2, \cdots, l_2)}(f) \tag{73}$$

$$q \mapsto q(x_{im}) \cdot n_i, \ \text{ for a set of sampling points } \{x_{im}\}_m \text{ on the face } f_i. \tag{74}$$

Note that in the case $l_2 < l_1$, the third set is merged with the second one.

One then just has to extract the number of required degrees of freedom $d(l_1 + 1)^{d-1} + (l_2 + 1)^{d-1}$ taking care that the integrands $p$ are free from each others. There, it reduces to choosing the multi - indices $\alpha$ for which the projection on the monomial $x^\alpha$ will be component wise and the other ones for which the projection will be done homogeneously thorough the components.

***Remark 6.***

• Once this decision is done for a each multi - index $\alpha$, one can go a bit weaker and accept a mixture in the nature of the degrees of freedom for given $\alpha$ where we decided a coordinate wise determination. We can by example set $(d - 1)$ component wise degrees of freedom and take the last one as the global one, or more generally take $c$ component wise, $e$ pointwise and possibly a global one so that $e + c + 1 = d$. However this complexifies the definition of the set of degrees of freedom for – up to the authors knowledge – no particular gain.

• One could also set one global degrees of freedom and only $c < d - 1$ component wise. This would tune identically $d - c$ components and provide an independent tuning for the other $c$ ones. However, this prevents the discretisation to be legible and will be disregarded.

• More global components are designed, more normal components of basis functions degenerate into internal basis functions in the dual space. Indeed, on the boundary the order of the polynomial is $\max l_1,\, l_2$. Therefore, the term $p \cdot n|_f$ requires only $\max l_1,\, l_2 + 1$ basis functions to be decomposed on. Thus, all the other basis functions required to represent totally the space $\mathbb{H}_k(K)|_f$ may see their global normal component vanish. Their component

- wise normal basis functions are however not vanishing as they are taking care of the coordinate - wise behaviours that are invisible through the global expression of $p \cdot n|_f$.

Therefore, more global degrees of freedom are designed, more the representation of $p \cdot n$ is completed globally, forcing the other basis functions to see their global normal component vanishing. Reversely, if the focus is made on the component - wise behaviour $p_{x_i} n_{ix}$, the expression of $p \cdot n$ will be decomposed on the coordinates and less basis functions would be forced to have a vanishing global component. One then need to think of the type of basis function he wants to design while selecting the degrees of freedom. ▲

***Remark.*** In all the following, we denote by $\int_{\partial K} \cdot \,\mathrm{d}x$ the group of moments $\int_f \cdot \,\mathrm{d}x$ for all faces $f \in \partial K$. ▲

At this point, any admissible definition would lead to $H(\mathrm{div},\, K)$ – conformal elements. Let us present some examples of sets of degrees of freedom and derive an admissibility condition opening further possibilities.

### 5.3.3 A possible selection of normal degrees of freedom

We detail here four possibilities as examples for the case $l_2 \geq l_1$ where every function in $\mathbb{H}_k(K)|_f$ reads under the form (68). The case $l_2 < l_1$ is more restrictive and will only be shortly evoked.

#### 5.3.3.1 First approach

A first natural definition emerging when reading out the structure of (68) provides the following degrees of freedom.

$$\sigma\colon p \longmapsto p_{x_i}(x_{im}) n_{ix_j}, \qquad \substack{j \in [\![1,\, d]\!],\, i \in [\![1,\, n]\!] \text{ and} \\ x_{im} \text{ middle or barycentric point in } f_i} \tag{75a}$$

$$\sigma\colon p \longmapsto \int\limits_{\partial K} p_{x_i}\, n_{x_i}\, x^{\alpha_j}\, \mathrm{d}x_i, \qquad \substack{i \in [\![1,\, d]\!], \\ \text{for any mutli - index } 1 \leq |\alpha_j| \leq l_1} \tag{75b}$$

$$\sigma\colon p \longmapsto \int\limits_{\partial K} p \cdot n x_i x^{\alpha_j}\, \mathrm{d}\gamma(x), \quad l_1 \leq |\alpha_j| \leq l_2,\, \forall i \in [\![1,\, d-1]\!] \tag{75c}$$

There, the value of the projection against the monomials furnishes the coefficients of each monomial part of the tested quantity for any coordinate, and the pointwise values drive the offset of polynomials in every direction.

***Remark 7.***

• Note that by construction, this configuration does not allow the use of more than $d$ degrees of freedom involving only constants, $d+1$ degrees of freedom involving a same monomial for those that belong both to $\mathcal{H}_{l_1}$ and $\mathcal{H}_{l_2}$, and one degree of freedom per monomial involved in the space $\mathcal{H}_{l_2} \setminus \mathcal{H}_{l_1}$, either as an integrand or as point values. Therefore, there is no over-determination possible when tuning the coefficients $\{\{a_{im}\}_{im}, \{b_i\}_i\}$ at fixed $i$.

• In practice and for unisolvence reasons, one may rather use any vector $v \neq n$ in place of the normal $n$ in the moment (75c). Indeed, for any $\alpha_i$ such that $|\alpha_i| = l_1$, the moment (75c) can be reconstructed from the degrees of freedom (75b). The order of any polynomial $p$ is reduced by one when tested against an edge's normal. Thus, the term $p \cdot n\, x_i$ reconstructs a behaviour equivalent to the $i^{\text{th}}$ coordinate of $p$, and the global moments (75c) when $\alpha_i$ is such that $|\alpha_i| = l_1$ duplicates the information that (75a) provides coordinate - wise. Therefore, in order to characterise completely the elements living in $xB_k$ in this layout, we need to define degrees of freedom that do not involve a loose of polynomial order. Using a vector different from the normal is a simple way. Even though this would not belong to the type of admissible degrees of freedom we described before, the $H(\text{div}, K)$ – conformity is preserved (though not enforced). Indeed, it reduces here to the continuity of the normal component, and on the boundary the order of the normal component is lowered. Furthermore, the lower order parts of the normal component is fully described by (75a), (75b) and by (75c) for any $\alpha_j$ such that $\alpha_j > l_1$. ▲

The full set $\Sigma :=$(75a) – (75c) has the following property.

**Property 5.13** Dimension

$$\dim(\Sigma) = n[d(l_1+1)^{d-1} + (d-1)(l_2+1)^{d-1} - (d-1)l_1^{d-1}]$$

***Proof.*** By construction, we directly obtain that on every face $f_i \in \partial K$,

$$\begin{aligned}
\dim\{(75a)\}|_{f_i} &= d \\
\dim\{(75b)\}|_{f_i} &= d((l_1+1)^{d-1} - (0+1)^{d-1}) \\
\dim\{(75c)\}|_{f_i} &= (d-1)((l_2+1)^{d-1} - (l_1-1+1)^{d-1}) \\
&= (d-1)((l_2+1)^{d-1} - l_1^{d-1}).
\end{aligned}$$

Indeed, the number of multi - indices $\alpha_j$ such that $1 \leq |\alpha_j| \leq l_1$ is

$|\{\alpha_j\}_j| = (l_1+1)^{d-1} - (0+1)^{d-1}$. Thus, it comes

$$\begin{aligned}\dim(\Sigma) &= d + d(l_1+1)^{d-1} - d + (d-1)((l_2+1)^{d-1} - l_1^{d-1}) \\ &= d(l_1+1)^{d-1} + (d-1)(l_2+1)^{d-1} - (d-1)l_1^{d-1}.\end{aligned}$$

As it holds similarly for every face, the claim follows.

■

We see immediately a criterion on the relationship between $l_1$ and $l_2$. Indeed, we need to set $d(l_1+1)^{d-1} + (l_2+1)^{d-1}$ degrees of freedom on each face. In this configuration, we can only set a maximum of $d(l_1+1)^{d-1} + (d-1)(l_2+1)^{d-1} - (d-1)l_1^{d-1}$ on each face. A straightforward feasibility condition then reads

$$\begin{aligned}&d(l_1+1)^{d-1} + (l_2+1)^{d-1} \leq d(l_1+1)^{d-1} + (d-1)(l_2+1)^{d-1} - (d-1)l_1^{d-1} \\ &\Leftrightarrow (1-d+1)(l_2+1)^{d-1} \leq -(d-1)l_{1,}^{d-1} \qquad (76)\end{aligned}$$

finally giving

$$\boldsymbol{(d-1)l_1^{d-1} \leq (d-2)(l_2+1)^{(d-1)}_.} \qquad (77)$$

***Remark.***
• By the condition 5.11 one already requires $l_1 = 0$. Therefore, in all dimensions bigger than two the above criterion is always fulfilled.
• In dimension two, the above relation turns to be an equality for any $l_2$. ▲

**A possible selection** As the condition (77) is always fulfilled, we can design the elements by considering the set $(\Sigma)$ and the relation (76). Indeed, we can set the degrees of freedom by considering fully the sets (75a) and (75b) while extracting $(l_2+1)^{d-1}$ elements from (75c). To select them, one can simply consider to extract them starting from the order $l_1$ and increasing with respect to the lexicographical order, though any other definition is equally admissible. We will denote this construction $Ia$.

**Assumption 5.14** Extraction of the $(l_2+1)^{d-1}$ degrees of freedom

Connecting with the third point of the *Remark 7*, we assume that when extracting the $(l_2+1)^{d-1}$ degrees of freedom, all the integrands of the moments and the polynomials $q \mapsto \sigma(q(x_{im}))$ that define the point values are chosen linearly independent.

***Remark.*** In the two dimensional case, the condition (77) turns to be an equality. Indeed, for $d = 2$, we have

$$(d-1)l_1^{d-1} \leq (d-2)(l_2+1)^{(d-1)} \Leftrightarrow l_1^1 \leq 0,$$

which together with the previous remark leads to a null identity. Therefore, we have as many liberties as possibilities to tune the polynomials, and thus have to consider the set $(\Sigma)$ in its totality. Its dimension will then match $2(0+1)+(l_2+1)-0 = l_2+3$, which corresponds to expected number of degrees of freedom on every face of the boundary. ▲

***Remark.*** An alternative to this definition keeping the same outer structure is to prescribe projections on constants instead of using moments in (75a). It would then read

$$p \mapsto \int_{f_i} p_{x_j} n_{ix_j} \,\mathrm{d}x_j, \quad j \in [\![1,\, d]\!],\; i \in [\![1,\, n]\!], \tag{78}$$

leading to a configuration that we denote by *Ib*. ▲

To conclude the first possibility of defining the normal degrees of freedom, we can derive the following property ensuring that the elements are well defined.

**Proposition 5.15** Unisolvence

Let us assume that the internal degrees of freedomare properly defined. Then, the two sets of normal degrees of freedom $\{\sigma\} = \{(75a), (75b), (75c)\}$ and $\{\sigma\} = \{(78), (75b), (75c)\}$ combined with the internal degrees of freedom are unisolvent for $\mathbb{H}_k(K)$ and define $H(\mathrm{div},\, K)$ – conformal elements.

***Proof.*** The proof is immediate given the *Proposition 5.21* and is succinctly given in the *Proof 5.3.6.*

■

#### 5.3.3.2 An alternative approach

Another definition is possible when observing the following relation.

$$d(l_1+1)^{d-1}+(l_2+1)^{d-1} = (d-1)(l_1+1)^{d-1}+(l_2+1)^{d-1}+(l_1+1)^{d-1} \tag{79}$$

There, keeping in mind that $\dim \mathbb{Q}_{l_2} = (l_2+1)^{d-1}$ and $\dim \mathbb{Q}_{l_1} = (l_1+1)^{d-1}$, one can define the first normal degrees of freedom by using the expressions of the two first types of degrees of freedom described in the *Paragraph 5.12.* Indeed, taking into consideration that there exists non - constants polynomial in any space $\mathbb{H}_k(K)$ for any order, we can split the first order projection space coordinate - wise and derive

$$\sigma\colon q \mapsto \int_{\partial K} q_{x_i}\, n_{x_i} p_k \,\mathrm{d}x_i, \qquad \text{for all } i \in [\![1,\, d]\!],\ \text{ all } p_k \in \mathcal{H}_{l_1} \setminus \mathcal{H}_0(\partial K) \tag{80a}$$

$$\sigma\colon q \mapsto \int_{\partial K} q_{x_i}\, n_{x_i}\, x_i^{l_1+1} \,\mathrm{d}x_i,\ \text{ for all } i \in [\![1,\, d]\!]. \tag{80b}$$

Note that the second line represents the component - wise projections of elements living in $xB_k$. Then, to retrieve the higher degree representations, we can use

$$\sigma(q) \mapsto \int_{\partial K} q \cdot n\, p_k \,\mathrm{d}\gamma(x), \qquad \text{for all } p_k \in \mathcal{H}_{l_2}(\partial K) \setminus \mathcal{H}_{l_1}(\partial K). \tag{80c}$$

Note that for all $m \in [\![l_1,\, l_2]\!]$ and any $i \in [\![1,\, d]\!]$, the monomial $x \mapsto x_i x^{\alpha_i}$ such that $|\alpha_i| = m$ belongs to $\mathcal{H}_{l_2}(\partial K)$. Furthermore, for all monomial $p \in \mathcal{H}_{l_2}(\partial K) \setminus \mathcal{H}_{l_1}(\partial K)$, there exists a $j \in [\![1,\, d]\!]$ and a $\alpha_i$ such that $|\alpha_i| \leq l_2 - 1$ such that $p = x_j x^{\alpha_i}$. Thus, no projection on (80c) will be null and those degrees of freedom are admissible. There, the lowest monomials are tuned coordinate wise while the highest ones are tuned globally.

As we can notice, the sets (80a) and (80c) together form a set of $d(l_1 + 1)^{d-1} + (l_2+1)^{d-1} - (l_1+1)^{d-1}$ degrees of freedom. We then have to choose $(l_1+1)^{d-1}$ other degrees of freedom. Reading out the nature of those two sets, it turns out that the only vectorial monomials lying in $\mathbb{H}_k(K)|_{\partial K}$ that are not described yet are of the form

$$x \mapsto (x_1 x_1^{l_2}\tilde{x},\, \cdots,\, x_d x_d^{l_2}\tilde{x})_j^T \in \bigtimes_{i=1}^{d} x_i \mathbb{P}_{\zeta_i([l_2],\, l_2,\, \cdots,\, l_2)} \tag{81}$$

with $\tilde{x}$ a monomial defining from the $(d-2)$ other variables whose degree is less than $l_2$. Those monomials correspond to the second term in the right hand side of the dimension relation (68), restricted to the boundary.

Therefore, any normal component can be tested against any of the coordinates of the vector. As we are on the boundary, the interdependency of the variables by the face parametrization implies similarities between some components of the vector (81) and lowers the number of degrees of freedom that we can set by one dimension. Thus, we can set at most

$$(d-1)(l_2+1)^{d-2}$$

degrees of freedom. The relationship between $l_1$ and $l_2$ then reads

$$(l_1+1)^{d-1} \leq (d-1)(l_2+1)^{d-2}. \tag{82}$$

***Remark.*** As in the previous case, for any $d \geq 2$ this condition is always fulfilled by the condition on $l_1$ itself. ▲

We can now set the left $(l_1+1)^{d-1}$ degrees of freedom either as moment based degrees of freedom projecting onto the monomial vectors (81) or by prescribing pointwise values of normal components. Those two categories respectively read

$$q \mapsto \int_{\partial K} q \cdot n x_j x_j^{l_2} \tilde{x} \, \mathrm{d}\gamma(x), \qquad \text{for any } \tilde{x} \in \mathbb{Q}_{l_2}(\partial_j K) \tag{83}$$

$$\text{and} \qquad q \mapsto q(x_{im}) \cdot n, \quad \substack{\forall i \in [\![1, d]\!],\, m \in [\![1, (1+l_1)^{d-1}]\!] \\ \text{and } x_{im} \text{ a sampling point}} \tag{84}$$

where $\partial_j K \subset \mathbb{R}^{d-2}$ emphases the parametrisation of any face $f \in \partial K$ from the variable $x_j$.

Selecting two by two different degrees of freedom among the above sets is enough to ensure that the corresponding element will be admissible. Mixing the type of degrees of freedom is also possible though not recommended for the sake of legibility.

Note that those degrees of freedom only describe functions of $\mathbb{H}_k(K)|_{\partial K}$ when their projection onto $x\mathbb{Q}_{l_2} \backslash \mathbb{Q}_{l_1}$ is not null. They are therefore dedicated to represent high variations. As in those cases the functions have a uniform polynomial definition as a core, the degrees of freedom are only tested globally without discrimination between the coordinates.

***Remark.***

• As before, in dimension two the inequality (82) reduces to an equality and when sticking to the same type of degrees of freedom we have to consider the full set (83) or (84).

• In this configuration, no identical monomial is shared between the different types of used degrees of freedom (80) – (83)/(84). Contrarily to the previous construction, one is therefore free to choose any subset of degrees of freedom he likes from any of the above sets without further assumption. ▲

***Note.***

• The liberty expressed by those degrees of freedom comes from the functions in $\mathbb{H}_k(K)|_{\partial K}$ whose projection is not null. Since $l_2 \geq l_1$ we have a direct sum and can state that they are are part of $B_k$.

• In the quadrilateral case, the lack of determination for high order polynomials was not described on the faces. Indeed, as it was not a block construction the determination was global and the freedom was reuptaken by the internal moments, preventing any supplementary moment to be required on the boundary. ▲

In both constructions, if either

$$(l_1+1)^{d-1} < (d-1)(l_2+1)^{d-2}$$

$$\text{or}$$

$$(d-1)l_1^{d-1} < (d-2)(l_2+1)^{d-1}$$

then the definition of $E_k(K)$ becomes sensitive to the choice of degrees of freedom if they are chosen as moments rather than point values. Indeed, the choice of the projection space for the $(l_1+1)^{d-1}$ components will define the degrees of freedom and tune the basis functions towards specific behaviours, and therefore emphasize specific properties.

As before, we can derive the following property ensuring that the elements are well defined.

**Proposition 5.16** Unisolvence

Let us assume that the internal degrees of freedom are properly defined. Then, the two sets of normal degrees of freedom $\{\sigma\} = \{(80a), (80c), (83)\}$ and $\{\sigma\} = \{(80a), (80c), (84)\}$ combined with the internal degrees of freedom are unisolvent for $\mathbb{H}_k(K)$ and define $H(\text{div}, K)$ – conformal elements.

***Proof.*** The proof is immediate given the *Proposition 5.21* and is succinctly given in the *Proof 5.3.6*

■

#### 5.3.3.3 Further definitions

The case discussed here is only an example the moments' selection. Any other extraction of the right number of degrees of freedom among the presented options (5.12) fulfilling the condition 5.14 would be admissible.

In particular, in dimension strictly bigger than two one can design elements in which no coordinate - wise degrees of freedom are used. Indeed, connecting with the *Remark 7* to tune $d+1$ coefficients associated to one monomial order whom $d$ are coordinate-wise and one is global, one can either define $d$ coordinate wise and 1 global, $d-1$ coordinate - wise and 2 global, or any other linearly independent combination of $d+1$ relations involving a same polynomial order. Setting even $d+1$ global coordinates would be admissible. Furthermore, as we work on the $\mathbb{Q}_k$ spaces, in dimension bigger than two there are enough moments to design $d+1$ global relations on a same polynomial order.

The case $l_1 > l_2$ is much more restrictive. Indeed, all the degrees of freedom have to be split coordinate wise in (80a), (80c) and we are left with $(l_2+1)^{d-1}$ degrees of freedom to set. However, every function is already fully characterized on the boundary. This is impossible unless $l_2 = -1$, which copes with the condition (5.11). In this case, any of the definitions presented in the above case can be used.

#### 5.3.3.4 Summary of the investigated configurations

Let us summarize the admissible constructions of degrees of freedom that were built in the previous paragraph. We could design four main types of $H(\text{div}, K)$ – conformal elements depending on the types of desired degrees of freedom and the layout of their distribution.

In all cases, the sets $(a)$ and $(b)$ have to be considered fully. Then, for the first configuration, only $(l_2+1)^{d-1}$ degrees of freedom have to be selected in the subset $(c)$ while for the second one $(l_1+1)^{d-1}$ degrees of freedom needs to be picked up. The selection is a matter of taste, possibly directed by properties of the discretised quantities that are known *a - priori*. For the sake of concision we will only summarize the four types of degrees of freedom that are natural in their construction. According to the above discussion, they are split into two main configurations $I$ and $II$, within which only slight changes occurs in the definition of degrees of freedom.

**First configuration, denoted by** $I$ By the relation (77), this configuration is available whenever $l_1 \leq 0$ and $(d-1)l_1^{d-1} \leq (d-2)(l_2+1)^{d-1}$.

**Element 5.17** Element $Ia$

$$\boldsymbol{p \mapsto \int_{f_i} p_{x_j} n_{ix_j}\, \mathrm{d}x_j,} \qquad \boldsymbol{j \in [\![1,\, d]\!],\, i \in [\![1,\, n]\!]} \tag{85a}$$

$$p \mapsto \int_{\partial K} p_{x_i}\, n_{x_i}\, x^{\alpha_j}\, \mathrm{d}x_i, \qquad \text{for any } 1 \leq |\alpha_j| \leq l_1, i \in [\![1,\, d]\!] \tag{85b}$$

$$p \mapsto \int_{\partial K} p(x)\cdot v\, x^{\alpha_j}\, \mathrm{d}\gamma(x), \qquad |\alpha_j| = l_1 + 1 \tag{85c}$$

$$p \mapsto \int_{\partial K} p(x)\cdot n\, x_i\, x^{\alpha_j}\, \mathrm{d}\gamma(x), \quad l_1 < |\alpha_j| \leq l_2,\, \forall j \in [\![1,\, d-1]\!], \tag{85d}$$

where for each face $f_i$, $v \neq n_i$.

**Element 5.18** Element $Ib$:

$$\boldsymbol{p \mapsto p_{x_i}(x_{im}) n_{ix_j},} \qquad \underset{\boldsymbol{x_{im}} \text{ middle or barycentric point on } \boldsymbol{f_i}}{\boldsymbol{j \in [\![1,\, d]\!],\, i \in [\![1,\, n]\!]} \text{ and}} \tag{86a}$$

$$p \mapsto \int_{\partial K} p_{x_i}\, n_{x_i}\, x^{\alpha_j}\, \mathrm{d}x_i, \qquad \text{for any } 1 \leq |\alpha_j| \leq l_1, i \in [\![1,\, d]\!] \tag{86b}$$

$$p \mapsto \int_{\partial K} p(x)\cdot v\, x^{\alpha_j}\, \mathrm{d}\gamma(x), \qquad |\alpha_j| = l_1 + 1 \tag{86c}$$

$$p \mapsto \int_{\partial K} p(x)\cdot n\, x_i\, x^{\alpha_j}\, \mathrm{d}\gamma(x), \quad l_1 < |\alpha_j| \leq l_2,\, \forall j \in [\![1,\, d-1]\!] \tag{86d}$$

where for each face $f_i$, $v \neq n_i$.

**Second configuration, denoted by** $II$ By the relation (82), this configuration is available when $l_1 \leq 0$ and $(l_1+1)^{d-1} < (d-1)(l_2+1)^{d-2}$.

**Element 5.19** Element $IIa$:

$$q \mapsto \int_{\partial K} q_{x_i}\, n_{x_i} p_k \,\mathrm{d}x_i, \qquad \substack{\text{for all } i \in [\![1,d]\!],\\ \text{for all } p_k \in \mathcal{H}_{l_1}(\partial K) \setminus \mathcal{H}_0(\partial K)} \tag{87a}$$

$$q \mapsto \int_{f_j} q_{x_i}\, n_{x_i} x_i^{l_1+1} \,\mathrm{d}x_i, \qquad \text{for all } i \in [\![1,\, d]\!] \text{ and all } j \in [\![1,\, n]\!] \tag{87b}$$

$$q \mapsto \int_{\partial K} q \cdot n\, p_k \,\mathrm{d}\gamma(x), \qquad \text{for all } p_k \in \mathcal{H}_{l_2}(\partial K) \setminus \mathcal{H}_{l_1}(\partial K) \tag{87c}$$

$$\boldsymbol{q \mapsto \int_{\partial K} q \cdot n x_j x_j^{l_2} \tilde{x} \,\mathrm{d}\gamma(x), \qquad \text{for any } \tilde{x} \in \mathbb{Q}_{l_2}(\partial_j K)} \tag{87d}$$

**Element 5.20** Element $IIb$:

$$q \mapsto \int_{\partial K} q_{x_i}\, n_{x_i} p_k \,\mathrm{d}x_i, \qquad \substack{\text{for all } i \in [\![1,d]\!],\\ \text{for all } p_k \in \mathcal{H}_{l_1}(\partial K) \setminus \mathcal{H}_0(\partial K)} \tag{88a}$$

$$q \mapsto \int_{f_j} q_{x_i}\, n_{x_i} x_i^{l_1+1} \,\mathrm{d}x_i, \qquad \text{for all } i \in [\![1,\, d]\!] \text{ and all } j \in [\![1,\, n]\!] \tag{88b}$$

$$q \mapsto \int_{\partial K} q \cdot n\, p_k \,\mathrm{d}\gamma(x), \qquad \text{for all } p_k \in \mathcal{H}_{l_2}(\partial K) \setminus \mathcal{H}_{l_1}(\partial K) \tag{88c}$$

$$\boldsymbol{q \mapsto q(x_{im}) \cdot n, \qquad \substack{\forall i \in [\![1,d]\!],\, m \in [\![1,\,(1+l_1)^{d-1}]\!]\\ \text{and } x_{im} \text{ a sampling point}}} \tag{88d}$$

***Remark 8.***

• The constraint (87d) projecting functions of $\mathbb{H}_k(K)|_{\partial K}$ onto high order polynomials tunes automatically the constant part of the tested quantities that are not explicitly considered in the second configuration. Indeed, on the boundary any element of $\mathbb{H}_k(K)$ read $q = (x, y)^T p_{l_2} + (1, 0)^T p_{l_{1,1}} + (0, 1)^T p_{l_{1,2}}$ for some polynomials $p_{l_2} \in \mathbb{Q}_{l_2}(\partial K)$ and $p_{l_{1,1}}, p_{l_{1,2}} \in \mathbb{Q}_{l_1}(\partial K)$. Thus, the projection against (87d) furnishes the value of $\int_{f_i}(Cp_{l_2} + n_{ix}p_{l_{1,1}} + n_{iy}p_{l_{1,2}})x_j^{l_2+1}\tilde{x}\mathrm{d}x = \int_{f_i}(C\tilde{p_{l_2}} + n_{ix}\tilde{p_{l_{1,1}}} + n_{iy}\tilde{p_{l_{1,2}}} + C + C_1 n_{ix} + C_2 n_{iy})x_j^{l_2+1}\tilde{x}\mathrm{d}x$, where $C = x \cdot n$ and the tilde over the polynomials represents the polynomial from which the constant part has been taken out. Therefore, as the polynomial order of the normal component is lowered on the boundary, this value can be combined with the degrees of freedom (5.19) – (87c), and therefore determines the constant part (remember that the monomials constants corresponding to the orders comprised between $l_1 + 1$ and $l_2 + 1$ are identical for all of the coordinates of the tested quantity). However, the constant value is not accessible by a direct lecture of the degrees of freedom.

• For the same reason, and by example in two dimensions where only one such degrees of freedom has to be designed, one can simply consider the global $\int_{\partial K} p \cdot n$ moment instead. The element remains unchanged. The implementation used for obtaining the results of the *Section 7* uses this equivalent formulation.

• In two dimension, every set has to be considered fully.

• When one distribution layout is selected, one is free to select the desired degrees of freedom within the given framework. ▲

The configuration (5.18) is the most intuitive one and allows to prescribe directly the offset of the polynomials in each direction. However, it is not fully moment based. The version (5.17) tackles this issue but we loose the direct setting of the offset.

The first type of the second configuration (5.19) is fully moment based. The functions are determined blindly only through average values on the faces, and is thus more suited to represent flux quantities. However, due to possibly high order of the projection polynomials, the coefficients of the inverse of the transfer matrix may be huge and scaling of the tuned basis functions may be required to obtain intelligible results. Furthermore, its definition may be counterintuitive.

The last presented type (5.20) gathers the inconvenient of the counter-intuitive definition and the delicate usage of pointwise normal values which here does not allow to tune the offset. It is then *a - priori* not advised.

***Note.*** In all of the above definitions, instead of projecting against the monomials $p_j$ one could project on any polynomial basis of $\text{span}\{p_j\}_j$. ▲

***Remark.*** We point out some last remarks.

• Note that for the normal degrees of freedom , it comes down to the same thing to consider projections onto Poisson's functions in $\mathbb{H}_k(K)$ or just the corresponding polynomials serving as boundary conditions as those spaces merges on the boundary.

• We saw that we had to set $l_1 \leq 0$ for any order and any dimension. If it may seem bothersome, as it corresponds to the part which is not conserved through the divergence it is actually welcomed that it vanishes.

• The mandatory $l_1 \leq 0$ also has a further advantage. Indeed, the two only types of degrees of freedom that enforces conformity are midpoint normal values, and testing the normal component against functions, and midpoints are delicate to handle. Only one shifts the behaviour, but having multiple prescribed losses in legibility *a - priori*. Furthermore, the configuration *IIb* only requires one pointwise value per face, and directly tunes the global offset.

• Any value of $l_2$ is allowed for any dimension. Therefore, the desired refinement sequence induced by $l_2 = a_2 k + b_2$ is admissible as any natural number $k$ would generate an admissible element.

• The sets of moments $(a)$ and $(b)$ of every configuration may also integrate with respect to the path $\gamma(x)$ rather than the coordinates $x_i$, as evoked in the first point of the *Box 5.12*.

• Lastly, note that the functions that are dual to those normal degrees of freedom do not vanish on every face. Indeed, by the unisolvence property they would be reduced to the identically null function. ▲

### 5.3.4 Definition of admissible Internal degrees of freedom

The definition of the internal degrees of freedom requires a little bit more work. Indeed, we have to make sure that the corresponding internal basis functions vanishes on every face. Given a generic basis vanishing on the faces, it is enough to ensure that the definition of degrees of freedom preserves it through the change of basis. We therefore stick to the idea of Raviart – Thomas and define moment based degrees of freedom that reads for any

$q \in \mathbb{H}_k(K)$;

$$\sigma(q) \mapsto \int_K q \cdot p_k \, \mathrm{d}x, \qquad \text{for all } p_k \in \mathcal{P}_k.$$

It is then left to define the space $\mathcal{P}_k$ of dimension $\left((m_1+1)^d + (m_2+1)^d\right)$ on which no projection of functions in $\mathbb{H}_k(K)$ is identically null or identically identical. Furthermore, it should not allow the projection on polynomials of degree at least $k+1$ when the order of the space $\mathbb{H}_k(K)$ is $k$.

As a test space, the space $\mathcal{P}_k$ can then simply gather functions if we want to emphasize the polynomial approximation of the discretised quantity. One may rather consider testing against Poisson's functions as in the definition of the space to stay coherent with the definition of the quantities' approximation.

### 5.3.5 Admissible definitions of $\mathcal{P}_k$

Depending on the space $\mathbb{H}_k(K)$ one would like to work with, the set $\mathcal{P}_k$ will enjoy different prescribed dimensions. We then have to design any projection space coping with any desired dimension and furnishing a decent quantification of the discretised quantity.

For a general polytope with $n$ faces and for any order $k \geq 0$, defining $\mathcal{P}_k$ is not a straightforward task. We detail some possible constructions preserving the $H(\mathrm{div},\, K)$ – conformity.

#### 5.3.5.1 The role of $\mathcal{P}_k$

As one can see in the *Section 2.3*, the set $\mathcal{P}_k$ is used to generate the set of internal moments of the element $E_k(K)$. Therefore, its definition plays a crucial role in the construction of the element and in the subsequent discretisation properties. In particular, it will highlight some specificities of the discretised quantity in the inner cell while totally neglecting others. The choice of $\mathcal{P}_k$ has thus to be done carefully.

One can notice that any function $p_k \in \mathcal{P}_k$ can be interpreted as a kernel for the moment $q \mapsto \int_K q \cdot p_k \, \mathrm{d}x$. Reversely, any moment can be expressed from a function belonging to $\mathcal{P}_k$. Therefore, the dimension of $\mathcal{P}_k$ should match the one of the internal moments' set, that is $\dim \mathcal{P}_k = \dim(\mathbb{H}_k(K) \setminus \mathbb{H}_k(K)|_{\partial K}$. As we enjoy a block construction, it reduces to $\dim \mathcal{P}_k = \dim \mathbb{H}_k(K) - \dim \mathbb{H}_k(K)|_{\partial K} = d(m_1+1)^d + (m_2+1)^d$. A general definition of $\mathcal{P}_k$ is then hard to set up. However, we can use the properties that have to enjoy the internal moments to construct $\mathcal{P}_k$ on a general spirit.

Meaningly, we know that the set of internal degrees of freedom should be defined such that the total set of degrees of freedom $\{\sigma\}$ is unisolvent for $\mathbb{H}_k(K)$, and that each function of the set should be a linear forms $K$.

***Remark 9.*** To ensure the unisolvence, it is enough to consider $\mathcal{P}_k$ as a free subspace of $\mathbb{H}_k(K)$. Indeed, the projections to this space would never be null or identical for every functions lying within $\mathbb{H}_k(K)$. Another way is to define $\mathcal{P}_k$ as a polynomial space. Indeed, for any $u \in \mathbb{H}_k(K)$, we have $\Delta u \in \mathbb{Q}_{m_1}$ or $\Delta u = 2\nabla \cdot (u) + x\Delta u$ with $\Delta u \in \mathbb{Q}_{m_2}$. Therefore, the projection of its Laplacian to $\mathbb{Q}_{\max\{m_1,\, m_2+1\}}$ is never null neither identical for every function of $\mathbb{H}_k(K)$, and we get non - null projection of the function itself onto polynomial of degree $\max m_1,\, m_2 + 1 + 2$. ▲

We present various ways to retrieve $\mathcal{P}_k$ depending on its dimension and of the user wishes.

#### 5.3.5.2 $\mathcal{P}_k$ as a polynomial projection space

When working with polynomial spaces, $\mathcal{P}_k$ has to fulfil an additionnal property coming from the fact that we work with vectorial basis functions. The chosen set $\mathcal{P}_k$ should then be defined, or at least could be recast as a Cartesian product of scalar polynomial spaces $\mathcal{P}_{k,\,i}$ whose degree is less than some integer $k_i$ in the $i^{\text{th}}$ variable. Meaningly, we should be able to write:

$$\begin{aligned} \mathcal{P}_k &= \bigtimes_{i=1}^{d} \mathcal{P}_{k,\,i} \\ \dim \mathcal{P}_k &= \dim \mathbb{H}_k(K) - \dim \mathbb{H}_k(K)|_{\partial K}. \end{aligned}$$

***Note.*** In practice, we ask $k_i \leq \overline{k_i} \colon = k$ where $k$ represents the order of the space. There, it prevents the projections onto $x^{k_i}$ to be identically identical for any polynomial function of the space. Indeed, their order would be lower and we would loose the unisolvence. Another way would be to ask for a $\overline{k_i} \leq \max\{m_1,\, m_2\}$ which is enough to cope with the dimension that guarantees the existence of a set of points of non-zero measure where the projections of any function of $\mathbb{H}_k(K)$ on some test basis functions differ from each other and are not vanishing. ▲

Furthermore, since we do not want to totally disregard one component of the discretisation in one spatial direction arbitrarily, we forbid the use of $\mathcal{P}_{k_i,\,i} = \varnothing$ for any $i \in [\![1,\, d]\!]$ unless $\dim \mathcal{P}_k = 0$. Thus, we have three possibilities to define $\mathcal{P}_k$;

$$(a) \quad \mathcal{P}_k = \bigtimes_{i=1}^{d} \mathbb{P}_{k_i}, \qquad \text{with } k_i \leq \overline{k_i}, \quad \forall i \in [\![1,\, d]\!] \tag{89}$$

$$(b) \quad \mathcal{P}_k = \bigtimes_{i=1}^{d} \mathbb{P}_{k_i^1, \cdots, k_i^d}, \quad \text{with } \substack{k_i^i \leq \overline{k_i}+1 \\ k_i^j \leq \overline{k_i},\, i \neq j} \quad \forall i,\, j \in [\![1,\, d]\!] \tag{90}$$

$$(c) \quad \mathcal{P}_k = \bigtimes_{i=1}^{d} \mathcal{P}_{k,\, i}, \qquad \text{with } \mathcal{P}_{k,\, i} = \begin{cases} \mathbb{P}_{k_i^1, \cdots, k_i^d}, & \substack{k_i^i \leq \overline{k_i}+1 \\ k_i^j \leq \overline{k_i},\, j \neq i} \\ \text{or} & \\ \mathbb{P}_{k_i}, & k_i \leq k. \end{cases} \tag{91}$$

***Remark.***

• Note that depending on the dimension of $\mathcal{P}_k$ it may happen that no space built on either (a) or (c) expressions can naturally have the same dimension as $\mathcal{P}_k$. Therefore, in the definition of (a) and (c), subspaces of $\mathbb{P}_{k_i}$ and $\mathbb{P}_{k_i^1, \cdots, k_i^d}$ may also be considered.

• The polynomial spaces $\mathbb{P}_{k_i}$, $k_i \leq \overline{k_i}$ for all $i \in [\![1,\, d]\!]$ used in the expressions (a) and (c) have been restricted from a more general definition for the sake of concision. One could in fact use slightly larger sets instead by asking that $\mathcal{P}_k \subset \mathbb{Q}_{\max\{m_1,\, m_2\}} \cup (x\mathbb{Q}_{[m_2]})$. There, it would mean that for any space $\mathcal{P}_{k,\, i}$, all the variables $x_j$, $j \neq i$ can be raised up to the power of $\max\{m_1,\, m_2\}$, and only the variable $x_i$ may be be raised to the power $\max\{m_1,\, m_2 + 1\}$.

• We only set $\mathcal{P}_k = \varnothing$ when $\dim \mathcal{P}_k = 0$. ▲

***Example.*** *We* consider the case $d = 2$, $n = 6$, $l_1 = 0$, $l_2 = 0$, $m_1 = -1$ and $m_2 = -1$. There, $\mathcal{P}_k = \varnothing$ (*see Figure 56*).

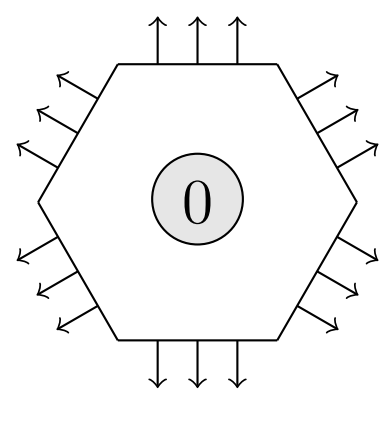

Fig. 56: Example of a configuration leading to an empty $\mathcal{P}_k$.

Here, we have $\dim \mathbb{H}_k(K) = 18$. We distribute $6(0 + 3) = 18$ normal moments on the edges. Therefore, we have a perfect distribution without any remaining internal moments. Thus, we are conformal and can set $\mathcal{P}_0 = \varnothing$.

Note that given the above definition of the indices $l_1$, $l_2$, $m_1$ and $m_2$, the case where $\dim \mathcal{P}_0 = 0$ appears only once, when $k = 0$. Thus, $\mathcal{P}_0 = \varnothing$ only if $k = 0$. ♦

For the sake of the explanation and with respect to the *Paragraph 2.3*, we will consider without loss of generality the spaces $\mathcal{P}_k$ endowed with their canonical basis in all what follows.

***Example.*** Let us assume that $\mathcal{P}_k$ is of the shape of $\mathcal{P}_k = \mathbb{P}_{2,1} \times \mathbb{P}_{1,2}$. Then, taking the canonical basis $\{1,\, x,\, xy,\, y,\, y^2\}$ of $\mathbb{P}_{1,2}$, the set of corresponding internal moments can be simply defined as

$$\sigma_i \colon q \mapsto \int_K q \cdot \begin{pmatrix} b_i(y,\, x) \\ 0 \end{pmatrix} \,\mathrm{d}x \quad \forall i \in [\![1,\, \dim \mathbb{P}_{2,1}]\!]$$

and

$$\sigma_{i+\dim \mathbb{P}_{2,1}} \colon q \mapsto \int_K q \cdot \begin{pmatrix} 0 \\ b_i(x,\, y) \end{pmatrix} \,\mathrm{d}x \quad \forall i \in [\![1,\, \dim \mathbb{P}_{1,2}]\!]$$

with $b_i(x,\, y)$ the $i^{\text{th}}$ basis function of $\mathcal{P}_k$. ♦

We now have a general sketch to define and generate $\mathcal{P}_k$. However, one can notice that the above designed spaces are not all symmetric with respect to the discretisation quality in each space direction. Therefore, allowing or not the use of non - symmetric $\mathcal{P}_k$ space is a choice one has to make before looking for a suitable definition of $\mathcal{P}_k$.

Please note that this choice have a huge impact on the possibilities one has for defining $\mathcal{P}_k$. Indeed, if for a non - symmetric choice it is always possible to find a classical polynomial space contained in $RT_k(K)$ having the required dimension, the same does not apply for the symmetric choice. Indeed, it is most unlikely that the dimension of the set of internal moments matches exactly with the one of a symmetry preserving space. It is especially true when $\mathcal{P}_k$ is a Cartesian product of several scalar polynomial spaces. In those cases, surgical treatments in the definition of $\mathcal{P}_k$ are required to preserve the symmetric properties. We briefly explore those two cases.

**Symmetry preserving spaces** Choosing to define the set of internal moments from a symmetric space $\mathcal{P}_k$ is a safe option. The solution will be represented homogeneously across the spatial coordinates. However, it brings severe difficulties when building the set. Indeed, it is usually not possible to match the dimension of a symmetric Cartesian product of classical polynomial spaces with the dimension of the set of internal moments one would like to generate.

***Example.*** We give several example for constructing the space $\mathcal{P}_k$.

• Choosing $\mathcal{P}_k$ in the shape of $(c)$ of (89) never leads to a symmetric space. Indeed, if we exclude the one dimensional case, no space $\mathbb{P}_l \times \mathbb{P}_{m,n}$ can be symmetric for any integers $l$, $m$ and $n$ non - equally null. Indeed, the definition of those spaces prevent any permutation to recover one or the other (*i.e.* $\mathbb{P}_l$ and $\mathbb{P}_{m,n}$ are not identical, even up to any permutation of the coordinates).

• In two dimensions, one can never match a symmetric space dimension with the number of internal functions when this one is odd. Indeed, the space $\mathcal{P}_k$ would read $\mathcal{P}_k = \mathcal{P}_{k,1} \times \mathcal{P}_{k,2}$. Thus, we would like to have $\dim \mathcal{P}_{k_1} \times \mathcal{P}_{k_2} = 2q+1$ for some $q \in \mathbb{N}$, where both $\mathcal{P}_{k,1}$ and $\mathcal{P}_{k,2}$ are non - empty sets.

If we consider the case $\mathcal{P}_k = \mathcal{P}_{k_1} \times \mathcal{P}_{k_2} = \mathbb{P}_l \times \mathbb{P}_m$ for some integers $l$ and $m$, we would ask $(l+1)(l+2)/2 + (m+1)(m+2)/2 = 2q+1$, which is not solvable over the integers for any $(m,\, l) \in \mathbb{N}$.

Similarly, considering the case $\mathcal{P}_k = \mathbb{P}_{k_1,k_2} \times \mathbb{P}_{k_3,k_4}$ implies that we ask $(k_1+1)(k_2+1) + (k_3+1)(k_4+1)$ to be odd. Therefore, it would imply that $(k_1+1)(k_2+1)$ is odd and $(k_3+1)(k_4+1)$ even (or reversely). But since for $(k_1+1)(k_2+1)$ to be odd you need $k_1$ and $k_2$ to be even, and for $(k_3+1)(k_4+1)$ to be even you need at least one of them to be even, the set $\mathbb{P}_{k_1,k_2} \times \mathbb{P}_{k_3,k_4} = \mathbb{P}_{odd,\,odd} \times \mathbb{P}_{even,\,odd}$ or $\mathbb{P}_{odd,\,odd} \times \mathbb{P}_{odd,\,even}$ cannot be symmetric.

• To give a more computational example, let us take $n = 6$, $l_1 = 0$, $l_2 = 1$, $m_1 = 0$ and $m_2 = 0$ in two dimensions. The total dimension of $RT_k(K)$ is $\dim RT_k(K) = n(2(l_1+1)^2 + (l_2+1)^d) + 2(m_1+1)^2 + (m_2+1)^2 - m_2 + 1 = 6(1+3) + 3 = 27$. The number of internal functions is $2 \times 1 \times 2 - 1 = 3$. We then have to find three internal moments, and thus a space $\mathcal{P}_k$ whose dimension is three. Let us show that this is not possible for $\mathcal{P}_k = \mathbb{P}_l \times \mathbb{P}_m$.

We would like $(l_1+1)(l_2+2)/2 + (m+1)(m+2)/2 = 3$, leading to $(m^2+3m+2)(l^2+3l+2) = 6$. Thus, the couple of integers $(l,\, m)$ has to fulfil the relation $m^2 + l^2 + 3(m+l) = 2$. Spanning computationally the integers to test the feasibility of this relation, we get

$$
\begin{aligned}
m = 0 &\Rightarrow l^2 + 3l = 2 \Rightarrow \text{ impossible over integers.}\\
m = 1 &\Rightarrow l^2 + 3l = -2 \Rightarrow \text{ impossible over integers.}\\
m = 2 &\Rightarrow l^2 + 3l = -8 \Rightarrow \text{ impossible over integers.}
\end{aligned}
$$

The second member is only decreasing as $m$ increases, making the relation impossible to achieve over non-negative integers. Therefore, by recursion no space of the shape $\mathbb{P}_l \times \mathbb{P}_m$ can be designed. ♦

More generally, one can easily notice that a necessary condition for a classical $\mathcal{P}_k$ built as a Cartesian product of general $\mathbb{P}_{k,\,i}$ of spaces defined through (a), (b) or (c) in (89) to preserve the symmetry is that the dimension of $\mathcal{P}_{k,\,i}$ should be the same for every $i \in [\![1,\, d]\!]$. Indeed, otherwise no permutation between the discretised components can be set, and it breaks the definition of symmetry given in the *Section* 2.2.

Based on this observation, one can infer that no classical symmetric space can have a dimension matching the required number of internal degrees of freedom. Indeed, a classical definition of symmetric space reads

$$\mathcal{P}_k = \bigtimes_{i=1}^{d} \mathcal{P}_{k,\,i}$$

with $\mathcal{P}_{k,\,i}$ taking the form of either (a), (b) or (c) in (89). Therefore, by the previous remark the $\mathcal{P}_{k,\,i}$ should enjoy the same dimension for any $i \in [\![1,\, d]\!]$. Thus, dim $\mathcal{P}_k = d \dim \mathcal{P}_{k,\,1}$. However, we need to have

$$\dim\, \mathcal{P}_{k,\,i} = (m_1+1)^d + \frac{1}{d}((m_2+1)^d - m_2^d),$$

leading to the condition

$$(m_2+1)^d - m_2^d \equiv 0[d].$$

Developing the term $(m_2+1)^d$ one get

$$\sum_{i=1}^{d} \binom{d}{i} (m_2^i) - m_2^d \equiv 0[d].$$

Thus, one has to find a constant $c \in \mathbb{N}$ such that

$$\sum_{i=1}^{d} \binom{d}{i} (m_2^i) - m_2^d = cd. \tag{92}$$

Note that in (92), all the terms $\binom{d}{i}$ but the ones corresponding to $i = 0$ and $i = 2$ can be factorized by $d$. Thus, as the term $m_2^d$ cancels with $-m_2^d$, we only have to make sure that the term $\binom{d}{0} m_2^0 = 1$ is a multiple of $d$. However, since we work in the case $d \geq 2$, this never happens and (92) does not have integer solutions. Therefore, one can never use classically defined symmetric spaces of the shapes (a), (b) and (c) where the natural subspaces share the same dimension.

We then propose a surgical construction of symmetric spaces that will fulfil the requirements of the degrees of freedom. This method can adapted on a case by case basis depending on the user wishes.

In short, one selects the smallest symmetric space $\mathcal{P}_k$ defined as (a) or (b) that have at least a dimension bigger of equal to $\dim \mathbb{H}_k(K) - \dim \mathbb{H}_k(K)|_{\partial K}$. Then, if required, one deletes some basis functions so that the dimension of the reduced space matches the expected one. There, one has to pay attention while performing this step that this process will not break the symmetry.

For the sake of clarity, we will explain the procedure step by step while considering the following example all along this paragraph.

***Example.*** We consider the two dimensional case with $m_1 = 2$ and $m_2 = 3$. By construction, the number of edges $n$ and the values of $l_1$ , $l_2$ are left free as it does not impact our discussion. There, as detailed generally before, no classical symmetric polynomial space can be designed for $\mathcal{P}_3$. Indeed, the dimension of the degrees of freedom describing the inner cell is 25, and we have to design 25 internal moments. Gathering the two previous assertions, we know that the expression (c) cannot be used to design a symmetric space, and that no space enjoying symmetry properties can be designed in their classical sense when the number of remaining degrees of freedom is odd. Therefore, no classical definition can be set here. ♦

We start by looking for the smallest classical symmetric space (a), (b) or (c) that has at least the dimension of the number of internal degrees of freedom.

***Example.*** *Continued* We look for any classical symmetric space whose dimension is the closest to 25 from above. The closest even value is 26. Therefore, we look for space whose dimension is 26 and on which projections of functions in $\mathbb{H}_k(K)$ are neither null nor identical. For safety and with the *Remark 9* in mind, we limit the polynomial order to $\max\{m_1,\, m_2+1\}$ though higher orders may be admissible too. As it is our preferred space, we start by looking for a competitor based on the $\mathbb{Q}_k$, namely $\mathcal{P}_k = \mathbb{P}_{k_{1,1},k_{2,1}} \times \mathbb{P}_{k_{2,1},k_{2,2}}$. We then have to look for $\{k_{1,1},\, k_{2,1}\}$ such that $2(k_{1,1}+1)(k_{2,1}+1) = 26$. The indices $k_{1,2}$ and $k_{2,2}$ will then be retrieved by symmetry.

However, 13 is prime. And since the maximum order of the Laplacian of functions in the space is $m_1,\, m_2 + 1 \leq 5$, the maximum admissible polyno-mial projection order is $7 < 13$. Thus, the desired order cannot be reached and no such combination can be found. Therefore, no space based on $\mathbb{Q}_k$ can be designed.

We then turn to the space built on $\mathbb{P}_k$ spaces. We want to look for $l$ such that $(l+1)(l+2) = 26$. The other index $m$ in $\mathcal{P}_3 = \mathbb{P}_l \times \mathbb{P}_m$ will be set with $l = m$ to guarantee symmetry. However here, this equation cannot be solved over the integers. Thus, we have to look for even higher dimension space.

We can show that the space of dimension 30 is the first case where an admissible space can be defined. Indeed, focusing our hunt on spaces based on $\mathbb{Q}_k$, we retrieve that $\mathcal{P}_3 = \mathbb{P}_{4,2} \times \mathbb{P}_{2,4}$ is an admissible space. We will now reduce the space so that its dimension comes down to 25, preserving its symmetry. ♦

During this process, if several spaces of the same dimension can be designed one should select the type of the space which is the most suited to the shape of the element. Indeed, classically for triangles the choice of $\mathcal{P}_k$ based on $\mathbb{P}_k$ would be made. Here, we work here with more complex polygonal shapes for whose the use of $\mathbb{Q}_k$ spaces is more recommended (*see e.g. [4]*). Furthermore, $\mathbb{H}_k(K)$ gets its finite dimension from the $\mathbb{Q}_k$ space. Therefore, if there is a choice, one would naturally prefer using $\mathcal{P}_k$ based on the $\mathbb{P}_{k_1,k_2}$ spaces. It is more compliant with the discretisation space $RT_k(K)$.

More generally, the hunt reads as follows. First, we try to find the closest largest set based on $\mathbb{Q}_k$ whose dimension is a multiple of $d$ and is the closest from above to $\mathcal{P}_k$. If impossible, we look in the same spirit, but for spaces based on $\mathbb{P}_k$. If it is still impossible, we enlarge again the dimension paying attention of its divisibility by $d$ and start again the hunt. One could also consider only $\mathbb{Q}_k$ and increase the admissible order to preserve the coherence of the type of the space rather than the coherence of the distribution of the order of the functions.

Then, we need to remove the relevant number of basis function from the found space taking care of the symmetry. In general, if the number of functions to remove is a multiple of the dimension, it can be easily done by taking some sample function and removing it as well as all the ones that can be obtained by a circular variable permutation. If one needs to remove only $d-1$ functions, the one could merge $d$ identical motives into a single vector. By example, the set $\left\{\begin{pmatrix}1\\0\end{pmatrix}, \begin{pmatrix}0\\1\end{pmatrix}\right\}$ can be reduced to $\left\{\begin{pmatrix}1\\1\end{pmatrix}\right\}$, or the set $\left\{\begin{pmatrix}xyz\\x\\zy^2\end{pmatrix}, \begin{pmatrix}xz^2\\xyz\\y\end{pmatrix}, \begin{pmatrix}z\\x^2y\\xyz\end{pmatrix}\right\}$ can be reduced to $\left\{\begin{pmatrix}xyz+xz^2+z\\xyz+yx^2+x\\xyz+zy^2+y\end{pmatrix}\right\}$.

That way, we lower the dimension while preserving the symmetry. It is particularly efficient when we know *a priori* that we are solving an isotropic problem.

***Example.*** *Continued* Here we want to reduce the dimension by five in dimension two. That means that we have to remove two functions with their permutations and merge two other ones into one symmetric vector.

We detail the process by assuming that a basis of $\mathcal{P}_4 = \mathbb{P}_{4,2} \times \mathbb{P}_{2,4}$ is the canonical one, without loss of generality. The canonical basis of $\mathcal{P}_4$ reads:

$$\left\{ \begin{pmatrix}1\\0\end{pmatrix}, \begin{pmatrix}x\\0\end{pmatrix}, \begin{pmatrix}x^2\\0\end{pmatrix}, \begin{pmatrix}x^3\\0\end{pmatrix}, \begin{pmatrix}x^4\\0\end{pmatrix}, \begin{pmatrix}xy\\0\end{pmatrix}, \begin{pmatrix}xy^2\\0\end{pmatrix}, \begin{pmatrix}x^2y\\0\end{pmatrix}, \begin{pmatrix}x^2y^2\\0\end{pmatrix}, \begin{pmatrix}y\\0\end{pmatrix}, \begin{pmatrix}x^3y\\0\end{pmatrix}, \begin{pmatrix}x^3y^2\\0\end{pmatrix}, \begin{pmatrix}x^4y\\0\end{pmatrix}, \begin{pmatrix}x^4y^2\\0\end{pmatrix}, \begin{pmatrix}y^2\\0\end{pmatrix}, \begin{pmatrix}0\\1\end{pmatrix}, \begin{pmatrix}0\\x\end{pmatrix}, \begin{pmatrix}0\\x^2\end{pmatrix}, \begin{pmatrix}0\\x^3\end{pmatrix}, \begin{pmatrix}0\\x^4\end{pmatrix}, \begin{pmatrix}0\\y\end{pmatrix}, \begin{pmatrix}0\\xy\end{pmatrix}, \begin{pmatrix}0\\xy^2\end{pmatrix}, \begin{pmatrix}0\\x^2y\end{pmatrix}, \begin{pmatrix}0\\x^2y^2\end{pmatrix}, \begin{pmatrix}0\\x^3y\end{pmatrix}, \begin{pmatrix}0\\x^3y^2\end{pmatrix}, \begin{pmatrix}0\\x^4y\end{pmatrix}, \begin{pmatrix}0\\x^4y^2\end{pmatrix}, \begin{pmatrix}0\\y^2\end{pmatrix} \right\}$$

We can merge any similar motive, as by example

$$\left\{ \begin{pmatrix}x^2\\0\end{pmatrix}, \begin{pmatrix}0\\y^2\end{pmatrix} \right\} \quad \text{to} \quad \left\{ \begin{pmatrix}x^2\\y^2\end{pmatrix} \right\},$$

which enforces an isotropic representation of quadratic behaviour. We can also forget some high order terms, as by example

$$\left\{ \begin{pmatrix}x^4y\\0\end{pmatrix}, \begin{pmatrix}0\\y^4x\end{pmatrix}, \begin{pmatrix}0\\x^3y\end{pmatrix}, \begin{pmatrix}xy^3\\0\end{pmatrix}, \right\} \to \varnothing.$$

Notice that when removing totally functions without merging them, one has to get rid of permutation couples to preserve the symmetry. After those changes, $\mathcal{P}_3$ can be defined and enjoys the desired dimension. ♦

***Remark.*** It is worth to notice that the set of functions that will be removed characterise the set that will remain, and therefore the moments' definition. Thus, the question of the choice of functions to omit is crucial. Furthermore, even if it is linked to the features one wants to instil to the element, this choice is arbitrary. ▲

Though arbitrary, the choice may be helped by keeping in mind that the set of behaviours that can be accurately represented by the space $\mathbb{P}_k$ is sensitive to this choice. In particular, if one knows *a priori* some qualitative behaviour of the represented solution, then the choice can be oriented so that the quality loss will be at its minimum.

***Example.*** *Two dimensional case*
• If one knows that the constant part of the solution is isotropic, on can merge without risk $\{(0,\, 1)^T,\, (1,\, 0)^T\}$ to $\{(1,\, 1)^T\}$.

• If one knows that the solution is close to a polynomial of order four, on can remove without so much impact the set $\{(x^3,\, 0)^T,\, (0,\, y^3)^T\} \to \varnothing$.

• If one knows that the solution is independent of $x^3$ for example, one could want to omit $\{(x^3,\, 0)^T,\, (0,\, x^3)^T\}$. However, the obtained set will not be symmetric anymore. ⧫

***Remark.***
• In general, the functions to remove are chosen on a case by case basis, with respect to the problem one wants to solve.

• This method is not suited to give a general framework. ▲

Using symmetric spaces is always guaranteeing a discretisation that does not gather troubles in itself. However, their definition is subject to arbitrary choices and one has to look deep into their structure to know their discretisation properties. If one does not really care of the homogeneity of the discretisation across the coordinates, then the use of non - symmetric space can be a good alternative as they can always be classically defined.

**Non - symmetric spaces** Allowing non - symmetric spaces can be dangerous when one has no idea of the behaviour that will be represented. Indeed, the representation quality will not be homogeneous across the space. However, in case of an anisotropic behaviour for which we know the dominant space variable, it can be very convenient.

One nice specificity of non - symmetric spaces is that for any number of desired internal functions, a classical construction of the form (a), (b) or (c) will always be possible. Indeed, the shape (b) will always furnish admissible spaces. One can see it by the onto property of its dimension onto the set of natural numbers.

***Example.*** *Two dimensional case* The problem of finding two couples $(l,\, m)$ and $(r,\, s)$ such that

$$\begin{cases} \dim \mathbb{P}_{l,\, m} \times \mathbb{P}_{r,\, s} = (l+1)(m+1) + (r+1)(s+1) = \alpha, \\ \quad \alpha \leq 2(\max\{m_1,\, m_2\}+1)^2 \in \mathbb{N} \\ \quad m \leq \max\{m_1,\, m_2\},\, l \leq \max\{m_1,\, m_2\}+1, \\ \quad r \leq \max\{m_1,\, m_2\},\, s \leq \max\{m_1,\, m_2\}+1 \end{cases}$$

has always a solution (true for any permutation in $\{r,\, s\}$ and $\{m,\, l\}$). ♦

In particular, as shown previously it may happen that non-symmetric spaces are the only possibility to build $\mathcal{P}_k$ through classical polynomial spaces.

***Example.*** *L*et us take the previous example where no classical symmetric spaces could have been used. We have $d = 2$, $m_1 = 2$, $m_2 = 3$ and are looking for non - symmetric $\mathcal{P}_3$ of the shape (a), (b) or (c) that would fit in our setting. As we prefer working with $\mathbb{Q}_k$ space, and since with non - symmetric spaces it is always possible, we start by looking for $\mathcal{P}_k$ in the form $\mathbb{P}_{l,\, m} \times \mathbb{P}_{r,\, s}$. We would like

$$\begin{cases} (l+1)(m+1) + (r+1)(s+1) = 25 \\ l \leq 5, \quad r \leq 4 \\ m \leq 4, \quad s \leq 5. \end{cases}$$

Spanning over the available cases shows that we have one solution;

$$l = 5, \quad r = 2, \quad m = 3, \quad s = 5 \tag{93}$$

The space $\mathbb{P}_{5,\, 3} \times \mathbb{P}_{2,\, 5}$ is the only available as its permutations within each subspace do not allow its monomials to define a Laplacian admissible in the space $\mathbb{H}_k(K)$. ♦

***Example.*** *F*or the case $m_1 = 1$ and $m_2 = 2$, allowing non - symmetric spaces provides more possibilities for $\mathcal{P}_k$. Indeed, the choice (c) could be considered (though not recommended in general). We would look for $\mathcal{P}_2$ in the shape $\mathcal{P}_2 = \mathbb{P}_l \times \mathbb{P}_{r,\, s}$, where several solutions can be derived. Indeed, we have to set up, up to possibly permutations between $r$ and $s$, triplets $(l,\, r,\, s)$ such that $(l+1)(l+2)/2 + (r+1)(s+1) = 21$. In particular, we can derive some working cases.

| l | r | s |
|---|---|---|
| 0 | 4 | 3 |
| 0 | 3 | 4 |
| 2 | 4 | 2 |
| 2 | 2 | 4 |

One could also use $l = 1$, $r = 2$, $s = 5$ with $\mathcal{P}_2 = \mathbb{P}_l \times \mathbb{P}_{r,s}$ and $l = 1$, $r = 5$, $s = 2$ with $\mathcal{P}_2 = \mathbb{P}_{r,s} \times \mathbb{P}_l$. Indeed, the shift in $x$ or $y$ in $B_k$ allow to project further in one order of the direction. ♦

Here again we are facing arbitrary choices in the selection of $\mathcal{P}_k$, and we fall back in the same discussion that we had when deleting functions from symmetric spaces. It is very unlikely that we will obtain uniqueness for non - symmetric $\mathcal{P}_k$, even for those build exclusively on $\mathbb{Q}_k$ spaces.

***Remark.*** An equivalent way to build admissible non-symmetric spaces would be to use the closest symmetric space as above, and to forget functions without taking care of preserving the symmetry of the space. ▲

### Latest remarks

• In conclusion, when one does not know *a - priori* the qualitative behaviour of the solution he wants to approximate, symmetric spaces should be used. Indeed, the structure of the space would preserve any isotropic behaviour and furnish a homogeneous discretisation in every direction.

However, it can happen that one has to reduce hardly the dimension of a classical symmetric space by leaving out some basis function in order to cope with the dimension of the inner discretisation liberty of $\mathbb{H}_k(K)$. This choice is arbitrary and the approximation is very sensitive on this choice.

If one prefers to use classically defined spaces in any case, or wants to represent an anisotropic behaviour, the use of non - symmetric spaces is a good choice. A classical definition can always be set up. However, there is no uniqueness of admissible $\mathcal{P}_k$ and its selection impact the way the solution will be represented.

• So far, we only spoke about the dimension to define admissible $\mathcal{P}_k$ spaces. We should not forget that the achieved moments should also be unisolvent for $\mathbb{H}_k(K)$. If by construction, the set $\{\sigma\}$ of internal degrees of freedom will be free and if the total number of degrees of freedom matches the dimension of $\mathbb{H}_k(K)$, it is not said that the total set of degrees of freedom will be unisolvent for $\mathbb{H}_k(K)$. One way to enforce it is to set $\mathcal{P}_k$ so that every function in the set of dual functions to internal degrees of freedom vanishes on the boundaries of $K$.

One could also use projection spaces involving the corresponding Poisson's solutions to stay more coherent with the shape of $\mathbb{H}_k(K)$.

#### 5.3.5.3 Projection spaces based on Poisson's problems

As we work on an space based on Poisson's problems, it is more convenient to have projections on Poisson's functions. Then, we can define $\mathcal{P}_k$ as a set of Poisson's functions solutions of

$$\begin{aligned} \Delta u &= p_k \\ u|_{\partial K} &= 0, \end{aligned}$$

where $p_k$ is a function belonging to the polynomial space $\mathcal{P}_k$ presented in the above subsection. The homogeneous Dirichlet conditions are not necessary here, but help to have smoother projections.

***Note.*** Any other Dirichlet problem would be admissible, but one should not use Neumann or Robin conditions as it would not fit in the space $\mathbb{H}_k(K)$ and result in unisolvence problems. ▲

### 5.3.6 Validity of designed elements

The element above - defined, we have the following property for any choice of a free set $\mathcal{P}_k \subset \mathbb{H}_k(K)$.

**Property 5.21**

For any $q \in \mathbb{H}_k(K)$, we have;

$$\begin{cases} (Ia) = 0 \text{ or } (Ib) = 0 \text{ or } (IIa) = 0 \text{ or } (IIb) = 0 \\ \displaystyle\int_K q \cdot p_k \, \mathrm{d}x = 0, \qquad \text{for all } p_k \in \mathcal{P}_k \end{cases}$$

$$\Rightarrow \quad q = 0,$$

where by $(Ia) = 0$ we mean that all the degrees of freedom considered in the configuration $(Ia)$ vanish.

***Remark.*** In the following proof, we refer to the functions $p_k$ by the term "kernel", while using the term "integrand" to represent the term $q \cdot p_k$. Immediate transfer of this designation apply to the normal moment based degrees of freedom. ▲

In a nutshell, the proof holds by the *Assumptions 5.14* ensuring the linear independence of the set of point-wise values and moment's integrands. The linearity of the integral operators transfers then this independence to the moments themselves, characterising any function of $\mathbb{H}_k(K)$ independently on the boundary and within the cell.

On the boundary, as the kernels are scalar polynomials and the functions of $\mathbb{H}_k(K)$ are vector polynomials, selecting the appropriate number of degrees of freedom in any of the sets $(Ia)$, $(Ib)$, $(IIa)$ or $(IIb)$ ensures a unique characterisation.

Within the element, the characterisation is done through projections over linearly independent sets on which projections of functions in $\mathbb{H}_k(K)$ are not identically identical (*i.e.* they differ at least on a subset of non-zero measure). We detail this determination process within the following proof.

***Proof.***
Let us first recall that the space $\mathbb{H}_k(K)$ is built from blocks of independent functions. In particular, the boundary behaviour of functions living in $\mathbb{H}_k(K)$ is independent of their behaviour within the inner cell. Therefore, reading out the structure of $\mathbb{H}_k(K)$ and making use of the superposition theorem, any function $q \in \mathbb{H}_k(K)$ can be recast as

$$q = f\mathbb{1}_{\partial K} + g\mathbb{1}_{\mathring{K}} \tag{94}$$

for two functions $f$ and $g$ belonging to $\mathbb{H}_k(K)$. As a consequence, characterising a function $q \in \mathbb{H}_k(K)$ comes down to characterising the independent functions $f$ and $g$ on the distinct supports $\partial K$ and $\mathring{K}$, respectively. Note also that necessarily, $f|_{f_j} \in \bigtimes_{i=1}^{d} \mathbb{Q}_{\max\{l_1, l_2+1\}}(\mathbb{R}^{d-1})$ for any face $f_j \in \partial K$.

We demonstrate that any admissible extraction (in the sense of the *Requirements 5.14*) from either of the four sets of degrees of freedom ($Ia$, Internal), ($Ib$, Internal), ($IIa$, Internal), ($IIb$, Internal) fully characterises the functions $f$ and $g$, independently.

• We first show that any above defined set of degrees of freedom preserve the independence of the boundary and inner characterisations. To this aim, we combine the relation (94) with the all possible definitions of the degrees of freedom encountered in the definition of the four elements. It comes the following relations.

All global normal moments lead to an expression of the form

$$\sigma(q) = \int_{\partial K} q \cdot n \, p_k = \int_{\partial K} (f \mathbb{1}_{\partial K} + g \mathbb{1}_{\mathring{K}}) \cdot n \, p_k = \int_{\partial K} f \cdot n \, p_k$$

for some polynomial function $p_k$ living on $\partial K$. On the other side, as $x_{im} \in \partial K$, the global degrees of freedom that are built from pointwise values read

$$\sigma(q) = q(x_{im}) \cdot n = f(x_{im}) \cdot n \mathbb{1}_{\partial K}(x_{im}) + g(x_{im}) \cdot n \mathbb{1}_{\mathring{K}}(x_{im}) = f(x_{im}) \cdot n.$$

Similar relations for coordinate - wise degrees of freedom can be derived, that is;

$$\sigma(q) = \int_{\partial K} q_{x_i} n_{x_i} \, p_k = \int_{\partial K} (f_{x_i} \mathbb{1}_{\partial K} + g_{x_i} \mathbb{1}_{\mathring{K}}) n_{x_i} \, p_k = \int_{\partial K} f_{x_i} n_{x_i} \, p_k$$

and $$\begin{aligned}\sigma(q) &= q_{x_i}(x_{im}) n_{x_i} = f_{x_i}(x_{im}) n_{x_i} \mathbb{1}_{\partial K}(x_{im}) + g_{x_i}(x_{im}) n_{x_i} \mathbb{1}_{\mathring{K}}(x_{im}) \\ &= f_{x_i}(x_{im}) n_{x_i}.\end{aligned}$$

Therefore, in any of the configurations $(Ia)$, $(Ib)$, $(IIa)$ and $(IIb)$ no contribution of the function $g$ representing the inner part of the cell is involved in the normal degrees of freedom.

The mirror case is achieved with the internal moments, leading to

$$\sigma(q) = \int_{K} q \cdot p_k = \int_{K} (f \mathbb{1}_{\partial K} + g \mathbb{1}_{\mathring{K}}) \cdot p_k = \int_{K} g \cdot p_k,$$

where $p_k$ stands for any Poisson's function living in $\mathbb{H}_k(K)$ or any polynomial function defining the second member of a Poisson's problem involved in the definition of $\mathbb{H}_k(K)$. There, the function $f$ representing the boundary part of the function $q$ is not involved, that for any definition of the space $\mathcal{P}_k$ generating the internal moments.

Thus, by linearity we can decompose the degrees of freedom $\{q \mapsto \sigma_i(q)\}_i$ in the following matrix form.

$$\begin{matrix}\text{Normal Dofs values}\left\{ \right. \\ \\ \text{Internal Dofs values}\left\{\right.\end{matrix}\begin{pmatrix}\sigma_1 \\ \vdots \\ \sigma_{N_N} \\ \sigma_{N_N+1} \\ \vdots \\ \sigma_{N_I}\end{pmatrix} = \begin{pmatrix}\text{Normal moments applied to } f & 0 \\ 0 & \text{Internal moments applied to } g\end{pmatrix}\begin{pmatrix} f \\ g\end{pmatrix}$$

Clearly, there is no interconnection between the function's characterisation on the boundaries and the one performed within the element. Thus, showing the *Property 5.21* reduces to show independently that

$$\{(Ia) = 0 \text{ or } (Ib) = 0 \text{ or } (IIa) = 0 \text{ or } (IIb) = 0\} \quad \Rightarrow \quad f|_{\partial K} = 0$$

and

$$\int_K g \cdot p_k \,\mathrm{d}x = 0, \quad \text{for all } p_k \in \mathcal{P}_k \quad \Rightarrow \quad g|_{\mathring{K}} = 0.$$

• Let us first consider the boundary characterisation. There, by definition of the spaces $\mathcal{H}_{l_1}$ and $\mathcal{H}_{l_2}$, the function $f|_{\partial K}$ is discontinuous at the polytope's vertices and can be decomposed into $n$ vectorial polynomial functions $\{f_j\}_j$ with distinct supports, each of them matching one particular face of the polytope. Thus, we can write

$$f|_{\partial K} = \sum_{j=1}^{n} r_j \mathbb{1}_{f_j}$$

with $r_j \in \bigtimes_{i=1}^{d} \mathbb{Q}_{\max\{l_1, l_2+1\}}(f_j)$ and $f_j$ any face belonging to $\partial K$. With a similar argument than in the previous point, the characterisation of $f|_{\partial K}$ can therefore be done edge-wise, and the determination matrix becomes block - diagonal. We discuss here the characterisation on one particular edge $f_j$ by showing the invertibility of the corresponding matrix block. The arguments naturally transpose to the other ones.

In this perspective, let us show that for any $r_j \in \bigtimes_{i=1}^{d} \mathbb{Q}_{\max\{l_1, l_2+1\}}(f_j)$, it holds $\{(Ia)|_{f_j} = 0 \text{ or } (Ib)|_{f_j} = 0 \text{ or } (IIa)|_{f_j} = 0 \text{ or } (IIb)|_{f_j} = 0\} \Rightarrow r_j = 0$, where $(\cdot)|_{f_j}$ represents the subset of the degrees of freedom $(\cdot)$ whose support (or evaluation point for point-values) matches (or lies on) $f_j$.

First of all, we recall that on the face $f_j$ the function $r_j$ is a multi-valued polynomial of the form

$$r_j|_{f_j} = \begin{pmatrix} a_{0,1} \\ \vdots \\ a_{0,d} \end{pmatrix} + \sum_{i=\dim(\mathcal{H}_0)}^{\dim(\mathcal{H}_{l_1}\cap\mathcal{H}_{l_2})} \begin{pmatrix} b_{\xi_1(i)} + a_{i,1} \\ \vdots \\ b_{\xi_d(i)} + a_{i,d} \end{pmatrix} m_{\alpha_i(x)} + \sum_{i=\dim(\mathcal{H}_{l_1}\cap\mathcal{H}_{l_2})}^{\dim(\mathcal{H}_{l_2})} \begin{pmatrix} x_1 b_i \\ \vdots \\ x_d b_i \end{pmatrix} m_{\alpha_i(x)},$$

where $m_{\alpha_i}$ represents a monomial of $\mathbb{Q}_{\max\{l_1,l_2\}}$ of multi-index degree $\alpha_i$ such that the set $\{m_{\alpha_i}\}_{i=\dim(\mathcal{H}_l)}^{\dim(\mathcal{H}_m)}$ forms a base of $\mathcal{H}_m \setminus \mathcal{H}_l$. The terms $\{\xi_i\}_{i\in[\![1,d]\!]}$ are the permutation operators detailed in the equation (68).

Note that the coefficients $a$ are defined coordinate - wise while the coefficients $b$ are identical for all the components. The function $r_j$ is therefore determined by $\dim(\{\{a_{i,m}\}_{\substack{i\in[\![0,\dim(\mathcal{H}_{l_1}\setminus\mathcal{H}_0)]\!],\\ m\in[\![1,d]\!]}}, \{b_i\}_{\substack{i\in[\![\dim(\mathcal{H}_{l_1}\setminus\mathcal{H}_0),\\ \dim(\mathcal{H}_{l_2})]\!]}}\})$ coefficients.

As in all configurations the function $r_j$ is determined only through its normal components, let us use the above expression to derive them more specifically. With the normal $n_j = (n_{jx_1}, \cdots, n_{jx_d})$ to the face $f_j$, it comes

$$r_j \cdot n_j|_{f_j} = \sum_{m=1}^{d} a_{0,m} n_{jx_m} + \sum_{i=\dim(\mathcal{H}_0)}^{\dim(\mathcal{H}_{l_1}\cap\mathcal{H}_{l_2})} \sum_{m=1}^{d} (b_{\xi_m(i)} + a_{i,m}) n_{jx_m} m_{\alpha_i}(x) + \sum_{i=\dim(\mathcal{H}_{l_1}\cap\mathcal{H}_{l_2})}^{\dim(\mathcal{H}_{l_2})} c_j\, b_i\, m_{\alpha_i}(x),$$

where $c_j = x \cdot n_j$ is a constant term on the face $f_j$. Reordering the terms, we end up with the formulation

$$\begin{aligned} r_j \cdot n|_{f_j} = &\sum_{m=1}^{d} \left( \left( a_{0,m} + \sum_{i=\dim(\mathcal{H}_0)}^{\dim(\mathcal{H}_{l_1}\cap\mathcal{H}_{l_2})} a_{i,m} m_{\alpha_i(x)} \right) n_{jx_m} \right) \\ &+ \left( \sum_{m=1}^{d} n_{jx_m} \right) \sum_{i=1}^{\dim(\mathcal{H}_{l_1}\cap\mathcal{H}_{l_2})} b_{\xi_m(i)} m_{\alpha_i}(x) + c_j \left( \sum_{i=\dim(\mathcal{H}_{l_1}\cap\mathcal{H}_{l_2})}^{\dim(\mathcal{H}_{l_2})} b_i\, m_{\alpha_i}(x) \right). \end{aligned}$$

The structure of the retrieved form makes clearly emerge the coefficients that should be retrieved with respect to a coordinate-wise behaviour and the ones that depend on a global behaviour.

In addition, as all the coefficients determining $r_j$ appear in this expression, using degrees of freedom defined only from the normal components of tested functions is admissible. Thus, the four configurations fitting this framework, we only have to make sure that the set of extracted degrees of freedom are uniquely characterising each of the involved coefficients. To this aim, we explicit all the possible degrees of freedom when applied to the function $r_j$.

For the sake of clarity, we denote by $\{\sigma_{M_{i,l}}\}_{il}$ the moments designed coordinate-wise, being of the form

$$\begin{aligned} &\sigma_{M_{i,l}}: q \mapsto q_{x_i}(x_{jl})n_{jx_i} \\ \text{or} \qquad &\sigma_{M_{i,l}}: q \mapsto \int_{f_j} q_{x_i} n_{jx_i} p_l \end{aligned}$$

for some scalar polynomial $p_l$, and by $\{\sigma_{T_l}\}_l$ the ones acting globally, being of the form

$$\begin{aligned} &\sigma_{T_l}: q \mapsto q(x_l) \cdot n \\ \text{or} \qquad &\sigma_{T_l}: q \mapsto \int_{f_j} q \cdot n \, p_l \end{aligned}$$

for some scalar polynomial $p_l$. Further, for convenience we denote by $\{\sigma_{V_l}\}_l$ the global degrees of freedom that comes into play to determining the coordinate - wise coefficients, whose expressions are done in the same way as $\{\sigma_{T_l}\}_l$.

We now express those degrees of freedom depending on the coefficients $\{b_{i,m}\}_{im}$ and $\{a_i\}_i$. Using the linearity of the degrees of freedom , plugging the expression (5.3.6) in place of $q$ and setting the permutation operator directly on the multi - indices $\alpha_i$ instead of the coefficients $a_i$, we can rewrite the moments as follows.

$$\begin{aligned} \sigma_{M_{m,l}}: (\{a_{i,m}\}, \{b_i\}) \longmapsto b_{0,m} &\int_{f_j} n_{jx_m} \, p_l \\ &+ \sum_{i=\mathcal{H}_0}^{\dim(\mathcal{H}_{l_1} \cap \mathcal{H}_{l_2})} a_{i,m} \int_{f_j} m_{\alpha_i}(x) n_{jx_m} \, p_l \\ &+ \sum_{i=\mathcal{H}_0}^{\dim(\mathcal{H}_{l_1} \cap \mathcal{H}_{l_2})} b_i \int_{f_j} n_{jx_m} m_{\xi_m(\alpha_i)}(x) p_l \\ &+ \sum_{i=\dim(\mathcal{H}_{l_1} \cap \mathcal{H}_{l_2})}^{\dim(\mathcal{H}_{l_2})} b_i \int_{f_j} (x_m \, n_{jx_m} m_{\alpha_i}(x)) p_l \end{aligned}$$

$$
\begin{aligned}
\sigma_{T_l} : (\{a_{i,m}\},\, \{b_i\}) \longmapsto & \sum_{m=1}^{d} a_{0,m} \int\limits_{f_j} n_{jx_m} p_l \\
& + \sum_{m=1}^{d} \sum_{i=\dim(\mathcal{H}_0)}^{\dim(\mathcal{H}_{l_1} \cap \mathcal{H}_{l_2})} a_{i,m} \int\limits_{f_j} m_{\alpha_i}(x) n_{jx_m} p_l \\
& + \sum_{i=\mathcal{H}_0}^{\dim(\mathcal{H}_{l_1} \cap \mathcal{H}_{l_2})} b_i \sum_{m=1}^{d} \left( \int\limits_{f_j} n_{jx_m} m_{\xi_m(\alpha_i)}(x) p_l \right) \\
& + \sum_{i=\dim(\mathcal{H}_{l_1} \cap \mathcal{H}_{l_2})}^{\dim(\mathcal{H}_{l_2})} b_i \int\limits_{f_j} (c_j m_{\alpha_i}(x) p_l)
\end{aligned}
$$

Thus, defining the component - wise parts of the global moments $\sigma_{T_l}$ by $\sigma_{T_{m,l}}(q) = \int_{f_j} n_{jx_m} q\, p_l$ such that $\sigma_{T_l} = \sum_{m=1}^{d} \sigma_{T_{m,l}}$, one can express any degrees of freedom of the four considered sets as

$$
\begin{aligned}
\sigma_{M_{m,l}} : q \longmapsto & \, a_{0,m} \sigma_{M_{m,l}}(1) + \sum_{i=\mathcal{H}_0}^{\dim(\mathcal{H}_{l_1} \cap \mathcal{H}_{l_2})} a_{i,m} \sigma_{M_{m,l}}(m_{\alpha_i}) \\
& + \sum_{i=\mathcal{H}_0}^{\dim(\mathcal{H}_{l_1} \cap \mathcal{H}_{l_2})} b_i \, \sigma_{M_{m,l}}(m_{\xi_m(\alpha_i)}) + \sum_{i=\dim(\mathcal{H}_{l_1} \cap \mathcal{H}_{l_2})}^{\dim(\mathcal{H}_{l_2})} b_i \, \sigma_{M_{m,l}}(x_m m_{\alpha_i})
\end{aligned}
$$

and
$$
\begin{aligned}
\sigma_{T_l} : q \longmapsto & \sum_{m=1}^{d} a_{0,m} \sigma_{T_{m,l}}(1) + \sum_{i=\mathcal{H}_0}^{\dim(\mathcal{H}_{l_1} \cap \mathcal{H}_{l_2})} \sum_{m=1}^{d} a_{i,m} \sigma_{T_{m,l}}(m_{\alpha_i}) \\
& + \sum_{i=\mathcal{H}_0}^{\dim(\mathcal{H}_{l_1} \cap \mathcal{H}_{l_2})} b_i \sum_{m=1}^{d} \sigma_{T_{m,l}}(m_{\xi_m(\alpha_i)}) + \sum_{i=\dim(\mathcal{H}_{l_1} \cap \mathcal{H}_{l_2})}^{\dim(\mathcal{H}_{l_2})} a_i c_j \, \sigma_{T_l}(m_{\alpha_i}).
\end{aligned}
$$

Note that in view of deriving the determination matrix, the last term can also be decomposed as follows.

$$
\sum_{i=\dim(\mathcal{H}_{l_1} \cap \mathcal{H}_{l_2})}^{\dim(\mathcal{H}_{l_2})} b_i c_j \, \sigma_{T_l}(m_{\alpha_i}) = \sum_{i=\dim(\mathcal{H}_{l_1} \cap \mathcal{H}_{l_2})}^{\dim(\mathcal{H}_{l_2})} b_i \sum_{m=1}^{d} \sigma_{T_{l,m}}(x_m m_{\alpha_i}).
$$

Similar relations for $\sigma_V$ can ve derived from the expression of $\sigma_T$. Thus, we can rewrite the degrees of freedom as a dot product and derive the characterisation matrix $\Sigma$

$$(\sigma_{M_{1,1}}, \cdots, \sigma_{M_{d,l}}, \sigma_{V_1}, \cdots, \sigma_{V_l}, \sigma_{T_1}, \cdots, \sigma_{T_l}, )^T = \Sigma(\{a_{i,m}\}_{im}, \{b_i\}_i)^T$$

whose shape is given in the following page. We now investigate its structure.

First of all, as the number of extracted degrees of freedom from the four sets $(Ia)$, $(Ib)$, $(IIa)$ and $(IIb)$ matches the number of coefficients determining $r_j$, the matrix $\Sigma$ is squared.

Let us focus on the top two-by-two left blocks, surrounded in blue. They correspond to the coefficients that should be determined coordinate - wise. Thus, by construction, there are $\dim(\{a_{im}\}_{im}) = d(l_1+1)^{d-1}$ columns. And by the definition of the configurations $I$ and $II$, we have $\dim(\{\sigma_{M_{i,j}}\}_{ij}) = d + d((l_1+1)^{d-1} - 1) = d(l_1+1)^{d-1}$. Therefore, this submatrix is squared.

Furthermore, each subblock corresponds to one member of the decomposition of $q$ tested through coordinate - wise degrees of freedom whose kernels are built on the same monomial. Therefore, as the degrees of freedom $\{\sigma_{M_{i,j}}\}$ consider one normal component only, the coefficients $\{a_{im}\}_{im}$ for $m \neq j$ are not involved, and the subblocks are diagonal. Thus, those submatrices are invertible and in particular their columns and rows are linearly independent.

On the other side of the matrix, the last bottom block surrounded in deep red matches the Raviart – Thomas moments tuning members of $\mathbb{H}_k(K)|_{\partial K}$ living exclusively in $xB_k|_{f_j}$. It is then a submatrix of the classical Raviart – Thomas's one, and is necessary invertible. Thus, its rows and columns are linearly independent.

The extended bottom right submatrix highlighted in dashed red corresponds to the previously described high - order submatrix of the Raviart – Thomas's setting, enriched by the moments $\{\sigma_V\}$ tuning the behaviour of members of $\mathbb{H}_k(K)|_{f_j}$ falling in the intersection $\mathcal{H}_{l_1} \cap \mathcal{H}_{l_2}(f_j)$.

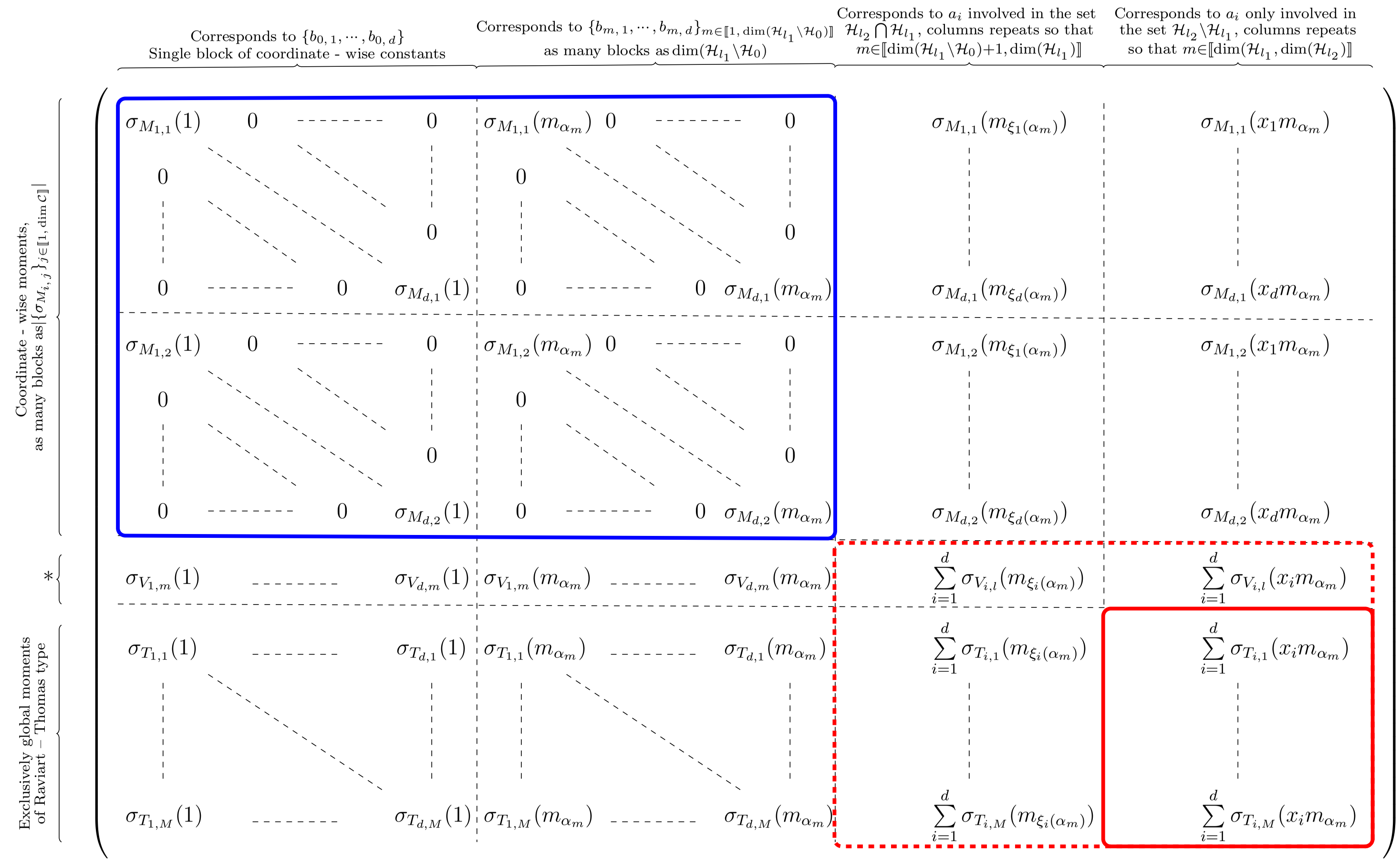

$\mathcal{C}$: Group of coordinate - wise moments associated to a same quantification process (same tesing monomial in moments or evaluation at a same same point-value).
(*) Tuning globally coefficients that have repercussion on the coordinate - wise coefficients $a_i$(misc moments).
As many rows as moments $\sigma_{V_m}$, with $m\in[\![\dim(\mathcal{H}_{l_1}\setminus\mathcal{H}_0)+1,\dim(\mathcal{H}_{l_1})]\!]$

This matrix is equivalent to the full Raviart – Thomas setting. Indeed, even if the moments $\{\sigma_V\}$ have to be slightly modified from the Raviart – Thomas setting in the configuration $I$, this modification leaves the projection order unchanged and the integrand still belongs to $\mathbb{H}_k(K)|_{f_j}$. Therefore, the dashed line block is invertible and its columns and rows are linearly independent.

There is only left to show that there is no linear dependence between rows of different row blocks. As the degrees of freedom are linear forms, it is enough to show that the integrand of moments (or polynomials constructing the point - wise values) that involve the same monomial are linearly independent.

Indeed, being linear forms whose kernels are polynomials, the degrees of freedom can combine each other only if their integrand ($q$ tested against the kernel) involve – up to constants – the same monomials. We then have to show that in both configurations, the rows involving terms whose projection onto the kernel can be expressed from a same monomial are linearly independent.

In the second configuration, $II$, this property comes automatically. Indeed, the only interaction between degrees of freedom having integrands sharing the same monomial order (and then possibly being based on the same monomial) is possible between (88c) and (88b) when $|p_k| = l_1+1$ (*resp.* (87c) *and* (87b)). Indeed, by definition of $\mathbb{H}_k(K)$, the polynomial $p_k \cdot n$ in (88b) is only of order $l_1$. However, no combination of (88c) can form the moments (88b). Indeed, for any real coefficients $c_i$ it holds

$$\sum c_i q_{x_i} n_{x_i} x_i^{l_1+1} \not\equiv q \cdot n p_k$$

for any monomial $p_k$ such that $|p_k| = l_1 + 1$. Note that in the left hand side, all the $c_i$ should be non - null to reconstruct the term $p \cdot n$. However, doing so no factorisation by a single monomial such that

$$\sum_{i=1}^{d} c_i q_{x_i} n_{x_i} x_i^{l_1+1} = \left( \sum_{i=1}^{d} c_i q_{x_i} n_{x_i} \right) p_k$$

is possible. Thus, the designed moments are linearly independent, and not row combination can occur for any tested polynomial belonging to $\mathbb{H}_k(K)|_{\partial K}$.

In the first configuration, $I$, the integrands involving the same monomials are found in the definitions (85c) and (85a) when $l_1 = 0$ (*resp.*

(86c) and (86a)). As the problem arises on constants, linear combinations would be possible in a classical Raviart – Thomas setting as the terms $p_{x_i}n_{ix}$ and $p_{y_i}n_{iy}$ can be combined to form $p \cdot n$. However, the use of the vector $v$ not collinear to $n$ in place of $n$ in the equation (85c) (*resp.* (86c)) makes the setting fulfil the *Assumptions 5.14*. Therefore, those degrees of freedom are not linearly dependent.

All in all, for both configurations all the rows are linearly independent. As by construction we have as many relations as unknowns, the matrix is invertible. Thus, we get a null kernel, meaningly

$$(Ia) = 0 \text{ or } (Ib) = 0 \text{ or } (IIa) = 0 \text{ or } (IIb) = 0 \qquad \Rightarrow \qquad f|_{\partial K} = 0.$$

• Let us now consider the internal characterisation of functions living in $\mathbb{H}_k(K)$. From the first point of the proof, it is enough to study the characterisation of $g$ within the inner cell.

By definition of $\mathbb{H}_k(K)$, any function $g \in \mathbb{H}_k(K)|_{\mathring{K}}$ can be decomposed over a set of Poisson's solutions as follows.

$$g = \sum_{i=1}^{d} \sum_{i=1}^{\dim A_k} a_{i,j} e_j u_i + \sum_{i=1}^{\dim B_k} b_j x \tilde{u}_i$$

Here, the vector $e_j$ stands for $e_j = (0, \cdots, 1, 0, \cdots, 0)^T$ where the 1 is in the $j^{\text{th}}$ position. The functions $u_i$ and $\tilde{u}_i$ represents the Poisson's solutions of the problems

$$\begin{cases} \Delta u_i = p_i,\ p_i \in \mathbb{Q}_{m_1}(K) \\ u_i|_{\partial K} = 0 \end{cases} \quad \text{and} \quad \begin{cases} \Delta \tilde{u}_i = \tilde{p}_i,\ \tilde{p}_i \in \mathbb{Q}_{[m_2]}(K) \\ \tilde{u}_i|_{\partial K} = 0 \end{cases} \tag{95}$$

where $\{p_i\}_i$ and $\{\tilde{p}_i\}_i$ form respectively a basis of $\mathbb{Q}_{m_1}(K)$ and $\mathbb{Q}_{[m_2]}(K)$. In any presented definition of the degrees of freedom, the internal characterisation is done through moment - based degrees of freedom of the form

$$\sigma_{I_k} : q \mapsto \int_K q \cdot p_k \, \mathrm{d}x$$

where the kernels $p_k \in \mathcal{P}_k$ consist of linearly independent polynomials belonging to $\mathbb{Q}_{\max\{m_1, m_2+1\}}(K)$, or of the solution of their corresponding problems of the form (95). Therefore, we can derive a characterisation matrix in the same spirit than in the case of the normal characterisation.

$$\begin{pmatrix} \sigma_{I_1} \\ \vdots \\ \sigma_{I_P} \end{pmatrix} = \begin{pmatrix} \overbrace{\begin{matrix} \int e_i u_1 \cdot p_{1,i} & \cdots & \int e_i u_A \cdot p_{1,i} \\ \vdots & & \vdots \\ \int e_i u_1 \cdot p_{P,i} & \cdots & \int e_i u_A \cdot p_{P,i} \end{matrix}}^{\substack{\text{Block repats as many time as coordinates;} \\ i \in [\![1, d]\!]}} & \Bigg| & \overbrace{\begin{matrix} \int x\tilde{u_1} \cdot p_1 & \cdots & \int x\tilde{u_B} \cdot p_1 \\ \vdots & & \vdots \\ \int x\tilde{u_1} \cdot p_P & \cdots & \int x\tilde{u_B} \cdot p_P \end{matrix}}^{\text{Single block}} \end{pmatrix} \begin{pmatrix} \left.\begin{matrix} a_{i,1} \\ \vdots \\ a_{i,A} \end{matrix}\right\} \text{Repeats } d \text{ times} \\ \hline b_1 \\ \vdots \\ b_P \end{pmatrix}$$

Let us consider the case where $\mathcal{P}_k$ forms a polynomial projection space. There, none of the $p_k \in \mathcal{P}_k$ is the zero function. In the same time, the functions $\{\{u_i\}_i, \{x\tilde{u_i}\}_i\}$ are linearly independent, and being solutions to some Poisson's problem with non-zero second member, they are by construction not identically vanishing on $K$. Indeed, even when $m_2 < m_1$ where second members of the problems (95) lives both in $\mathbb{Q}_{m_1}$ and $\mathbb{Q}_{[m_2]}$, it holds $\Delta(x\,\tilde{u_i}) = 2\nabla \cdot \tilde{u_i} + \Delta(\tilde{u_i})$. Thus, it is impossible to combine linearly the function $x\,\tilde{u_i}$ with functions of the set $\{u_i\}_i$.

Furthermore, the degrees of the polynomials belonging to the space $\mathcal{P}_k$ are lower or equal than the highest degree of the second members of the Poisson's problem defining the space $\mathbb{H}_k(K)$. Thus, every projection of function of $\mathbb{H}_k(K)$ onto the space $\mathcal{P}_k$ is not null. And as the internal moments are linear forms, any linear combination of those moments at fixed $p_k$ could have its integrand factorised by the kernel $p_k$ for any $p_k \in \mathcal{P}_k$, transferring the linear independency of the set $\{\{u_i\}_i, \{x\tilde{u_i}\}_i\}$ to the terms $\int u_i \cdot p_k\}_i$ for any fixed $p_k \in \mathcal{P}_k$.

Lastly, as the space $\mathcal{P}_k$ contains only linearly independent functions the previous argument can be repeated for each row of the above matrix. And as by construction the number of internal degrees of freedom matches the dimension of the space $\mathbb{H}_k(K)|_{\mathring{K}}$, the linear independence of functions of $\mathcal{P}_k$ combined with the linear independence of the tested functions transfer automatically to the moments tested against a basis of $\mathbb{H}_k(K)|_{\mathring{K}}$. Thus, the internal submatrix is invertible.

The same reasoning can be applied when $\mathcal{P}_k$ is built from the Poisson's solutions themselves, as the projections of functions would decompose the functions directly.

- Merging the above two points together, we retrieve

$$\{(Ia) = 0 \text{ or } (Ib) = 0 \text{ or } (IIa) = 0 \text{ or } (IIb) = 0\} \quad \Rightarrow \quad f|_{\partial K} = 0$$

and
$$\int_K g \cdot p_k \,\mathrm{d}x = 0, \qquad \text{for all } p_k \in \mathcal{P}_k \quad \Rightarrow \quad g|_{\mathring{K}} = 0,$$

implying by definition of any $q \in \mathbb{H}_k(K)$

$$\begin{cases} \{(Ia) = 0 \text{ or } (Ib) = 0 \text{ or } (IIa) = 0 \text{ or } (IIb) = 0 \\ \int\limits_K q \cdot p_k \, \mathrm{d}x = 0, \quad \text{for all } p_k \in \mathcal{P}_k \end{cases} \Rightarrow \quad q = 0,$$

reducing to the claim.

■

***Remark.*** We presented the proof for the four main elements under consideration. However, as expressed in the *Paragraph 5.3.3.3*, one could have chosen to extract another layout of the degrees of freedom that are also admissible under the *Assumptions 5.14*. The only change in the proof concerns the top left two by two blocks describing the coordinate - wise behaviours.

There, in dimension strictly bigger than two one may choose instead to take care of coordinate - wise components by using further global moments in the shape of $\sigma_V$, breaking the diagonal structure by full line(s) of non - vanishing contributions. However, as by the definitions of the degrees of freedom in the configurations $I$ and $II$, there cannot be more than $d+1$ polynomials having the same order being moments or point-values quantifiers. Thus, as we have $d+1$ coordinate - wise moments to tune per decomposed monomials, there is no over-determination. Further, by the *Assumption 5.14*, the linear independence of the sub-matrix's columns is ensured here too. Indeed, those polynomials cannot be linearly dependent. Thus, by linearity of the degrees of freedom, the independence of the moments is transferred to the moments and there cannot be any row dependency. The submatrix block corresponding to one specific order is therefore invertible. ▲

The unisolvence of the elements presented in the *Propositions 5.15 and 5.16* is a consequence of the above proposition.

***Proof.*** **5.3.6** *of the Propositions 5.15 and 5.16* By construction, the number of degrees of freedom matches the dimension of the space $\mathbb{H}_k(K)$. Furthermore, by the *Proposition 5.21* the kernel of the linear operator defined by the set of degrees of freedom has a null kernel. Thus, the four presented elements $(Ia)$, $(Ib)$, $(IIa)$ and $(IIb)$ are unisolvent for the space $\mathbb{H}_k(K)$.

■

***Example.*** Let us provide an example of the normal unisolvence property in the two dimensional case for the element $IIb$ with $l_1 = 0$ and $l_2 = 3$. We derive its normal submatrix corresponding to the characterisation of $f$ on one single edge $f_j \in \partial K$ which is not parallel to the axis $y$.

Here, as the edge can be seen as a subset of $\mathbb{R}$, any function $f|_{f_j} \in \mathbb{H}_k(K)|_{f_j}$ can be decomposed on a basis of $(\mathbb{P}(\mathbb{R}))^d$. In particular, it reads on the face $f_j$

$$\begin{aligned} f|_{f_j} = &\begin{pmatrix} a_{0,1} \\ a_{0,2} \end{pmatrix} + \begin{pmatrix} b_0 \\ b_0 \end{pmatrix} \circ \begin{pmatrix} x \\ y \end{pmatrix} + \begin{pmatrix} b_1 \\ b_1 \end{pmatrix} \circ \begin{pmatrix} x \\ y \end{pmatrix} \circ \begin{pmatrix} x \\ x \end{pmatrix} + \begin{pmatrix} b_2 \\ b_2 \end{pmatrix} \circ \begin{pmatrix} x \\ y \end{pmatrix} \circ \begin{pmatrix} x^2 \\ x^2 \end{pmatrix} \\ &+ \begin{pmatrix} b_3 \\ b_3 \end{pmatrix} \circ \begin{pmatrix} x \\ y \end{pmatrix} \circ \begin{pmatrix} x^3 \\ x^3 \end{pmatrix} \end{aligned}$$

for some coefficients $\{a_0, a_1\}$ and $\{b_m\}_{m=1,2,3}$, and where the coordinate $y$ can be parametrised by $x$ along the edge $f_j$. We can also directly express the normal components of the tested functions $f$ as

$$\begin{aligned} f_x\, n_x &= a_{0,1}\, n_x + b_0\, x\, n_x + b_1\, x^2\, n_x + b_2\, x^3\, n_x + b_3\, x^4\, n_x, \\ f_y\, n_y &= a_{0,2}\, n_y + b_0\, y\, n_y + b_1\, xy\, n_y + b_2\, x^2 y\, n_y + b_3\, x^3 y\, n_y, \\ f \cdot n &= a_{0,1}\, n_x + a_{0,2}\, n_y + b_0\, x\, n_x + b_0\, y\, n_y + b_1\, x^2\, n_x + b_1\, xy\, n_y \\ &\quad + b_2\, x^3\, n_x + b_2\, x^2 y\, n_y + b_3\, x^4\, n_x + b_3\, x^3 y\, n_y. \end{aligned}$$

On the side of the elements construction, the considered degrees of freedom read as follows.

Coordinate - wise moments $\sigma_{M_{ij}}$ :

$$q \mapsto \int q_x n_{jx}\, x \mathrm{d}\gamma(\mathrm{x}) \tag{96a}$$

$$q \mapsto \int q_y n_{jy}\, y \mathrm{d}\gamma(\mathrm{x}) \tag{96b}$$

Misc global moment $\sigma_{V_j}$ :

$$q \mapsto q(\mathrm{x}_{im}) \cdot n \tag{96c}$$

Global moments $\sigma_{T_j}$ :

$$q \mapsto \int q \cdot n\, p_k \mathrm{d}\gamma(\mathrm{x}) \quad \text{for } p_k \in \mathcal{H}_3(f_j) \setminus \mathcal{H}_0(f_j), \tag{96d}$$

There, $\mathrm{x}_{im}$ is the midpoint of the edge $f_j$, and the basis of $\mathcal{H}_3(f_j) \setminus \mathcal{H}_0(f_j)$ is taken as the canonical one; $\mathcal{H}_3(f_j) \setminus \mathcal{H}_0(f_j) = \{x \mapsto x,\, x \mapsto x^2,\, x \mapsto x^3\}$.

Thus, the characterisation of the function $f|_{f_j}$ reads under the matrix form;

$$\begin{pmatrix}\sigma_{M_{1,1}}\\ \sigma_{M_{1,2}}\\ \sigma_{V_1}\\ \sigma_{T_1}\\ \sigma_{T_2}\\ \sigma_{T_3}\end{pmatrix} = \left(\begin{array}{cc|c|ccc} \int x\,n_x & 0 & \int x^2 n_x & \int x^3 n_x & \int x^4\,n_x & \int x^5\,n_x \\ 0 & \int y\,n_y & \int y^2\,n_y & \int x\,y^2 n_y & \int x^2\,y^2\,n_y & \int x^3 y^2\,n_y \\ \hline n_x & n_y & C & Cx_{im} & Cx_{im}^2 & Cx_{im}^3 \\ \hline \int x\,n_x & \int x\,n_y & \int Cx & \int Cx^2 & \int Cx^3 & \int Cx^4 \\ \int x^2\,n_x & \int x^2\,n_y & \int Cx^2 & \int Cx^3 & \int Cx^4 & \int Cx^5 \\ \int x^3\,n_x & \int x^3\,n_y & \int Cx^3 & \int Cx^4 & \int Cx^5 & \int Cx^6 \end{array}\right) \begin{pmatrix}a_{0,1}\\ a_{0,2}\\ b_0\\ b_1\\ b_2\\ b_3\end{pmatrix}$$

where $C$ represents the constant $\mathrm{x} \cdot n$ for any $\mathrm{x} \in f_j$.

Let us also assume in this example that the edge $f_j$ can be parametrised by the couple $(x,\, 1-x)$, with the coordinate $x$ spanning $[0,\, 1]$. The variable $y$ can then be substituted above by $1-x$, and one can integrate from 0 to 1 the variable $x$ to skim the edge's path. The above matrix reduces then to

$$(\Sigma) = \left(\begin{array}{cc|c|ccc} \frac{n_x}{2} & 0 & \frac{n_x}{3} & \frac{n_x}{4} & \frac{n_x}{5} & \frac{n_x}{6} \\ 0 & \frac{n_y}{2} & \frac{n_y}{3} & \frac{n_x}{12} & \frac{n_x}{30} & \frac{n_x}{60} \\ \hline n_x & n_y & C & Cx_{im} & Cx_{im}^2 & Cx_{im}^3 \\ \hline \frac{n_x}{2} & \frac{n_y}{2} & \frac{C}{2} & \frac{C}{3} & \frac{C}{4} & \frac{C}{5} \\ \frac{n_x}{3} & \frac{n_y}{3} & \frac{C}{3} & \frac{C}{4} & \frac{C}{5} & \frac{C}{6} \\ \frac{n_x}{4} & \frac{n_x}{4} & \frac{C}{4} & \frac{C}{5} & \frac{C}{6} & \frac{C}{7} \end{array}\right).$$

Assuming further in this example that $n_f = +(\frac{\sqrt{2}}{2},\, \frac{\sqrt{2}}{2})$, $C = \frac{\sqrt{2}}{2}$.
The matrix thus becomes

$$(\Sigma) = \frac{\sqrt{2}}{2} \left(\begin{array}{cc|c|ccc} \frac{1}{2} & 0 & \frac{1}{3} & \frac{1}{4} & \frac{1}{5} & \frac{1}{6} \\ 0 & \frac{1}{2} & \frac{1}{3} & \frac{1}{12} & \frac{1}{30} & \frac{1}{60} \\ \hline 1 & 1 & 1 & \frac{1}{2} & \frac{1}{4} & \frac{1}{8} \\ \hline \frac{1}{2} & \frac{1}{2} & \frac{1}{2} & \frac{1}{3} & \frac{1}{4} & \frac{1}{5} \\ \frac{1}{3} & \frac{1}{3} & \frac{1}{3} & \frac{1}{4} & \frac{1}{5} & \frac{1}{6} \\ \frac{1}{4} & \frac{1}{4} & \frac{1}{4} & \frac{1}{5} & \frac{1}{6} & \frac{1}{7} \end{array}\right),$$

which is invertible. However, as in this example we decomposed the function $f$ on the canonical basis, and as the degrees of freedom' kernels of $\{\sigma_T\}$ form also a canonical basis, $(\Sigma)$ is close to the Vandermonde matrix and has a bad conditioning value of 17479. ♦

***Note.*** Let us assume that we have a generic basis of $\mathbb{H}_k(K)$ whose functions can be split into internal and normal classification. When tuning it against the above defined degrees of freedom, the following behaviour can be observed.

By definition of the normal degrees of freedom and the boundary vanishing property of internal basis functions, the contribution of any internal basis function through a normal moment vanishes. Therefore, no contribution of internal basis function will be contained in the transfer matrix and the split between internal and normal degrees of freedom is preserved.

However, a normal function given through the degrees of freedommay contribute to the definition of internal functions. Indeed, the normal basis classification does not impose that the normal functions should vanish within the cell (and neither should their internal degrees of freedom). Therefore, by definition of the internal basis functions, non - vanishing degrees of freedom values corresponding to some contribution of normal degrees of freedom will appears in the rows defining the tuned internal basis functions.

As a consequence, the split preservation is only guaranteed by the definition of generic basis functions that vanishes on the boundaries. Indeed, the normal contributions acting on the faces will be shut down on when applying the transfer matrix to the initial raw basis, whose internal functions have a vanishing normal component on the boundary. Therefore, the tuned internal basis functions vanishes on the faces of the elements. This gives the last condition on the construction given in the *Requirements* (5.4). ▲

We detail in the next section the way to have an intelligible definition of $\mathcal{P}_k$ in order to fulfil the second point of the *Requirements 5.3*.

## 5.4 Summary of the construction

For the sake of convenience, we summarize the construction of the conformal elements on general polytopes. All the elements presented here are $H(\text{div}, K)$ – conformal.

The class of discretisation spaces is given by the following definition.

$$\begin{aligned}\mathbb{H}_k(K) =& \{u \in H^1(K),\, u|_{\partial K} \in \mathcal{H}_{l_1}(\partial K),\, \Delta u \in \mathbb{Q}_{m_1}(K)\}^d \\ &+ x\,\{u \in H^1(K),\, u|_{\partial K} \in \mathcal{H}_{l_2}(\partial K),\, \Delta u \in \mathbb{Q}_{[m_2]}(K)\},\end{aligned} \tag{97}$$

with the convention that $\mathbb{Q}_{-1} = \{0\}$ and where the integers $l_1$, $l_2$, $m_1$ and $m_2$ verify:

$$\begin{aligned} m_1,\, m_2 &\geq -1 \\ l_2 &\geq -1 \\ -1 \leq l_1 &\leq 0 \\ (d-2)(l_2+1)^{d-1} \geq (d-1)l_1^{d-1} \quad &\text{or} \quad (d-1)(l_2+1)^{d-2} \geq (l_1+1)^{d-1} \end{aligned}$$

So defined, the space enjoys the following properties:

- $\dim \mathbb{H}_k(K) = n(d(l_1+1)^{d-1} + (l_2+1)^{d-1}) + ((m_1+1)^d + (m_2+1)^d - m_2^2)$.
- For all $q \in \mathbb{H}_k(K)$, $q \cdot n|_{\partial K} \in \mathcal{H}_{\max\{l_1, l_2\}}(\partial K)$.
- The two natural subspaces are in direct sum.

Thus, $H(\text{div}, K)$ – conformal elements can be set up by using the following definition of degrees of freedom.

**Internal degrees of freedom:**

$$\sigma(q) \mapsto \int_K q \cdot p_k \,\mathrm{d}x, \qquad \text{for all } p_k \in \mathcal{P}_k, \tag{98}$$

for any space $\mathcal{P}_k$ that can be defined either as any finite polynomial or subspace of Poisson's functions having for dimension $(m_1+1)^d + (m_2+1)^d - m_2^d$. We only ask any projection of a function living in $\mathbb{H}_k(K)$ to be not identically null and the projection values to be by two different. For this, it is enough to make sure that if the order of the space is $k$, there is no projection onto a polynomial of degree bigger than $k+1$. The same goes analogously for any Poisson's subspace $\mathcal{P}_k \subset \mathbb{H}_k(K)$.

**Normal degrees of freedom:**

| Configuration | Ia | Ib | IIa | IIb |
|---|---|---|---|---|
| Low order representation | $\int_{f_i} q_{x_j} n_{ix_j}$<br>$\forall j \in [\![1,\, d]\!]$<br>$\forall i \in [\![1,\, n]\!]$ | $q_{x_i}(x_{im}) n_{ix_j}$<br>$\forall j \in [\![1,\, d]\!]$,<br>$\forall i \in [\![1,\, n]\!]$,<br>$x_m$ midpoint | Inherited from the highest order representation | Inherited from the highest order representation |
| Representation of elements in $A_k \cap B_k\|_{\partial K}$ | $\int_{\partial K} q_{x_i}\, n_{x_i}\, x^{\beta_j}$<br>$\boldsymbol{\int_{\partial K} q \cdot v\, x^{\alpha_j}}$<br>$\forall i \in [\![1,\, d]\!]$<br>$\forall 1 \leq \|\beta_j\| \leq l_1$,<br>$\|\alpha_j\| = l_1 + 1$,<br>$v \neq n_i$ | $\int_{\partial K} q_{x_i}\, n_{x_i}\, x^{\beta_j}$<br>$\boldsymbol{\int_{\partial K} q(x) \cdot v\, x^{\alpha_j}}$<br>$\forall i \in [\![1,\, d]\!]$<br>$\forall 1 \leq \|\beta_j\| \leq l_1$,<br>$\|\alpha_j\| = l_1 + 1$,<br>$v \neq n_i$ | $\int_{\partial K} q_{x_i}\, n_{x_i} p_k$<br>$\int_{f_j} q_{x_i}\, n_{jx_i} x_i^{l_1+1}$<br>$\forall i \in [\![1,\, d]\!]$,<br>$\forall j \in [\![1,\, n]\!]$,<br>$p_k \in \mathcal{H}_{l_1}(\partial K) \setminus \mathcal{H}_0(\partial K)$ | $\int_{\partial K} q_{x_i}\, n_{x_i} p_k$<br>$\int_{f_j} q_{x_i}\, n_{x_i} x_i^{l_1+1}$<br>$\forall i \in [\![1,\, d]\!]$,<br>$\forall j \in [\![1,\, n]\!]$,<br>$p_k \in \mathcal{H}_{l_1}(\partial K) \setminus \mathcal{H}_0$ |
| Higher orders representation: | $\boldsymbol{\int_{\partial K} q \cdot n x_j x^{\alpha_j}}$<br>$\forall l_1 < \|\alpha_j\| \leq l_2$,<br>$\forall j \in [\![1,\, d-1]\!]$ | $\boldsymbol{\int_{\partial K} q \cdot n x_j x^{\alpha_j}}$<br>$l_1 < \|\alpha_j\| \leq l_2$,<br>$\forall j \in [\![1,\, d-1]\!]$ | $\int_{\partial K} q \cdot\, n\, p_k$<br>$\boldsymbol{\int_{\partial K} q \cdot n x_j x_j^{l_2} \tilde{x}}$<br>$\forall p_k \in \mathcal{H}_{l_2}(\partial K) \setminus \mathcal{H}_{l_1}(\partial K)$<br>$\forall \tilde{x} \in \mathbb{Q}_{l_2}(\partial_j K)$ | $\int_{\partial K} q \cdot\, n\, p_k$<br>$\boldsymbol{q(x_{im}) \cdot n}$<br>$\forall p_k \in \mathcal{H}_{l_2}(\partial K) \setminus \mathcal{H}_{l_1}(\partial K)$<br>$\forall i \in [\![1,\, n]\!]$,<br>$\{x_{im}\}_m$ sampling points |
| Pick $N_I$ elements per face from the bold elements | $N_I = (l_2+1)^{d-1}$ | $N_I = (l_2+1)^{d-1}$ | $N_I = (l_1+1)^{d-1}$ | $N_I = (l_1+1)^{d-1}$ |

where the component wise line integrals sum up on the corresponding coordinate $x_i$ and the integrals considering the normal components $q \cdot n$ sum up along the path $\gamma(x)$.

***Remark.*** The indication "inherited from the highest order representation" connects with the first and second points of the *Remark 8.* ▲

We adopt the same representation convention than with the classical $RT_k(K)$ elements, that is:

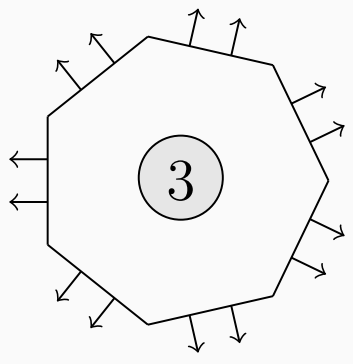

Fig. 57: Representation of the element for a septagon.

***Remark.*** *Disclaimer* Due to the definition of component - wise degrees of freedom, the construction method presented here only works for shapes whose edges vertices are not aligned with the origin and whose edges are not collinear with any of the axis. However, once a reference element is built, a transformation to those critical elements is always possible. ▲

## 5.5 Meaning of the choice of the orders

**Meaning of the orders** In order to help with the choice of the space $\mathbb{H}_k(K)$, we detail the impact of the coefficients $l_1$, $l_2$, $m_1$ and $m_2$ on the obtained discretisation.

• To begin with, $l_1$ describes the part of polynomial projection on the boundary which is not conserved through the gradient operator. Thus, setting as required $l_1 = 0$ makes then its gradient contribution vanish, which is welcomed. It also allows a simple definition of normal degrees of freedom.

• The value of $m_1$ drives the discretisation quality of the inner cell. As it does not link to any particular feature, one is free to choose it. The variable $m_2$ describes the highest inner cell discretisation quality. It will be preserved through the gradient operator. Note that setting $m_2 < m_1$ is not recommended since the benefit obtained by the subspace $xB_k$ would be absorbed by $A_k$ when applying the gradient.

• Lastly, the value of $l_2$ describes the boundary regularity that will be preserved through the gradient and can be chosen freely.

**Selection of combination of suitable orders** In practice, we define the four coefficients with respect to a single $k$ to be able to define several discretisation qualities of a space having a predefined structure. We therefore set the four coefficients through an affine relation of the type $l_0 = ak + b$ where $a$ and $b$ are two constant integers. Let us now derive the impact of their simultaneous definition.

• *Combination between $m_1$ and $l_1$* When prescribing in $A_k$ values such that $l_1 >> m_1$, one prescribes severe oscillations on the boundary while preserving a complete smoothness of the Laplacian within the domain. Indeed, there will be functions satisfying $\Delta u = 0$ with $u|_{\partial K} \in \mathcal{H}_{l_1}$, $l_1 >> 1$. This may end up to weaken the discretisation in a neighbourhood of the boundary. The same general discussion applies for the relation between $m_2$ and $l_2$.

• *Further combination between $m_2$ and $l_2$* Reversely, prescribing $m_2 >> l_2$ in $B_k$ allows functions to have reckless Laplacian within the domain and be totally vanishing on the boundary. If it is convenient to design the internal moments, the constraints on the faces appear to be brutal.

• *Combination of $m_1$ and $m_2$* It is strongly advised to set $m_1 = m_2$ to avoid holes in the projection space due to the lack of intermediate order polynomials.

• *A note on $l_1$* If $l_1$is such that $b_1 > 0$, than the lowest order space will have $db$ supplementary normal degrees of freedom per face that will require special definition. Therefore, a high value of b is not recommended.

So far, no particular restriction on the interdependency between $l_1$ and $l_2$. They will come when constructing the degrees of freedom . However, some important remark has to be done in the case $l_1 > l_2$. We assume for the sake of convenience that we are in the case $m_1 \leq m_2$. The block construction implies that the case $m_2 < m_1$ is analogous as their internal dimension is the same. There, one clearly has $x\mathbb{Q}_{l_2} \subset (\mathbb{Q}_{l_1})^d$. Indeed, let us denote by $\zeta_i$ the circular permutation of the vector $(l_2 + 1, l_2, \cdots, l_2)$. Then, all the elements in $x\mathbb{Q}_{l_2}$ are part of $\{\mathbb{P}_{\zeta_i(l_2+1,l_2,\cdots,l_2)}\}_{i\in[\![1,d]\!]} \subset \mathbb{Q}_{l_2+1} \subset \mathbb{Q}_{l_1}$. Therefore, since $l_2 < l_1$, for any element $q$ of $x\mathbb{Q}_{l_2}$ one can find vectors in $(\mathbb{Q}_{l_1})^d$ for which some combination will reduce to $q$. Thus, there is no added contribution from the constraint $\mathbb{Q}_{l_2}$ in the space $B_k$. If the solutions are still properly generated within the cell thanks to the structure of $B_k$, the definition of the normal degrees of freedom 6.1.2 will be, when possible, counterintuitive.

# 6 An example in two dimensions

We detail two possible definitions of spaces. First, we detail an example of a space lying in the previously presented framework, enjoying a parallel with the Raviart – Thomas setting from the order $k = 1$ on. In second time, we present an example of a reduced space enjoying a parallel with the Raviart – Thomas setting at any order.

## 6.1 An example of a classical space

Let us first detail a possible definition of elements living within the space $\mathbb{H}_k(K)$ built through the coefficients $l_1 = 0$, $l_2 = k$, $m_1 = k-1$ and $m_2 = k-1$ for any $k \in \mathbb{N}$ before building their corresponding basis functions. The two dimensional shape $K$ of the element is left free.

The considered discretisation space then reads

$$\begin{aligned}\mathbb{H}_k(K) =&\{u \in H^1(K),\, u|_{\partial K} \in \mathcal{H}_0(\partial K),\, \Delta u \in \mathbb{Q}_{k-1}(K)\}^2 \\ &+ \mathrm{x}\,\{u \in H^1(K),\, u|_{\partial K} \in \mathcal{H}_k(\partial K),\, \Delta u \in \mathbb{Q}_{[k-1]}(K)\},\end{aligned} \tag{99}$$

while the elements of interest $E_k = (K, \mathbb{H}_k(K), \{\sigma\})$ are defined by the following sets of degrees of freedom.

**Internal degrees of freedom:**

$$\sigma(q) \mapsto \int_K q \cdot p_k \,\mathrm{d}x, \qquad \text{for all } p_k \in \mathcal{P}_k \tag{100}$$

We chose identically the symmetric space

$$\mathcal{P}_k = \mathbb{P}_{k,\,k-1} \times \mathbb{P}_{k-1,\,k} \setminus \left(\mathbb{P}_{[k],\,[k-1]} \times \mathbb{P}_{[k-1],\,[k]}\right) \cup \left\{ (x,\, y)^T \mapsto \begin{pmatrix} x^k y^{k-1} \\ x^{k-1} y^k \end{pmatrix} \right\}$$

for every configuration. One can observe that the internal projection space is less refined than the one set on the edges. This is not bothersome as the impact of the divergence within the cell is less dramatic. Note also that in practice, for defining the projections (100) one works with the canonical basis of $\mathcal{P}_k$.

**Normal degrees of freedom:**

| Configuration | Ia | Ib | IIa | IIb |
|---|---|---|---|---|
| Low order representation | $q_{x_i}(x_{im})n_{ix_j}$ <br> $\forall j \in [\![1, 2]\!]$, <br> $\forall i \in [\![1, n]\!]$ | $\int_{f_j} q_{x_i} n_{jx_i}$ <br> $\forall i \in [\![1, 2]\!]$ <br> $\forall j \in [\![1, n]\!]$ | $\int_{\partial K} q \cdot n$ | $q(x_{im}) \cdot n$ <br> $\forall i \in [\![1, n]\!]$ |
| Representation of elements in $A_k \cap B_k\|_{\partial K}$ | $\int_{\partial K} q \cdot v\, x$ | $\int_{\partial K} q(x) \cdot v\, x$ | $\int_{f_j} q_{x_i}\, n_{x_i} x_i$ <br> $\forall i \in [\![1, 2]\!]$, <br> $\forall j \in [\![1, n]\!]$ | $\int_{f_j} q_{x_i}\, n_{x_i} x_i$ <br> $\forall i \in [\![1, 2]\!]$, <br> $\forall j \in [\![1, n]\!]$ |
| Higher orders representation | $\int_{\partial K} q \cdot n x^{\alpha+1}$ <br> $\forall 0 < \alpha \leq k$ | $\int_{\partial K} q \cdot n x^{\alpha+1}$ <br> $0 < \alpha \leq k$ | $\int_{\partial K} q \cdot n\, p_k$ <br> $\forall p_k \in \mathcal{H}_k(\partial K) \setminus \mathcal{H}_0(\partial K)$ | $\int_{\partial K} q \cdot n\, p_k$ <br> $\forall p_k \in \mathcal{H}_k(\partial K) \setminus \mathcal{H}_0(\partial K)$ |

Here, $x_m$ stands for the middle point on each edge. Furthermore, in those examples both the component wise line integrals and the integrals considering the full normal component $q \cdot n$ sum up along the path $\gamma(x)$.

***Remark.*** As pointed out in the *Remark 8* of the previous section, the term $\int_{\partial K} q \cdot n$ can be replaced by $\int_{\partial K} q \cdot n\, x^{k+1} \mathrm{d}\gamma(x)$ without changing the element. ▲

### 6.1.1 Properties

By a straightforward application of the previous section, we derive the dimension of the space being

**Property 6.1** Dimension

$$\dim \mathbb{H}_k(K) = n(k+3) + 2k(k+1) - \mathbb{1}_{k>0} \tag{101}$$

***Proof.*** By a direct application of the formula obtained in (66), we get

$$\begin{aligned}
\dim \mathbb{H}_k(K) &= n\ \left(d(l_1+1)^{d-1} + (l_2+1)^{d-1}\right) \\
&\quad + \left(d(m_1+1)^d + (m_2+1)^d - m_2^d\right) \\
&= n(2(0+1)^1 + (k+1)^1) + 2(k-1+1)^2 \\
&\quad + (k-1+1)^2 - (k-1)^2 \\
&= n(k+3) + 2k^2 + k^2 - (k^2 - 2k + 1) \\
\dim \mathbb{H}_k(K) &= n(k+3) + 2k(k+1) - \mathbb{1}_{k>0}.
\end{aligned}$$ ■

### 6.1.2 Construction of the corresponding basis functions

Let us present a possible construction for the basis of our sample space $\mathbb{H}_k(K)$ presented in the *Section 6*. The method is detailed specifically for the two - dimensional case, although an extension to general dimension is straightforward. It can also be simply adapted to any element of the class built on any space $\mathbb{H}_k(K)$. The required changes will be pointed out for the sake of completeness.

As in the two previous cases of simplicial and quadrilateral elements, this canonical basis of $\mathbb{H}_k(K)$ can be later adapted towards the designed degrees of freedom (6.1) and (100). Endowed in $\mathbb{H}_k(K)$, the resulting basis will then define the corresponding element $E_k(K)$. This adjustment on the degrees of freedom applies in the exact same way as before and will therefore not be repeated. We refer to the *Section 2.4* for more details.

The only major difference from the framework of the two *Sections 3 and 4* concerns the shape of the reference element. Here, no reference element is prescribed and one is free to use any shape fulfilling the three following criterion.

- No vertices that lie on a same edge are aligned with the origin.
- No edge is collinear with any of the axes.
- In view of using the elements $Ia$ and $Ib$, no edge's normal should be collinear with the vector $v$.

If the former two restrictions are constraints emerging only from our construction method, the ground of the later one lies in the definition of the degrees of freedom and thus cannot be avoided.

Any reference element will be denoted $K$. We adopt the notations represented in the *Figure 58*, which is adaptive towards the number of edges one wants to consider.

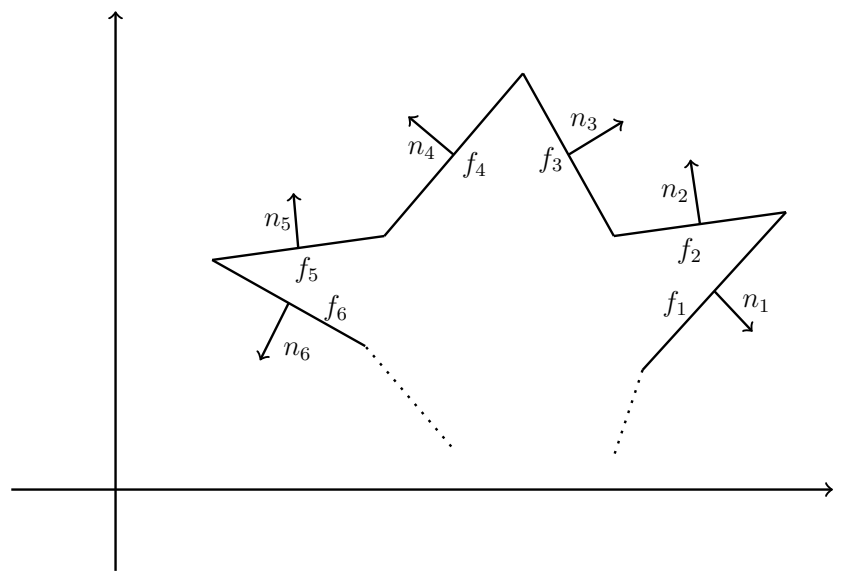

Fig. 58: Notations for any reference shape with $n$ edges, direct orientation.

***Note.*** As we do not restrict ourselves to the convex case, the transformation from the reference element to a generic one having a similar shape is not straightforward and is treated separately in the *Section 2.5.3.* ▲

Now that the framework is set, we define $\dim \mathbb{H}_k(K)$ functions that lie in $\mathbb{H}_k(K)$ and that are two by two free. As always, the basis functions are split into two categories.

- $\mathcal{N}$, a free set of functions living in $\mathbb{H}_k(K)$ whose dimension value $\dim \mathbb{H}_k(K)|_{\partial K} = n(2 \dim \mathbb{Q}_0(f) + \dim \mathbb{Q}_k(f))$. They should enforce the $H(\text{div}, K)$ – conformity.
- $\mathcal{I}$, a free set of function living in $\mathbb{H}_k(K)$ whose dimension is $\dim \mathbb{P}_k = (k+1)^2+2k+1 = 2k(k+1) - \mathbb{1}_{k>0}$. They should preserve the $H(\text{div}, K)$ – conformity.

Note that here as well, we ask the functions falling in each of the two categories to belong to the full space $\mathbb{H}_k(K)$ and not either to one of the two subspaces that are in direct sum. The constraints given by the definition of the Poisson's problems defining the subspaces prevent to define the split directly on the direct sum. Indeed, the obtained divergence conformity of the normal component on the boundary is made possible through smooth boundary conditions that are not identical for both subspaces. The freeness of $\mathcal{N} \cup \mathcal{I}$ will then come as a by product of our definitions.

**Let us first generate $\mathcal{N}$.** We want to generate a set of $\dim \mathcal{H}_k(\partial K) = n(2 \dim \mathbb{Q}_0(f) + \dim \mathbb{Q}_k(f)) = n(k+3)$ functions that are free from each other and that lie in $\mathbb{H}_k(K)$. They will later be tuned to be dual to the set of degrees of freedom (6.1). Therefore, we follow their layout and set $k+3$ functions on each edge $f_i$, $i \in [\![1, n]\!]$.

***Note.*** As $\mathcal{T}_k(\partial K)$, the space $\mathcal{H}_k(\partial K)$ is based on $\mathbb{Q}_k(f)$. Thus, and as we are on the boundary of a two dimensional element the edge - wise discretisation space $\mathcal{P}_{b,k} = \mathbb{Q}_k(f)$ coincides with $\mathbb{P}_k(f)$ for any $k \in \mathbb{N}$. ▲

We keep the same spirit as for the construction of the basis functions in the simplicial and quadrilateral cases. However, as here the two scalar subspaces $A_k$ and $B_k$ that generate the $\mathbb{H}_k(K)$ space are different, we need to define two vectors per edge belonging respectively to $\mathbb{P}_{1,0} \times \mathbb{P}_{0,1}$ and $\mathbb{Q}_0 \times \mathbb{Q}_0$. We will then stitch both to functions belonging respectively to the space $A_k$ or $B_k$ before summing the resulting vectorial functions up.

In view of later use, we recall that their Laplacian respectively matches $(x,\, y) \mapsto x^{k-1}y^{k-1}$ or have the possibility to vanish within the cell.

For any edge $f_i$, $i \in [\![1,\, n]\!]$, we define the two vectors

$$e_{1,\,i} = \begin{pmatrix} x \\ y \end{pmatrix} \quad \text{and} \quad e_{2,\,i} = \begin{pmatrix} n_{ix} \\ n_{iy} \end{pmatrix} \tag{102}$$

where $n_i = (n_{ix},\, n_{iy})^T$ is the vector normal to the edge $f_i$.

***Note.*** Those vectors can be seen as the split of the one defined in the two previous *Sections 3.5 and 4.3.3.* ▲

However, since on the boundary we are constrained by the order of the polynomial and can therefore only create $k+1$ functions in $A_k$ and one function in $B_k$, we need to use two more vectors in order to create the remaining functions. We set

$$e_{3,\,i} = \begin{pmatrix} n_{ix} + n_{iy}^2/n_{ix} \\ 0 \end{pmatrix} \quad \text{and} \quad e_{4,\,i} = \begin{pmatrix} 0 \\ n_{iy} + n_{ix}^2/n_{iy} \end{pmatrix}. \tag{103}$$

As we excluded the cases of edges that are collinear to the axes, those vectors are always defined.

***Note.***

• One could also have used other coefficients in their definition, but doing so ensures that $e_{3,\,i} \cdot n_i$ and $e_{4,\,i} \cdot n_i$ scales to one.

• The extension of (103) to a general dimension is made by defining $d$ vectors whose only non - zero coordinates are in position $j$, $j \in [\![1,\, d]\!]$ and value $n_{ix_j} + \sum_{m \neq j} n_{ix_m}^2 / n_{ix_j}$. ▲

The $(k+1)$ functions per edge belonging to $A_k$ or x$B_k$ are then generated by Poisson's functions whose boundary conditions are prescribed edge - wise through a set of Lagrangian functions $\{l_{i,\,m}\}_{\substack{m \in [\![1,\, k+1]\!] \\ i \in [\![1,\, n]\!]}}$. To construct them, we define $f_{i,\,m} \in L^2(K)$ as the solution of

$$\begin{cases} \Delta u = x^{k-1}y^{k-1} \\ u|_{\partial K} = l_{i,\,m} \mathbb{1}_{f_i} \end{cases} \tag{104}$$

and $g_i \in L^2(K)$ as the harmonic solution of

$$\begin{cases} \Delta u = 0 \\ u|_{\partial K} = 2\mathbb{1}_{f_i}. \end{cases} \tag{105}$$

For more details regarding the construction of the set $\{l_{i,m}\}_{i,m}$, we refer to the *Section 4.3.3* from which a generalisation to $n$ edges immediately follows.

***Note.*** One can also naturally extend the generation method for the Lagrangian functions from a single Lagrangian set $\{l_m\}_{m\in[\![1,\,k+1]\!]}$ built on the segment [0, 1] as it was done in the *Paragraph 9.* ▲

We then define the normal basis functions as being:

$$(x,\, y) \mapsto e_{i,1} f_{i,m} + e_{i,2} g \tag{106a}$$

$$(x,\, y) \mapsto e_{i,1} f_{i,1} - e_{i,3} g \tag{106b}$$

$$(x,\, y) \mapsto e_{i,1} f_{i,k+1} - e_{i,4} g \tag{106c}$$

***Remark.*** The constant assigned on the boundary in the problem (105) can be anything except one. Indeed, in the case $k = 0$ the Lagrangian functions reduces to the identically one function. There, it would exist for each edge $f_i$ an equation of the set (104) for which the solution matches the one of (105) on the boundary. The constructed functions would then differ only by the factor $(x,\, y) \mapsto (x{-}n_i,\, x,\, y{-}n_{iy})^T$, and since we saw in the *Section 4* that the set (49) is not free, the set of basis functions (106) would not be free either. ▲

Finally, we set:

$$\mathcal{N} = \{e_{i,1} f_{i,m} + e_{i,2} g\}_{\substack{i\in[\![1,\,n]\!] \\ m\in[\![1,\,k+1]\!]}} \cup \{e_{i,1} f_{i,m} - e_{i,l} g\}_{\substack{i\in[\![1,\,n]\!] \\ (m,\,l)\in\{(1,3) \\ ,(k+1,\,4)\}}}. \tag{107}$$

***Remark.*** *Generalisation to other spaces* The only modifications that are required to derive a general construction for any $\mathbb{H}_k(K)$ space in two dimension read:

- Define (104) with a second member as $(x,\, y) \mapsto x^{m_2} y^{m_2}$.
- Define the boundary condition of (104) through Lagrangian functions of degree $l_2$.
- By definition of $RT_k(K)$, we have $d(l_1{+}1)^{d-1}$ misc basis functions to set up. As here $d = 2$, define the $2(l_1{+}1)$ misc basis functions as (106b) and (106c) by extending the set $m \in \{1,\, k+1\}$ to $\{m \in [\![1,\, k+1]\!]\}$ such that its dimension is $2(l_1+1)^{2-1}/2$. As we already have $l_2+1$ functions in $A_k$ and since we duplicate them, we can define up to $2(l_2+1)$ functions, which is more than we need. If we have an odd number of functions, the choice of either $e_3$ or $e_4$ in the definitions (106b) of the misc function which will not be symmetrised is up to the user.

- The extension to $d$ dimensions is doable by defining multidimensional functions enjoying a Lagrangian property on each edge. ▲

***Remark.*** Let us point out the reasons of the second shape restriction we made in the definition of the reference element, namely that no edge of the element has its vertices aligned with the origin.

We assume that the reference element comprises some edge whose verticies are aligned with the origin. Then, for the same reasons as presented in the simplicial case and emphasized in the *Figure 16* page 68, the term $(x,\, y) \mapsto (x,\, y)^T f_{i,j}$ would vanish for any $j \in [\![1,\, k+1]\!]$. We would then be left with the term driven solely by $e_2$, which is not dependent on $j$. The set of functions would then not be free on the edges for any $k \geq 1$. Furthermore, it would shut down the desired order on that edge and we would be left on this edge with only a constant projection of the discretised quantity. ▲

Just before detailing some properties of $\mathcal{N}$, let us derive the following straightforward property that will be needed in later proofs.

**Lemma 6.2**

Let $u$ be solution to the following Poisson's problem.

$$\begin{cases} \Delta u = p_l \\ u|_{\partial K} = q_m, \end{cases}$$

where $p_l \in \mathbb{Q}_l(K)$ and $q_m \in \mathbb{Q}_m(f)$ for some $l,\, m \in \mathbb{N}$. Then, for any $\alpha \in \mathbb{R}$, the function $(x,\, y) \mapsto \alpha u(x,\, y)$ fulfils

$$\begin{cases} \Delta(\alpha u) \in \mathbb{Q}_l(K) \\ (\alpha u)|_{\partial K} \in \mathbb{Q}_m(f). \end{cases}$$

***Proof.*** It comes directly from the linearity of the Laplacian. Indeed,

as $\alpha \in \mathbb{R}$ is a constant, we can write

$$\begin{cases} \Delta(\alpha u) = \alpha\Delta(u) = \alpha p_l \in \mathbb{Q}_l(K) \\ (\alpha u)|_{\partial K} = \alpha q_m \in \mathbb{Q}_m(f). \end{cases}$$

■

Let us now state some properties of $\mathcal{N}$.

**Property 6.3** Dimension of $\mathcal{N}$

$$\dim \mathcal{N} = n(k+3)$$

***Proof.*** By simply reading out the construction (107), we obtain

$$\begin{aligned}\dim \mathcal{N} &= \dim\{\{e_{i,1}f_{i,m} + e_{i,2}g\}_{\substack{i\in[\![1,n]\!]\\ m\in[\![1,k+1]\!]}} \cup \{e_{i,1}f_{i,m} - e_{i,l}g\}_{\substack{i\in[\![1,n]\!]\\ (m,l)\in\{(1,3)\\ ,(k+1,4)\}}}\} \\ &= n(k+1) + n(1+1) \\ &= n(k+3).\end{aligned}$$

Indeed, by construction the Lagrangian functions are two by two different. Therefore, as the indicator function prevents any interaction between the basis functions emerging from different edges, there is no identical functions in the set $\mathcal{N}$ and its dimension is preserved.

■

**Proposition 6.4**

For any $q$ belonging to $\mathcal{N}$, $q$ belongs to $\mathbb{H}_k(K)$.

***Proof.*** Let $q$ belong to $\mathcal{N}$. Then, by construction $q$ takes one of the following forms.

$$q = e_{i,1}f_{i,m} + e_{i,2}g_i \tag{108a}$$

$$\text{or} \quad q = e_{i,1}f_{i,1} - e_{i,3}g_i \tag{108b}$$

$$\text{or} \quad q = e_{i,1}f_{i,k+1} - e_{i,4}g_i \tag{108c}$$

for some index $i \in [\![1,\, n]\!]$ and some index $m \in [\![1,\, k+1]\!]$.

- We first assume that $q$ is of the former form, that is

$$q = \begin{pmatrix} x \\ y \end{pmatrix} f_{i,m} + \begin{pmatrix} n_{ix} \\ n_{iy} \end{pmatrix} g_i.$$

By construction, we know that $f_{i,m}$ satisfies

$$\begin{cases} \Delta f_{i,m} = x^{k-1}y^{k-1} \in \mathbb{Q}_{[k-1]}(K) \\ f_{i,m}|_{\partial K} = l_{i,m}\mathbb{1}_{f_i} = \begin{cases} l_{i,m} \in \mathbb{Q}_k(f) & \text{if x} \in f_i \\ 0 \in \mathbb{Q}_0(f) \subset \mathbb{Q}_k(f) & \text{if x} \notin f_i \end{cases} \in \mathcal{H}_k(\partial K) \end{cases}$$

Thus, $f_{i,m}$ belongs naturally to $B_k$ and $(x,\, y)^T f_{i,m} \in \mathrm{x}B_k$. Similarly, one has by construction that $g$ satisfies

$$\begin{cases} \Delta g = 0 \in \mathbb{Q}_{-1}(K) \subset \mathbb{Q}_{k-1}(K) \\ g|_{\partial K} = 2\mathbb{1}_{f_i} = \begin{cases} l_{i,m} \in \mathbb{Q}_0(f) & \text{if } x \in f_i \\ 0 \in \mathbb{Q}_0(f_j) & \text{if } x \in f_j \neq f_i \end{cases} \in \mathcal{H}_0(\partial K) \end{cases}$$

for any $k \geq 0$, which fits our setting. Therefore, $g \in A_k$. Lastly, as $n_{ix}$ and $n_{iy}$ are two real constants, we have by the *Property 6.2* that $n_{ix}g$ and $n_{iy}g$ both belong to $A_k$. Thus, $(n_{ix},\, n_{iy})^T g$ belongs to $(A_k)^2$, and $q \in \mathbb{H}_k(K)$.

• Let us now treat the case where $q$ writes under the form

$$q = \begin{pmatrix} x \\ y \end{pmatrix} f_{i,1} - \begin{pmatrix} n_{ix} + n_{iy}^2/n_{ix} \\ 0 \end{pmatrix} g_i$$

We get in the exact same way that $\mathrm{x}f_{i,1} \in B_k$. Furthermore, as the term $(n_{ix} + n_{iy}^2/n_{ix})$ belongs to $\mathbb{R}$, we also retrieve similarly $-(n_{ix} + n_{iy}^2/n_{ix})g_i \in A_k$. Lastly, we notice that the identically null function $(x,\, y) \mapsto 0$ belongs to $A_k$ for any $k \geq -1$as being solution to

$$\begin{cases} \Delta u \equiv 0 \in \mathbb{Q}_{k-1}(K) \\ u|_{\partial K} \equiv 0 \in \mathcal{H}_0(\partial K). \end{cases}$$

Thus, $-(n_{ix} + n_{iy}^2/n_{ix},\, 0)^T g_i \in A_k$, and $q \in \mathbb{H}_k(K)$.

• The last case where $q$ writes $q = e_{i,1}f_{i,k+1} - e_{i,4}g_i$ is treated in the exact same way as the second point, and concludes the proof.

■

**Property 6.5** Pointwise vanishing property

Let us consider any function $q \in \mathcal{N}$ defined from any edge $f_i$, $i \in [\![1,\, n]\!]$. Then, for any sampling point $x_l \in f_j$, $l \in [\![1,\, k+1]\!]$ that was used to generate the $k+1$ Lagrangian functions $\{l_{j,m}\}_{m\in[\![1,\,k+1]\!]}$ from some edge $f_j$, $j \in [\![1,\, n]\!]$, we obtain the following property:

$$q(x_l) \cdot n = \delta_{ij}(c\,\delta_{ml} \pm 2),$$

where $c$ represents a constant depending on the angle formed between the edge $i$ and the axes.

***Proof.*** Let $q$ belong to $\mathcal{N}$. We recall that by construction $q$ takes one of the forms presented in (108).

• Le us first assume that $q$ is of the form (106a):

$$q = \begin{pmatrix} x \\ y \end{pmatrix} f_{i,\,m} + \begin{pmatrix} n_{ix} \\ n_{iy} \end{pmatrix} g_i.$$

Then, one can write on any edge $f_j$ belonging to $\partial K$:

$$\forall (x,\, y) \in f_j,\; q \cdot n|_{\partial K} = (x n_{jx} + y n_{jy}) f_{i,\,m} + (n_{ix} n_{jx} + n_{iy} n_{jy}) g_i.$$

Plugging the values of $f_{i,\,m}$ and $g_{i,\,m}$ on the boundary into the previous equation, we get

$$q \cdot n|_{\partial K} = (x n_{jx} + y n_{jy}) l_{i,\,m} \mathbb{1}_{f_i} + 2(n_{ix} n_{jx} + n_{iy} n_{jy}) \mathbb{1}_{f_i}.$$

Recalling that $l_{i,\,m}$ is a Lagrangian function generated from a set of points $\{\mathrm{x}_{i,\,l}\}_{\substack{l\in[\![1,\,k+1]\!] \\ i\in[\![1,\,n]\!]}}$ all belonging to the edge $f_i$, we have for any edge $f_j \in \partial K$:

$$\begin{aligned}
q(\mathrm{x}_{j,\,l}) \cdot n|_{f_j} &= (x_{j,\,l} n_{jx} + y_{j,\,l} n_{jy}) l_{i,\,m}(\mathrm{x}_{j,\,l}) \mathbb{1}_{f_i} + 2(n_{ix} n_{jx} + n_{iy} n_{jy}) \mathbb{1}_{f_i} \\
&= \begin{cases} (x_{j,\,l} n_{jx} + y_{j,\,l} n_{jy}) l_{i,\,m}(\mathrm{x}_{j,\,l}) \times 0 + 2(n_{ix} n_{jx} + n_{iy} n_{jy}) \times 0 & \text{if } i \neq j \\ \underbrace{(x_{j,\,l} n_{jx} + y_{j,\,l} n_{jy})}_{=C} \delta_{ml} \times 1 + 2\, \underbrace{(n_{ix} n_{jx} + n_{iy} n_{jy})}_{=1} & \text{if } i = j \end{cases} \\
&= \begin{cases} 0 & \text{if } i \neq j \\ C\delta_{ml} + 2 & \text{if } i = j \end{cases}
\end{aligned}$$

for some constant $C$ dependent on the position and angle of the edge $f_j$ with respect to the axes. Thus,

$$q(\mathrm{x}_{i,\,l}) \cdot n|_{f_j} = \delta_{ij}(c\delta_{ml} + 2).$$

• The two other cases are treated equally. Indeed, the only change to make in the developments above arises in the definition of the constant stitched to $g_i$. Meaningly, one has to replace the coefficient $(n_{iy},\, n_{iy})^T$ by $-(n_{ix}+n_{iy}^2/n_{ix},\, 0)^T$ or $-(0,\, n_{iy}+n_{ix}^2/n_{iy})^T$ in all of the above equations. And since by definition those terms still scales to one with respect to the normal vector to the edge $f_j$, we get

$$\begin{aligned} q\cdot n|_{f_j}(\mathrm{x}_{j,l}) &= (x_{j,l}n_{jx} + y_{j,l}n_{jy})l_{i,m}(\mathrm{x}_{j,l})\mathbb{1}_{f_i} \\ &\quad - 2(n_{ix}^2 + n_{ix}n_{iy}^2/n_{ix} + 0)\mathbb{1}_{f_i} \\ &= \begin{cases} 0 & \text{if } i \neq j \\ Cl_{i,m}(\mathrm{x}_{i,l}) - 2 & \text{if } i = j \end{cases} \\ &= \begin{cases} 0 & \text{if } i \neq j \\ C\delta_{ml} - 2 & \text{if } i = j \end{cases} \\ &= \delta_{ij}(c\,\delta_{ml} - 2) \end{aligned}$$

• The last case being analogous to the previous one, we have:

$$\forall q \in \mathcal{N}, q(\mathrm{x}_{j,l})\cdot n|_{f_i} = \delta_{ij}(C\delta_{ml} \pm 2).$$ ■

**Property 6.6**

The functions of $\mathcal{N}$ are linearly independent.

***Proof.***

• Let us show first that the set of functions (106a) is free. We know that on each edge $f_i$, the functions $\{f_{i,m}\}_{m\in[\![1,\,k+1]\!]}$ are linearly independent as they are forming a Lagrangian set $\{l_{i,m}\}_{m\in[\![1,\,k+1]\!]}$ on the boundary. Furthermore, we have by construction that the set $\{f_{i,m}\}_{m\in[\![1,\,k+1]\!]}$ comprises functions that are null on all the edges but $f_i$. Thus, as those observations hold for any edge $f_i$, $i \in [\![1,\, n]\!]$, the set $\{f_{i,m},\, i\in[\![1,\, n]\!],\, m\in[\![1,\, k+1]\!]\}$ is naturally free on the boundary. And as the considered elements are not reduced to a point, the vector $\mathrm{x}\in K$ is not equally null. Thus, the set $\{(x,\, y)^T f_{i,m}\}_{\substack{i\in[\![1,\,n]\!] \\ m\in[\![1,\,k+1]\!]}}$ is free on the boundary.

Similarly, the functions $\{g_i\}_{i\in[\![1,\,n]\!]}$ are free from each other as each of them are vanishing on all the edges but $f_i$. And since $e_{i,2} \not\equiv 0$ lies in $\mathbb{R}^2$ for any $i\in[\![1,\, n]\!]$, the set $\{e_{i,2}g_i\}_{i\in[\![1,\,n]\!]}$ is also free. Furthermore, as in the expression (106a) the term $e_{i,2}g_i$ is only translating the expression of $f_{i,m}$ up or down on the boundary, it does not impact the previously described behaviour and the set of functions (106a) is free.

• Let us now show that (106b) is free with the set (106a). To see this, it is enough to show that no linear combination of (106a) can generate (106b). By the Lagrangian properties the function $f_{i,1}$ enjoys on a boundary, the only way to generate the function $f_{i,1}$ is to use the function $f_{i,1}$ itself. Therefore, the term $e_{i,1}f_{i,1}$ can be generated only using the term $e_{i,1}f_{i,1} + e_{i,2}g_i$.

Thus, to generate (106b) on has to remove $e_{i,2}g_i + e_{i,3}g_i$ from the previous quantity. However, using functions of (106a) to remove them would in the same time make use of some functions $f_{m,1}$ for $m \neq i$ and thus break the relation on the $f_{i,j}$. Thus, we only have to check if $e_{i,1}f_{i,1} + e_{i,2}g_i$ is not collinear with $e_{i,1}f_{i,1} - e_{i,3}g_i$. As $e_{i,3}g_i \neq e_{i,1}$, no relation satisfying $Ae_{1,i}f_{i,1} + Ae_{2,i} = Bf_{i,1} - Be_{i,3}g_i$ can be found for any $A$, $B \in \mathbb{R}$. Thus, {(106a), (106b)} is free.

• Lastly, we can see in a similar way that the only mean to construct (106c) is to make the use of $f_{i,k+1}$ itself without any other contribution, whose only appearance is in (106a). For the same reason as above, the vector $e_{i,j}f_{i,k+1} + e_{i,2}g_i$ if not collinear with (106c), and the set $\mathcal{N}$ is free.

■

**Let us now generate $\mathcal{I}$** We now want to generate $2k(k+1) - \mathbb{1}_{k>0}$ functions lying in $\mathbb{H}_k(K)$ that are free from each other and free with $\mathcal{N}$. They will later be tuned to be dual to the degrees of freedom (98).

***Note.*** For $k = 0$, $2k(k+1) - \mathbb{1}_{k>0} = 0$. The set $\mathcal{I}$ is then designed as an empty set, and there exists no internal functions. ▲

As the space is built as a direct sum and as our construction makes use of Poisson's functions, we can simply define the internal functions by reading out the structure of the space. Indeed, making use of the superposition theorem, one can extract the constraints on the Laplacian of function belonging to either $A_k$ or $B_k$ while imposing homogeneous Dirichlet boundary condition. The corresponding basis functions can then be constructed with those functions as a core. In this perspective, we define $f_l$ as being the solution to

$$\begin{cases} \Delta u = x^l y^{k-1} \\ u|_{\partial K} = 0, \end{cases}$$

$g_l$ being the solution to

$$\begin{cases} \Delta u = x^{k-1}y^l \\ u|_{\partial K} = 0, \end{cases}$$

and $h_{l,m}$ being the solution to

$$\begin{cases} \Delta u = x^l y^m \\ u|_{\partial K} = 0. \end{cases}$$

Then, we define the functions

$$F_l \colon (x,\, y)^T \mapsto \begin{pmatrix} x \\ y \end{pmatrix} f_l(x,\, y) \qquad \text{for all } l \in [\![0,\, k-2]\!] \tag{109}$$

$$G_l \colon (x,\, y)^T \mapsto \begin{pmatrix} x \\ y \end{pmatrix} g_l(x,\, y) \qquad \text{for all } l \in [\![0,\, k-1]\!] \tag{110}$$

$$H^x_{l,m} \colon (x,\, y)^T \mapsto \begin{pmatrix} 1 \\ 0 \end{pmatrix} h_{l,m}(x,\, y) \qquad \text{for all } l,\, m \in [\![0,\, k-1]\!] \tag{111}$$

$$H^y_{l,m} \colon (x,\, y)^T \mapsto \begin{pmatrix} 0 \\ 1 \end{pmatrix} h_{l,m}(x,\, y) \qquad \text{for all } l,\, m \in [\![0,\, k-1]\!] \tag{112}$$

to finally set

$$\mathcal{I} = \{F_l\}_{l\in[\![0,\,k-2]\!]} \cup \{G_l\}_{l\in[\![0,\,k-1]\!]} \cup \{H^x_{l,m}\}_{l,\,m\in[\![0,\,k-1]\!]} \cup \{H^y_{l,m}\}_{l,\,m\in[\![0,\,k-1]\!]}.$$

We can immediately derive the following properties.

**Property 6.7** Dimension

For any $k > 0$, $\mathcal{I}$ enjoys the following dimension.

$$\dim \mathcal{I} = 2k(k+1) - 1$$

***Proof.*** We directly have by construction that $\mathcal{I}$ is made of the sets of functions $\{F_l\}_l$, $\{G_l\}_l$, $\{H^x_{l,m}\}_{l,m}$ and $\{H^y_{l,m}\}_{l,m}$ providing respectively $k-1$, $k$ and $k^2$ functions. The dimension of the set is then directly given by $k-1+k+k^2+k^2 = 2k(k+1)-1$.

■

**Proposition 6.8**

For any $q$ belonging to $\mathcal{I}$, $q$ belongs to $\mathbb{H}_k(K)$.

***Proof.*** Any function $q$ of $\mathcal{I}$ can be written either on the form of $F_l$, $G_l$ or $H_l$. We derive the above property case by case.

- If $q$ belongs to the subset $\{F_l\}$, then $q$ is of the form

$$(x,\, y)^T \mapsto \begin{pmatrix} x \\ y \end{pmatrix} f_l(x,\, y) \qquad \text{for all } l \in [\![0,\, k-2]\!]$$

By construction, $f_l$ belongs to the space $B_k$. Therefore, $q \in \mathrm{x}B_k$ and $q \in \mathbb{H}_k(K)$.

- As $g_l$ also belongs to the space $B_k$, $q$ belongs to the subset $B_k$ in the exact same way than when constructed from $G_l$.

- If $q$ belongs to the subset $\{H_{lx}\}$, then $q$ is of the form

$$(x,\, y)^T \mapsto \begin{pmatrix} 1 \\ 0 \end{pmatrix} h_{l\,m}(x,\, y)$$

with $h_{l,m}$ belonging by construction to $A_k$. And since $(1,\, 0)^T \in \mathbb{R}^2$, $q$ belongs to $(A_k)^2$, and thus to $\mathbb{H}_k(K)$. The exact same applies when $p$ has the form of $\{H_{ly}\}$.

■

**Proposition 6.9**

For any $q$ belonging to $\mathcal{I}$, one has

$$q \cdot n|_{\partial K} = 0$$

***Proof.*** By construction, one has $h_{l,m}$, $f_l$ and $g_l$ vanishing on the boundary of $K$. Therefore, any function of the shape $Af_l + bg_l$ will vanish on the boundary identically for any constant vectors $A, b, \in \mathbb{R}$. As by construction $q \cdot n = (an_{ix} + bn_{iy})h_l$, or $(xn_{ix} + yn_{iy})f_l$ or $(xn_{ix} + yn_{iy})g_l$, and since on any boundary the term $(xn_{ix} + yn_{iy})$ is constant, we get $q \cdot n|_{\partial K} = 0$.

■

**Proposition 6.10**

The set $\mathcal{I}$ is free strictly within the polygon $K$

***Proof.*** We want to show that the set $\mathcal{I}$ is free.

• Let us first show that the set $\{H^x_{l,m},\, H^y_{l,m}\}_{lm}$ is free. To see it, we start by recalling that the fundamental function $h_{l,m}$ of this set verifies

$$\begin{cases} \Delta h_{l,m} = x^l y^m \\ h_{l,m}|_{\partial K} = 0 \end{cases}$$

for any $l,\, m \in [\![0,\, k-1]\!]$ such that any pair $(l,\, l)_{l\in[\![0,\, k-1]\!]}$ only appear once. Therefore, the polynomial degrees appearing in the right - hand - side of the Laplacian equalities of this set are two by two different. Since we have by *Proposition 6.2* that if $h_{l,m} \in \mathbb{Q}_{l,m}$, then $\alpha h_{l,m} \in \mathbb{Q}_{l,m}$, no linear combination of the form

$$p_{i,j} = \sum_{m,\, l=1}^{k-1} \alpha_{l,m} h_{l,m} + \sum_{l=1}^{k-1} \beta_{l,l} h_{l,l} - \alpha_{i,j} h_{i,j}$$

can have the property

$$\Delta(p_{i,j}) = x^i y^j$$

for any $i,\, j \in [\![0,\, k-1]\!]$. In other words, since the degrees of the polynomials being right hand sides of the Laplacian equations are two by two different, there cannot be any linear combination between them that would generate a polynomial having a degree not already described in the set. Thus, by linearity of the Laplacian this impossibility gets back to $h_{l,m}$. Therefore, the set $\{h_{l,m}\}_{l,m}$ is free. Thus, the set

$$\left\{ \begin{pmatrix} 1 \\ 0 \end{pmatrix} h_{l,m},\, \begin{pmatrix} 0 \\ 1 \end{pmatrix} h_{l,m} \right\}_{lm}$$

is free by construction.

• Then, for the same reason as evoked above, the sets $\{f_l\}_l$ and $\{g_l\}_l$ are also free. And since the orders present in those two subsets are different from each others, the set $\{f_l,\, g_l\}$ is also free. Furthermore, since x $\not\equiv 0$ on $K$, the set

$$\{F_l,\, G_l\}_l = \left\{ \begin{pmatrix} x \\ y \end{pmatrix} f_l,\, \begin{pmatrix} x \\ y \end{pmatrix} g_l \right\}_l$$

is equally free.

• Let us now bring those two sets together. We start by considering any function of $\{F_l,\, G_l\}_l$. As we saw before, we have

$$\Delta(x f_l) = 2\nabla \cdot f_l + x\Delta(f_l) \in x\mathbb{Q}_{[l],\,[k-1]} \subset \mathbb{Q}_{l+1,\,k-1} \subset \mathbb{Q}_{k-1}$$

and $\quad \Delta(yf_l) = 2\nabla \cdot f_l + y\Delta(f_l) \in y\mathbb{Q}_{[l],\,[k-1]} \subset \mathbb{Q}_{l,\,k} \not\subset \mathbb{Q}_{k-1}.$

Thus, $(x,\, y)^T f_l$ cannot be recast as $(xf_l,\, yf_l)^T$ such that both $(x,\, y) \mapsto xf_l$ and $(x,\, y) \mapsto yf_l$ belong to $A_k$. Furthermore, as the degree of $f_l$ are two by two different, no combination of such constructions can compensate this problem. Therefore, no function of $\{F_l,\, G_l\}_l$ can generate a function of $\{H^x_{l,m},\, H^y_{l,m}\}$, not even with the help of other functions of $\{H^x_{l,m},\, H^y_{l,m}\}$.

Reversely, no function of $\{H^x_{l,m},\, H^y_{l,m}\}$ can generate functions in $\{F_l,\, G_l\}_l$. Indeed, we have $\Delta(h^x_{l,m}) \in \mathbb{Q}_{k-1}$. Therefore, assuming that we could find a function $v$ in $\mathrm{span}\{H^x_{l,m},\, H^y_{l,m}\}$ such that

$$v = \begin{pmatrix} v_1 \\ v_2 \end{pmatrix} = \begin{pmatrix} x \\ y \end{pmatrix} u$$

for some function $u$, then

$$\begin{aligned} \Delta(v_1) &= \Delta(xu) = 2\nabla \cdot u + x\Delta(u) \in \mathbb{Q}_{k-1} \\ \Delta(v_2) &= \Delta(yu) = 2\nabla \cdot u + y\Delta(u) \in \mathbb{Q}_{k-1}. \end{aligned}$$

Thus, $\Delta u \in \mathbb{Q}_{k-2} \not\subset \mathbb{Q}_{[k-1]}$ and $u \notin B_k$. And since in $B_k$ the orders of the degrees defining the right hand sides of the Laplacian equations are also two by two different, the functions cannot compensate each other this lack of order. Therefore, the set

$$\{F_l,\, G_l,\, H^x_{l,m},\, H^y_{l,m}\}_{\substack{l_1 \in [\![0,\, k-2]\!] \\ (l,\, m) \text{ unique pairs in } [\![0,\, k-1]\!]^2}}$$

is linearly independent.

■

**Let us now gather the sets $\mathcal{N}$ and $\mathcal{I}$** Let us show that the sets $\mathcal{N}$ and $\mathcal{I}$ together generate $\mathbb{H}_k(K)$.

**Proposition 6.11**

It holds:

$$\mathbb{H}_k(K) = \mathrm{span}\,\{\mathcal{N},\, \mathcal{I}\}$$

The proof is immediate through the following assertions.

**Proposition 6.12**

The set $\{\mathcal{N}, \mathcal{I}\}$ is free.

***Proof.*** We know that the sets $\mathcal{N}$ and $\mathcal{I}$ are free. Therefore, to show that the set $\{\mathcal{I}, \mathcal{N}\}$ is free, it is enough to show that no function of $\mathcal{N}$ combined with functions of $\mathcal{I}$ cannot create another function of $\mathcal{I}$, and *vice - versa.*

• Let us show first that no function of $\mathcal{N}$ combined with functions of $\mathcal{I}$ can create an other function of $\mathcal{I}$.

By definition, every function of $\mathcal{I}$ vanishes on the boundary. However, no function of $\mathcal{N}$ vanishes identically on the boundary, and as we saw in the previous *Section 4.3.1*, the only constants one can build with function in $\mathcal{N}$ without every function to be zero happens when every function is in use with an equal weight. As functions of $\mathcal{I}$ are all vanishing on the boundary, one may not use some to drag the constant down. Thus, no other function of $\mathcal{I}$ can be created.

• Reversely, functions in $\mathcal{I}$ have their normal components vanishing on every edge. Therefore, when combined with $\dim \mathcal{N} - 1$ functions in $\mathcal{N}$, the last function of $\mathcal{N}$ cannot be created. Indeed, to create the last function one would need to create some polynomial profile on the boundary. However as the function in $\mathcal{I}$ are vanishing on the boundary, one can only use functions in $\mathcal{N}$ to create a behaviour on the boundary. Furthermore, $\mathcal{N}$ being a free set, no function of $\mathcal{N}$ can create the last function in $\mathcal{N}$.

Thus, combining the two cases we have directly that the set $\{\mathcal{I}, \mathcal{N}\}$ is free.

■

**Property 6.13** Dimension

$$\dim \operatorname{span} \{\mathcal{N}, \mathcal{I}\} = \dim \mathbb{H}_k(K)$$

***Proof.*** We saw before that the set $\{\mathcal{I}, \mathcal{N}\}$ is free. Thus,

$$\begin{aligned}\dim \operatorname{span} \{\mathcal{N}, \mathcal{I}\} &= \dim\{\mathcal{N}, \mathcal{I}\} \\ &= \dim \mathcal{N} + \dim \mathcal{I}\end{aligned}$$

where the last line holds as we saw before that the sets $\mathcal{N}$ and $\mathcal{I}$ are free from each other and therefore their functions are two by two different. Gathering their previously computed dimensions, we get that

$$\dim \operatorname{span}\{\mathcal{N}, \mathcal{I}\} = n(k+3) + (2k(k-1)-1)\mathbb{1}_{k>0}$$

where the indicator function in the first line is there only to emphasise that there exist internal moments only for the cases $k > 0$. Thus, we derive

$$\begin{aligned}\dim \operatorname{span}\{\mathcal{N}, \mathcal{I}\} &= \underbrace{n(k+1) + k^2 - (k-1)^2}_{} + \underbrace{2n + 2k^2 \mathbb{1}_{k>0}}_{} \\ &= \dim B_k + 2\dim A_k \\ &= \dim \mathbb{H}_k(K).\end{aligned}$$

■

**Property 6.14**

If $p$ belongs to $\operatorname{span}\{\mathcal{N}, \mathcal{I}\}$, then $p$ belongs to $\mathbb{H}_k(K)$.

***Proof.*** If $p$ belongs to $\operatorname{span}\{\mathcal{N}, \mathcal{I}\}$, then $p$ can be written as a linear combination of elements in $\mathcal{N}$ and $\mathcal{I}$, that is

$$p = \sum_i \alpha_i q_{i,1} + \sum_i \beta_i q_{i,2}$$

with $\{\alpha_i\}_i$, $\{\beta_i\}_i$ sets of real constants, $q_{i,1}$ a function belonging to $\mathcal{N}$ and $q_{i,2}$ a function belonging to $\mathcal{I}$. Thus, since for any $i$ the terms $\{\alpha_i\}_i$ and $\{\beta_i\}_i$ are constants, we can infer by the *Property 6.2* with $l = k-1$, $m = 0$ and $l = k-1$, $m = k$ that the functions $\alpha_i q_{1,i}$ and $\beta_i q_{i,2}$ also belong respectively to $\mathcal{N}$ and $\mathcal{I}$. Therefore, by linearity of the Lagrangian we have by extension that

$$\sum_i \alpha_i q_{i,1} \in \mathcal{N} \quad \text{and} \quad \sum_i \alpha_i q_{i,2} \in \mathcal{I}. \tag{113}$$

Thus, $p$ belongs to $\mathbb{H}_k(K)$.

■

**Conclusion** The following functions form a basis of $\mathbb{H}_k(K)$ for polygonal elements in two dimensions.

*Normal basis functions*

$$\phi_j = e_{i,1} f_{i,m} + e_{i,2} g \quad {}_{\substack{i\in[\![1,\,n]\!], m\in[\![1,\,k+1]\!] \\ m\in[\![1,\,n(k+1)]\!]}}$$
$$\phi_j = e_{i,1} f_{i,m} - e_{i,l} g \quad {}_{\substack{i\in[\![1,\,n]\!], (m,l)\in\{(1,3),\,(k+1,4)\} \\ m\in[\![n(k+1),\,n(k+3)]\!]}}$$

*Internal basis functions*

$$\phi_j = F_l \quad l \in [\![0,\, k-2]\!],\, j \in [\![1,\, k-1]\!]$$
$$\phi_j = G_l \quad l \in [\![0,\, k-1]\!],\, j \in [\![k,\, 2k]\!]$$
$$\phi_j = H_l \quad l \in [\![0,\, k-1]\!],\, j \in [\![2k+1,\, k^2+2k-1]\!]$$

The vectors $e_{i,1}$, $e_{i,2}$ are defined in (102), $e_{i,3}$, $e_{i,4}$ in (103) and the functions $F_l$, $G_l$ and $H_l$ are defined in (109).

## 6.2 An example of a reduced space

As quickly addressed in the *Paragraph 6* of the *Section 5.3.2* and as it will be shown in the numerical results, a classical construction of the space $\mathbb{H}_k(K)$ implies degeneration of the misc normal functions into internal ones when tuned against the sets of degrees of freedom defining the elements $Ia$, $Ib$, $IIa$ and $IIb$.

This is a consequence of determining the boundary part of functions living in $\mathbb{H}_k(K)$ only by its normal component until completely describing its global polynomial behaviour, before using misc normal degrees of freedom to characterising information that is only accessible component -wise. There, the global normal component of the dual basis corresponding to those misc moments will vanish.

Thus, to allow a parallel with the Raviart–Thomas elements from the lowest order, one can consider the last paragraph of the *Remark 5* and set by example the following space up.

As this is a pure reduction from the general setting, we only gather quickly its main properties and basis functions. A choice of reduced elements will also be presented.

### 6.2.1 Definition of the reduced space

**Definition 6.15** Reduced space

$$\mathbb{H}_k(K) = \left\{u \in H^1(K),\ u|_{\partial K} \equiv 1,\ \Delta u \in \mathbb{Q}_{k-1}(K)\right\}^2 + \begin{pmatrix} x \\ y \end{pmatrix} \left\{u \in H^1(K),\ u|_{\partial K} \in \mathcal{H}_k(\partial K),\ \Delta u \in \mathbb{Q}_{[k-1]}(K)\right\}$$

***Note.***
• Here, any definition of $l_2$, $m_1$ and $m_2$ would be fine. The important modification is the prescription of a uniformly constant boundary condition within the subspace $A_k$, not to create further liberty where it is not relevant.

• One could also have prescribed $u|_{\partial K} \equiv 1$ in $B_k$ instead, but the control of the normal component on the boundary would be more delicate as the normal component of functions living in $A_k$ is not driven with respect to the edge position in space. ▲

### 6.2.2 Properties of the reduced space

The $H(\text{div})$ – conformity ready and the direct sum construction properties are inherited from the general setting. However, its dimension is reduced.

**Property 6.16** Dimension of the reduced space

$$\dim \mathbb{H}_k(K) = n(k+1) + 2k(k-1) - \mathbb{1}_{k>0}$$

***Note.*** Here, we retrieve a setting allowing the definition of exactly $k+1$ normal functions per edge as in the Raviart – Thomas elements. As this is the dimension of $\mathbb{P}_k(f)$, all the liberty will be required to totally determine the global normal component and no misc degenerating moments will be defined. ▲

### 6.2.3 Construction of basis functions

The construction of a set of basis functions for the reduced space is done similarly as in the general case, removing the two supplementary basis functions. It then comes the following definitions.

**Lowest order case ($k = 0$)** (only normal functions, for any edge $f_i \in \partial K$)**:**

$$f_1 := (x,\, y) \mapsto \begin{pmatrix} x \\ y \end{pmatrix} \begin{Bmatrix} \Delta u = 0 \\ u|_{\partial K} = \mathbb{1}_{f_i} \end{Bmatrix} + \begin{pmatrix} n_{ix} \\ n_{iy} \end{pmatrix} \begin{Bmatrix} \Delta u = 0 \\ u|_{\partial K} \equiv 2 \end{Bmatrix}$$

**General case ($k > 0$):** By denoting $h_{l,\,m}$ a polynomial of order $l$ on the $x$ component and $m$ on the $y$ component, the basis functions of the reduced space read as follows.

*Normal functions*: for any $j \in [\![0,\, k]\!]$,

$$f_j := (x,\, y) \mapsto \begin{pmatrix} x \\ y \end{pmatrix} \begin{Bmatrix} \Delta u = h_{k-1,\,k-1}(x, y) \\ u|_{\partial K} = g_j(x, y)\mathbb{1}_{f_i} \end{Bmatrix} + \begin{pmatrix} n_{ix} \\ n_{iy} \end{pmatrix} \begin{Bmatrix} \Delta u = 0 \\ u|_{\partial K} \equiv 2 \end{Bmatrix}$$

*Internal functions:*

$$f_{xl} := (x,\, y) \mapsto \begin{pmatrix} x \\ y \end{pmatrix} \begin{Bmatrix} \Delta u = h_{k-1,\,l}(x,\, y) \\ u|_{\partial K} = 0 \end{Bmatrix},\ 0 \leq l \leq k-1$$

$$f_{yl} := (x,\, y) \mapsto \begin{pmatrix} x \\ y \end{pmatrix} \begin{Bmatrix} \Delta u = h_{l,\,k-1}(x,\, y) \\ u|_{\partial K} = 0 \end{Bmatrix},\ 0 \leq l < k-1$$

$$f_{xlm} := (x,\, y) \mapsto \begin{pmatrix} 1 \\ 0 \end{pmatrix} \begin{Bmatrix} \Delta u = h_{l,\,m}(x,\, y) \\ u|_{\partial K} = 0 \end{Bmatrix},\ 0 \leq l,\, m \leq k-1$$

$$f_{ylm} := (x,\, y) \mapsto \begin{pmatrix} 0 \\ 1 \end{pmatrix} \begin{Bmatrix} \Delta u = h_{l,\,m}(x,\, y) \\ u|_{\partial K} = 0 \end{Bmatrix},\ 0 \leq l,\, m \leq k-1$$

***Note.***

• In practice, the space where the functions $\{h_{lm}\}_{lm}$ is either taken as the space of Hermite polynomials or Legendre polynomials, offering an acceptable conditioning of the transfer matrix to the designed elements.

• The more sensitive choice relates to the selection of the two normal moment one has to kick out per edge. Typically, one has to remove the one that do not characterise the functions anymore since we have less liberty. As we constricted the lowest dimensional space $A_k$, it is a safe choice for the first configuration to remove the two supplementary moments tuning the constant part on the boundary edge and coordinate wise. The elements $Ia$ and $Ib$ thus merge.

• For the configuration $II$, we keep the lowest order representation but remove the the two coordinate wise first order projection. Indeed, as the normal component is tuned, all the determination comes to lowest order. The difference between elements $IIa$ and $IIb$ is preserved. ▲

### 6.2.4 Construction of reduced elements

One can notice that the internal part of the space is left unchanged from the general setting. Therefore, the only change occurring in the element's definition lies in the selected normal degrees of freedom. We summarise a possible choice of normal moments to disregard from each classical configuration presented in the previous section.
*Normal moments*: for any edge $f_j \in \partial K$:

| Element | Core Moments | Misc Moment |
|---|---|---|
| Ia | $\int_{f_j} p_i(v)\,\varphi\cdot n\,\mathrm{d}\gamma,\ i\in[\![1,\,k]\!]$ | $\int_{f_j} v\,\varphi\cdot\begin{pmatrix}1\\1\end{pmatrix}\,\mathrm{d}\gamma$ |
| Ib | $\int_{f_j} p_i(v)\,\varphi\cdot n\,\mathrm{d}\gamma,\ i\in[\![1,\,k]\!]$ | $\int_{f_j} v\,\varphi\cdot\begin{pmatrix}1\\1\end{pmatrix}\,\mathrm{d}\gamma$ |
| IIa | $\int_{f_j} p_i(v)\,\varphi\cdot n\,\mathrm{d}\gamma,\ i\in[\![1,\,k]\!]$ | $\int_{f_j} \varphi\cdot n\,\mathrm{d}\gamma$ |
| IIb | $\int_{f_j} p_i(v)\,\varphi\cdot n\,\mathrm{d}\gamma,\ i\in[\![1,\,k]\!]$ | $\varphi(x_m)\cdot n$ |

*Internal moments:*

| | Core Moments | Misc Moment |
|---|---|---|
| All elements | $\int_K \begin{pmatrix}q_{l,m}(x,\,y)\\0\end{pmatrix}\cdot\varphi\,\mathrm{d}x\mathrm{d}y$ <br> $\int_K \begin{pmatrix}0\\q_{m,l}(x,\,y)\end{pmatrix}\cdot\varphi\,\mathrm{d}x\mathrm{d}y,$ <br> $l\in[\![0,k]\!],\,m\in[\![0,k-1]\!]$ <br> $(l,m)\neq(k,k-1)$ | $\int_K \begin{pmatrix}q_{k,\,k-1}(x,\,y)\\q_{k-1,\,k}(x,\,y)\end{pmatrix}\cdot\varphi\,\mathrm{d}x\mathrm{d}y$ |

***Note.*** In all configurations, there is no supplementary normal moment anymore. Furthermore, the elements $Ia$ and $Ib$ are now identical. ▲

# 7 Numerical Results

Let us now test our constructions. In a first study, we simply generate the canonical basis of the space $\mathbb{H}_k(K)$ presented in the *Section 6.1.2*. Then, after confirming its properties we investigate its tuning towards the sets of degrees of freedom $Ia$, $Ib$, $IIa$ and $IIb$. Thereupon, the sensitivity of the transfer matrix towards the shape of the element will be tested and the reasons of the geometrical restrictions we made in the *Section 5.2.1* will be witnessed. In a third study, we record the properties of each element through the shapes and layout of their corresponding basis functions. Further tests over the choice of the basis functions used in the definition of the degrees of freedom will also be performed. Finally, we show a parallel with the classical Raviart – Thomas setting for elements constructed from the degrees of freedom $IIb$ endowed with the reduced space presented in the *Section 6.2*.

***Note.***

- In all of the presented experiments, the *Fenics* $\mathbb{P}_2$ finite element solver have been used to compute the underneath Poisson's solutions.
- Unless mentioned, the basis functions presented in this section have been generated through the the Poisson's problems by using the Hermite polynomial space in the definition of the second members and the Lagrangian space for defining the boundary conditions.
- The computations of the Lagrangian sets $\{l_{i,j}\}_{i,j}$ have been performed using Gauss - Legendre sampling points.
- When building the basis functions, the last expression of the basis function (108) has been replaced by $e_{1,i}f_{i,k} - e_{i,4}g_i$ for the sake of coding convenience. Remark that when $k > 1$, as the functions $f_{i,1}$,$f_{i,k+1}$ and $f_{i,k}$ are both part of the set of solutions to the problem (104) built from a same Lagrangian set, this change does not impact the freeness of the generated set of functions. Furthermore, when we consider the lowest and first order spaces, the functions $f_{i,1}$ and $f_{i,k}$ are identical and the freeness of the functional set is only ensured by the definitions of $e_{i,3}$ and $e_{i,4}$.
- During the visualisation process, the basis functions are plotted through a nearest interpolation where the boundary points are interpolated only from the sampled neighbours that are also lying on the boundary.
- All the presented conditionings of transfer matrices are given truncated whenever there is no particular specification.

▲

## 7.1 Canonical Basis functions of $\mathbb{H}_k(K)$

We start by testing our construction of the canonical basis of $\mathbb{H}_k(K)$ on both convex and non - convex elements by considering the three shapes represented in the *Figure 59*. As they are planned to be mapped onto general elements of potentially fine meshes, we designed them lying within the unit circle in order to reduce the impact of a possible change in the order of magnitude c n

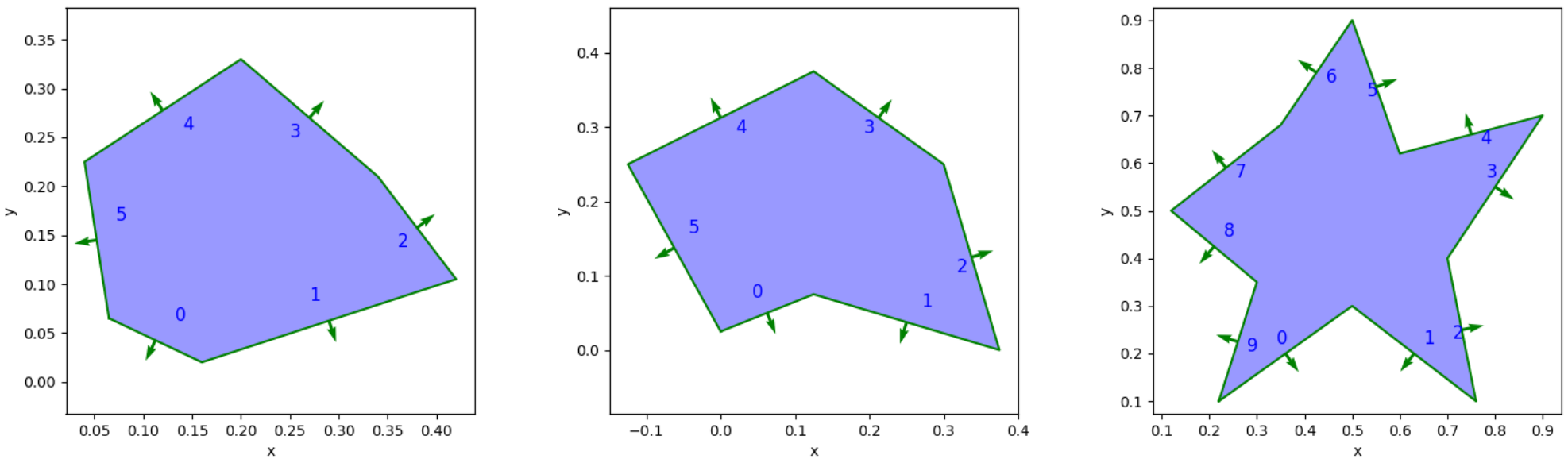

(a) Convex hexagon (b) Non - convex hexagon (c) Star shaped decagon

Fig. 59: Reference elements used for testing the construction of the canonical basis functions of $\mathbb{H}_k(K)$.

### 7.1.1 Normal Basis functions

To begin with, let us investigate the obtained normal basis functions through qualitative considerations.

On the boundary, any normal function $p$ lying in $\mathbb{H}_k(K)$ should fulfil the properties

$$\begin{cases} p \cdot n|_{\partial K} \in \mathcal{T}_k(\partial K) \\ p \cdot n|_{\partial K} \not\equiv 0. \end{cases}$$

Therefore, we aim to see our computed normal basis function $\phi$ fulfilling them. In particular, we would like to confirm that the boundary support of the basis function is restricted to a single edge. Indeed, the boundary conditions generating the core of the normal functions only consists of possibly shifted Lagrangian functions to which is stitched an edge indicator *(see the construction* (107)*)*. This vanishing property is then directly transmitted to $\phi$ component wise by the vectors $e_{i,1}$, $e_{i,2}$, $e_{i,3}$ and $e_{i,4}$. Thus, as the dot product is linear, the quantity $\phi \cdot n|_{\partial K}$ also enjoys a support which is restricted to a single edge. This can be directly observed on the *Figure 60*. Please also note that those normal functions will therefore never be continuous at the vertices.

| Convex hexagon | Non - convex hexagon | Star shaped decagon |
|---|---|---|

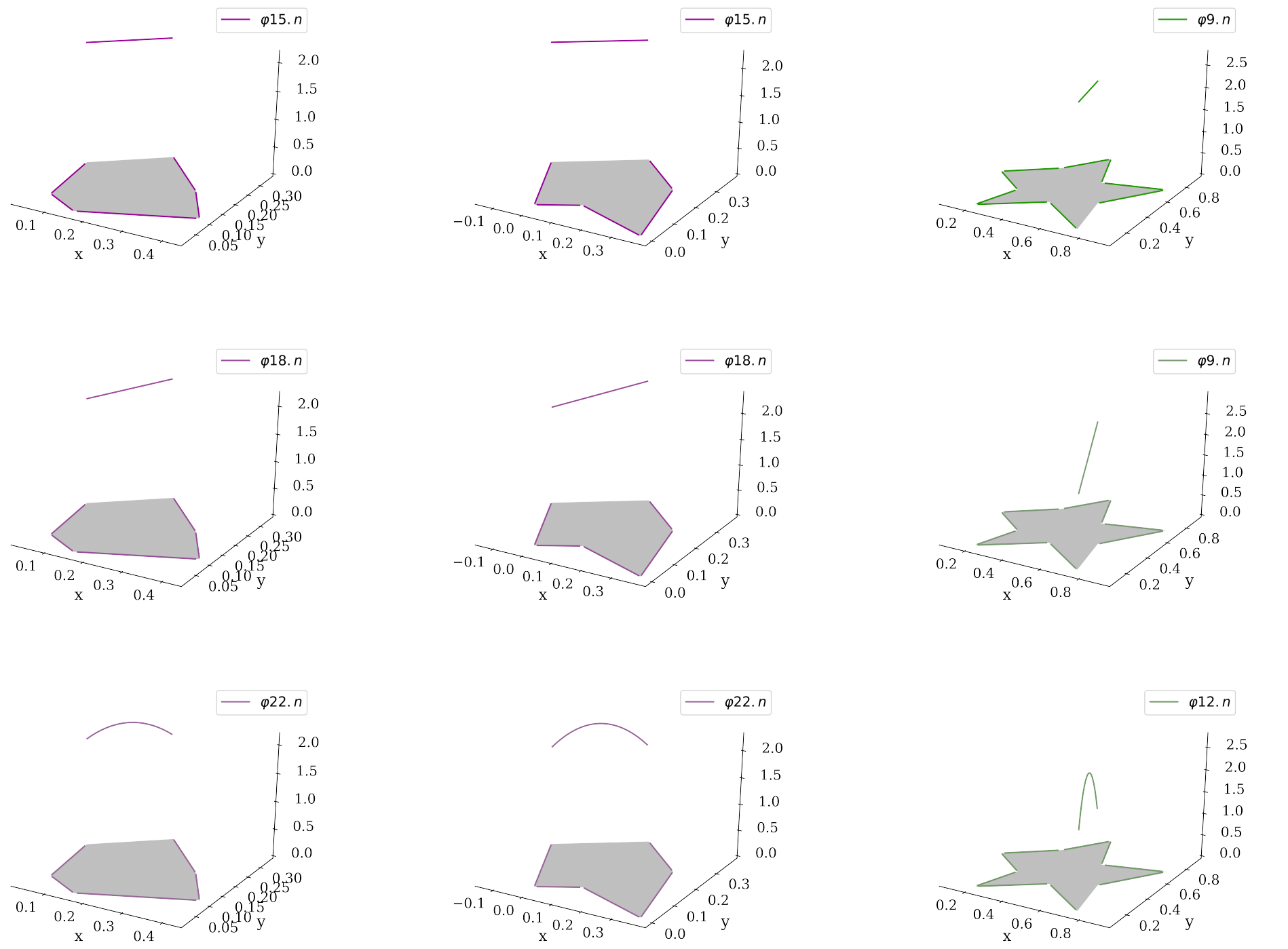

Fig. 60: Representation of some $\phi \cdot n|_{\partial K}$ for various orders and elements' shapes. First row: $k = 0$. Second row: $k = 1$, last row: $k = 2$.

As a by product, the boundary restrictions being built from Lagrangian sets, no normal function is identically vanishing on every edge. The second property is then naturally fulfilled. It also comes immediately that wherever the functions are not vanishing, they belong to $\mathbb{P}_k$. Indeed, the first row of the cases presented in *Figure 60* shows constant behaviours, the second linear behaviours and the last quadratic ones. Thus, the first property is respected for any order and any shape, verifying our theoretical results. In particular, the non - convexity of the thirdly presented shape has absolutely no impact on the construction of the normal functions.

Within the element, we expect the normal basis functions to enjoy some smoothness component wise, coming from the regularity of their Laplacian. And indeed, all the representative functions plotted in the *Figure 61* have decently regular coordinate wise components. There again, neither the shape of the element or the order of the discretisation space impact their regularity. In particular, the inner neighbourhood of the discontinuities present on the edges are handled smoothly.

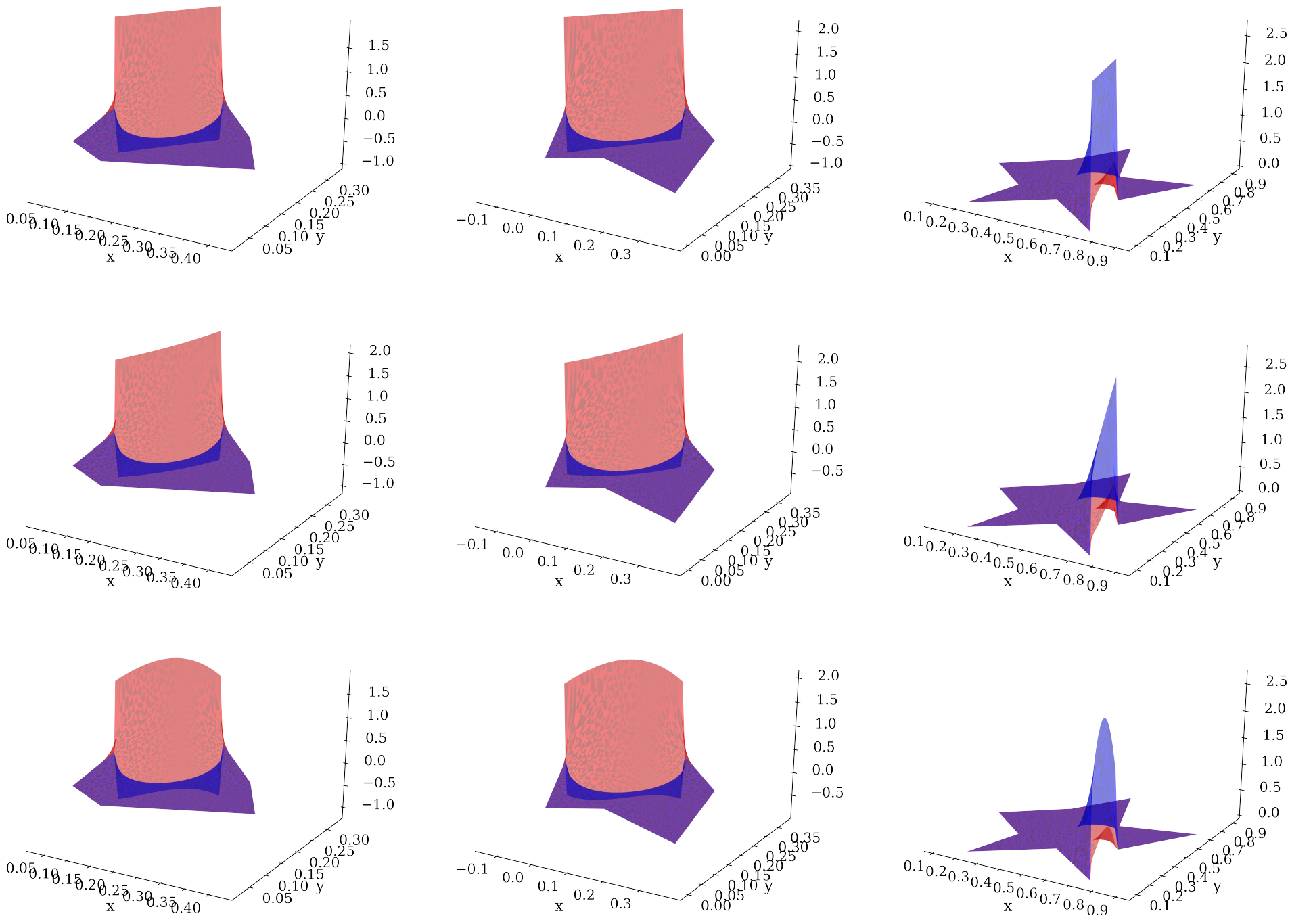

Fig. 61: Component wise internal functions within the element.
Blue: x component. Red: y component.

Furthermore, it can be seen that the polynomial variations of the basis functions' components on the boundary can be higher by one degree than their corresponding global normal component $\varphi \cdot n$. Indeed, by construction the polynomials emerging from the space $\mathrm{x}B_k$ are of degree $k+1$ on each component, degree that reduces when considering the normal component $\varphi \cdot n$ by the degeneration of the term $\mathrm{x} \cdot n$ into a constant.

This can be observed by example in the case of the quadrilateral element presented on the *Figure 62c* for the order $k = 2$. There, the $x$ component of the thirteenth basis function pictured in the *Figure 62a* clearly belongs to $\mathbb{P}_3(f_3)$. However, its global normal component plotted in the *Figure 62b* is a polynomial of order two, as expected.

***Remark.*** As noticeable in the *Figure 62b*, all the generated basis function do not share the same positivity. In fact, by construction the supplementary basis functions corresponding to the further freedom furnished by the space $A_k$ are shifted below, and become therefore most of the time negative. ▲

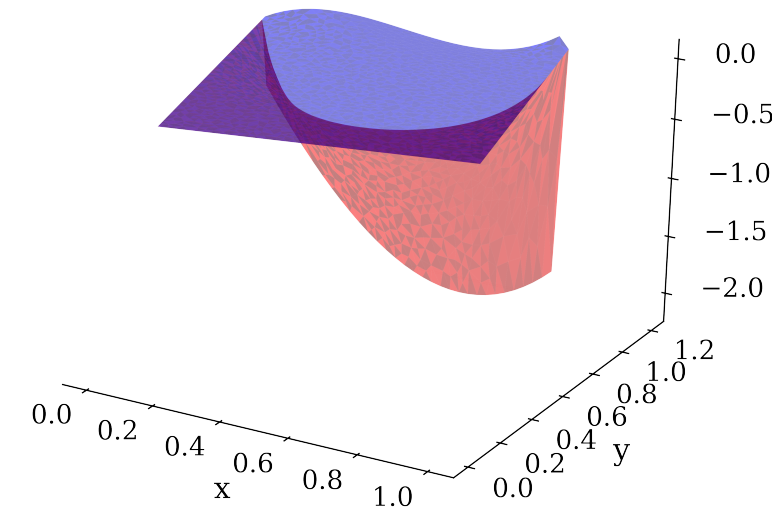

(a) Components of the 13th basis function. Blue: $x$ component. Red: $y$ component.

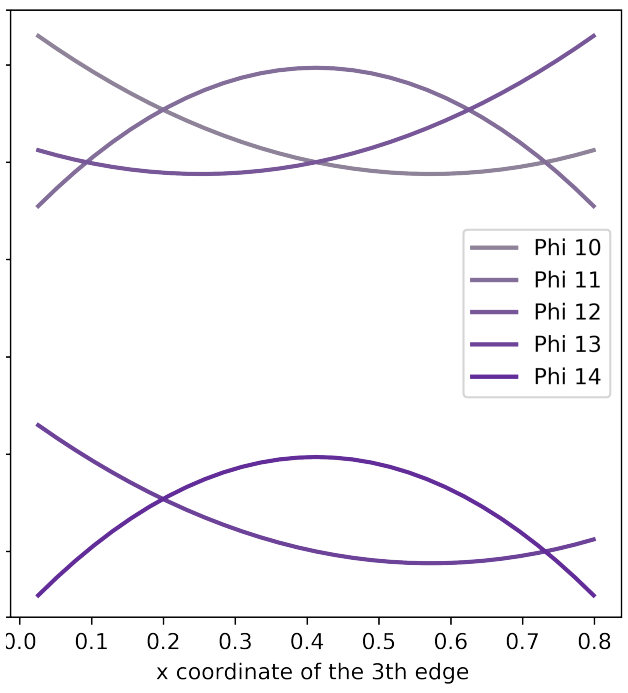

(b) Basis functions from the third edge

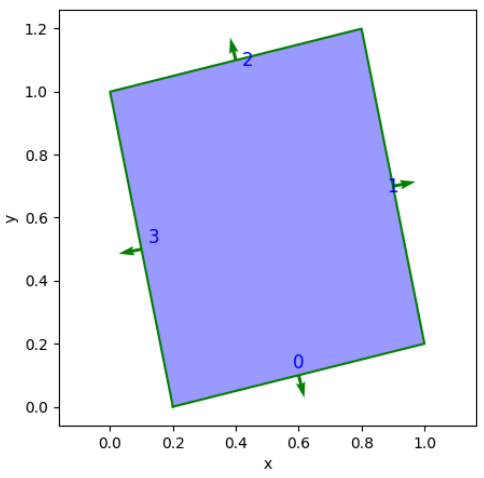

(c) Quadrilateral shape

Fig. 62: Example of an higher polynomial order component wise than when testing the normal component

The last property of the normal basis functions we would like to witness in the experiments involves the local property (6.5). We investigate it through the shift and local variations of $\varphi \cdot n$.

In the lowest order case represented in the *Figure 63*, the values of $\phi \cdot n$ are always constant for any of the presented shapes. By reading out the constant values, one can confirm a shift of $4 = 2-(-2)$ between the functions $\{e_{i,1} f_{3,0} + e_{3,2} g_i\}_i$ and the functions $\{e_{i,1} f_{3,0} - e_{3,3} g_0\}_i$ or $\{e_{i,1} f_{3,0} - e_{3,4} g_i\}_i$. In addition, one can notice that the shift raising the constant $-2$ to lower than $-1.5$ or from 2 to a higher than 2.5 in the decagon case corresponds to the value of $\mathrm{x} \cdot n_3$ for any x belonging to $f_3$. Indeed, the $\mathrm{x} \cdot n_3$ value can be computed by taking by example the value of one of the two edge's vertices and testing it against the normal; $\mathrm{x} \cdot n_3 = (0.76,\, 0.10)^T \cdot (0.98,\, 0.2)^T \approx 0.76$, corresponding to the shift witnessed on the graph. The complete layout of the third edge can be found in the *Appendix B*.

In the hexagonal cases, the considered edge implies a value of $\mathrm{x} \cdot n_5$ almost vanishing, resulting in a corresponding shift which is almost non - existent. Indeed, *e.g.* for the convex element, $\mathrm{x} \cdot n_5 = (0.20,\, 0.33)^T \cdot (-0.55,\, 0.84)^T \approx 0.17$, which is not noticeable at the scale of our picture. Thus, we can confirm the relation $\phi \cdot n = \mathrm{x} \cdot n \pm 2$.

In the first order case represented in the *Figure 64*, we observe linear variations of the normal components. In particular for the decagon, their variation rate match in absolute value the constant $\mathrm{x} \cdot n_3$ for any x belonging to $f_3$. Furthermore, they are shifted by two. Thus, the behaviour of those normal functions is comparable to

$$\phi \cdot n \; p(\mathrm{x} \cdot n_3) \pm 2$$

for a given function $p$ belonging to $\mathbb{P}_1$. Observing that the two functions built from the vector $e_{i,2}$ sum up to the constant function $\mathrm{x} \mapsto \mathrm{x} \cdot n_3 + 2$, one directly retrieve the pointwise local property (6.5) for two sampling points distributed on the edge $f_3$. Note however that the two functions built from the vectors $e_{3,3}$ and $e_{3,4}$ do not sum to $\mathrm{x} \cdot n - 2$ due to their identical direction of variation, but still enjoy the pointwise property (6.5) as begin on the edges shifted copies of $e_{3,1} f_{3,1} + e_{3,2} g_3$. The same observation can be drawn for the two hexagons.

For the last case $k = 2$ that we describe here *(see Figure 65)*, one can directly observe the Lagrangian property of the three functions that are constructed from the vector $e_{3,2}$. The two last ones corresponds to the first and second functions that have been shifted below by four, and only inherit the pointwise property (6.5). Note that in this case, and more generally for any order bigger or equal to two, the constant $\mathrm{x} \cdot n$ impacts every coefficients of the monomials entering in the definition of the normal component $\phi \cdot n$ that have a degree greater or equal to one.

Lastly, let us remark that for the three different orders $k$ detailed before, no normal function is vanishing on all edges, preserving in essence the dimension of the split into normal and internal subspaces. However, for the two lowest orders the normal functions generated from $e_{i,3}$ and $e_{i,4}$ merge (only) on the edges. For the first order space, it is a consequence of the computational simplification made by taking $f_{i,k}$ instead of $f_{i,k+1}$ in their definition. For the lowest order space, this side effect is unavoidable as there is only one Lagrangian function per edge.

***Remark.*** This merge directly connects to the fact that the space $\mathbb{H}_k(K)$ restricted to the boundary is not in direct sum anymore. ▲

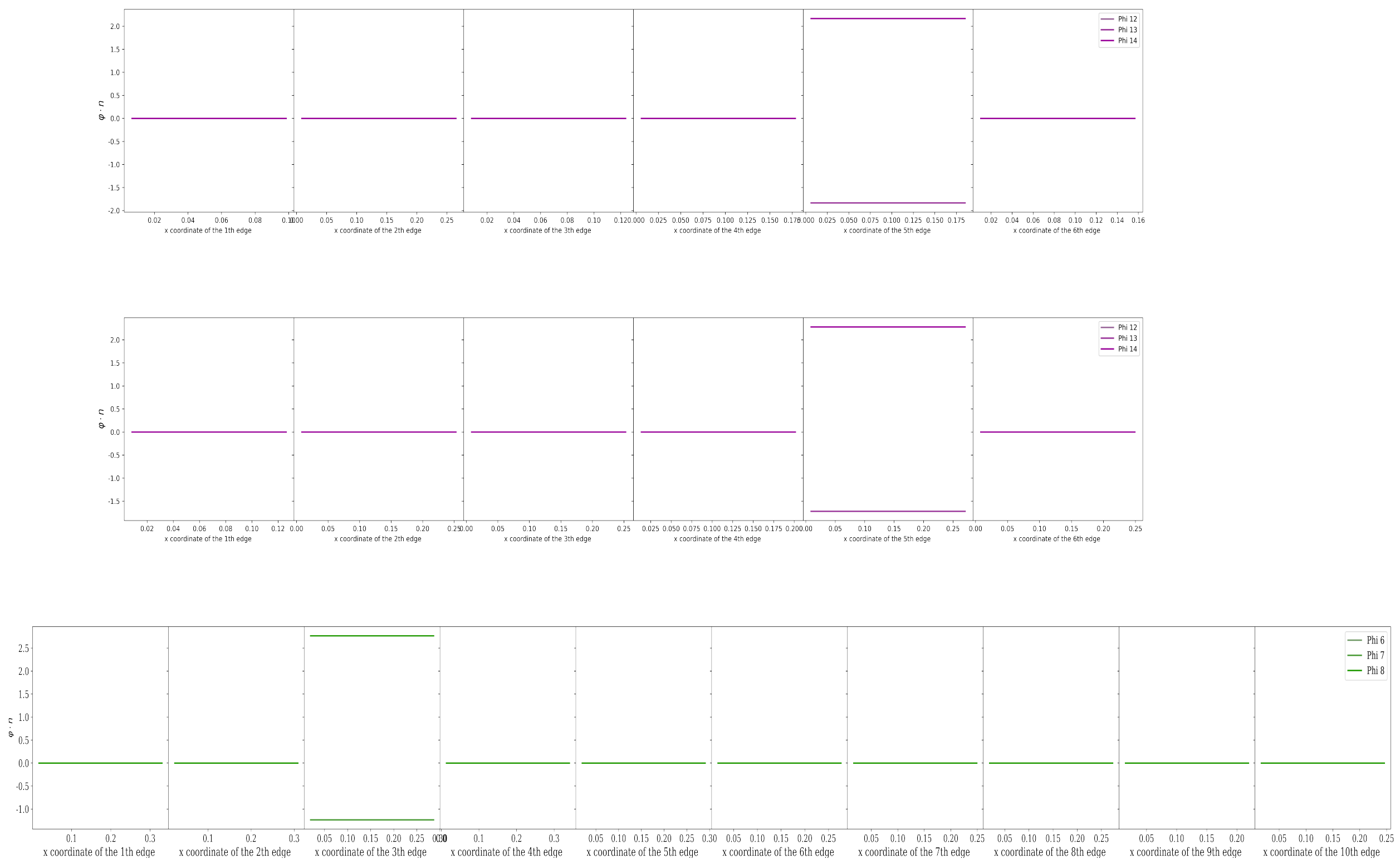

Fig. 63: Split view of all the normal components of basis functions on all the boundaries. Case for $k = 0$ and all basis functions that are generated from the boundary conditions set that are not vanishing on the fourth edge for the two first elements and second one for the last presented element. Top: convex hexagon, middle: not convex hexagon, bottom: star shaped decagon.

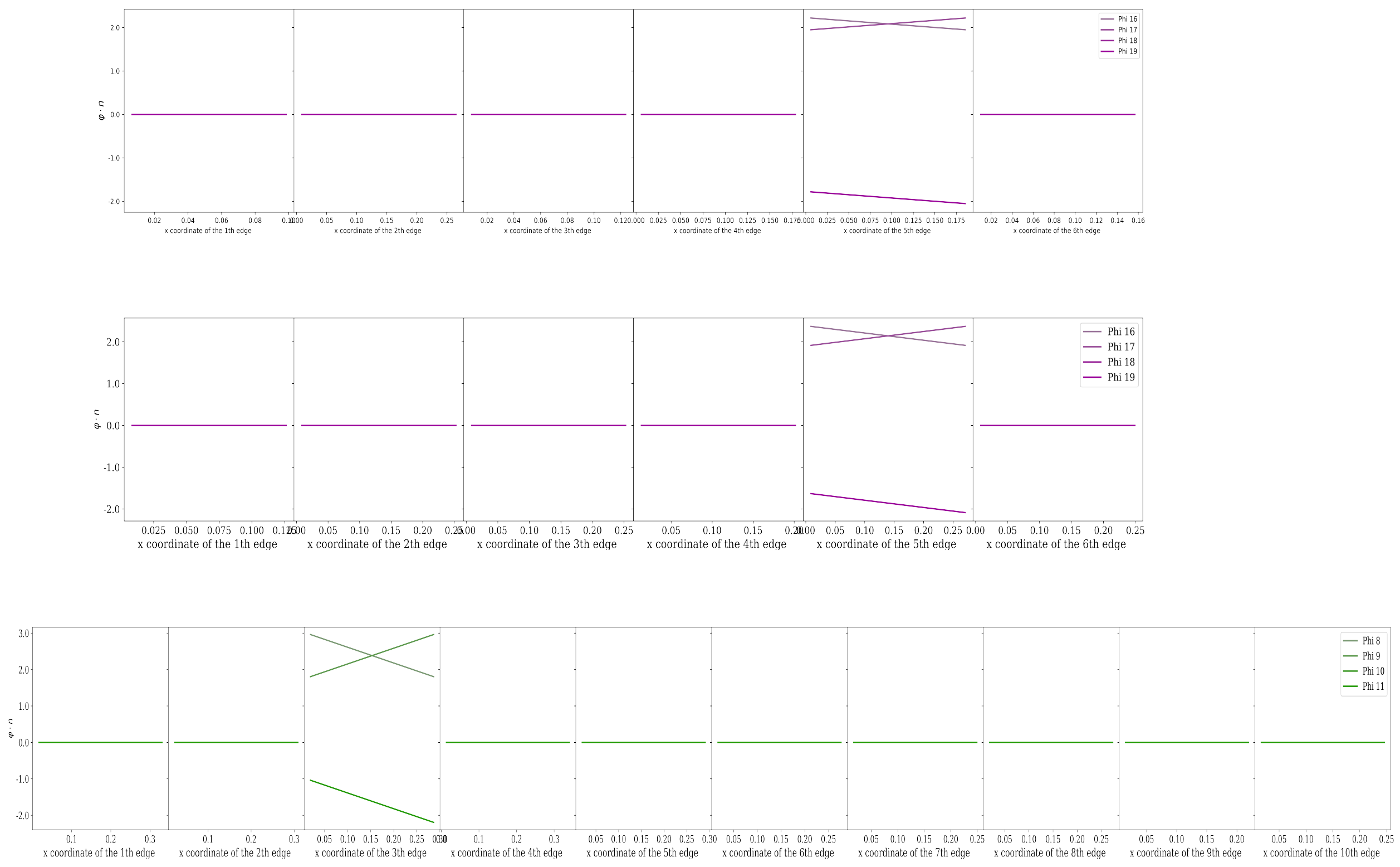

Fig. 64: Split view of all the normal components of basis functions on all the boundaries. Case for $k = 1$ and all basis functions that are generated from the boundary conditions set that are not vanishing on the fifth edge for the hexagons and third for the decagon. Top: convex hexagon, middle: not convex hexagon, bottom: star shaped decagon.

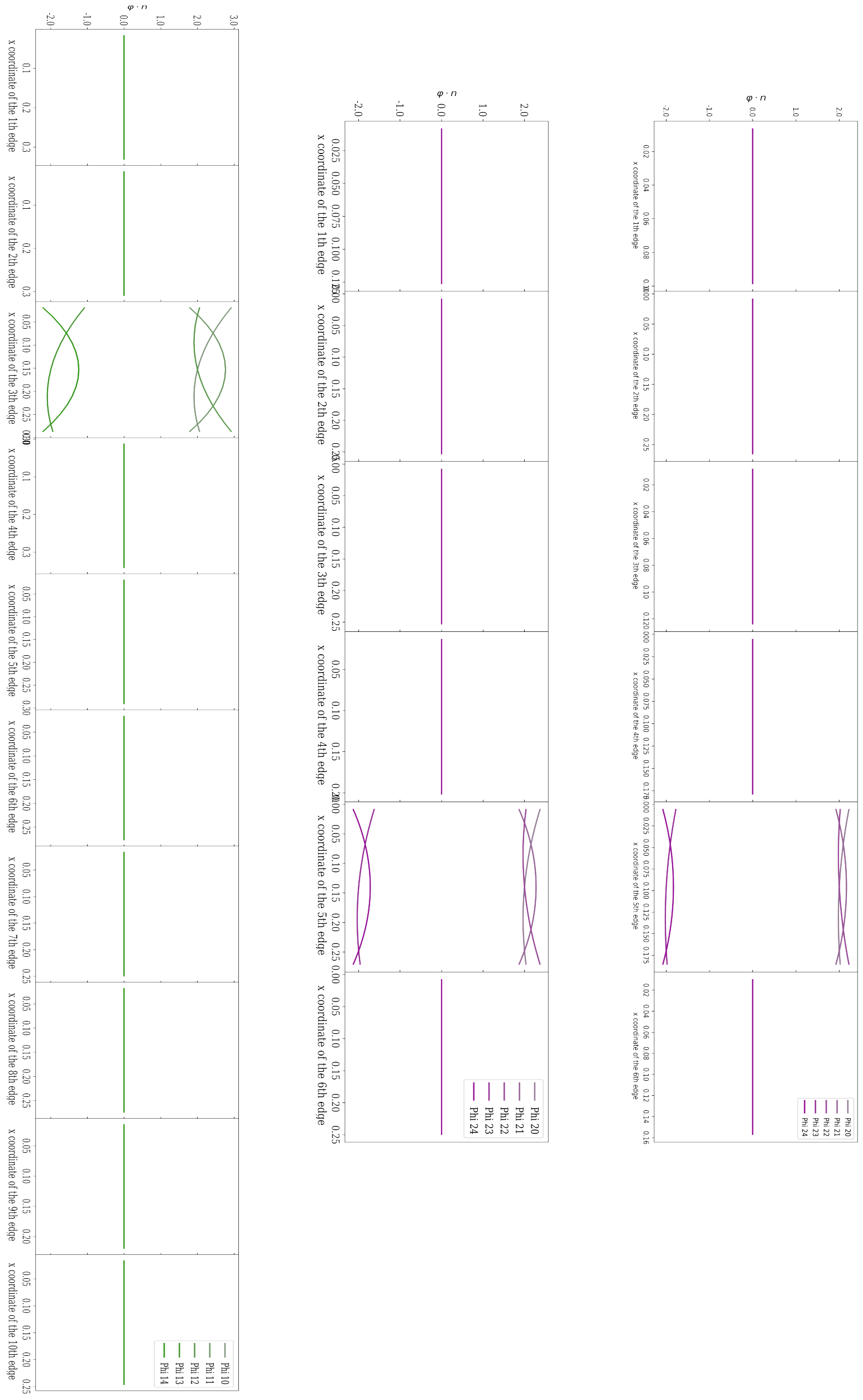

Fig. 65: Split view of all the normal components of basis functions on all the boundaries. Case for $k = 2$ and all basis functions that are generated from the boundary conditions set that are not vanishing on the second edge. Top: convex hexagon, middle: not convex hexagon, bottom: star shaped decagon.

So far, the qualitative behaviour is similar for any of the edges. Let us now investigate how the canonical basis functions are sensitive to the slope and position of the edges with respect to the axes. To this sake, we define the element presented in *Figure 66* and focus on its three edges that share similarities.

| Edge | Normal $(n_x, n_y)$ | Norm |
|---|---|---|
| 1 | (0.24, -0.97) | 0.29 |
| 3 | (0.24, 0.97) | 0.21 |
| 4 | (-0.24, 0.97) | 0.21 |

Tab. 2: Properties of the similar edges of the element.

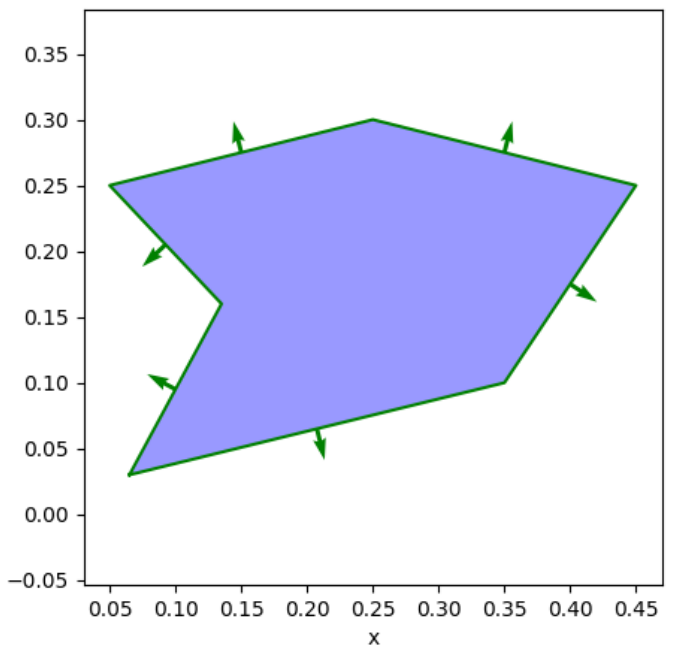

Fig. 66: Non - convex hexagon with edge similarities

Let us focus on the case $k = 0$ and $k = 1$, being sufficient to draw general conclusions. The behavior of the normal components are represented in the *Figure 67* where for the sake of concision only their non - vanishing support is considered.

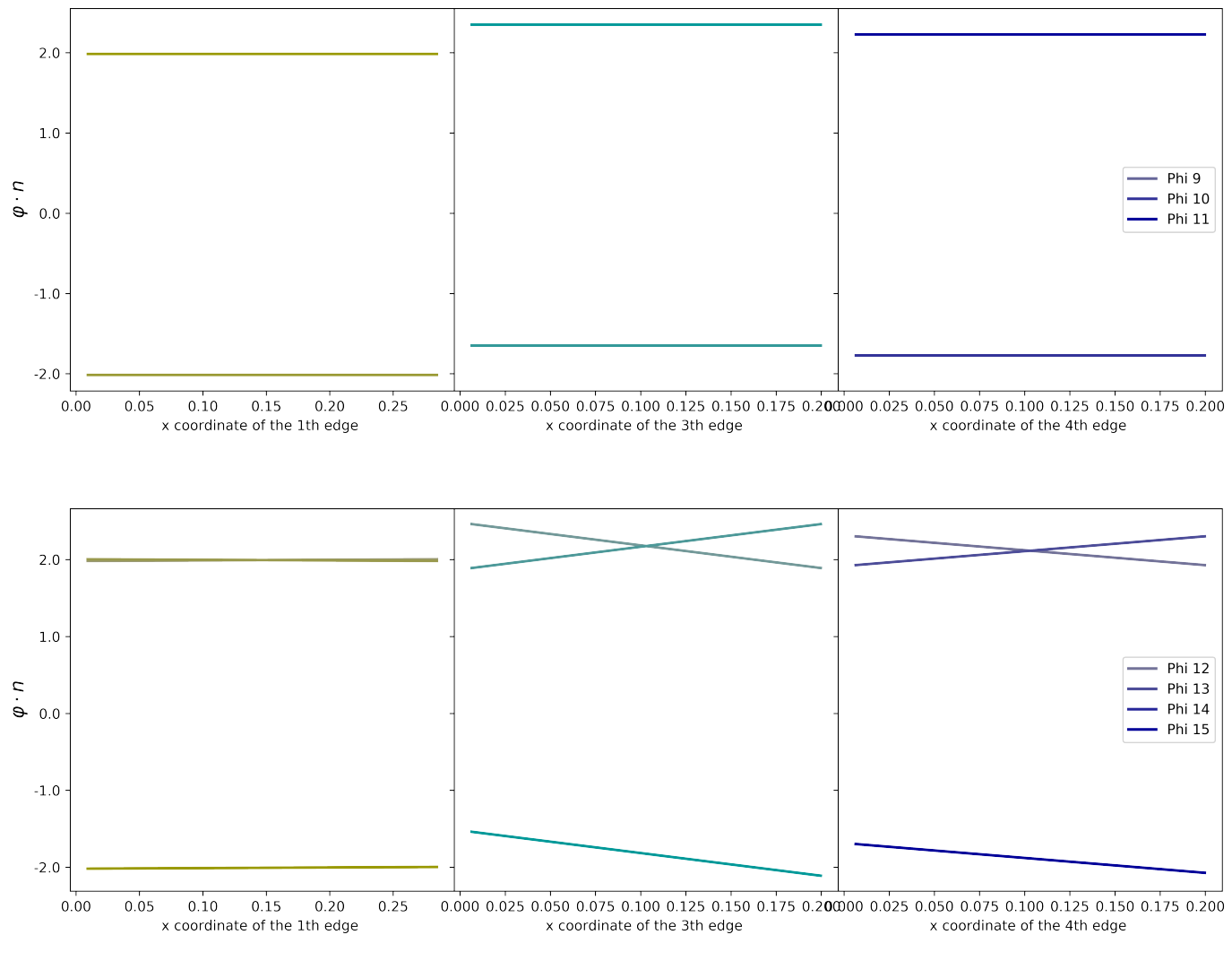

Fig. 67: Normal components of the basis functions built from the edges 1, 3 and 4 that are plotted respectively on the edges 1, 3 and 4.. Top: $k = 1$. Bottom: $k = 2$.

We first study the local variation rate of those basis functions, vanishing only in the case $k = 0$. There, in the case $k = 1$ and as quickly mentioned before, the growing rate of the basis function is given by the coefficient $x \cdot n$. As this value only depends on the angle the edge forms with with the $x$ - axis *see Figure 68*, it only depends on the angle $\alpha$ and on the position of the edge in the space.

However, this dependence is smooth thorough the plane by the linearity of the dot product. Therefore, slight variations in the angle $\alpha$ or of the position of the edge only implies a slight variation in the growing rate of the basis function. As a consequence, in order to avoid similarities in the basis functions (even if it would be up to their boundary support), it is advised to define edges having clearly different shapes.

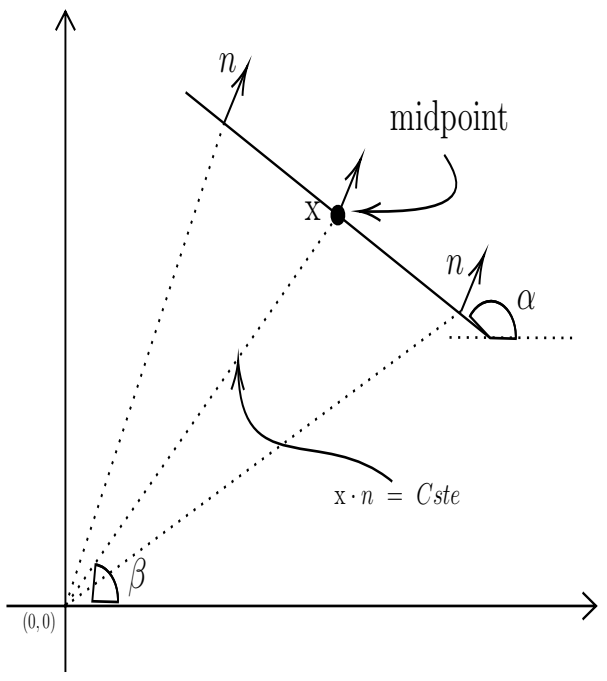

Fig. 68: Determination of $x \cdot n$

Please also note that since the normal is normalized and that the element lies within the unit circle, the slope of the basis function will be bounded by $\max|x \cdot n| = 1$. This bound is realized whenever the edge is perpendicular to the vector $(x - 0,\, y - 0)^T$ for any x belonging to the edge. Similarly, the minimum value, 0, is realized when the normal of the edge is perpendicular to the vector $(x - 0,\, y - 0)^T$ for any x belonging to the edge. Thus, the variations of the normal components of the basis functions are always lower than the ones of the classical Lagrangian functions.
Those observations generalises to higher orders when considering the local variations of Lagrangian polynomials squeezed by the value $x \cdot n$.

The impact of the length of the edge on the definition of the basis functions is much more restricted.
Indeed, as one can see for the first edge in the case $k = 0$ plotted in the top left corner of the *Table 67*, the offset of the basis function, 2, is not impacted by the length of the edge as it does not take part in the determination of the constant part of the function. However, in conjunction with the value $x \cdot n$, the length does determine the highest value the basis function can take on the edge.

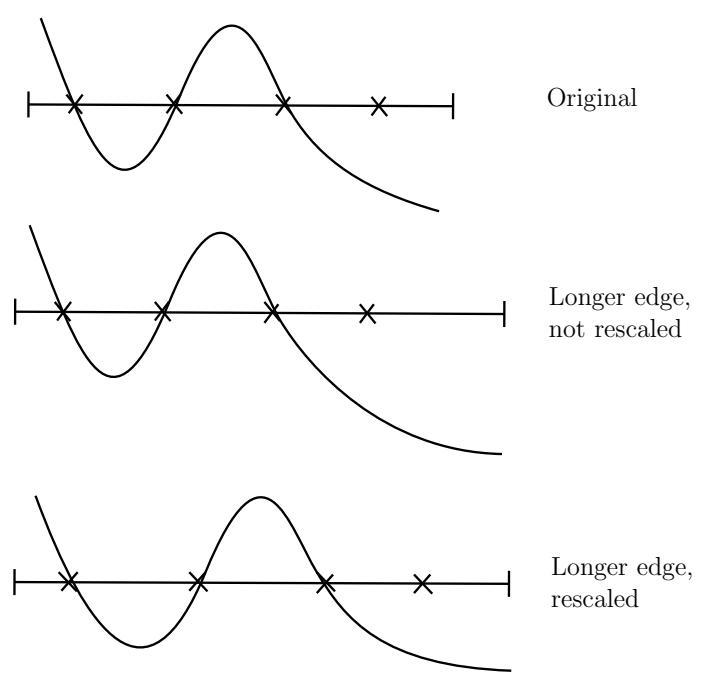

Fig. 69: Impact of uncontrolled variations

Indeed, since the sampling points of the Lagrangian functions have to be distributed strictly within the edge, its length impacts the amount of uncontrolled variation between the two sampling points closest to the vertices and the vertices themselves. Therefore, if the points are not homogeneously rescaled with the modification of the edge's length, then one may have an higher variation *(see Figure 69)*. However, assuming an homogeneous sampling method of the Gauss - Legendre points on the edges thorough the edges, the impact of the edge's length is limited to the definition of the boundary support of the basis function.

***Remark.*** In a perspective of lowering the sensitivity of the transfer matrix to small perturbations, we can also confirm by the above examples that the amplitude of the basis functions only depends on the length of the edge and of the distribution of the sampling points. Therefore, it is advised to choose them homogeneously thorough the edges while prescribing some close enough to the verticies. ▲

### 7.1.2 Internal Basis functions

Let us now investigate the behaviour of the internal basis functions. First of all, we can confirm that we enjoy some regularity within the cell component wise. Indeed, as the internal basis functions are built as regular Laplacian solutions with homogeneous Dirichlet boundary conditions to which have been stitched vectors belonging to at most $(\mathbb{P}_1)^2$, we obtain results from which some are presented in the *Figures 70* and *113*. There, one can directly observe a decent regularity regardless the order or shape that is in use.

One can also confirm on those two set of examples that due to the definition of $\mathcal{P}_k$ used in the *Section 6*, the last internal basis function is the only one whose both components are not vanishing. Indeed, all the other internal basis functions see either their $x$ or $y$ component vanishing *(see Figure 113)*. Furthermore, the smallness of their amplitude will turn to be an advantage when tuning them towards the sets of degrees of freedom.

Lastly, we also tested the global vanishing property of their normal component on the boundary of $K$, and according to the results represented in the *(see Figure 72)* it is preserved.

| Convex hexagon | Non - convex hexagon | Star shaped decagon |
|---|---|---|

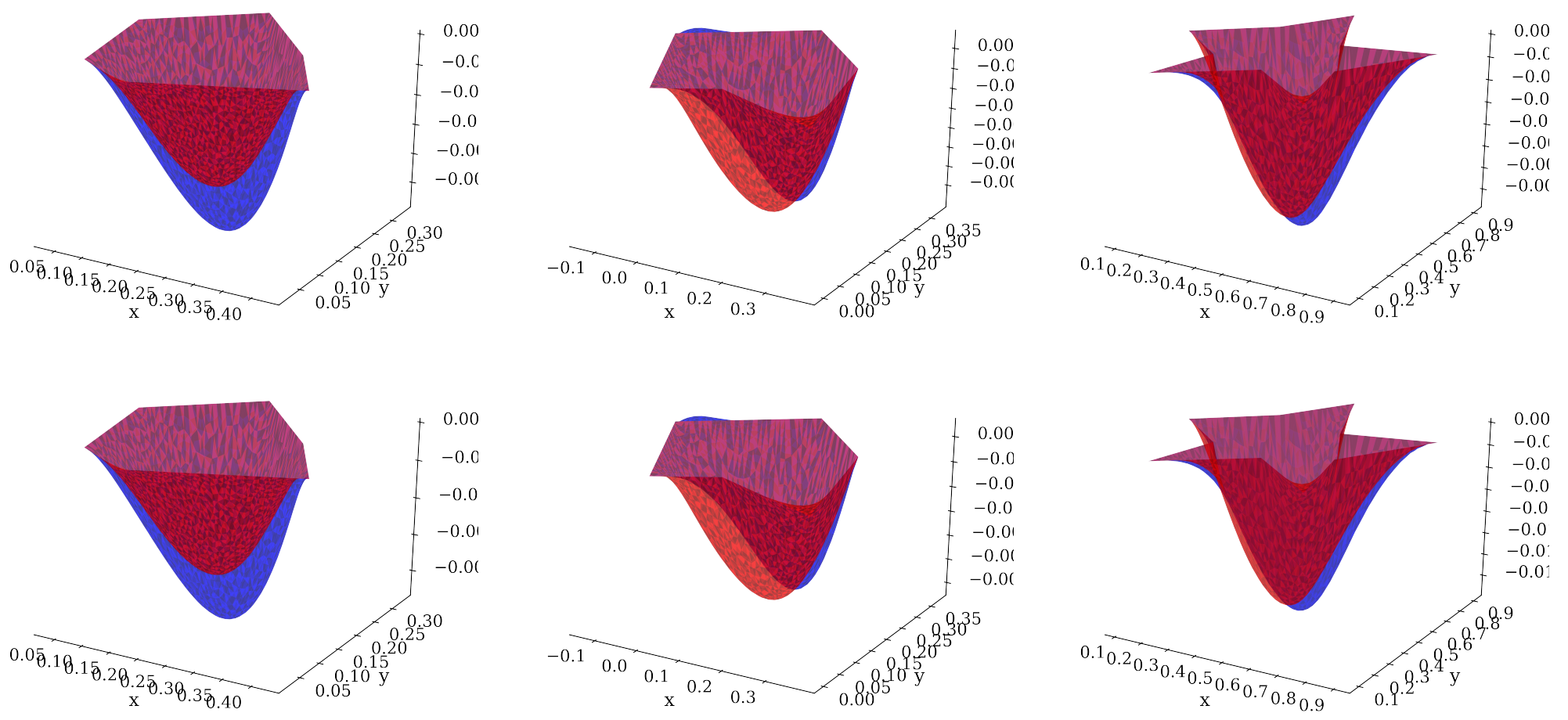

Fig. 70: Example of the last internal basis function within the elements. Top: $k = 1$. Bottom: $k = 2$. Blue: x component, Red: y component.

| Convex hexagon | Not convex hexagon | Star shaped decagon |
|---|---|---|

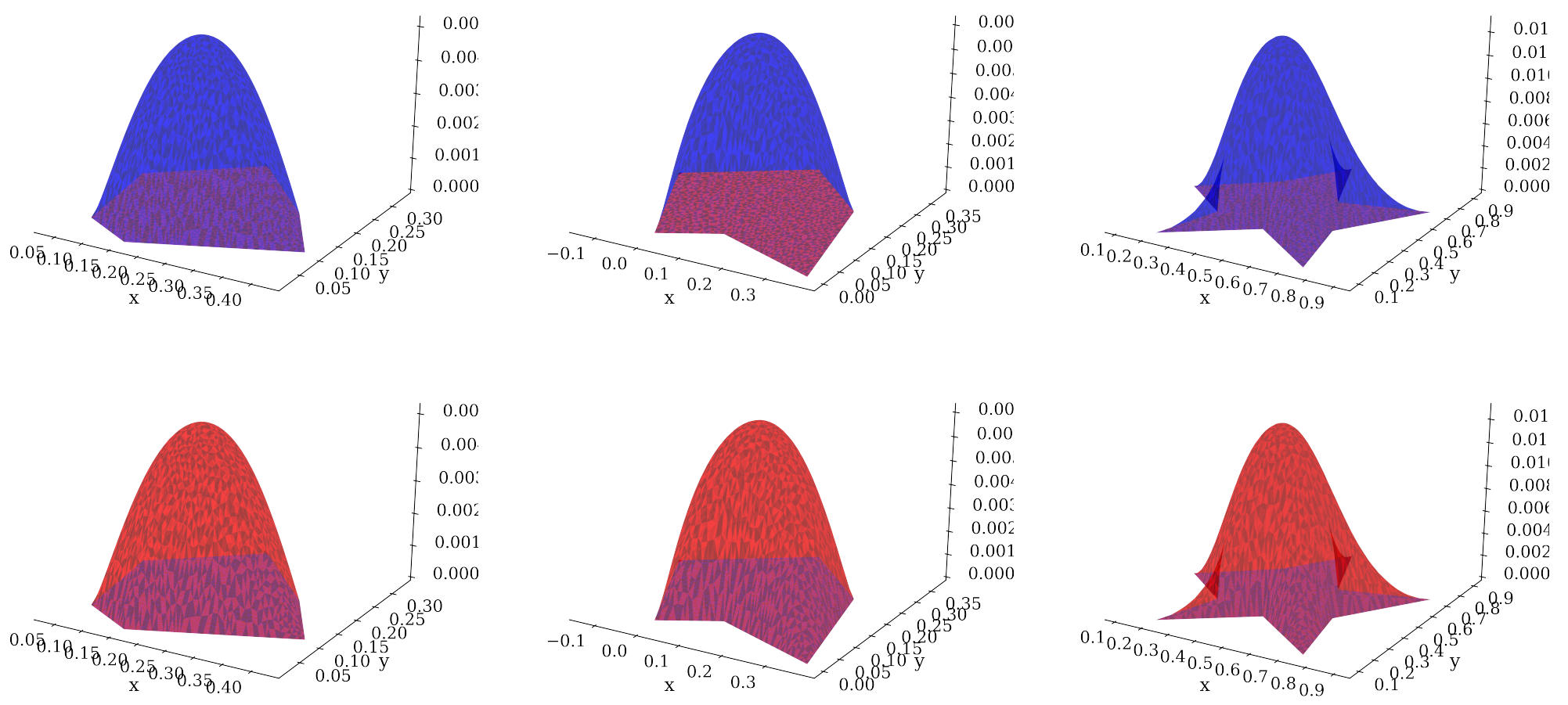

Fig. 71: Some of coordinate wise internal basis function in the case $k = 2$. Blue: x component. Red: y component.

In conclusion, our construction of the canonical basis fulfils our wishes of smooth functions within the elements that are preserving the split into internal and normal classifications. However, the reliability of the basis functions depends on the quality of the numerical approximations of the solution to the Laplacians problem. In particular, even though their settings are quite classical, it can quickly turn to be a problem for heavily squeezed elements. By example, star shaped element with high number of edges may be delicate to design.

Convex hexagon | Non - convex hexagon | Star shaped decagon

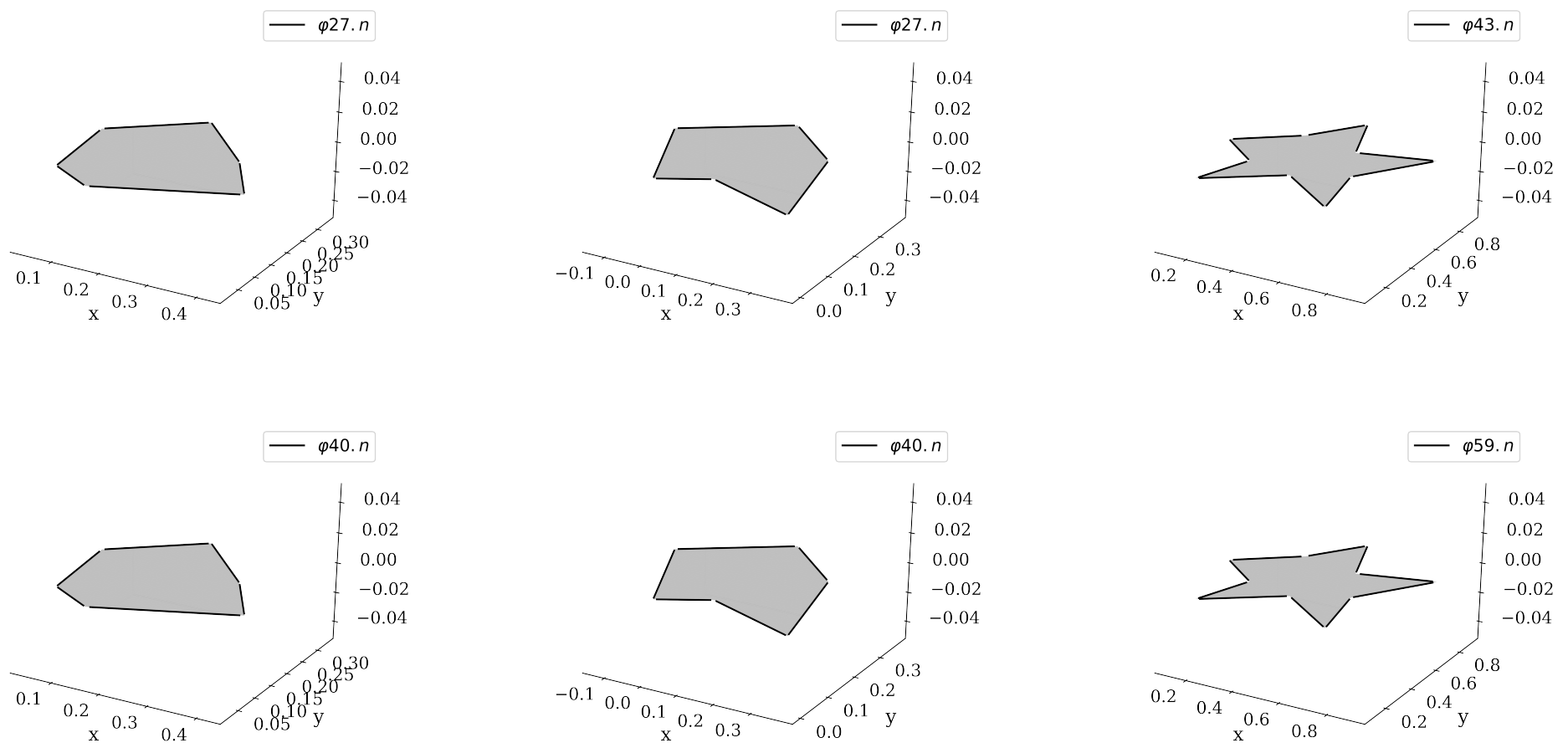

Fig. 72: Normal component of the last internal basis function on the boundary of the elements. Top: $k = 1$. Bottom: $k = 2$.

As our canonical basis construction is reliable, let us tune them against the previously described degrees of freedom (6.1) and (100) to form the desired elements.

## 7.2 Tuning to the basis functions corresponding to the said elements

In order to define the $H(\text{div}, K)$ – conformal elements presented in the *Section 5*, we combine the previously built basis functions and the sets of degrees of freedom $Ia$, $Ib$, $IIa$ and $IIb$ with the method given in the *Paragraph 2.4*. Thus, in all this section the term "transfer matrix" refers to the matrix $\Lambda$ built from the considered set of degrees of freedom. Unless explicitly mentioned, the elements considered in all this section have been defined via a Hermite polynomial space for the boundary projections. Within the elements, the space $\mathcal{P}_k$ has also been taken as a Hermite polynomial space in dimension two.

We start by retrieving the structure of this matrix before detailing its invertibility and robustness towards the shape the elements are built on.

### 7.2.1 Structure of the transfer matrix

Let us show first the transfer matrices corresponding to the two main types of degrees of freedom, *Ib* and *IIb*, built on the previous basis function set. The relation of their structure with the definition of the degrees of freedom will show up.

For the sake of concision we consider the simplicial element presented in the *Figure 73* and use the first order elements. The corresponding matrix is presented in the following page.

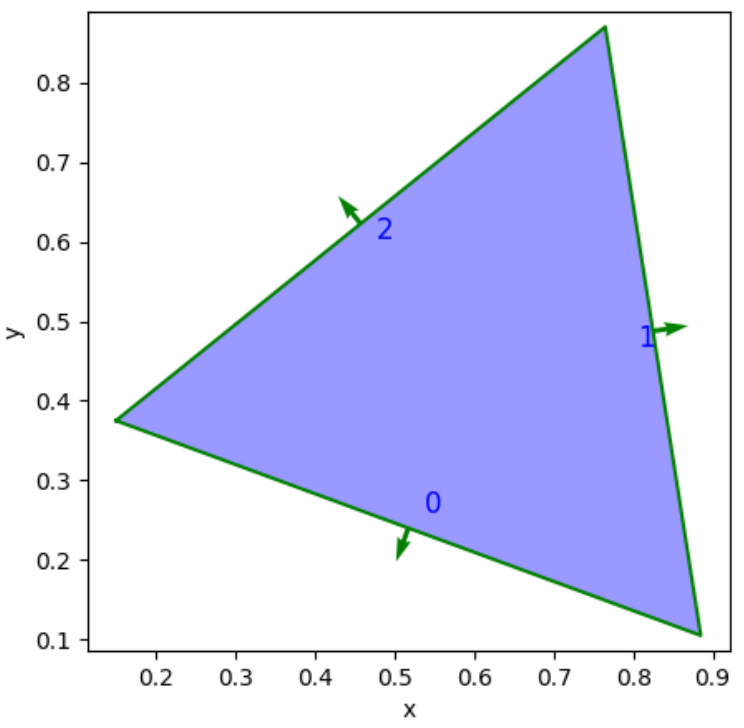

Fig. 73: Triangle of reference.

In both configurations, the twelve first rows and twelve first columns form a matrix which is diagonal by block. Knowing that each column corresponds to a basis function and each row to a degree of freedom, it clearly emphases the edge wise boundary support of the basis functions and edge wise definition of the degrees of freedom. One block then corresponds to the informations furnished by the normal degrees of freedom on one specific edge. Therefore, all the informations retrieved from a specific edge are only used to tune the subset of normal basis functions whose support is not vanishing on this edge.

The last three rows and last three columns corresponds respectively to internal degrees of freedom and internal basis functions. Therefore, only by reading out the structure of the matrix we can confirm that no internal function is involved in the construction of the normal functions. The reverse however is not true, as detailed previously when constructing the basis functions. Furthermore, having a closer look to the internal submatrix (the intersection of the last three rows and last three columns), we can infer that the internal basis functions see their $x$ or $y$ contribution vanish for the first or second internal moment, which is compliant with their definition.

Lastly, one can observe that the internal submatrix is identical for any of the four proposed elements, confirming that our definition of internal degrees of freedom is identical for any element considered in this scope. The internal basis functions will then always be tuned in a similar way. Furthermore, by definition the internal degrees of freedom all characterizes different aspects of functions in $\mathbb{H}_k(K)$ within the element.

**Element IIb**

| | | | | | | | | | | | | | | |
|---|---|---|---|---|---|---|---|---|---|---|---|---|---|---|
| 0.082 | 0.025 | -0.013 | -0.811 | 0 | 0 | 0 | 0 | 0 | 0 | 0 | 0 | 0 | 0 | 0 |
| 0.290 | 0.316 | -0.406 | -0.036 | 0 | 0 | 0 | 0 | 0 | 0 | 0 | 0 | 0 | 0 | 0 |
| 0.348 | -0.348 | 0.348 | 0.348 | 0 | 0 | 0 | 0 | 0 | 0 | 0 | 0 | 0 | 0 | 0 |
| 1.798 | 1.798 | -2.202 | -2.202 | 0 | 0 | 0 | 0 | 0 | 0 | 0 | 0 | 0 | 0 | 0 |
| 0 | 0 | 0 | 0 | 1.463 | 1.505 | 0.236 | -1.022 | 0 | 0 | 0 | 0 | 0 | 0 | 0 |
| 0 | 0 | 0 | 0 | 0.047 | 0.022 | -0.714 | 0.029 | 0 | 0 | 0 | 0 | 0 | 0 | 0 |
| 0 | 0 | 0 | 0 | -0.760 | 0.760 | -0.760 | -0.760 | 0 | 0 | 0 | 0 | 0 | 0 | 0 |
| 0 | 0 | 0 | 0 | 2.481 | 2.409 | -1.519 | -1.519 | 0 | 0 | 0 | 0 | 0 | 0 | 0 |
| 0 | 0 | 0 | 0 | 0 | 0 | 0 | 0 | 0.259 | 0.183 | -0.020 | -0.731 | 0 | 0 | 0 |
| 0 | 0 | 0 | 0 | 0 | 0 | 0 | 0 | 0.658 | 0.763 | -0.896 | 0.071 | 0 | 0 | 0 |
| 0 | 0 | 0 | 0 | 0 | 0 | 0 | 0 | -0.173 | 0.171 | -0.171 | -0.171 | 0 | 0 | 0 |
| 0 | 0 | 0 | 0 | 0 | 0 | 0 | 0 | 2.100 | 2.098 | -1.900 | -1.900 | 0 | 0 | 0 |
| -0.040 | -0.031 | 0.021 | 0.537 | 0.202 | 0.205 | 0.031 | -0.144 | -0.091 | -0.083 | 0.021 | 0.306 | 0.002 | 0 | -0.001 |
| -0.149 | -0.152 | 0.207 | 0.018 | 0.052 | 0.044 | -1.093 | 0.025 | 0.161 | 0.168 | -0.208 | 0.021 | 0 | 0.002 | -0.001 |
| 0.116 | 0.135 | -0.174 | 0.182 | 0.170 | 0.165 | 0.535 | -0.132 | -0.065 | -0.038 | 0.041 | 0.072 | 0.001 | -0.001 | -0 |

**Element Ib**

| | | | | | | | | | | | | | | |
|---|---|---|---|---|---|---|---|---|---|---|---|---|---|---|
| 0.149 | 0.149 | -0.089 | -2.089 | -0 | -0 | -0 | -0 | -0 | -0 | -0 | -0 | -0 | -0 | -0 |
| 1.650 | 1.650 | -2.113 | -0.113 | -0 | -0 | -0 | -0 | -0 | -0 | -0 | -0 | -0 | -0 | -0 |
| -0.733 | -0.565 | 0.685 | 1.792 | 0 | 0 | 0 | 0 | 0 | 0 | 0 | 0 | 0 | 0 | 0 |
| 0.348 | -0.348 | 0.348 | 0.348 | 0 | 0 | 0 | 0 | 0 | 0 | 0 | 0 | 0 | 0 | 0 |
| 0 | 0 | 0 | 0 | 2.391 | 2.325 | 0.439 | -1.561 | 0 | 0 | 0 | 0 | 0 | 0 | 0 |
| 0 | 0 | 0 | 0 | 0.090 | 0.084 | -1.958 | 0.042 | 0 | 0 | 0 | 0 | 0 | 0 | 0 |
| 0 | 0 | 0 | 0 | 0.771 | 0.935 | -3.712 | -0.501 | 0 | 0 | 0 | 0 | 0 | 0 | 0 |
| 0 | 0 | 0 | 0 | -0.760 | 0.760 | -0.760 | -0.760 | 0 | 0 | 0 | 0 | 0 | 0 | 0 |
| -0 | -0 | -0 | -0 | -0 | -0 | -0 | -0 | 0.642 | 0.646 | -0.145 | -2.145 | -0 | -0 | -0 |
| 0 | 0 | 0 | 0 | 0 | 0 | 0 | 0 | 1.459 | 1.452 | -1.755 | 0.245 | 0 | 0 | 0 |
| 0 | 0 | 0 | 0 | 0 | 0 | 0 | 0 | 0.146 | 0.427 | -0.735 | 1.031 | 0 | 0 | 0 |
| 0 | 0 | 0 | 0 | 0 | 0 | 0 | 0 | -0.173 | 0.171 | -0.171 | -0.171 | 0 | 0 | 0 |
| -0.040 | -0.031 | 0.021 | 0.537 | 0.202 | 0.205 | 0.031 | -0.144 | -0.091 | -0.083 | 0.021 | 0.306 | 0.002 | 0 | -0.001 |
| -0.149 | -0.152 | 0.207 | 0.018 | 0.052 | 0.044 | -1.093 | 0.025 | 0.161 | 0.168 | -0.208 | 0.021 | 0 | 0.002 | -0.001 |
| 0.116 | 0.135 | -0.174 | 0.182 | 0.170 | 0.165 | 0.535 | -0.132 | -0.065 | -0.038 | 0.041 | 0.072 | 0.001 | -0.001 | -0 |

Therefore, the invertibility of the matrix (and thus the unisolvence of the element corresponding to the used degrees of freedom) comes down to the invertibility of each of the normal submatrices. Note that in the case of non - invertibility, if one block contains some row that is a linear combinations of the others, it means that the set of degrees of freedom is not robust towards the layout of the corresponding edge. If several columns are related by some linear combination, then the initial basis is not free when built on the considered edge.

However, due to this block construction, if the matrix is not invertible it is rather easy to identify the faulting edge and use the sensitivity of the degrees of freedom towards the layout of the edges to surgically modify the reference element and resolve the issue. Furthermore, we have the insurance that those local unisolvence problems never spread over the matrix but are always limited to their respective blocks. Therefore, even in the case of several problematic edges, all of the non - invertible submatrices are easily recognizable, and stitches can easily be performed. Let us detail more precisely those limitations.

### 7.2.2 Confirmation of failing cases

Let us start by confirming the necessity of the restrictions on the elements' shape we derived in the *Paragraph 6.1.2*. As in this framework the arising problems are limited to the boundary, it is enough to restrict the study to the case $k = 0$ where there is no internal moment.

**Edge whose normal is collinear with the vector $v$** For elements corresponding to the sets of degrees of freedom $Ia$ and $Ib$ (see the *Definitions 5.17* and *5.18* ), the case of edges whose normal is collinear with the vector $v$ used in the definition of the degrees of freedom breaks the unisolvence. Indeed, any function $q$ of $\mathbb{H}_0(K)$ read

$$q\colon \begin{pmatrix} x \\ y \end{pmatrix} \mapsto \begin{pmatrix} x \\ y \end{pmatrix} C + \begin{pmatrix} A \\ B \end{pmatrix}$$

for three real constants $A, B$ and $C$. Therefore, assuming without loss of gen-

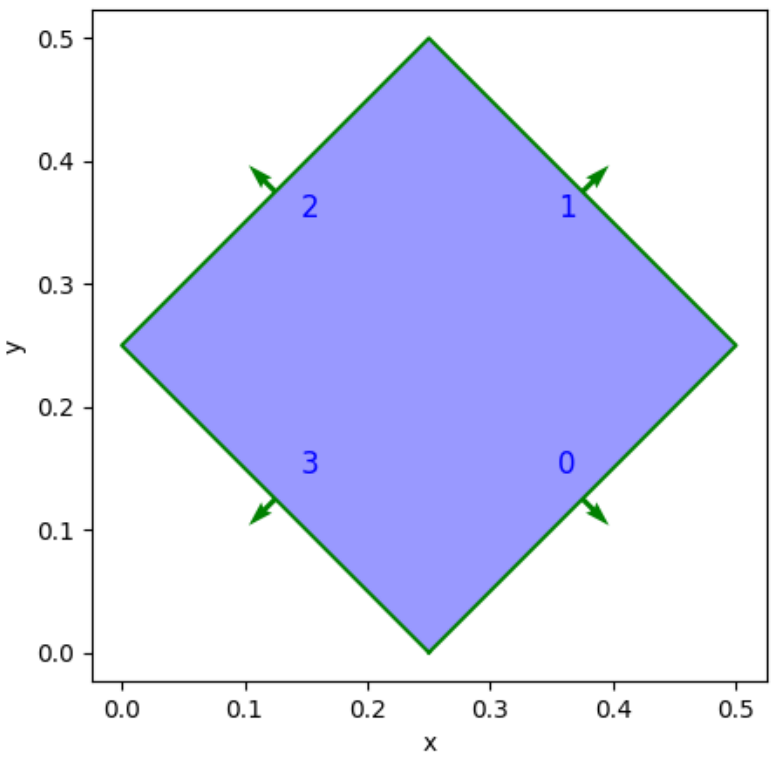

Fig. 74: Reference element where $v = n_3$

erality that $n = v$, the three degrees of freedom read

$$\int_f q \cdot v \, x \, \mathrm{d}\gamma(x) = \int_f \underbrace{(xn_x + yn_y)}_{=1} Cx + An_x + Bn_y \, \mathrm{d}\gamma(x)$$

$$\int_f q_x \, n_x \, \mathrm{d}\gamma(x) = \int_f (xn_x)C + An_x \, \mathrm{d}\gamma(x)$$

$$\int_f q_y \, n_y \, \mathrm{d}\gamma(x) = \int_f (yn_y)C + Bn_y \, \mathrm{d}\gamma(x).$$

By definition, the first degree of freedom quantifies globally the aspects of functions living in $\mathbb{H}_0(K)$, involving their two components tested against $v$. The two last moments are quantifying the functions of $\mathbb{H}_0(K)$ coordinate wise, tested against the corresponding normal component $n_{ix}$ or $n_{iy}$. Thus, as $v$ is collinear to the normal, the global moment can be reconstructed from the two coordinate - wise moments, and is therefore not unisolvent. However, it is worth to notice that the two elements $IIa$ and $IIb$ do not have this limitation.

A example of this failure case can be obtained by considering the element presented in the *Figure 74* while taking $v = (1, 1)$. There, we would obtain the conditionings

| Element | $Ia$ | $Ib$ | $IIa$ | $IIb$ |
|---|---|---|---|---|
| Conditioning | 1.9e+17 | 1.7e+17 | 1097 | 3225, |

basically meaning that the two first matrices are not invertible, and that the basis functions corresponding to the two last elements will be numerically reliable. For reference, we give the transfer matrix of the element $Ib$ where a combination between the three **rows** of the block corresponding to the second and fourth edge can be observed.

| | | | | | | | | | | | |
|---|---|---|---|---|---|---|---|---|---|---|---|
| 0.265 | -1.735 | 1.265 | 0 | 0 | 0 | 0 | 0 | 0 | 0 | 0 | 0 |
| -2.088 | -0.088 | 0.912 | -0 | -0 | -0 | -0 | -0 | -0 | -0 | -0 | -0 |
| 0.199 | -0.132 | 0.034 | 0 | 0 | 0 | 0 | 0 | 0 | 0 | 0 | 0 |
| 0 | 0 | 0 | 0.265 | -1.735 | 1.265 | 0 | 0 | 0 | 0 | 0 | 0 |
| 0 | 0 | 0 | -1.735 | 0.265 | 1.265 | 0 | 0 | 0 | 0 | 0 | 0 |
| 0 | 0 | 0 | -0.122 | -0.122 | 0.210 | 0 | 0 | 0 | 0 | 0 | 0 |
| -0 | -0 | -0 | -0 | -0 | -0 | -0.088 | -2.088 | 0.912 | -0 | -0 | -0 |
| 0 | 0 | 0 | 0 | 0 | 0 | -1.735 | 0.265 | 1.265 | 0 | 0 | 0 |
| 0 | 0 | 0 | 0 | 0 | 0 | -0.132 | 0.199 | 0.034 | 0 | 0 | 0 |
| -0 | -0 | -0 | -0 | -0 | -0 | -0 | -0 | -0 | -0.083 | -2.083 | 0.917 |
| -0 | -0 | -0 | -0 | -0 | -0 | -0 | -0 | -0 | -2.094 | -0.094 | 0.906 |
| 0 | 0 | 0 | 0 | 0 | 0 | 0 | 0 | 0 | 0.181 | 0.181 | -0.151 |

Indeed, let us consider by example the fourth edge.

There, the sum of the two first rows give us

$$-2.177 \quad -2.77 \quad 1.823$$

which divided by $-12.02$ leads to

$$0.181 \quad 0.181 \quad -0.151,$$

which equals the last row.

| | | |
|---|---|---|
| -0.083 | -2.083 | 0.917 |
| -2.094 | -0.094 | 0.906 |
| 0.181 | 0.181 | -0.151 |

Tab. 3: Fourth block of the matrix corresponding to the element $Ib$ for $k = 1$.

This observation also holds computationally up to the machine's precision. The restriction either on the shape of the element or on the definition of $v$ is therefore justified.

**Edges aligned with the axes** The second restriction applies to all the considered elements. Indeed, in any case the collinearity of some edge with one of the axes forces at least one component - wise quantifier to reduce to the zero function. In particular, any tested canonical basis function will reduce to a zero contribution through those degrees of freedom, and the transfer matrix will comprise empty **rows**.

Thus, no matrix then being full ranked, none of the previously defined sets of degrees of freedom will be unisolvent for $\mathbb{H}_k(K)$. As a consequence, no such element can be defined on those shapes.

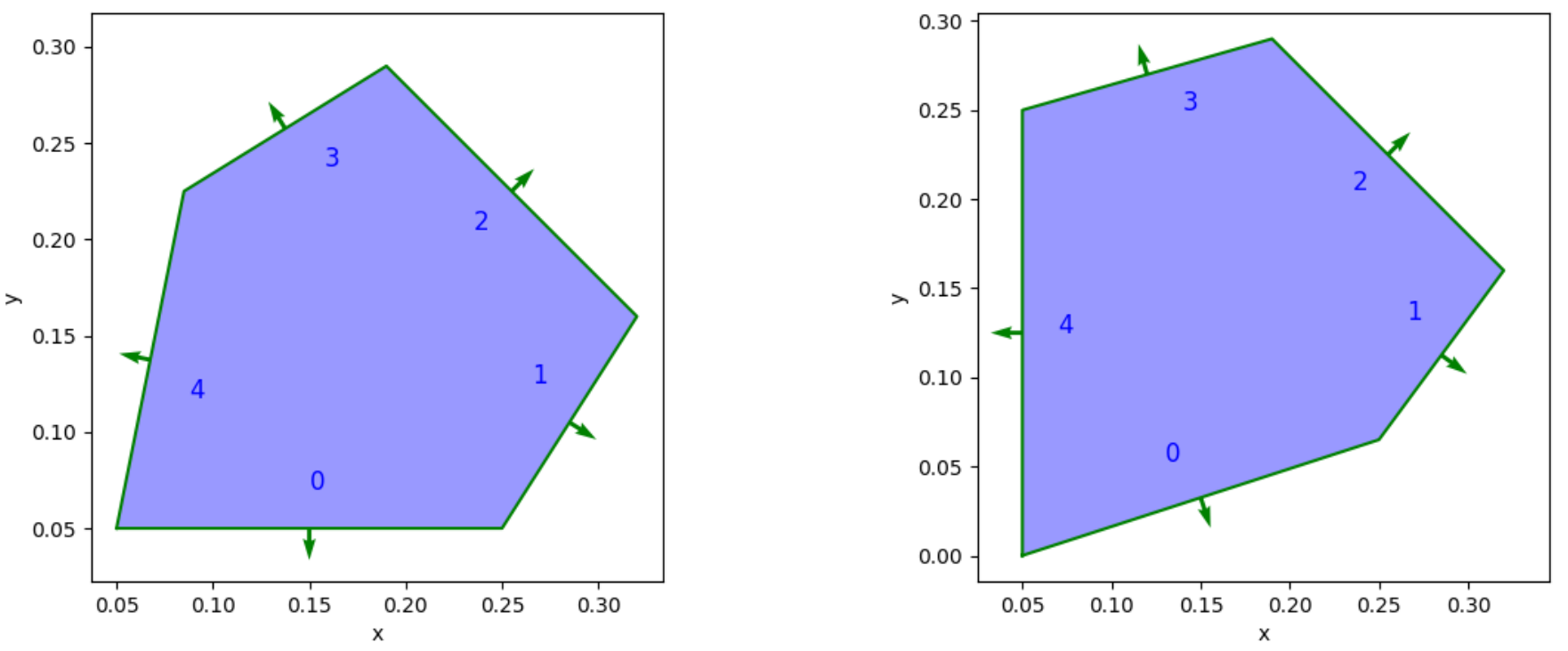

(a) Shape within one edge is aligned with the axis $x$

(b) Shape within one edge is aligned with the axis $y$

Fig. 75: Examples of elements have edges collinear to one of the axes.

By example, considering the shapes given in the *Figures 75a* and *75b*, computing the degrees of freedom that corresponds to the element $IIb$ leads respectively to the following transfer matrices.

| | | | | | | | | | | | | | | |
|---|---|---|---|---|---|---|---|---|---|---|---|---|---|---|
| **0** | **0** | **0** | **0** | **0** | **0** | **0** | **0** | **0** | **0** | **0** | **0** | **0** | **0** | **0** |
| -0.020 | -0.000 | 0.019 | 0 | 0 | 0 | 0 | 0 | 0 | 0 | 0 | 0 | 0 | 0 | 0 |
| -2.050 | -0.050 | 1.950 | 0 | 0 | 0 | 0 | 0 | 0 | 0 | 0 | 0 | 0 | 0 | 0 |
| 0 | -0 | 0 | 0.008 | -0.060 | 0.060 | 0 | 0 | 0 | -0 | 0 | 0 | -0 | -0 | 0 |
| -0 | -0 | 0 | 0.001 | -0.000 | 0.006 | 0 | 0 | 0 | -0 | 0 | 0 | -0 | 0 | 0 |
| 0 | 0 | 0 | 0.371 | -1.725 | 2.184 | 0 | 0 | 0 | 0 | 0 | 0 | 0 | 0 | 0 |
| 0 | 0 | 0 | 0 | 0 | 0 | 0.008 | -0.070 | 0.054 | 0 | 0 | 0 | 0 | 0 | 0 |
| 0 | 0 | 0 | 0 | 0 | 0 | -0.062 | 0.006 | 0.047 | 0 | 0 | 0 | 0 | 0 | 0 |
| 0 | 0 | 0 | 0 | 0 | 0 | -1.367 | -1.367 | 2.339 | 0 | 0 | 0 | 0 | 0 | 0 |
| 0 | 0 | 0 | 0 | 0 | 0 | 0 | 0 | 0 | -0.001 | 0.002 | 0.008 | 0 | 0 | 0 |
| 0 | 0 | 0 | 0 | 0 | 0 | 0 | 0 | 0 | -0.052 | 0.006 | 0.052 | 0 | 0 | 0 |
| 0 | 0 | 0 | 0 | 0 | 0 | 0 | 0 | 0 | -1.770 | 0.353 | 2.146 | 0 | 0 | 0 |
| 0 | 0 | 0 | 0 | 0 | 0 | 0 | 0 | 0 | 0 | 0 | 0 | -0.0008 | -0.021 | 0.021 |
| 0 | 0 | 0 | 0 | 0 | 0 | 0 | 0 | 0 | 0 | 0 | 0 | -0.0098 | 0.000 | 0.002 |
| 0 | 0 | 0 | 0 | 0 | 0 | 0 | 0 | 0 | 0 | 0 | 0 | -0.4933 | -1.886 | 1.960 |
| 0.001 | -0.021 | 0.007 | 0 | 0 | 0 | 0 | 0 | 0 | 0 | 0 | 0 | 0 | 0 | 0 |
| -0.011 | -0.000 | 0.012 | 0 | 0 | 0 | 0 | 0 | 0 | 0 | 0 | 0 | 0 | 0 | 0 |
| -1.611 | -0.734 | 2.015 | 0 | 0 | 0 | 0 | 0 | 0 | 0 | 0 | 0 | 0 | 0 | 0 |
| 0 | 0 | 0 | 0.007 | -0.053 | 0.049 | 0 | 0 | 0 | 0 | 0 | 0 | 0 | 0 | 0 |
| 0 | 0 | 0 | -0.000 | -0.000 | 0.008 | 0 | 0 | 0 | 0 | 0 | 0 | 0 | 0 | 0 |
| 0 | 0 | 0 | 0.227 | -1.700 | 2.162 | 0 | 0 | 0 | 0 | 0 | 0 | 0 | 0 | 0 |
| 0 | 0 | 0 | 0 | 0 | 0 | 0.008 | -0.070 | 0.054 | 0 | 0 | 0 | 0 | 0 | 0 |
| 0 | 0 | 0 | 0 | 0 | 0 | -0.062 | 0.006 | 0.047 | 0 | 0 | 0 | 0 | 0 | 0 |
| 0 | 0 | 0 | 0 | 0 | 0 | -1.367 | -1.367 | 2.339 | 0 | 0 | 0 | 0 | 0 | 0 |
| -0 | 0 | -0 | -0 | 0 | -0 | -0 | 0 | -0 | -0.000 | 0.005 | 0.002 | -0 | -0.000 | 0.000 |
| 0 | 0 | -0 | 0 | 0 | -0 | -0 | 0 | 0 | -0.067 | 0.010 | 0.081 | -0.006 | 0.001 | 0.001 |
| 0 | 0 | 0 | 0 | 0 | 0 | 0 | 0 | 0 | -1.767 | 0.583 | 2.226 | 0 | 0 | 0 |
| -0.0001 | 0.0060 | -0.0017 | 0 | 0 | 0 | 0 | 0 | 0 | -0 | 0.001 | 0.000 | -0.000 | -0.025 | 0.024 |
| **0** | **0** | **0** | **0** | **0** | **0** | **0** | **0** | **0** | **0** | **0** | **0** | **0** | **0** | **0** |
| 0 | 0 | 0 | 0 | 0 | 0 | 0 | 0 | 0 | 0 | 0 | 0 | -0.050 | -2.050 | 1.950 |

***Note.***

• In practice, the vectors $e_{1,3}$ and $e_{1,4}$ cannot be designed in such cases. For the sake of the example we computed the basis functions with the equivalent vectors (1, 0) and (0, 1) instead.

• The non - vanishing terms in the bottom left corner of the second presented matrix, corresponding to the element (75b) only corresponds to numerical defects in the integration method due to the span of a vertical edge. ▲

**Vertices aligned with the origin**

The third shape restriction that also impacts all the presented elements involves the alignment of the edges' vertices with the origin of the axes. Indeed, for all edges $f_i$ whose vertices line up with the origin, the constant $\mathrm{x} \cdot n$ vanishes. Therefore, all the normal components of the basis functions that are built from the vectors $e_{i,2}$, $e_{i,3}$ and $e_{i,4}$ coincide.

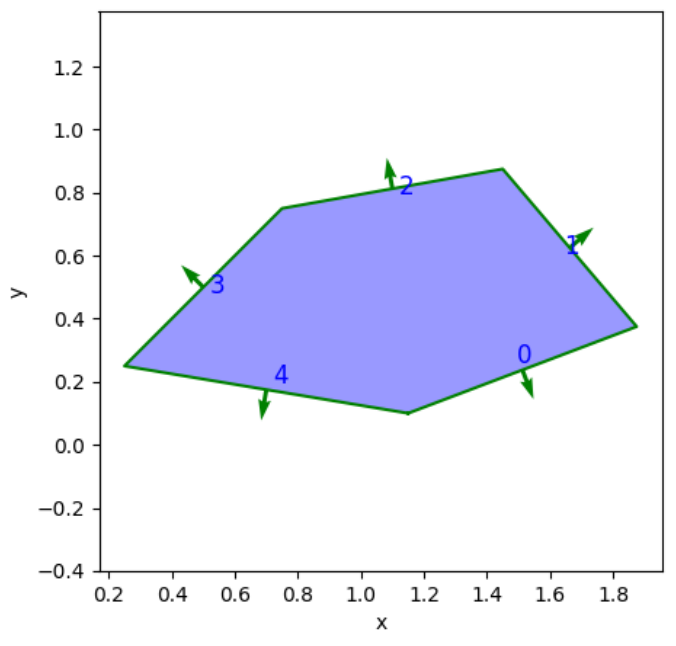

Fig. 76: Aligned origin

Furthermore, as this issue comes from the co - occurrence of our construction of canonical basis functions and degrees of freedom involving normal components, no $H(\mathrm{div}, K)$ – conformal element can be defined through the technique given in the *Section 2.4*. To give an example, let us consider the element pictured in the *Figure 76*.

First of all, one can notice that the term $\mathrm{x} \cdot n$ vanishes for any x belonging to the problematic edge. Therefore, the basis functions built from Lagrangian functions whose support match this edge naturally loose one polynomial order on the boundary (*see the Figure 77*).

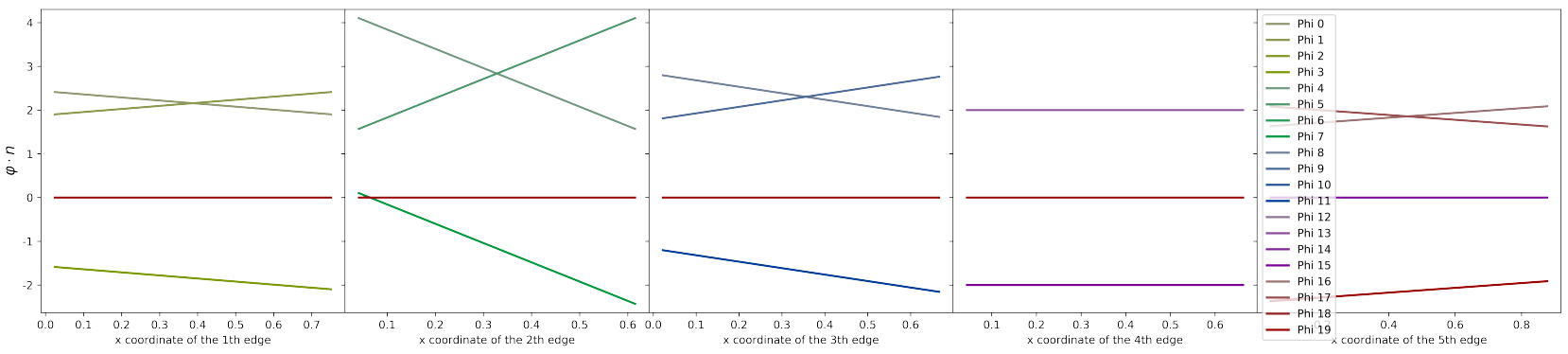

Fig. 77: Normal component of the canonical basis functions when the fourth edge has his vertices collinear with the origin.

As a consequence, the basis functions are not free and no element can be built in the way that we presented here. Indeed, in that case several **rows** will be linearly dependent and some **columns** will share important similarities. To emphasise it, let us consider the first block of the transfer matrix corresponding to the element $IIb$ for $k = 1$, computed here from the canonical polynomial projection spaces. We retrieve

| | | | |
|---|---|---|---|
| 0.0141 | 0.0.147 | -0.0003 | -0.0275 |
| 0.136 | 0.0131 | -0.0280 | 0.0003 |
| 0.0331 | 0.0331 | -0.0331 | -0.0331 |
| 0.2772 | 0.2772 | -0.2772 | -0.2772. |

There, the similarity between the first and second column as well as the linear dependence between the third and fourth row can be clearly established. When considering the same element $IIb$ generated from the projection space discussed in all of this section, based on Hermite polynomials, the following matrix is retrieved instead.

| | | | |
|---|---|---|---|
| 0.275 | 0.227 | -0.034 | -0.653 |
| 0.343 | 0.392 | -0.585 | 0.034 |
| -0 | -0 | 0 | 0 |
| 1.237 | 1.237 | -1.237 | -1.237 |

The lack of unisolvence is then emphasised through a vanishing row, being a consequence of the definition of the Hermite's polynomials.

This behaviour is proper to each edge and impacts only the respective matrix blocks. Thus, even when more that one problematic edges are encountered, this problem does not get global. By example, taking care of the element given in the *Figure 78*, we end up with the following behaviour where the problematic is restricted to the two mis - designed edges.

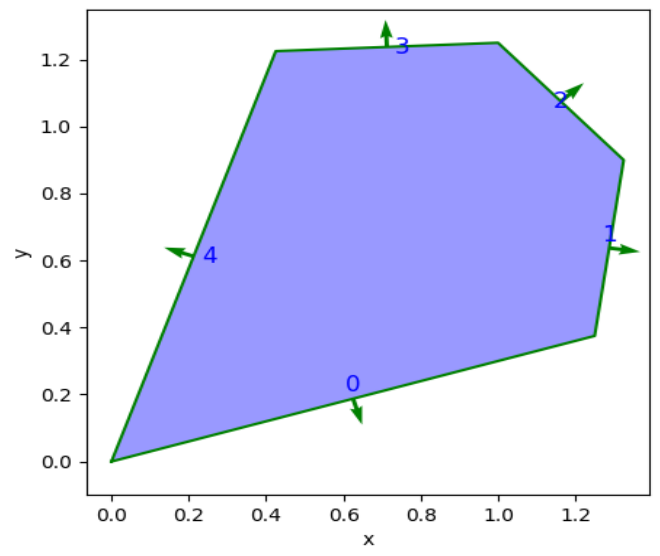

Fig. 78: Element with two problematic edges

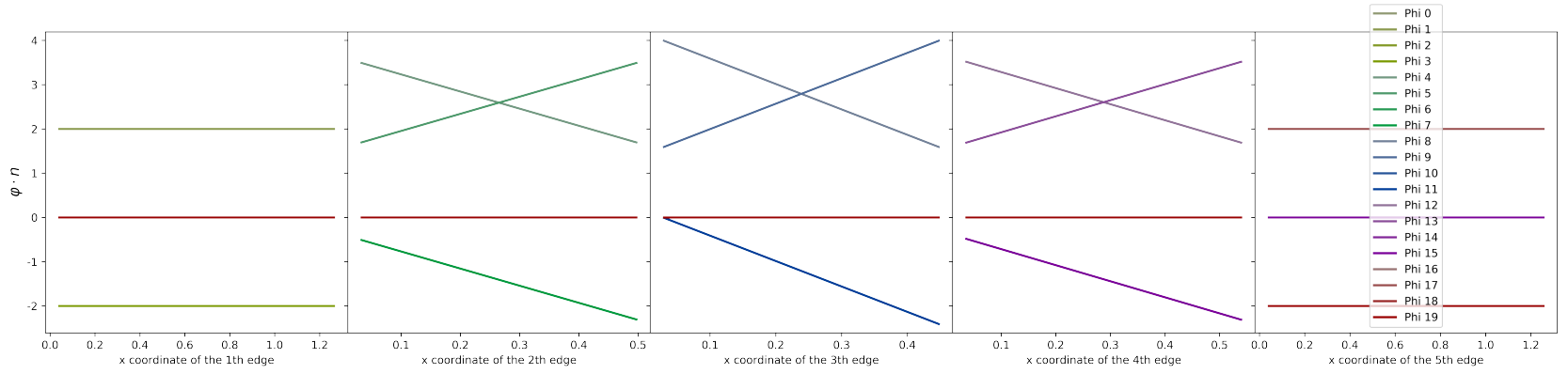

Fig. 79: Canonical basis function when two edges are collinear with the origin

***Remark.*** *Important remark* The value of the term $\mathrm{x} \cdot n$ is also related to the conditioning of the transfer matrix. Indeed, connecting with the end of the *Paragraph 7.1.1*, even if the coefficient $\mathrm{x} \cdot n$ is not vanishing, it being small is enough to create similarities between the basis functions. Thus, the smaller this coefficient is, worse is the conditioning of the transfer submatrix associated to the corresponding edge. ▲

Let us now detail the limit cases where elements are definable, but where the matrices share some similarities in between the normal blocks.

### 7.2.3 Limit cases

As the transfer matrix benefits from a block layout, the presence of similar edges in terms of length and orientation (collinear normals) in the polygonal shape does not prevent the definition of any of the elements $Ia$, $Ib$, $IIa$ or $IIb$. Indeed, each block being independent from another, the presence of similar diagonal blocks does not impact the invertibility of the full matrix. Taking into consideration the shape limitations of the previous section is then enough to ensure the unisolvence of all the discussed elements.

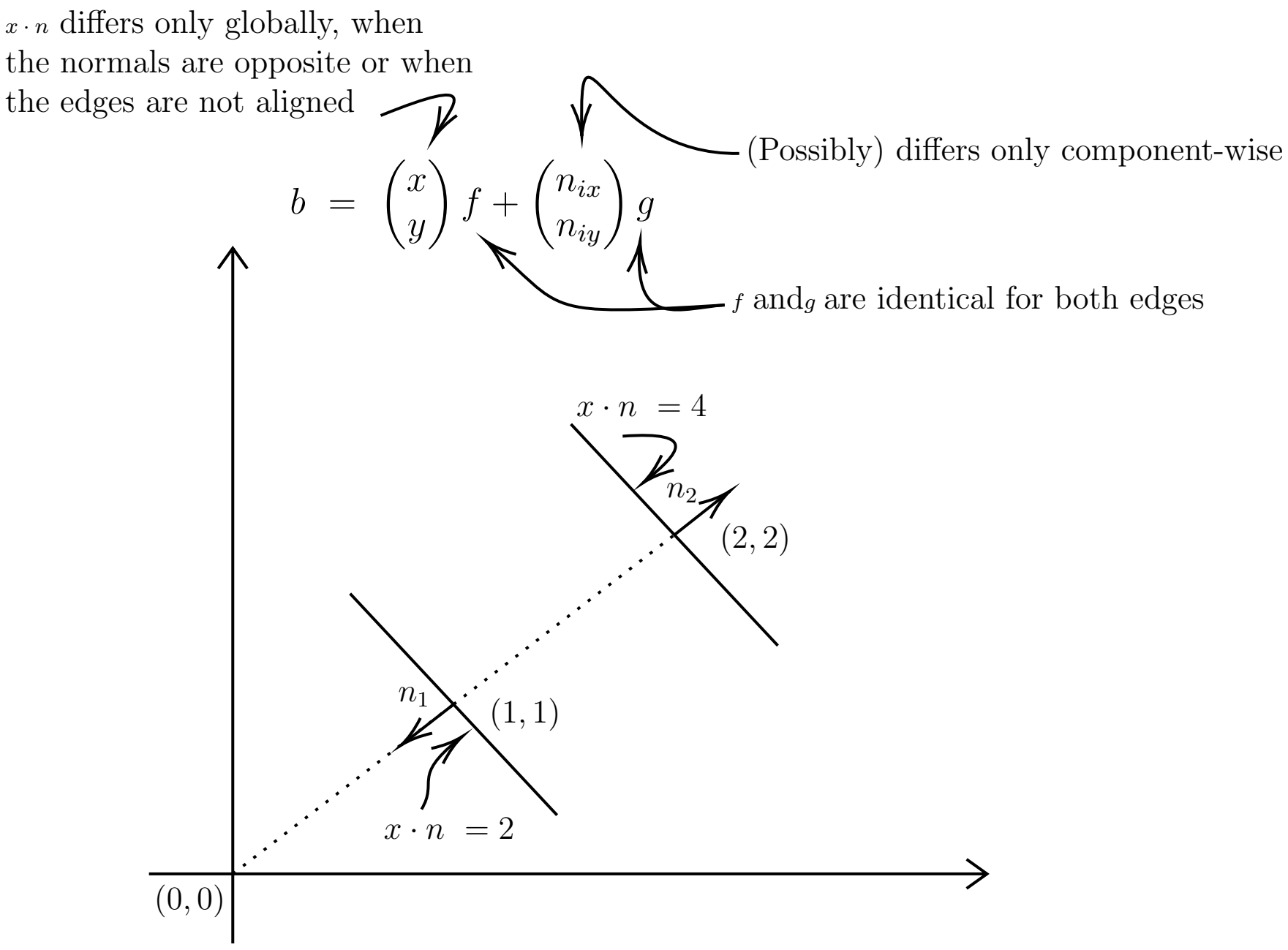

Fig. 80: Sketch of the basis functions layout for two parallel edges having the same length.

However, the restriction of the underneath canonical basis functions to the edges differ only from their support and from the values of $n$ and $\mathrm{x} \cdot n$ (*see the Figure 80*). Thus, if the polygonal shape contains short similar edges that are close to each other, there will be some similarity between the two groups of basis functions emerging from the similar edges. Furthermore, the dissimilarities between the basis functions when seen through the sets of degrees of freedom can be tiny.

**Parallel similar edges** In the case of elements having parallel edges whose length is identical, the only difference in the normal component of the basis functions reside as always in the coefficient $\mathrm{x} \cdot n$ and possibly in the orientation of the term $(n_{ix}, n_{iy})^T$. Therefore, if the constant $\mathrm{x} \cdot n$ is small, similarities in the matrix blocks will appear.

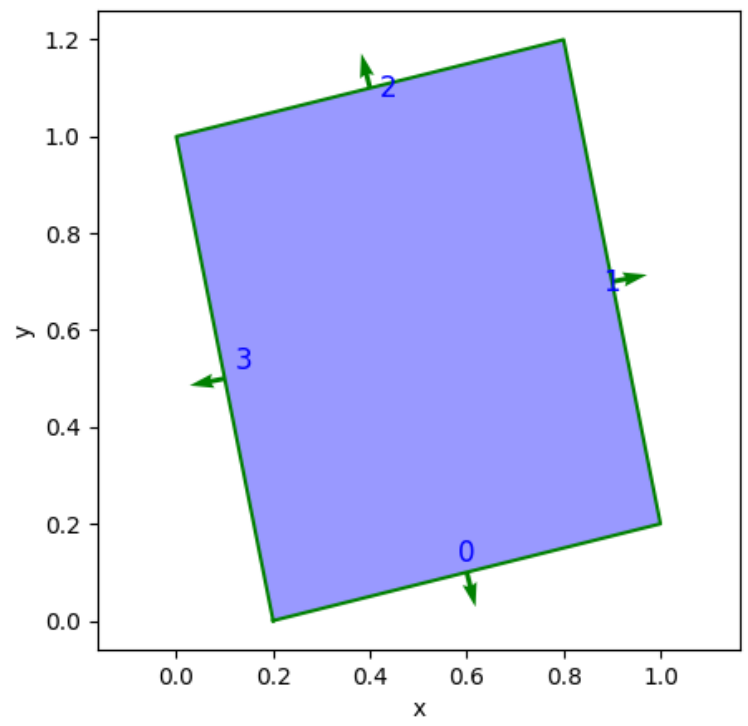

Fig. 81: Element of reference having parallel edges

By example, let us consider the shape represented in the *Figure 81*, where the four elements can be defined. The conditionings of their respective transfer matrix read

| Element | $Ia$ | $Ib$ | $IIa$ | $IIb$ |
|---|---|---|---|---|
| Conditioning | 4000 | 234 | 3373 | 3375, |

which are decent values indicating that the created elements will be reliable. Furthermore, when deriving the transfer matrix corresponding to the element $Ib$ (*see below*), one can observe that the small differences in the coefficients $\mathrm{x} \cdot n$ allow the achievement of two by two different blocks, despite the similar couples of edges (1, 3) and (2, 4). Indeed, the shift in the $\mathrm{x} \cdot n$ value combined with the normals having opposite sign hides the edge's similarities.

| | | | | | | | | | | | |
|---|---|---|---|---|---|---|---|---|---|---|---|
| 0.146 | -1.854 | 0.263 | 0 | 0 | 0 | 0 | 0 | 0 | 0 | 0 | 0 |
| -2.097 | -0.097 | 1.785 | -0 | -0 | -0 | -0 | -0 | -0 | -0 | -0 | -0 |
| 0.961 | -2.434 | -0.197 | 0 | 0 | 0 | 0 | 0 | 0 | 0 | 0 | 0 |
| 0 | 0 | 0 | 0.879 | -1.121 | 2.803 | 0 | 0 | 0 | 0 | 0 | 0 |
| 0 | 0 | 0 | -1.860 | 0.140 | 0.217 | 0 | 0 | 0 | 0 | 0 | 0 |
| 0 | 0 | 0 | -4.466 | -0.291 | 1.957 | 0 | 0 | 0 | 0 | 0 | 0 |
| -0 | -0 | -0 | -0 | -0 | -0 | -0.097 | -2.097 | 0.021 | -0 | -0 | -0 |
| 0 | 0 | 0 | 0 | 0 | 0 | -0.933 | 1.067 | 2.950 | 0 | 0 | 0 |
| 0 | 0 | 0 | 0 | 0 | 0 | -0.133 | 3.262 | 1.025 | 0 | 0 | 0 |
| -0 | -0 | -0 | -0 | -0 | -0 | -0 | -0 | -0 | -0.104 | -2.104 | 1.819 |
| -0 | -0 | -0 | -0 | -0 | -0 | -0 | -0 | -0 | -2.092 | -0.092 | -0.015 |
| 0 | 0 | 0 | 0 | 0 | 0 | 0 | 0 | 0 | 5.460 | 1.285 | -0.964 |

However, the intrinsic basis functions emerging from the coupled edges are similar and the order of magnitude of the matrix coefficients are comparable to each other. Therefore, the tuning of the basis functions towards the degrees of freedom defining the four elements will not enforce a clearer distinction between the couples of initially similar canonical basis functions. This can be bothersome for the solutions' representation.

**Similar aligned edges** Another limit case corresponds to two aligned similar edges within the same element. There, the obtained rows corresponding to global normal components of blocks linked to similar edges will be identical.

By example, we can consider the element presented in the *Figure 82*, having the edges one, five and seven that are aligned and similar, sharing the following properties.

| Edge | Normal ($n_x$, $n_y$) | Norm |
|---|---|---|
| 3,5,7 | (-0.24, 0.97) | 0.82 |

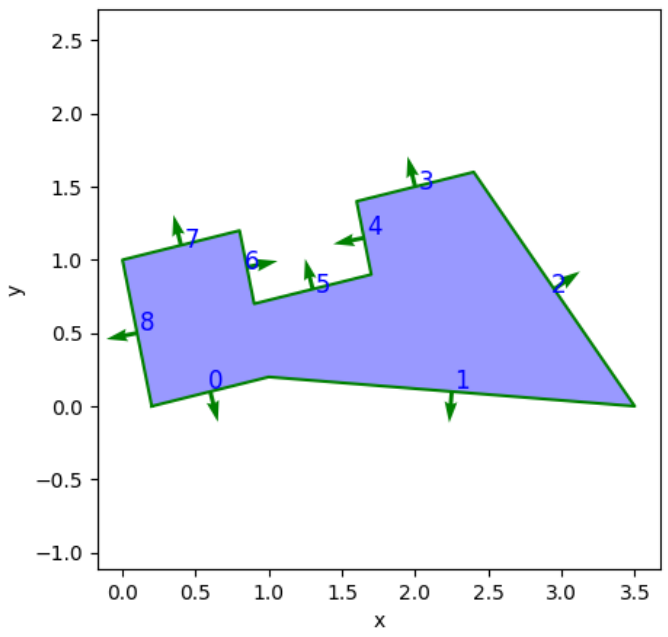

Fig. 82: Element with similar aligned edges.

In this case, as the basis functions have their boundary support which is distinct from each other, the block layout of the transfer matrix is preserved and the conditioning is equally decent.

| Element | Ia | Ib | IIa | IIb |
|---|---|---|---|---|
| Conditioning | 32345 | 25075 | 42144 | 31312 |

For reference, we transcript hereafter the third, fifth and seventh blocks of the transfer matrix corresponding to the element *IIb*.

| | | |
|---|---|---|
| -0.707 | -3.593 | -0.537 |
| -0.588 | 1.577 | 3.614 |
| **-1.030** | **-1.030** | **2.970** |

Block corresponding to the third edge

| | | |
|---|---|---|
| -0.303 | -2.179 | -0.193 |
| -0.705 | 0.450 | 1.536 |
| -1.539 | -1.539 | 2.461 |

Block corresponding to the fifth edge

| | | |
|---|---|---|
| -0.035 | -0.612 | -0.001 |
| -0.739 | 0.849 | 2.343 |
| **-1.030** | **-1.030** | **2.970** |

Block corresponding to the seventh edge

The similarity between the third rows – expressing the global normal degree of freedom – of the matrix's block corresponding to the edges three and seven can be observed. Furthermore, notice that as in the previous case, the block corresponding to the fifth edge, which is similar but not aligned with the other two, do not share the same row.

**Hanging node** A natural consequence of the case presented above concerns the hanging node. Indeed, a hanging node can be seen as a vertex delimiting two consecutively aligned edges having a same outer normal.

To point out this link, let us consider the polygon represented in the *Figure 83*, where the third and fourth edges are consecutive, aligned and of the same length. As previously, the four elements under consideration can be defined and the conditionings of their corresponding matrices allow to define the tuned basis functions with reliability.

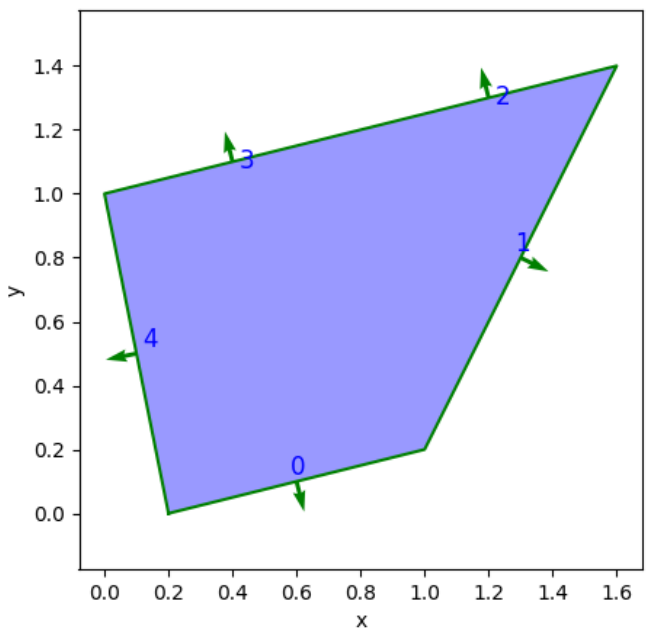

Fig. 83: Element with an hanging node.

| Element | *Ia* | *Ib* | *IIa* | *IIb* |
|---|---|---|---|---|
| Conditioning | 8563 | 714 | 6700 | 6265 |

Deriving by example the transfer matrix corresponding to the element *IIb* one can witness the similarities between the two last rows corresponding to the global normal basis degree of freedom of the two consecutive edges.

| | | | | | | | | | | | | | | |
|---|---|---|---|---|---|---|---|---|---|---|---|---|---|---|
| 0.070 | -0.796 | 0.121 | 0 | 0 | 0 | 0 | 0 | 0 | 0 | 0 | 0 | 0 | 0 | 0 |
| -0.153 | -0.009 | 0.127 | 0 | 0 | 0 | 0 | 0 | 0 | 0 | 0 | 0 | 0 | 0 | 0 |
| -1.951 | -1.951 | 2.049 | 0 | 0 | 0 | 0 | 0 | 0 | 0 | 0 | 0 | 0 | 0 | 0 |
| 0 | 0 | 0 | 1.918 | -1.341 | 4.526 | 0 | 0 | 0 | 0 | 0 | 0 | 0 | 0 | 0 |
| 0 | 0 | 0 | -2.402 | -0.412 | -0.013 | 0 | 0 | 0 | 0 | 0 | 0 | 0 | 0 | 0 |
| 0 | 0 | 0 | -1.195 | -1.195 | 2.805 | 0 | 0 | 0 | 0 | 0 | 0 | 0 | 0 | 0 |
| 0 | 0 | 0 | 0 | 0 | 0 | -0.259 | -1.991 | -0.157 | 0 | 0 | 0 | 0 | 0 | 0 |
| 0 | 0 | 0 | 0 | 0 | 0 | -0.691 | 1.185 | 2.950 | 0 | 0 | 0 | 0 | 0 | 0 |
| 0 | 0 | 0 | 0 | 0 | 0 | **-1.030** | **-1.030** | **2.970** | 0 | 0 | 0 | 0 | 0 | 0 |
| 0 | 0 | 0 | 0 | 0 | 0 | 0 | 0 | 0 | -0.039 | -0.657 | -0.002 | 0 | 0 | 0 |
| 0 | 0 | 0 | 0 | 0 | 0 | 0 | 0 | 0 | -0.791 | 0.910 | 2.510 | 0 | 0 | 0 |
| 0 | 0 | 0 | 0 | 0 | 0 | 0 | 0 | 0 | **-1.030** | **-1.030** | **2.970** | 0 | 0 | 0 |
| 0 | 0 | 0 | 0 | 0 | 0 | 0 | 0 | 0 | 0 | 0 | 0 | -0.012 | -0.203 | 0.172 |
| 0 | 0 | 0 | 0 | 0 | 0 | 0 | 0 | 0 | 0 | 0 | 0 | -1.017 | -0.061 | -0.024 |
| 0 | 0 | 0 | 0 | 0 | 0 | 0 | 0 | 0 | 0 | 0 | 0 | -2.196 | -2.196 | 1.804 |

Thus, having a hanging node has for impact to double the information that would have been furnished from a single edge. However, the price to pay is the possible creation of a discontinuity at the hanging node. Indeed, the introduction of a vertex generates a new virtual edge, leading to the development of further freedom in the space $\mathcal{H}_k(\partial K)$.

***Note.*** The case where a hanging node separates two edges of different length do not impact those last observations. Only the similarities of the last block rows would be broken, as the edges would not be similar anymore in the sense of the previous paragraph. ▲

### Edge whose normal is collinear to the vector $(x - 0,\, y - 0)^T$

Though being disconnected from the layout of the transfer matrix, let us derive an important last remark which is useful when the tuned basis functions are wished to enjoy the same order of magnitude. Indeed, we can note that having one edge whose normal is collinear with the vector $(x - 0,\, y - 0)^T$ for some x lying on that edge creates a disequilibrium in the amplitude of the tuned basis function.

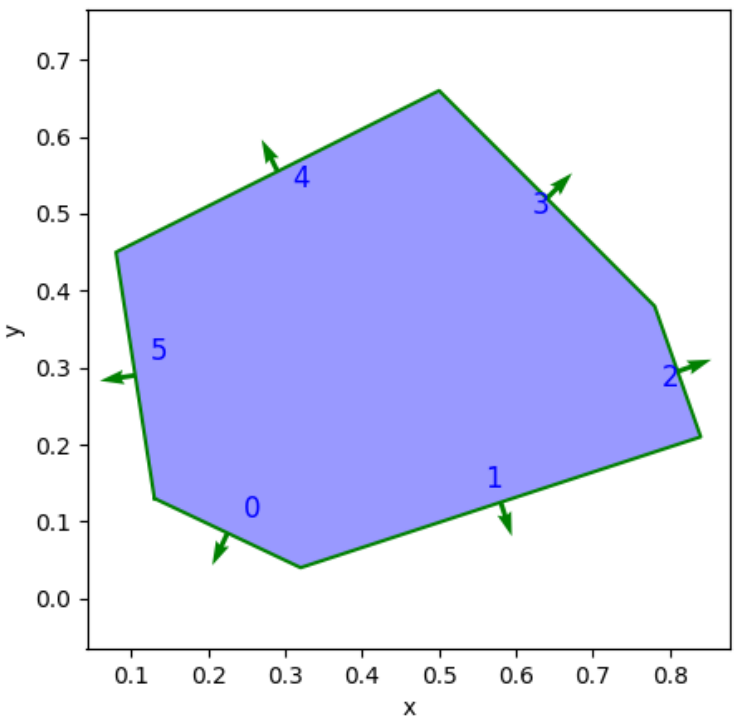

Fig. 84: Element with a problematic edge

By example, when one takes the element *Ib* constructed on the shape given in the *Figure 84*, the coefficients $\{\mathrm{x} \cdot n_i\}_{i \in [\![1,\, 6]\!]}$ are very unbalanced and one ends up having very high variations of the tuned basis functions.

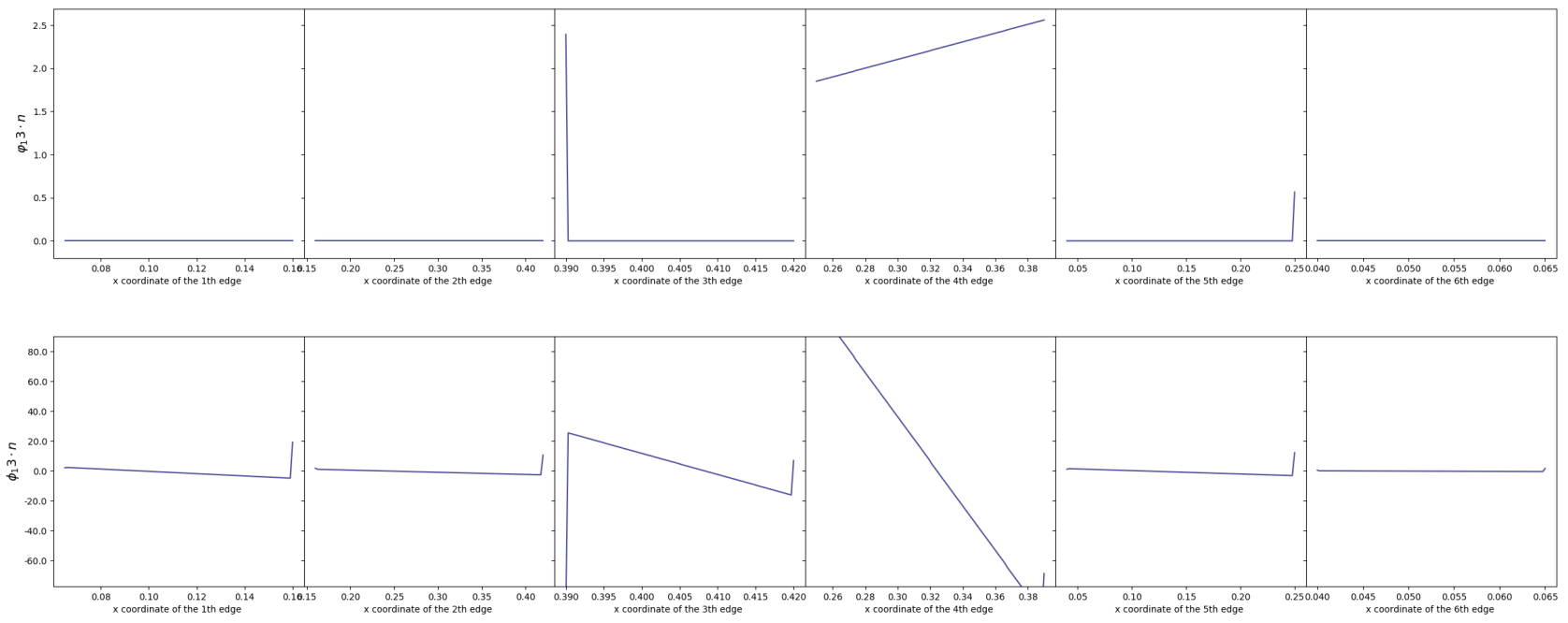

Fig. 85: Amplitude of the basis functions through Dofs when the face normal is collinear with the x vector. Top: canonical. Bottom: Tuned against Ib. On the third edge, the very high neighbouring constant is even numerically creating noisy unwanted variations.

This problematic arises especially when using polynomial projectors that form a basis which is known for its bad conditioning. We represent such a case in the *Figure 85* when the canonical polynomial basis is used as the normal projectors set.

We now assume that all the considered shapes do not fall in the failing cases and investigate the sensitivity of the degrees of freedom towards the geometry.

### 7.2.4 Conditioning of the transfer matrix with respect to the shape of the element

We study here the sensitivity of the transfer matrix towards the shape, scale and order of the considered elements. In particular, both convex and non - convex star shaped cases are under our scope.

**Impact of the shape of the element on the transfer matrix** Let us first investigate the impact of the elements' shape on the transfer matrix. To this sake, we select four shapes having an increasing number of edges. The details of the test cases can be found in the tables given in the *Appendix B*. Then, using comparable edge lengths we test the conditionings of the transfer matrices corresponding to the four elements $Ia$, $Ib$, $IIa$ and $IIb$ built on the convex shapes presented in the *Figure 86* and on non - convex shapes given in the *Figure 87*.

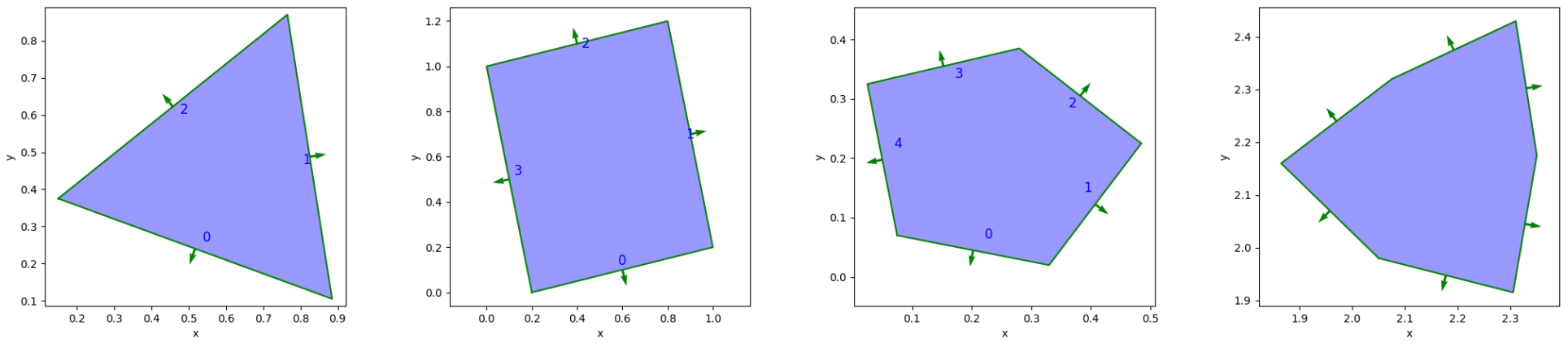

(a) Triangle (b) Quadrilateral (c) Pentagon (d) Hexagon

Fig. 86: Considered convex shapes

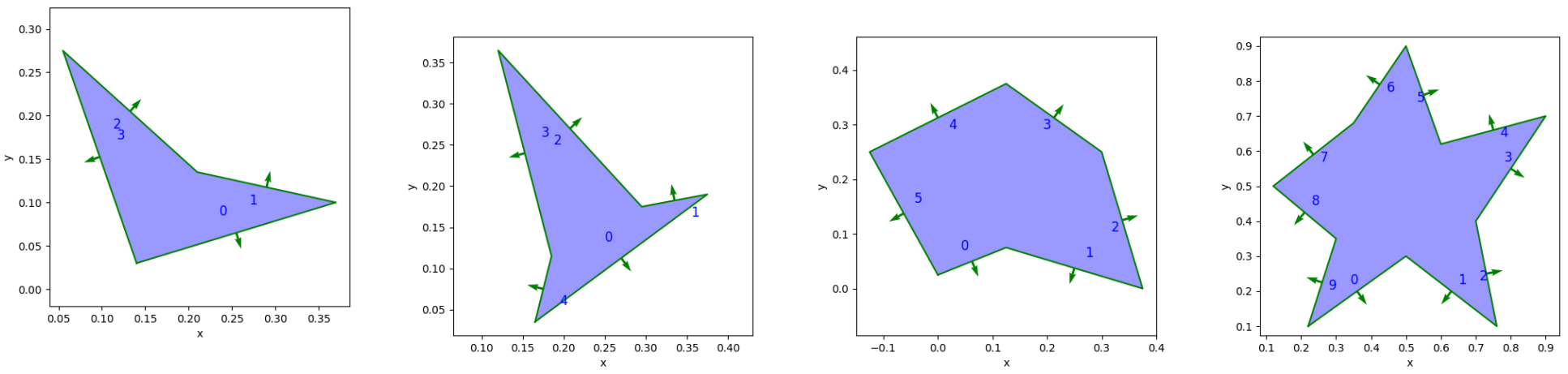

(a) Non - convex quadrilateral (b) Non - convex pentagon (c) Non - convex hexagon (d) Star shaped decagon

Fig. 87: Considered non - convex shapes

Note that as the shape impacts the internal degrees of freedom only by its size, squeezeness and by the quality of the underneath Laplacian solver, the obtained variations in the conditionings of the internal submatrix between two different shapes of comparable areas are negligible for the most common applications. We then consider only the case $k = 0$ and derive the results in the *Table 4* for the convex cases and in the *Table 5* for the star - shaped elements.

| Element | Triangle | Quadrilateral | Pentagon | Hexagon |
|---|---|---|---|---|
| Ia | 637 | 4000 | 5836 | 3805501 |
| Ib | 185 | 234 | 7905 | 245708 |
| IIa | 462 | 3373 | 6693 | 54473 |
| IIb | 556 | 3375 | 27371 | 78111 |

Tab. 4: Conditionings of the elements constructed on the shapes given in the *Figure 86*.

First of all, one can observe that all of the obtained conditionings are reasonable and allow to define reliable elements.

In the convex case, we mainly observe a worsening of the conditioning with the increasing number of edges. However, shapes leading to particularly good conditionings can be found even for relatively high number of edges, as for example the pentagon for the elements $IIa$ and $IIb$ (comparatively to the other shapes).

This behaviour reflects that the conditioning is only driven by the collection of edges' layouts within the space, and that an increasing number of edges rises the probability of the shape to encompass edges whose coefficients $\mathrm{x} \cdot n$ is relatively small. Indeed, as the boundary of the shapes are closed, more we have edges higher is the probability to encounter an edge whose coefficient $\mathrm{x} \cdot n$ is low. Furthermore, as evoked previously the conditioning of the normal sub - matrices increases as the coefficient $\mathrm{x} \cdot n$ decreases. Therefore, by the diagonal block construction this bad conditioning will be transferred directly to the transfer matrix itself. Thus, even if no direct correlation can be drawn between the number of edges and the conditioning of the transfer matrix, the likelihood of having bad conditioning for transfer matrices built on shapes having lot of edges is high.

| Element | Quadrilateral | Pentagon | Hexagon | Decagon |
|---|---|---|---|---|
| Ia | 12414 | 135880 | 6836 | 32747 |
| Ib | 36646 | 57292 | 25086 | 22406 |
| IIa | 44730 | 381198 | 42920 | 15849 |
| IIb | 179466 | 1541968 | 170260 | 47571 |

Tab. 5: Conditionings of the elements constructed on the non - convex shapes given in the *Figure 87*

In the non - convex case, the same observations can be drown. However, one may get better results for high number of edges than in the convex case, in particular for the elements $Ia$ and $Ib$ that are more sensitive component - wise. Indeed, at fixed number of edges it is easier for a non - convex element to avoid critical edges layout than for a convex one.

Thus, one can see that our sets of degrees of freedom are acceptable, but not robust enough to our canonical basis functions in order to prevent the conditioning of the transfer matrix to get high. Therefore, if theoretically any number of edges is admissible, numerically we advise to refrain the use until twelve edges in order to prevent too much similarities between the basis functions emerging from the edges whose coefficient $\mathrm{x} \cdot n << 1$.

**Impact of the scale of the element's shape** Let us now investigate how the scale of the element's shape impacts the conditioning of the transfer matrix. For the sake of convenience, we consider the simplicial case and use the two lowest order elements. To allow a comparison with respect to the size of the element's shape, we consider an original triangle motive $T1$ and half, double, quadruple and sextuple its size to form the shapes $T0$, $T2$, $T3$ and $T4$ (*see the Figure 88*).

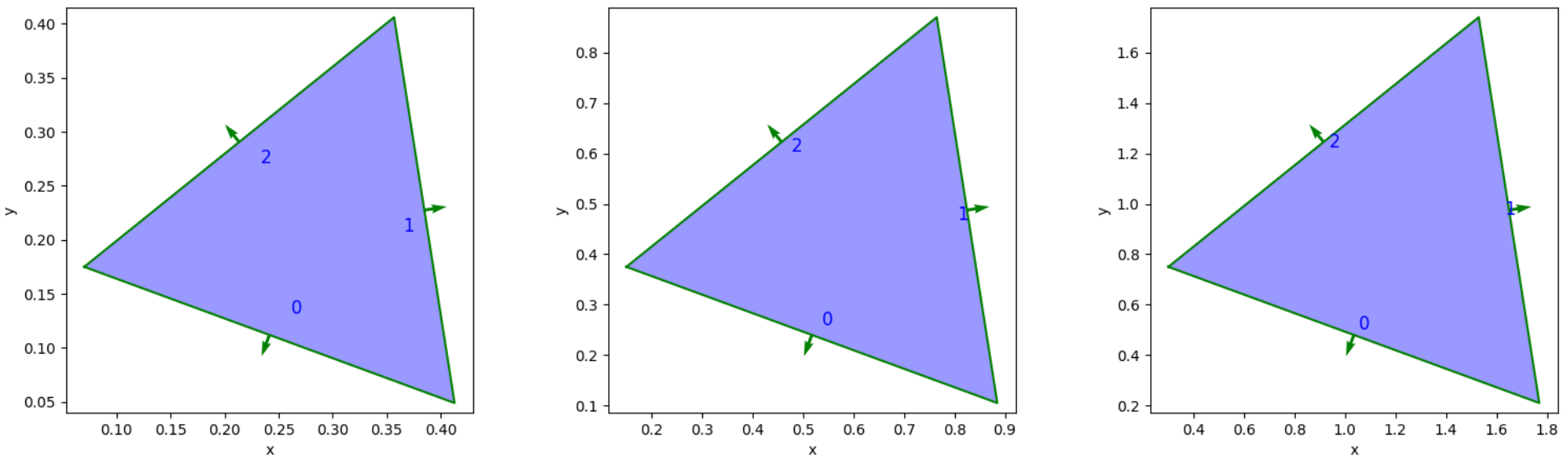

(a) Original triangle, $T1$ (b) Original triangle, $T1$ (c) Doubled size, $T2$

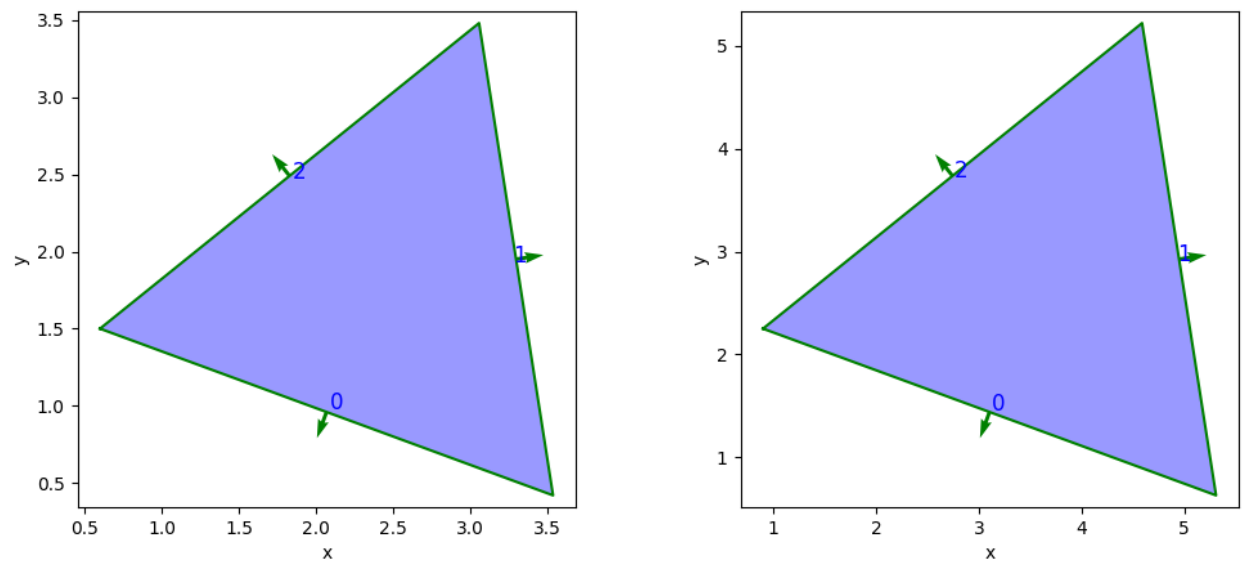

(d) Quadrupled size, $T3$ (e) Sextupled size, $T4$

Fig. 88: Considered triangular shapes of different scales

***Note.*** The polygons $T0$, $T1$ and lie within the unit circle while the polygon $T3$ and $T4$ are out of it. The case $T2$ represents a limit case where its shape mostly lies within the unit circle, but where one edge is totally out of it. ▲

Let us derive the conditioning of the transfer matrices for the lowest order elements in the *Table 6*. There, it can be seen that for the polygons contained within the unit circle, bigger the element is better is the conditioning. And indeed, by construction of the basis functions on the boundary, wider an edge is more the sampling points of the Lagrangian functions are spaced. Therefore, the basis functions are stretcher and the sensitivity of the degrees of freedom reflects better the enhanced dissimilarities of the basis functions.

In particular for the lowest elements case, as all the degrees of freedom are normal ones, longer the edges are more the normal sub - matrices comprise distinguishable rows. Thus, the conditioning of every submatrix gets better, leading to a better conditioning of the transfer matrix itself.

| Element | $T0$ | $T1$ | $T2$ | $T3$ | $T4$ |
|---|---|---|---|---|---|
| Ia | 1238 | 637 | 566 | 2683 | 52763 |
| Ib | 930 | 185 | 331 | 1096 | 9115 |
| IIa | 1791 | 462 | 312 | 1002 | 16532 |
| IIb | 4998 | 566 | 358 | 1986 | 18335 |

Tab. 6: Conditioning of the transfer matrix for $k = 0$

However, when at least an element's edge falls outside of the unit circle, bigger the element is worse is the conditioning. Indeed, in this case the testing polynomials dominates the Lagrangian functions generating the canonical basis functions on the boundary, preventing a clear differentiation between the basis functions of a same block through their set of degrees of freedom' values. This phenomena gets worse as the edges being outside of the unit circle gets wider.

Let us also remark that for the elements $Ib$ and $IIb$, the gain of conditioning when doubling the triangle size is better than the gain obtained for the elements $Ia$ and $IIa$, respectively. Indeed, some of their degrees of freedom are pointwise values that are not averaged with respect to the length of the edge. Thus, provided a good selection of the evaluation point, a clear differentiation between the basis function is easier to perform only through the pointwise values than from an averaged value.

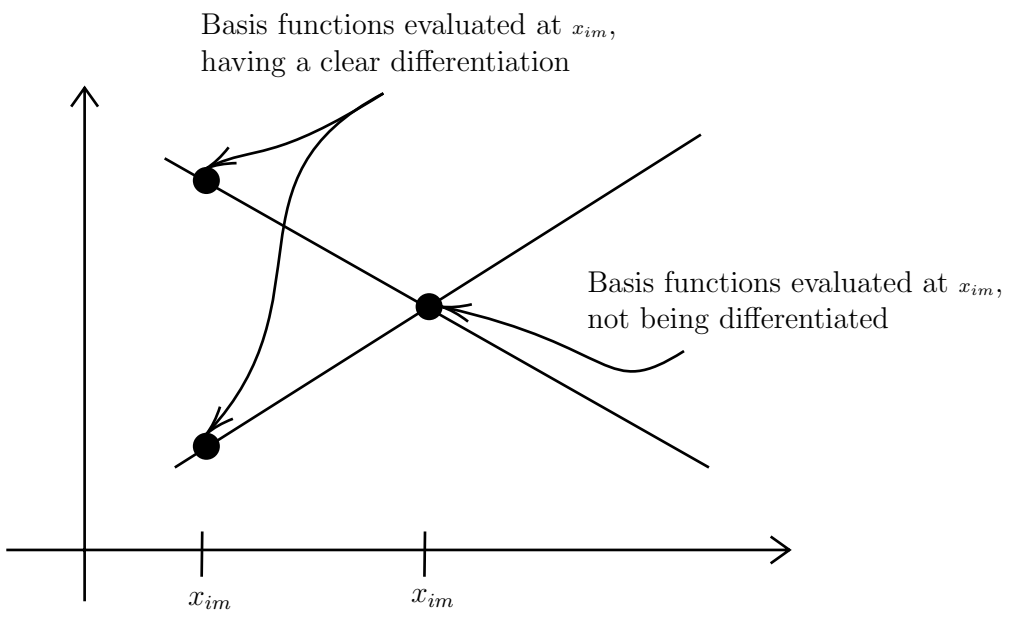

Fig. 89: Example of point wise evaluations of the normal components of two basis functions

However, if the considered point falls in a part of the edge where the normal component of the basis functions are close to each other (typically when they are crossing), then one may loose this advantage or even loose the gain of conditioning one would expect when expanding the element from a smaller shape (*see the Figure 89*).

This is the case by example for the element *Ib* where the characterization performed by one of the two component - wise pointwise values do not allow a clear differentiation between the basis functions. Indeed, for $k = 0$ the pointwise value is evaluated at the mid-point, where the two Lagrangian functions generating the basis functions crosses. Note that this drawback is not easy to stitch as it is very hard to know a - priori which is the best evaluation point to consider. In particular, it depends on the edge's layout, on the used sampling points when creating the Lagrangian functions and on the order of the element.

Same conclusions can be drawn for higher orders. To give an example, we derive the conditioning of the above experiment for the case $k = 1$ in the *Table 7.*

| Element | $T0$ | $T1$ | $T2$ | $T3$ | $T4$ |
|---|---|---|---|---|---|
| Ia | 521340 | 64041 | 16688 | 7156 | 96927 |
| Ib | 1096845 | 40738 | 12681 | 5063 | 8381 |
| IIa | 979645 | 67006 | 17783 | 4110 | 24102 |
| IIb | 2631593 | 78954 | 13148 | 9460 | 19899 |

Tab. 7: Considered elements of reference and respective global conditioning for $k = 1$

There, the conditioning is homogeneously lowering with the increasing size of the element's shape for polynomials contained within the unit circle while increasing for elements lying mainly outside the unit circle, this for any element type. Let us determine if this homogeneous gain/loss is already visible at the levels of internal and normal sub matrices.

| Element | $T0$ | $T1$ | $T2$ | $T3$ | $T4$ |
|---|---|---|---|---|---|
| Ia | 2435 | 1215 | 1198 | 3715 | 77531 |
| Ib | 2191 | 439 | 782 | 1434 | 5754 |
| IIa | 4335 | 1231 | 663 | 1625 | 17148 |
| IIb | 11619 | 1438 | 759 | 3585 | 13901 |

Tab. 8: Considered elements of reference and respective conditioning of the normal submatrix for $k = 1$.

When considering the normal submatrix given in the *Table 8*, one can derive similar observations as in the case $k = 0$. Considering the internal submatrix (neglecting the contributions of normal basis functions) lead to the conditionings given in the *Table 9.*

| Element | $T0$ | $T1$ | $T2$ | $T3$ | $T4$ |
|---|---|---|---|---|---|
| All configurations | 35.0 | 20.3 | 20.8 | 32.3 | 47.6 |

Tab. 9: Considered elements of reference and respective conditioning of the internal submatrix for $k = 1$.

There, we can observe that for shapes contained within the unit circle the conditioning gets better as the shape gets larger. This comes from the construction's method of the canonical basis functions. Indeed, as built from Laplacian solutions with homogeneous Dirichlet conditions, larger the shape is higher the amplitude of the solution is. Thus, by their construction, the values of the internal degrees of freedom retrieved from the basis functions enjoy a larger variance for larger element's shapes. In addition, as the order of magnitude of the degrees of freedom values increase they are less similar and emphases the differences between the considered basis functions with more ease. Thus, the conditioning values are decreasing.

For elements that are not contained within the unit circle, the conditioning is increasing with the size of the element for the same reason as evoked in the case $k = 0$ for the normal degrees of freedom. Indeed, the values of the testing polynomials supersede the ones of the canonical basis functions and prevent from achieving a clear distinction of the basis functions only through their corresponding degrees of freedom values.

Therefore, for the shapes $T0$, $T1$, $T2$ and $T4$ the global observations hold already at the sub-matrices level. The shape $T3$, needs a special attention. If the two submatrices indicates that the conditioning is already worsening, the impact of the normal canonical basis function onto the internal ones lowers the impact of the birth of lack of dissimilarities. Indeed, most of the polygon is still contained within the unit circle and the dissimilarities can be enhanced by combining boundary and inner projections during the tuning process. Therefore, the combination of the two sub-matrices lead to an improvement of the conditioning of the global matrix.

### 7.2.5 Conditioning of the transfer matrix with respect to the order

Let us now discuss the impact that the element's order has on the transfer matrix. As a case study, we consider the non - convex decagon presented in the *Figure 90* on which elements of various orders are built. We derive the conditionings associated with each of the elements in the *Table 10.*

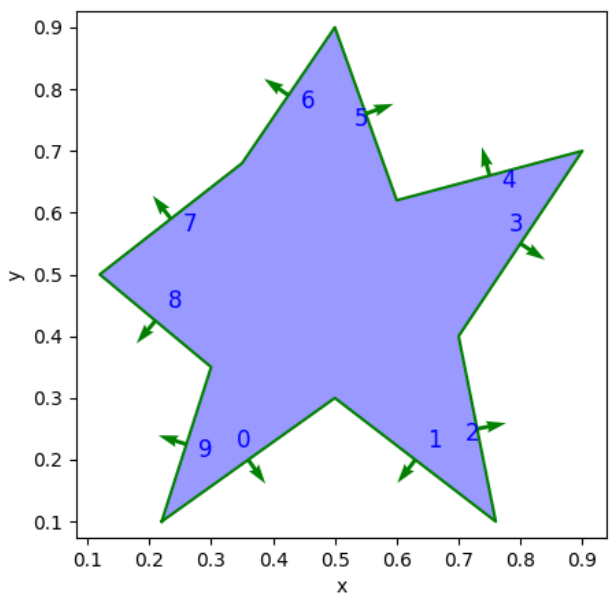

Fig. 90: Non - convex hexagon

A dramatic increase can be observed. Let us determine which of the normal or internal sets of degrees of freedom weakens the most the transfer matrix across the order's increments.

| k | 0 | 1 | 2 | 3 |
|---|---|---|---|---|
| Ia | 32747 | 168806 | 85753753 | 450045582314 |
| Ib | 22406 | 134478 | 84440957 | 480280015315 |
| IIa | 15849 | 99636 | 174815839 | 448778665091 |
| IIb | 47571 | 253343 | 246522510 | 522414830048 |

Tab. 10: Conditionings of the transfer matrix for the several order elements built on the shape 90.

We consider first the conditionings of the normal submatrix presented in the *Table 11*. There, even if the values are increasing, it only gains one order of magnitude per order increment. The weakening of the normal submatrix is then intelligible and high order element could be used provided a full awareness of the fact.

| k | 0 | 1 | 2 | 3 |
|---|---|---|---|---|
| Ia | 32747 | 54276 | 147210 | 2168033 |
| Ib | 22406 | 47126 | 78949 | 117673 |
| IIa | 15849 | 38580 | 126504 | 188478 |
| IIb | 47571 | 99555 | 172541 | 258340 |

Tab. 11: Condition number of the normal submatrix of each element.

However, when considering the conditionings of internal submatrix, it results a dramatical increase presented in the *Table 12*. This extreme weakening comes from the decrease of the latest internal degrees of freedom values within the cell for increasing orders. Indeed, for high orders the values of the degrees of freedom are very small and – for a computer – very similar. Therefore, the conditioning gets high.

| k | 1 | 2 | 3 |
|---|---|---|---|
| Ia | 17.8 | 130224 | 1762489008 |

Tab. 12: Conditioning of the internal submatrix for various orders.

***Note.*** On the point of view of the discretised quantity, this behaviour can be seen as a side effect of the definition of the discretisation space itself. Indeed, the improvement of the inner discretisation quality along the increment of the order $k$ simply matches the refinement of the projection space when adding projections on higher order monomials. It then forces the latests quantifiers to furnish much smaller values than the ones corresponding to the first quantifiers grasping the rough behaviours. ▲

This observation is particularly strong for elements whose shape is contained within the unit circle, as the powers of monomials encountered in the definition of the degrees of freedom will get close to zero. Therefore the degrees of freedom values will get similarly faster. Furthermore, note that due to the contribution of normal basis functions to internal basis functions, the bad conditioning of the internal submatrix is linked back with the normal sub - matrices in the inversion process.

In practice, one can consider safely elements until the fourth order. Note however that it is possible to reduce the conditioning by choosing another projection basis for $\mathcal{P}_k$ than the one based on Hermite polynomials that is considered here.

### 7.2.6 Conditioning of the transfer matrix with respect to the projection bases

We end the discussion on the conditioning by testing its value with respect to the chosen sets of polynomials, either defining the moments or generating the canonical basis functions. To this end, we plotted conditioning graphs representing the conditioning's values of the transfer matrix mapping the canonical basis functions to the ones defining the elements $Ia$, $Ib$, $IIa$ or $IIb$.

There, the types of employed polynomials are seen as parameters and are represented by a colour whose meaning is given in the *Table 13.* The corresponding conditioning value is depicted by a blue line.

More precisely, each of the four considered elements is using at least one definition of a moment based degrees of freedom of the form

$$\sigma_N : q \mapsto \int_{\partial K} p \cdot n\, r_i \,\mathrm{d}\gamma(x) \quad \text{for } \{r_i\}_i \text{ in } \mathcal{H}_k(K)$$

or

$$\sigma_I : q \mapsto \int_{K} p \cdot m_i \,\mathrm{d}x \quad \text{for } \{m_i\}_i \text{ base of } \mathcal{P}_k(K).$$

Its associated basis functions is built by adapting the canonical ones, taking by definition a form similar to

$$f_j := (x,\, y) \mapsto \begin{pmatrix} x \\ y \end{pmatrix} \left\{ \begin{aligned} &\Delta u = h_{k-1,\, k-1}(x, y) \\ &u|_{\partial K} = g_j(x, y)\mathbb{1}_{f_i} \end{aligned} \right\} + \begin{pmatrix} n_{ix} \\ n_{iy} \end{pmatrix} \left\{ \begin{aligned} &\Delta u = 0 \\ &u|_{\partial K} \equiv 2 \end{aligned} \right\}$$

and

$$f_{xlm} := (x,\, y) \mapsto \begin{pmatrix} 1 \\ 0 \end{pmatrix} \left\{ \begin{aligned} &\Delta u = h_{l,m}(x,\, y) \\ &u|_{\partial K} = 0 \end{aligned} \right\},\; 0 \le l,\, m \le k-1$$

where the functions $r$, $m$, $g$ and $h$ are polynomials (see the *Paragraph 6.1.2*).

In order to test the elements' construction depending on the used polynomial bases in the definition of the functions $r$, $m$, $g$ and $h$, we considered the possibilities presented in the *Table 13*.

| Colour/Function | $m$ (2D) and $r$ (1D) | $g$ (1D) | $h$ (2D) |
|---|---|---|---|
| Blue | Canonical (C + S) | Lagrangian | Chebyshev (S) |
| Red | Chebyshev (S) | Canonical (C+S) | Hermite (S) |
| Green | Hermite (S) | Canonical (C) | Legendre (S) |
| Yellow | Legendre (S) | — | Canonical (C+S) |
| Magenta | Laguerre (S) | — | Canonical (C) |
| Cyan | Canonical (C) | — | — |
| Taupe | Canonical (-) | — | — |

Tab. 13: Type of used polynomial basis and their corresponding colour indexation. C: centred with respect to the element's barycentre or the edge's midpoint. S: scaled with respect to the area of the element or the edge's length.

Then, for each tested parameter set we plotted the types of used basis functions $\{r_i\}_i$, $\{m_i\}_i$, $\{h_{l,m}\}_{l,m}$ and $\{g_j\}_j$ through a colour code, from the top to the bottom (*see the Figure 91*). For convenience, the parameter sets have been sorted so that the conditioning values are given in an ascending order.

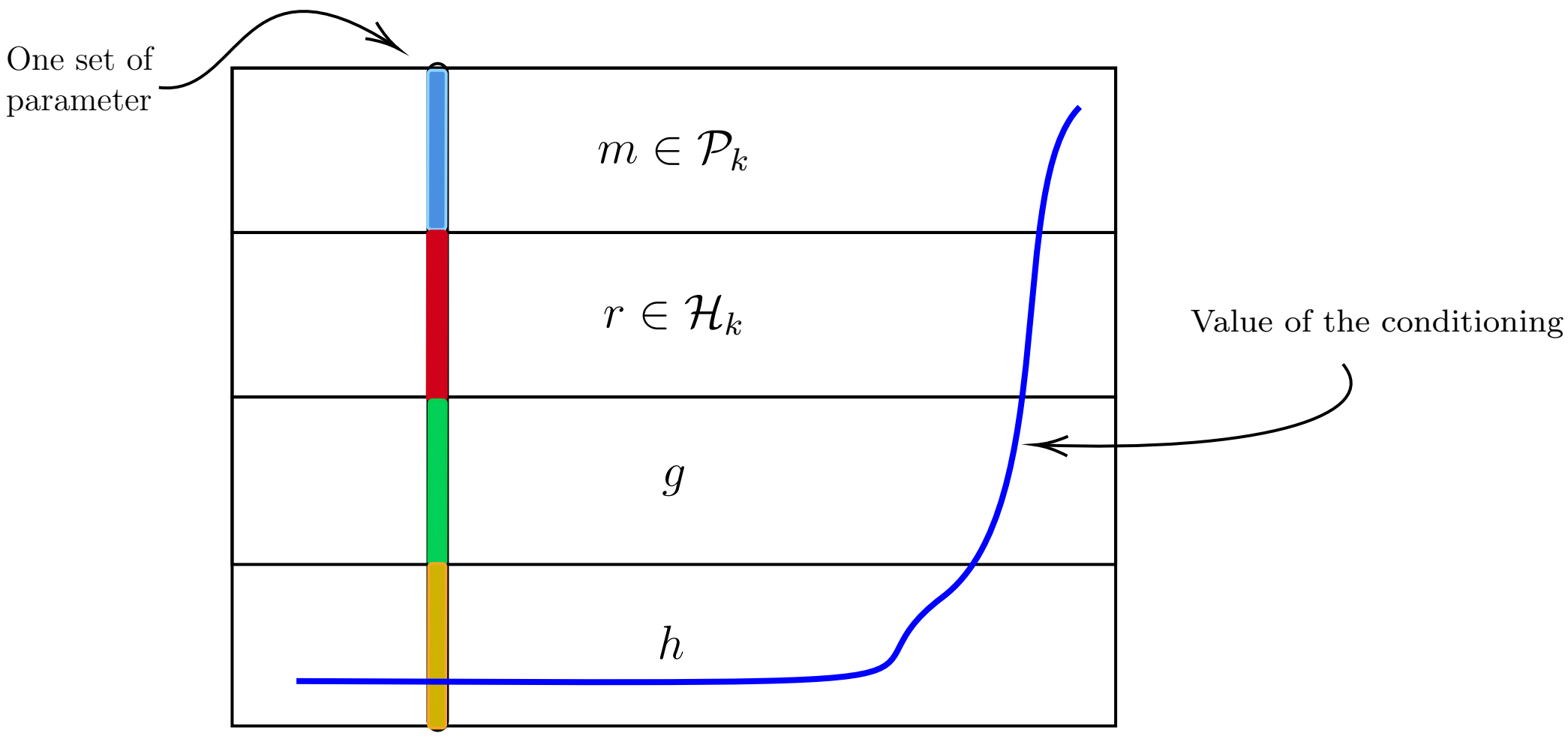

Fig. 91: Chart representing the layout of the conditioning graphs.

***Note.*** The $x$ - axis of the conditioning graphs only corresponds to the number of tested cases and has no particular meaning. ▲

***Remark.***
• In the case $k = 0$, the polynomials constructing the basis functions are either identically null or constant. Furthermore, there is no internal basis function nor internal moments. Thus, all the investigated types of basis function merge and there is no point of discussing the benefits from one configuration over another. Therefore, only the cases $k > 0$ will be considered.

• Unless mentioned, all the investigated polynomials lie within the unit circle. ▲

#### 7.2.6.1 A first test case

We start the discussion by analysing the configurations of parameters leading to the best conditionings in the particular case of the element *IIb* for $k = 2$, constructed over the non - convex hexagonal shape presented in the *Figure 87*.

As one can observe on the *Figures 92*, the parameter mostly influencing the conditioning drives the definition of the polynomials $\{m_i\}_i$, defining the projections within the element. Indeed, the parameter set driving the inner projectors are presented in the first parameter row, where large blocks of the same colour emerge. Observing those blocks while keeping in mind the conditioning values they map with, it can be observed that any conditioning's value contained within some contiguous range has been obtained from a parameter set using a specific parameter as an internal projector. Furthermore, this contiguous range of values is not wide, typically sharing the same order of magnitude.

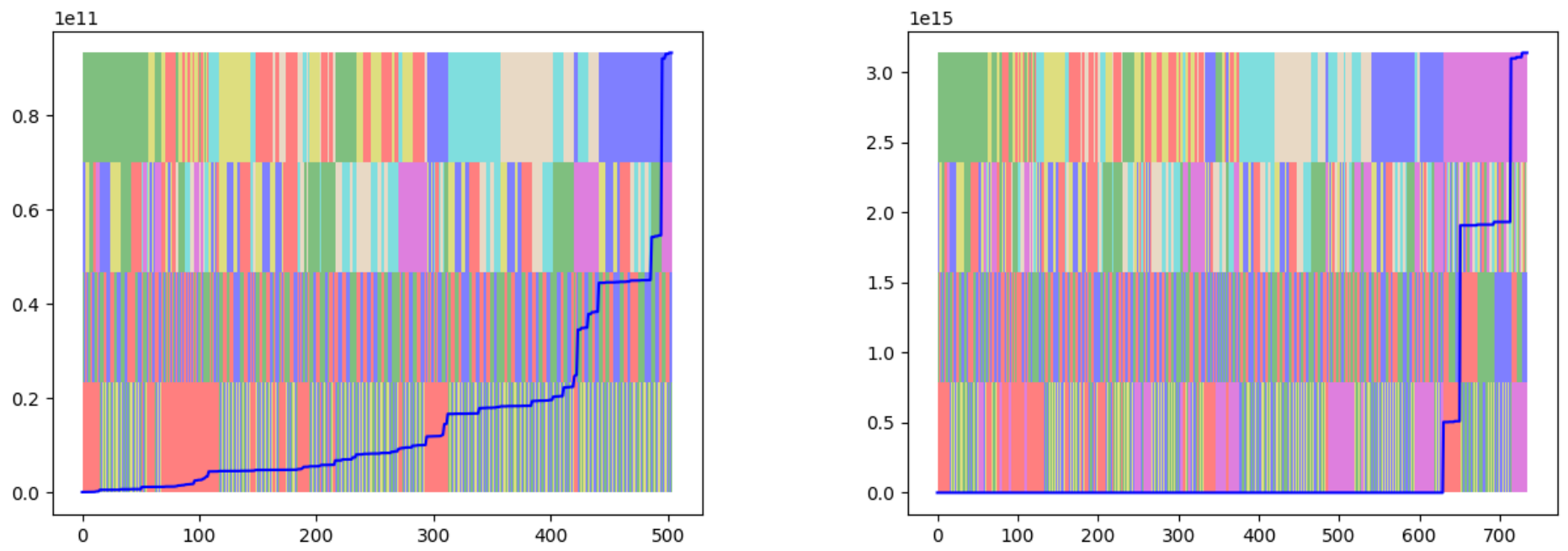

Fig. 92: Conditioning graph of the element $Ib$ for $k = 2$ on a non - convex hexagon. Left: without considering the Laguerre cases. Right: full considered test.

However, when the internal parameter changes, the variations of the conditioning values are steep. This can be witnessed at the bounds of each main colour block of the first row. Thus, the internal parameter leads to an almost piecewise constant conditioning line whose jumps follow its changes. Therefore, more than dominating the set of parameters, it also drives the order of magnitude of the conditioning value itself.

***Note.*** Qualitatively, the worse conditionings are achieved when the Laguerre polynomials are in use as internal projectors, followed by the use of the unscaled canonical basis functions. On the opposite side, the Hermite case represented in green shows a significant improvement of the conditioning's value. ▲

The secondary most impacting parameter turns out to be the choice of the polynomial space defining the second members $\{h_{l,m}\}_{l,m}$. Indeed, depicted in the fourth row, several blocks of colours can still be distinguished under each dominant internal projector block. In particular, we observe the same qualitative behaviour as above, where the magenta blocks (Laguerre case) corresponds to the highest conditioning values and the red blocks (Hermite case) to the lowest ones.

When the Laguerre case is excluded (*see the left part of the Figure 92*), one can also witness the secondary impact of the boundary projectors. Indeed, qualitatively they appear to be of the same importance than the inner constructors by defining clearly distinct groups of parameters. However, quantitatively, for any boundary projector the conditioning values are contained within the same order of magnitude. Therefore, any choice of boundary projector would be admissible. Lastly, it can be observed that boundary constructors do not show any significant grouping, meaning that their choice is the least impacting on the conditioning value.

However, those boundary constructors severely impact the conditioning values when Laguerre internal projectors are used (*see the right part of the Figure 92*). Indeed, more than defining clear groups, the boundary constructors are secondary impacting parameters, dominating the impact of the internal constructors. One exception has to be made when Laguerre internal constructors are used, where they are shifted to a tertiary important level. Surprisingly, the blocks of boundary constructors read from the worse to the better case; Lagrangian, canonical, and Centered scaled canonical polynomials.

In conclusion, provided that the canonical and Laguerre cases are not used, the normal projectors and normal constructors play only a very little role in the conditioning values. And indeed, it connects with the observations made in the previous paragraph where we saw that the worsening of the conditioning was mostly due to the internal submatrix.

Thus, at a first glance, it is recommended to use a combination of Hermite polynomials for both internal projectors and constructors. The choice of the normal projectors and constructors is not relevant at this point and will be discussed later on, with respect to further investigations.

#### 7.2.6.2 Sensitivity with respect to the polygonal shape

Let us now test the sensitivity of the best parameter set towards the element's shape. To illustrate our conclusions, we select the particular case of the element $Ib$ for $k = 2$, applied on a triangle, quadrilateral and on the non - convex hexagon presented above.

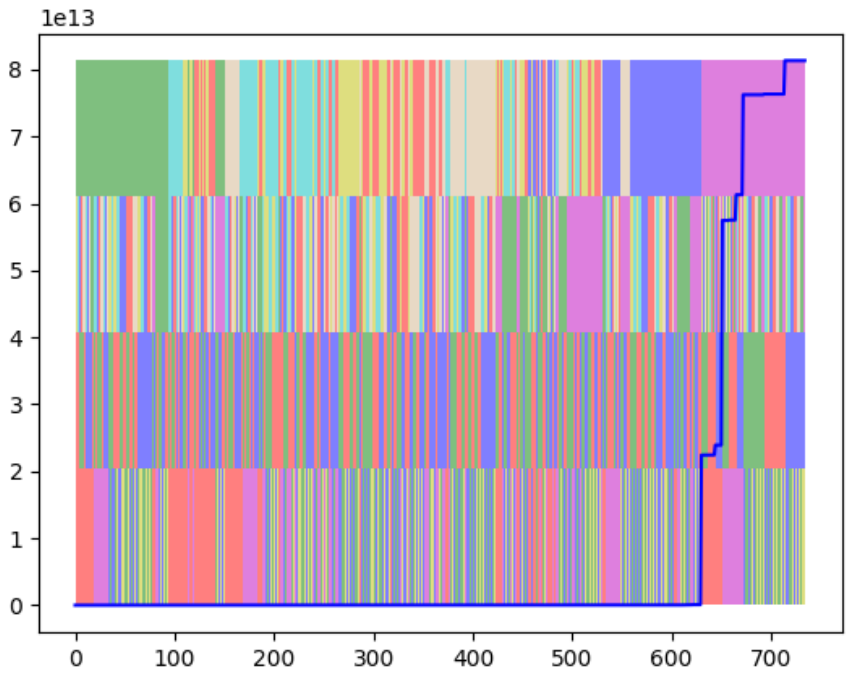

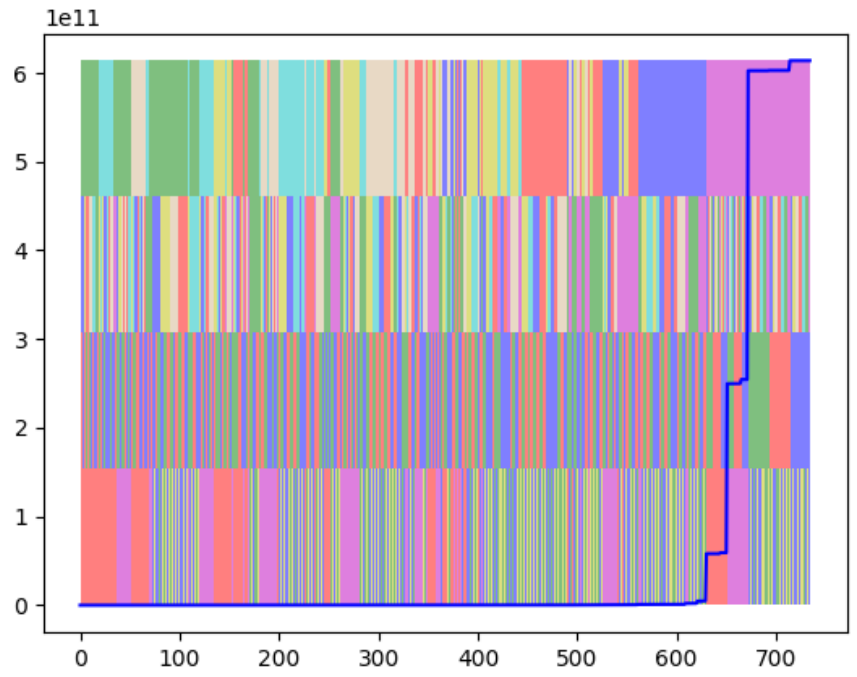

First of all, we can observe that for any chosen parameter set the order of magnitude of the conditionings are greater for the non - convex hexagon than for the quadrilateral, which are in their turn greater than for the triangular shape. This matches the observations that were done in the previous subsection.

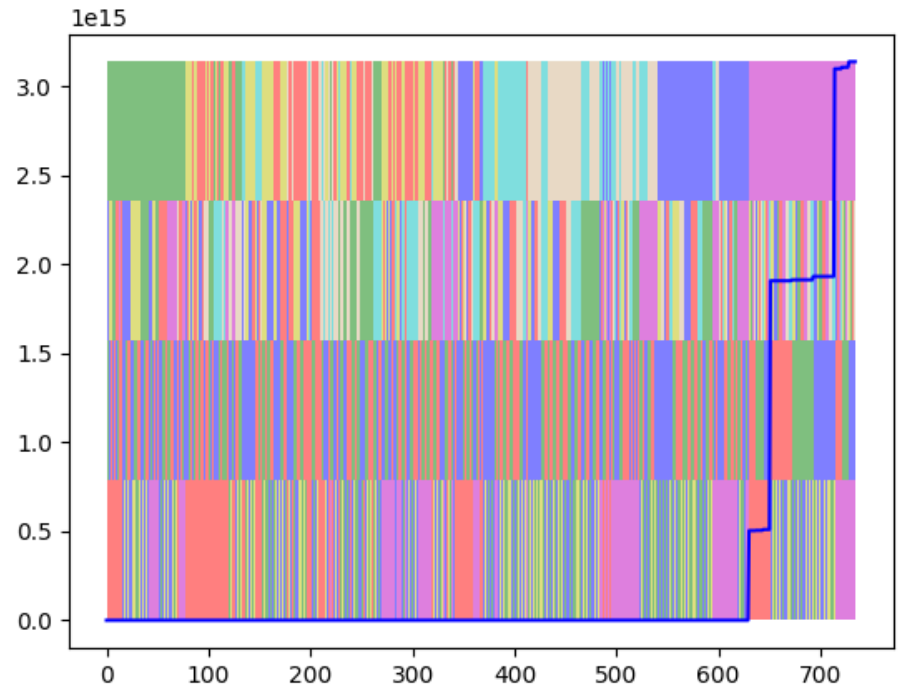

Fig. 93: Conditioning graphs for the element *Ib* at $k = 2$. Top left: Triangle case. Top right: Quadrangular case.

Then, it comes out that the layout of the best parameter sets are very similar from one case to another. Especially, the three extreme cases of the Hermite, canonical and Laguerre polynomials are preserved regardless the polygonal shape.

The main change consists in the permutations of the inner projectors (depicted in the first row) between the red (Chebyshev), cyan (centred canonical), yellow (Legendre) and taupe (canonical) cases. In particular, when taking the triangle case as a reference, the Chebyshev and canonical cases flip their precedence for the quadrilateral, while a combined flip yellow/red to cyan occur when going to the non - convex case.

Those changes emphasise the impact of the inner properties of the element. Indeed, the canonical basis functions are heavily dependent on the constructors' types, especially when the polygon is squeezed. Thus, the best projectors depends on the squeezeness of the polygon, as the projection of shape - sensitive Poisson's solution onto the functions $\{m_i\}_i$ impacts the characterisation of any function living in $\mathbb{H}_k(K)|_{\mathring{K}}$. Furthermore, not knowing the squeezeness of the polygon prevents any *a - priori* control over the parameter sets.

In particular, one can compare the quadrilateral case, which leads to very regular Poisson's solutions, with the triangular case. There, it can be observed a clarification of the best parameter sets. Indeed, a condensation of the red (Hermite) cases for both internal projectors and generators can be noticed. The Hermite projectors are not interspersed anymore by Legendre cases, and the Hermite constructors become of primary importance when considering only the low conditioning cases (rather than being of secondary importance after the boundary projectors in the triangle case). This is a simple consequence of the regularity of the Poissons' problems, where the definition of the second member functions can be transferred to the canonical basis function itself without being altered by the squeezeness of the element. The prevalence of one parameter over another thus becomes clearer.

A similar phenomena of precedence reversal is witnessed in case of the non - convex hexagon. There, the squeezeness of the polygon amplifies the similarities between the canonical basis functions generated from the Laguerre cases. In particular for the sets leading to the worse conditionings, the Laguerre inner constructors prevail on the boundary constructors, which is not the case *e.g.* for the triangle and quadrilateral cases.

However, the conditioning values corresponding to the cases impacted by a permutation of the best parameter sets remains in the same order of magnitude. Furthermore, the boundary cases still have only a few and negligible impact for any polygonal shape. Therefore, we can conclude that the behaviour is globally similar, preserving the decreasing pattern of the conditioning and three extreme cases detailed in the previous paragraph. In all the next, the non - convex hexagon will be taken as a model case.

#### 7.2.6.3 Sensitivity with respect to the discretisation order

Let us now discuss the impact of the discretisation order on the best parameter sets. To this sake, we represent in the *Figure 94* the conditioning graphs of the element *Ib* for various orders, built on the non - convex hexagon presented in the *Figure 87*.

Let us first point out that in the case $k = 1$, the inner constructors are constants. Therefore, the last row of the top left graph has no meaning. However, we can already witness the primary importance of the inner projectors, already defining the same parameter precedence as witnessed in the previous paragraphs in the case $k = 2$. Furthermore, as the inner constructor types are of no importance here, it is possible to observe a criterion over the boundary constructor and projector sets.

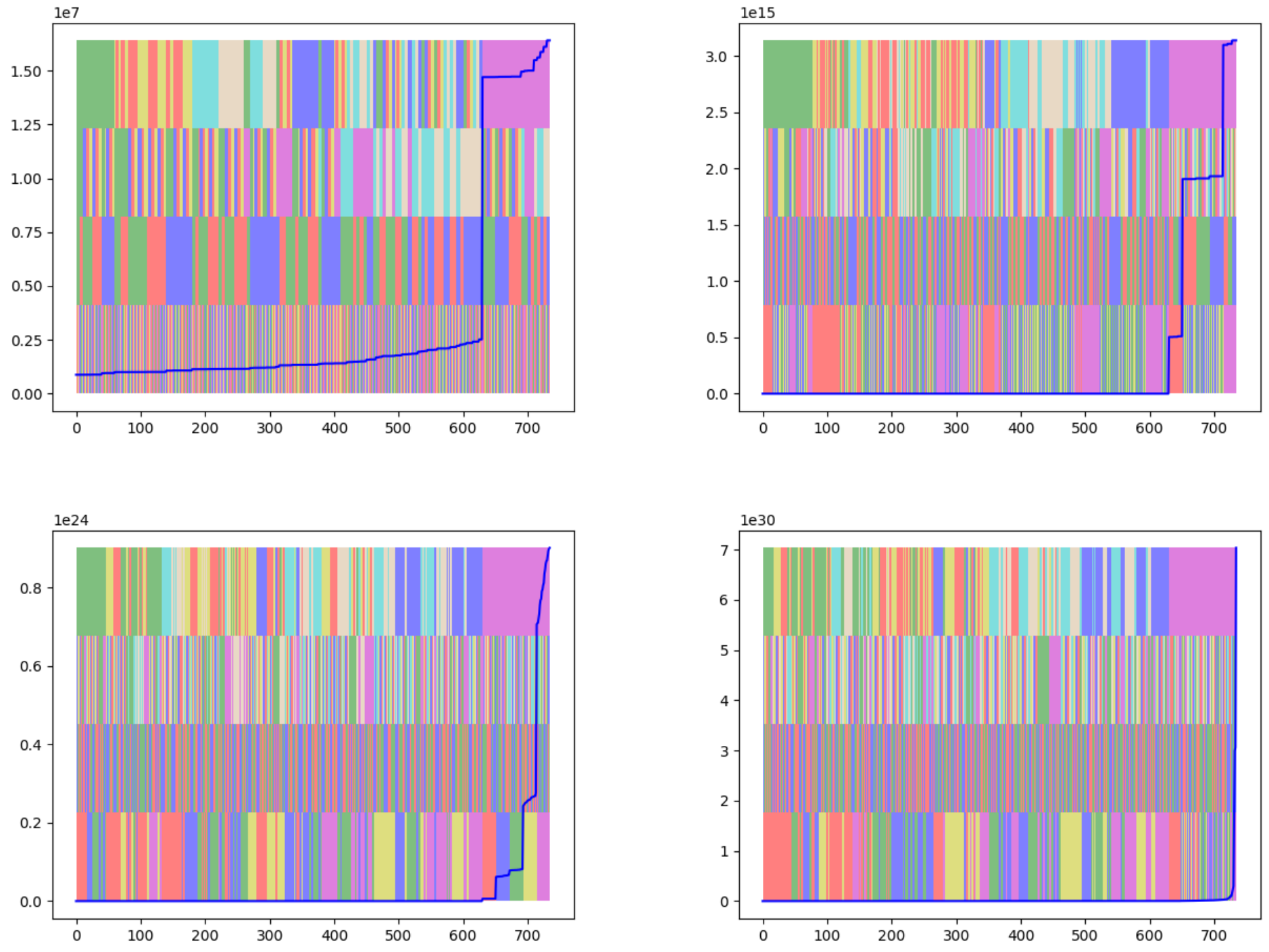

Fig. 94: Conditioning graph for the element $Ib$ constructed on a non - convex hexagon. Top left: $k = 1$. Top right: $k = 2$. Bottom left: $k = 3$. Bottom right: $k = 4$.

Indeed, the choice of the boundary projector appears to be mainly of secondary importance, where a pattern driving the improvement of the conditioning number can be clearly established. Indeed, a main repetition sequence blue (Lagrangian) – red (centred scaled canonical) – green (canonical), from the worst to the best, emerges. Thereafter, the boundary constructors play mainly the role of a tertiary important parameter, especially in the low conditioning area. Indeed, under each block of same boundary constructor a pattern yellow (Legendre) – red (Chebyshev) – violet(Laguerre) – green (Hermite), from the worst to the best, can be witnessed.

However, as discussed before, the changes in those conditioning values are contained within a negligible range. Furthermore, note that all the cases that are not represented in this pattern are concentrated the right area of the graph where no clear precedence rule between the boundary constructors and projectors can be established. It corresponds to the case happening after the elimination of the inner Laguerre projectors, where both boundary projectors and constructors appears to be conjointly operating a classification of primary importance.

For higher orders, the choice over the inner constructors comes into play, and at a first glance, the sorting of the parameters are similar for any order. However, it can be seen that higher the order is, the most important the inner characterisation becomes. Indeed, the small clusters driven by a boundary projectors in the case $k = 1$ are vanishing in favour of the inner generators in the case $k = 2$. From the order 2 to the order 3, the precedence of the inner projector over the inner constructor is also challenged, especially for the cases leading to low conditionings. Here again, it is due to the definition of the basis functions whose construction through a Poisson's problem amplifies the impact of the second member's definition.

#### 7.2.6.4 Sensitivity on the element

So far, all the discussions were conducted by taking the example of a particular element, either $Ib$ or $IIb$. Let us now test the sensitivity of the best parameter sets with respect to the considered element. As for the discussion on the order, we separate the case $k = 1$ from the higher orders, as there the last parameter row is not significant.

In the case $k = 1$ represented in the *Figure 95*, the conditionings of each element lead to two main ordering patterns of the parameter sets, as well as two main layouts of the conditioning decay.

Indeed, one can observe that for both configurations $I$ and $II$, the shape of the conditioning's decay is impacted by the type of the element, $a$ or $b$. After excluding the internal Laguerre projector case, if the type $a$ generates a plateau corresponding to the boundary Laguerre projector case, the types $b$ never stagnates and enjoy a smooth decay. This is due to the use of point values rather than further moments in the definition of the lowest order characterisation. As the type $b$ does not use lowest order moments, the impact of the projectors is less accentuated and the variations between two parameter sets are accordingly decreased.

This observation especially holds for the configuration $I$ where this impact can also be slightly witnessed in the ordering of the parameter sets. Indeed, there the lowest order elements are only tuned coordinate wise, which makes the element's definition more sensitive to the choice of using specific point values or moments. And in particular, choosing points value entirely breaks the distinction that can be made between boundary projectors at the lowest orders, as the element would be less sensitive to their definition. This phenomenon can be seen in the fracturing of the colour block corresponding to the Laguerre boundary projection case, as well its permutation with the cou-ple of taupe and blue blocks. There, all those parameters tend to compete as the corresponding conditioning values are concentrated within a small range.

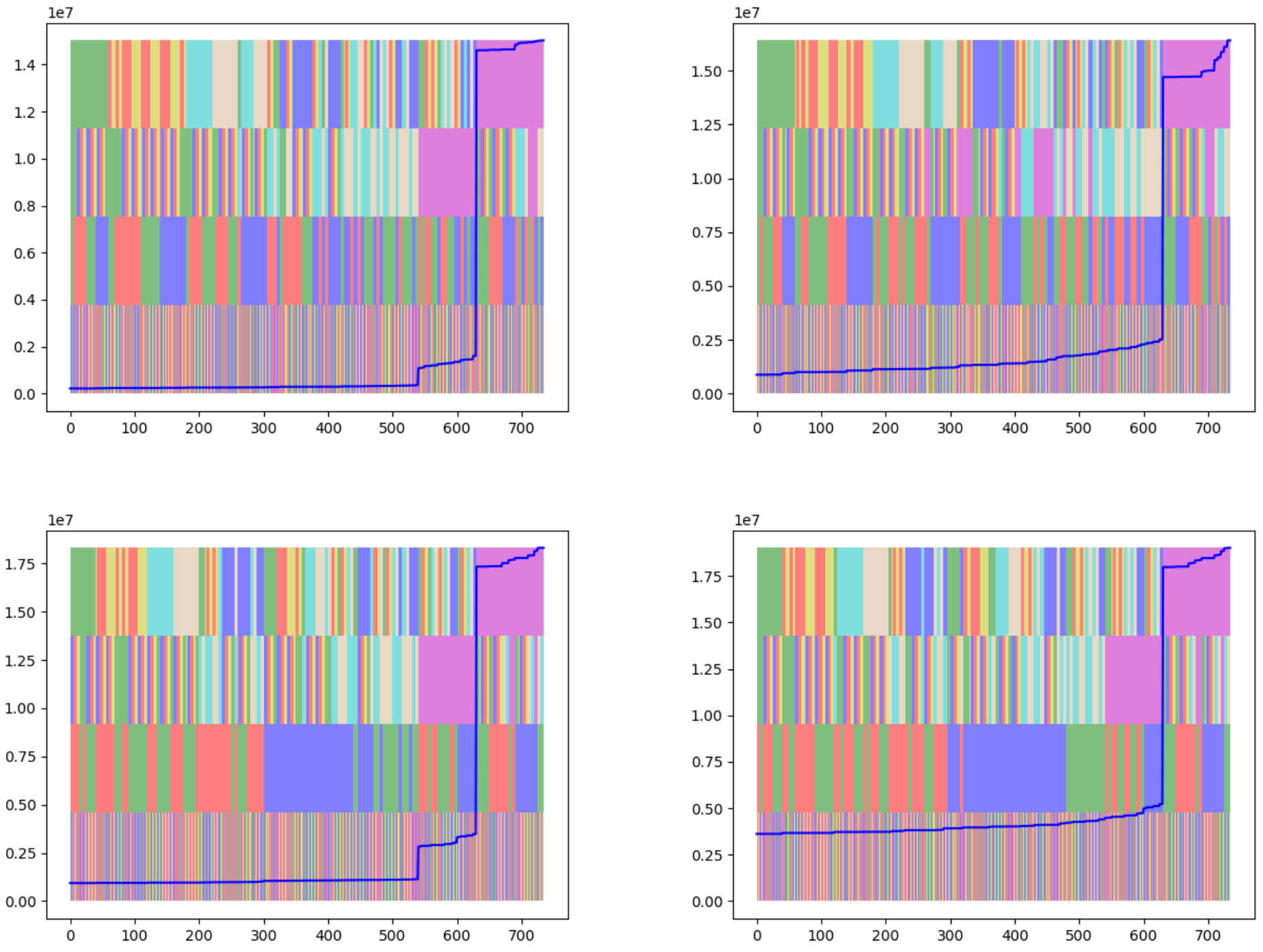

Fig. 95: Conditioning graphs for elements built on a non - convex hexagon for $k = 1$. Top left: $Ia$. Top right: $Ib$. Bottom left: $IIa$. Bottom right: $IIb$.

As the case $II$ is always considering at least first order projections coordinate wise, the ordering of the best parameters are changed not significantly by using point values. However, note that the decay of the conditioning is still smoothed.

For a fixed element type, the choice of the configuration alters the precedence of the boundary constructors. Indeed, the difference between the configurations $I$ and $II$ lies in the grouping of the centred scaled canonical polynomials, becoming of first importance in the precedence of the parameter sets. This change is compliant with the elements' definitions. The configuration $II$ is characterising the elements of $\mathbb{H}_k(K)$ with a more global approach than the configuration $I$. Thus, the boundary constructors are more globally solicited, especially at the lowest order level, allowing a clearer ordering of the bests boundary constructors.

***Remark.*** Despite the above comments, please note that in any case the extreme cases remain unchanged. ▲

Let us now consider the higher order cases, *i.e.* from $k = 2$ on. We pictured the case $k = 2$ in the *Figure 96* for reference.

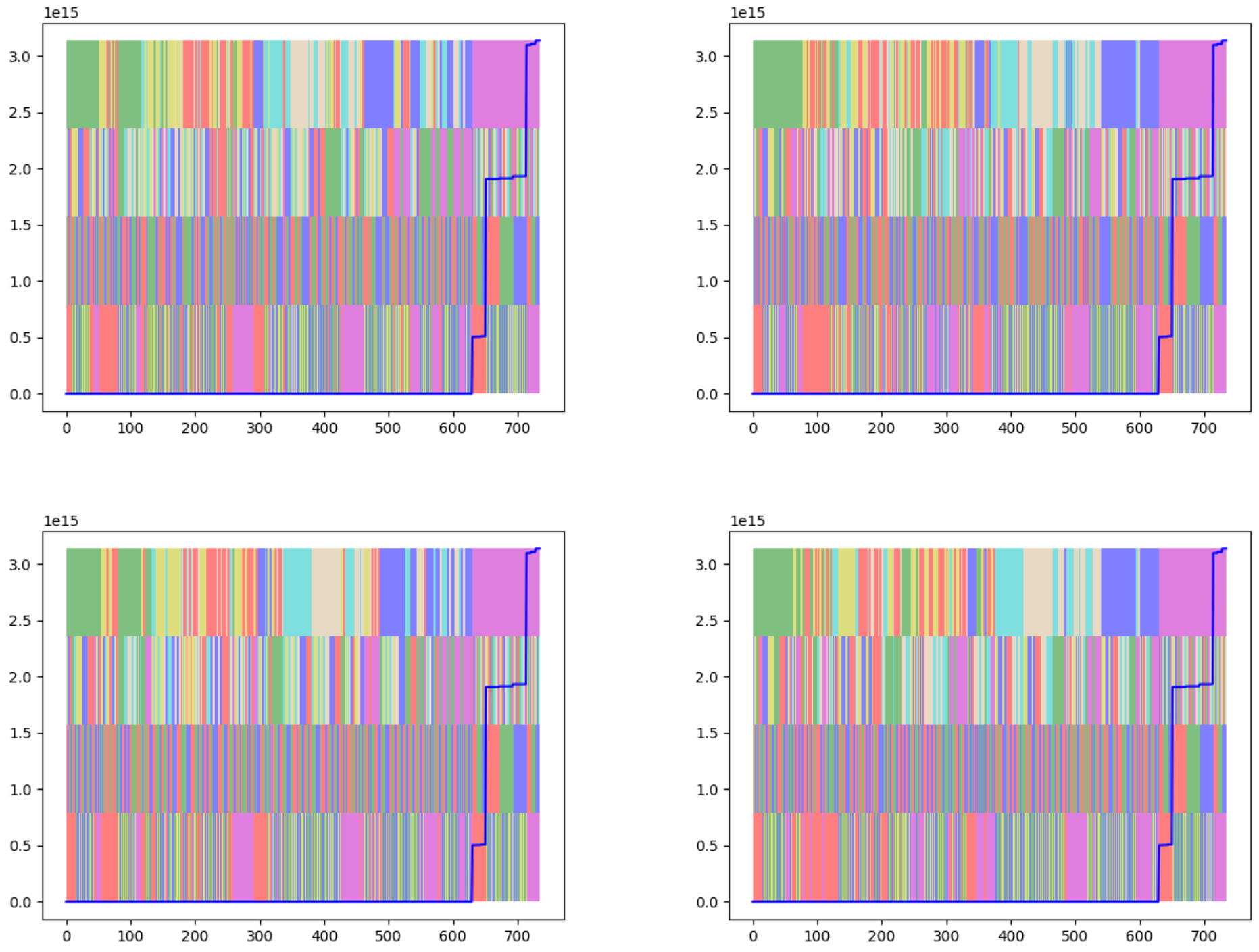

Fig. 96: Conditioning graphs for elements built on a non - convex hexagon for $k = 2$. Top left: $Ia$. Top right: $Ib$. Bottom left: $IIa$. Bottom right: $IIb$.

There, as detailed in the previous paragraph, the impact of the boundary constructor is disregarded for the benefit of the inner constructor. Thus, as the inner characterisation of functions living in $\mathbb{H}_k(K)|_{\mathring{K}}$ is identical for any considered element, no significant change occurs when comparing one element's definition with an other.

Indeed, no major difference emerges when comparing the ordering of parameter sets corresponding to two different elements configurations $I$ and $II$ having the same type $a$ or $b$. Furthermore, the only difference between the types $a$ and $b$ impacts the grouping of the best projectors. Indeed, the use of point values do not involve lowest order moment on the boundary, allowing the inner projectors to take over the precedence. However, as in the case $k = 1$, this grouping has no major influence on the conditioning's value.

To conclude, let us recall that higher the order is, the most impacting the internal projectors and constructors are. Thus, as the internal construction of the four considered elements is the same there will be less and less difference between their conditioning graphs. In particular, the conditioning decay will be of the same shape for any of them, as it can already be seen on the *Figure 96*.

Therefore, we can state that all the discussion detailed on particular elements in the above paragraphs still hold for any other considered element in the cases $k > 1$.

#### 7.2.6.5 Sensitivity with respect to the size of element

Let us now end this discussion by testing the sensitivity of the conditioning graphs with respect to the size of the element. We present in the *Figure 97* the results for the element $IIb$ constructed on the second and third triangular shapes of the *Figures 88*, for the orders $k = 1$ and $k = 2$.

For both orders, one can observe that increasing the size of the triangle makes little changes in the ordering of the boundary constructors and allow a clearer separation of the inner projector blocks. Indeed, the values of the Poisson's solutions are amplified with the size of the element, allowing a better comparison of the projectors. However, as always in this section those small changes do not severely impact the corresponding conditioning values. Therefore, for polygons contained within the unit circle their scaling is not severely impacting the precedence of the parameters sets.

However, some noticeable changes can be observed when the size increase forces the polygon to be mostly out of the unit circle. We plotted the conditioning graphs for the third triangular shape presented in the *Figure 88* on the *Figure 98*.

In the first order case, one can notice a severe change in the precedence order of the parameters, where the outcome does not allow us to derive tangible conclusions elsewhere from the limit cases. Furthermore, considering the boundary projectors, the blocks corresponding to the canonical and Laguerre cases appears to be the best choice, while defining the worst cases for the inner projectors. Thus, all the contributions become severely interconnected and it is hard to draw a general rule.

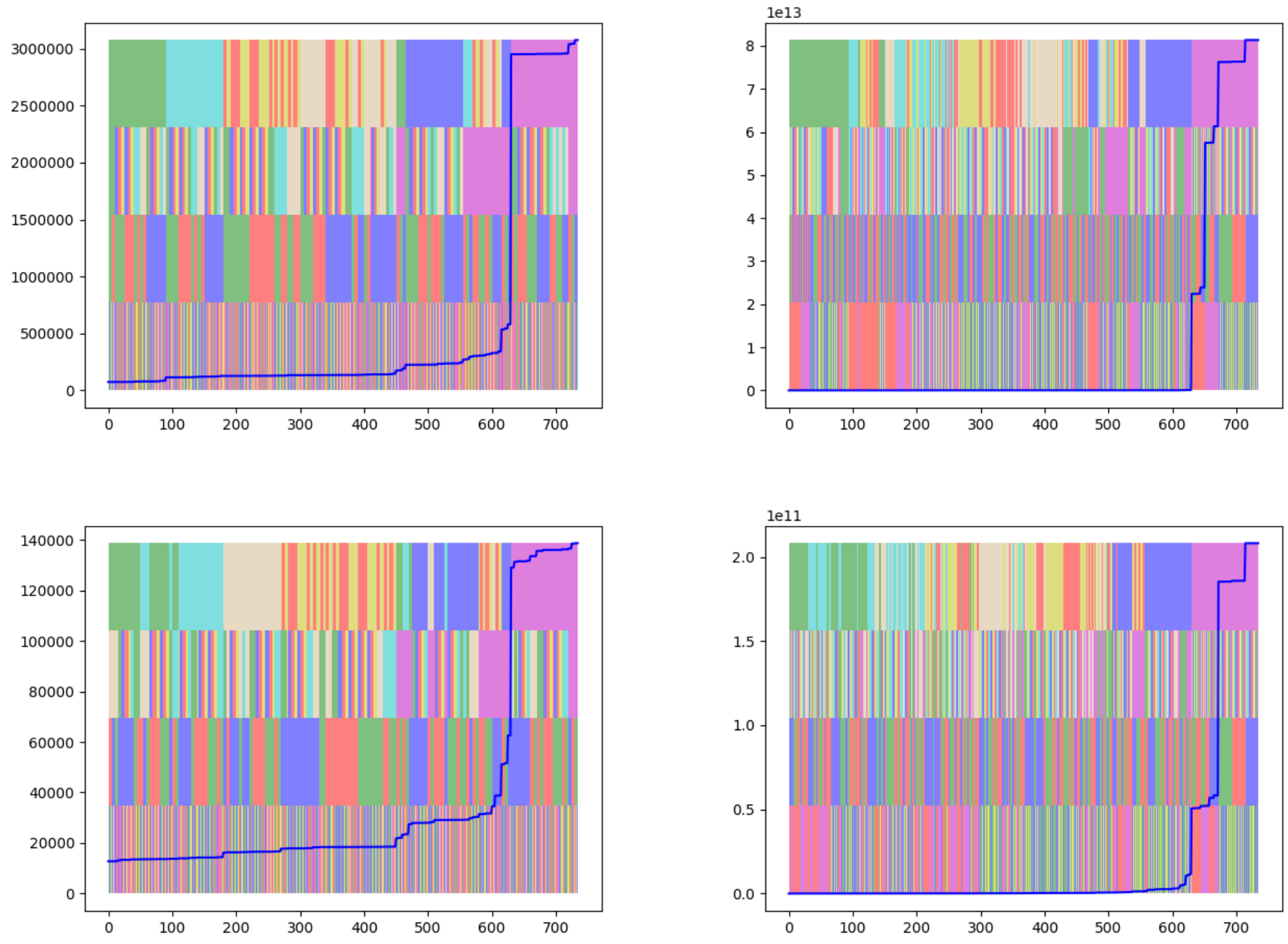

Fig. 97: Conditioning graphs for the element $IIb$ constructed on a triangles. Left column: $k = 1$. Right column: $k = 2$. Top: Small triangle. Bottom: Bigger triangle.

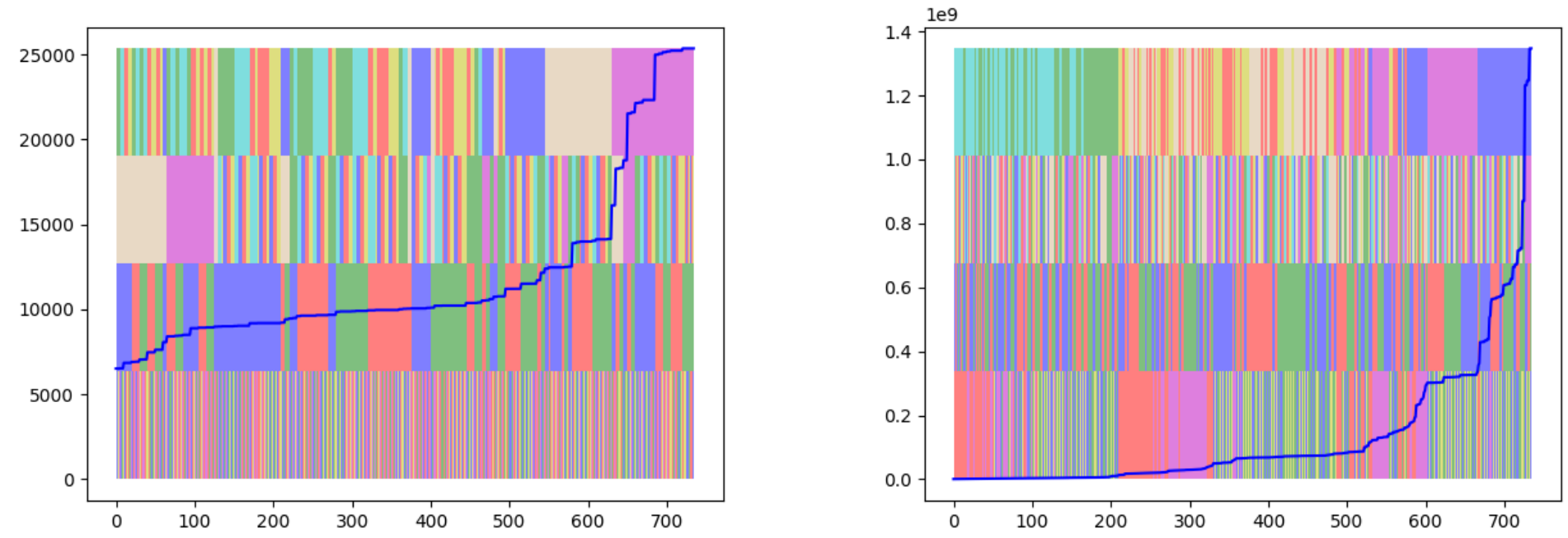

Fig. 98: Conditioning graphs for the element $IIb$ constructed on a triangle lying outside of the unit circle. Left column: $k = 1$. Right column: $k = 2$.

If for $k = 2$ the results are more compliant with the conclusions drawn for elements living within the unit circle, one has to notice the permutation between the colour blocks corresponding to the boundary Laguerre and centered scaled projectors. A similar flip can be observed as well between the non - scaled canonical and Hermite polynomials. Therefore, the conclusions of above may change for elements outside the outer circle.

#### 7.2.6.6 Conclusion

As a practical conclusion of the conditioning tests, we can state that overall, the Hermite internal projectors and constructors are leading to the most reliable elements. The selection of polynomials generating of characterising the boundary behaviour are only of secondary importance, except in the first order case.

Let us now observe the shapes of the tuned basis functions. In all the following results, the elements have been constructed by using Hermite polynomials as boundary projectors and Lagrange polynomials as boundary constructors, unless explicitly mentioned.

## 7.3 Study of the basis functions of the three elements

Let us now explore more in details the properties of all the discussed elements through their basis functions. In all this paragraph, we consider the non - convex hexagon given in the *Figure 90.* As the internal basis functions are tuned in the same manner for all of the elements, they will be detailed only once in the context of the element $Ia$.

### 7.3.1 Layout of the basis functions

We first represent the basis functions corresponding to the four elements under consideration when the degrees of freedom are built from the Hermite polynomial space.

#### 7.3.1.1 Element Ia

The element $Ia$ is built from the intuitive set-up of normal degrees of freedom whose definitions are only making use of moments. Let us investigate how the definition of those degrees of freedom tune the canonical basis functions.

**Normal basis functions** We start by analysing the basis functions that are generated from the canonical normal basis functions. As an example, we plot in the *Figure 99* the internal behaviour of one of the basis functions enjoying the highest variations, associated to the second edge of the polygon. Its boundary behaviour is plotted in the *Figure 100*.

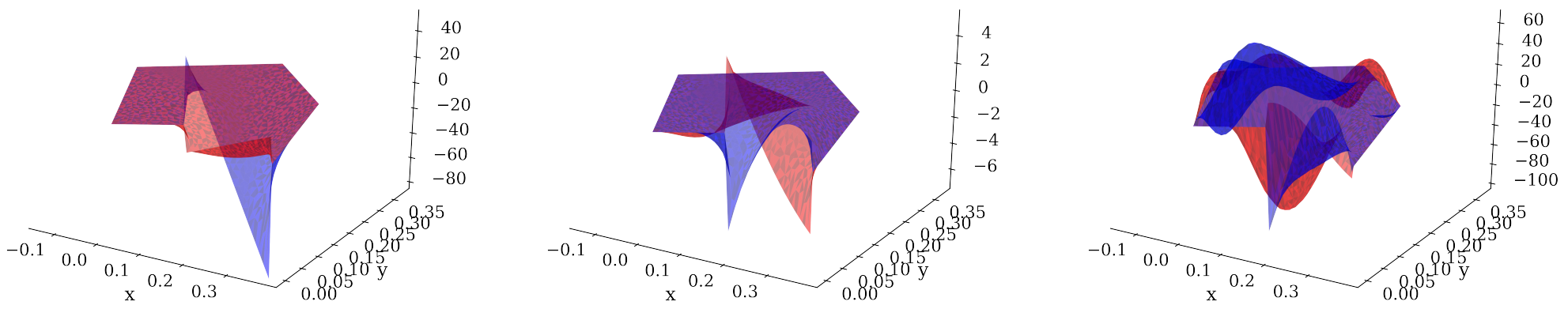

Fig. 99: Internal components of a basis function representative of the order of the space. Blue: $x$ component. Red: $y$ component. From left to right: $k = 0$, $k = 1$, $k = 2$.

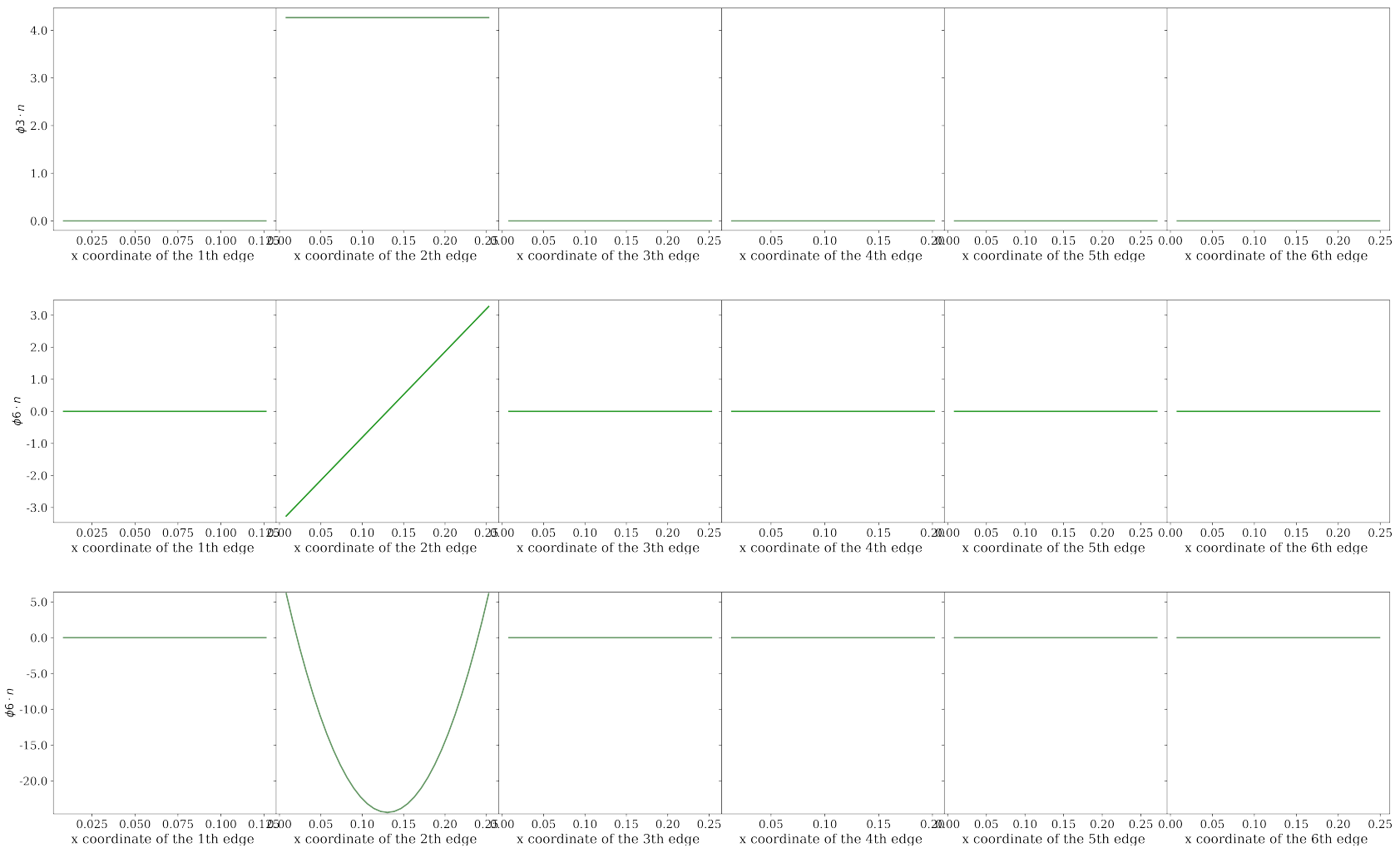

Fig. 100: Representative normal component of normal tuned basis functions that have a non - identically vanishing normal component. Eclated view on every edges. Top to bottom: $k = 0$, $k = 1$, $k = 2$.

It can be observed that the two components of the normal basis function still enjoy some smoothness within the element for any order $k$. And indeed, the new basis functions is only a linear combinations of the smooth canonical basis functions.

In addition, we observe that the support of the normal component of the normal basis function is still limited to a single edge, which makes easier the representation of the discretised quantity across the boundaries. The observed variations on its support are then still compliant with the order of the element.

However, we do not enjoy anymore the local point - wise Lagrangian property. Indeed, as pictured in the *Figure 101* for the second edge, the tuning method does not respect the assigned pointwise values nor the Lagrangian layout that the $k+1$ first functions of each edge were enjoying.

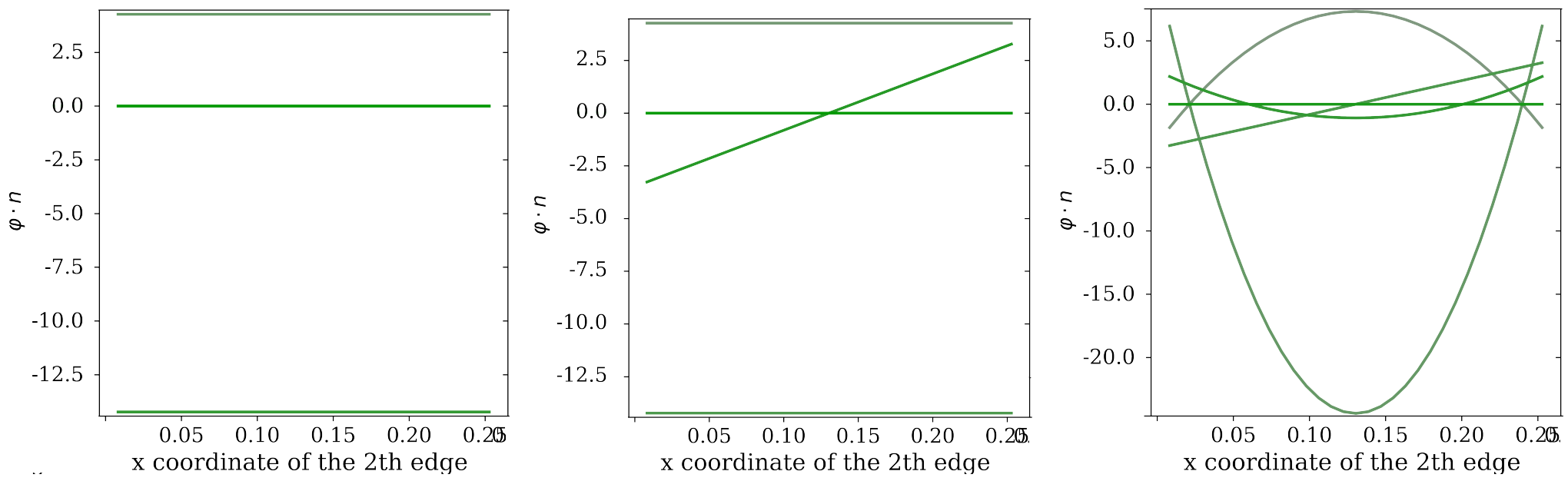

Fig. 101: Plot of all the functions emerging from the second edge on the second edge of the element. From left to right: $k=0$, $k=1$, $k=2$

As a consequence, in the lowest order case the normal component of the normal functions do not scale to one. We can notice that their normal component is not merging either. More importantly, one function seems to be missing. In fact, during the tuning process one normal basis function is degenerating into an internal basis function, without totally vanishing. In particular, its components are both not vanishing on the boundary, making it slightly different from the originally designed internal basis functions. It in fact corresponds to the global normal behaviour being tuned directly from the two coordinate wise ones. We represent the behaviour of the three normal basis functions of the first edge on the *Figures 102* and *103* in the case $k=0$.

There, one can observe that though its normal component is vanishing, the its components are not vanishing coordinate - wise on the boundary. Furthermore, within the element the variations are smooth, as being a linear combination of functions belonging to either to the spaces $A_k$ or x$B_k$. Thus, its regularity allows the degenerating function to be reuptaken by set of internal basis functions. Note that for higher order elements, always only one normal basis function will be vanishing.

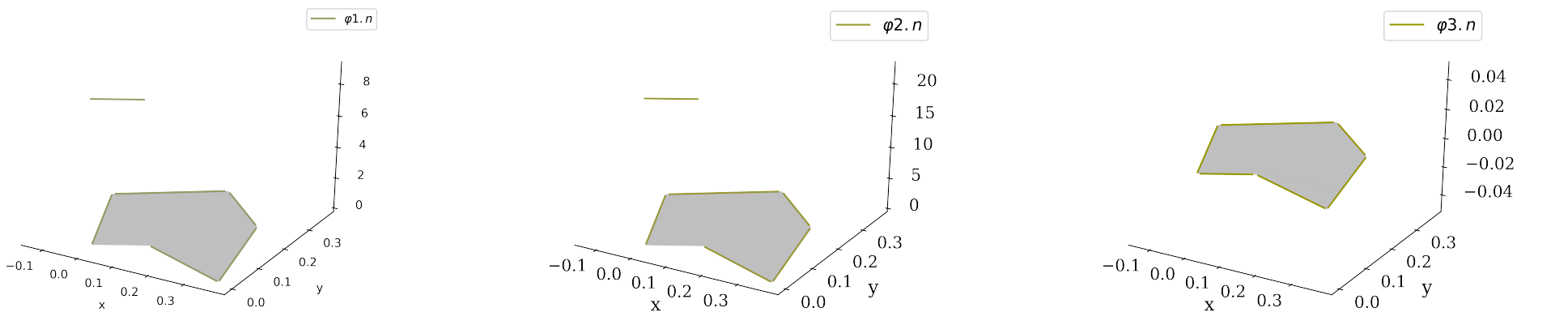

Fig. 102: Plot of the normal component of the basis functions emerging from the former canonical normal basis functions in the case $k = 1$. From left to right: first normal basis function, second normal basis function, third degenerating basis function generated from a canonical normal basis function.

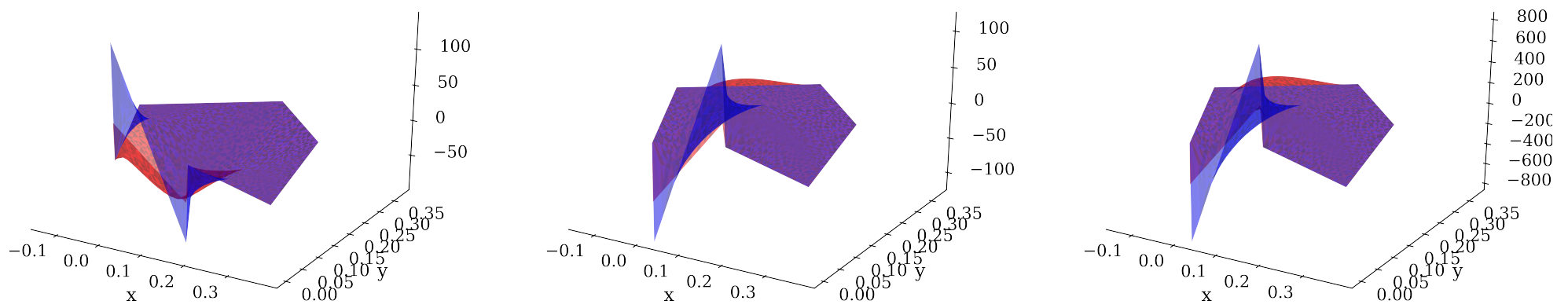

Fig. 103: Plot of the internal behaviour of the basis functions emerging from the former canonical normal basis functions in the case $k = 1$. From left to right: first normal basis function, second normal basis function, third degenerating basis function generated from a canonical normal basis function. None of the coordinate - wise components are vanishing identically on all the boundaries.

**Internal basis functions** Lastly, we can confirm that the properties of the internal functions are well preserved. Let us plot the last internal basis function in the *Figure 104*.

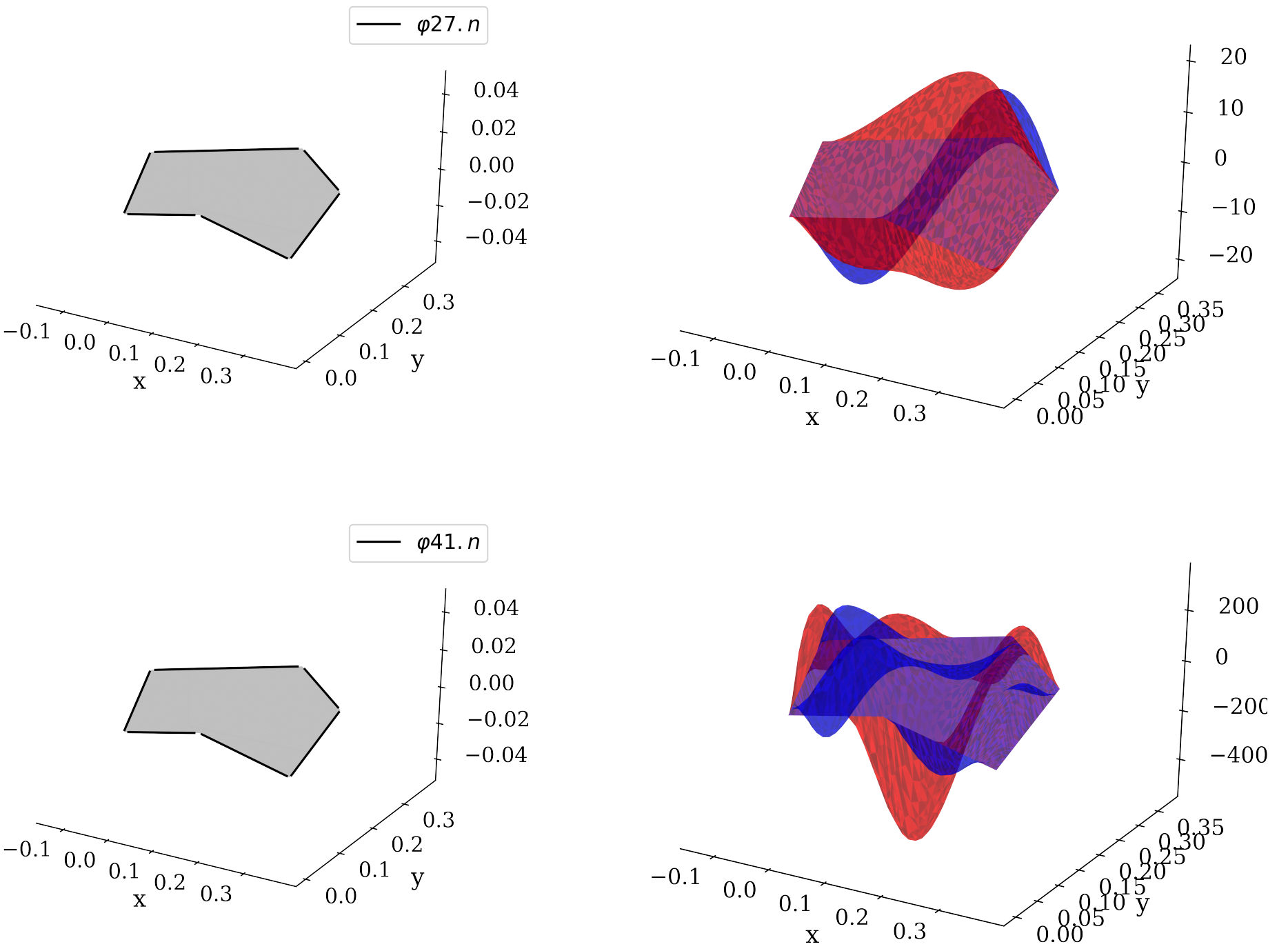

Fig. 104: Plot of the last internal basis function. Top: $k = 1$. Bottom: $k = 2$. Left: Normal component of the basis functions. Right: internal behaviour of the basis function.

There, we can draw the same conclusions as for the normal basis functions when considering the smoothness and the amplitude within the element. We also note that the vanishing property of the normal component on the boundary is well preserved. Therefore, the split between internal and normal basis functions is preserved, though the internal subgroup is enlarged by $n$ basis functions.

#### 7.3.1.2 Element Ib

The element $Ib$ is built in the same spirit than the element $Ia$, but using point - wise values of the normal components for the lowest order projections. Let us see how this little change in definition of the degrees of freedom impacts the tuned canonical basis functions.

When considering the normal basis functions, the same qualitative conclusions as for the previous element can be derived for both the variations of the normal components along the edges and the smoothness of the basis functions' components within the element. The only change is quantitative, as the amplitude of the functions tend to be higher. For the sake of completeness, we plot in the *Figure 105* the normal component of one basis function representative of the space order, emerging from the second edge. Its behaviour within the element is pictured in the *Figure 106*.

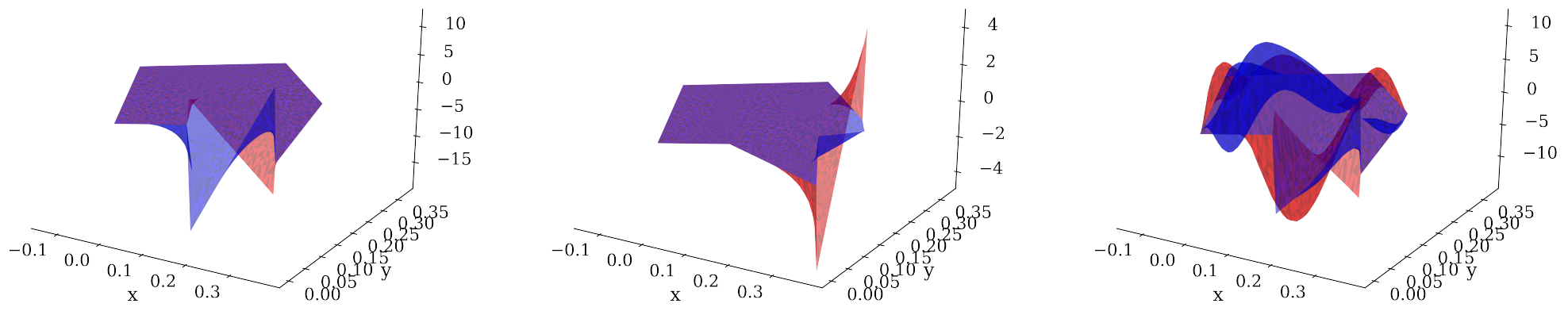

Fig. 105: Internal components of a basis function representative of the order of the space. Blue: $x$ component. Red: $y$ component. From left to right: $k = 0$, $k = 1$, $k = 2$.

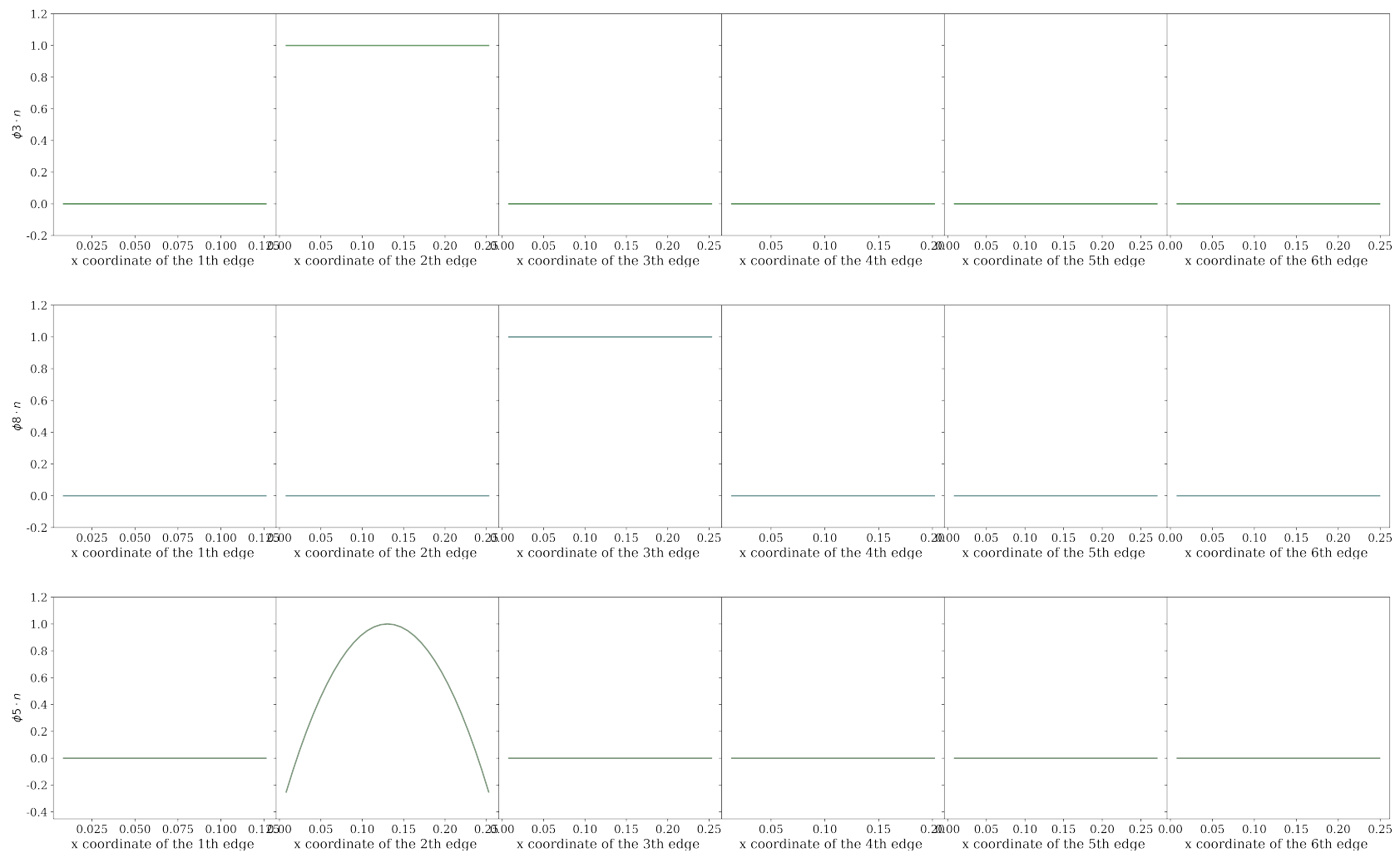

Fig. 106: Representative normal component of normal tuned basis functions that have a non - identically vanishing normal component. Eclated view on every edges. Top to bottom: $k = 0$, $k = 1$, $k = 2$.

However, slight changes can be noticed when paying a closer attention to the layout of the normal basis function, represented in the *Figure 107*.

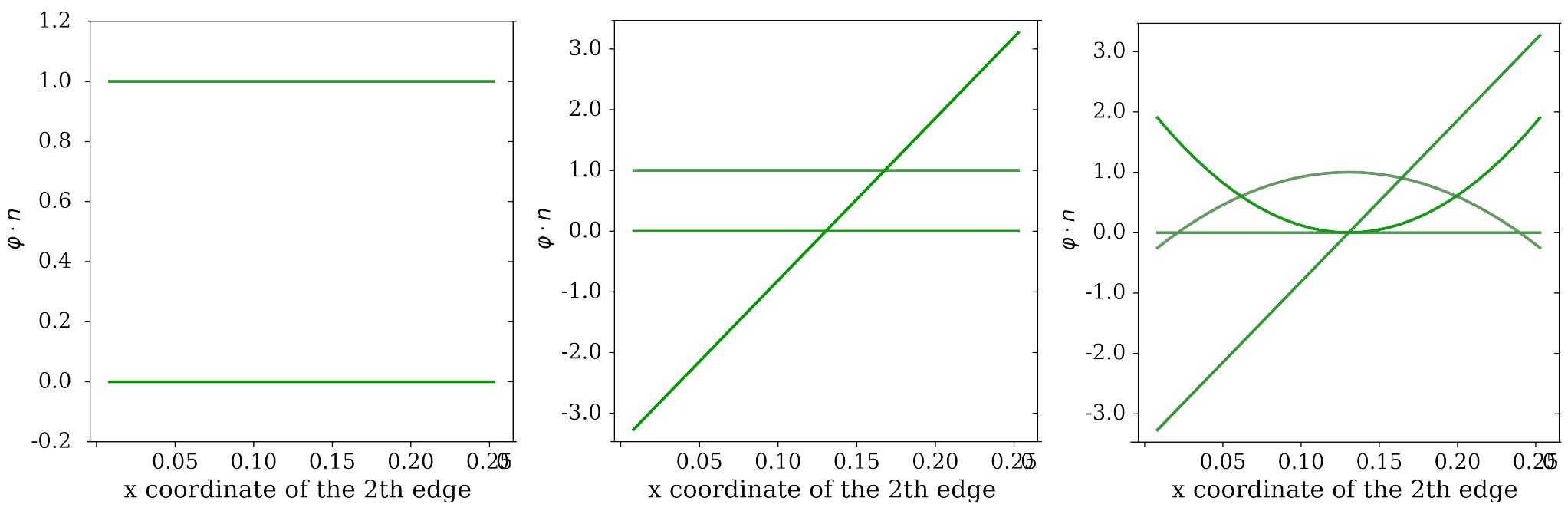

Fig. 107: Plot of all the functions emerging from the second edge on all the second edge of the element. From left to right: $k = 0$, $k = 1$, $k = 2$

Indeed, we can observe that the normal components of the normal basis function scales to one in the lowest order case, regardless the conditioning of the matrix. For the highest orders however, the Lagrangian behaviour is still not respected between the normal components of the first $k + 1$ normal functions.

Furthermore, not one but two basis functions seems to be missing. In fact, as it can be seen in the *Figure 108*, still only one normal basis function is degenerating into an internal function, but the normal component of the regular normal basis function merge. The first auxiliary normal basis function see its normal component vanishing on every edge and the second auxiliary normal basis function adopts the other direction of variation while possibly inheriting a much higher amplitude than the other normal basis functions.

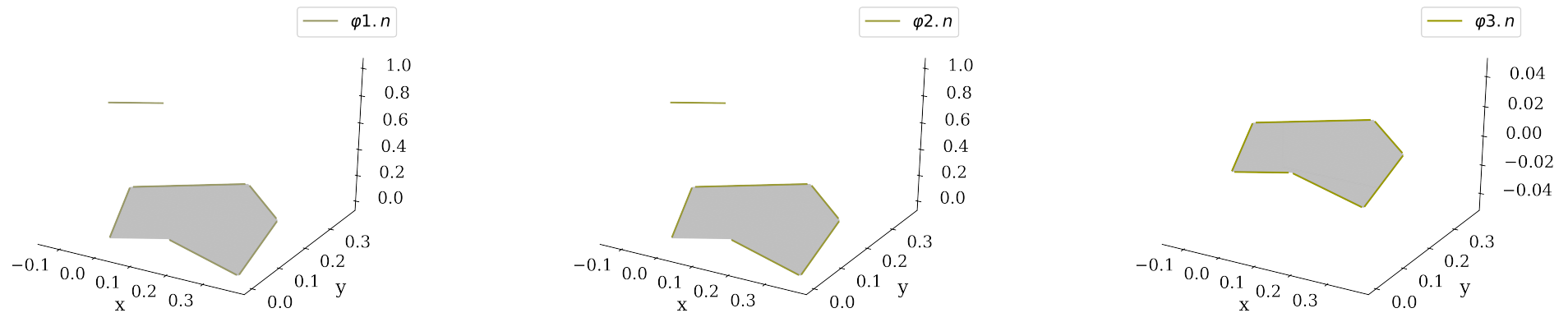

Fig. 108: Plot of the normal component of the basis functions emerging from the former canonical normal basis functions in the case $k = 1$. From left to right: first normal basis function, second normal basis function, third degenerating basis function generated from a canonical normal basis function.

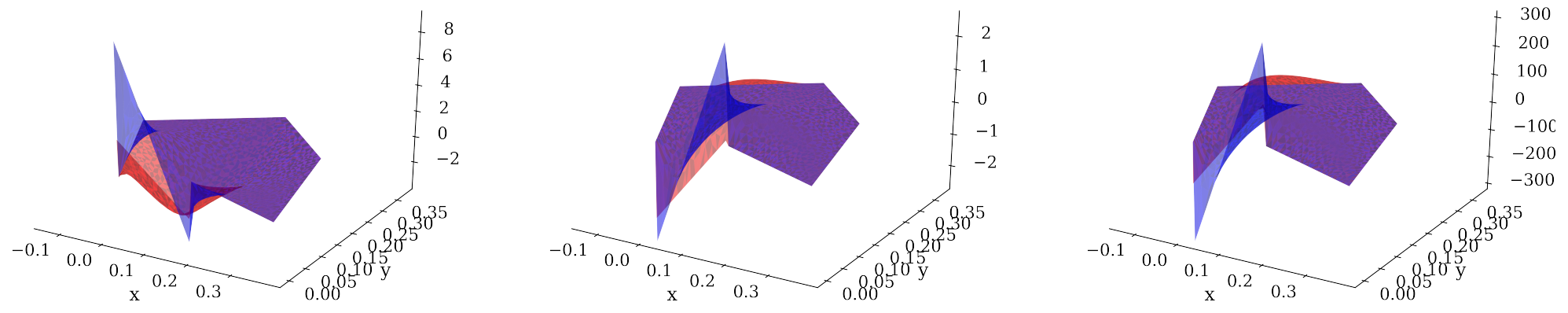

Fig. 109: Plot of the internal behaviour of the basis functions emerging from the former canonical normal basis functions in the case $k = 1$. From left to right: first normal basis function, second normal basis function, third degenerating basis function generated from a canonical normal basis function. None of the coordinate - wise components are vanishing identically on all the boundaries.

#### 7.3.1.3 Element IIa

*A - contrario* from the two previous elements, the element $IIa$ is built from the counter - intuitive definition of the normal degrees of freedom whose definitions are only making use of moments.

We start as before by studying the tuned basis functions generated from the canonical normal basis functions. As an example, we plot in the *Figure 110* the internal component of one of their representative having the highest variations, whose support matches the second edge. Its boundary behaviour is represented in the *Figure 111*.

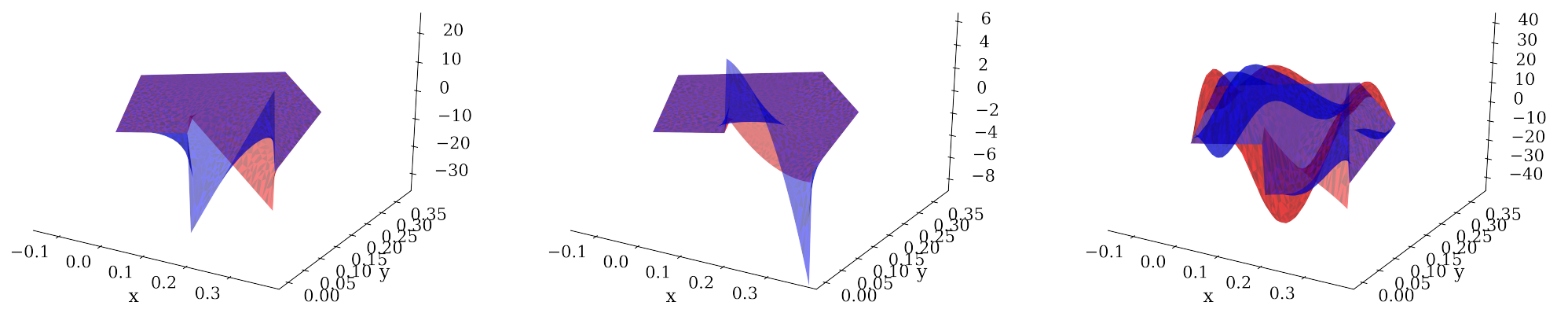

Fig. 110: Internal components of a basis function representative of the order of the space. Blue: $x$ component. Red: $y$ component. From left to right: $k = 0$, $k = 1$, $k = 2$.

Similarly as in the previous cases, the two components of the normal basis function are smooth within the element for any order $k$. Their amplitudes are also comparable to the one found in the case of the element $Ia$.

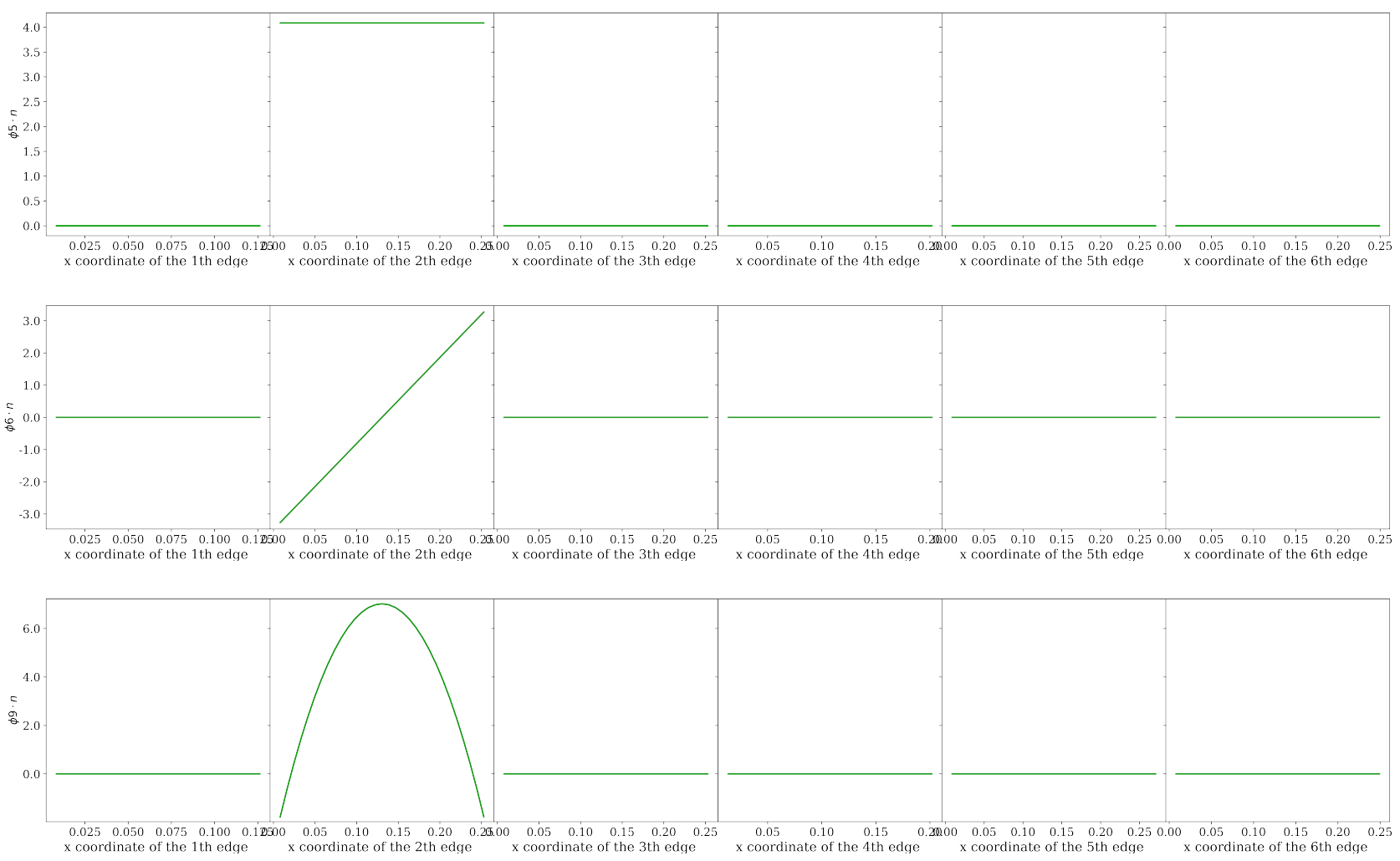

Fig. 111: Representative normal component of normal tuned basis functions that have a non - identically vanishing normal component. Eclated view on every edges. Top to bottom: $k = 0$, $k = 1$, $k = 2$.

Furthermore, the support of the normal component of the normal basis functions is also limited to a single edge and still allows an easier representation of the discretized quantity across the boundaries. The observed variations on its support are also still compliant with the order of the element. However, as for the element $Ia$, we do not enjoy anymore the local point - wise Lagrangian property, even for the lowest element. Indeed, the tuning method is not changed and the degrees of freedom still do not allow the assigned pointwise values to be kept during the tuning process. We picture the results in the *Figure 112* for the basis functions whose support matches the second edge.

There, the normal basis functions of the lower elements do not scale to one either. However, and this is the important change compared to the previously detailed elements, not one but two functions see their normal component vanishing on every edge. As pictured in the *Figures 113 and 114* they are not identically vanishing within the element. Thus, even if they are slightly different from the originally designed internal basis functions by their non - vanishing $x$ or $y$ components on the boundary, they can still be considered as such.

Note that so far, this element is the most compliant with the variations of the normal component of the basis functions on the boundary as only $k + 1$ polynomials are generated on the boundary.

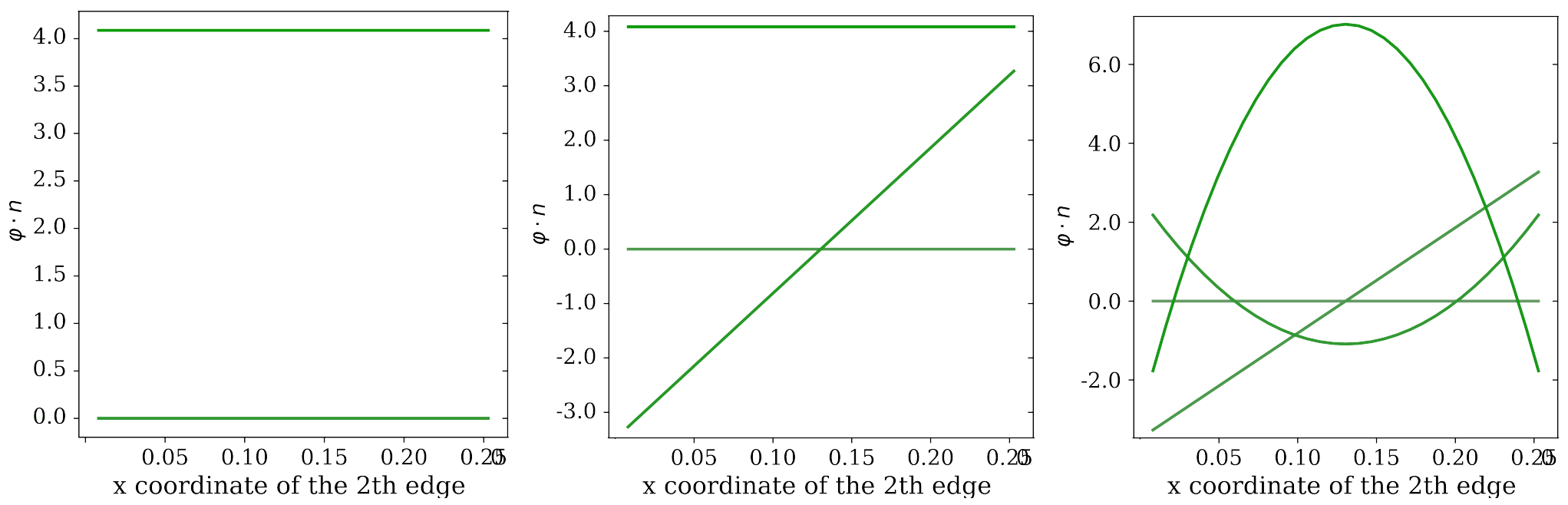

Fig. 112: Plot of all the functions emerging from the second edge on all the second edge of the element. From left to right: $k = 0$, $k = 1$, $k = 2$.

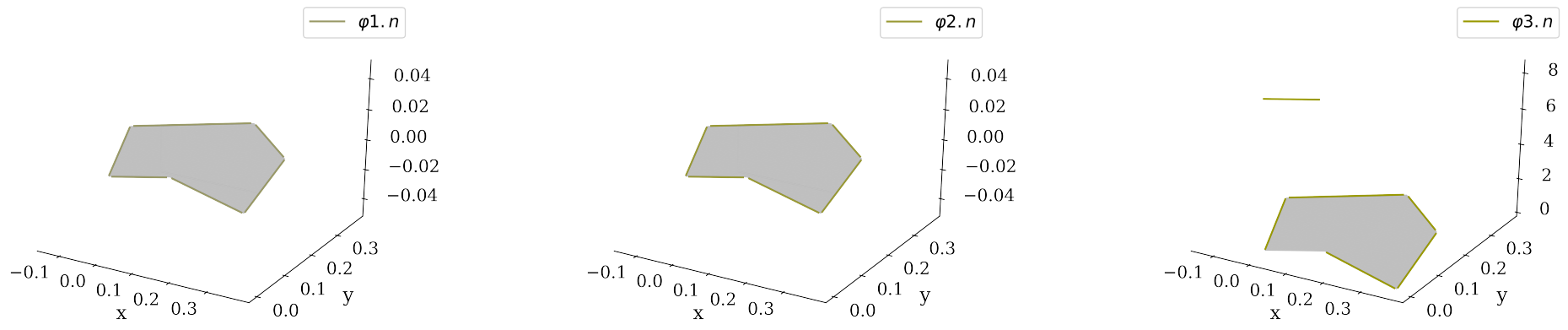

Fig. 113: Plot of the normal component of the basis functions emerging from the former canonical normal basis functions in the case $k = 1$. From left to right: first normal basis function, second normal basis function, third degenerating basis function generated from a canonical normal basis function.

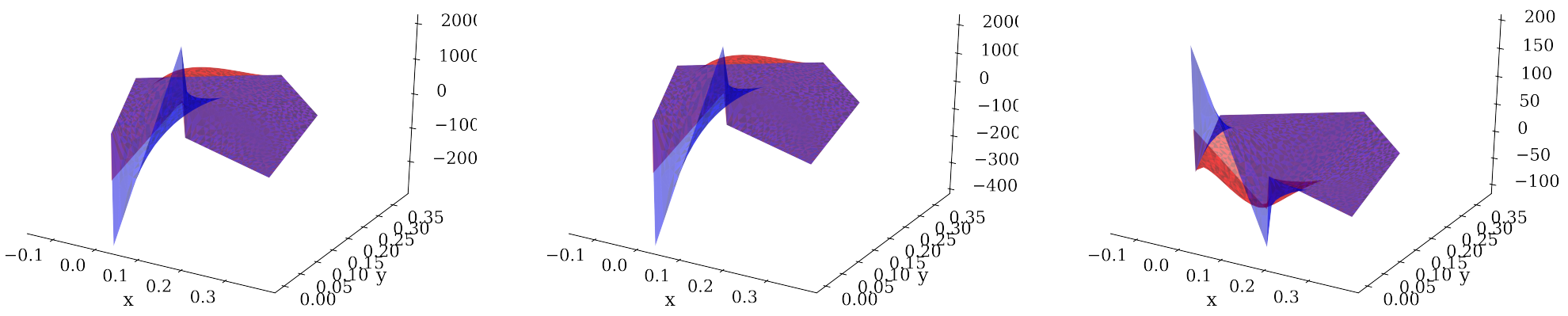

Fig. 114: Plot of the internal behaviour of the basis functions emerging from the former canonical normal basis functions in the case $k = 1$. From left to right: first normal basis function, second normal basis function, third degenerating basis function generated from a canonical normal basis function. None of the coordinate - wise components are vanishing identically on all the boundaries.

### 7.3.1.4 Element IIb

The last element combines the scaling advantages of the element $Ib$ with the advantage of the degeneration of the two auxiliary normal basis functions as internal basis functions that we encountered in the case of the element $IIa$.

First of all, let us point out that as all of our presented elements, the functions tuned from canonical normal basis functions enjoy smoothness within the element and when not vanishing, enjoy also a polynomial regularity of their normal component on the boundary. The internal behaviour of a representative normal functions enjoying the highest variations is given in the *Figure 115*, while its normal behaviour is pictured in the *Figure 116*.

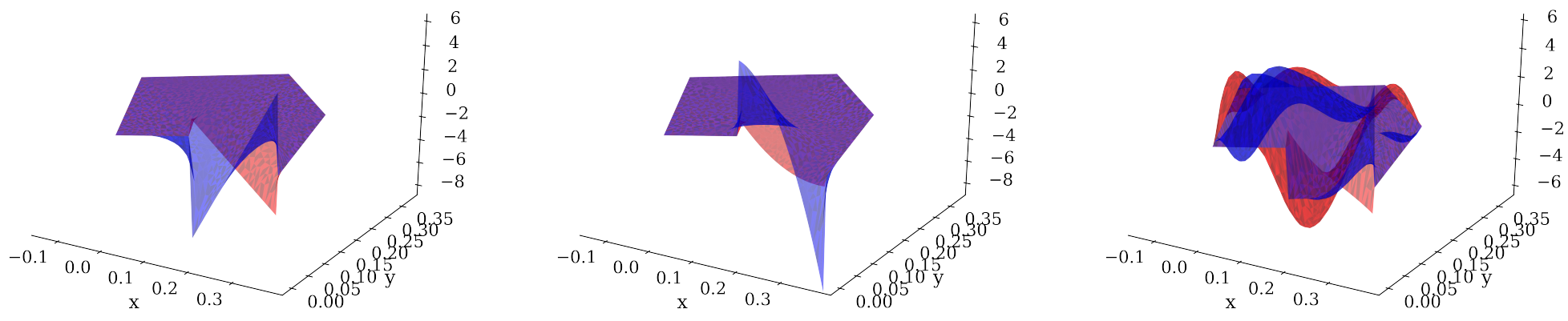

Fig. 115: Internal components of a basis function representative of the order of the space. Blue: $x$ component. Red: $y$ component. From left to right: $k = 0$, $k = 1$, $k = 2$.

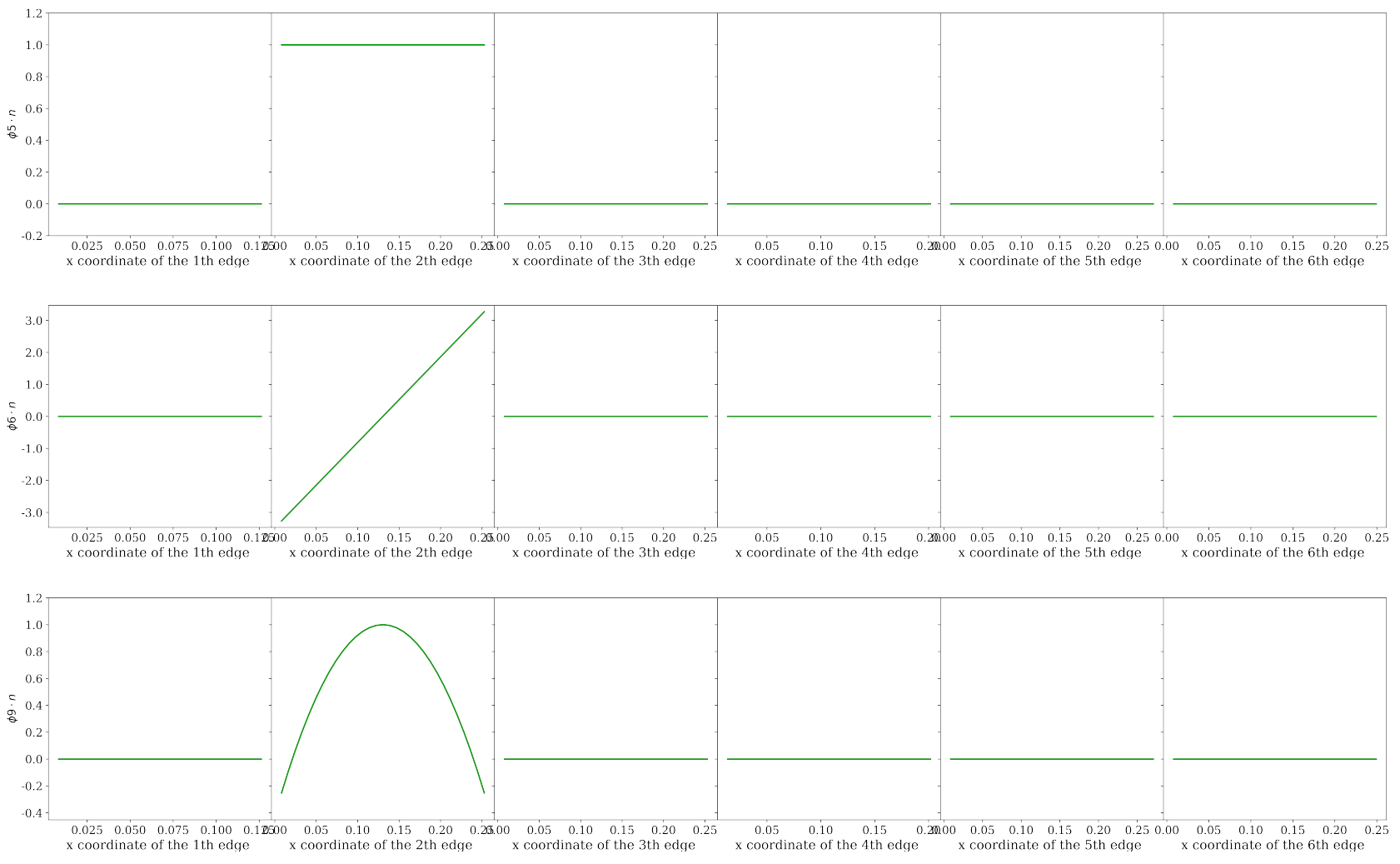

Fig. 116: Representative normal component of normal tuned basis functions that have a non - identically vanishing normal component. Eclated view on every edges. Top to bottom: $k = 0$, $k = 1$, $k = 2$.

On a local point of view, we still do not enjoy any Lagrangian property. Indeed, for the same reason as above the tuning process do not allow to keep the pointwise values, and the degrees of freedom are not enforcing a Lagrangian behaviour for the $k+1$ functions having a non - vanishing normal component. We plot the results of functions living on the second edge in the *Table 117*. However, for the lowest order element the normal components of the basis functions do scale to one, when not vanishing.

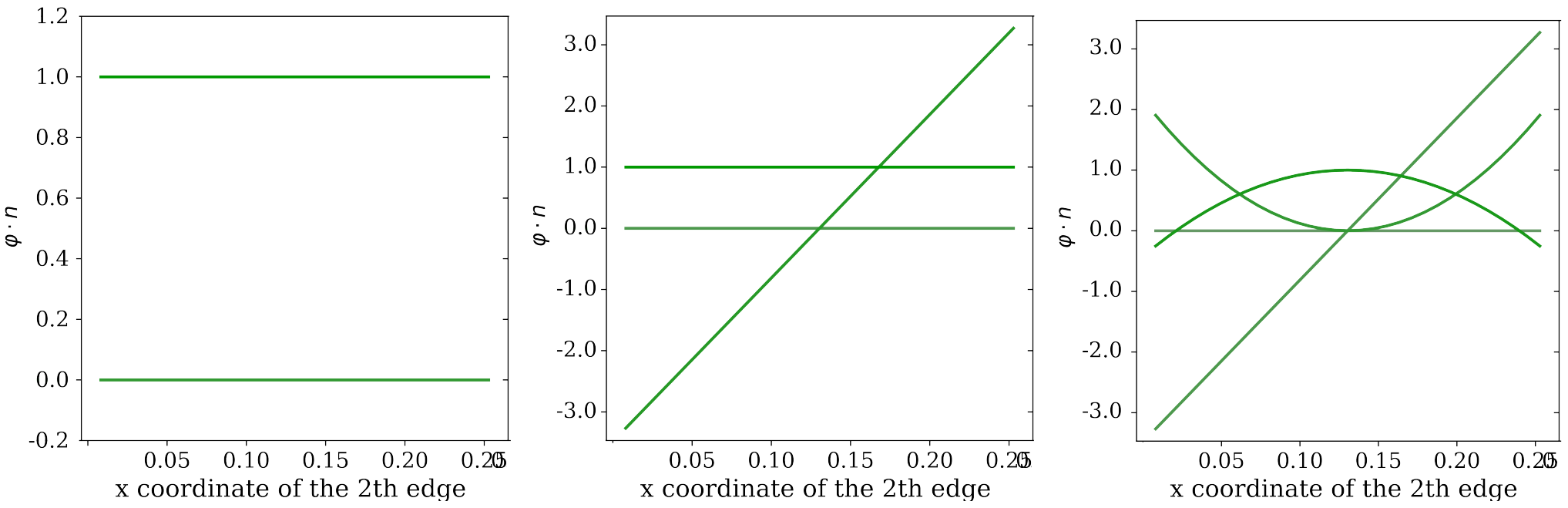

Fig. 117: Normal functions emerging from the second edge.

Lastly, we can witness that similarly as before, the functions whose normal component vanishes on the edges are not identically vanishing *(see the Figures 118 and 119)*. Thus, due to their regularity they can be reuptaken as an internal basis function, even if component wise they are not vanishing on the boundary.

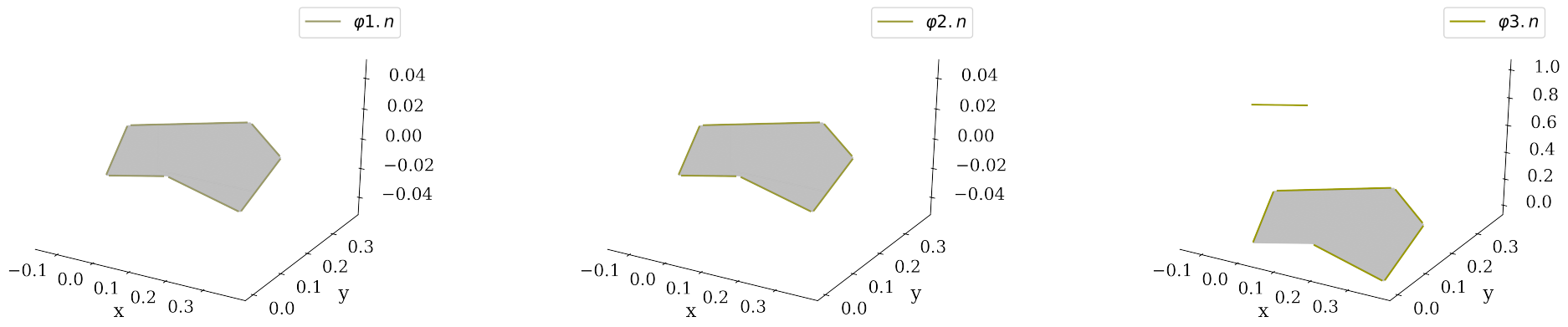

Fig. 118: Plot of the normal component of the basis functions emerging from the former canonical normal basis functions in the case $k=0$. From left to right: first normal basis function, second normal basis function, third degenerating basis function generated from a canonical normal basis function.

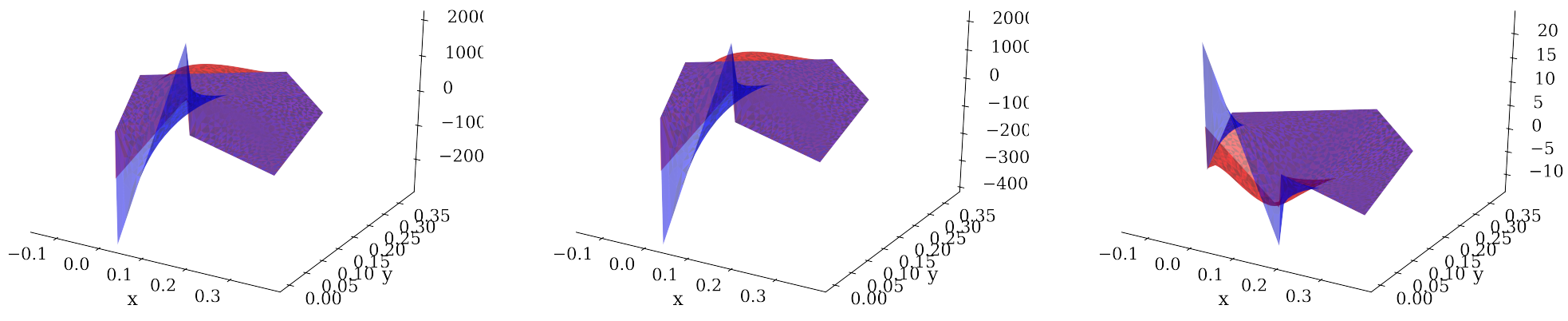

Fig. 119: Plot of the internal behaviour of the basis functions emerging from the former canonical normal basis functions in the case $k = 0$. From left to right: first normal basis function, second normal basis function, third degenerating basis function generated from a canonical normal basis function. None of the coordinate - wise components are vanishing identically on all the boundaries.

### 7.3.2 Impact of the shift in the pointwise value

Let us do a miscellaneous study for the two elements that comprises point - wise values of normal components. We would like to see the impact of the offset prescribed in the pointwise value. To this sake, we select a pointwise normal value and shift it by one below. By example, instead of prescribing $q \mapsto q(\mathrm{x}_m) \cdot n$, we define the corresponding degree of freedom by $q \mapsto q(\mathrm{x}_m) \cdot n - 1$. We then look at the behaviour of the tuned normal basis functions.

**For the element $Ib$,** after shifting in the same way both component wise pointwise values, it comes the results presented in the *Figure 120.*

Then, when flattening the previous graph, we retrieve boundary behaviour presented in the *Figure 121.* We can here observe that as one would expect, this shift impacts all the normal basis functions of all the edges homogeneously. All the normal components that were previously vanishing are lifted up. The unitary value of functions whose components were not vanishing is however preserved.

One can also see that the transfer matrix is weaken. Indeed, the conditioning of the matrix considering a shifted point value is about 58818 while the conditioning of the transfer matrix for the regular element $Ib$ value 9188.

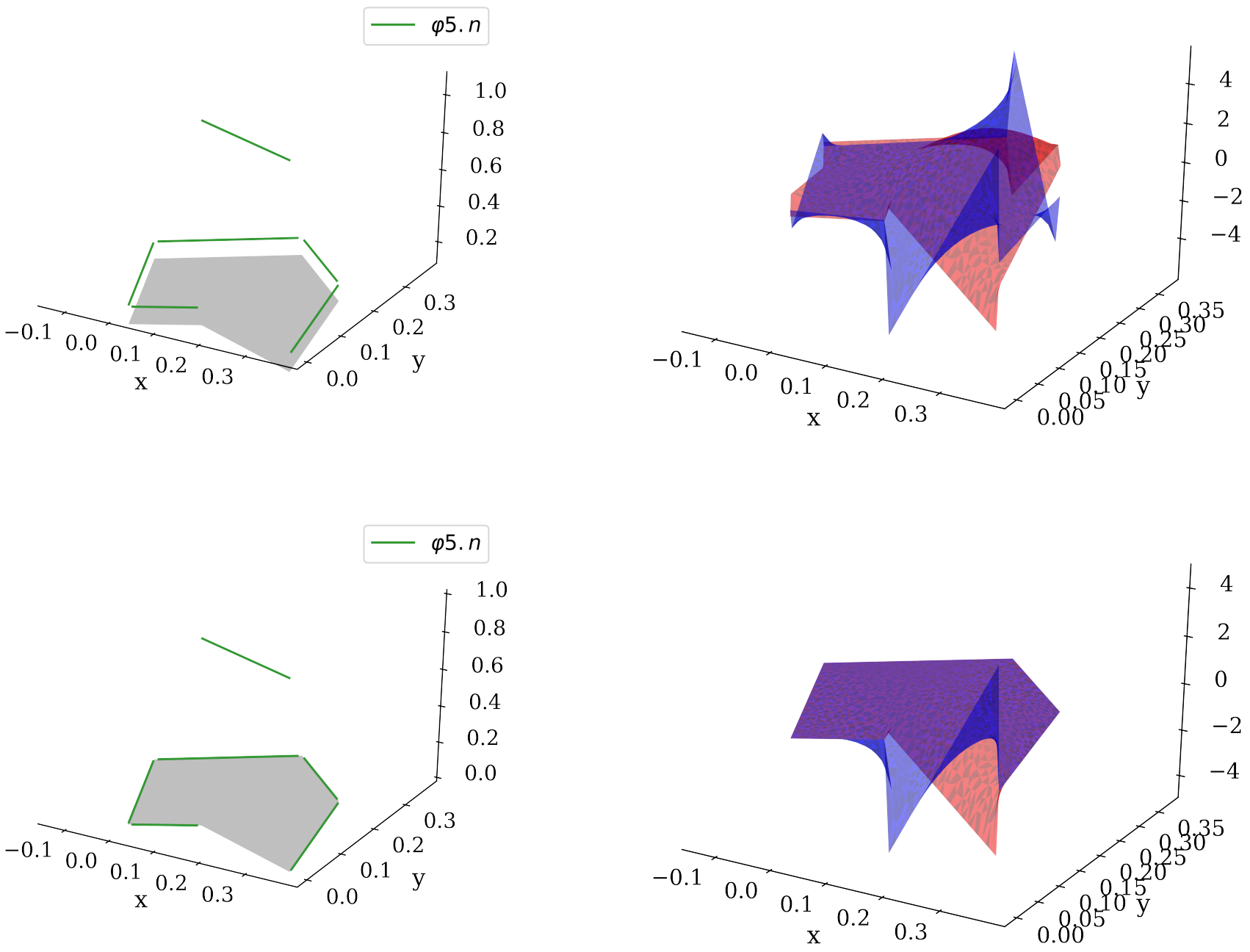

Fig. 120: First normal moment built from the first edge in the case $k = 0$. Left: shifted pointwise value. Right: classical setting.

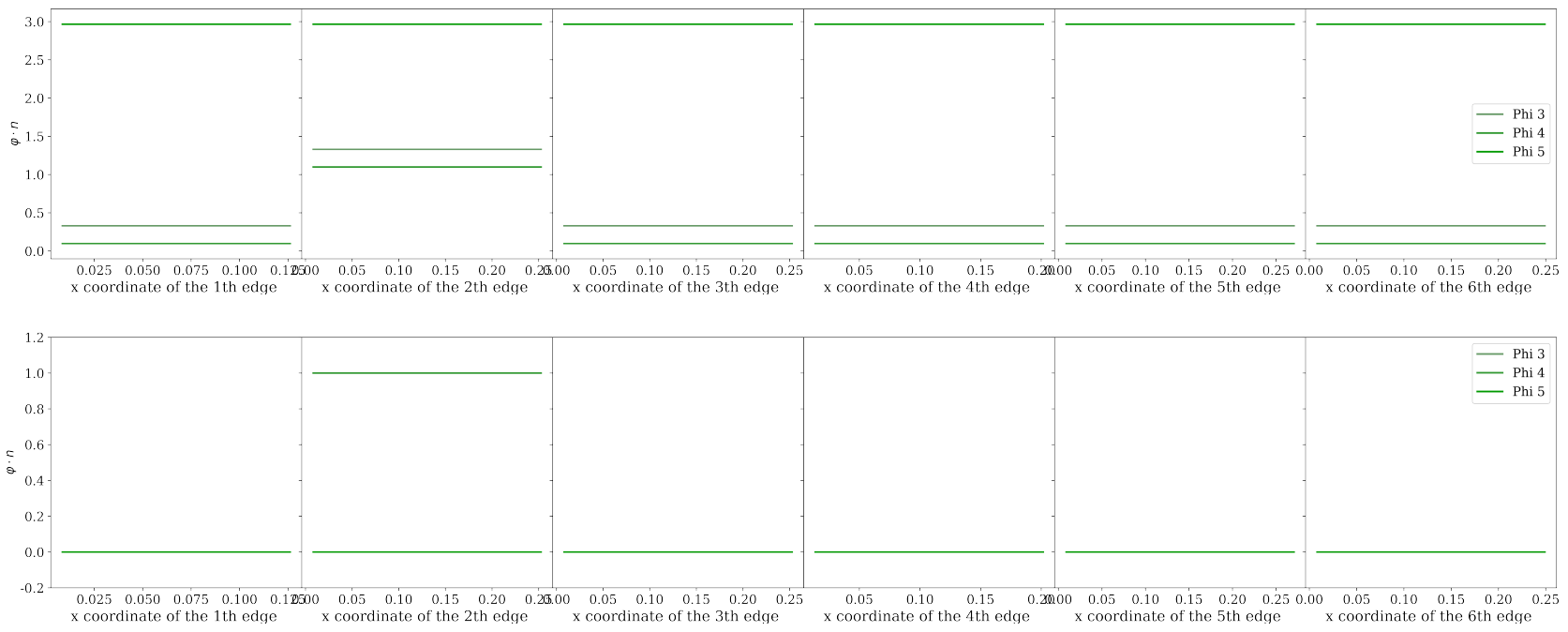

Fig. 121: First normal moment built from the first edge in the case $k = 0$: Top Ib shifted, Bottom Ib

**For the element** $IIb$**,** the conclusions are identical though here only one pointwise value is used. Indeed, by example for $k = 0$, it comes the following results.

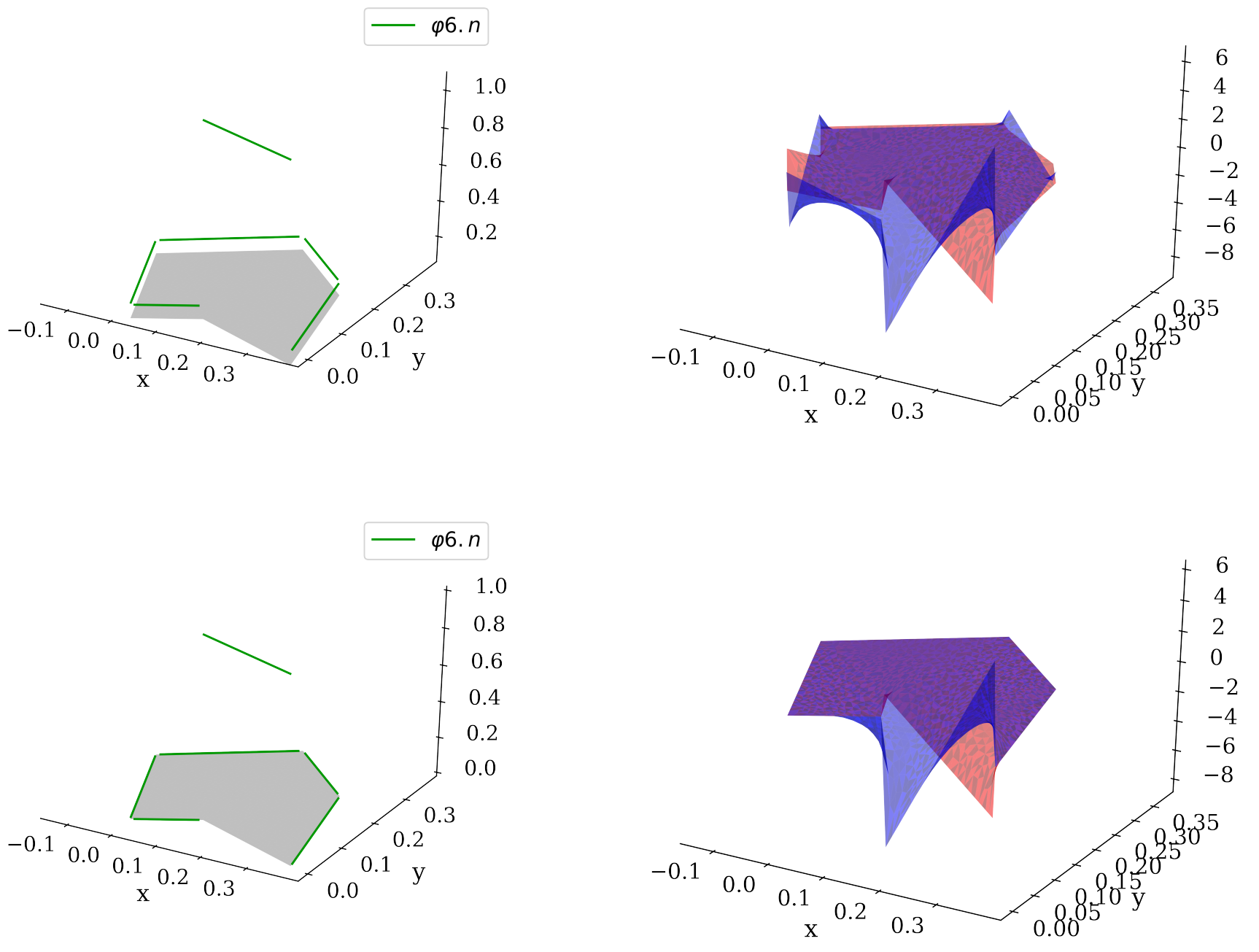

Fig. 122: First normal moment built from the first edge in the case $k = 0$. Left: shifted pointwise value. Right: classical setting.

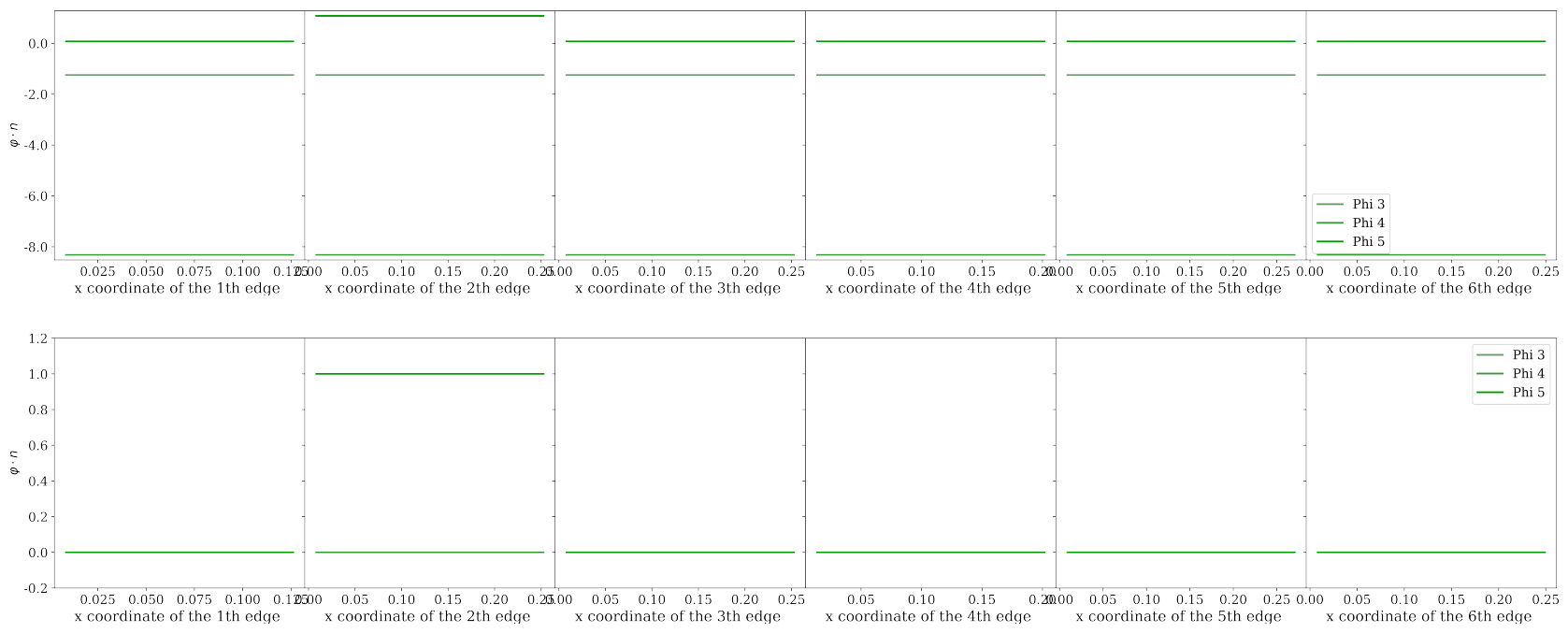

Fig. 123: First normal moment built from the first edge in the case $k = 0$: Top Ib shifted, Bottom Ib

As before, when flattening the previous graph in the above figure, we get where same conclusions as in the case of the element *Ib* can be drawn. Especially, we get a conditioning of 328924 for the shifted set of degrees of freedom against 119664 for the classical setting.

Therefore, if the behaviour of the shifted element is as one would expect, it is better to stick to the classical setting.

### 7.3.3 Impact of the projector on the shape of the normal component of the tuned basis functions

For the sake of completeness, let us now study the impact of the choice of the polynomial projectors on the shape of the normal components of the basis functions. Note that the choice of the polynomial constructors has no impact there, provided that they lead to decent matrices conditionings.

In all the cases presented below, we only consider the case of the non - convex hexagon presented in the *Figure 90* and plot all the normal component of the basis functions whose support matches the second edge.

***Note.*** As the Laguerre boundary projector leads to unreliable results (see the discussion on the conditionings), this case has been excluded from the discussion. ▲

***Remark.*** The layout of the pictures presented below corresponds to the following sketch.

| | | |
|---|---|---|
| Centred - scaled canonical polynomials | Chebyshev polynomials | Hermite polynomials |
| Legendre polynomials | Centred - unscaled canonical polynomials | Not centred - unscaled canonical polynomials |

Fig. 124: Normal components of functions emerging from the second edge.

▲

#### 7.3.3.1 First configuration, $I$

Let us start with the first configuration. We plotted the first order case in the *Figures 125* for the element $Ia$ and *127* for the element $Ib$. The second order case is reported in the *Figures 126* for the element $Ia$ and *128* for the element $Ib$.

First of all, we can observe for the type $a$ that the conditioning of the transfer matrix impacts on the amplitude of the basis functions. Indeed, for both orders $k = 1$ and $k = 2$, the two cases of non - centred unscaled canonical polynomials and centred unscaled canonical polynomials, whose respective colour codes are both located in the high conditioning area of the conditioning graphs, lead to amplitudes being at least ten times higher than the one obtained in the other cases. In particular, one to two functions are dominating the other ones, which creates disequilibrium in the representation of functions living in the space $\mathbb{H}_k(K)$. Reversely, the Hermite polynomial case see its corresponding basis functions enjoy the smallest obtained amplitude.

One can also observe that for all cases but the ones involving non - centred canonical polynomials, two basis functions are constant in the case $k = 1$, and one function is linear in the case $k = 2$. Thus, despite the bad conditioning one can still want to use the unscaled canonical functions so that any basis functions has the same order, corresponding to the maximal one the normal component of functions of $\mathbb{H}_k(K)$ can achieve.

Similar observations can be derived for the element $Ib$. However, note that in the cases of centred projectors, the two constant functions in the case $k = 1$ merge. Furthermore, the impact of the definition of the point value shows up when using non - centred polynomials. Indeed, when the projectors are unscaled, the point value contributes to the tuning for all the monomial. Its contribution to the shape of the tuned basis functions is then diffuse. However, the above discussion on the shift would still hold. This can be observed on the bottom right of the *Figure 128* where the tuned basis functions are shifted not only in its value but also along its corresponding edge.

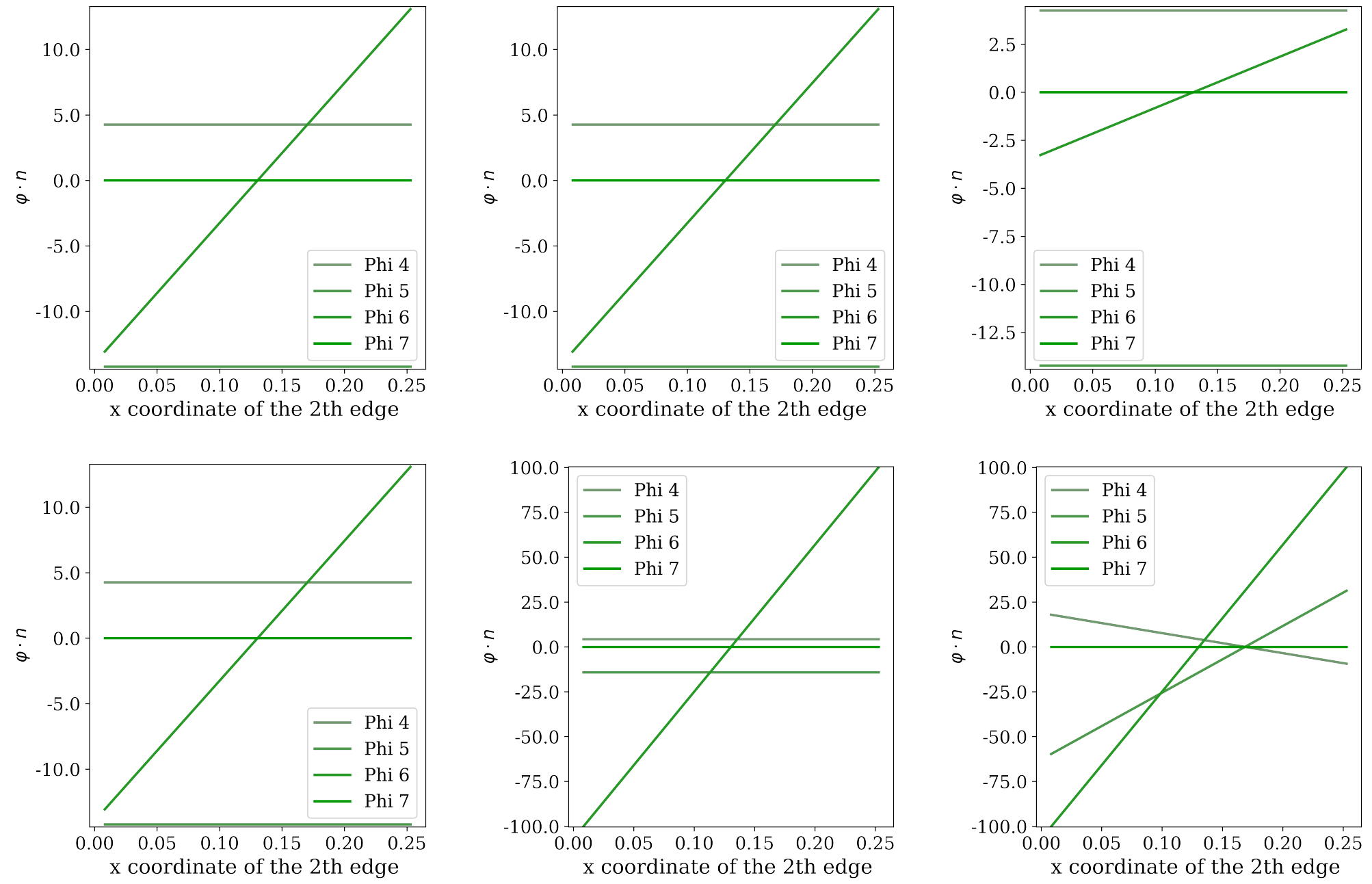

Fig. 125: Normal components of basis functions living on the second edge. Computed for the element $Ia$ on a non - convex hexagon, $k = 1$.

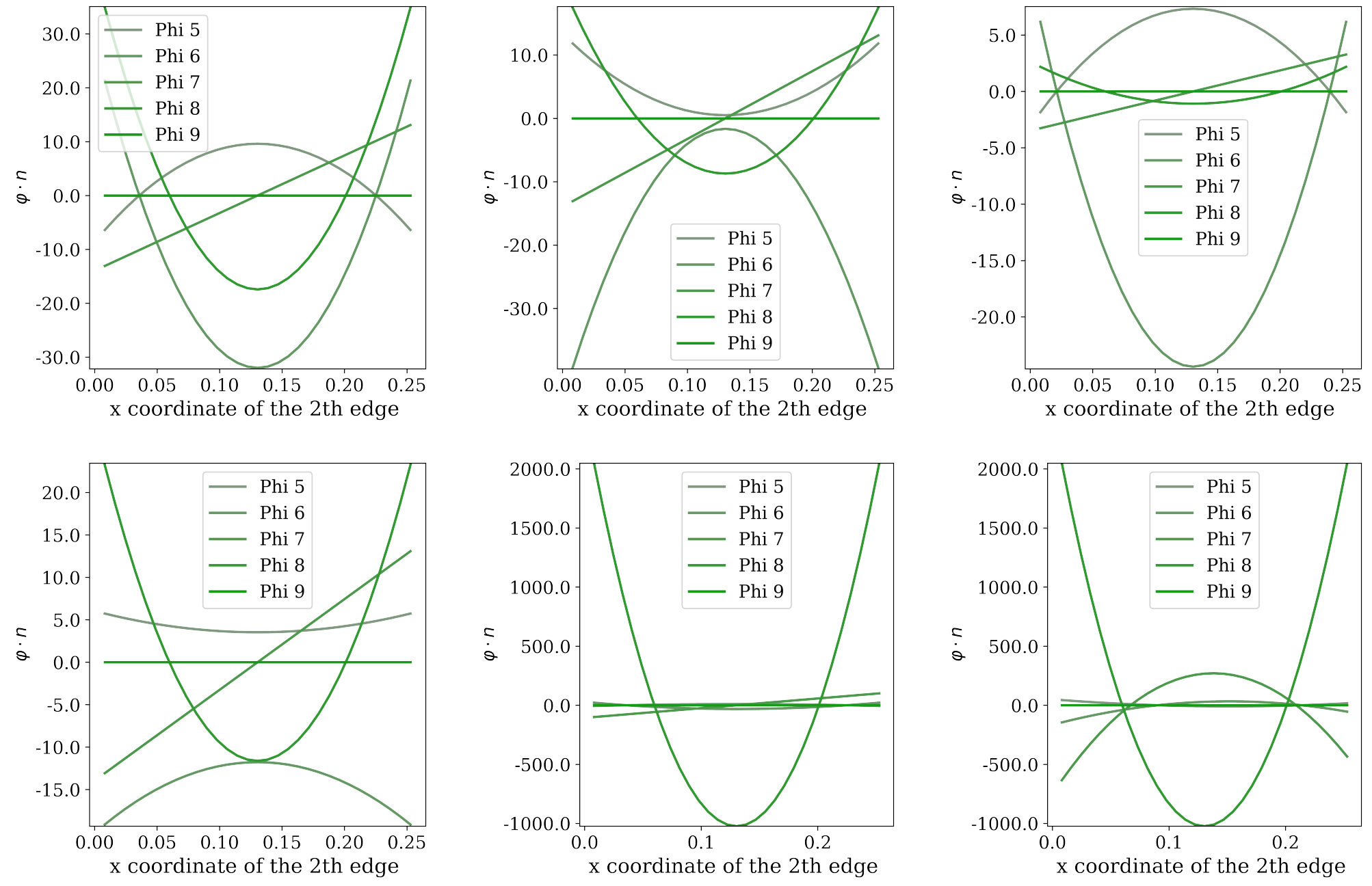

Fig. 126: Normal components of basis functions living on the second edge. Computed for the element $Ia$ on a non - convex hexagon, $k = 2$.

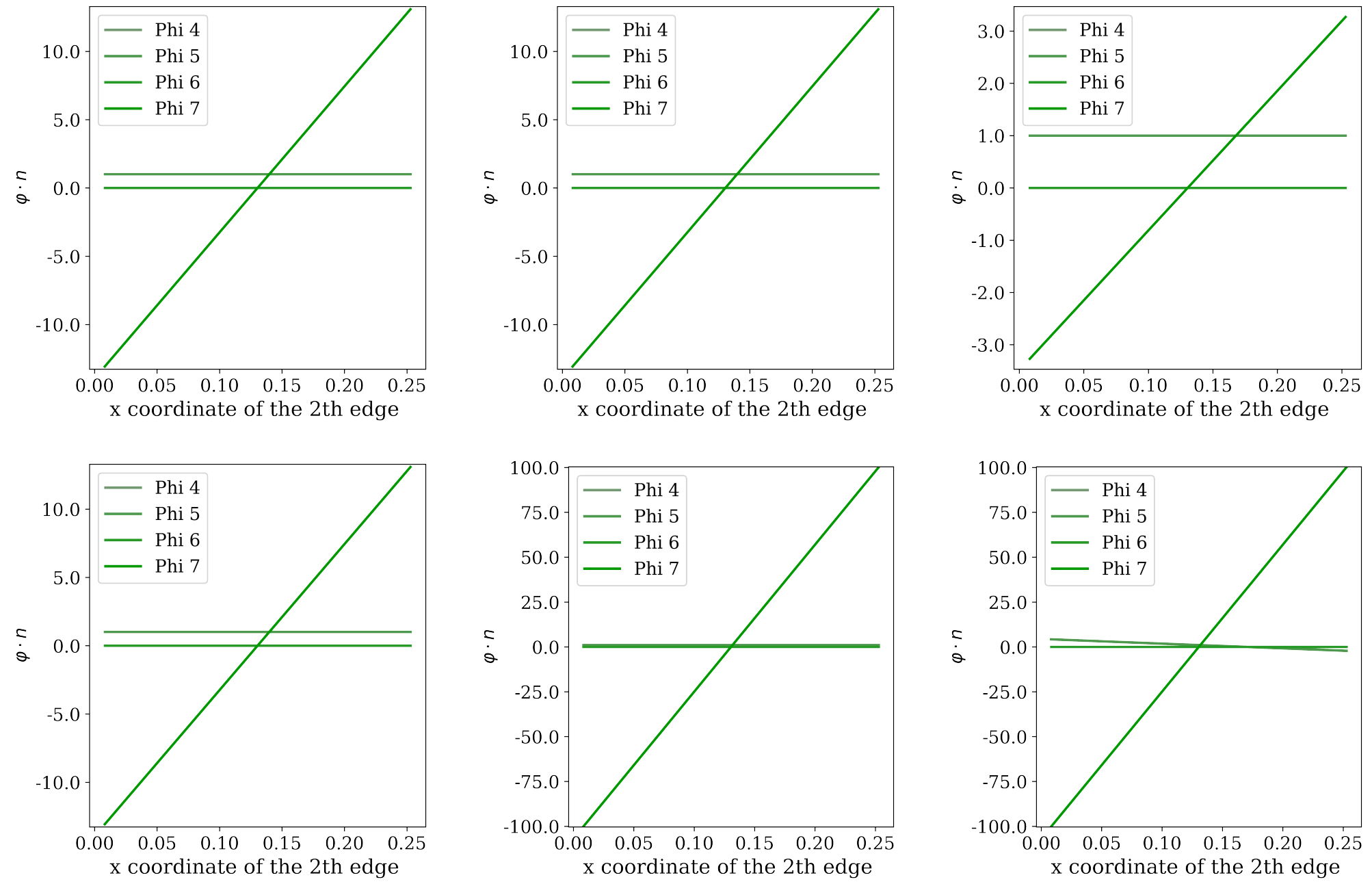

Fig. 127: Normal components of basis functions living on the second edge. Computed for the element $Ib$ on a non - convex hexagon, $k = 1$.

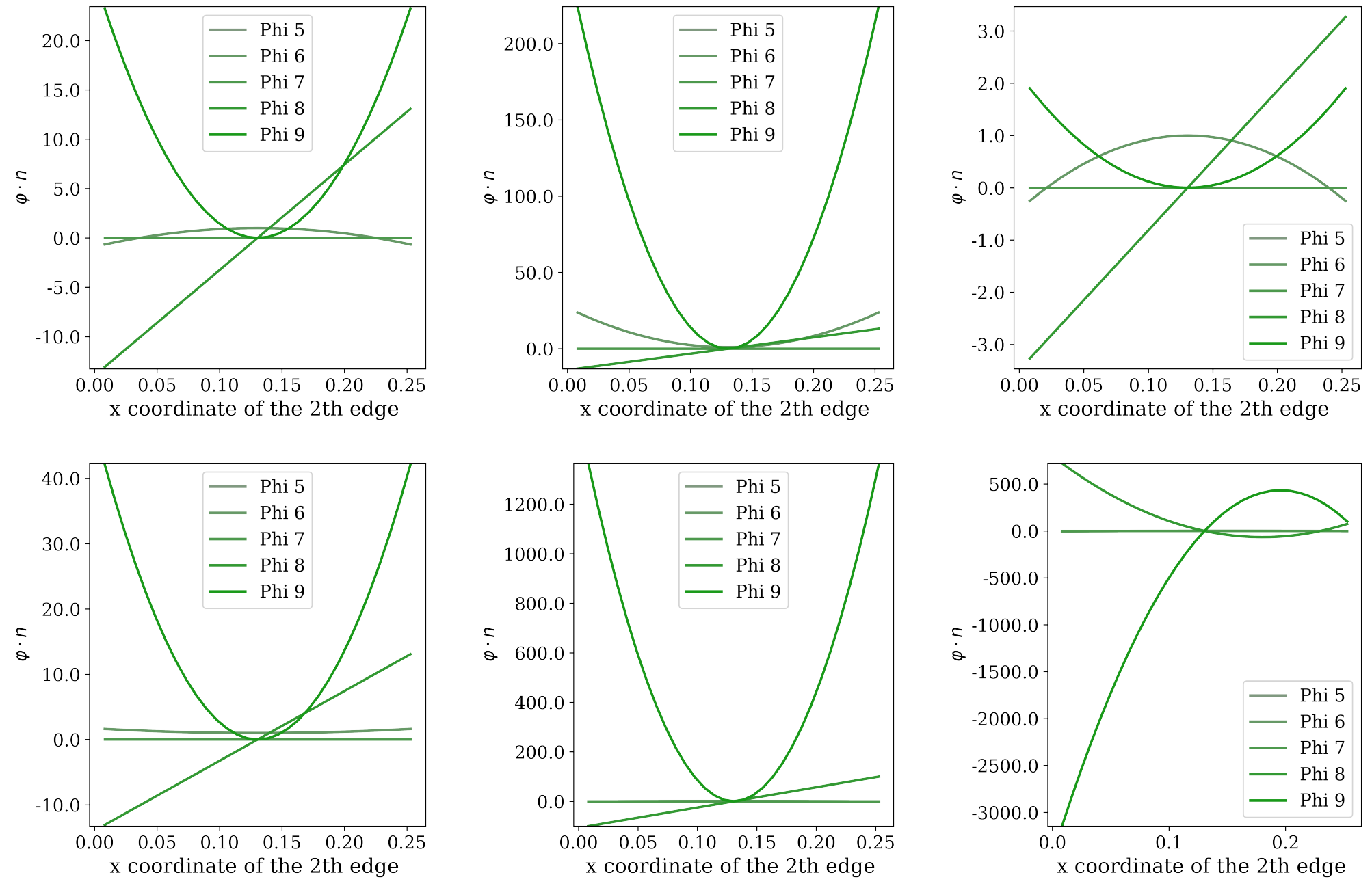

Fig. 128: Normal components of basis functions living on the second edge. Computed for the element $Ib$ on a non - convex hexagon, $k = 2$.

#### 7.3.3.2 Second configuration, $II$

The second configuration $II$ achieves a similar behaviour. We plotted the first order case in the *Figures 129* for the element $IIa$ and *131* for the element $IIb$. The second order case is reported in the *Figures 130* for the element $IIa$ and *132* for the element $IIb$.

First of all, for any type we observe the same behaviour in the basis functions' order as in the previous case, especially regarding the amplitude of the functions. In particular, in the case $k = 1$ we witness one basis function of order zero and one basis functions of order one for any centred projector. The only non - centred case is also still presenting two functions with linear variations. However, note that as here two misc basis functions degenerates, there is only one constant function rather than two in the previous case.

This observation can be extended to the second order case, where among the three non - vanishing normal components of basis functions are two quadratic functions and one function linearly increasing.

Lastly, let us remark that in the case of the element $IIb$, any centred projectors scales the constant function to one for the order $k = 1$. Furthermore, there too the Hermite projectors lead to the most reliable results in terms of the boundary characterisation of functions living in $\mathbb{H}_k(K)$.

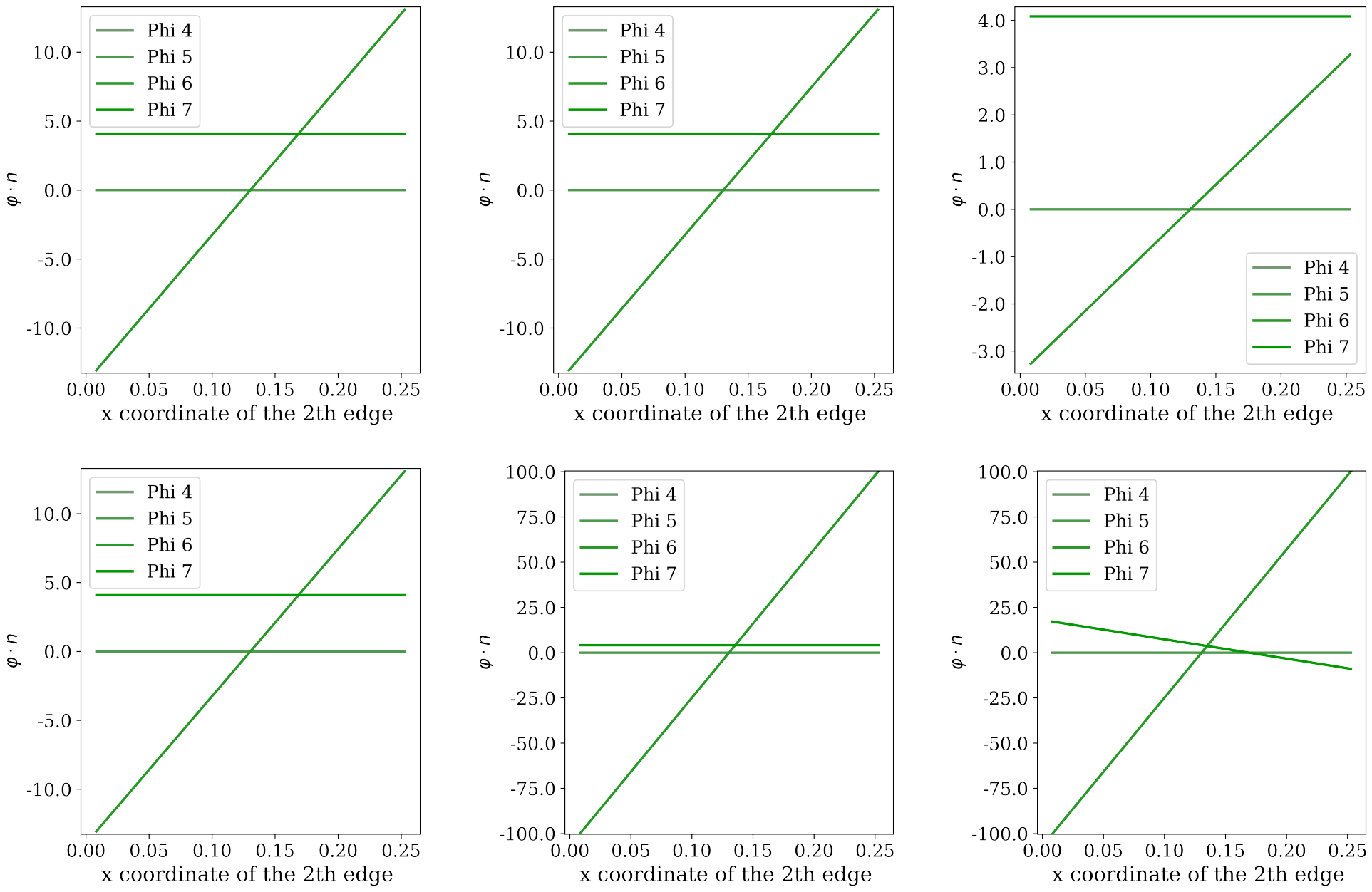

Fig. 129: Normal components of basis functions living on the second edge. Computed for the element $IIa$ on a non - convex hexagon, $k = 1$.

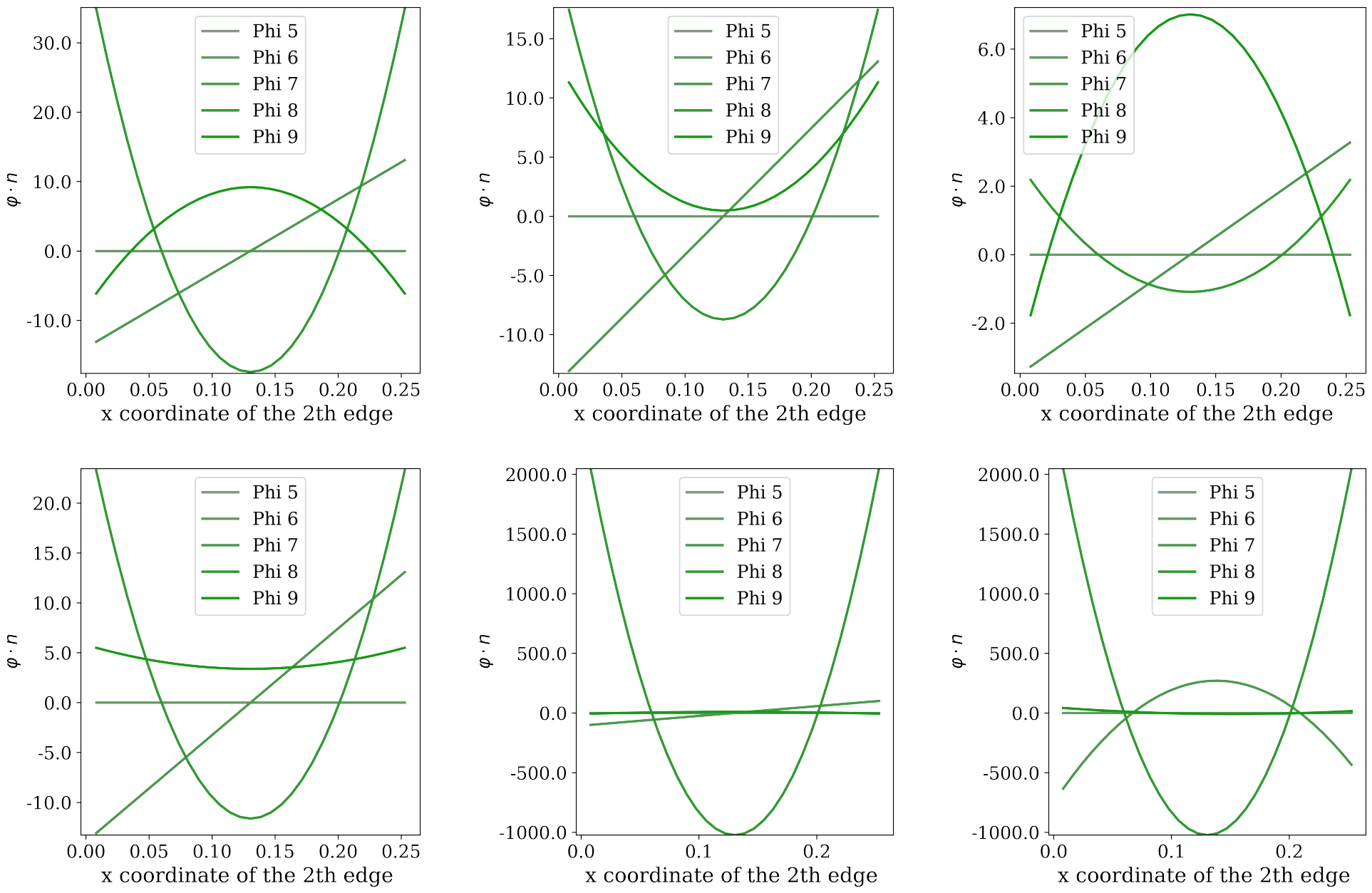

Fig. 130: Normal components of basis functions living on the second edge. Computed for the element $IIa$ on a non - convex hexagon, $k = 2$.

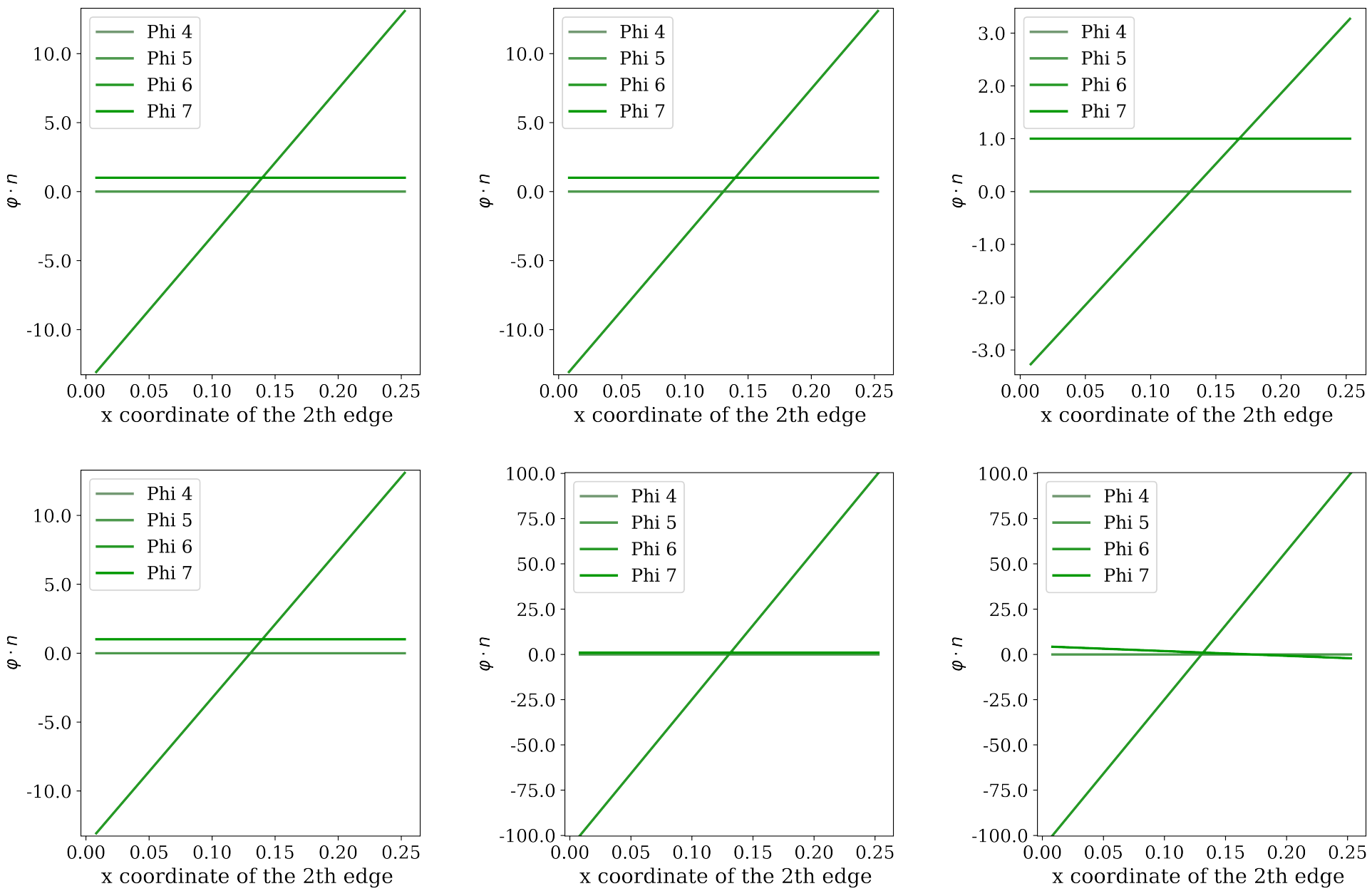

Fig. 131: Normal components of basis functions living on the second edge. Computed for the element $IIb$ on a non - convex hexagon, $k = 1$.

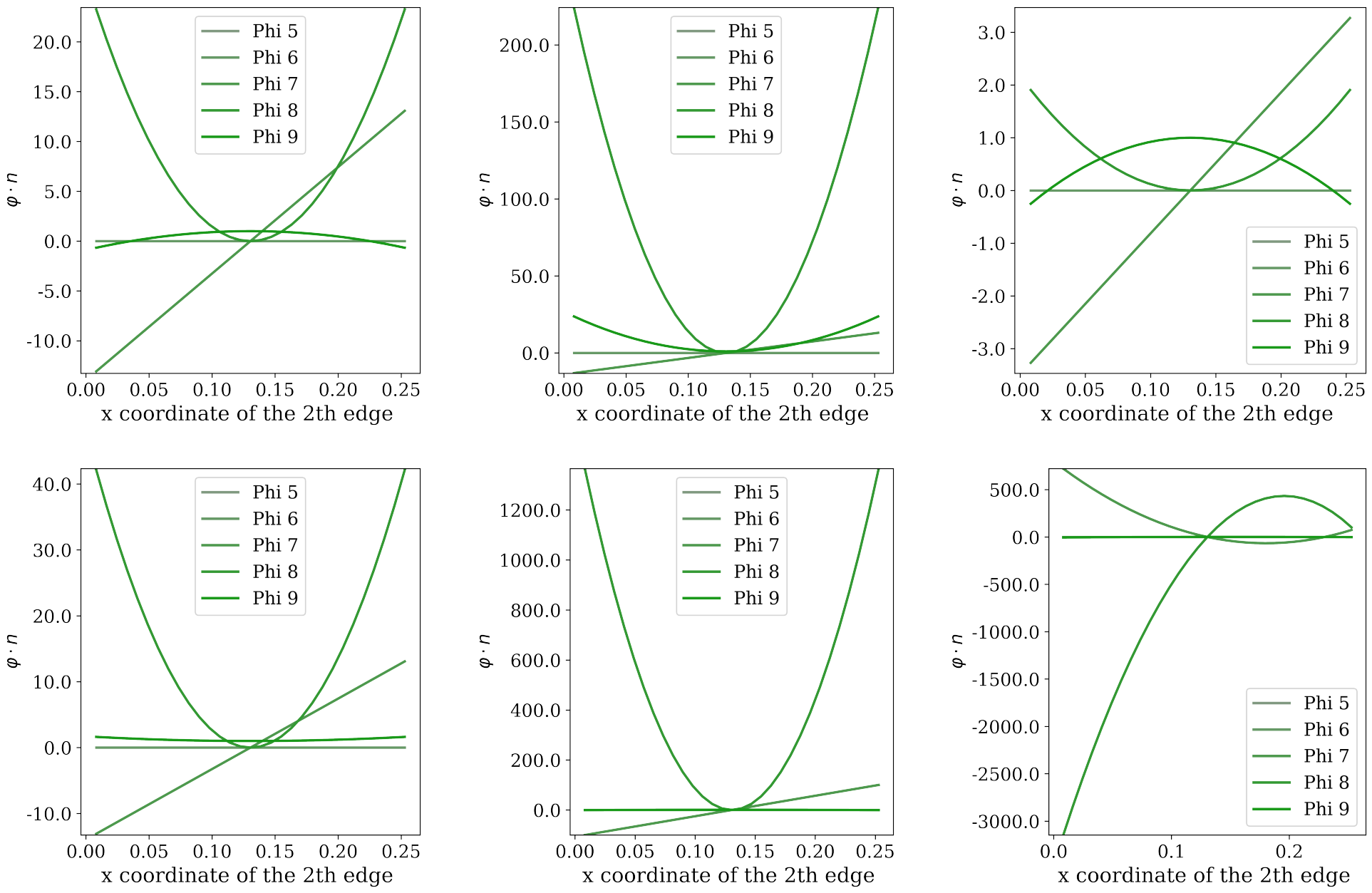

Fig. 132: Normal components of basis functions living on the second edge. Computed for the element $IIb$ on a non - convex hexagon, $k = 2$.

## 7.4 Comparison with classical RT

As a last investigation, let us compare the results obtained from the reduced element $IIb$ presented in the *Section 6.2* with the classical Raviart – Thomas elements built on simplicial and quadrangular shapes.

### 7.4.1 Triangular shape

We start by considering the triangular shape represented in the *Figure 133*. Note that because of our shape restrictions we could not study the exact same triangle as the one presented in the Raviart – Thomas section, but had to use a rotated one.

For the sake of concision we will only represent the normal basis functions on the first edge. The basis functions constructed from two other edges observe the same behaviour.

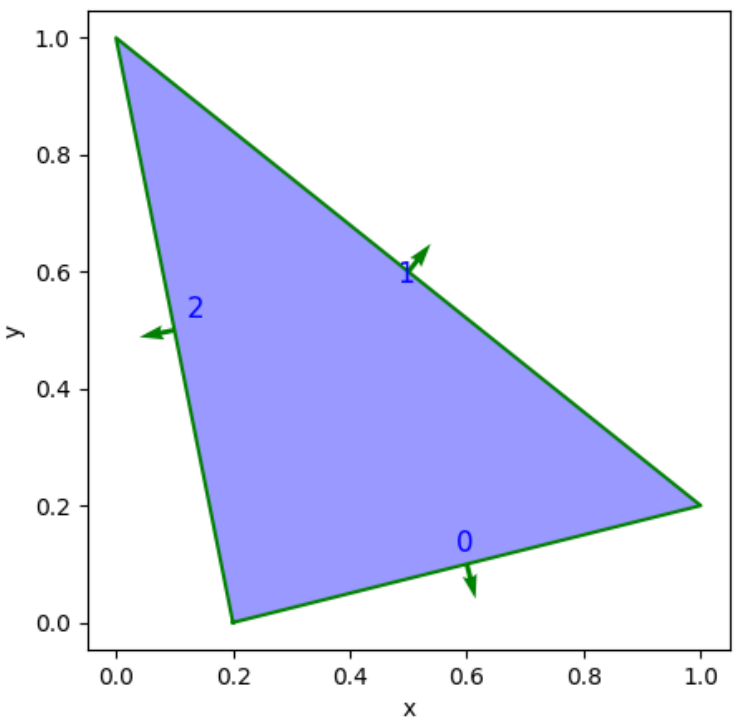

Fig. 133: Triangle

Let us start by deriving the behaviour of the normal basis functions of the lowest order space. There, the representative normal basis function on the first edge reads as follows.

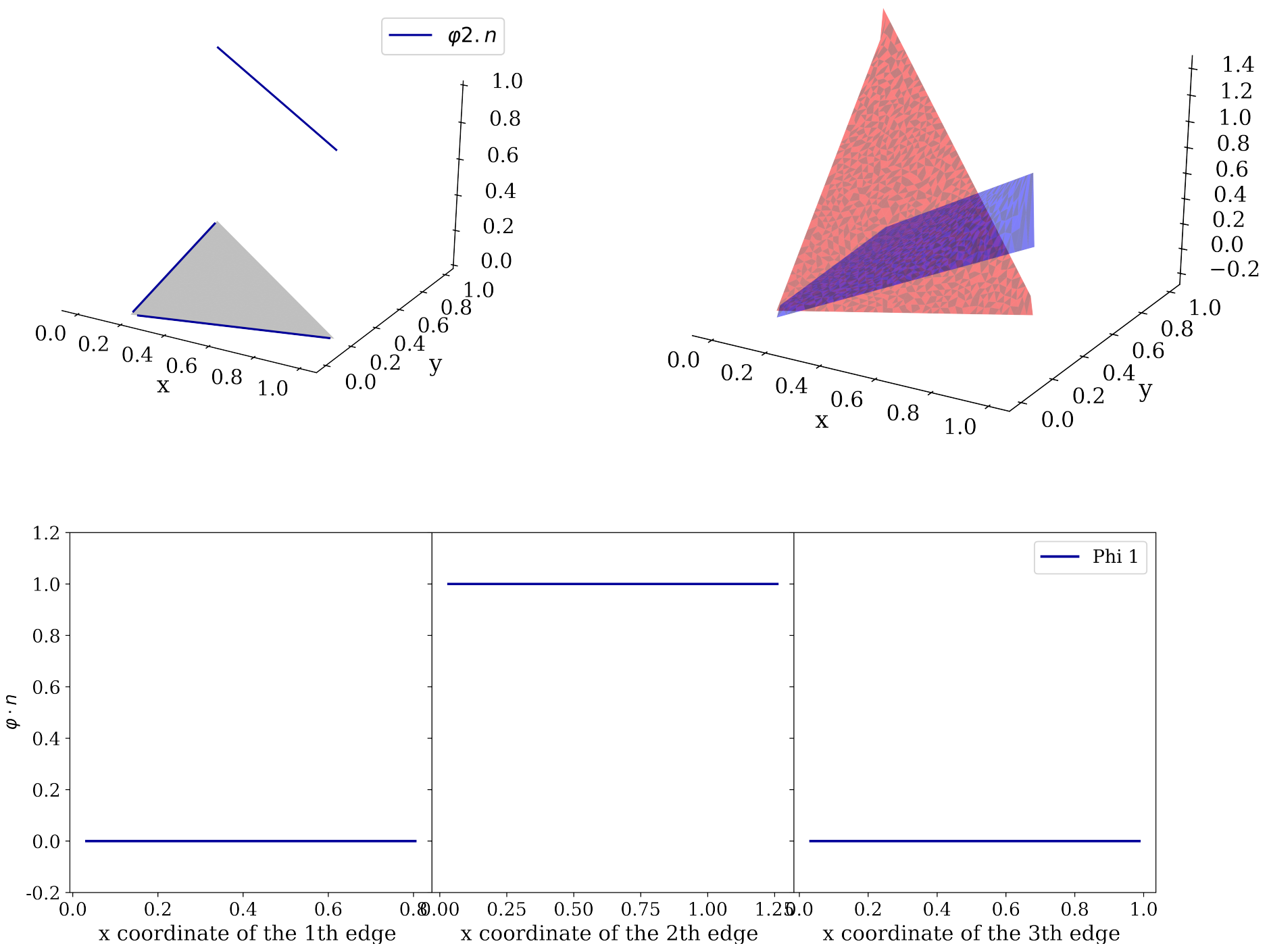

Fig. 134: First basis function in the case $k = 0$. Top left: normal component along the boundaries. Top right: internal behaviour of the basis function. Blue: $x$ component. Red: $y$ component.

As in the Raviart – Thomas setting, the only normal basis functions whose normal component is non - vanishing scales to one. However, note that within the element its behaviour is not polynomial.

Similarly, for a first order element presented in the *Figure 135*, the comparison with the Raviart – Thomas element on the boundary is immediate. Especially, there is no degeneration of canonical normal basis function into tuned internal ones. Furthermore, the amplitude of the basis functions is still contained within a range that allows reliable a decomposition of functions living in $\mathbb{H}_k(K)|_{\partial K}$. The last case presented in the *Figure 136*, for $k = 2$, drives to similar conclusions on the edges.

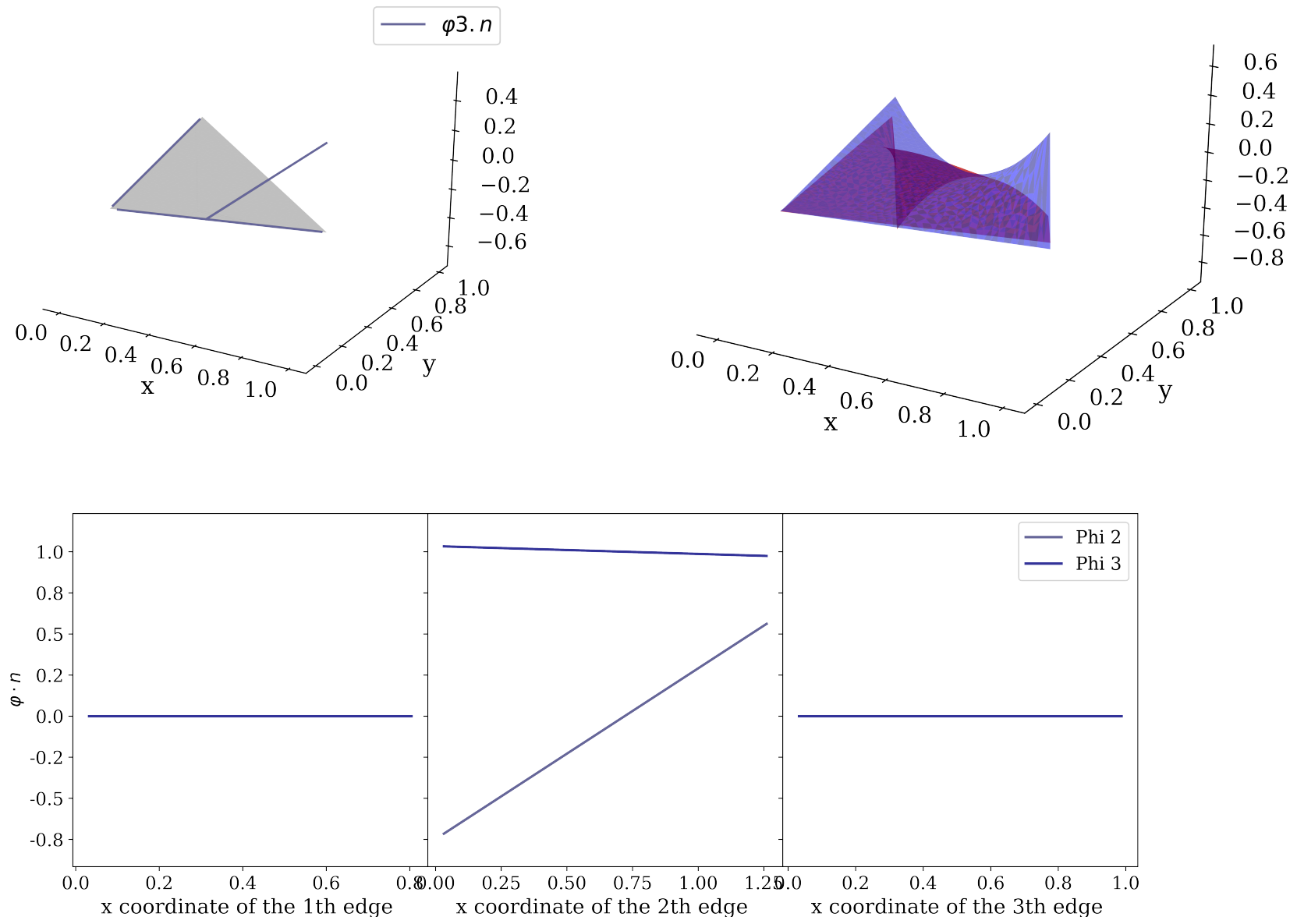

Fig. 135: Second basis function in the case $k = 1$.

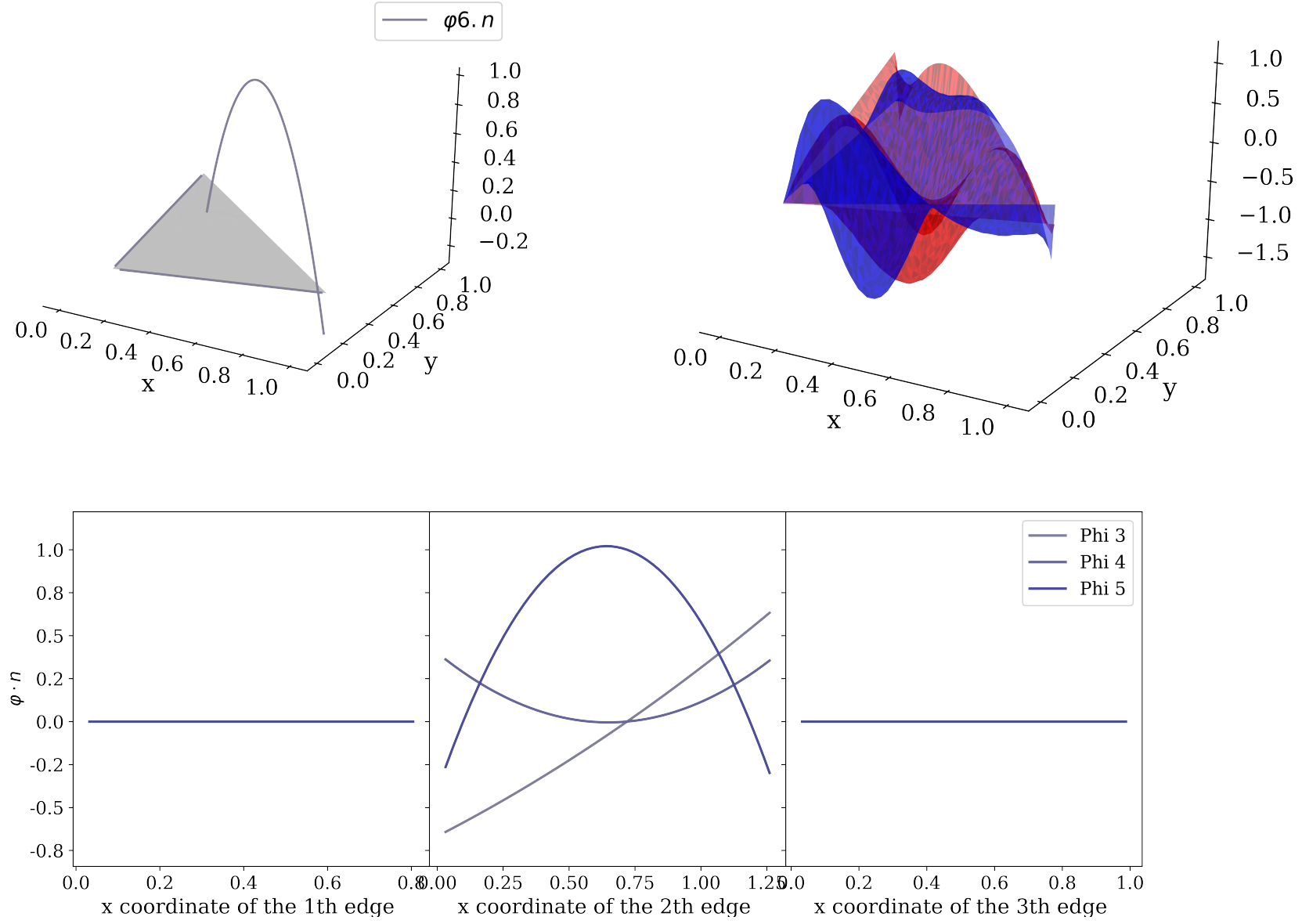

Fig. 136: Fifth basis function in the case $k = 2$. Top left: normal component along the boundaries. Top right: internal behaviour of the basis function. Blue: $x$ component. Red: $y$ component.

Furthermore, as shown on the *Figure 137* the amplitude of the retrieved basis functions is comparable for any edge.

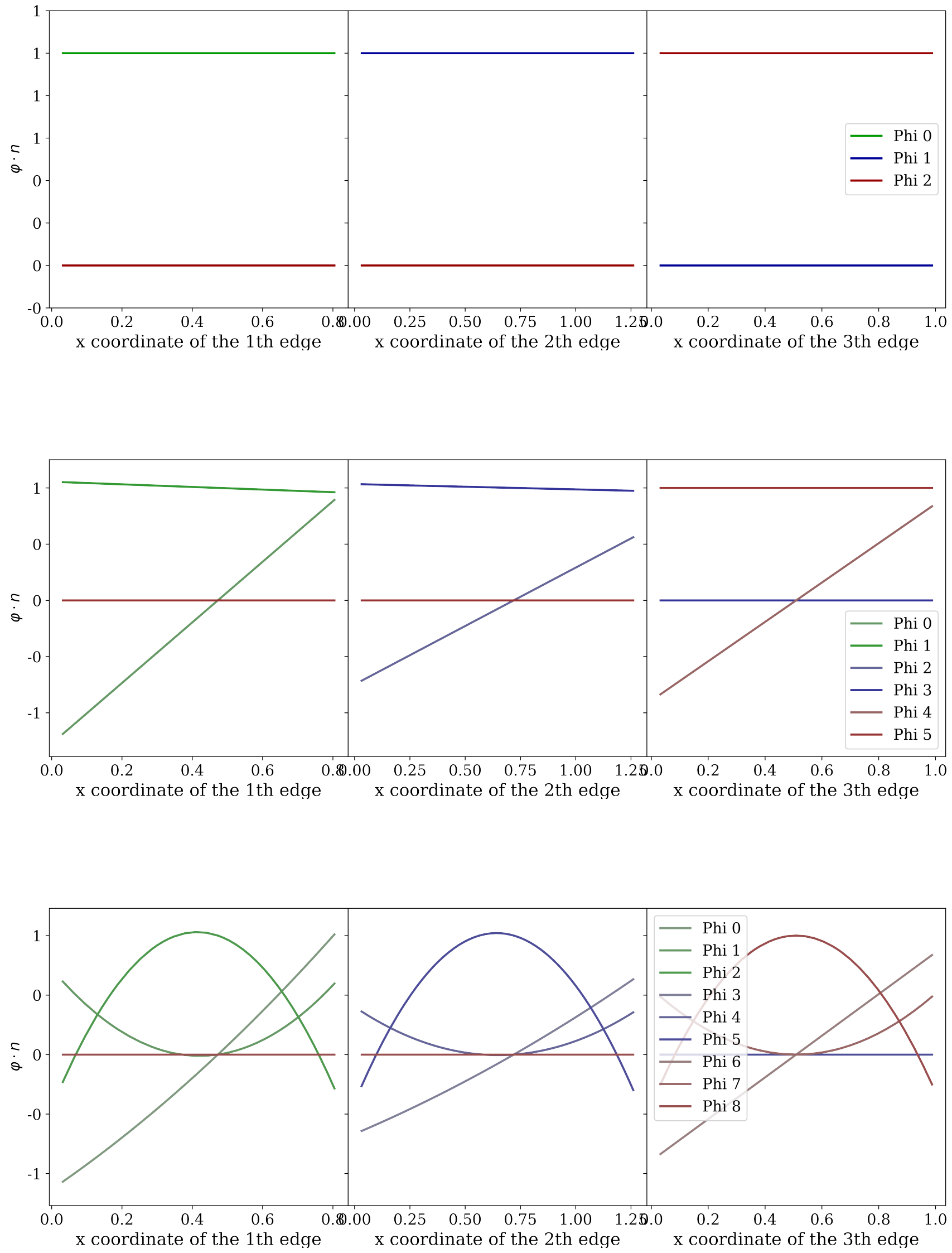

Fig. 137: All the basis functions plotted on all the edges. From top to bottom: $k = 0$, $k = 1$, $k = 2$.

On the side of internal functions, as expected, for any order $k$ bigger or equal than one, the normal components of the internal moments are vanishing on the edges. The parallel with the Raviart – Thomas setting is here also immediate.

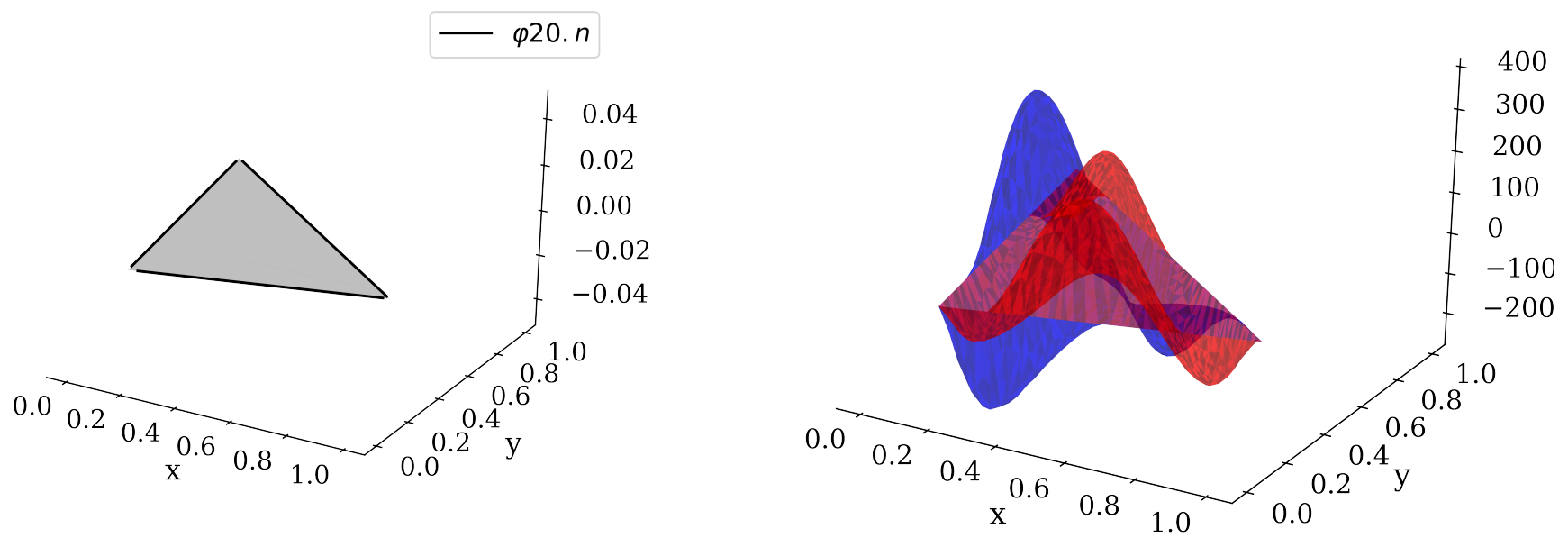

Fig. 138: Last internal basis function in the case $k = 2$.

Lastly, let us derive a comment on the conditionings. As we work on a reduced space where no misc basis functions nor misc moments are considered, the conditioning is further reduced compared to the general setting. In particular, using the Hermite polynomials as projectors we retrieve the following truncated values.

| Order | 1 | 2 | 3 |
|---|---|---|---|
| Conditionings | 15 | 12536 | 670252578. |

Tab. 14: Conditioning values depending on the order

### 7.4.2 Quadrangular shape

We perform the same test for a quadrangular shape. Here too, we considered the reference shape of the quadrilateral Raviart – Thomas elements and rotated it not to fall in our shape limitations. The used element is depicted in the *Figure 139*.

As before, for the sake of concision only one representative of the normal basis functions whose normal component is not vanishing on the boundary will be represented.

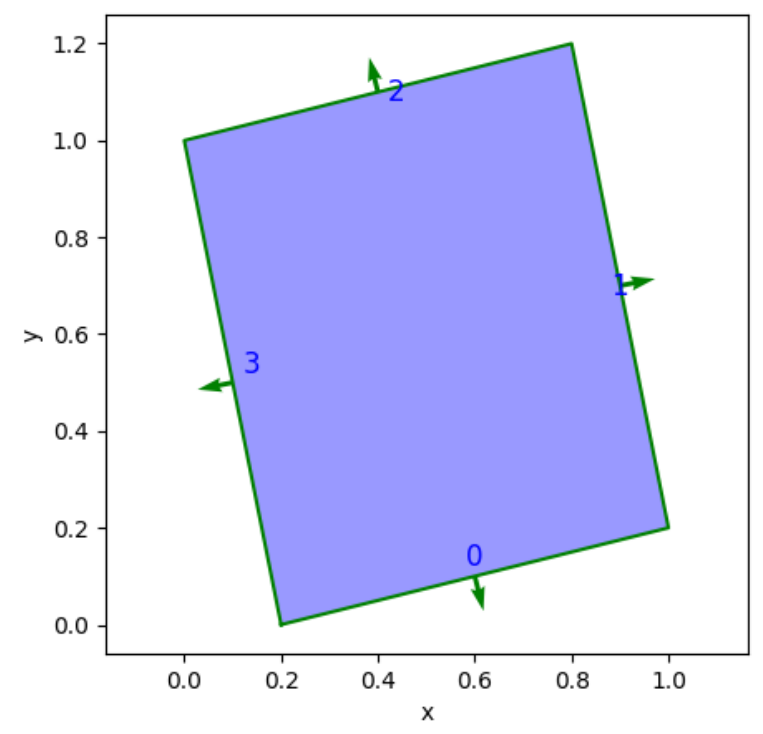

Fig. 139: Quadrangle of reference

For the lowest order space the setting of the element $IIb$ leads to the behaviour represented in the *Figure 140*.

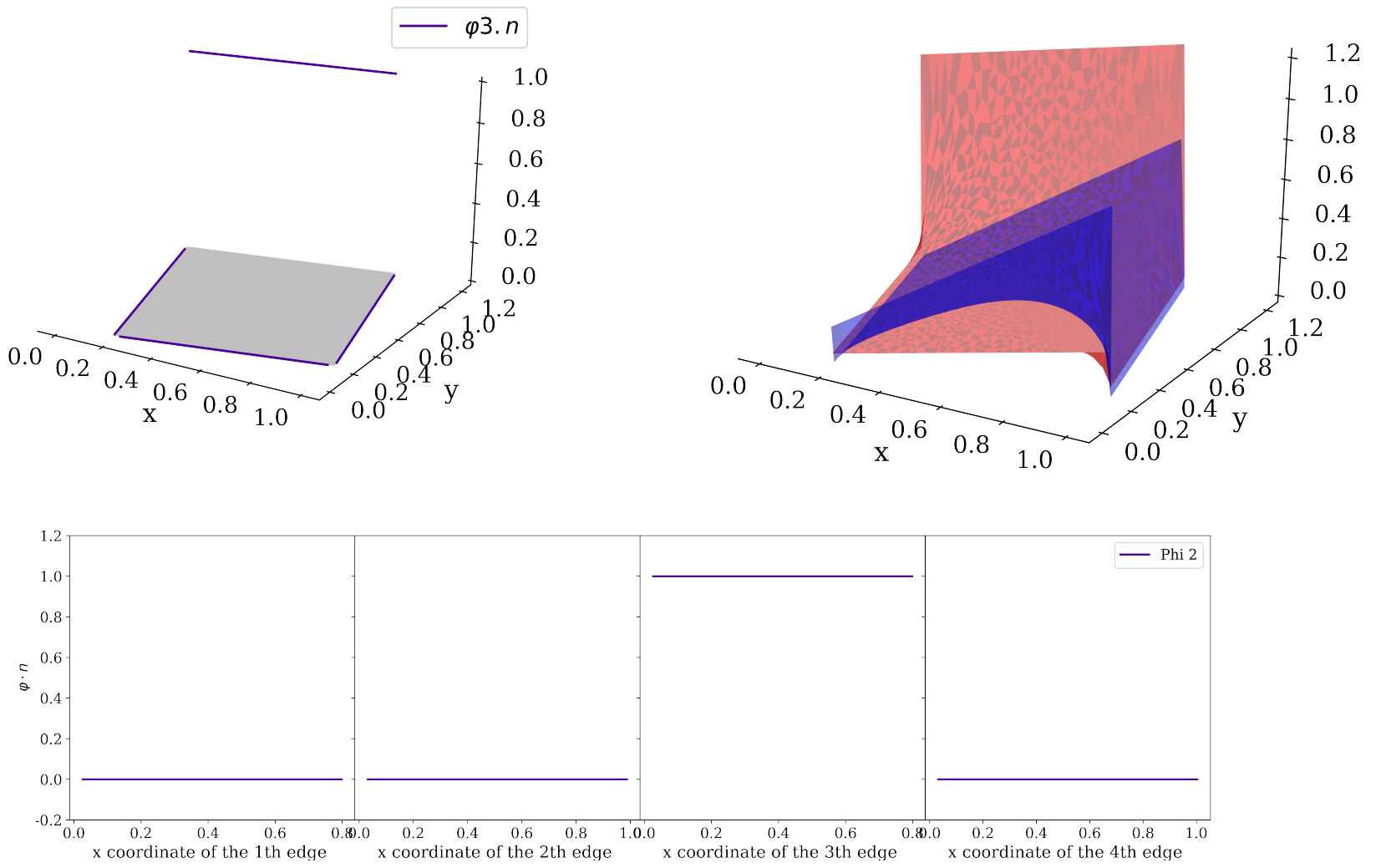

Fig. 140: Second basis function in the case $k = 0$. Top left: normal component along the boundaries. Top right: internal behaviour of the basis function. Blue: $x$ component. Red: $y$ component.

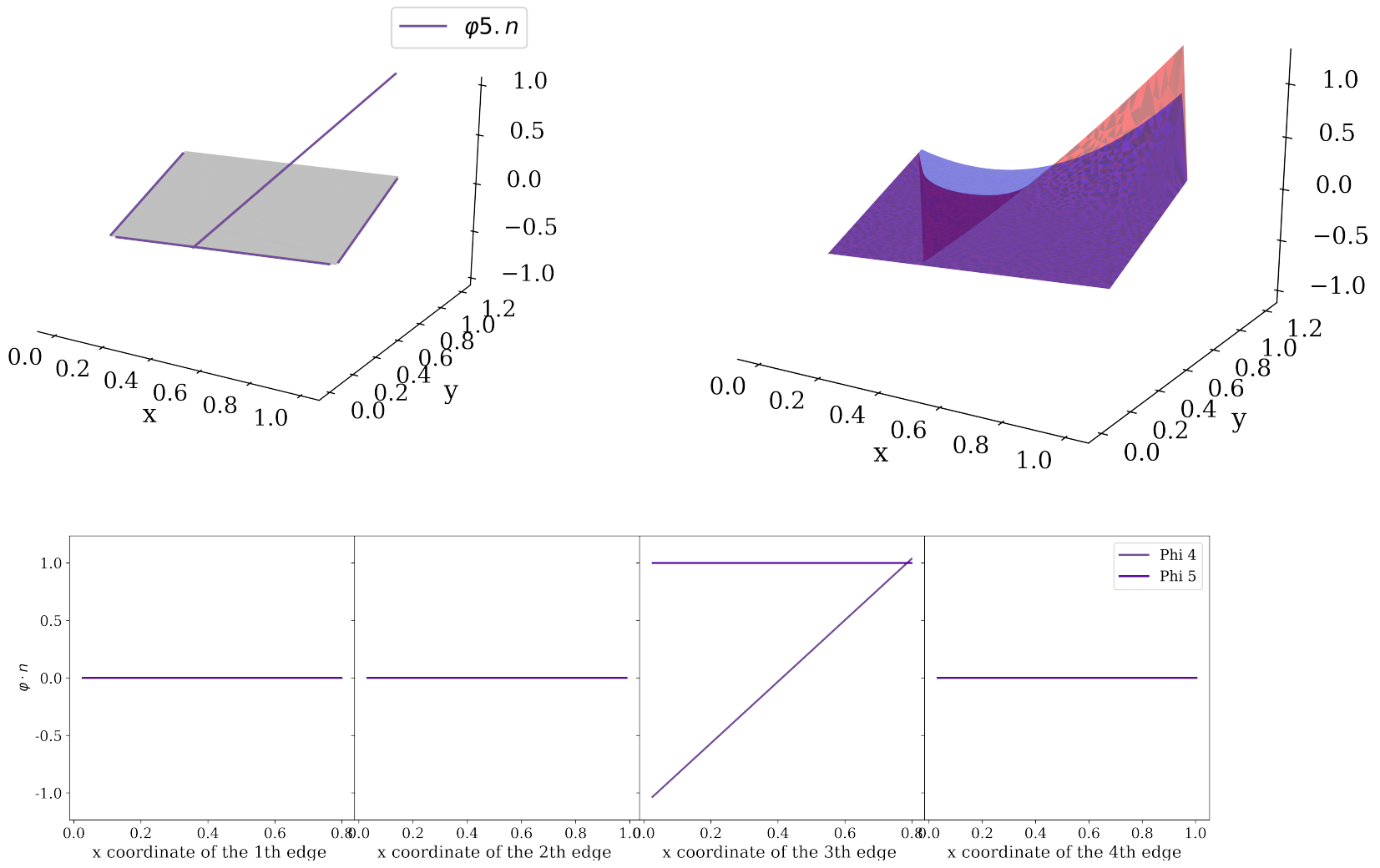

Fig. 141: Fourth basis function in the case $k = 1$. Top left: normal component along the boundaries. Top right: internal behaviour of the basis function. Blue: $x$ component. Red: $y$ component.

As in the simplicial case, one gets immediately the parallel with the classical Raviart – Thomas setting on the boundary. In a similar way, we can derive the same conclusions for the first order element given in the *Figure 141*.

The case of the second order element presented in the *Figure 142* also leads to a parallel with the Raviart – Thomas element. In particular, by the choice of the Hermite polynomials as projectors, the amplitude of the non - vanishing normal components are still contained in a range allowing to project any function living on $\mathbb{H}_k(K)|_{\partial K}$ with reliability.

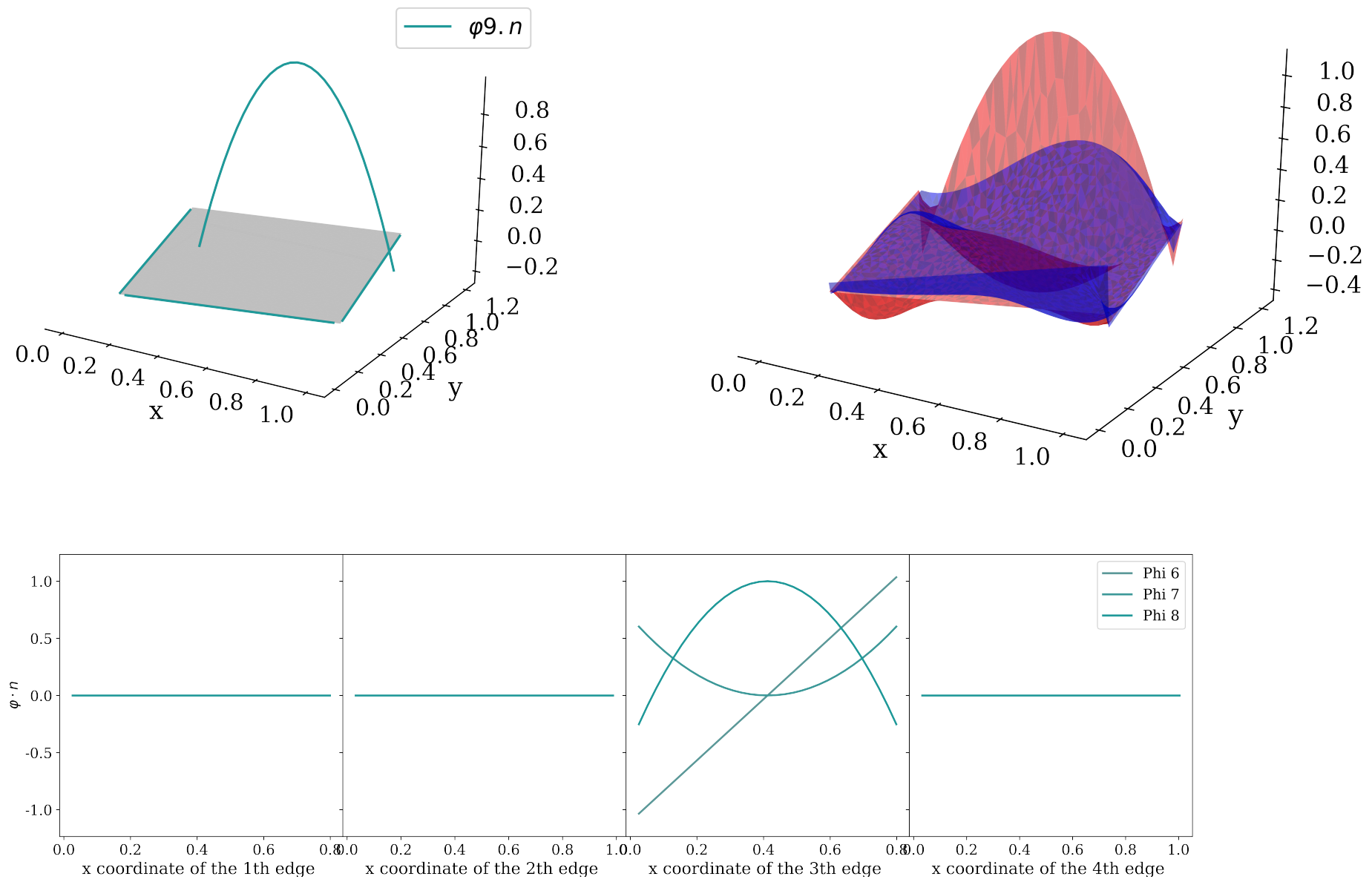

Fig. 142: Eight basis function in the case $k = 2$. Top left: normal component along the boundaries. Top right: internal behaviour of the basis function. Blue: $x$ component. Red: $y$ component.

Furthermore, as shown on the *Figure 143*, the amplitude of the retrieved basis functions is comparable for any edge. More importantly, we can witness that as each function see its support restricted to one single edge, the discretisation of functions living in $\mathbb{H}_k(K)|_{\partial K}$ is similar for every edge.

On the side of the internal basis functions, one can also derive the same results on as in the simplicial case. More precisely, for any order $k$ bigger or equal than one, the normal components of the internal moments are vanishing on the edges. The parallel with the Raviart – Thomas setting is here again immediate.

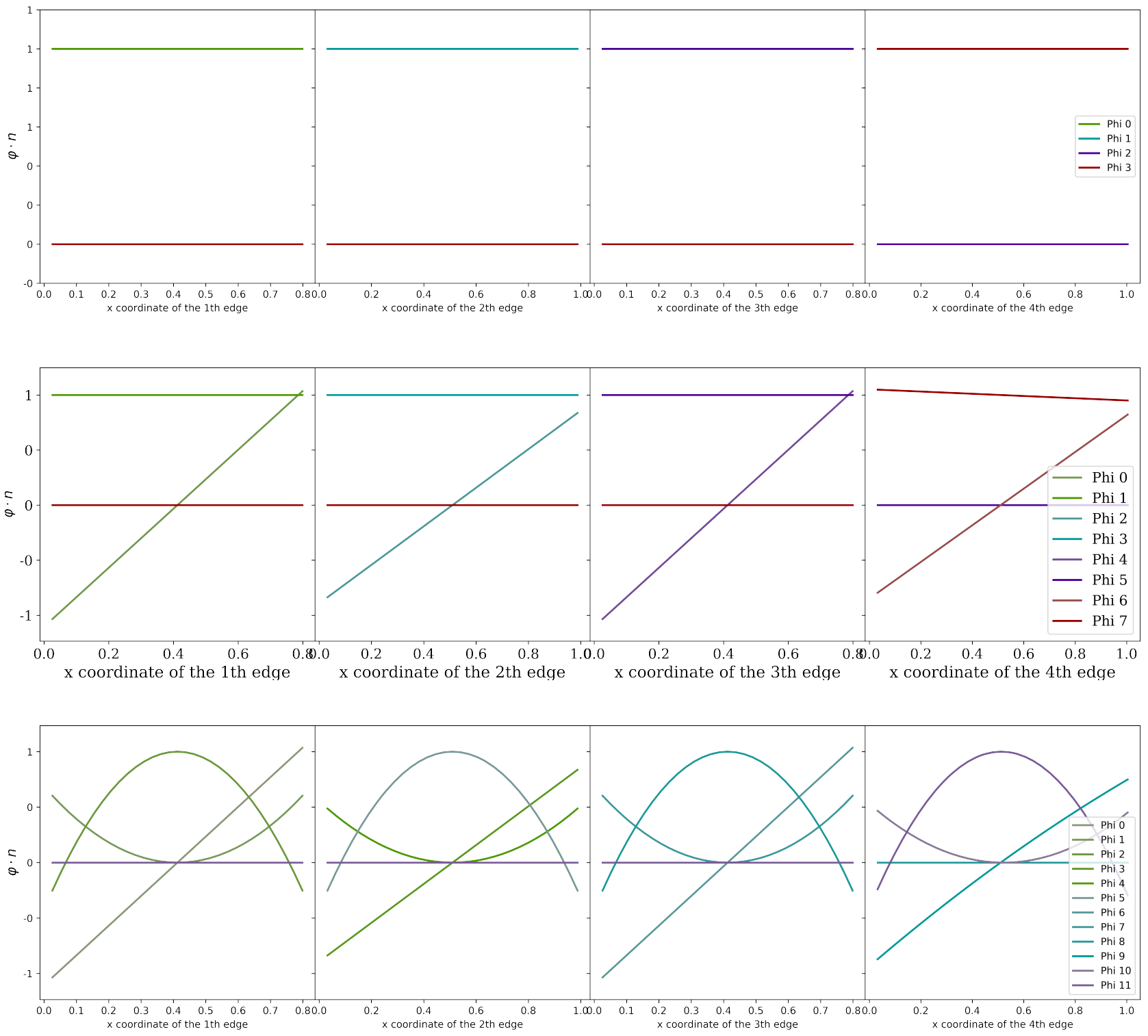

Fig. 143: All the basis functions plotted on all the edges. From top to bottom: $k = 0$, $k = 1$, $k = 2$.

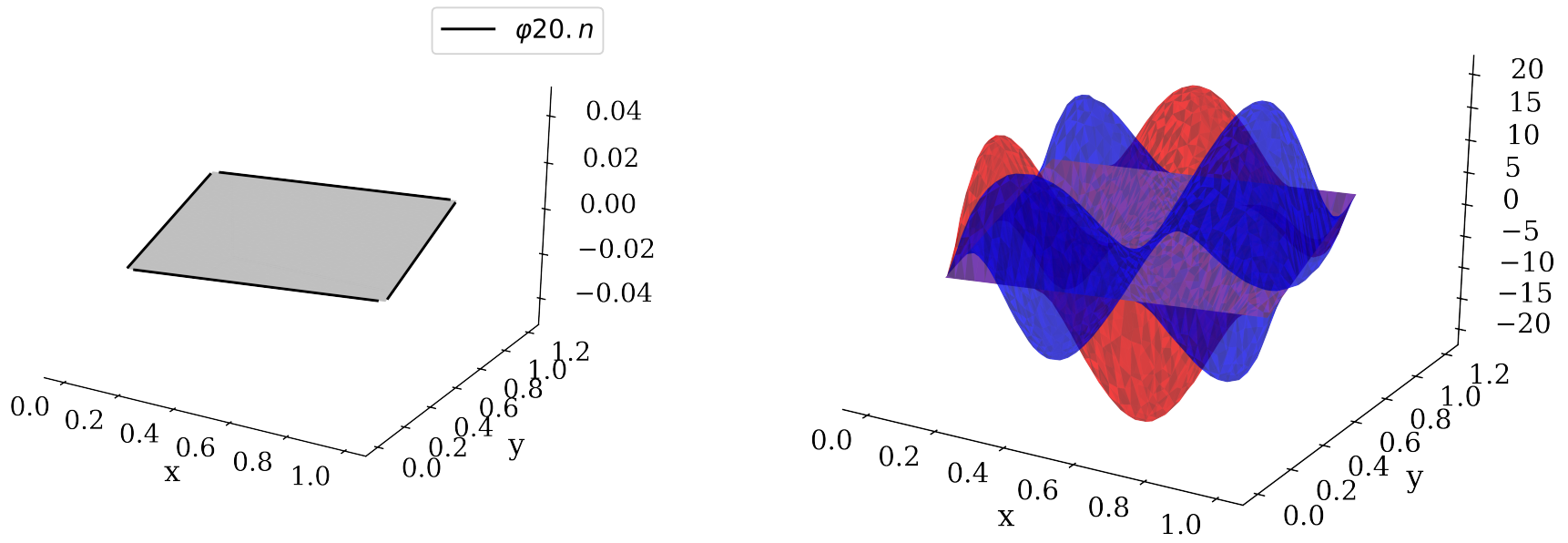

Fig. 144: Internal study

Lastly, let us comment the obtained conditionings in the quadrilateral case. As we work on a reduced space where no misc basis functions nor misc moments are considered, the conditioning is by far reduced compared to the general setting. In particular, using the Hermite polynomials as projectors we retrieve the following truncated values. One can observe that though both elements mainly enjoy the same orientation, the regularity of the quadrilateral shape leads to better conditionings than in the triangular case.

| Order | 1 | 2 | 3 |
|---|---|---|---|
| Conditionings | 14 | 2276 | 15108654. |

Tab. 15: Conditioning values depending on the order

### 7.4.3 Hexagonal shape

Let us end the discussion on the numerical results by constructing the alternative element *IIb* on a non - convex hexagon given in the *Figure 145*. By performing the same tests as for the classical polygons that are used to build the Raviart – Thomas elements on, we show that we achieved an extension of the Raviart – Thomas element on any polygonal shape.

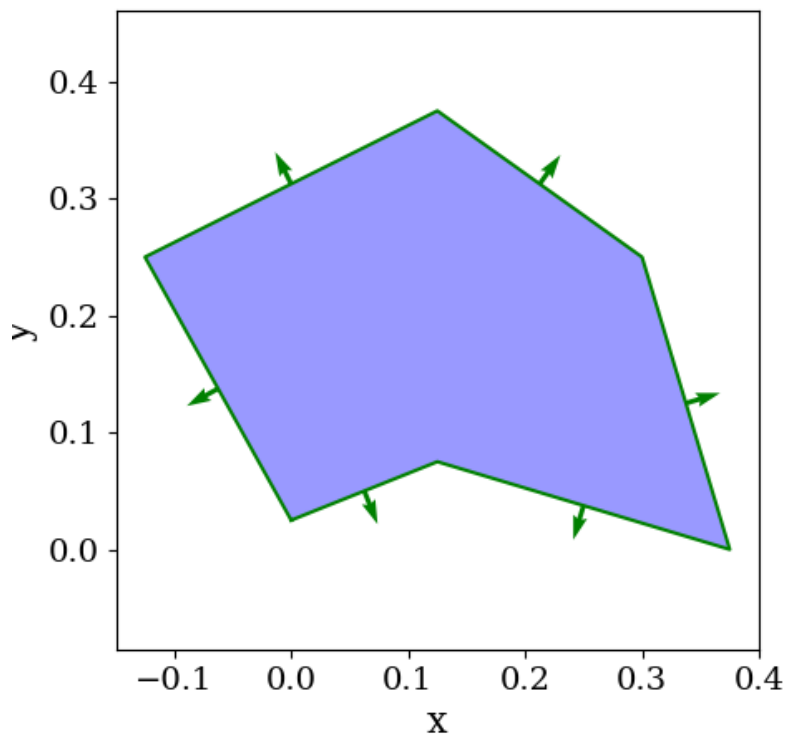

Fig. 145: Quadrangle of reference

As before, for the sake of concision only one representative of the normal basis functions whose normal component is not vanishing on the boundary will be represented.

For the lowest order space, the setting of the element *IIb* leads to the results presented in the *Figure 146*. As in the simplicial and quadrilateral cases, one immediately observe a parallel with the classical Raviart – Thomas setting on the boundary. In a similar way, we can derive the same conclusions for the first order element whose outcome is pictured in the *Figure 147*.

The case of the second order element presented in the *Figure 148* also leads to a parallel with the Raviart – Thomas element. Especially, one can notice that the amplitude of the basis functions are contained in a range close to $[0, 1]$, which eases the decomposition of functions living in $\mathbb{H}_k(K)|_{\partial K}$ over the basis functions. It also emphasises the reliability of the dual set of degrees of freedom in the characterisation of those functions.

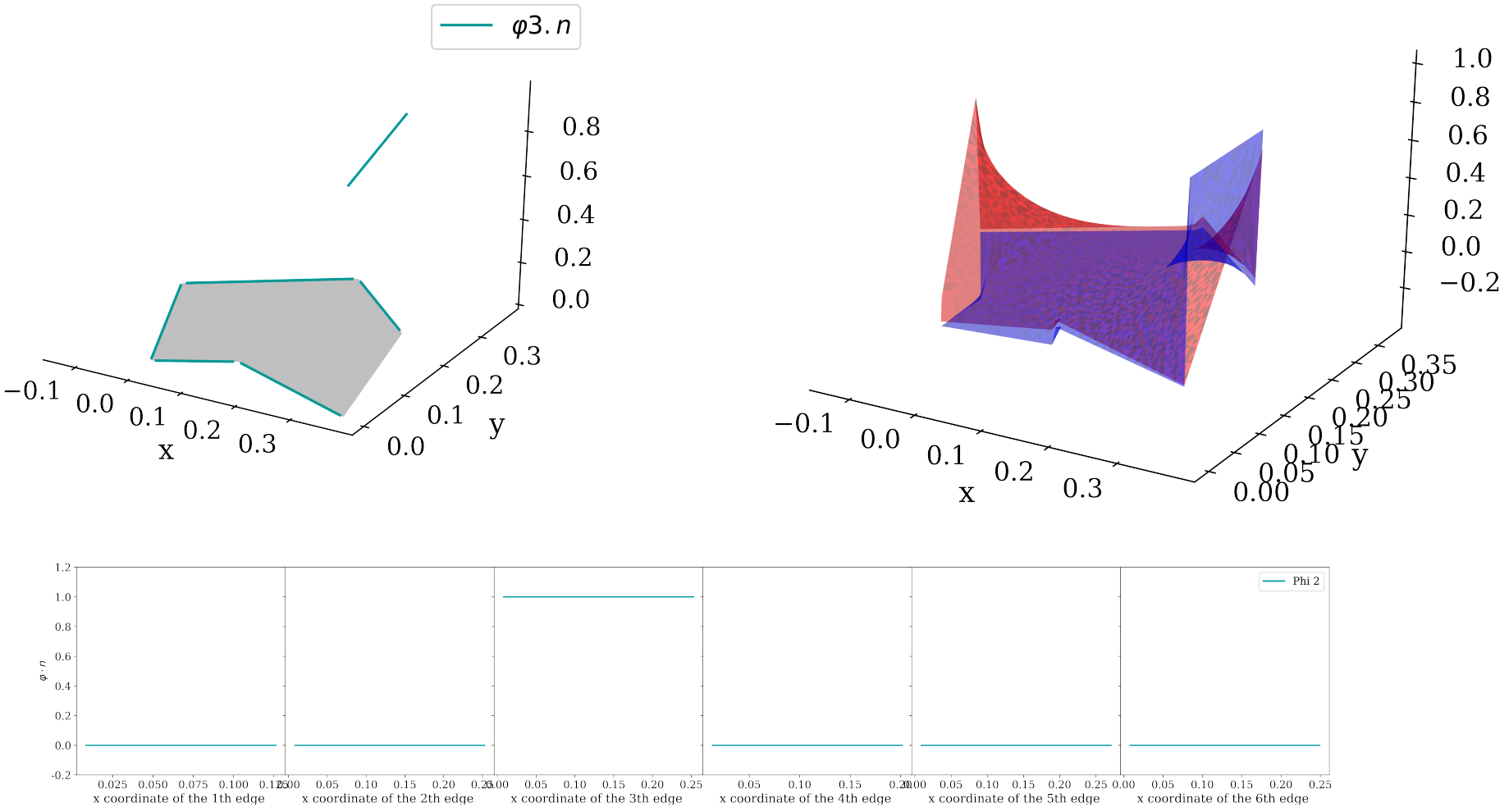

Fig. 146: Second basis function in the case $k = 0$. Top left: normal component along the boundaries. Top right: internal behaviour of the basis function. Blue: $x$ component. Red: $y$ component.

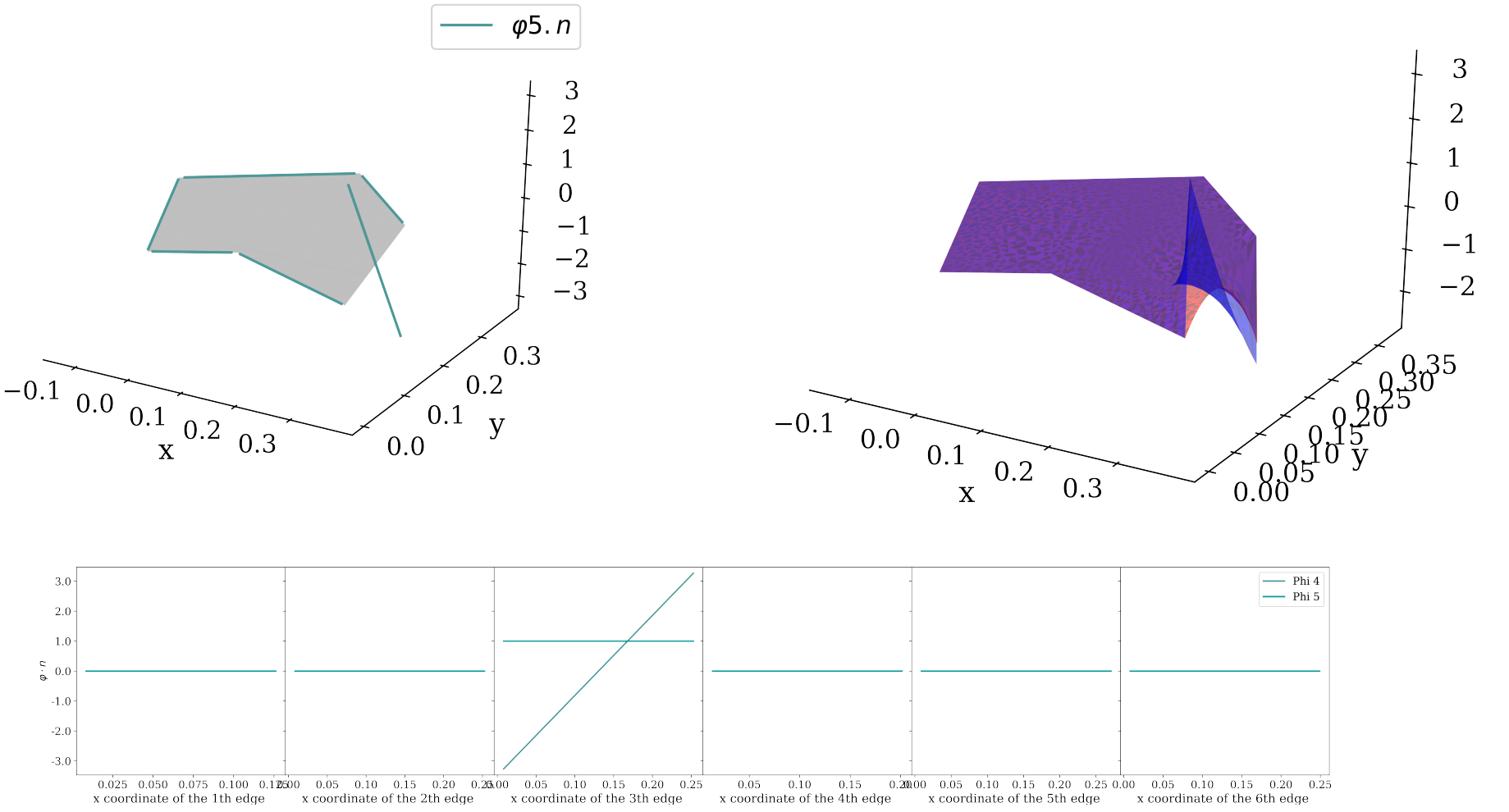

Fig. 147: Fourth basis function in the case $k = 1$. Top left: normal component along the boundaries. Top right: internal behaviour of the basis function. Blue: $x$ component. Red: $y$ component.

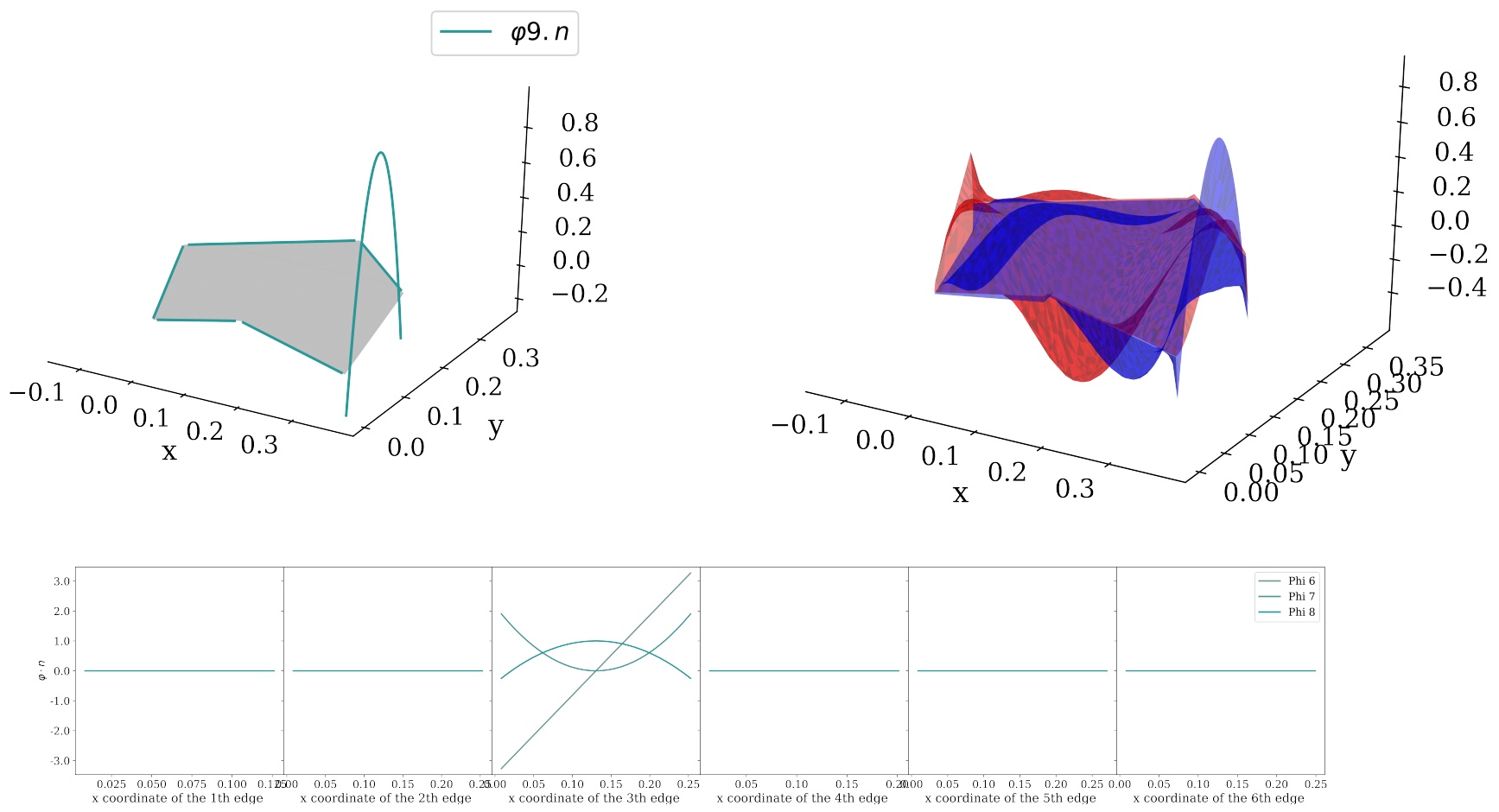

Fig. 148: Eight basis function in the case $k = 2$. Top left: normal component along the boundaries. Top right: internal behaviour of the basis function. Blue: $x$ component. Red: $y$ component.

Furthermore, as shown on the *Figure 149*, the amplitude of the retrieved basis functions is comparable for any edge, as in the simplicial and quadrangular cases. And here again, we can witness that as each function see its support restricted to one single edge, the discretisation of functions living in $\mathbb{H}_k(K)|_{\partial K}$ is similar for every edge.

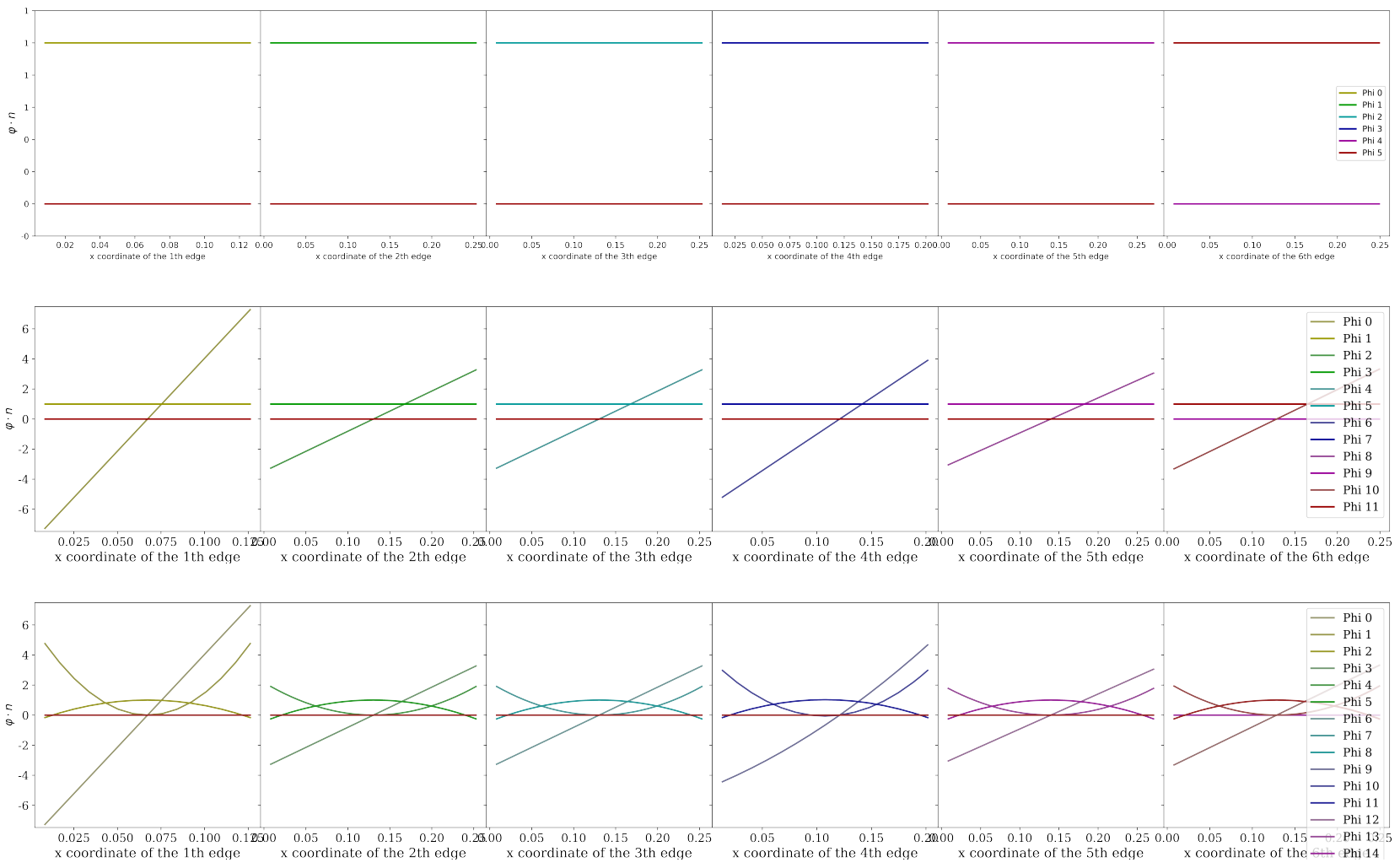

Fig. 149: All the basis functions plotted on all the edges. From top to bottom: $k = 0$, $k = 1$, $k = 2$.

One can also derive same results as before for the internal functions. Similarly to the simplicial and quadrangular cases their normal components are vanishing on the boundary (see the *Figure 150*).

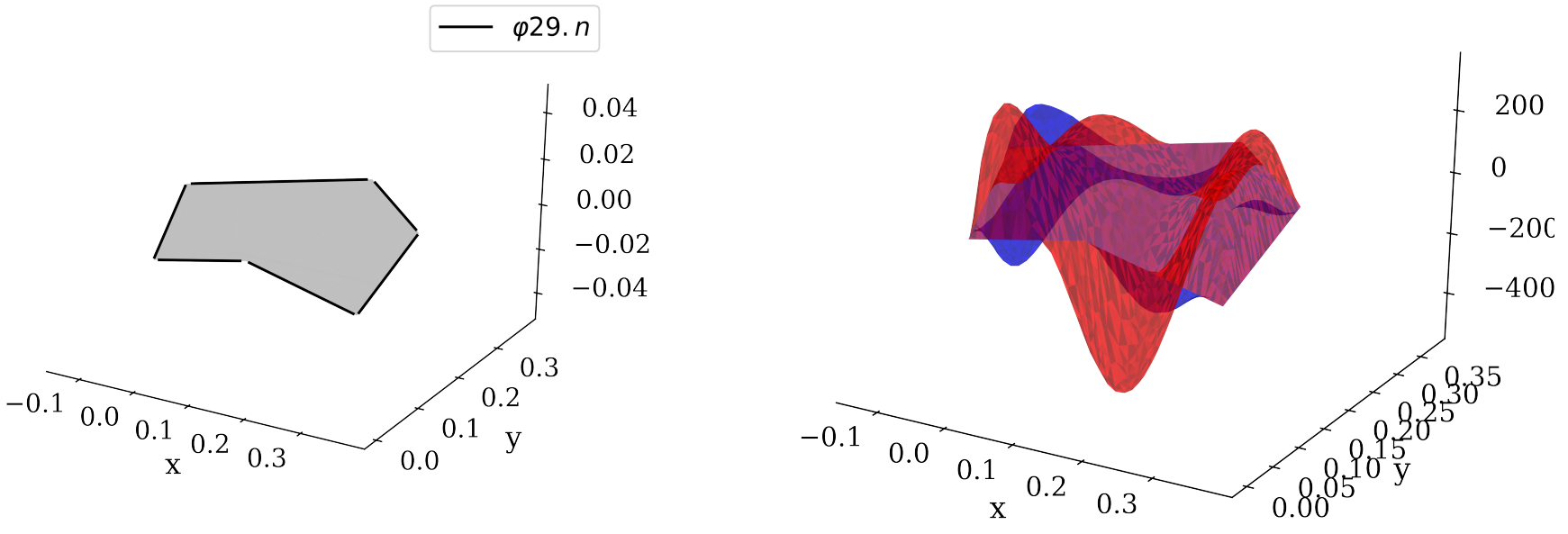

Fig. 150: Internal study

Lastly, we can observe that the conditioning is increased compared to the two classical cases. Indeed, connecting with the discussion on the conditioning for the classical space, combining the squezeness of the element with an increased number of edges let us expect an increase of the conditioning. However, please note that we still retrieve lower values with the reduced space than with the general one.

| Order | 1 | 2 | 3 |
|---|---|---|---|
| Conditionings | 113 | 479311 | 939600583. |

Tab. 16: Conditioning values depending on the order

### 7.4.4 Partial conclusion

The general space $\mathbb{H}_k(K)$ allowed us to define $H(\text{div})$ – conformal elements on any shape. Especially, any of the considered four configurations $Ia$, $Ib$, $IIa$ and $IIb$ offers reliable settings on any polygonal shape fulfilling the two restrictions we derived in the *Section 6.1.2*.

In particular, the retrieved tuned basis functions enhanced the properties instilled to the discretisations, depending on the definition of the corresponding degrees of freedom. Indeed, for the configuration $I$, where the characterisation of the boundary part of functions living in $\mathbb{H}_k(K)$ is wished to be the more coordinate wise as possible, only one normal canonical basis function is degenerating into an internal one. Thus, the $d = 2$ constants of the functions are tuned coordinate - wise. In the second configuration where the tuning is wished to be more global, two basis functions are degenerating, meaning that even at the level of constants the tuning is done globally. The determination of coordinate - wise coefficients is done through higher order projections.

Similarly, we could test the scaling of the basis functions depending on the type of the element. The type $a$, which is entirely defined from moments did not prescribe any scale, though the type $b$ forced the normal component of the constant basis functions to scale to one, while keeping the highest order normal components of basis functions within a reasonable amplitude.

Lastly, it was shown that the reduced space furnishes a more reliable element through a lower conditioning of its transfer matrices. Its dimension matching directly the target order on the boundary, no degeneration of canonical basis functions is witnessed and a parallel with the classical Raviart – Thomas setting can be derived on the boundary of any polygonal element from the lowest order. Furthermore, though not polynomials, the vectorial basis functions still preserve a smooth behaviour within the polygon.

## Acknowledgements.

E. Le Mélédo and P. Öffner have been funded by SNF project 200020_175784 "Solving advection dominated problems with high order schemes with polygonal meshes: application to compressible and incompressible flow problems".

# A Implemented spaces

In all the elements presented below, the functions $g$, $h$, $p$, $q$ have to be chosen given the list of implemented functions provided in the last subsection.

## A.1 Classical space

This space corresponds to the one presented in the *Section 6*, implemented in the code: "All_Combined".

**Space:**

$$\mathbb{H}_k(K) = \{u \in H^1(K),\, u|_{\partial K} \in \mathcal{H}_0(\partial K),\, \Delta u \in \mathbb{Q}_{k-1}(K)\}^2 + \begin{pmatrix} x \\ y \end{pmatrix} \{u \in H^1(K),\, u|_{\partial K} \in \mathcal{H}_k(\partial K),\, \Delta u \in \mathbb{Q}_{[k-1]}(K)\}$$

**Dimension:**

$$\dim \mathbb{H}_k(K) = n(k+3) + 2k(k-1) - \mathbb{1}_{k>0}$$

### A.1.1 Implemented canonical basis functions

**Case** $k = 0$ (Only normal functions, for any edge $f_i \in \partial K$. One to two of them will degenerate during the tuning, for any edge $f_i \in \partial K$)**:**

$$f_1 := (x,\, y) \mapsto \begin{pmatrix} x \\ y \end{pmatrix} \begin{Bmatrix} \Delta u = 0 \\ u|_{\partial K} = \mathbb{1}_{f_i} \end{Bmatrix} - \begin{pmatrix} n_{ix} + n_{iy}^2/n_{ix} \\ 0 \end{pmatrix} \begin{Bmatrix} \Delta u = 0 \\ u|_{\partial K} = 2\mathbb{1}_{f_i} \end{Bmatrix}$$

$$f_2 := (x,\, y) \mapsto \begin{pmatrix} x \\ y \end{pmatrix} \begin{Bmatrix} \Delta u = 0 \\ u|_{\partial K} = \mathbb{1}_{f_i} \end{Bmatrix} - \begin{pmatrix} 0 \\ n_{iy} + n_{ix}^2/n_{iy} \end{pmatrix} \begin{Bmatrix} \Delta u = 0 \\ u|_{\partial K} = 2\mathbb{1}_{f_i} \end{Bmatrix}$$

$$f_3 := (x,\, y) \mapsto \begin{pmatrix} x \\ y \end{pmatrix} \begin{Bmatrix} \Delta u = 0 \\ u|_{\partial K} = \mathbb{1}_{f_i} \end{Bmatrix} + \begin{pmatrix} n_{ix} \\ n_{iy} \end{pmatrix} \begin{Bmatrix} \Delta u = 0 \\ u|_{\partial K} = 2\mathbb{1}_{f_i} \end{Bmatrix}$$

**General case** $k > 0$**:**

*Normal functions* (One to two of them will degenerate during the tuning, for any edge $f_i \in \partial K$):

$$f_1 := (x,\, y) \mapsto \begin{pmatrix} x \\ y \end{pmatrix} \begin{Bmatrix} \Delta u = h_{k-1,\,k-1}(x,y) \\ u|_{\partial K} = g_0(x,y)\mathbb{1}_{f_i} \end{Bmatrix} - \begin{pmatrix} n_{ix} + n_{iy}^2/n_{ix} \\ 0 \end{pmatrix} \begin{Bmatrix} \Delta u = 0 \\ u|_{\partial K} = 2\mathbb{1}_{f_i} \end{Bmatrix}$$

$$f_2 := (x,\, y) \mapsto \begin{pmatrix} x \\ y \end{pmatrix} \begin{Bmatrix} \Delta u = h_{k-1,\,k-1}(x,y) \\ u|_{\partial K} = g_{k-1}(x,y)\mathbb{1}_{f_i} \end{Bmatrix} - \begin{pmatrix} 0 \\ n_{iy} + n_{ix}^2/n_{iy} \end{pmatrix} \begin{Bmatrix} \Delta u = 0 \\ u|_{\partial K} = 2\mathbb{1}_{f_i} \end{Bmatrix}$$

$$f_{2+j} := (x,\, y) \mapsto \begin{pmatrix} x \\ y \end{pmatrix} \left\{ \begin{aligned} &\Delta u = h_{k-1,\,k-1}(x,y) \\ &u|_{\partial K} = g_j(x,y)\mathbb{1}_{f_i} \end{aligned} \right\} + \begin{pmatrix} n_{ix} \\ n_{iy} \end{pmatrix} \left\{ \begin{aligned} &\Delta u = 0 \\ &u|_{\partial K} = 2\mathbb{1}_{f_i} \end{aligned} \right\},\ j \in [\![0,\, k]\!]$$

*Internal functions:*

$$f_{xl} := (x,\, y) \mapsto \begin{pmatrix} x \\ y \end{pmatrix} \left\{ \begin{aligned} &\Delta u = h_{k-1,\,l}(x,\, y) \\ &u|_{\partial K} = 0 \end{aligned} \right\},\ 0 \le l \le k-1$$

$$f_{yl} := (x,\, y) \mapsto \begin{pmatrix} x \\ y \end{pmatrix} \left\{ \begin{aligned} &\Delta u = h_{l,\,k-1}(x,\, y) \\ &u|_{\partial K} = 0 \end{aligned} \right\},\ 0 \le l < k-1$$

$$f_{xlm} := (x,\, y) \mapsto \begin{pmatrix} 1 \\ 0 \end{pmatrix} \left\{ \begin{aligned} &\Delta u = h_{l,\,m}(x,\, y) \\ &u|_{\partial K} = 0 \end{aligned} \right\},\ 0 \le l,\, m \le k-1$$

$$f_{ylm} := (x,\, y) \mapsto \begin{pmatrix} 0 \\ 1 \end{pmatrix} \left\{ \begin{aligned} &\Delta u = h_{l,\,m}(x,\, y) \\ &u|_{\partial K} = 0 \end{aligned} \right\},\ 0 \le l,\, m \le k-1$$

### A.1.2 Implemented Elements (degrees of freedom)

*Normal moments*, for any edge $f_j \in \partial K$:

| | Core Moments | Misc Moment | Supplementary moments |
|---|---|---|---|
| Ia | $\int_{f_j} p_i(v)\,\varphi\cdot n\,\mathrm{d}\gamma,\ i \in [\![1,\, k]\!]$ | $\int_{f_j} v\,\varphi \cdot \begin{pmatrix} 1 \\ 1 \end{pmatrix}\,\mathrm{d}\gamma$ | $\int_{f_j} \varphi_x n_{ix}\,\mathrm{d}\gamma$ and $\int_{f_j} \varphi_y\, n_{iy}\,\mathrm{d}\gamma$ |
| Ib | $\int_{f_j} p_i(v)\,\varphi\cdot n\,\mathrm{d}\gamma,\ i \in [\![1,\, k]\!]$ | $\int_{f_j} v\,\varphi \cdot \begin{pmatrix} 1 \\ 1 \end{pmatrix}\,\mathrm{d}\gamma$ | $\varphi_x(x_m,\, y_m)n_{ix}$ and $\varphi_y(x_m,\, y_m)n_{iy}$ |
| IbShifted | $\int_{f_j} p_i(v)\,\varphi\cdot n\,\mathrm{d}\gamma,\ i \in [\![1,\, k]\!]$ | $\int_{f_j} v\,\varphi \cdot \begin{pmatrix} 1 \\ 1 \end{pmatrix}\,\mathrm{d}\gamma$ | $\varphi_x(x_m,\, y_m)n_{ix} - 1$ and $\varphi_y(x_m,\, y_m)n_{iy} - 1$ |
| IIa | $\int_{f_j} p_i(v)\,\varphi\cdot n\,\mathrm{d}\gamma,\ i \in [\![1,\, k]\!]$ | $\int_{f_j} \varphi\cdot n\,\mathrm{d}\gamma$ | $\int_{f_j} x\,\varphi_x n_{ix}\,\mathrm{d}\gamma$ and $\int_{f_j} y\,\varphi_y\, n_{iy}\,\mathrm{d}\gamma$ |
| IIb | $\int_{f_j} p_i(v)\,\varphi\cdot n\,\mathrm{d}\gamma,\ i \in [\![1,\, k]\!]$ | $\varphi(x_m)\cdot n$ | $\int_{f_j} x\,\varphi_x n_{ix}\,\mathrm{d}\gamma$ and $\int_{f_j} y\,\varphi_y\, n_{iy}\,\mathrm{d}\gamma$ |
| IIbShifted | $\int_{f_j} p_i(v)\,\varphi\cdot n\,\mathrm{d}\gamma,\ i \in [\![1,\, k]\!]$ | $\varphi(x_m)\cdot n - 1$ | $\int_{f_j} x\,\varphi_x n_{ix}\,\mathrm{d}\gamma$ and $\int_{f_j} y\,\varphi_y\, n_{iy}\,\mathrm{d}\gamma$ |

*Internal moments*:

| | Core Moments | Misc Moment |
|---|---|---|
| All elements | $\int_K \begin{pmatrix} q_{l,\,m}(x,\, y) \\ 0 \end{pmatrix} \cdot \varphi\,\mathrm{d}\gamma$ and $\int_K \begin{pmatrix} 0 \\ q_{m,\,l}(x,\, y) \end{pmatrix} \cdot \varphi\,\mathrm{d}\gamma,\ \substack{l\in[\![0,\,k]\!] \\ m\in[\![0,\,k-1]\!] \\ (l,\,m)\neq(k,\,k-1)}$ | $\int_K \begin{pmatrix} q_{k,\,k-1}(x,\, y) \\ q_{k-1,\,k}(x,\, y) \end{pmatrix} \cdot \varphi\,\mathrm{d}\gamma$ |

## A.2 Reduced space

This space corresponds to the one discussed in the *Section 6.2*, but where the boundary conditions of functions are left on $\mathbb{P}_k(f)$ instead of prescribing them to the constant. This is not natural but do not make necessary the tuning of the basis functions of the space $\mathbb{H}_k(K)$ towards the designed elements for the basis functions to have the property $p \cdot n \in \mathbb{P}_k(f)$. The corresponding code is named " All_Combined_Alternative_Space".

***Note.*** If one wants nevertheless to retrieve the basis functions corresponding to the presented elements, this space is not recommended as the definition boundary conditions are not natural for the space. Furthermore, after tuning it reduces to the All_Combined_Alternative_Space_Uniform below, that is constructed from a natural definition. ▲

**Space:**

$$\begin{aligned}\mathbb{H}_k(K) = &\big\{u \in H^1(K),\, u|_{\partial K} \equiv 1,\, \Delta u \in \mathbb{Q}_{k-1}(K)\big\}^2 \\ &+ \begin{pmatrix} x \\ y \end{pmatrix} \big\{u \in H^1(K),\, u|_{\partial K} \in \mathcal{H}_k(\partial K),\, \Delta u \in \mathbb{Q}_{[k-1]}(K)\big\}\end{aligned}$$

**Dimension:**

$$\dim \mathbb{H}_k(K) = n(k+1) + 2k(k-1) - \mathbb{1}_{k>0}$$

### A.2.1 Implemented canonical basis functions

**Case** $k = 0$ (only normal functions, for any edge $f_i \in \partial K$)**:**

$$f_1 := (x,\, y) \mapsto \begin{pmatrix} x \\ y \end{pmatrix} \begin{Bmatrix} \Delta u = 0 \\ u|_{\partial K} = \mathbb{1}_{f_i} \end{Bmatrix} + \begin{pmatrix} n_{ix} \\ n_{iy} \end{pmatrix} \begin{Bmatrix} \Delta u = 0 \\ u|_{\partial K} = 2\mathbb{1}_{f_i} \end{Bmatrix}$$

**General case** $k > 0$**:**

*Normal functions:*

$$f_j := (x,\, y) \mapsto \begin{pmatrix} x \\ y \end{pmatrix} \begin{Bmatrix} \Delta u = h_{k-1,\, k-1}(x, y) \\ u|_{\partial K} = g_j(x, y)\mathbb{1}_{f_i} \end{Bmatrix} + \begin{pmatrix} n_{ix} \\ n_{iy} \end{pmatrix} \begin{Bmatrix} \Delta u = 0 \\ u|_{\partial K} = 2\mathbb{1}_{f_i} \end{Bmatrix},\ j \in [\![0,\, k]\!]$$

*Internal functions:*

$$f_{xl} := (x,\, y) \mapsto \begin{pmatrix} x \\ y \end{pmatrix} \begin{Bmatrix} \Delta u = h_{k-1,\, l}(x,\, y) \\ u|_{\partial K} = 0 \end{Bmatrix},\ 0 \le l \le k-1$$

$$f_{yl} := (x,\, y) \mapsto \begin{pmatrix} x \\ y \end{pmatrix} \left\{ \begin{aligned} &\Delta u = h_{l,\,k-1}(x,\, y) \\ &u|_{\partial K} = 0 \end{aligned} \right\},\ 0 \leq l < k-1$$

$$f_{xlm} := (x,\, y) \mapsto \begin{pmatrix} 1 \\ 0 \end{pmatrix} \left\{ \begin{aligned} &\Delta u = h_{l,\,m}(x,\, y) \\ &u|_{\partial K} = 0 \end{aligned} \right\},\ 0 \leq l,\, m \leq k-1$$

$$f_{ylm} := (x,\, y) \mapsto \begin{pmatrix} 0 \\ 1 \end{pmatrix} \left\{ \begin{aligned} &\Delta u = h_{l,\,m}(x,\, y) \\ &u|_{\partial K} = 0 \end{aligned} \right\},\ 0 \leq l,\, m \leq k-1$$

### A.2.2 Implemented Elements (degrees of freedom)

Normal moments, for any edge $f_j \in \partial K$:

| | Core Moments | Misc Moment | Supplementary moments |
|---|---|---|---|
| Ia | $\int_{f_j} p_i(v)\,\varphi\cdot n\,\mathrm{d}\gamma,\ i \in [\![1,\,k]\!]$ | $\int_{f_j} v\,\varphi\cdot \begin{pmatrix} 1 \\ 1 \end{pmatrix}\,\mathrm{d}\gamma$ | None |
| Ib | $\int_{f_j} p_i(v)\,\varphi\cdot n\,\mathrm{d}\gamma,\ i \in [\![1,\,k]\!]$ | $\int_{f_j} v\,\varphi\cdot \begin{pmatrix} 1 \\ 1 \end{pmatrix}\,\mathrm{d}\gamma$ | None |
| IbShifted | $\int_{f_j} p_i(v)\,\varphi\cdot n\,\mathrm{d}\gamma,\ i \in [\![1,\,k]\!]$ | $\int_{f_j} v\,\varphi\cdot \begin{pmatrix} 1 \\ 1 \end{pmatrix}\,\mathrm{d}\gamma$ | None |
| IIa | $\int_{f_j} p_i(v)\,\varphi\cdot n\,\mathrm{d}\gamma,\ i \in [\![1,\,k]\!]$ | $\int_{f_j} \varphi\cdot n\,\mathrm{d}\gamma$ | None |
| IIb | $\int_{f_j} p_i(v)\,\varphi\cdot n\,\mathrm{d}\gamma,\ i \in [\![1,\,k]\!]$ | $\varphi(x_m)\cdot n$ | None |
| IIbShifted | $\int_{f_j} p_i(v)\,\varphi\cdot n\,\mathrm{d}\gamma,\ i \in [\![1,\,k]\!]$ | $\varphi(x_m)\cdot n - 1$ | None |

Internal moments:

| | Core Moments | Misc Moment |
|---|---|---|
| All elements | $\int_K \begin{pmatrix} q_{l,m}(x,\, y) \\ 0 \end{pmatrix}\cdot\varphi\,\mathrm{d}x\mathrm{d}y$ and $\int_K \begin{pmatrix} 0 \\ q_{m,l}(x,\, y) \end{pmatrix}\cdot\varphi\,\mathrm{d}x\mathrm{d}y,\ \substack{l\in[\![0,\,k]\!] \\ m\in[\![0,\,k-1]\!] \\ (l,m)\neq(k,\,k-1)}$ | $\int_K \begin{pmatrix} q_{k,\,k-1}(x,\, y) \\ q_{k-1,\,k}(x,\, y) \end{pmatrix}\cdot\varphi\,\mathrm{d}x\mathrm{d}y$ |

***Note.*** Due to the selected moments to leave out from the general setting (as explained in the *Section 6.2*), the elements Ia, Ib and IbShifted are Identical.
▲

## A.3 Reduced space endowed with natural basis functions

This space corresponds to the one discussed in the *Section 6.2*, and its basis functions are constructed naturally. The corresponding code is named " All_Combined_Alternative_Space_Uniform".

**Space:**

$$\mathbb{H}_k(K) = \left\{u \in H^1(K),\, u|_{\partial K} \equiv 1,\, \Delta u \in \mathbb{Q}_{k-1}(K)\right\}^2 + \begin{pmatrix} x \\ y \end{pmatrix} \left\{u \in H^1(K),\, u|_{\partial K} \in \mathcal{H}_k(\partial K),\, \Delta u \in \mathbb{Q}_{[k-1]}(K)\right\}$$

**Dimension:**

$$\dim \mathbb{H}_k(K) = n(k+1) + 2k(k-1) - \mathbb{1}_{k>0}$$

### A.3.1 Implemented canonical basis functions

**Case** $k = 0$ (only normal functions, for any edge $f_i \in \partial K$)**:**

$$f_1 := (x,\, y) \mapsto \begin{pmatrix} x \\ y \end{pmatrix} \begin{Bmatrix} \Delta u = 0 \\ u|_{\partial K} = \mathbb{1}_{f_i} \end{Bmatrix} + \begin{pmatrix} n_{ix} \\ n_{iy} \end{pmatrix} \begin{Bmatrix} \Delta u = 0 \\ u|_{\partial K} \equiv 2 \end{Bmatrix}$$

**General case** $k > 0$**:**

*Normal functions:*

$$f_j := (x,\, y) \mapsto \begin{pmatrix} x \\ y \end{pmatrix} \begin{Bmatrix} \Delta u = h_{k-1,\,k-1}(x, y) \\ u|_{\partial K} = g_j(x, y)\mathbb{1}_{f_i} \end{Bmatrix} + \begin{pmatrix} n_{ix} \\ n_{iy} \end{pmatrix} \begin{Bmatrix} \Delta u = 0 \\ u|_{\partial K} \equiv 2 \end{Bmatrix},\; j \in [\![0,\, k]\!]$$

*Internal functions:*

$$f_{xl} := (x,\, y) \mapsto \begin{pmatrix} x \\ y \end{pmatrix} \begin{Bmatrix} \Delta u = h_{k-1,\,l}(x,\, y) \\ u|_{\partial K} = 0 \end{Bmatrix},\; 0 \le l \le k-1$$

$$f_{yl} := (x,\, y) \mapsto \begin{pmatrix} x \\ y \end{pmatrix} \begin{Bmatrix} \Delta u = h_{l,\,k-1}(x,\, y) \\ u|_{\partial K} = 0 \end{Bmatrix},\; 0 \le l < k-1$$

$$f_{xlm} := (x,\, y) \mapsto \begin{pmatrix} 1 \\ 0 \end{pmatrix} \begin{Bmatrix} \Delta u = h_{l,\,m}(x,\, y) \\ u|_{\partial K} = 0 \end{Bmatrix},\; 0 \le l,\, m \le k-1$$

$$f_{ylm} := (x,\, y) \mapsto \begin{pmatrix} 0 \\ 1 \end{pmatrix} \begin{Bmatrix} \Delta u = h_{l,\,m}(x,\, y) \\ u|_{\partial K} = 0 \end{Bmatrix},\; 0 \le l,\, m \le k-1$$

### A.3.2 Implemented Elements (degrees of freedom)

*Normal moments*, for any edge $f_j \in \partial K$:

| | Core Moments | Misc Moment | Supplementary moments |
|---|---|---|---|
| Ia | $\int_{f_j} p_i(v)\, \varphi \cdot n \,\mathrm{d}\gamma,\ i \in [\![1,\, k]\!]$ | $\int_{f_j} v\, \varphi \cdot \begin{pmatrix} 1 \\ 1 \end{pmatrix} \,\mathrm{d}\gamma$ | None |
| Ib | $\int_{f_j} p_i(v)\, \varphi \cdot n \,\mathrm{d}\gamma,\ i \in [\![1,\, k]\!]$ | $\int_{f_j} v\, \varphi \cdot \begin{pmatrix} 1 \\ 1 \end{pmatrix} \,\mathrm{d}\gamma$ | None |
| IbShifted | $\int_{f_j} p_i(v)\, \varphi \cdot n \,\mathrm{d}\gamma,\ i \in [\![1,\, k]\!]$ | $\int_{f_j} v\, \varphi \cdot \begin{pmatrix} 1 \\ 1 \end{pmatrix} \,\mathrm{d}\gamma$ | None |
| IIa | $\int_{f_j} p_i(v)\, \varphi \cdot n \,\mathrm{d}\gamma,\ i \in [\![1,\, k]\!]$ | $\int_{f_j} \varphi \cdot n \,\mathrm{d}\gamma$ | None |
| IIb | $\int_{f_j} p_i(v)\, \varphi \cdot n \,\mathrm{d}\gamma,\ i \in [\![1,\, k]\!]$ | $\varphi(x_m) \cdot n$ | None |
| IIbShifted | $\int_{f_j} p_i(v)\, \varphi \cdot n \,\mathrm{d}\gamma,\ i \in [\![1,\, k]\!]$ | $\varphi(x_m) \cdot n - 1$ | None |

*Internal moments:*

| | Core Moments | Misc Moment |
|---|---|---|
| All elements | $\int_K \begin{pmatrix} q_{l,m}(x,\, y) \\ 0 \end{pmatrix} \cdot \varphi \,\mathrm{d}x\mathrm{d}y$ and $\int_K \begin{pmatrix} 0 \\ q_{m,l}(x,\, y) \end{pmatrix} \cdot \varphi \,\mathrm{d}x\mathrm{d}y,\ {\scriptstyle l \in [\![0,\, k]\!],\ m \in [\![0,\, k-1]\!],\ (l,\, m) \neq (k,\, k-1)}$ | $\int_K \begin{pmatrix} q_{k,\, k-1}(x,\, y) \\ q_{k-1,\, k}(x,\, y) \end{pmatrix} \cdot \varphi \,\mathrm{d}x\mathrm{d}y$ |

***Note.*** Due to the selected moments to leave out from the general setting (as explained in the *Section 6.2*), the elements Ia, Ib and IbShifted are Identical.

▲

## A.4 Implemented functions

### A.4.1 Functions used in the Poisson's problems

#### A.4.1.1 Boundary functions for Poisson's problems (g)

The functions $g_j := (x, y) \mapsto g_j(x, y)$ are designed to furnish a 1D function that is defined along the considered boundary. Here, we only propose polynomial functions. The term $(x, y)$ represents the running point along this boundary and $j$ the polynomial order. The proposed boundary functions are:

- 1: **Lagrangian function (recommended).**

  1D $j$th Lagrangian function (of degree $k$) defined from the $x$ coordinate only (one can see it as a projection). It is possible to change it with respect to the norm of the vector line 180-181 of the file Problem.py, without much impact.

- 2: **Canonical polynomial (centred, scaled).**

  Scaled and centred unscaled 1D canonical polynomial, determined depending on the running norm along the edge. Practically, it computes: $l = length(edge); u = 2(x-xa)/l; w = 2(y-ya)/l; z = norm((u, w), 2) - 1; result = z^j$, where $(xa, ya)$ are the coordinates of one of the two vertices attached to the edge.

- 3: **Canonical polynomial (centred, unscaled):**

  Centred unscaled 1D canonical polynomial, determined depending on the running norm along the edge. Practically, it computes:$u = x - xa; w = y - ya; z = norm((u, w), 2) - 0.5length(edge); result = z^j$, where $(xa, ya)$ are the coordinates of one of the two vertices attached to the edge.

***Note.*** The function type is assumed to be chosen identically for all the element's boundaries.

▲

#### A.4.1.2 Second member functions for Poisson's problems (h):

The functions $h_{i,j} := (x, y) \mapsto h_{i,j}(x, y)$ are designed to furnish a 2D function that is defined within the considered polygon. Here, we only consider polynomial cases. Thus, $(x, y)$ represent the evaluation point, $i$ the polynomial order on the $x$ coordinate and $j$ the polynomial order on the $y$ coordinate. The proposed functions are:

- 1: **Chebyshev polynomial (scaled).**

  Chebyshev polynomial evaluated at a rescaled and centred point with respect to the element's shape. $result = Chebyshev(2(x - bary_x)/Area, i)chebval(2(y - bary_y)/Area, j)$.

- 2: **Hermite polynomial (scaled).**

  Hermite polynomial evaluated at a rescaled and centred point with respect to the element's shape, with coefficient one on the term involving the $i$th power of $x$ and $j$th power of $y$: $result = Hermite(4(x - bary_x)/Area, 4(y - bary_y)/Area, (i, j))$.

- 3: **Legendre polynomial (scaled).**

  Legendre polynomial evaluated at a rescaled and centred point with respect to the element's shape, with coefficient one on the term involving the $i$th power of $x$ and $j$th power of $y$: $result = Legendre(2(x - bary_x)/Area, 2(y - bary_y)/Area, (i, j))$.

- 4: **Canonical polynomial (centred, scaled).**

  Returns the usual canonically evaluated polynomial, centred w.r.t. the barycentre of the convex hull of the element's shape and scale with respect to the area of the convex hull: $result = (2(x - bary_x)/Area)^i(2(y - bary_y)/Area)^j$.

- 5: **Canonical polynomial (centred, unscaled).**

  Returns the usual canonically evaluated polynomial, centred w.r.t. the barycentre of the convex hull of the element's shape: $result = (x - bary_x)^i(y - bary_y)^j$.

***Note.*** As before, the function type is assumed to be chosen identically for all the element's boundaries.

▲

### A.4.2 Projection functions

#### A.4.2.1 Boundary projectors as integrand functions (p):

The functions $g_i := v \mapsto g_i(v)$ are designed to furnish a 1D function that is defined along the considered boundary. Here, $v$ represents the running norm along this boundary. The proposed boundary functions are all polynomials. Thus, $i$ represent the degree of the wished polynomial, being of the following types:

- 1: **Canonical polynomial (centred, scaled).**

  Scaled and centred 1D canonical polynomial, determined depending on the running norm along the edge. Practically, it computes: $z = 2v/length(edge) - 1; result = z^j$.

- 2: **Chebyshev polynomial (scaled).**

  Chebyshev polynomial evaluated at a rescaled and centred running point with respect to the edge. Practically, it computes: $z = 2v/length(edge) - 1; result = Chebyshev(z, i)$.

- 3: **Hermite polynomial (scaled).**

  Hermite polynomial evaluated at a rescaled and centred running point with respect to the edge. Practically, it computes: $z = 4v/length(edge) - 2; result = Hermite(z, i)$.

- 4: **Legendre polynomial (scaled).**

  Legendre polynomial evaluated at a rescaled and centred running point with respect to the edge. Practically, it computes: $z = 2v/length(edge) - 1; result = Legendre(z, i)$.

- 5: **Laguerre Polynomial (scaled).**

  Legendre polynomial evaluated at a rescaled and centred running point with respect to the edge. Practically, it computes: $z = 12v/length(edge) - 2; result = Laguerre(z, i)$.

- 6: **Canonical polynomial (centred, unscaled).**

  Centred unscaled 1D canonical polynomial, determined depending on the running norm along the edge. Practically, it computes: $z = v - 0.5 * length(edge); result = v^j$.

- 7: **Canonical polynomial (not - centred, unscaled).**

  Not-centred unscaled 1D canonical polynomial, determined depending on the running norm along the edge. Practically, it computes: $z = v; result = v^j$.

#### A.4.2.2 Inner cell projectors as integrand functions (q):

The functions are defined in the exact same way as for the Harmonic problems. We refer to the above section for the description. However, the indexing is slightly different. Here are available:

- 1: **Canonical polynomial (centred, scaled).**
- 2: **Chebyshev polynomial (scaled).**
- 3: **Hermite polynomial (scaled).**
- 4: **Legendre polynomial (scaled).**
- 5: **Laguerre Polynomial (scaled).**

  Laguerre polynomial evaluated at a rescaled and centred point with respect to the element's shape, with coefficient one on the term involving the $i$th power of $x$ and $j$th power of $y$: $result = Laaguerre(12(x - bary_x + 4)/Area, 12(y - bary_y + 4)/Area, (i, j))$.
- 6: **Canonical polynomial (centred, unscaled).**
- 7: **Canonical polynomial (not - centred, unscaled).**

  Returns the usual canonically evaluated polynomial at the point whose coordinates are taken as absolute, without taking into account the size or location of the element with respect to the axes: $result = (x)^i(y)^j$.

# B Details of some elements used in the numerical experiments

We collect here the informations corresponding to the definition of the elements used in the numerical results section.

## B.1 Elements for a comparison with the Raviart – Thomas setting

| Edge | Vertex1 | Vertex2 | Normal | Norm |
|---|---|---|---|---|
| 0 | [0.20, 0.00] | [1.00, 0.20] | [0.24, -0.97] | [0.82] |
| 1 | [1.00, 0.20] | [0.00, 1.00] | [0.62, 0.78] | [1.28] |
| 2 | [0.00, 1.00] | [0.20, 0.00] | [-0.98, -0.20] | [1.02] |

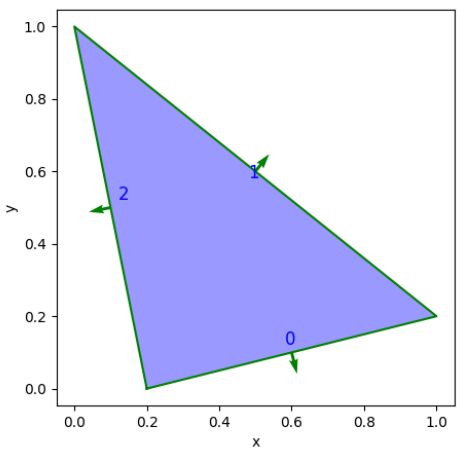

Fig. 151: Triangular shape

| Edge | Vertex1 | Vertex2 | Normal | Norm |
|---|---|---|---|---|
| 0 | [0.20, 0.00] | [1.00, 0.20] | [0.24, -0.97] | [0.82] |
| 1 | [1.00, 0.20] | [0.80, 1.20] | [0.98, 0.20] | [1.02] |
| 2 | [0.80, 1.20] | [0.00, 1.00] | [-0.24, 0.97] | [0.82] |
| 3 | [0.00, 1.00] | [0.20, 0.00] | [-0.98, -0.20] | [1.02] |

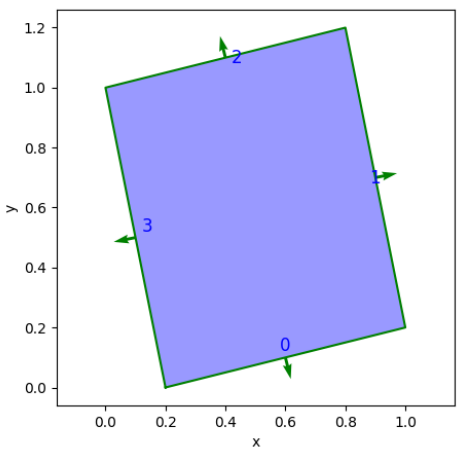

Fig. 152: Quadrangular shape

## B.2 Convex elements

| Edge | Vertex1 | Vertex2 | Normal | Norm |
|---|---|---|---|---|
| 0 | [0.07, 0.18] | [0.41, 0.05] | [-0.34, -0.94] | [0.37] |
| 1 | [0.41, 0.05] | [0.36, 0.41] | [0.99, 0.15] | [0.36] |
| 2 | [0.36, 0.41] | [0.07, 0.18] | [-0.63, 0.78] | [0.37] |

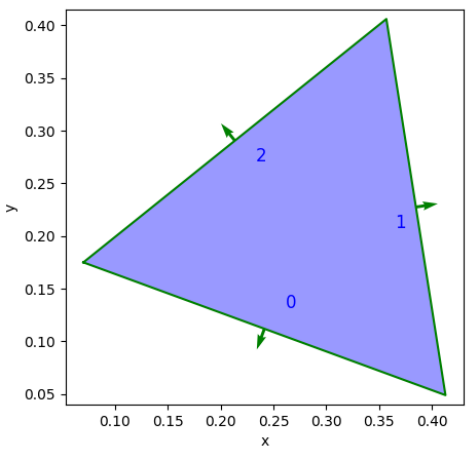

Fig. 153: Smallest considered triangular shape

| Edge | Vertex1 | Vertex2 | Normal | Norm |
|---|---|---|---|---|
| 0 | [0.15, 0.38] | [0.89, 0.10] | [-0.34, -0.94] | [0.78] |
| 1 | [0.89, 0.10] | [0.76, 0.87] | [0.99, 0.15] | [0.77] |
| 2 | [0.76, 0.87] | [0.15, 0.38] | [-0.63, 0.78] | [0.79] |

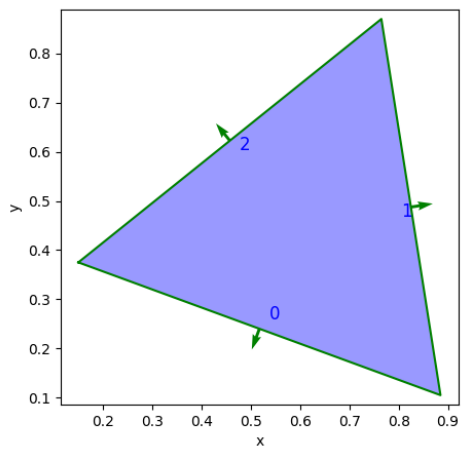

Fig. 154: Small considered triangular shape

| Edge | Vertex1 | Vertex2 | Normal | Norm |
|---|---|---|---|---|
| 0 | [0.30, 0.75] | [1.77, 0.21] | [-0.34, -0.94] | [1.57] |
| 1 | [1.77, 0.21] | [1.53, 1.74] | [0.99, 0.15] | [1.55] |
| 2 | [1.53, 1.74] | [0.30, 0.75] | [-0.63, 0.78] | [1.58] |

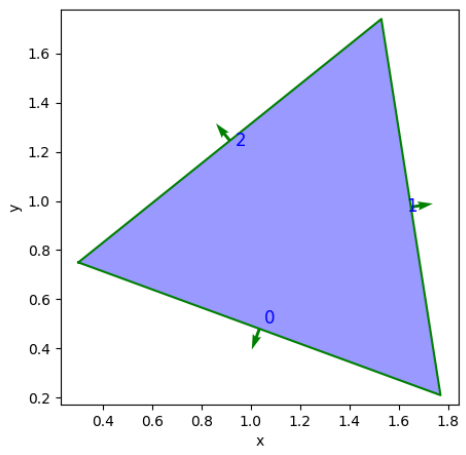

Fig. 155: Considered triangular shape

| Edge | Vertex1 | Vertex2 | Normal | Norm |
|---|---|---|---|---|
| 0 | [0.60, 1.50] | [3.54, 0.42] | [-0.34, -0.94] | [3.13] |
| 1 | [3.54, 0.42] | [3.06, 3.48] | [0.99, 0.15] | [3.10] |
| 2 | [3.06, 3.48] | [0.60, 1.50] | [-0.63, 0.78] | [3.16] |

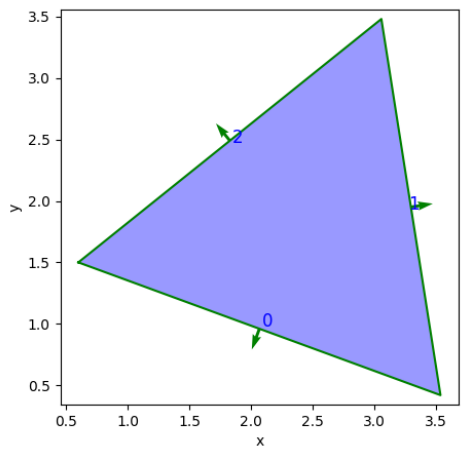

Fig. 156: Big considered triangular shape

| Edge | Vertex1 | Vertex2 | Normal | Norm |
|---|---|---|---|---|
| 0 | [0.90, 2.25] | [5.31, 0.63] | [-0.34, -0.94] | [4.70] |
| 1 | [5.31, 0.63] | [4.59, 5.22] | [0.99, 0.15] | [4.65] |
| 2 | [4.59, 5.22] | [0.90, 2.25] | [-0.63, 0.78] | [4.74] |

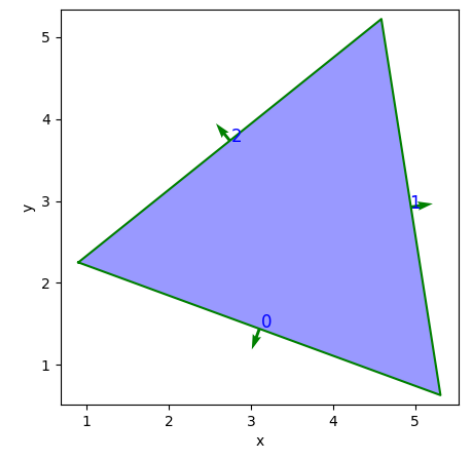

Fig. 157: Biggest considered triangular shape

| Edge | Vertex1 | Vertex2 | Normal | Norm |
|---|---|---|---|---|
| 0 | [0.25, 0.00] | [0.50, 0.25] | [0.71, -0.71] | [0.35] |
| 1 | [0.50, 0.25] | [0.25, 0.50] | [0.71, 0.71] | [0.35] |
| 2 | [0.25, 0.50] | [0.00, 0.25] | [-0.71, 0.71] | [0.35] |
| 3 | [0.00, 0.25] | [0.25, 0.00] | [-0.71, -0.71] | [0.35] |

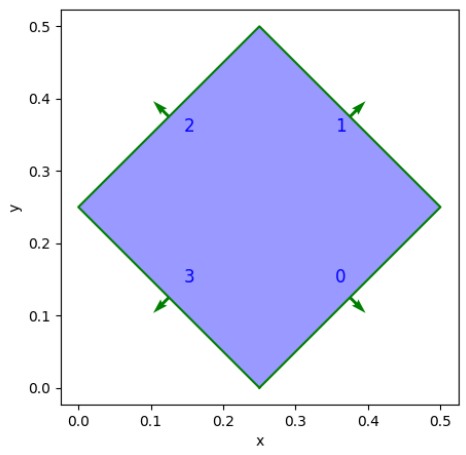

Fig. 158: Quadrangular shape with parallel edges

| Edge | Vertex1 | Vertex2 | Normal | Norm |
|---|---|---|---|---|
| 0 | [0.08, 0.07] | [0.33, 0.02] | [-0.19, -0.98] | [0.26] |
| 1 | [0.33, 0.02] | [0.48, 0.23] | [0.80, -0.60] | [0.26] |
| 2 | [0.48, 0.23] | [0.28, 0.39] | [0.62, 0.79] | [0.26] |
| 3 | [0.28, 0.39] | [0.03, 0.33] | [-0.23, 0.97] | [0.26] |
| 4 | [0.03, 0.33] | [0.08, 0.07] | [-0.98, -0.19] | [0.26] |

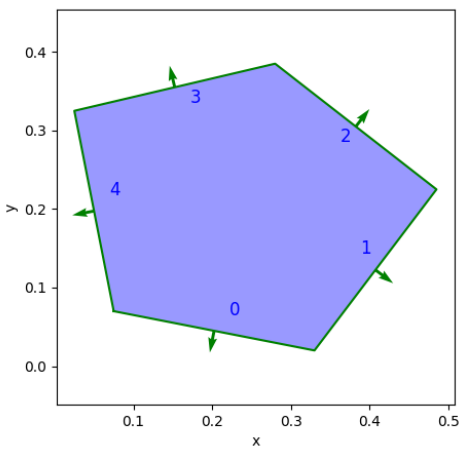

Fig. 159: Pentagonal shape

| Edge | Vertex1 | Vertex2 | Normal | Norm |
|---|---|---|---|---|
| 0 | [2.05, 1.98] | [2.31, 1.92] | [-0.25, -0.97] | [0.26] |
| 1 | [2.31, 1.92] | [2.35, 2.17] | [0.99, -0.17] | [0.26] |
| 2 | [2.35, 2.17] | [2.31, 2.43] | [0.99, 0.15] | [0.26] |
| 3 | [2.31, 2.43] | [2.08, 2.32] | [-0.42, 0.91] | [0.26] |
| 4 | [2.08, 2.32] | [1.86, 2.16] | [-0.61, 0.80] | [0.26] |
| 5 | [1.86, 2.16] | [2.05, 1.98] | [-0.70, -0.72] | [0.26] |

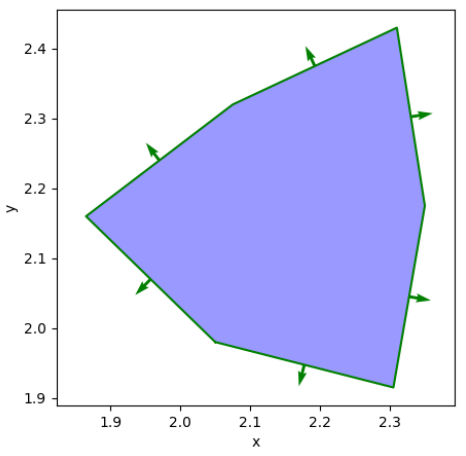

Fig. 160: Regular hexagonal shape

| Edge | Vertex1 | Vertex2 | Normal | Norm |
|---|---|---|---|---|
| 0 | [0.07, 0.07] | [0.16, 0.02] | [-0.43, -0.90] | [0.11] |
| 1 | [0.16, 0.02] | [0.42, 0.11] | [0.31, -0.95] | [0.27] |
| 2 | [0.42, 0.11] | [0.34, 0.21] | [0.80, 0.61] | [0.13] |
| 3 | [0.34, 0.21] | [0.20, 0.33] | [0.65, 0.76] | [0.18] |
| 4 | [0.20, 0.33] | [0.04, 0.23] | [-0.55, 0.84] | [0.19] |
| 5 | [0.04, 0.23] | [0.07, 0.07] | [-0.99, -0.15] | [0.16] |

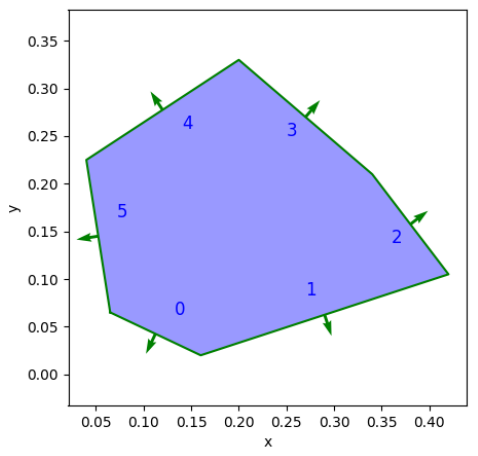

Fig. 161: Hexagonal shape

| Edge | Vertex1 | Vertex2 | Normal | Norm |
|---|---|---|---|---|
| 0 | [0.13, 0.13] | [0.32, 0.04] | [-0.43, -0.90] | [0.21] |
| 1 | [0.32, 0.04] | [0.84, 0.21] | [0.31, -0.95] | [0.55] |
| 2 | [0.84, 0.21] | [0.78, 0.38] | [0.94, 0.33] | [0.18] |
| 3 | [0.78, 0.38] | [0.50, 0.66] | [0.71, 0.71] | [0.40] |
| 4 | [0.50, 0.66] | [0.08, 0.45] | [-0.45, 0.89] | [0.47] |
| 5 | [0.08, 0.45] | [0.13, 0.13] | [-0.99, -0.15] | [0.32] |

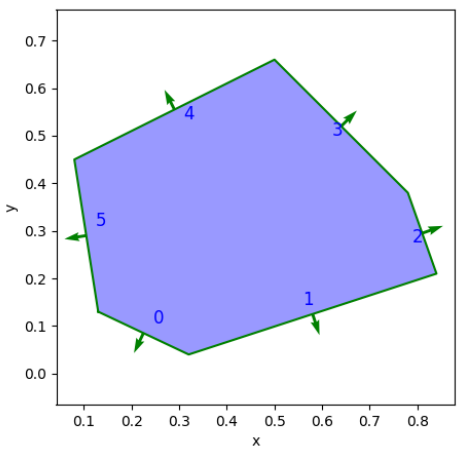

Fig. 162: Alternative hexagonal shape

## B.3 Non - convex elements

| Edge | Vertex1 | Vertex2 | Normal | Norm |
|---|---|---|---|---|
| 0 | [0.14, 0.03] | [0.37, 0.10] | [0.29, -0.96] | [0.24] |
| 1 | [0.37, 0.10] | [0.21, 0.14] | [0.21, 0.98] | [0.16] |
| 2 | [0.21, 0.14] | [0.06, 0.28] | [0.67, 0.74] | [0.21] |
| 3 | [0.06, 0.28] | [0.14, 0.03] | [-0.94, -0.33] | [0.26] |

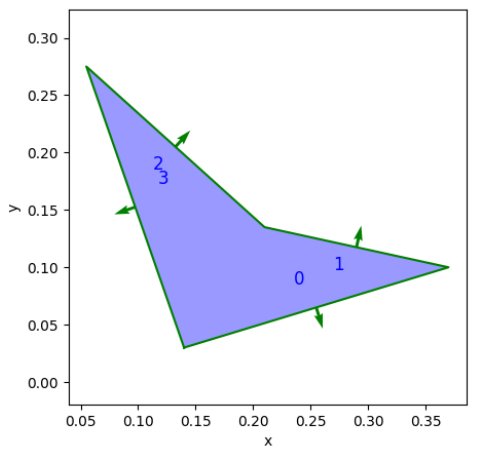

Fig. 163: Non - convex quadrilateral

| Edge | Vertex1 | Vertex2 | Normal | Norm |
|---|---|---|---|---|
| 0 | [0.17, 0.03] | [0.38, 0.19] | [0.59, -0.80] | [0.26] |
| 1 | [0.38, 0.19] | [0.30, 0.18] | [-0.18, 0.98] | [0.08] |
| 2 | [0.30, 0.18] | [0.12, 0.36] | [0.74, 0.68] | [0.26] |
| 3 | [0.12, 0.36] | [0.19, 0.11] | [-0.97, -0.25] | [0.26] |
| 4 | [0.19, 0.11] | [0.17, 0.03] | [-0.97, 0.24] | [0.08] |

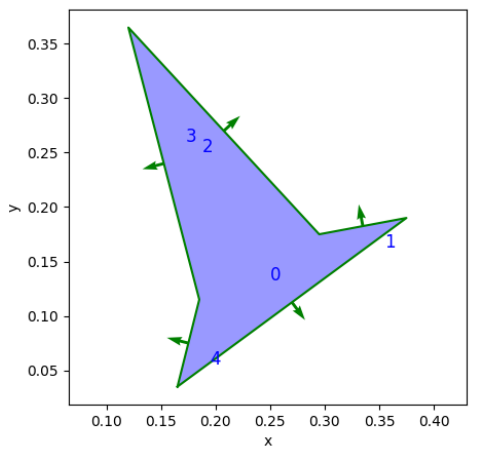

Fig. 164: Non - convex pentagon

| Edge | Vertex1 | Vertex2 | Normal | Norm |
|---|---|---|---|---|
| 0 | [0.00, 0.03] | [0.12, 0.07] | [0.37, -0.93] | [0.13] |
| 1 | [0.12, 0.07] | [0.38, 0.00] | [-0.29, -0.96] | [0.26] |
| 2 | [0.38, 0.00] | [0.30, 0.25] | [0.96, 0.29] | [0.26] |
| 3 | [0.30, 0.25] | [0.12, 0.38] | [0.58, 0.81] | [0.22] |
| 4 | [0.12, 0.38] | [-0.12, 0.25] | [-0.45, 0.89] | [0.28] |
| 5 | [-0.12, 0.25] | [0.00, 0.03] | [-0.87, -0.49] | [0.26] |

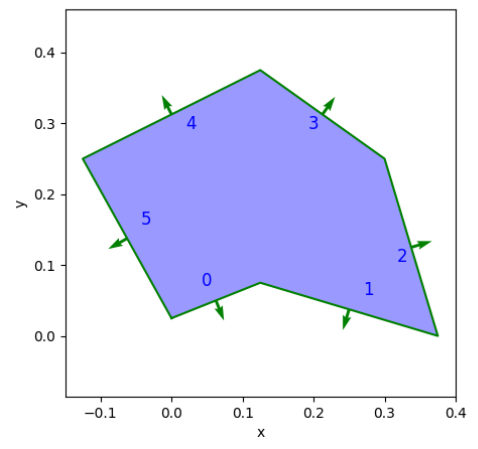

Fig. 165: Non - convex hexagon

| Edge | Vertex1 | Vertex2 | Normal | Norm |
|---|---|---|---|---|
| 0 | [0.07, 0.03] | [0.35, 0.10] | [0.24, -0.97] | [0.29] |
| 1 | [0.35, 0.10] | [0.45, 0.25] | [0.83, -0.55] | [0.18] |
| 2 | [0.45, 0.25] | [0.25, 0.30] | [0.24, 0.97] | [0.21] |
| 3 | [0.25, 0.30] | [0.05, 0.25] | [-0.24, 0.97] | [0.21] |
| 4 | [0.05, 0.25] | [0.14, 0.16] | [-0.73, -0.69] | [0.12] |
| 5 | [0.14, 0.16] | [0.07, 0.03] | [-0.88, 0.47] | [0.15] |

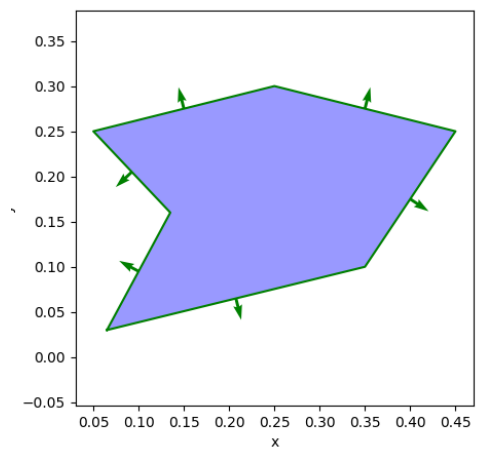

Fig. 166: Alternative non - convex hexagon

| Edge | Vertex1 | Vertex2 | Normal | Norm |
|---|---|---|---|---|
| 0 | [0.22, 0.10] | [0.50, 0.30] | [0.58, -0.81] | [0.34] |
| 1 | [0.50, 0.30] | [0.76, 0.10] | [-0.61, -0.79] | [0.33] |
| 2 | [0.76, 0.10] | [0.70, 0.40] | [0.98, 0.20] | [0.31] |
| 3 | [0.70, 0.40] | [0.90, 0.70] | [0.83, -0.55] | [0.36] |
| 4 | [0.90, 0.70] | [0.60, 0.62] | [-0.26, 0.97] | [0.31] |
| 5 | [0.60, 0.62] | [0.50, 0.90] | [0.94, 0.34] | [0.30] |
| 6 | [0.50, 0.90] | [0.35, 0.68] | [-0.83, 0.56] | [0.27] |
| 7 | [0.35, 0.68] | [0.12, 0.50] | [-0.62, 0.79] | [0.29] |
| 8 | [0.12, 0.50] | [0.30, 0.35] | [-0.64, -0.77] | [0.23] |
| 9 | [0.30, 0.35] | [0.22, 0.10] | [-0.95, 0.30] | [0.26] |

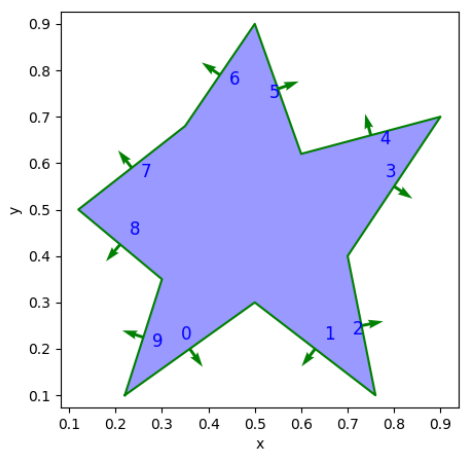

Fig. 167: Non - convex decagon

## B.4 Elements used to test the failing and limit cases

| Edge | Vertex1 | Vertex2 | Normal | Norm |
|---|---|---|---|---|
| 0 | [0.20, 0.00] | [1.00, 0.20] | [0.24, -0.97] | [0.82] |
| 1 | [1.00, 0.20] | [1.60, 1.40] | [0.89, -0.45] | [1.34] |
| 2 | [1.60, 1.40] | [0.80, 1.20] | [-0.24, 0.97] | [0.82] |
| 3 | [0.80, 1.20] | [0.00, 1.00] | [-0.24, 0.97] | [0.82] |
| 4 | [0.00, 1.00] | [0.20, 0.00] | [-0.98, -0.20] | [1.02] |

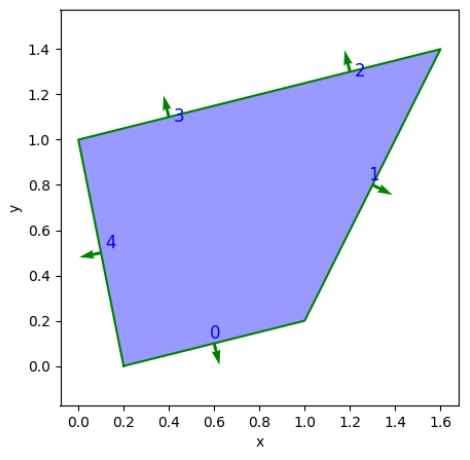

Fig. 168: Hanging node

| Edge | Vertex1 | Vertex2 | Normal | Norm |
|---|---|---|---|---|
| 0 | [0.20, 0.00] | [1.00, 0.20] | [0.24, -0.97] | [0.82] |
| 1 | [1.00, 0.20] | [3.50, 0.00] | [-0.08, -1.00] | [2.51] |
| 2 | [3.50, 0.00] | [2.40, 1.60] | [0.82, 0.57] | [1.94] |
| 3 | [2.40, 1.60] | [1.60, 1.40] | [-0.24, 0.97] | [0.82] |
| 4 | [1.60, 1.40] | [1.70, 0.90] | [-0.98, -0.20] | [0.51] |
| 5 | [1.70, 0.90] | [0.90, 0.70] | [-0.24, 0.97] | [0.82] |
| 6 | [0.90, 0.70] | [0.80, 1.20] | [0.98, 0.20] | [0.51] |
| 7 | [0.80, 1.20] | [0.00, 1.00] | [-0.24, 0.97] | [0.82] |
| 8 | [0.00, 1.00] | [0.20, 0.00] | [-0.98, -0.20] | [1.02] |

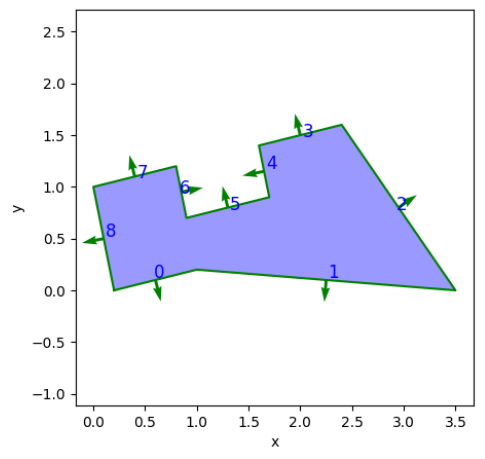

Fig. 169: Polygon having similar edges

| Edge | Vertex1 | Vertex2 | Normal | Norm |
|---|---|---|---|---|
| 0 | [1.15, 0.10] | [1.88, 0.38] | [0.35, -0.93] | [0.78] |
| 1 | [1.88, 0.38] | [1.45, 0.88] | [0.76, 0.65] | [0.66] |
| 2 | [1.45, 0.88] | [0.75, 0.75] | [-0.18, 0.98] | [0.71] |
| 3 | [0.75, 0.75] | [0.25, 0.25] | [-0.71, 0.71] | [0.71] |
| 4 | [0.25, 0.25] | [1.15, 0.10] | [-0.16, -0.99] | [0.91] |

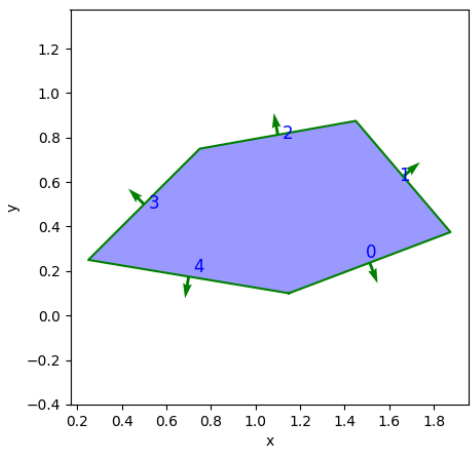

Fig. 170: Two vertices aligned with the origin

| Edge | Vertex1 | Vertex2 | Normal | Norm |
|---|---|---|---|---|
| 0 | [0.00, 0.00] | [1.25, 0.38] | [0.29, -0.96] | [1.31] |
| 1 | [1.25, 0.38] | [1.32, 0.90] | [0.99, -0.14] | [0.53] |
| 2 | [1.32, 0.90] | [1.00, 1.25] | [0.73, 0.68] | [0.48] |
| 3 | [1.00, 1.25] | [0.42, 1.23] | [-0.04, 1.00] | [0.58] |
| 4 | [0.42, 1.23] | [0.00, 0.00] | [-0.94, 0.33] | [1.30] |

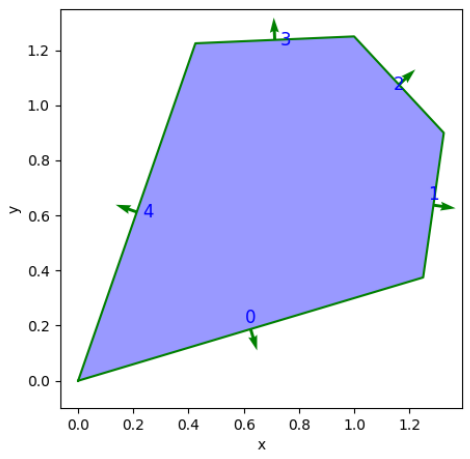

Fig. 171: Two vertices aligned with the origin

| Edge | Vertex1 | Vertex2 | Normal | Norm |
|---|---|---|---|---|
| 0 | [0.05, 0.05] | [0.25, 0.05] | [0.00, -1.00] | [0.20] |
| 1 | [0.25, 0.05] | [0.32, 0.16] | [0.84, -0.54] | [0.13] |
| 2 | [0.32, 0.16] | [0.19, 0.29] | [0.71, 0.71] | [0.18] |
| 3 | [0.19, 0.29] | [0.09, 0.23] | [-0.53, 0.85] | [0.12] |
| 4 | [0.09, 0.23] | [0.05, 0.05] | [-0.98, 0.20] | [0.18] |

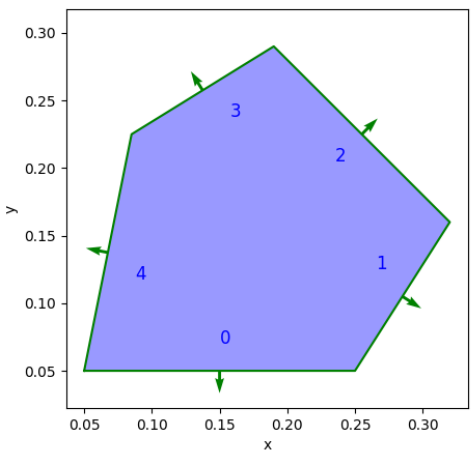

Fig. 172: Edge parallel to the $x$ axis

| Edge | Vertex1 | Vertex2 | Normal | Norm |
|---|---|---|---|---|
| 0 | [0.05, 0.00] | [0.25, 0.07] | [0.31, -0.95] | [0.21] |
| 1 | [0.25, 0.07] | [0.32, 0.16] | [0.81, -0.59] | [0.12] |
| 2 | [0.32, 0.16] | [0.19, 0.29] | [0.71, 0.71] | [0.18] |
| 3 | [0.19, 0.29] | [0.05, 0.25] | [-0.27, 0.96] | [0.15] |
| 4 | [0.05, 0.25] | [0.05, 0.00] | [-1.00, -0.00] | [0.25] |

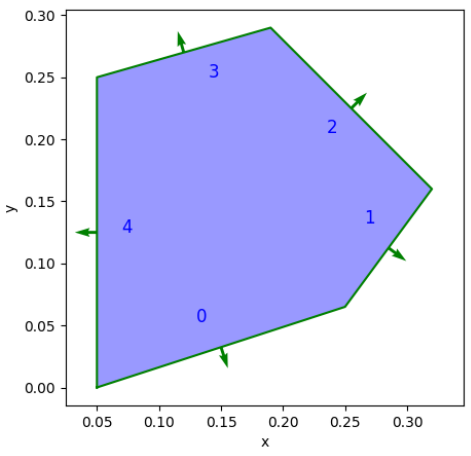

Fig. 173: Edge parallel to the $y$ axis

## References

[1] R. Abgrall, E. le Mélédo, and P. Öffner. On the Connection between Residual Distribution Schemes and Flux Reconstruction. `https://arxiv.org/abs/1807.01261`, July 2018.

[2] D. Boffi, F. Brezzi, M. Fortin, et al. *Mixed finite element methods and applications*, volume 44. Springer, 2013.

[3] F. Brezzi, J. Douglas, and L. D. Marini. Two families of mixed finite elements for second order elliptic problems. *Numerische Mathematik*, 47(2):217–235, Jun 1985.

[4] F. Brezzi and M. Fortin, editors. *Mixed and Hybrid Finite Element Methods*. Springer New York, 1991.

[5] J. Chabrowski. *The Dirichlet problem with L2-boundary data for elliptic linear equations*. Springer, 2006.

[6] B. Cockburn, G. E. Karniadakis, and C.-W. Shu. *Discontinuous Galerkin methods: theory, computation and applications*, volume 11. Springer Science & Business Media, 2012.

[7] D. A. Di Pietro and S. Lemaire. An extension of the Crouzeix-Raviart space to general meshes with application to quasi-incompressible linear elasticity and Stokes flow. *Mathematics of Computation*, 84(291):1–31, 2015.

[8] F. Dubois, I. Greff, and C. Pierre. Raviart-thomas finite elements of petrov-galerkin type. *arXiv preprint arXiv:1710.04395*, 2017.

[9] V. Ervin. Computational bases for rtk and bdmk on triangles. *Computers & Mathematics with Applications*, 64(8):2765–2774, 2012.

[10] A. Gillette, A. Rand, and C. Bajaj. Construction of scalar and vector finite element families on polygonal and polyhedral meshes.

[11] J. Gopalakrishnan, L. E. García-Castillo, and L. F. Demkowicz. Nédélec spaces in affine coordinates. *Computers & Mathematics with Applications*, 49(7-8):1285–1294, 2005.

[12] H. T. Huynh. A flux reconstruction approach to high-order schemes including discontinuous galerkin methods. In *18th AIAA Computational Fluid Dynamics Conference*, page 4079, 2007.

[13] J.C. Nédélec. Mixed finite elements in $\mathbb{R}^3$. *Numerische Mathematik*, 35(3):315–341, 1980.

[14] M. V. Kondratieva, A. B. Levin, A. V. Mikhalev, and E. V. Pankratiev. *Differential and Difference Dimension Polynomials*. Springer Netherlands, 1999.

[15] H. Ranocha, P. Öffner, and T. Sonar. Summation-by-parts operators for correction procedure via reconstruction. *Journal of Computational Physics*, 311:299–328, 2016.

[16] P. A. Raviart and J. M. Thomas. A mixed finite element method for 2-nd order elliptic problems. In I. Galligani and E. Magenes, editors, *Mathematical Aspects of Finite Element Methods*, pages 292–315, Berlin, Heidelberg, 1977. Springer Berlin Heidelberg.

[17] N. Sukumar and E. Malsch. Recent advances in the construction of polygonal finite element interpolants. *Archives of Computational Methods in Engineering*, 13(1):129, 2006.

[18] C. Talischi, A. Pereira, G. H. Paulino, I. F. Menezes, and M. S. Carvalho. Polygonal finite elements for incompressible fluid flow. *International Journal for Numerical Methods in Fluids*, 74(2):134–151, oct 2013.

[19] J.-M. Thomas. Méthode des éléments finis hybrides duaux pour les problèmes elliptiques du second ordre. *Revue française d'automatique, informatique, recherche opérationnelle. Analyse numérique*, 10(R3):51–79, 1976.

[20] L. B. D. Veiga, F. Brezzi, A. Cangiani, G. Manzini, L. D. Marini, and A. Russo. Basic principles of virtual element methods. *Mathematical Models and Methods in Applied Sciences*, 23(01):199–214, jan 2013.

[21] P. E. Vincent, P. Castonguay, and A. Jameson. A new class of high-order energy stable flux reconstruction schemes. *Journal of Scientific Computing*, 47(1):50–72, 2011.